\documentclass[12pt, oneside]{book}

\usepackage[sort,numbers]{natbib}

\usepackage[svgnames]{xcolor}

\usepackage{amsmath,amsfonts,amsthm,tikz,comment,xfrac,amssymb}
\usetikzlibrary{decorations.markings}
\usepackage[T1]{fontenc}
\usepackage[font={scriptsize}]{caption}

\usepackage[utf8]{inputenc}
\usepackage[T1]{fontenc}
\usepackage{lmodern}          
\usepackage{microtype}        
\usepackage[
            a4paper,            
            includeheadfoot,    
            textwidth=16cm,
            hmarginratio=1:1,   
            top=2.5cm,          
            bottom=2.7cm]{geometry}  

\definecolor{ThmColor}{rgb}{1.00, 0.95, 0.75}    
\definecolor{CorColor}{rgb}{0.93, 0.86, 0.72}    
\definecolor{PropColor}{rgb}{1.00, 0.88, 0.88}   
\definecolor{LemColor}{rgb}{0.88, 0.96, 0.88}    
\definecolor{DefColor}{rgb}{0.85, 0.94, 0.980}   
\definecolor{NotColor}{rgb}{1.00, 0.94, 0.86}    
\definecolor{RemColor}{rgb}{0.95, 0.90, 0.950}   
\definecolor{ExoColor}{rgb}{0.89, 0.993, 0.890}  
\definecolor{AssColor}{rgb}{0.72, 0.92, 0.90}    

\usepackage{amsmath}
\usepackage{amsfonts}
\usepackage{amssymb}
\usepackage{amsthm}
\usepackage{mathtools}

\usepackage{mdframed}

\mdfdefinestyle{thmstyle}{backgroundcolor=ThmColor, nobreak=true,
    linewidth=0.5pt, linecolor=Goldenrod!60}
\mdfdefinestyle{corstyle}{backgroundcolor=CorColor, nobreak=true,
    linewidth=0.5pt, linecolor=Brown!50}
\mdfdefinestyle{propstyle}{backgroundcolor=PropColor, nobreak=true,
    linewidth=0.5pt, linecolor=RosyBrown!50}
\mdfdefinestyle{lemstyle}{backgroundcolor=LemColor, nobreak=true,
    linewidth=0.5pt, linecolor=SeaGreen!40}
\mdfdefinestyle{defstyle}{backgroundcolor=DefColor, nobreak=true,
    linewidth=0.5pt, linecolor=SteelBlue!50}
\mdfdefinestyle{notstyle}{backgroundcolor=NotColor, nobreak=true,
    linewidth=0.5pt, linecolor=Orange!40}
\mdfdefinestyle{remstyle}{backgroundcolor=RemColor,
    linewidth=0.5pt, linecolor=MediumOrchid!40}
\mdfdefinestyle{exostyle}{backgroundcolor=ExoColor,
    linewidth=0.5pt, linecolor=SeaGreen!50}
\mdfdefinestyle{assstyle}{backgroundcolor=AssColor, nobreak=true,
    linewidth=0.5pt, linecolor=Teal!55}

\mdtheorem[style=thmstyle]{theorem}{Theorem}[section]
\mdtheorem[style=propstyle]{proposition}[theorem]{Proposition}
\mdtheorem[style=corstyle]{corollary}[theorem]{Corollary}
\mdtheorem[style=lemstyle]{lemma}[theorem]{Lemma}

\mdtheorem[style=defstyle]{definition}[theorem]{Definition}
\mdtheorem[style=notstyle]{notation}[theorem]{Notation}
\mdtheorem[style=assstyle]{assumption}[theorem]{Assumption}

\mdtheorem[style=remstyle]{example}[theorem]{Example}
\newtheorem{remark}[theorem]{Remark}

\mdtheorem[style=exostyle]{exercise}[theorem]{Exercise}

\usepackage{needspace}
\makeatletter
\newcommand\mdf@avoidtopclip{%
  \Needspace*{6\baselineskip}%
  \ifdim\pagetotal<\baselineskip
    \vspace*{\dimexpr\topskip+0.5\baselineskip\relax}%
  \fi
}
\@for\@mdf@env:=theorem,proposition,corollary,lemma,definition,notation,assumption,example,exercise\do{%
  \expandafter\BeforeBeginEnvironment\expandafter{\@mdf@env}{\mdf@avoidtopclip}%
}
\makeatother

\usepackage[Sonny]{fncychap}   

\usepackage{fancyhdr}
\usepackage[bookmarks=true,
            colorlinks=true,
            urlcolor=blue,
            citecolor=Purple,
            linkcolor=blue,
            hypertexnames=false]{hyperref}

\usepackage{enumitem}
\setlist{parsep=0pt, leftmargin=5mm}
\setlist[itemize,enumerate]{nolistsep, itemsep=1pt, topsep=1pt}
\setlist[enumerate,1]{wide, label=\textnormal{\textbf{(\roman*)}}, labelindent=0pt, itemsep=.1in}

\usepackage[margin=10pt, font=small, labelfont=bf, labelsep=endash]{caption}

\usepackage{graphicx}
\usepackage{float}            

\usepackage{cleveref}
\crefname{exercise}{exercise}{exercises}
\crefname{theorem}{Theorem}{Theorems}
\crefname{lemma}{Lemma}{Lemmas}
\crefname{definition}{Definition}{Definitions}

\newcommand{\R}{\mathbb{R}}
\newcommand{\N}{\mathbb{N}}
\newcommand{\Z}{\mathbb{Z}}
\newcommand{\C}{\mathbb{C}}
\newcommand{\E}{\mathbb{E}}                 
\newcommand{\Prob}{\mathbb{P}}              

\graphicspath{{pics/}}

\begin{document}

\frontmatter

\begin{titlepage}
  \thispagestyle{empty}
  \vspace*{2cm}
  \begin{center}
    {\Huge\bfseries\scshape Singular parabolic Anderson models\\[4mm]}
    {\Large\itshape From flat space to geometric contexts}\\[2cm]
    {\large Samy Tindel}\\[0.5cm]
    {\normalsize Purdue University}\\[0.5cm]
    {\normalsize \today}
  \end{center}
  \vfill
\end{titlepage}
\setcounter{page}{2}   

\chapter*{Acknowledgments}
\addcontentsline{toc}{chapter}{Acknowledgments}
\markboth{ACKNOWLEDGMENTS}{ACKNOWLEDGMENTS}

This book grew out of the lecture notes I prepared for the first Aarhus Summer
School in Analysis and Probability, held in 2026. I wish to thank the organizers, Fabrice
Baudoin and Li Chen, for running it so smoothly and in such a relaxed
atmosphere. My thanks also go
to all the participants, who made the summer school an engaging and rewarding
experience.

My journey into the world of parabolic Anderson models has been shared with
many collaborators. I am grateful to Fabrice Baudoin, Le Chen, Li Chen,
Xia Chen, Aur\'elien Deya, Yaozhong Hu, Jingyu Huang, Khoa L\^e, David Nualart,
Cheng Ouyang, Jian Song, Frederi Viens, Jing Wang and Panqiu Xia, for the
mathematics we did together and for the pleasure of doing it.

Finally, I thank my family, Ana, Lea and Marco, for their constant support.
Living with a mathematician obsessed (among other things) with the parabolic Anderson model requires
a good deal of patience. They have shown it with good humor, even when the
conversation at the dinner table drifted toward intermittency.

\tableofcontents

\mainmatter

\chapter{Introduction}
\label{chap:intro}

This book is an introduction to the \emph{parabolic Anderson model} (PAM),
a stochastic partial differential equation lying at the interface of
mathematical physics, probability theory and analysis. We pursue two goals.
Our first aim is to describe the model on the Euclidean space~$\R^{d}$. We
shall explain what the equation is, why it has attracted attention in physics
and in probability, and which phenomena it accounts for. The two phenomena
under consideration are \emph{intermittency} and its polymer counterpart,
called \emph{localization}. Our second aim, which is the more original part
of this book, is to analyze the effect of the ambient geometry on the model.
Namely we shall replace the flat space~$\R^{d}$ by a sub-Riemannian manifold
such as a Heisenberg group, or by a fractal such as the Sierpi\'nski gasket,
and see how the exponents attached to the equation are modified.

This introductory chapter contains no proof. We first fix the equation once
and for all in Section~\ref{sec:pam-model}. Then we describe in
Sections~\ref{sec:intermittency} and~\ref{sec:polymers} the two phenomena
which motivate every estimate proved later. The analytic tools we shall need
are listed in Section~\ref{sec:techniques}. We close in
Section~\ref{sec:geometry} by explaining why the underlying geometry deserves
a treatment of its own.

\begin{notation}\label{not:proba}
  Throughout the book $\left(\Omega,\mathcal{F},\Prob\right)$ denotes an
  underlying probability space, and $\E$ the expectation with respect
  to~$\Prob$. We write $a\lesssim b$ when there is an unspecified constant
  $C>0$ such that $a\le C\,b$, and $a\asymp b$ when both $a\lesssim b$ and
  $b\lesssim a$ hold.
\end{notation}

\section{The parabolic Anderson model}
\label{sec:pam-model}

The name of the model goes back to P.~W.~Anderson's work on the quantum
mechanics of disordered media. In this context an electron evolves in a
random potential, and it may become spatially \emph{localized} instead of
delocalized. The stationary object behind this picture is a random
Schr\"odinger operator, called the \emph{Anderson Hamiltonian}
\begin{equation}\label{eq:anderson-hamiltonian}
  \mathcal{H} \,=\, \tfrac{1}{2}\Delta \,+\, V,
\end{equation}
where $\Delta$ is the Laplacian on~$\R^{d}$ and $V=V(x)$ is a random
potential. The \emph{parabolic} model is the heat flow generated
by~\eqref{eq:anderson-hamiltonian}. That is, one studies the large-time
behavior of $\exp(t\,\mathcal{H})$, or equivalently of the solution to the
Cauchy problem $\partial_{t}u=\mathcal{H}u$. In this book we shall focus on
the situation where the static potential $V$ is replaced by a
\emph{space-time} noise. Writing $\dot{W}$ for this noise, the parabolic
Anderson model on $\R_{+}\times\R^{d}$ reads
\begin{equation}\label{eq:pam}
  \partial_{t} u(t,x)
  \,=\, \tfrac{1}{2}\Delta\, u(t,x) \,+\, \beta\, u(t,x)\,\dot{W}(t,x),
  \qquad u(0,x)=u_{0}(x),
\end{equation}
where $\beta>0$ is a parameter tuning the intensity of the noise (the
inverse temperature, in the language of statistical mechanics).

Let us highlight three features of equation~\eqref{eq:pam}.

\begin{enumerate}
\item
{\it The noise is multiplicative.} In the right-hand side
of~\eqref{eq:pam}, the unknown $u$ multiplies the noise~$\dot{W}$. The medium
thus acts on the solution proportionally to the solution itself. This is the
mechanism which produces the tall and rare peaks described in
Section~\ref{sec:intermittency}. If one replaces $u\,\dot{W}$ by $\dot{W}$
alone, one gets the \emph{additive} stochastic heat equation, whose behavior
is much tamer. We shall nevertheless meet this equation in
Chapter~\ref{chap:existence}, as an intermediate step.

\item
{\it The noise is rough.} The driving noise $\dot{W}$ is white in time,
and either white or colored in space. Such a $\dot{W}$ is a random
distribution and not a function. In this sense equation~\eqref{eq:pam} is a
\emph{singular} stochastic partial differential equation. This singularity of
the potential~$\dot{W}$ is the feature which gives its title to the book. The product $u\,\dot{W}$ in~\eqref{eq:pam}
has thus no classical meaning, and a stochastic integration theory is needed
in order to give a rigorous sense to~\eqref{eq:pam}. We shall adopt the
It\^o--Walsh--Dalang point of view. Unless otherwise stated, our solutions
are therefore understood in the It\^o (or It\^o--Skorohod) sense. The
competing Stratonovich interpretation is the natural one for pathwise and
polymer questions, and we will comment on it when polymers enter the picture.
Let us also mention that for a rougher noise, like the space-time white noise
in spatial dimension $d\ge 2$, the stochastic integration approach breaks
down. The product $u\,\dot{W}$ becomes genuinely ill-posed, and the equation
has to be \emph{renormalized}. Making sense of the renormalized model is one
of the achievements of the theory of regularity
structures~\cite{hairer:14:theory}. We refer to~\cite{hairer:15:simple} for
the construction of the continuum parabolic Anderson model on~$\R^{2}$, and
to~\cite{chen:21:moment} for moment estimates in this renormalized setting.
Those pathwise techniques fall outside the scope of this book. We shall stay
within the range of regularity for which the It\^o--Walsh--Dalang theory
applies.

\item
{\it The regularity of the noise is a parameter.} The spatial regularity
of $\dot{W}$ will be encoded by a parameter $\alpha$. This parameter decides
\emph{whether} equation~\eqref{eq:pam} can be solved at all. It also governs
\emph{how} the solution behaves in large time. The interplay between this
regularity and the geometry of the underlying space is the main subject of
this book.
\end{enumerate}

\section{Intermittency}
\label{sec:intermittency}

The most striking feature of the parabolic Anderson model is
\emph{intermittency}. Roughly speaking, a field is intermittent when its mass
does not spread out as it would under the plain heat equation. Instead, as
time grows, the mass concentrates on a few small islands which are widely
separated, and on which the field reaches huge heights. This behavior is
observed for magnetic fields in turbulent astrophysical plasmas. This is the
context in which intermittency in random media was popularized by Zel'dovich
and his co-authors~\cite{zeldovich:87:intermittency}. The mathematical theory
of the parabolic Anderson model was then founded by Carmona and
Molchanov~\cite{carmona:94:parabolic}, and surveys can be found
in~\cite{gartner:05:parabolic,konig:16:parabolic}.

\begin{figure}[H]
  \centering
  \begin{minipage}[b]{0.42\textwidth}
    \centering
    \includegraphics[width=\linewidth, trim={18 12 18 32}, clip]{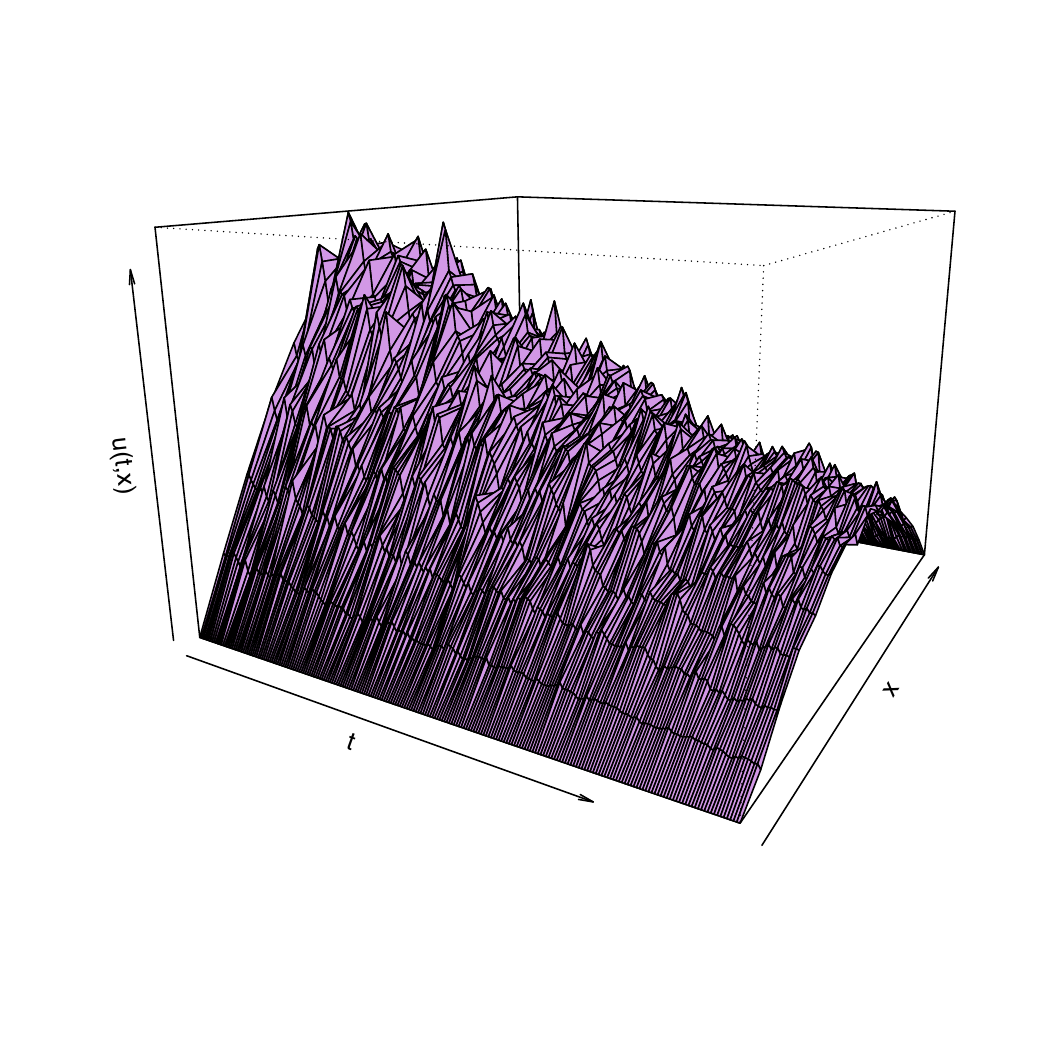}\\[-1mm]
    {\small $\lambda = 0.1$}
  \end{minipage}\hfill
  \begin{minipage}[b]{0.42\textwidth}
    \centering
    \includegraphics[width=\linewidth, trim={18 12 18 32}, clip]{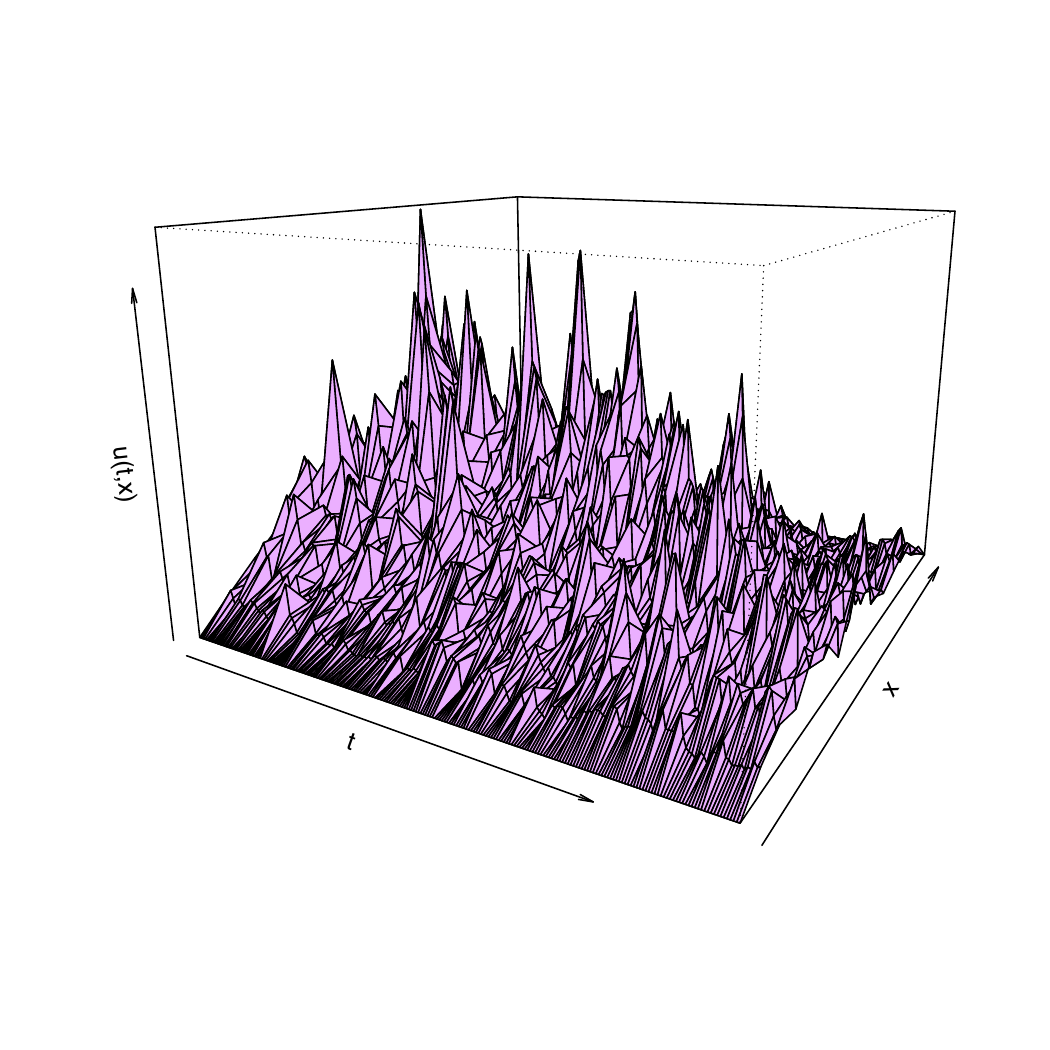}\\[-1mm]
    {\small $\lambda = 0.5$}
  \end{minipage}
  \caption[The weak-noise regime of the parabolic Anderson model]{The
    solution $u(t,x)$ of~\eqref{eq:pam} for small noise intensity $\lambda$:
    only a mild random perturbation of the smooth heat-equation profile.}
  \label{fig:intermittency-low}
\end{figure}

\begin{figure}[H]
  \centering
  \begin{minipage}[b]{0.42\textwidth}
    \centering
    \includegraphics[width=\linewidth, trim={18 12 18 32}, clip]{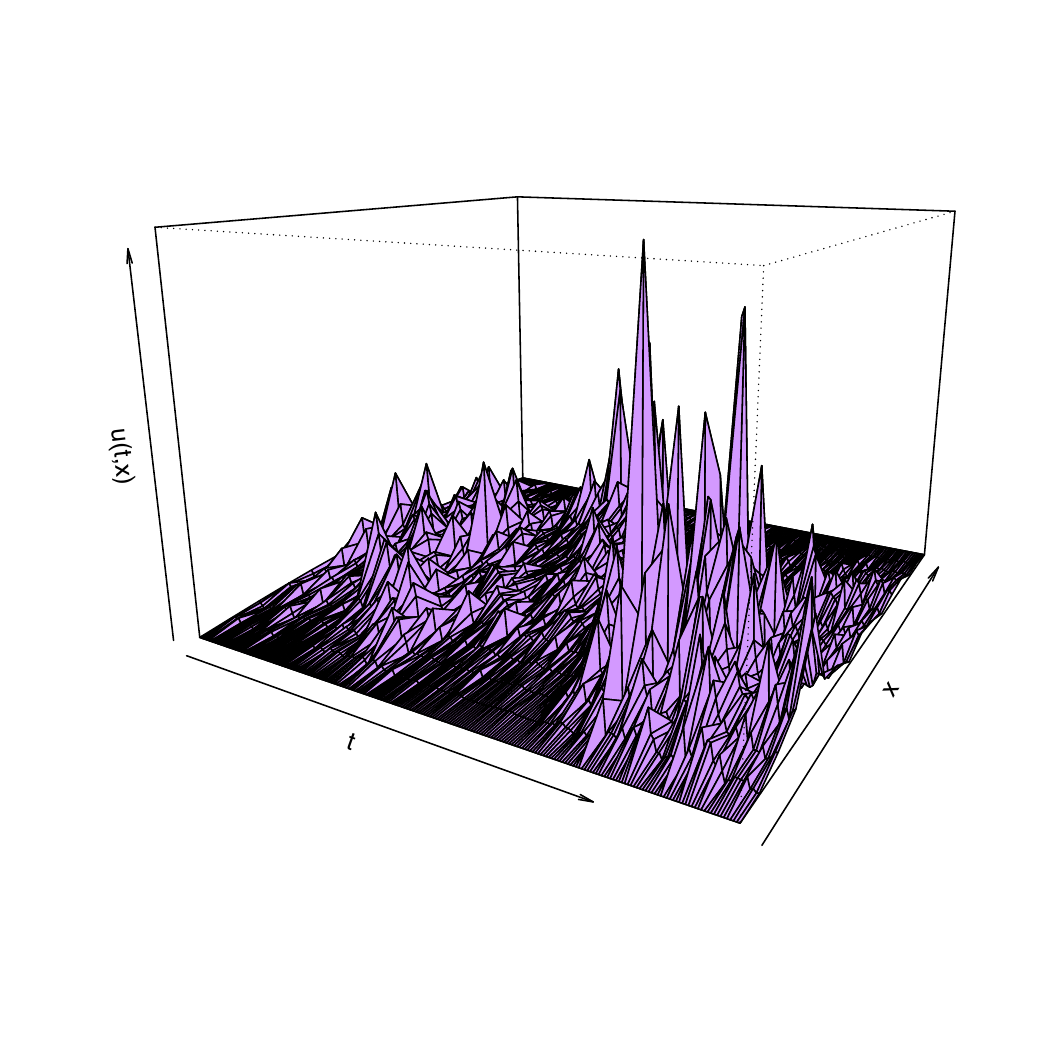}\\[-1mm]
    {\small $\lambda = 1$}
  \end{minipage}\hfill
  \begin{minipage}[b]{0.42\textwidth}
    \centering
    \includegraphics[width=\linewidth, trim={18 12 18 32}, clip]{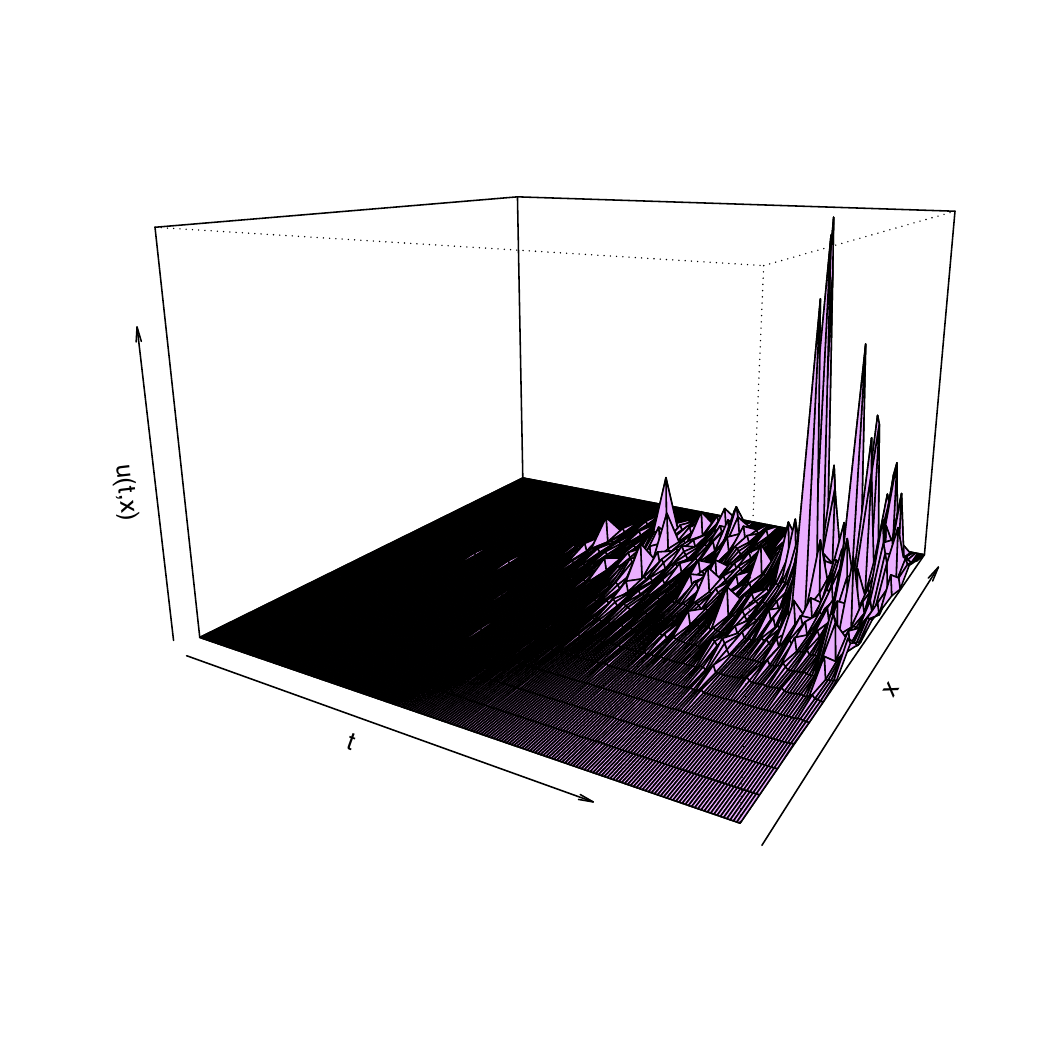}\\[-1mm]
    {\small $\lambda = 2$}
  \end{minipage}
  \caption[The onset of intermittency in the parabolic Anderson model]{The
    same solution for larger noise intensity $\lambda$: the mass concentrates
    on a few tall, well-separated peaks (the onset of intermittency).}
  \label{fig:intermittency-high}
\end{figure}

Figures~\ref{fig:intermittency-low} and~\ref{fig:intermittency-high}
illustrate this phenomenon. They display simulations of the solution $u(t,x)$
to~\eqref{eq:pam}, for increasing values of the noise intensity $\lambda$.
For a small intensity, the surface is a mild random perturbation of the
smooth heat-equation profile (see Figure~\ref{fig:intermittency-low}). The
profile is then progressively reshaped as $\lambda$ grows. For $\lambda=2$,
almost all the mass has collapsed onto a few sharp and isolated spikes,
standing against an essentially flat background (see
Figure~\ref{fig:intermittency-high}). This is a clear manifestation of
intermittency.

The quantitative counterpart of this picture is expressed in terms of the
moments of the solution. Namely, for $p\ge 1$ one defines the $p$-th
\emph{moment Lyapunov exponent} (also called upper Lyapunov exponent) as
\begin{equation}\label{eq:lyapunov}
  \gamma(p) \,:=\, \limsup_{t\to+\infty}\, \frac{1}{t}\,
  \ln \E\!\left[\,u(t,x)^{\,p}\,\right].
\end{equation}
Whenever the relevant quantities exist and are finite, the model is said to
be (fully) \emph{intermittent} if the normalized exponents are strictly
increasing, that is
\begin{equation}\label{eq:intermittency}
  \frac{\gamma(1)}{1} \,<\, \frac{\gamma(2)}{2} \,<\, \cdots \,<\,
  \frac{\gamma(p)}{p} \,<\, \cdots.
\end{equation}

\noindent
The strict inequalities in~\eqref{eq:intermittency} are the analytic
counterpart of the spikes observed in Figure~\ref{fig:intermittency-high}.
Indeed, higher moments are dominated by the rare and very high peaks. On the
exponential scale they thus grow strictly faster than what the lower moments
would predict. Establishing two-sided bounds on the
exponents~\eqref{eq:lyapunov} is one of the central goals of this book, and
this is the subject of Chapter~\ref{chap:moments}.

\section[Polymers]{Directed polymers and localization}
\label{sec:polymers}

Equation~\eqref{eq:pam} admits a second reading, which is borrowed from
statistical mechanics and gives a geometric meaning to intermittency. Indeed,
by the Feynman--Kac formula, the solution of~\eqref{eq:pam} can be
represented as
\begin{equation}\label{eq:feynman-kac}
  u(t,x) \,=\, \E_{B}\!\left[\, u_{0}(B_{t})\,
  \exp\Big(\textstyle\int_{0}^{t} \beta\,\dot{W}(t-s,B_{s})\,ds\Big)\right],
\end{equation}
where $B$ is a Brownian motion started at $x$ and $\E_{B}$ is the
expectation over $B$ only, the noise $\dot{W}$ being frozen. The right-hand
side of~\eqref{eq:feynman-kac} is the \emph{partition function} of a
\emph{directed polymer}. That is, a Brownian path is reweighted by the
exponential of the energy it collects along the random
environment~$\dot{W}$. Normalizing this weight, one defines the
\emph{polymer measure}. This is a random probability measure on paths, under
which the favorable regions of the environment are visited more often. In
this language, the parameter $\beta$ in~\eqref{eq:pam} is an inverse
temperature.

\begin{figure}[H]
  \centering
  \begin{minipage}[b]{0.42\textwidth}
    \centering
    \includegraphics[width=\linewidth, trim={70 55 55 40}, clip]{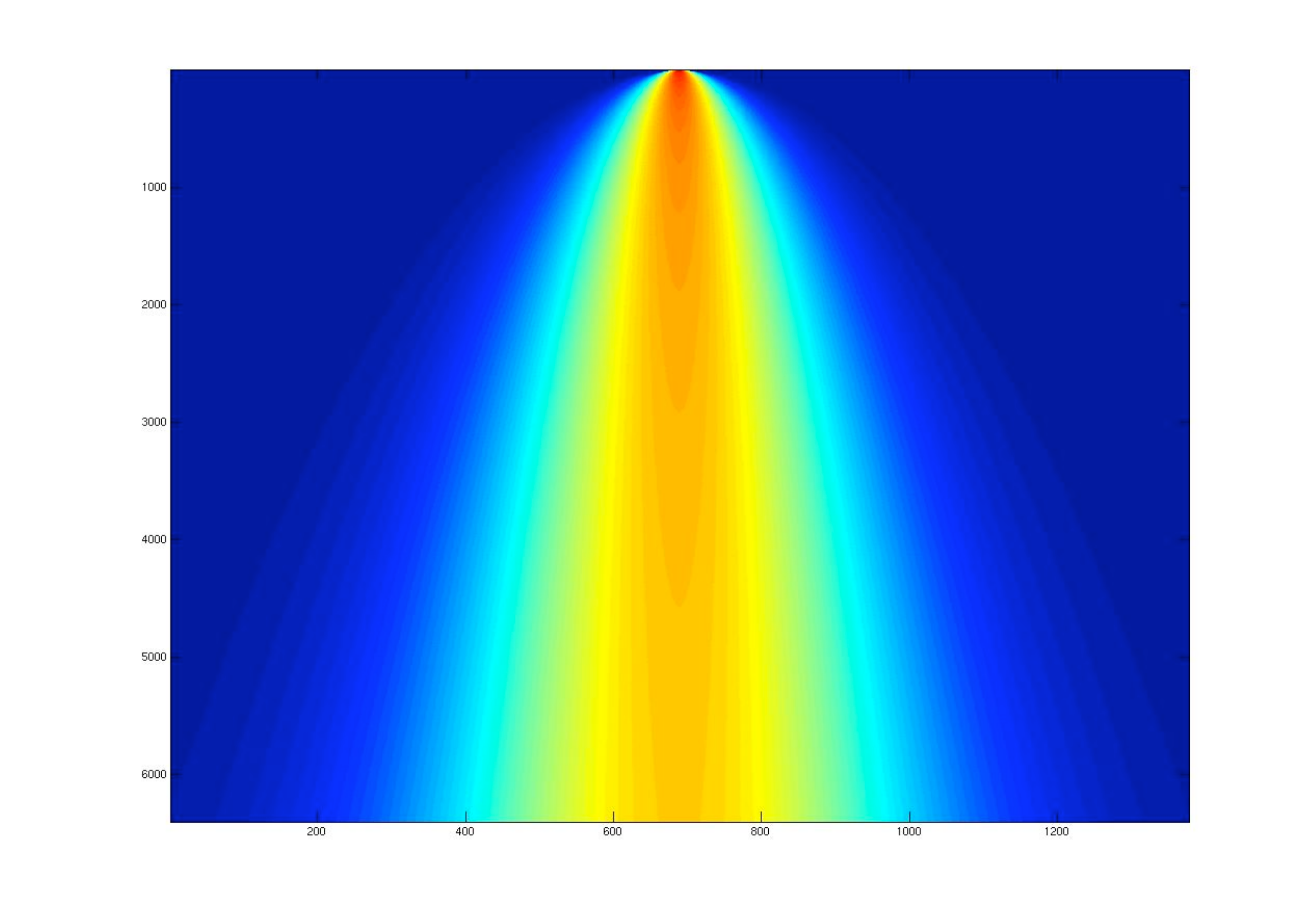}\\[1mm]
    {\small (a) diffusive behavior}
  \end{minipage}\hfill
  \begin{minipage}[b]{0.42\textwidth}
    \centering
    \includegraphics[width=\linewidth, trim={70 55 55 40}, clip]{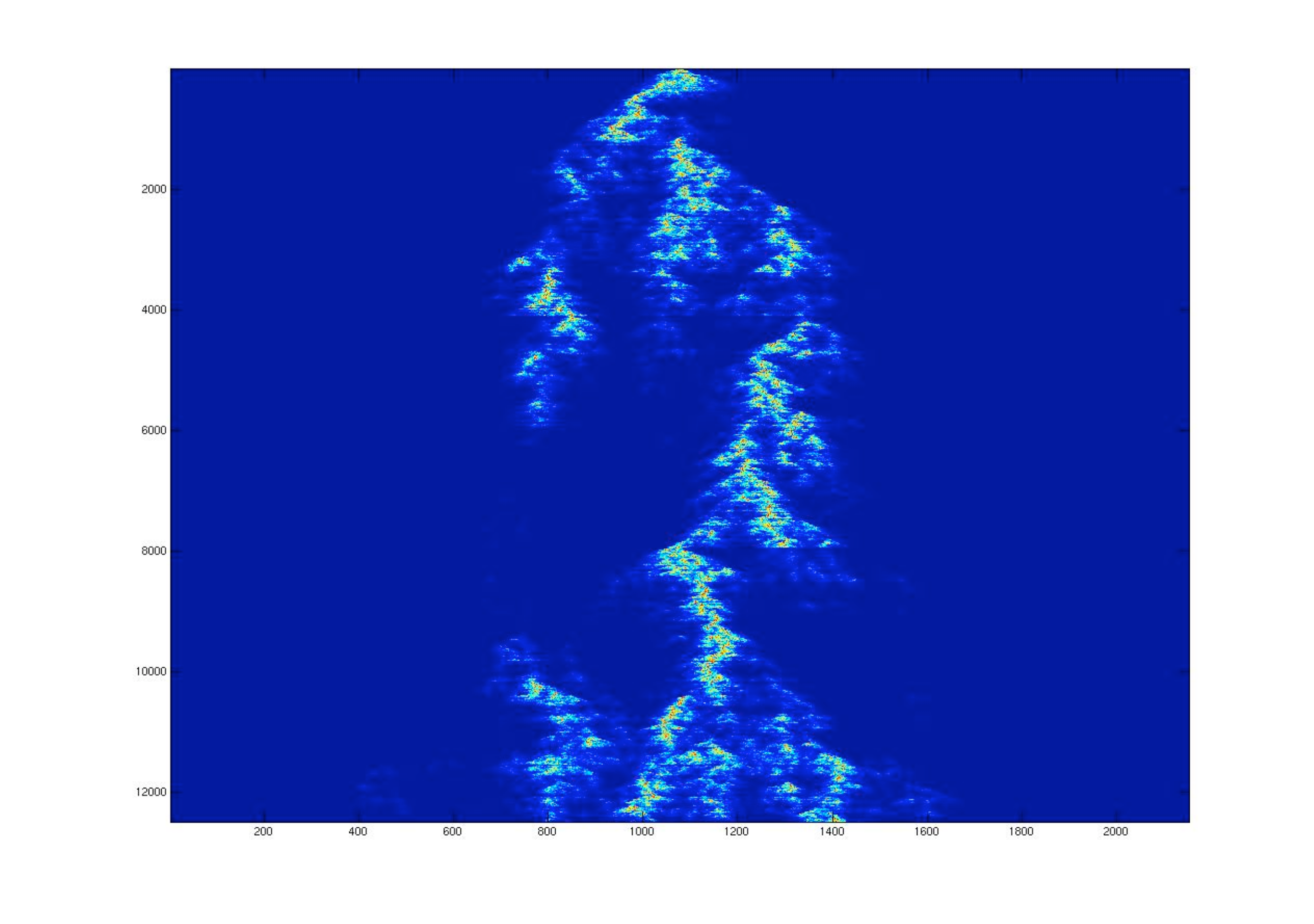}\\[1mm]
    {\small (b) localized behavior}
  \end{minipage}
  \caption[Diffusive versus localized polymer]{The weight of the polymer
    endpoint in space (horizontal) and time (vertical). \emph{(a)} In the
    weak-disorder regime the weight spreads parabolically: diffusive
    behavior. \emph{(b)} In the strong-disorder regime the weight
    concentrates on a thin wandering filament: localization.}
  \label{fig:polymer}
\end{figure}

With this correspondence in hand, the analytic question of intermittency
becomes a pathwise question about the polymer. It is organized around the
classical dichotomy between \emph{weak} and \emph{strong} disorder. In the
weak-disorder regime, the energy collected along the path is too weak to
overcome the entropy of the Brownian motion. The polymer measure stays close
to the Brownian one, and the path remains \emph{diffusive}, spreading
like~$\sqrt{t}$. In the strong-disorder regime the energy prevails over the
entropy. The path then \emph{localizes} and concentrates on a thin corridor
of exceptionally favorable sites. Those two regimes are displayed in
Figure~\ref{fig:polymer}, which shows the weight of the endpoint of the
polymer as a function of space (horizontal) and time (vertical).

Directed polymers have a long and rich history, both in the discrete setting
of random walks on $\Z^{d}$ and in the continuous setting underlying
equation~\eqref{eq:pam}. We refer to the Saint-Flour lecture notes of
Comets~\cite{comets:17:directed} for a panorama, and
to~\cite{rovira:05:brownian,marquez:11:model,lacoin:10:new,cosco:21:directed}
for results of the type which will guide us. Let us also stress that the
polymer point of view is one of the main reasons to study the parabolic
Anderson model beyond the flat case. Indeed, the localization corridor of
Figure~\ref{fig:polymer}\,(b) is a path in the underlying space. Its behavior
should therefore depend on the geometry of that space.

\section[Toolkit]{Overview of the toolkit}
\label{sec:techniques}

The analysis of equation~\eqref{eq:pam} relies on a small number of tools,
which we now list. Each of them will be developed in its own right later in
the book. Our aim here is only to present those methods and to indicate where
they will be used.

\begin{enumerate}
\item
{\it Stochastic calculus in infinite dimensions.} In order to give a
meaning to the product $u\,\dot{W}$ in~\eqref{eq:pam}, one needs a theory of
stochastic integration with respect to a space-time noise. We shall use the
random-field approach of Walsh and
Dalang~\cite{walsh:86:introduction,dalang:99:extending}, together with the
Malliavin calculus~\cite{nualart:06:malliavin} when the It\^o--Skorohod
interpretation is at stake. This is the content of
Chapters~\ref{chap:noise} and~\ref{chap:existence}.

\item
{\it The Feynman--Kac formula.} The
representation~\eqref{eq:feynman-kac} relates the equation to the polymer. It
is also our main tool for moment estimates, as we shall see in
Chapter~\ref{chap:moments}.

\item
{\it Intersection local times.} Computing the moments $\E[u(t,x)^{p}]$
through~\eqref{eq:feynman-kac}, one is led to expressions involving the
mutual intersection of $p$ independent copies of the underlying Brownian
motion. The corresponding intersection local times control the
exponents~\eqref{eq:lyapunov}.

\item
{\it Wiener chaos expansions.} Since the noise is Gaussian, the solution
of~\eqref{eq:pam} can be expanded in a series of multiple stochastic
integrals. This chaos expansion yields existence and uniqueness results, as
well as moment formulas (see Chapter~\ref{chap:existence}).

\item
{\it Heat semigroup estimates.} All the quantities above are governed by
the heat kernel $p_{t}(x,y)$ of the underlying Laplacian. Sharp upper and
lower bounds on $p_{t}$ are therefore essential for the whole theory. As
explained in the next section, this is also the point where the geometry of
the space enters.
\end{enumerate}

\section[Geometric structures]{From Euclidean space to geometric structures}
\label{sec:geometry}

On the flat space $\R^{d}$, the model~\eqref{eq:pam} is now fairly well
understood. The foundational analytic results we shall rely on go back
to~\cite{bertini:95:stochastic,carmona:94:parabolic,gartner:00:moment}, and
to~\cite{hu:09:stochastic,hu:15:stochastic} for the stochastic PDE
formulation. In this book we wish to understand what survives, and what
changes, when $\R^{d}$ is replaced by a space with a non-trivial geometry.

As mentioned above, the geometry enters through the heat kernel. On $\R^{d}$
this kernel is Gaussian, that is $p_{t}(x,y)\asymp
t^{-d/2}\exp(-|x-y|^{2}/2t)$, and everything is controlled by the single
number $d$. On a geometric space the Gaussian bound is replaced by a
sub-Gaussian one,
\begin{equation}\label{eq:subgaussian}
  p_{t}(x,y) \,\asymp\, t^{-d_{h}/d_{w}}
  \exp\left(-C\Big(\frac{d(x,y)^{d_{w}}}{t}\Big)^{\!\frac{1}{d_{w}-1}}\right),
\end{equation}
where \emph{two} exponents are now involved. The Hausdorff dimension
$d_{h}$ governs the growth of the volume of balls, while the walk dimension
$d_{w}\ge 2$ governs how fast the Brownian motion escapes them. One recovers
the Euclidean case for $d_{h}=d$ and $d_{w}=2$. On the Sierpi\'nski gasket,
for instance, one has $d_{h}=\ln 3/\ln 2$ and $d_{w}=\ln 5/\ln 2$. The
combination which governs the parabolic Anderson model is the \emph{spectral
dimension}
\begin{equation}\label{eq:spectral-dimension}
  d_{s} \,:=\, \frac{2\,d_{h}}{d_{w}},
\end{equation}
which measures the growth of the eigenvalues of the Laplacian.

Replacing $d$ by the pair $(d_{h},d_{w})$ has a real effect on the model.
Two concrete instances will recur throughout the book.

\begin{enumerate}
\item
{\it Sub-Riemannian spaces.} On the Heisenberg group $\mathbb{H}^{n}$ the
sub-Riemannian structure is degenerate, and the relevant notion of dimension
is no longer the topological dimension $d=2n+1$, but the Hausdorff
dimension $Q=2n+2$. The threshold of noise regularity for which
equation~\eqref{eq:pam} can be solved shifts accordingly. Written in terms of
$Q$ it coincides with the Euclidean condition, while written in terms of $d$
it does not. This effect can be traced back to the degeneracy of the
geometry~\cite{baudoin:23:parabolic-b,baudoin:23:parabolic-a}.

\item
{\it Fractals.} On post-critically finite fractals, the exponents $d_{h}$
and $d_{w}$ take a whole range of non-classical values. The Lyapunov
exponents~\eqref{eq:lyapunov} then acquire a rich dependence on the noise
regularity and on the boundary conditions~\cite{baudoin:24:parabolic}. In
particular the condition $d_{s}<2$, that is $d_{h}<d_{w}$, makes even the
space-time white noise admissible in~\eqref{eq:pam}. This phenomenon has no
Euclidean analogue in dimension $d\ge 2$.
\end{enumerate}

Those two examples are our main motivation for this book.

\section[Outline]{Outline of the book}
\label{sec:outline}

The book is organized as follows.

\begin{enumerate}
\item
{\it Chapter~\ref{chap:noise}, Gaussian noises and their regularity.} We
introduce the Wiener space and the Malliavin calculus related to a noise which
is white in time and colored in space. We then measure the spatial regularity
of this noise in a Besov scale. The Riesz, Bessel and fractional kernels are
our main examples.

\item
{\it Chapter~\ref{chap:existence}, existence and uniqueness.} We construct
the stochastic integral, following Walsh and then Dalang, and we solve the
equation. We first handle the additive case, for which the solution is
explicit and the question is the regularity of the stochastic convolution. The
multiplicative case is then treated by a contraction argument.

\item
{\it Chapter~\ref{chap:moments}, moment estimates.} We resort to chaos
expansions and to a Feynman-Kac representation of the moments of the solution.
This yields matching upper and lower bounds on the Lyapunov
exponents~\eqref{eq:lyapunov}. The intermittency described above follows.

\item
{\it Chapter~\ref{chap:geom-setting}, a general geometric setting.} We
introduce the metric measure spaces on which the second half of the book is
set. Namely we focus on the spaces carrying a sub-Gaussian heat kernel, with a
Hausdorff dimension $d_{h}$ and a walk dimension $d_{w}$. We also describe the
Dirichlet and Neumann problems on a domain. The Sierpi\'nski gasket and the
metric graphs are the examples to keep in mind.

\item
{\it Chapter~\ref{chap:geom-pam}, PAM on metric spaces.} We carry out the
program of the three Euclidean chapters in this geometric setting. The noise is
now produced by an inverse fractional Laplacian. Existence is obtained through
the chaos expansion instead of a contraction argument. The moment estimates are
read off a Feynman-Kac formula written with the underlying diffusion.
\end{enumerate}

Let us close this outline with a remark about the Euclidean part of the
book. The special form of the multiplicative term in~\eqref{eq:pam} is not
used in Chapters~\ref{chap:noise} and~\ref{chap:existence}. Indeed, the
integration theory of Section~\ref{sec:integration} and the fixed point
argument of Section~\ref{sec:contraction} are performed for a general
Lipschitz coefficient $\sigma(u)$, the parabolic Anderson model corresponding
to the case $\sigma(r)=\beta r$ (see Corollary~\ref{cor:nonlinear-wp}). Read
in this way, those two chapters can be used as a basic introduction to
stochastic partial differential equations in dimension $d\ge1$. 

\chapter[Gaussian noises]{Gaussian noises and their regularity}
\label{chap:noise}

\section[Malliavin calculus]{Wiener space and Malliavin calculus}\label{sec:wiener-malliavin}

In this book, the noise $\dot{W}$ driving the PAM equation~\eqref{eq:pam} will be a random
distribution. We will thus see it as a random linear functional on
suitable test functions. The Gaussian Hilbert space structure of this noise,
together with the related Malliavin calculus, is at the basis of all our
subsequent analysis. This section is devoted to the corresponding framework.
We refer to Nualart~\cite{nualart:06:malliavin} for the Malliavin calculus,
and to Hu, Huang, Nualart and Tindel~\cite{hu:15:stochastic} for the specific
Gaussian noises we work with.

\subsection{The driving noise}\label{ssec:noise}

On the probability space $(\Omega,\mathcal{F},\Prob)$ of
Notation~\ref{not:proba}, the space-time noise $W$ is encoded by a
centered Gaussian family
$\{W(\varphi);\,\varphi\in\mathcal{D}([0,\infty)\times\R^{d})\}$.
Throughout this book the noise is \emph{white in time} and correlated in
space only, with spatial covariance function $\Lambda:\R^{d}\to\R_{+}$:
\begin{equation}\label{eq:noise-cov}
  \E\bigl[W(\varphi)\,W(\psi)\bigr]
  = \int_{\R_{+}\times\R^{2d}}
    \varphi(t,x)\,\psi(t,y)\,\Lambda(x-y)\,
    dx\,dy\,dt.
\end{equation}
We assume throughout that $\Lambda$ is non-negative definite, so that its
Fourier transform $\mu:=\mathcal{F}\Lambda$ is a non-negative tempered
measure on $\R^{d}$, called the \emph{spectral measure} of the noise.
In order to define stochastic integrals with respect to $W$ in
$L^{2}(\Omega)$, the spectral measure should not charge high frequencies too
heavily. Namely the heat kernel $p_{t}$, whose Fourier transform is
$e^{-t|\xi|^{2}/2}$, has to belong to $\mathcal{H}$ for each $t>0$. This
amounts to the condition $\int_{\R^{d}}e^{-t|\xi|^{2}}\mu(d\xi)<\infty$. A
uniform version of this requirement, valid for all $t\ge 0$, is Dalang's
condition~\cite{dalang:99:extending}.

\begin{assumption}\label{hyp:spectral}
  The spectral measure $\mu$ satisfies Dalang's condition:
  \begin{equation}\label{eq:dalang-cond}
    \int_{\R^{d}}\frac{\mu(d\xi)}{1+|\xi|^{2}} < \infty.
  \end{equation}
\end{assumption}

\noindent
As we shall see in Chapter~\ref{chap:existence}, this integrability
condition is exactly what is needed in order to define the stochastic
integral with respect to $W$.

The Hilbert space on which $W$ acts as an isometry is obtained by completing
the space of test functions, with respect to the inner product induced by the
covariance~\eqref{eq:noise-cov}. Let thus $\mathcal{H}$ be the completion of
$\mathcal{D}([0,\infty)\times\R^{d})$ with respect to the inner product
\begin{equation}\label{eq:H-innprod}
  \langle\varphi,\psi\rangle_{\mathcal{H}}
  = \int_{\R_{+}\times\R^{2d}}
    \varphi(t,x)\,\psi(t,y)\,\Lambda(x-y)\,
    dx\,dy\,dt.
\end{equation}
Equivalently, writing $\mathcal{F}_{x}$ for the Fourier transform in the
space variable, one can write~\eqref{eq:H-innprod} as
\begin{equation}\label{eq:H-innprod-fourier}
  \langle\varphi,\psi\rangle_{\mathcal{H}}
  = \int_{\R_{+}\times\R^{d}}
    \mathcal{F}_{x}\varphi(t,\xi)\,\overline{\mathcal{F}_{x}\psi(t,\xi)}\,
    \mu(d\xi)\,dt.
\end{equation}
The covariance identity $\E[W(\varphi)W(\psi)]=\langle\varphi,\psi\rangle_{\mathcal{H}}$
means that $\varphi\mapsto W(\varphi)$ extends to a linear isometry
\begin{equation}\label{eq:W-isometry}
  W:\mathcal{H}\longrightarrow L^{2}(\Omega),
  \qquad
  \E[W(h)\,W(k)] = \langle h,k\rangle_{\mathcal{H}},
\end{equation}
and throughout the book we will work with Wiener integrals written as 
\begin{equation*}
W(h)=\int_{0}^{\infty}\int_{\R^{d}}h(s,x)\,W(ds,dx),
\quad\text{for}\quad
h\in\mathcal{H}. 
\end{equation*}

\begin{example}[Riesz kernel]
  The primary example throughout the book is the \emph{$\eta$-Riesz
  spatial covariance}, defined for $0<\eta<d$ by
  \begin{equation}\label{eq:riesz-kernel}
    \Lambda(x) = |x|^{-\eta},
    \qquad
    \mu(d\xi) = c_{\eta,d}\,|\xi|^{-(d-\eta)}\,d\xi,
  \end{equation}
  where $c_{\eta,d}>0$ is a normalizing constant.
  Assumption~\ref{hyp:spectral} holds precisely when $0<\eta<2$
  (Lemma~\ref{lem:dalang-riesz} in Section~\ref{sec:noise-def}).
  The limiting case $\eta\to 0$ corresponds informally to white noise in space.
\end{example}

\subsection{Isonormal Gaussian processes and Wiener chaos}\label{ssec:isonormal}

The isometry~\eqref{eq:W-isometry} is the defining property of a general
class of Gaussian processes, which we now introduce.

\begin{definition}[Isonormal Gaussian process]\label{def:isonormal}
  Let $\mathcal{H}$ be a separable real Hilbert space. A centered Gaussian
  family $W=\{W(h);\,h\in\mathcal{H}\}$ on $(\Omega,\mathcal{F},\Prob)$
  is called an \emph{isonormal Gaussian process on $\mathcal{H}$} if
  \begin{equation}\label{eq:isonormal}
    \E[W(h_{1})\,W(h_{2})] = \langle h_{1},h_{2}\rangle_{\mathcal{H}}
    \qquad\text{for all }h_{1},h_{2}\in\mathcal{H}.
  \end{equation}
\end{definition}

The noise of Section~\ref{ssec:noise} is an isonormal Gaussian process on
$\mathcal{H}$. The Hilbert space structure of $\mathcal{H}$ yields a
canonical orthogonal decomposition of $L^{2}(\Omega)$, whose building blocks
are defined thanks to Hermite polynomials. Let us recall that the $n$-th
\emph{probabilists' Hermite polynomial} $H_{n}$ is characterized by
\begin{equation}\label{eq:hermite}
  \exp\Bigl(tz - \frac{t^{2}}{2}\Bigr)
  = \sum_{n=0}^{\infty}\frac{t^{n}}{n!}\,H_{n}(z),
  \qquad t,z\in\R,
\end{equation}
giving the first values
\begin{equation}\label{eq:hermite-values}
  H_{0}(z) = 1, \qquad
  H_{1}(z) = z, \qquad
  H_{2}(z) = z^{2}-1, \qquad
  H_{3}(z) = z^{3}-3z.
\end{equation}

For every unit vector $h\in\mathcal{H}$, the random variable $W(h)$ is
standard Gaussian. The composition $H_{n}(W(h))$ is thus a well-defined
element of $L^{2}(\Omega)$ for each $n\ge 0$. Let now $h$ vary over all unit
vectors, and take closures in $L^{2}(\Omega)$. For each $n$, one obtains in
this way a closed subspace which stands for the degree $n$ part of the
Gaussian structure.

\begin{definition}[Wiener chaos]\label{def:chaos}
  For $n\ge 1$, the \emph{$n$-th homogeneous Wiener chaos}
  $\mathbf{H}_{n}$ is the closed subspace of $L^{2}(\Omega)$ generated by
  \begin{equation}\label{eq:chaos-span}
    \bigl\{H_{n}(W(h));\;h\in\mathcal{H},\;\|h\|_{\mathcal{H}}=1\bigr\}.
  \end{equation}
  We set $\mathbf{H}_{0}:=\R$ (the constants). The subspaces
  $\{\mathbf{H}_{n}\}_{n\ge 0}$ are pairwise orthogonal in $L^{2}(\Omega)$.
\end{definition}

Each chaos $\mathbf{H}_{n}$ can be parametrized by symmetric kernels.
Namely, the \emph{$n$-th multiple stochastic integral} $I_{n}(f)$ of a kernel
$f\in\mathcal{H}^{\odot n}$ (the $n$-th symmetric tensor product of
$\mathcal{H}$) is the unique element of $\mathbf{H}_{n}$ extending, by
linearity and $L^{2}$-continuity, the relation
\begin{equation}\label{eq:In-def}
  I_{n}(h^{\otimes n}) := H_{n}(W(h)),
  \qquad h\in\mathcal{H},\;\|h\|_{\mathcal{H}}=1.
\end{equation}
The resulting maps $I_{n}:\mathcal{H}^{\odot n}\to\mathbf{H}_{n}$ are
isometries: for all $f\in\mathcal{H}^{\odot n}$ and $g\in\mathcal{H}^{\odot m}$,
\begin{equation}\label{eq:In-isometry}
  \E\bigl[I_{n}(f)\,I_{m}(g)\bigr]
  = \mathbf{1}_{n=m}\cdot n!\,\langle f,g\rangle_{\mathcal{H}^{\otimes n}}.
\end{equation}
The factor $n!$ is due to the symmetrization, since each symmetric kernel
$f$ is built from $n!$ permutations of its arguments. Notice also that the
pairwise orthogonality in Definition~\ref{def:chaos} is an immediate
consequence of~\eqref{eq:In-isometry} with $n\ne m$.

With the maps $I_{n}$ in hand, one may wonder whether the orthogonal direct
sum $\bigoplus_{n\ge 0}\mathbf{H}_{n}$ exhausts
$L^{2}(\Omega,\sigma(W),\Prob)$. The answer is given by the following
classical result.

\begin{theorem}[Wiener chaos decomposition; {\cite[Theorem~1.1.1]{nualart:06:malliavin}}]
  \label{thm:chaos-decomp}
  Every $F\in L^{2}(\Omega,\sigma(W),\Prob)$ admits a unique orthogonal
  decomposition
  \begin{equation}\label{eq:chaos-decomp}
    F = \E[F] + \sum_{n=1}^{\infty}I_{n}(f_{n}),
  \end{equation}
  where the series converges in $L^{2}(\Omega)$ and $f_{n}\in\mathcal{H}^{\odot n}$
  for each $n\ge 1$.
\end{theorem}

Decomposition~\eqref{eq:chaos-decomp} will be called the
\emph{Wiener chaos expansion} of $F$.

\subsection{Malliavin derivative and Skorohod integral}\label{ssec:malliavin}

Two canonical operators are related to the Wiener chaos structure, namely a
derivative $D$ and its $L^{2}$ adjoint $\delta$. To begin with, the
derivative $D$ can be seen as a gradient with respect to the underlying
Gaussian randomness. Given a direction $h\in\mathcal{H}$, one
can perturb the noise $W$ by the deterministic shift $\varepsilon h$: for
every $k\in\mathcal{H}$, replace $W(k)$ by $W(k)+\varepsilon\langle
k,h\rangle_{\mathcal{H}}$. For a smooth cylindrical random variable
$F=f(W(h_{1}),\dots,W(h_{n}))$, this transforms $F$ into
\begin{equation*}
f(W(h_{1})+\varepsilon\langle h_{1},h\rangle_{\mathcal{H}},\dots,
W(h_{n})+\varepsilon\langle h_{n},h\rangle_{\mathcal{H}}) .
\end{equation*}
Differentiating this expression at $\varepsilon=0$, one gets the \emph{directional
Malliavin derivative of $F$ in the direction $h$}:
\begin{eqnarray}\label{eq:dir-derivative}
  D_{h}F
  &=& \frac{d}{d\varepsilon}\bigg|_{\varepsilon=0}
  f\bigl(W(h_{1})+\varepsilon\langle h_{1},h\rangle_{\mathcal{H}},
  \dots,W(h_{n})+\varepsilon\langle h_{n},h\rangle_{\mathcal{H}}\bigr) \notag\\
  &=& \sum_{j=1}^{n}
  \frac{\partial f}{\partial x_{j}}(W(h_{1}),\dots,W(h_{n}))\,
  \langle h_{j},h\rangle_{\mathcal{H}}.
\end{eqnarray}
The map $h\mapsto D_{h}F$ is a bounded linear functional on $\mathcal{H}$.
Owing to the Riesz representation theorem, there exists thus a unique
$\mathcal{H}$-valued random variable $DF$ such that
\begin{equation}\label{eq:dir-riesz}
  D_{h}F = \langle DF, h\rangle_{\mathcal{H}}
  \qquad\text{for all }h\in\mathcal{H}.
\end{equation}
The random variable $DF$ is called the Malliavin derivative of $F$. Reading
the $\mathcal{H}$-representative in~\eqref{eq:dir-derivative}, one obtains the
following explicit formula.

\begin{definition}[Malliavin derivative]\label{def:malliavin-D}
  A random variable is \emph{smooth cylindrical} if it has the form
  $F=f(W(h_{1}),\dots,W(h_{n}))$ for some $h_{1},\dots,h_{n}\in\mathcal{H}$
  and $f\in C^{\infty}_{p}(\R^{n})$ (smooth with all derivatives of
  polynomial growth). The \emph{Malliavin derivative} $DF$ is the
  $\mathcal{H}$-valued random variable
  \begin{equation}\label{eq:malliavin-D}
    DF = \sum_{j=1}^{n}
    \frac{\partial f}{\partial x_{j}}(W(h_{1}),\dots,W(h_{n}))\,h_{j}.
  \end{equation}
\end{definition}

Formula~\eqref{eq:malliavin-D} defines $D$ on smooth cylindrical variables
only. In order to extend it to a larger class of random variables, one has to
check that limits of smooth cylindrical approximations admit a well-defined
derivative. This is the closability property below, whose proof can be found
in~\cite[Proposition~1.2.1]{nualart:06:malliavin}.

\begin{proposition}[Closability of $D$]\label{prop:D-closable}
  The Malliavin derivative $D$, defined on smooth cylindrical variables
  by~\eqref{eq:malliavin-D}, is closable as an operator from $L^{2}(\Omega)$
  to $L^{2}(\Omega;\mathcal{H})$: if $(F_{n})$ is a sequence of smooth
  cylindrical variables satisfying $F_{n}\to 0$ in $L^{2}(\Omega)$ and
  $DF_{n}\to G$ in $L^{2}(\Omega;\mathcal{H})$, then $G=0$.
\end{proposition}

Proposition~\ref{prop:D-closable} guarantees that the norm below is
well-defined. Indeed, if two approximating sequences of smooth cylindrical
variables converge to the same limit in $L^{2}(\Omega)$, their derivatives
converge to the same limit in $L^{2}(\Omega;\mathcal{H})$. The extension of
$D$ is thus unambiguous. The \emph{Malliavin--Sobolev space}
$\mathbb{D}^{1,2}$ is then defined as the closure of the smooth cylindrical
variables under
\begin{equation}\label{eq:D12-norm}
  \|F\|_{1,2}^{2}
  = \E[F^{2}]+\E\bigl[\|DF\|_{\mathcal{H}}^{2}\bigr].
\end{equation}
The action of the derivative on chaos expansions is also easily described.
Namely, if $F=\sum_{n\ge 0}I_{n}(f_{n})$ belongs to $\mathbb{D}^{1,2}$, then
\begin{equation}\label{eq:D-chaos}
  DF = \sum_{n=1}^{\infty}n\,I_{n-1}\bigl(f_{n}(\,\cdot\,,\star)\bigr),
\end{equation}
where $\star$ stands for the free variable in $\mathcal{H}$ and
$f_{n}(\,\cdot\,,\star)$ denotes $f_{n}$ with one argument singled out.

Since $D$ is a closed operator from $L^{2}(\Omega)$ into
$L^{2}(\Omega;\mathcal{H})$, it admits a well-defined $L^{2}$-adjoint. This
adjoint maps $\mathcal{H}$-valued square integrable processes to scalar
random variables. As we shall see, it plays the role of a stochastic integral
for non-adapted integrands.

\begin{definition}[Skorohod integral]\label{def:skorohod}
  The \emph{Skorohod integral} $\delta$, also called the \emph{divergence
  operator}, is the $L^{2}(\Omega;\mathcal{H})$-adjoint of $D$. An element
  $u\in L^{2}(\Omega;\mathcal{H})$ belongs to $\mathrm{Dom}(\delta)$ if
  the map $F\mapsto\E[\langle DF,u\rangle_{\mathcal{H}}]$ extends to a
  bounded linear functional on $L^{2}(\Omega)$. In that case, $\delta(u)$
  is the unique element of $L^{2}(\Omega)$ satisfying the duality relation
  \begin{equation}\label{eq:duality}
    \E\bigl[\delta(u)\,F\bigr]
    = \E\bigl[\langle DF,u\rangle_{\mathcal{H}}\bigr]
    \qquad\text{for all }F\in\mathbb{D}^{1,2}.
  \end{equation}
\end{definition}

\noindent
Expressing the inner product in~\eqref{eq:duality} thanks
to~\eqref{eq:H-innprod}, one gets an explicit space-time formula. Namely, for
$u\in\mathrm{Dom}(\delta)$ and $F\in\mathbb{D}^{1,2}$ we have
\begin{equation}\label{eq:duality-explicit}
  \E\bigl[\delta(u)\,F\bigr]
  = \E\!\int_{\R_{+}\times\R^{2d}}
    D_{t,x}F\;u_{t,y}\;\Lambda(x-y)\,dx\,dy\,dt,
\end{equation}
where $D_{t,x}F$ denotes the Malliavin derivative of $F$ evaluated at
$(t,x)\in\R_{+}\times\R^{d}$. When $u$ is adapted to the filtration
generated by $W$, the Skorohod integral $\delta(u)$ coincides with the
It\^o--Walsh integral $\int u\,W(dt,dx)$. In general, $\delta$ extends the
It\^o integral to non-adapted integrands.

We close this section with three identities that will be invoked repeatedly.
\begin{enumerate}
\item
For any $F\in\mathbb{D}^{1,2}$ and any $h\in\mathcal{H}$,
\begin{equation}\label{eq:product-rule}
  F\,W(h) = \delta(F\,h) + \langle DF,h\rangle_{\mathcal{H}}.
\end{equation}

\item
For any $F\in\mathbb{D}^{2,2}$ and $h,g\in\mathcal{H}$,
\begin{equation}\label{eq:second-duality}
  \E\bigl[F\,W(h)\,W(g)\bigr]
  = \E\bigl[\langle D^{2}F,h\otimes g\rangle_{\mathcal{H}^{\otimes 2}}\bigr]
  + \E[F]\,\langle h,g\rangle_{\mathcal{H}}.
\end{equation}

\item
For any $v\in\mathrm{Dom}(\delta)$ with chaos decomposition
$v = \sum_{n=0}^{\infty}I_{n}(g_{n}(\,\cdot\,,\star))$, where $\star$ denotes
the free variable in $\mathcal{H}$ (as in~\eqref{eq:D-chaos}) and each
$g_{n}(\,\cdot\,,\star)\in\mathcal{H}^{\odot n}$ is symmetric in its first $n$
variables, one can express $\delta(v)$ as
\begin{equation}\label{eq:delta-chaos}
  \delta(v) = \sum_{n=0}^{\infty} I_{n+1}\!\left(\widetilde{g}_{n}\right),
\end{equation}
where $\widetilde{g}_{n}$ denotes the symmetrization of
$g_{n}(\,\cdot\,,\star)$ in all $n+1$ variables of $\mathcal{H}$.

\end{enumerate}
Identities~\eqref{eq:product-rule} and~\eqref{eq:second-duality} follow from
the duality formula~\eqref{eq:duality}. Identity~\eqref{eq:delta-chaos}
follows from the action~\eqref{eq:D-chaos} of $D$ on chaos. We refer
to~\cite[Chapter~1]{nualart:06:malliavin} for the proofs.

We close this section with an important comparison of norms for Wiener chaos
elements. Namely the $L^{k}$ norm is bounded by the $L^{2}$ norm, at the cost
of an explicit factor in $n$ and $k$. This will be one of our main tools for
the moment estimates of Section~\ref{sec:moment-estimates}.

\begin{theorem}[Nelson hypercontractivity]\label{thm:hypercontractivity}
Let $k\ge 2$ be an integer. For any $n\ge 0$ and any $f_{n}\in\mathcal{H}^{\odot n}$,
the following holds true:
\begin{equation}\label{eq:hypercontractivity}
  \bigl\|I_{n}(f_{n})\bigr\|_{L^{k}(\Omega)}
  \;\le\; (k-1)^{n/2}\,\bigl\|I_{n}(f_{n})\bigr\|_{L^{2}(\Omega)}.
\end{equation}
Consequently, for any $F$ admitting the chaos decomposition $F=\sum_{n=0}^{\infty}I_{n}(f_{n})\in L^{2}(\Omega)$
with $f_{n}\in\mathcal{H}^{\odot n}$,
we have
\begin{equation}\label{eq:hypercontractivity-sum}
  \bigl\|F\bigr\|_{L^{k}(\Omega)}
  \;\le\;
  \sum_{n=0}^{\infty}(k-1)^{n/2}\,\sqrt{n!}\,\|f_{n}\|_{\mathcal{H}^{\otimes n}}.
\end{equation}
\end{theorem}

\begin{proof}
Inequality~\eqref{eq:hypercontractivity-sum} follows immediately
from~\eqref{eq:hypercontractivity} by the triangle inequality in $L^{k}(\Omega)$
and the isometry $\|I_{n}(f_{n})\|_{L^{2}}=\sqrt{n!}\,\|f_{n}\|_{\mathcal{H}^{\otimes n}}$
from~\eqref{eq:In-isometry}. It therefore suffices to
prove~\eqref{eq:hypercontractivity}.

In order to show~\eqref{eq:hypercontractivity}, we introduce the \emph{Ornstein-Uhlenbeck semigroup} $(P_{t}^{\mathrm{OU}})_{t\ge 0}$ on
$L^{2}(\Omega)$, which acts on each Wiener chaos level by
\begin{equation}\label{eq:OU-eigen}
  P_{t}^{\mathrm{OU}}\,I_{n}(f_{n}) \;=\; e^{-nt}\,I_{n}(f_{n}),
\end{equation}
so that $I_{n}(f_{n})$ is an eigenfunction of the generator with eigenvalue $n$.
Nelson's theorem asserts that $P_{t}^{\mathrm{OU}}:L^{2}(\Omega)\to L^{q}(\Omega)$
with operator norm at most $1$ whenever $q\le 1+e^{2t}$;
see~\cite[Theorem~1.4.1]{nualart:06:malliavin}.
Setting $t:=\frac{1}{2}\log(k-1)$, so that $1+e^{2t}=k$, and
using~\eqref{eq:OU-eigen}, we obtain
\begin{equation*}
  \bigl\|I_{n}(f_{n})\bigr\|_{L^{k}(\Omega)}
  \;=\; e^{nt}\,\bigl\|P_{t}^{\mathrm{OU}}\,I_{n}(f_{n})\bigr\|_{L^{k}(\Omega)}
  \;\le\; e^{nt}\,\bigl\|I_{n}(f_{n})\bigr\|_{L^{2}(\Omega)}
  \;=\; (k-1)^{n/2}\,\bigl\|I_{n}(f_{n})\bigr\|_{L^{2}(\Omega)},
\end{equation*}
which is~\eqref{eq:hypercontractivity}.
\end{proof}

\section{Besov spaces}\label{sec:besov}

Besov spaces form a flexible scale of function spaces, which interpolates
between the classical Sobolev and H\"{o}lder spaces. They are well adapted to
the regularity of distributions such as our noise $\dot{W}$. Their definition
is based on the \emph{Littlewood--Paley decomposition}. This decomposition
splits a tempered distribution into spectrally localized blocks, and the
regularity is then read off the dyadic decay of those blocks.

\subsection{Littlewood--Paley decomposition}\label{ssec:LP}

The decomposition is built from a smooth, dyadic partition of unity in
Fourier space.

\begin{notation}[Littlewood--Paley cutoffs]\label{not:LP-partition}
  Fix two smooth, radial, nonnegative functions $\chi,\varphi:\R^{d}\to[0,1]$
  with compact support satisfying:
  \begin{enumerate}
    \item the support of $\chi$ is contained in a ball $\{|\xi|\le b\}$ and
      the support of $\varphi$ is contained in an annulus
      $C=\{a\le|\xi|\le b\}$ with $a=3/4$ and $b=8/3$;
    \item $\chi(\xi)+\displaystyle\sum_{j\ge 0}\varphi(2^{-j}\xi)=1$
      for all $\xi\in\R^{d}$;
    \item $\mathrm{supp}(\chi)\cap\mathrm{supp}(\varphi(2^{-i}\cdot))
      =\varnothing$ for $i\ge 1$, and
      $\mathrm{supp}(\varphi(2^{-i}\cdot))\cap
      \mathrm{supp}(\varphi(2^{-j}\cdot))=\varnothing$ whenever $|i-j|>1$.
  \end{enumerate}
  The existence of such $\chi,\varphi$ is
  classical~\cite[Proposition~2.10]{bahouri:11:fourier}.
  We write $\varphi_{j}:=\varphi(2^{-j}\cdot)$ for $j\ge 0$.
\end{notation}

\noindent
Those functions are interpreted as follows. The function $\chi$ concentrates
near zero frequency, while $\varphi_{j}$ selects the frequencies of order
$2^{j}$. Together they tile the frequency space in dyadic annuli, and they
allow to localize any distribution spectrally. The resulting building blocks
are the Littlewood--Paley projections.

\begin{notation}[Schwartz space and tempered distributions]\label{not:schwartz}
  We denote by $\mathcal{S}(\R^{d})$ the \emph{Schwartz space} of rapidly
  decreasing smooth functions on $\R^{d}$, equipped with the usual
  Fr\'echet topology generated by the seminorms
  \begin{equation*}
\|f\|_{\alpha,\beta}:=\sup_{x\in\R^{d}}|x^{\alpha}\partial^{\beta}f(x)|.
\end{equation*}
Its topological dual $\mathcal{S}'(\R^{d})$ is the space of
  \emph{tempered distributions}. The Fourier transform $\mathcal{F}$
  extends continuously from $\mathcal{S}(\R^{d})$ to $\mathcal{S}'(\R^{d})$,
  so that the operations in Definition~\ref{def:LP-blocks} below are
  well-defined for every tempered distribution.
\end{notation}

\begin{definition}[Littlewood--Paley blocks]\label{def:LP-blocks}
  Let $u\in\mathcal{S}'(\R^{d})$. The \emph{Littlewood--Paley blocks} are
  \begin{equation}\label{eq:LP-blocks}
    \Delta_{-1}u\,:=\,\mathcal{F}^{-1}(\chi\,\mathcal{F}u),
    \qquad
    \Delta_{j}u\,:=\,\mathcal{F}^{-1}(\varphi_{j}\,\mathcal{F}u)
    \quad\text{for }j\ge 0.
  \end{equation}
  The \emph{low-frequency projection} at level $k\ge 0$ is
  $S_{k}u:=\sum_{j=-1}^{k-1}\Delta_{j}u$.
\end{definition}

\noindent
Observe that each block $\Delta_{j}u$ is smooth: it coincides with the convolution
$K_{j}*u$, where 
\begin{equation*}
K_{j}(z):=2^{jd}(\mathcal{F}^{-1}\varphi)(2^{j}z)
\end{equation*}
is a Schwartz-class kernel. The spectral disjointness condition in
Notation~\ref{not:LP-partition}-(iii) implies that $\Delta_{i}u$ and
$\Delta_{j}u$ have disjoint Fourier supports whenever $|i-j|>1$. In
addition, owing to Notation~\ref{not:LP-partition}-(ii), the blocks sum to $u$ in $\mathcal{S}'(\R^{d})$:
\begin{equation}\label{eq:LP-sum}
  u = \sum_{j\ge -1}\Delta_{j}u \, ,
  \qquad\text{in }\mathcal{S}'(\R^{d}).
\end{equation}

\subsection{Besov--H\"{o}lder spaces}\label{ssec:besov-def}

The Littlewood--Paley decomposition $u=\sum_{j\ge -1}\Delta_{j}u$
(equation~\eqref{eq:LP-sum}) provides a canonical way to measure the
regularity of a distribution. Namely a function is smooth when its
high-frequency blocks $\Delta_{j}u$ are small, and it is rough when those
blocks grow. Quantifying this idea leads to the Besov--H\"{o}lder scale. We
introduce below three levels of Besov spaces, each of them being adapted to a
different aspect of our problem.
\begin{enumerate}
  \item \emph{Unweighted Besov--H\"{o}lder spaces}
    $\mathcal{B}^{\kappa}(\R^{d})$ (Definition~\ref{def:besov}): the
    baseline scale, built from the $L^{\infty}$ norms of the $L^{p}$ blocks.
    The exponent $\kappa\in\R$ is the regularity index; for
    $\kappa\in(0,1)$ these spaces coincide with H\"{o}lder spaces, while
    negative $\kappa$ captures rough distributions such as the noise
    $\dot{W}$.
  \item \emph{Weighted Besov--H\"{o}lder spaces}
    $\mathcal{B}^{\kappa}_{w}(\R^{d})$ (Definition~\ref{def:besov-weighted}):
    the $L^{\infty}$ norm of each block is replaced by
    $\|w\,\Delta_{j}f\|_{L^{\infty}}$, where $w$ is a weight function.
    The polynomial weight $\rho_{\sigma}$
    (Notation~\ref{not:poly-weight}) is the only instance needed here;
    it controls a moderate spatial growth of $f$ at infinity.
  \item \emph{$L^{p}$-weighted Besov spaces}
    $\mathcal{B}^{\kappa}_{p,w}(\R^{d})$ (Definition~\ref{def:besov-Lp}):
    the $L^{\infty}$ norm of each block is replaced by the weighted
    $L^{p}_{w}$ norm, and the sup over $j$ is replaced by a
    $\ell^{p}$-type sum.  The main instance is the space
    $\mathcal{A}_{q}=\mathcal{B}^{-\kappa'}_{2q,\rho_{\sigma}}(\R^{d})$,
    which appears as an intermediate space in the proof that the noise
    lies in $\mathcal{B}^{-\kappa}_{\rho_{\sigma}}$.
\end{enumerate}

Let us start with the unweighted case.

\begin{definition}[Besov--H\"{o}lder space]\label{def:besov}
  Let $\kappa\in\R$. The \emph{Besov--H\"{o}lder space}
  $\mathcal{B}^{\kappa}(\R^{d}):=\mathcal{B}^{\kappa}_{\infty,\infty}(\R^{d})$
  consists of all tempered distributions $f\in\mathcal{S}'(\R^{d})$ for which
  \begin{equation}\label{eq:def-besov}
    \|f\|_{\mathcal{B}^{\kappa}}
    \,:=\,\sup_{j\ge -1}\,2^{j\kappa}\,\|\Delta_{j}f\|_{L^{\infty}(\R^{d})}
    \,<\,\infty.
  \end{equation}
\end{definition}

\noindent
Recall that the exponent $\kappa$ in Definition~\ref{def:besov} is the
regularity index. A large positive $\kappa$ forces the high-frequency blocks
to be small, which corresponds to smooth functions. A negative $\kappa$
allows a polynomial growth of those blocks, and thus rough distributions. The
connection with classical function spaces is transparent for
$\kappa\in(0,1)$. In this range, $\mathcal{B}^{\kappa}(\R^{d})$ coincides
with the space of $\kappa$-H\"{o}lder continuous functions
$\mathcal{C}^{\kappa}(\R^{d})$~\cite[Theorem~2.36]{bahouri:11:fourier}. We
shall see in Section~\ref{sec:noise-besov} that the noise $\dot{W}$ belongs
almost surely to $\mathcal{B}^{\kappa}(\R^{d})$ for a proper negative
$\kappa$, which quantifies its roughness.

When the solution $u$ to the PAM equation is not bounded in space, one has to
keep track of its growth at infinity. This is achieved by inserting a weight
which penalizes large values of $|x|$. The weight we shall use in this book is
the polynomial decay function $\rho_{\sigma}$ below.

\begin{notation}[Polynomial weight]\label{not:poly-weight}
  For $\sigma>0$, the \emph{polynomial weight} $\rho_{\sigma}:\R^{d}\to(0,\infty)$
  is defined by
  \begin{equation}\label{eq:def-rho}
    \rho_{\sigma}(x) \,:=\, \bigl(1+|x|^{2}\bigr)^{-\sigma/2}.
  \end{equation}
  Thus $\rho_{\sigma}(x)\sim|x|^{-\sigma}$ as $|x|\to\infty$. In
  particular, $\rho_{\sigma}\in L^{p}(\R^{d})$ whenever $p\sigma>d$, so
  by choosing $\sigma$ large enough one can ensure any desired
  integrability at infinity.
\end{notation}

\noindent
With the weight $\rho_{\sigma}$ at hand, one obtains a weighted
analogue of Definition~\ref{def:besov} simply by replacing the
$L^{\infty}$ norm of each LP block by a weighted $L^{\infty}$ norm.

\begin{definition}[Weighted Besov--H\"{o}lder space]\label{def:besov-weighted}
  Let $w:\R^{d}\to(0,\infty)$ be a measurable weight and $\kappa\in\R$.
  The \emph{weighted Besov--H\"{o}lder space} $\mathcal{B}^{\kappa}_{w}(\R^{d})$
  consists of all $f\in\mathcal{S}'(\R^{d})$ (Notation~\ref{not:schwartz})
  for which
  \begin{equation}\label{eq:def-besov-weighted}
    \|f\|_{\mathcal{B}^{\kappa}_{w}}
    \,:=\,\sup_{j\ge -1}\,2^{j\kappa}\,
    \bigl\|w\,\Delta_{j}f\bigr\|_{L^{\infty}(\R^{d})}
    \,<\,\infty.
  \end{equation}
  When $w\equiv 1$ this reduces to Definition~\ref{def:besov}:
  $\mathcal{B}^{\kappa}_{1}(\R^{d})=\mathcal{B}^{\kappa}(\R^{d})$.
  The primary instance in this book is $w=\rho_{\sigma}$
  (Notation~\ref{not:poly-weight}), giving the space
  $\mathcal{B}^{\kappa}_{\rho_{\sigma}}(\R^{d})$ of distributions
  whose LP blocks satisfy $\sup_{j}2^{j\kappa}\|\rho_{\sigma}\Delta_{j}f\|_{L^{\infty}}<\infty$.
\end{definition}

\noindent
In the space $\mathcal{B}^{\kappa}_{\rho_{\sigma}}$, each block is measured
with an $L^{\infty}$ norm. Replacing it by an $L^{p}_{\rho_{\sigma}}$ norm and
summing over $j$, one gets a finer scale, better suited to moment estimates.
The instance we shall need below is $\mathcal{A}_{q}$, which is an
intermediate space in the proof that the noise lies in
$\mathcal{B}^{-\kappa}_{\rho_{\sigma}}$.

\begin{definition}[$L^{p}$-weighted Besov space and the space $\mathcal{A}_{q}$]%
  \label{def:besov-Lp}
  Let $w:\R^{d}\to(0,\infty)$, $p\ge 1$, and $\kappa\in\R$.
  For $g\in L^{1}_{\mathrm{loc}}(\R^{d})$ set
  \begin{equation}\label{eq:def-Lp-w}
    \|g\|_{L^{p}_{w}}^{p}
    \,:=\,\int_{\R^{d}}|g(x)|^{p}\,w(x)\,dx.
  \end{equation}
  The \emph{$L^{p}$-weighted Besov space} $\mathcal{B}^{\kappa}_{p,w}(\R^{d})$
  consists of all $f\in\mathcal{S}'(\R^{d})$ for which
  \begin{equation}\label{eq:def-besov-Lp}
    \|f\|_{\mathcal{B}^{\kappa}_{p,w}}^{p}
    \,:=\,\sum_{j\ge -1}\,2^{jp\kappa}\,
    \|\Delta_{j}f\|_{L^{p}_{w}(\R^{d})}^{p}
    \,<\,\infty.
  \end{equation}
  Given $\kappa'>0$, $\sigma>0$, and $q\ge 1$, we define
  \begin{equation}\label{eq:def-Aq}
    \mathcal{A}_{q}
    \,:=\,\mathcal{B}^{-\kappa'}_{2q,\,\rho_{\sigma}}(\R^{d}),
  \end{equation}
  so that
  $\|f\|_{\mathcal{A}_{q}}^{2q}
  =\sum_{j\ge -1}2^{-2qj\kappa'}\|\Delta_{j}f\|_{L^{2q}_{\rho_{\sigma}}}^{2q}$.
  By Besov embedding, $\mathcal{A}_{q}\hookrightarrow\mathcal{B}^{-\kappa}_{\rho_{\sigma}}$
  for any $\kappa>\kappa'$, provided $q$ is large enough (see
  \cite[Chapter~6]{triebel:83:theory}).
\end{definition}

Two properties of Besov--H\"{o}lder spaces will be used repeatedly. The
first one quantifies the gain of regularity when the heat semigroup $p_{t}*f$
is applied to a distribution of regularity $\kappa$. The second one is a
product estimate.

\begin{theorem}[Heat semigroup on Besov spaces;
  {\cite[Lemma~2.4]{hu:15:stochastic}}]\label{prop:heat-besov}
  Let $\kappa\in\R$ and $f\in\mathcal{B}^{\kappa}(\R^{d})$. Then for all
  $t>0$ and $\hat{\kappa}>\kappa$,
  \begin{equation}\label{eq:heat-besov}
    \|p_{t}*f\|_{\mathcal{B}^{\hat{\kappa}}}
    \,\le\,c\,t^{-\frac{\hat{\kappa}-\kappa}{2}}\,\|f\|_{\mathcal{B}^{\kappa}},
  \end{equation}
  where $c>0$ depends only on $d$, $\kappa$, and $\hat{\kappa}$.
\end{theorem}

\begin{proposition}[Product in Besov spaces;
  {\cite[Theorem~2.52]{bahouri:11:fourier}}]\label{prop:product-besov}
  Let $\kappa_{1},\kappa_{2}\in\R$ with $\kappa_{2}<0<\kappa_{1}$ and
  $\kappa_{1}>|\kappa_{2}|$. Then the pointwise product extends to a
  continuous bilinear map
  \begin{equation}\label{eq:product-besov}
    \mathcal{B}^{\kappa_{1}}(\R^{d})\times\mathcal{B}^{\kappa_{2}}(\R^{d})
    \longrightarrow\mathcal{B}^{\kappa_{2}}(\R^{d}),
  \end{equation}
  with bound
  $\|f_{1}f_{2}\|_{\mathcal{B}^{\kappa_{2}}}
  \le c\,\|f_{1}\|_{\mathcal{B}^{\kappa_{1}}}\|f_{2}\|_{\mathcal{B}^{\kappa_{2}}}$.
\end{proposition}

\noindent
Proposition~\ref{prop:heat-besov} is one of the pillars of the existence and
uniqueness argument in Chapter~\ref{chap:existence}. It asserts that a
convolution with the heat kernel gains $(\hat{\kappa}-\kappa)/2$ derivatives,
at the cost of a power of $t$. As for Proposition~\ref{prop:product-besov}, it
will be invoked in order to give a meaning to the product $u\,\dot{W}$ in the
PAM equation, when $u$ is smooth and $\dot{W}$ is a rough distribution.

\section[Noises]{Definition of our classes of noises}\label{sec:noise-def}

The framework for the noise has been set up in
Subsection~\ref{ssec:noise}, that is the covariance
structure~\eqref{eq:noise-cov}, the spectral measure $\mu=\mathcal{F}\Lambda$,
the Hilbert space $\mathcal{H}$ with inner
product~\eqref{eq:H-innprod-fourier}, the isometry~\eqref{eq:W-isometry} and
Dalang's condition~\eqref{eq:dalang-cond}. In the present section we collect
concrete examples of spatial correlations $\Lambda$. For each of them we
specify $\Lambda(x)$ and the spectral measure $\mu(d\xi)$, and we check
Assumption~\ref{hyp:spectral} in a short lemma.

\begin{example}[Riesz noise]
  \label{ex:riesz}
  For $\eta\in(0,d)$, the \emph{Riesz correlation} is
  \begin{equation}\label{eq:riesz-corr}
    \Lambda_{\eta}(x)\,:=\,|x|^{-\eta},
    \qquad x\in\R^{d}\setminus\{0\}.
  \end{equation}
  It is nonneg-definite and locally integrable for $\eta\in(0,d)$.
  Its Fourier transform is
  \begin{equation}\label{eq:riesz-spectral}
    \mu_{\eta}(d\xi)
    \,:=\,c_{\eta,d}\,|\xi|^{-(d-\eta)}\,d\xi,
    \qquad
    c_{\eta,d}
    =\frac{2^{d-\eta}\pi^{d/2}\,
    \Gamma\!\left(\tfrac{d-\eta}{2}\right)}
    {\Gamma\!\left(\tfrac{\eta}{2}\right)}.
  \end{equation}
\end{example}

\begin{remark}\label{rmk:riesz}
The spectral density $|\xi|^{-(d-\eta)}$ is homogeneous of degree
  $\eta-d$. Also note that this spectral density is singular at the origin and non-integrable at
  infinity for $\eta\ge 0$.
\end{remark}

\noindent
The threshold $\eta<2$ is dictated by the high-frequency behavior of
$\mu_{\eta}$. Indeed, the singularity at the origin is always integrable. In
contrast, the power-law decay at infinity is not integrable against
$(1+|\xi|^{2})^{-1}$ once $\eta\ge 2$.

\begin{lemma}[Dalang's condition for Riesz noise]\label{lem:dalang-riesz}
  Assumption~\ref{hyp:spectral} holds for the noise of
  Example~\ref{ex:riesz} if and only if $\eta<2$.
\end{lemma}
\begin{proof}
  Inserting~\eqref{eq:riesz-spectral} into~\eqref{eq:dalang-cond}, the
  condition to verify is
  \begin{equation}\label{eq:dalang-riesz-integral}
    \int_{\R^{d}}\frac{c_{\eta,d}\,|\xi|^{-(d-\eta)}}{1+|\xi|^{2}}\,d\xi
    \,=\,c_{\eta,d}\left(
      \int_{|\xi|\le 1}\frac{|\xi|^{-(d-\eta)}}{1+|\xi|^{2}}\,d\xi
      +\int_{|\xi|>1}\frac{|\xi|^{-(d-\eta)}}{1+|\xi|^{2}}\,d\xi
    \right)<\infty.
  \end{equation}
  Each piece in the right hand side of \eqref{eq:dalang-riesz-integral} is handled separately.

  \smallskip
  \noindent
  {\it Low frequencies $|\xi|\le 1$.}
  Since $1+|\xi|^{2}\ge 1$, one has
  $\int_{|\xi|\le 1}|\xi|^{-(d-\eta)}\,d\xi
  \sim\int_{0}^{1}r^{\eta-1}\,dr<\infty$ for all $\eta>0$.

  \smallskip
  \noindent
  {\it High frequencies $|\xi|>1$.}
  Here $(1+|\xi|^{2})^{-1}\sim|\xi|^{-2}$, so
  \begin{equation*}
    \int_{|\xi|>1}|\xi|^{-(d-\eta)-2}\,d\xi
    \,\sim\,\int_{1}^{\infty}r^{\eta-3}\,dr,
  \end{equation*}
  which converges if and only if $\eta<2$.
  Combining both pieces, the integral in~\eqref{eq:dalang-riesz-integral}
  is finite if and only if $\eta<2$, which completes the proof.
\end{proof}

Before stating our second example, we recall how the operator
$(I-\Delta)^{-\eta/2}$ can be written as an integral of the heat
semigroup. We denote by $P_{t}=e^{t\frac{1}{2}\Delta}$ the heat semigroup
generated by $\tfrac{1}{2}\Delta$ on $\R^{d}$ (the operator driving the
PAM~\eqref{eq:pam}), whose integral kernel is the Gaussian heat
kernel~$p_{t}$ introduced in~\eqref{eq:dalang-cond}'s discussion, namely
\begin{equation}\label{eq:heat-kernel-explicit}
  p_{t}(x)\,=\,(2\pi t)^{-d/2}\, \exp\left(-\frac{|x|^{2}}{2t}\right),
  \qquad t>0,\ x\in\R^{d},
\end{equation}
with spatial Fourier transform $\widehat{p_{t}}(\xi)=e^{-t|\xi|^{2}/2}$. The
full Laplacian entering the Bessel operator generates $e^{t\Delta}=P_{2t}$,
whose Gaussian kernel is
$p_{2t}(x)=(4\pi t)^{-d/2}\,e^{-|x|^{2}/(4t)}$.
Our starting point is the elementary Laplace transform identity for the
Gamma function. Namely, for every $\lambda>0$ and $\eta>0$ we have
\begin{equation}\label{eq:gamma-laplace}
  \lambda^{-\eta/2}
  \,=\,\frac{1}{\Gamma(\eta/2)}
  \int_{0}^{\infty}t^{\eta/2-1}\,e^{-t\lambda}\,dt.
\end{equation}
Applying~\eqref{eq:gamma-laplace} through the spectral calculus with
$\lambda$ replaced by the (positive, self-adjoint) operator $I-\Delta$,
and using $e^{-t(I-\Delta)}=e^{-t}\,P_{2t}$, we obtain the subordination
formula
\begin{equation}\label{eq:bessel-subord}
  (I-\Delta)^{-\eta/2}
  \,=\,\frac{1}{\Gamma(\eta/2)}
  \int_{0}^{\infty}t^{\eta/2-1}\,e^{-t}\,P_{2t}\,dt.
\end{equation}
Reading the integral kernel of both sides of~\eqref{eq:bessel-subord},
thanks to the kernel $p_{2t}$ of $P_{2t}$ obtained above, one gets the kernel
$\Lambda^{\mathrm{B}}_{\eta}$ of the example below.

\begin{example}[Bessel noise.]
  \label{ex:bessel}
  For $\eta>0$, the \emph{Bessel correlation}
  $\Lambda^{\mathrm{B}}_{\eta}$ is the integral kernel of the operator
  $(I-\Delta)^{-\eta/2}$:
  \begin{equation}\label{eq:bessel-corr}
    \Lambda^{\mathrm{B}}_{\eta}(x)
    =c\int_{0}^{\infty}w^{(\eta-d)/2-1}\,e^{-w}\,
    e^{-\frac{|x|^{2}}{4w}}\,dw,
    \qquad x\in\R^{d},
  \end{equation}
  with normalizing constant $c=(4\pi)^{-d/2}\Gamma(\eta/2)^{-1}$, as
  follows from the subordination formula~\eqref{eq:bessel-subord}.
  Its Fourier transform is
  \begin{equation}\label{eq:bessel-spectral}
    \mu^{\mathrm{B}}_{\eta}(d\xi)
    \,:=\,(1+|\xi|^{2})^{-\eta/2}\,d\xi.
  \end{equation}
  For $0<\eta<d$, the kernel $\Lambda^{\mathrm{B}}_{\eta}$ has the same
  Riesz-type singularity $|x|^{\eta-d}$ at the origin, but (unlike the
  Riesz correlation) it decays exponentially at infinity (see~\cite[Section~1.2.1]{grafakos:14:modern}):
  \begin{equation}\label{eq:bessel-asymp}
    \Lambda^{\mathrm{B}}_{\eta}(x)\,\asymp\,
    \begin{cases}
      |x|^{\eta-d} & \text{as }|x|\to 0,\\
      |x|^{(\eta-d-1)/2}\,e^{-|x|} & \text{as }|x|\to\infty.
    \end{cases}
  \end{equation}
  For $\eta\ge d$ the kernel is bounded near the origin, the
  singularity in~\eqref{eq:bessel-asymp} being replaced by a finite
  value.  The spectral density $(1+|\xi|^{2})^{-\eta/2}$ has
  polynomial decay for large $|\xi|$, and no singularity at $\xi=0$.
\end{example}

\begin{remark}
Let us mention an interesting duality. The Bessel kernel
  $\Lambda^{\mathrm{B}}_{\eta}$ may be used either as the correlation
  $\Lambda$ (with spectral density~\eqref{eq:bessel-spectral}), or conversely
  as the spectral density $\hat{\Lambda}(\xi)\,d\xi$. In the latter case the
  correlation function is $(1+|x|^{2})^{-\eta/2}$, a polynomial-decay
  function. Both choices lead to well-defined noises, with different
  conditions on $\eta$.

\end{remark}

\begin{lemma}[Dalang's condition for Bessel noise]\label{lem:dalang-bessel}
  Assumption~\ref{hyp:spectral} holds for the noise of
  Example~\ref{ex:bessel} if and only if $\eta>d-2$.
\end{lemma}
\begin{proof}
  Substituting~\eqref{eq:bessel-spectral} into~\eqref{eq:dalang-cond}, we have to check the finiteness of the following integral:
  \begin{equation*}
    \int_{\R^{d}}(1+|\xi|^{2})^{-(\eta+2)/2}\,d\xi
    \,\sim\,\int_{0}^{\infty}r^{d-1}(1+r^{2})^{-(\eta+2)/2}\,dr.
  \end{equation*}
  The integrand near $r=0$ is $O(1)$, hence integrable.  For $r\to\infty$
  it behaves like $r^{d-\eta-3}$, and the integral converges if and only
  if $d-\eta-3<-1$, i.e., $\eta>d-2$.
\end{proof}

The two previous examples isolate the two features of a correlation
function that matter for the PAM, namely its singularity at the origin
and its decay at infinity. Each example controls these features by
a single parameter. On the one hand, the Bessel kernel of
Example~\ref{ex:bessel} produces, through~\eqref{eq:bessel-asymp}, a
prescribed power singularity $|x|^{\eta-d}$ at the origin (for
$0<\eta<d$), together with a rapid (exponential) decay at infinity. On
the other hand, the polynomial-decay function $(1+|x|^{2})^{-s/2}$ is
bounded near the origin and behaves like $|x|^{-s}$ as
$|x|\to\infty$. Adding one of each type therefore lets us tune the two
exponents independently. The Bessel term fixes the order of the singularity
at $0$, while the polynomial term fixes the power of the decay at $\infty$.
This is the content of our next example, in which any pair of power behaviors
at $0$ and at $\infty$ can be realized.

\begin{example}[Riesz-type noise]
  \label{ex:riesz-type}
  For $s_{1},s_{2}\in(0,d)$, the \emph{Riesz-type correlation with
  parameters $(s_{1},s_{2})$} is defined by
  \begin{equation}\label{eq:riesz-type-corr}
    \Lambda_{s_{1},s_{2}}(x)
    \,:=\,\Lambda^{\mathrm{B}}_{s_{1}}(x)+(1+|x|^{2})^{-s_{2}/2},
  \end{equation}
  combining the Bessel kernel of Example~\ref{ex:bessel} and a
  polynomial-decay tail.  Its spectral density is
  \begin{equation}\label{eq:riesz-type-spectral}
    \hat{\Lambda}_{s_{1},s_{2}}(\xi)
    \,:=\,(1+|\xi|^{2})^{-s_{1}/2}+\Lambda^{\mathrm{B}}_{s_{2}}(\xi),
  \end{equation}
  where $\Lambda^{\mathrm{B}}_{s_{2}}(\xi)$ denotes the Bessel kernel
  of Example~\ref{ex:bessel} evaluated at $\xi$.
  The kernel $\Lambda_{s_{1},s_{2}}$ interpolates between
  Riesz-type behavior near the origin and polynomial decay at infinity:
  \begin{equation}\label{eq:riesz-type-asymp}
    \Lambda_{s_{1},s_{2}}(x)\,\asymp\,
    \begin{cases}
      |x|^{s_{1}-d} & \text{as }|x|\to 0,\\
      |x|^{-s_{2}} & \text{as }|x|\to\infty.
    \end{cases}
  \end{equation}
  The parameter $s_{1}$ controls the singularity at the origin, and
  $s_{2}$ controls the spatial decay at infinity.
\end{example}

\begin{lemma}[Dalang's condition for Riesz-type noise]%
  \label{lem:dalang-riesz-type}
  For any $s_{2}\in(0,d)$, Assumption~\ref{hyp:spectral} holds for the
  noise of Example~\ref{ex:riesz-type} if and only if $s_{1}>d-2$.
\end{lemma}
\begin{proof}
The spectral measure $\mu_{s_{1},s_{2}}$ of $\Lambda_{s_{1},s_{2}}$ satisfies
\begin{equation*}
\mu_{s_{1},s_{2}}=\mu^{\mathrm{B}}_{s_{1}}+\nu ,
\quad\text{where}\quad
\nu(d\xi)=\Lambda^{\mathrm{B}}_{s_{2}}(\xi)\,d\xi
\end{equation*}
The contribution
  of $\mu^{\mathrm{B}}_{s_{1}}$ is handled by
  Lemma~\ref{lem:dalang-bessel}: it satisfies~\eqref{eq:dalang-cond} iff
  $s_{1}>d-2$.  For $\nu$, we split at $|\xi|=1$ in the following way:

  \smallskip
  \noindent
  {\it Low frequencies.}  Near $\xi=0$ the Bessel kernel satisfies
  $\Lambda^{\mathrm{B}}_{s_{2}}(\xi)\asymp|\xi|^{s_{2}-d}$ according to \eqref{eq:bessel-asymp}, so that
  $\int_{|\xi|\le 1}\Lambda^{\mathrm{B}}_{s_{2}}(\xi)\,d\xi
  \sim\int_{0}^{1}r^{s_{2}-1}\,dr<\infty$ for $s_{2}>0$.

  \smallskip
  \noindent
  {\it High frequencies.}  The Bessel kernel decays exponentially:
  $\Lambda^{\mathrm{B}}_{s_{2}}(\xi)\lesssim e^{-c|\xi|}$ for
  $|\xi|>1$, so $\int_{|\xi|>1}\Lambda^{\mathrm{B}}_{s_{2}}(\xi)/(1+|\xi|^{2})\,d\xi<\infty$.

  \smallskip
  \noindent
  Hence $\int_{\R^{d}}\nu(d\xi)/(1+|\xi|^{2})<\infty$ for all
  $s_{2}\in(0,d)$, and the condition reduces entirely to $s_{1}>d-2$.
\end{proof}

All the correlations considered so far (the Riesz kernel of
Example~\ref{ex:riesz}, the Bessel kernel of Example~\ref{ex:bessel},
and the Riesz-type kernel of Example~\ref{ex:riesz-type}) are
isotropic. Namely they depend on $x$ through the Euclidean norm $|x|$ only,
and they treat all the directions of $\R^{d}$ in the same way. Our last
example breaks this symmetry. It is built as a tensor product over the $d$
coordinates, with a Hurst parameter $H_{i}$ attached to each direction. The
roughness of the noise is thus allowed to differ from one coordinate axis to
the other. 

\begin{example}[Fractional spatial noise]
  \label{ex:frac-space}
  Let $H_{i}\in\bigl(\tfrac{1}{2},1\bigr)$ for $i=1,\ldots,d$.  The
  \emph{fractional spatial correlation} is
  \begin{equation}\label{eq:frac-space-corr}
    \Lambda_{H}(x)
    \,:=\,\prod_{i=1}^{d}H_{i}(2H_{i}-1)\,|x_{i}|^{2H_{i}-2},
    \qquad x=(x_{1},\ldots,x_{d})\in\R^{d}\setminus\{0\}.
  \end{equation}
  Each factor is (up to a constant) the covariance kernel of a
  one-dimensional fractional Brownian motion with Hurst parameter
  $H_{i}$.  The spectral measure is
  \begin{equation}\label{eq:frac-space-spectral}
    \mu_{H}(d\xi)
    \,:=\,C_{H}\prod_{i=1}^{d}|\xi_{i}|^{1-2H_{i}}\,d\xi,
  \end{equation}
  where $C_{H}>0$ is a constant depending on the $H_{i}$.  The
  spectral density is singular at each coordinate hyperplane
  $\{\xi_{i}=0\}$; the condition $H_{i}>1/2$ ensures that this
  singularity is integrable.
\end{example}

\begin{lemma}[Dalang's condition for fractional spatial noise]%
  \label{lem:dalang-frac-space}
  Assumption~\ref{hyp:spectral} holds for the noise of
  Example~\ref{ex:frac-space} if and only if
  $\sum_{i=1}^{d}H_{i}>d-1$.
\end{lemma}
\begin{proof}
  Set $\bar{H}:=\sum_{i=1}^{d}H_{i}$. We divide again the verification of \eqref{eq:dalang-cond} in two steps.

  \smallskip
  \noindent
  {\it Near the origin.}  Since $1-2H_{i}\in(-1,0)$, each factor
  $|\xi_{i}|^{1-2H_{i}}$ is integrable near $\xi_{i}=0$.  Evaluating
  separately in each coordinate, we get
  \begin{equation*}
\int_{[-1,1]^{d}} \mu_{H}(d\xi)
=\int_{[-1,1]^{d}}\prod_{i}|\xi_{i}|^{1-2H_{i}}\,d\xi
  =\prod_{i}\int_{-1}^{1}|\xi_{i}|^{1-2H_{i}}\,d\xi_{i}
<\infty .
\end{equation*}

  \smallskip
  \noindent
  {\it Near infinity.}  The spectral density
  $\prod_{i}|\xi_{i}|^{1-2H_{i}}$ is homogeneous of degree
  $\sum_{i}(1-2H_{i})=d-2\bar{H}$.  Passing to polar coordinates
  $\xi=r\omega$ with $r=|\xi|$ and $\omega\in S^{d-1}$, so that
  $d\xi=r^{d-1}\,dr\,d\sigma(\omega)$, we obtain
  \begin{equation*}
    \prod_{i}|\xi_{i}|^{1-2H_{i}}
    =r^{d-2\bar{H}}\prod_{i}|\omega_{i}|^{1-2H_{i}}.
  \end{equation*}
  The angular integral is finite, that is
  \begin{equation*}
    \int_{S^{d-1}}\prod_{i}|\omega_{i}|^{1-2H_{i}}\,d\sigma(\omega)
    <\infty,
  \end{equation*}
  because each exponent satisfies $1-2H_{i}>-1$. This makes the
  singularities on the coordinate hyperplanes $\{\omega_{i}=0\}$
  integrable, exactly as in the previous step.  Hence
  \begin{equation*}
    \int_{|\xi|>1}\frac{\prod_{i}|\xi_{i}|^{1-2H_{i}}}{1+|\xi|^{2}}\,d\xi
    \,\sim\,\int_{1}^{\infty}\frac{r^{d-2\bar{H}}}{1+r^{2}}\,r^{d-1}\,dr
    \,\sim\,\int_{1}^{\infty}r^{2d-2\bar{H}-3}\,dr,
  \end{equation*}
  which converges if and only if $2d-2\bar{H}-3<-1$. We thus obtain the desired condition
  $\bar{H}>d-1$.
\end{proof}

\begin{remark}[Further examples]\label{rem:noise-further}
The examples below are widely considered in the literature, although they will play a more modest role in this book.
\begin{enumerate}
\item
\emph{White noise in space.}  
Setting $\Lambda=\delta_{0}$ gives
  $\mu=\text{Lebesgue}$, and~\eqref{eq:dalang-cond} reduces to
  $\int_{\R^{d}}(1+|\xi|^{2})^{-1}\,d\xi<\infty$, which holds if
  and only if $d=1$.  Hence space-time white noise is only admissible
  in one spatial dimension.

\item
\emph{Fractional noise in time.}  
The framework extends to include time
  correlations by replacing the white-noise factor $dt$
  in~\eqref{eq:noise-cov} with a convolution kernel $\gamma(s-t)\,ds\,dt$,
  where $\gamma(t)=H_{0}(2H_{0}-1)|t|^{2H_{0}-2}$, $H_{0}\in(\tfrac{1}{2},1)$,
  is the covariance of a one-dimensional fractional Brownian motion with
  Hurst parameter $H_{0}$.  The full covariance is then
  \begin{equation}\label{eq:frac-time-cov}
    \mathbb{E}\bigl[W(\varphi)\,W(\psi)\bigr]
    =\int_{\R_{+}^{2}}\!\int_{\R^{2d}}
    \varphi(s,x)\,\psi(t,y)\,\gamma(s-t)\,\Lambda(x-y)\,
    dx\,dy\,ds\,dt,
  \end{equation}
  and Dalang's condition acquires an additional integrability factor
  from $\hat\gamma$; we refer to~\cite[Section~2.1]{hu:15:stochastic} for
  details.  The four spatial examples above all extend naturally to this
  setting.
\end{enumerate}
 \end{remark}

\section{Inclusion in Besov spaces}\label{sec:noise-besov}

We are now ready to make precise the statement, announced after
Definition~\ref{def:besov}, that the driving noise $\dot{W}$ belongs almost
surely to a Besov space with negative regularity. Since $W$ is white in time
(see~\eqref{eq:noise-cov}), the natural object here is not a single
distribution. We shall rather consider the whole trajectory $t\mapsto W_{t}$,
where $W_{t}$ is the spatial distribution obtained by integrating $W$ over the
time window $[0,t]$. The main result of this section is
Theorem~\ref{thm:noise-besov} below. It states that this trajectory is
H\"{o}lder continuous in time, with values in a weighted Besov space
$\mathcal{B}^{-\kappa}_{\rho_{\sigma}}$
(Definition~\ref{def:besov-weighted}). This type of regularity
is exploited when solving stochastic PDEs in the Stratonovich sense,
but will not be directly featured in the It\^o setting adopted in this book.
We still find it a valuable information, since it gives a precise idea of the noise singularity.
The result is taken from~\cite[Proposition~5.5]{hu:15:stochastic}, specialized here to the
white-in-time setting~\eqref{eq:noise-cov}.
We first record the distribution-valued process associated with $W$.

\begin{notation}[The noise as a distribution-valued process]\label{not:Wt}
  For $t\ge 0$, let $W_{t}$ denote the random tempered distribution
  acting on test functions $\varphi\in\mathcal{D}(\R^{d})$ by
  \begin{equation}\label{eq:def-Wt}
    W_{t}(\varphi)\,:=\,W\bigl(\mathbf{1}_{[0,t]}\otimes\varphi\bigr),
  \end{equation}
  with $W$ the Gaussian family of
  Subsection~\ref{ssec:noise}. For $0\le s\le t$ we write the increment
  \begin{equation}\label{eq:def-incr-W}
    \delta W_{st}\,:=\,W_{t}-W_{s},
    \qquad\text{so that}\qquad
    \delta W_{st}(\varphi)=W\bigl(\mathbf{1}_{[s,t]}\otimes\varphi\bigr).
  \end{equation}
\end{notation}

\noindent
In order to measure the joint space-time regularity of $t\mapsto W_{t}$, we
combine the weighted Besov norm of Definition~\ref{def:besov-weighted} with a
H\"{o}lder norm in time.

\begin{definition}[Time-H\"{o}lder, Besov-valued space]\label{def:holder-besov}
  Fix $T>0$, $\theta\in(0,1)$, $\kappa\in\R$, and $\sigma>0$, and let
  $\rho_{\sigma}$ be the polynomial weight of
  Notation~\ref{not:poly-weight}. The space
  $\mathcal{C}^{\theta,-\kappa}_{\rho_{\sigma}}([0,T])$ consists of all
  maps $f:[0,T]\to\mathcal{B}^{-\kappa}_{\rho_{\sigma}}(\R^{d})$ for which
  \begin{equation}\label{eq:def-holder-besov}
    \|f\|_{\mathcal{C}^{\theta,-\kappa}_{\rho_{\sigma}}}
    \,:=\,\sup_{0\le t\le T}\|f_{t}\|_{\mathcal{B}^{-\kappa}_{\rho_{\sigma}}}
    \,+\,\sup_{0\le s<t\le T}
    \frac{\|f_{t}-f_{s}\|_{\mathcal{B}^{-\kappa}_{\rho_{\sigma}}}}{(t-s)^{\theta}}
    \,<\,\infty.
  \end{equation}
\end{definition}

\noindent
The regularity of $W_{t}$ in space is dictated by the integrability of the
spectral measure $\mu$ at high frequencies. Recall that
Assumption~\ref{hyp:spectral} (Dalang's condition) only guarantees
$\int_{\R^{d}}(1+|\xi|^{2})^{-1}\,\mu(d\xi)<\infty$. In order to get a
genuine H\"{o}lder regularity in the Besov scale, a small amount of extra
integrability is needed. It is quantified by a parameter $\alpha$ in the
following assumption.

\begin{assumption}\label{hyp:spectral-holder}
  There exists $\alpha\in(0,1)$ such that
  \begin{equation}\label{eq:dalang-holder}
    \int_{\R^{d}}\frac{\mu(d\xi)}{1+|\xi|^{2(1-\alpha)}}\,<\,\infty.
  \end{equation}
\end{assumption}

\noindent
Since $2(1-\alpha)<2$, Assumption~\ref{hyp:spectral-holder} is stronger than
Assumption~\ref{hyp:spectral}, and it reduces to it in the limit $\alpha\to
0$. The gap between the two conditions is the price to pay for a pathwise
regularity of the noise, instead of a mere $L^{2}$ regularity. For our running
examples of Section~\ref{sec:noise-def}, condition~\eqref{eq:dalang-holder} is
easily checked.
\begin{enumerate}
\item For the Riesz noise of Example~\ref{ex:riesz}, with
  $\mu(d\xi)\asymp|\xi|^{-(d-\eta)}\,d\xi$ by~\eqref{eq:riesz-spectral},
  condition~\eqref{eq:dalang-holder} holds for every
  $\alpha<1-\tfrac{\eta}{2}$, hence (for some $\alpha>0$) under the same
  restriction $\eta<2$ as Dalang's condition.

\item For the Bessel noise of Example~\ref{ex:bessel}, with
  $\mu^{\mathrm{B}}_{\eta}(d\xi)=(1+|\xi|^{2})^{-\eta/2}\,d\xi$
  by~\eqref{eq:bessel-spectral},
  condition~\eqref{eq:dalang-holder} holds for every
  $\alpha<1-\tfrac{d-\eta}{2}$, i.e. whenever $\eta>d-2$.
\end{enumerate}

We can now state the main result of this section.

\begin{theorem}[Pathwise Besov regularity of the white-in-time noise.]\label{thm:noise-besov}
  Let $W$ be the white-in-time Gaussian noise with
  covariance~\eqref{eq:noise-cov} and spectral measure $\mu$, and suppose
  that Assumption~\ref{hyp:spectral-holder} holds for some
  $\alpha\in(0,1)$. Fix a horizon $T>0$ and a weight exponent
  $\sigma>d$. Then for every $\theta\in\bigl(0,\tfrac{1}{2}\bigr)$ and
  every $\kappa>1-\alpha$, the process $t\mapsto W_{t}$ of
  Notation~\ref{not:Wt} admits a modification belonging almost surely to
  $\mathcal{C}^{\theta,-\kappa}_{\rho_{\sigma}}([0,T])$. Moreover the
  random variable $\|W\|_{\mathcal{C}^{\theta,-\kappa}_{\rho_{\sigma}}}$
  has finite moments of all orders.
\end{theorem}

\noindent
The proof relies on a single second moment computation. The white-in-time
structure turns the time average into a factor $(t-s)$, while the spatial
integrability~\eqref{eq:dalang-holder} controls the growth of the
Littlewood--Paley blocks. The intermediate space~$\mathcal{A}_{q}$ of
Definition~\ref{def:besov-Lp} then converts those moment bounds into a
pathwise H\"{o}lder regularity, thanks to the Garsia--Rodemich--Rumsey
lemma~\cite{garsia:70:real}.

\begin{proof}
  Fix $\kappa>\kappa'>1-\alpha$ and recall from~\eqref{eq:def-Aq} the
  space $\mathcal{A}_{q}=\mathcal{B}^{-\kappa'}_{2q,\,\rho_{\sigma}}(\R^{d})$,
  whose norm is
  \begin{equation}\label{eq:Aq-norm-recall}
    \|f\|_{\mathcal{A}_{q}}^{2q}
    =\sum_{j\ge -1}2^{-2qj\kappa'}\,\|\Delta_{j}f\|_{L^{2q}_{\rho_{\sigma}}}^{2q}
    =\sum_{j\ge -1}2^{-2qj\kappa'}
    \int_{\R^{d}}|\Delta_{j}f(x)|^{2q}\,\rho_{\sigma}(x)\,dx,
  \end{equation}
  the last equality being the definition~\eqref{eq:def-Lp-w} of the
  weighted $L^{2q}$ norm. We will prove that, for $q$ large enough,
  \begin{equation}\label{eq:as-bnd-increments-W}
    \E\Bigl[\|\delta W_{st}\|_{\mathcal{A}_{q}}^{2q}\Bigr]
    \,\le\,c_{q}\,(t-s)^{q},
    \qquad 0\le s\le t\le T.
  \end{equation}
  Starting from~\eqref{eq:as-bnd-increments-W}, the proof is concluded as
  follows. The Besov embedding recorded after~\eqref{eq:def-Aq} gives
  $\mathcal{A}_{q}\hookrightarrow\mathcal{B}^{-\kappa}_{\rho_{\sigma}}$
  for $q$ large, so~\eqref{eq:as-bnd-increments-W} yields
  $$
  \E[\|\delta W_{st}\|_{\mathcal{B}^{-\kappa}_{\rho_{\sigma}}}^{2q}]\le
  c_{q}(t-s)^{q}.
  $$ 
  Hence by the Garsia--Rodemich--Rumsey
  lemma~\cite{garsia:70:real}, a moment bound of the form above
implies that $t\mapsto
  W_{t}$ has a modification which is H\"{o}lder continuous of any order
  $\theta<\tfrac{1}{2}-\tfrac{1}{2q}$ with values in
  $\mathcal{B}^{-\kappa}_{\rho_{\sigma}}$, and that the corresponding
  H\"{o}lder norm has moments of all orders. Letting $q\to\infty$ covers
  every $\theta<\tfrac{1}{2}$ and every $\kappa>1-\alpha$, which is the
  assertion of the theorem. It remains to
  prove~\eqref{eq:as-bnd-increments-W}, which we do in three steps.

  \smallskip

  \noindent
  {\it Step 1: Reduction to a second moment.}  By
  Notation~\ref{not:Wt}, $\Delta_{j}(\delta W_{st})$ is the random field
  obtained by applying the Littlewood--Paley block $\Delta_{j}$
  (Definition~\ref{def:LP-blocks}) to the spatial distribution
  $\delta W_{st}$. Writing $\Delta_{j}f=K_{j}*f$ and setting
  $K_{j,x}(y):=K_{j}(x-y)$, we have
  \begin{equation*}
\Delta_{j}(\delta W_{st})(x)=\delta W_{st}(K_{j,x})
  =W(\mathbf{1}_{[s,t]}\otimes K_{j,x}) ,
\end{equation*}
by~\eqref{eq:def-incr-W}.
  Each such quantity is a centered Gaussian random variable, so its even
  moments are controlled by its variance,
  \begin{equation}\label{eq:gauss-moments}
    \E\Bigl[\bigl|W(\mathbf{1}_{[s,t]}\otimes K_{j,x})\bigr|^{2q}\Bigr]
    \,=\,(2q-1)!! \
    \E\Bigl[\bigl|W(\mathbf{1}_{[s,t]}\otimes K_{j,x})\bigr|^{2}\Bigr]^{q}.
  \end{equation}
  Taking expectations in~\eqref{eq:Aq-norm-recall} applied to
  $f=\delta W_{st}$ and inserting~\eqref{eq:gauss-moments}, we obtain
  \begin{equation}\label{eq:bnd-EW-Aq}
    \E\Bigl[\|\delta W_{st}\|_{\mathcal{A}_{q}}^{2q}\Bigr]
    \,=\, c_{q}
    \sum_{j\ge -1}2^{-2qj\kappa'}
    \int_{\R^{d}}
    \E\Bigl[\bigl|W(\mathbf{1}_{[s,t]}\otimes K_{j,x})\bigr|^{2}\Bigr]^{q}
    \rho_{\sigma}(x)\,dx.
  \end{equation}

  \smallskip

  \noindent
  {\it Step 2: The second moment.}  Using the spectral form
  \eqref{eq:H-innprod-fourier} of the covariance and the white-in-time
  structure~\eqref{eq:noise-cov}, the time integration reduces to the
  diagonal and produces a factor $(t-s)$:
  \begin{equation}\label{eq:second-moment-W}
    \E\Bigl[\bigl|W(\mathbf{1}_{[s,t]}\otimes K_{j,x})\bigr|^{2}\Bigr]
    =\int_{s}^{t}\!\int_{\R^{d}}\bigl|\mathcal{F}K_{j,x}(\xi)\bigr|^{2}\,
    \mu(d\xi)\,dr
    =(t-s)\int_{\R^{d}}\bigl|\varphi_{j}(\xi)\bigr|^{2}\,\mu(d\xi),
  \end{equation}
  where we used $|\mathcal{F}K_{j,x}(\xi)|=|\varphi_{j}(\xi)|$ (the
  translation $x$ only contributes a unit-modulus phase), with
  $\varphi_{j}=\varphi(2^{-j}\cdot)$ from
  Notation~\ref{not:LP-partition} for $j\ge 0$. In~\eqref{eq:second-moment-W}, we also use 
  $\varphi_{-1}=\chi$.
  In particular the right-hand side of~\eqref{eq:second-moment-W} does
  not depend on $x$. To estimate it, introduce the finite measure
  \begin{equation}\label{eq:def-nu}
    \nu(d\xi)\,:=\,\frac{\mu(d\xi)}{1+|\xi|^{2(1-\alpha)}},
  \end{equation}
  which is finite by Assumption~\ref{hyp:spectral-holder}. For $j\ge 0$
  the function $\varphi_{j}$ is supported in the annulus
  $\{2^{j}a\le|\xi|\le 2^{j}b\}$ (Notation~\ref{not:LP-partition}), where
  $1+|\xi|^{2(1-\alpha)}\le c\,2^{2(1-\alpha)j}$; hence
  \begin{equation}\label{eq:bnd-Ij}
    \int_{\R^{d}}|\varphi_{j}(\xi)|^{2}\,\mu(d\xi)
    \,\le\,\int_{\{|\xi|\le 2^{j}b\}}
    \bigl[1+|\xi|^{2(1-\alpha)}\bigr]\,\nu(d\xi)
    \,\le\,c_{\mu}\,2^{2(1-\alpha)j}.
  \end{equation}
  The low-frequency block $j=-1$ contributes the finite constant
  $\int|\chi(\xi)|^{2}\mu(d\xi)$ and is consistent with~\eqref{eq:bnd-Ij}.

  \smallskip

  \noindent
  {\it Step 3: Summation.}  Insert~\eqref{eq:second-moment-W}
  and~\eqref{eq:bnd-Ij} into~\eqref{eq:bnd-EW-Aq}. Since the integrand no
  longer depends on $x$, the spatial integral contributes the finite
  factor $\int_{\R^{d}}\rho_{\sigma}(x)\,dx<\infty$, which is where the
  assumption $\sigma>d$ is used (Notation~\ref{not:poly-weight}). We are
  left with
  \begin{equation*}
    \E\Bigl[\|\delta W_{st}\|_{\mathcal{A}_{q}}^{2q}\Bigr]
    \,\le\,c_{q}\,(t-s)^{q}
    \sum_{j\ge -1}2^{-2qj\kappa'}\,2^{2q(1-\alpha)j}
    \,=\,c_{q}\,(t-s)^{q}
    \sum_{j\ge -1}2^{2qj(1-\alpha-\kappa')}.
  \end{equation*}
  Because $\kappa'>1-\alpha$, the geometric series converges, and we
  obtain~\eqref{eq:as-bnd-increments-W}. This completes the proof.
\end{proof}

\noindent
Theorem~\ref{thm:noise-besov} gives a precise sense to the roughness of the
noise. Indeed, for any $\alpha>0$ allowed by~\eqref{eq:dalang-holder}, the
spatial regularity index $-\kappa$ can be taken just below $\alpha-1<0$. This
confirms that $\dot{W}$ is a genuine distribution and not a function. The
negative index is also the reason why the product $u\,\dot{W}$
in~\eqref{eq:pam} is delicate.

\section*{Bibliographic notes}
\addcontentsline{toc}{section}{Bibliographic notes}
\markright{BIBLIOGRAPHIC NOTES}

\smallskip

\noindent
{\it Wiener space and Malliavin calculus.} Section~\ref{sec:wiener-malliavin}
follows Nualart's monograph~\cite{nualart:06:malliavin}, which is the standard
reference for the Malliavin calculus over an abstract Wiener space. Its
Chapter~1 contains, in the generality we need, the isonormal process, the Wiener
chaos decomposition and the pair of operators $D$ and $\delta$. Three further
classical accounts are worth quoting. Malliavin's own
book~\cite{malliavin:97:stochastic} presents the theory from the
differential-geometric angle, which was Malliavin's original point of view. Janson~\cite{janson:97:gaussian}
works entirely with Gaussian Hilbert spaces, and is the cleanest source for the
algebra of the chaoses (product formula, hypercontractivity, moment
comparisons). Nualart and Nualart~\cite{nualart:18:introduction} is a shorter
textbook treatment, convenient for a first reading. For the use of Malliavin calculus
in normal approximation (a direction not represented in this book), we refer to Nourdin
and Peccati~\cite{nourdin:12:normal}.

\smallskip

\noindent
{\it Besov spaces.} The Littlewood--Paley construction of Section~\ref{sec:besov}
is taken from Chapter~2 of Bahouri, Chemin and Danchin~\cite{bahouri:11:fourier},
whose normalizations and notation (the dyadic partition of unity, the blocks
$\Delta_{j}$, the Bernstein inequalities) we have kept throughout. The scale
$\mathcal{B}^{\kappa}$ itself is of course classical. The standard treatises are
those of Triebel~\cite{triebel:83:theory,triebel:92:theory}, with
Bergh and L\"{o}fstr\"{o}m~\cite{bergh:76:interpolation} for the
interpolation-theoretic construction and Runst and
Sickel~\cite{runst:96:sobolev} for nonlinear operations on these spaces. The
harmonic-analysis prerequisites can all be found in
Grafakos~\cite{grafakos:14:modern}.

The weighted spaces of Definition~\ref{def:besov-weighted} are used here for the
reason that brought them into the singular SPDE literature. On the unbounded
domain $\R^{d}$ the noise is not uniformly controlled at infinity, and the
polynomial weight $\rho_{\sigma}$ of Notation~\ref{not:poly-weight} is thus necessary. Our formulation is that of Mourrat and
Weber~\cite{mourrat:17:global}. A systematic theory of Besov spaces with
admissible weights is developed by Edmunds and
Triebel~\cite{edmunds:96:function} and in Triebel~\cite{triebel:06:theory}.

\smallskip

\noindent
{\it The classes of noises.} Dalang's condition (Assumption~\ref{hyp:spectral})
and the extension of Walsh's integral that goes with it are due to
Dalang~\cite{dalang:99:extending}. The same condition was obtained
independently, and at essentially the same time, by Peszat and
Zabczyk~\cite{peszat:97:stochastic,peszat:00:nonlinear}. Their route is a
different one. Instead of Walsh's martingale measures, they work in the Da
Prato and Zabczyk setting of stochastic evolution equations in Hilbert
spaces~\cite{daprato:14:stochastic}. In their framework, the noise is a spatially homogeneous
Wiener process, and~\eqref{eq:dalang-cond} is then what makes the stochastic
convolution a well-defined process with values in a suitable weighted $L^{2}$
space. Two useful surveys of the resulting integration theories are the
minicourse~\cite{dalang:09:minicourse} and the comparison of Dalang and
Quer-Sardanyons~\cite{dalang:11:stochastic}, where the Walsh and Da
Prato-Zabczyk approaches are shown to produce the same solutions. The four
families of Section~\ref{sec:noise-def} are the ones treated by Hu, Huang,
Nualart and Tindel~\cite{hu:15:stochastic}, from which
Examples~\ref{ex:riesz},~\ref{ex:bessel} and~\ref{ex:frac-space} are taken. The
Bessel and Riesz-type kernels of Examples~\ref{ex:bessel}
and~\ref{ex:riesz-type} are discussed in Chen, Ouyang, Tindel and
Xia~\cite{chen:26:ergodic}, where the exponential decay~\eqref{eq:bessel-asymp}
at infinity is what makes the ergodic theory work. The fractional spatial noise
of Example~\ref{ex:frac-space} first appears, in the context of the stochastic
heat equation, in Hu and Nualart~\cite{hu:09:stochastic}. Its rough counterpart
$H<\tfrac12$ falls outside Dalang's condition and calls for entirely different
tools, for which we refer to Hu, Huang, L\^{e}, Nualart and
Tindel~\cite{hu:17:stochastic}. Wider classes of correlations, where $\Lambda$ is
only assumed to be a nonneg-definite measure, are handled by Chen and
Huang~\cite{chen:19:comparison} and Chen and Kim~\cite{chen:19:nonlinear}. For a
noise which is also colored in time, see Balan and
Conus~\cite{balan:16:intermittency}.

\smallskip

\noindent
{\it Inclusion in Besov spaces.} Theorem~\ref{thm:noise-besov} is
\cite[Proposition~5.5]{hu:15:stochastic}, specialized to the white-in-time
setting. The passage from moment bounds to pathwise regularity uses the lemma of
Garsia, Rodemich and Rumsey~\cite{garsia:70:real}. The same
information about $\dot{W}$ is packaged with wavelets rather than
Littlewood--Paley blocks in the regularity structures literature, where the
spaces $\mathcal{C}^{\kappa}_{\mathfrak{s}}$ of~\cite[Section~3]{hairer:14:theory}
play the role of our $\mathcal{B}^{\kappa}$. The precise dictionary between the
two descriptions is the content of the Besov reconstruction theorem of Hairer and
Labb\'{e}~\cite{hairer:17:reconstruction}. Weighted versions on the whole space,
close to those of Section~\ref{sec:besov}, are constructed
in~\cite{hairer:18:multiplicative}, and the two-dimensional PAM is solved in that
framework in~\cite{hairer:15:simple}.

\chapter[Existence and uniqueness]{Existence and uniqueness for the stochastic heat equation}
\label{chap:existence}

\section[Integration theory]{Integration theory with respect to the Gaussian noise}\label{sec:integration}

The mild (Duhamel) formulation of the PAM equation~\eqref{eq:pam} reads
\begin{equation}\label{eq:mild-pam}
  u(t,x) = P_{t}u_{0}(x)
  + \beta\int_{0}^{t}\!\int_{\R^{d}}
    p_{t-s}(x-y)\,u(s,y)\,W(ds,dy),
\end{equation}
where $P_{t}=e^{t\frac{1}{2}\Delta}$ is the heat semigroup introduced at the
beginning of the Bessel-noise example (see~\eqref{eq:heat-kernel-explicit}),
whose action on $u_{0}$ is the ordinary convolution
$P_{t}u_{0}=p_{t}*u_{0}$ with the heat
kernel~$p_{t}$~\eqref{eq:heat-kernel-explicit}.  Let us mention two features
which make the stochastic integral in the right hand side non-trivial. First,
the deterministic factor $s\mapsto p_{t-s}(x-\cdot)$ does not obviously sit in
$\mathcal{H}$. It is not ill-defined, but checking that it belongs to
$\mathcal{H}$ requires a Fourier computation against the spectral measure
$\mu$. Second, the random factor $u(s,y)$ is adapted but not deterministic.
The Hilbert space isometry~\eqref{eq:W-isometry} does not apply directly to
the product $p_{t-s}(x-y)\,u(s,y)$. In this section we set up an integration
theory which handles both features. We follow
Dalang~\cite{dalang:99:extending}, who extended Walsh's martingale measure
integral~\cite{walsh:86:introduction} in order to allow distributional
integrands. One of Dalang's main ideas is to use a Fourier-based norm, which
is weaker than Walsh's.

\subsection{Martingale measure}
We have introduced $W$ as a centered Gaussian family indexed by test
functions in Subsection~\ref{ssec:noise}. However, the white-in-time
covariance~\eqref{eq:noise-cov} shows that $\varphi\mapsto W(\varphi)$ extends
$\sigma$-additively to a $\sigma$-finite $L^{2}(\Omega)$-valued measure,
defined on bounded Borel subsets of $[0,\infty)\times\R^{d}$. This object is
usually called a \emph{martingale measure}.

\begin{notation}[Worthy martingale measure associated with $W$]\label{not:mg-measure}
  For $B\in\mathcal{B}_{b}(\R^{d})$ (the bounded Borel subsets of $\R^{d}$), define
  \begin{equation}\label{eq:mg-process}
    M_{t}(B) \,:=\, W\!\left(\mathbf{1}_{[0,t]\times B}\right),
    \qquad t\ge 0,
  \end{equation}
  where the right-hand side is an isometric extension of the map
  $\varphi\mapsto W(\varphi)$ of Subsection~\ref{ssec:noise} to
  indicator functions.  The covariance~\eqref{eq:noise-cov} gives
  \begin{equation}\label{eq:mg-covariation}
    \E\bigl[M_{t}(A)\,M_{t}(B)\bigr]
    \,=\,t\int_{A}\!\int_{B}\Lambda(x-y)\,dx\,dy,
    \qquad A,B\in\mathcal{B}_{b}(\R^{d}).
  \end{equation}
  Under Assumption~\ref{hyp:spectral}, the process $t\mapsto M_{t}(B)$ is a
  continuous, square-integrable $(\mathcal{F}_{t})$-martingale, where
  \begin{equation}\label{eq:filtration}
    \mathcal{F}_{t}
    \,:=\,\sigma\bigl\{M_{s}(B);\;s\le t,\;B\in\mathcal{B}_{b}(\R^{d})\bigr\}\vee\mathcal{N},
  \end{equation}
  with $\mathcal{N}$ the $\sigma$-field of $\Prob$-null events.  In the language
  of Walsh~\cite{walsh:86:introduction}, the triple
  \begin{equation*}
(M_{t}(B),\mathcal{F}_{t},Q)
\end{equation*}
  is a \emph{worthy martingale measure}, where the covariation measure~$Q$ is
  given by~\eqref{eq:mg-covariation}.  We retain the differential notation
  $W(ds,dx)$ for the stochastic increment of $M$, writing
  \begin{equation}\label{eq:W-differential}
    M_{t}(B)\,=\,\int_{0}^{t}\!\int_{B}W(ds,dx).
  \end{equation}
\end{notation}

\begin{remark}\label{rmk:filtration}
  The filtration $(\mathcal{F}_{t})_{t\ge 0}$ defined in~\eqref{eq:filtration}
  is the natural filtration of the noise $W$. Namely it is the smallest
  filtration for which the process $t\mapsto M_{t}(B)$ is adapted, for every
  $B\in\mathcal{B}_{b}(\R^{d})$. It also coincides with the filtration alluded
  to in~\eqref{eq:duality-explicit} of Subsection~\ref{ssec:malliavin}. Indeed
  we have mentioned there that the Skorohod integral $\delta(u)$ reduces to
  the It\^o--Walsh integral $\int u\,W(dt,dx)$ when $u$ is adapted to the
  filtration generated by $W$. This filtration is exactly
  $(\mathcal{F}_{t})_{t\ge 0}$ as defined in~\eqref{eq:filtration}.
\end{remark}

\begin{remark}[The covariation measure $Q$]\label{rmk:covariation-measure}
  Let us make precise the third entry of the triple
  $(M_{t}(B),\mathcal{F}_{t},Q)$ of
  Notation~\ref{not:mg-measure}.  In Walsh's
  terminology~\cite{walsh:86:introduction}, a martingale measure $M$ is
  \emph{worthy} when there exists a (random) signed measure $Q(dx,dy,ds)$ on
  $\R^{d}\times\R^{d}\times[0,\infty)$, symmetric in $(x,y)$, positive definite,
  and $\sigma$-finite, that represents the covariation of the two
  martingales $t\mapsto M_{t}(A)$ and $t\mapsto M_{t}(B)$ in the sense
  \begin{equation}\label{eq:covariation-bracket}
    \langle M(A),M(B)\rangle_{t}
    \,=\,Q\bigl(A\times B\times(0,t]\bigr),
    \qquad A,B\in\mathcal{B}_{b}(\R^{d}),
  \end{equation}
  and that is controlled by a (positive) \emph{dominating measure} $K$ in the
  sense that $|Q(\Gamma)|\le K(\Gamma)$ for every
  $\Gamma\in\mathcal{B}(\R^{d}\times\R^{d}\times[0,\infty))$.  Here
  $\langle M(A),M(B)\rangle$ denotes the predictable covariation of the two
  square-integrable martingales of Notation~\ref{not:mg-measure}, i.e.\ the
  predictable process making
  $M_{t}(A)\,M_{t}(B)-\langle M(A),M(B)\rangle_{t}$ a martingale.

  For the noise $W$ considered here, the covariance
  computation~\eqref{eq:mg-covariation} identifies $Q$ explicitly: it is the
  \emph{deterministic} measure
  \begin{equation}\label{eq:Q-explicit}
    Q(dx,dy,ds)\,=\,\Lambda(x-y)\,dx\,dy\,ds,
  \end{equation}
  with $\Lambda$ the spatial covariance function of~\eqref{eq:noise-cov}.
  Indeed, evaluating~\eqref{eq:Q-explicit} on $A\times B\times(0,t]$ gives
  \begin{equation*}
    Q\bigl(A\times B\times(0,t]\bigr)
    \,=\,t\int_{A}\!\int_{B}\Lambda(x-y)\,dx\,dy,
  \end{equation*}
  which is exactly the right-hand side of~\eqref{eq:mg-covariation}.  Two
  comments are in order.
  \begin{enumerate}
  \item \emph{Worthiness is automatic here.}  Recall from the covariance
  structure~\eqref{eq:noise-cov} that the spatial covariance
  $\Lambda:\R^{d}\to\R_{+}$ is non-negative, and that it is non-negative
  definite (the sentence following~\eqref{eq:noise-cov}).  The first property
  gives $Q\ge 0$ and the second that $Q$ is positive definite, so that the
  measure $Q$ of~\eqref{eq:Q-explicit} dominates itself: one may take the
  dominating measure $K=Q$, which is what makes $M$ worthy.  Note also that the
  required symmetry of $Q$ in $(x,y)$ is automatic, since a real non-negative
  definite $\Lambda$ is even, $\Lambda(x-y)=\Lambda(y-x)$.

  \item \emph{The bracket is deterministic.}  Because $Q$ does not depend on
  $\omega$, the covariation~\eqref{eq:covariation-bracket} is deterministic as
  well, reflecting the Gaussian, spatially homogeneous and white-in-time nature
  of $W$.  This is the structural feature that will let us reduce, in
  Proposition~\ref{prop:walsh-isometry}, the second-moment of the stochastic
  integral to the explicit bilinear form~\eqref{eq:walsh-isometry}.
  \end{enumerate}
\end{remark}

\subsection{Walsh's stochastic integral}
In his landmark lecture notes~\cite{walsh:86:introduction}, Walsh has built
the stochastic integral $\int g\,W(ds,dx)$ for predictable processes $g$,
thanks to the following $L^{2}$ completion.

\begin{definition}[Elementary processes and Walsh's integral]\label{def:walsh-integral}
  A process $(g(s,x;\omega))$ is \emph{elementary} if
  \begin{equation}\label{eq:elementary-process}
    g(s,x;\omega) \,=\, X(\omega)\,\mathbf{1}_{(a,b]}(s)\,\mathbf{1}_{A}(x),
  \end{equation}
  where $0\le a<b$, $A\in\mathcal{B}_{b}(\R^{d})$, and $X$ is a bounded,
  $\mathcal{F}_{a}$-measurable random variable.  Denote by $\mathcal{E}$ the
  vector space of all finite linear combinations of elementary processes.  For
  $g\in\mathcal{E}$, the \emph{Walsh stochastic integral} is defined by linearity
  from
  \begin{equation}\label{eq:walsh-elementary}
    \int_{0}^{T}\!\int_{\R^{d}}
    X\,\mathbf{1}_{(a,b]}(s)\,\mathbf{1}_{A}(x)\,W(ds,dx)
    \,:=\, X\,\bigl(M_{T\wedge b}(A)-M_{T\wedge a}(A)\bigr).
  \end{equation}
\end{definition}

So far the integral is defined on the dense class $\mathcal{E}$ of
elementary processes~\eqref{eq:elementary-process} only. The next result turns
it into a genuine stochastic integral. It is an $L^{2}$ isometry, which
controls the second moment of~\eqref{eq:walsh-elementary} by a deterministic
bilinear form. The integral is then extended by closure to a large space of
predictable integrands.

\begin{proposition}[Walsh's isometry and extension]\label{prop:walsh-isometry}
  For a predictable process $g$, set
  \begin{equation}\label{eq:walsh-norm-plus}
    \|g\|_{+}^{2}
    \,:=\,\E\!\left[\int_{0}^{T}\!\int_{\R^{d}}\!\int_{\R^{d}}
        |g(s,x)|\,\Lambda(x-y)\,|g(s,y)|\,dx\,dy\,ds\right],
  \end{equation}
  and let $P_{+}$ be the completion of $(\mathcal{E},\|\cdot\|_{+})$.  For every
  $g\in\mathcal{E}$, the Walsh integral~\eqref{eq:walsh-elementary} satisfies
  \begin{equation}\label{eq:walsh-isometry}
    \E\!\left[\left|\int_{0}^{T}\!\int_{\R^{d}}g\,W(ds,dx)\right|^{2}\right]
    \,=\,\E\!\left[\int_{0}^{T}\!\int_{\R^{d}}\!\int_{\R^{d}}
        g(s,x)\,\Lambda(x-y)\,g(s,y)\,dx\,dy\,ds\right]
    \,\le\,\|g\|_{+}^{2}.
  \end{equation}
  Consequently $g\mapsto\int_{0}^{T}\!\int_{\R^{d}}g\,W(ds,dx)$ extends uniquely
  to an $L^{2}(\Omega)$-bounded linear map on all of $P_{+}$, the
  \emph{Walsh stochastic integral}.
\end{proposition}

\begin{proof}
  We prove the equality in~\eqref{eq:walsh-isometry} on the elementary
  integrands~\eqref{eq:elementary-process} and extend it by bilinearity; the
  inequality and the passage to $P_{+}$ are then immediate.

  \smallskip

  \noindent
  {\it Step 1: The identity on a single elementary process.}
  Fix one elementary process
  \begin{equation*}
    g(s,x;\omega)\,=\,X(\omega)\,\mathbf{1}_{(a,b]}(s)\,\mathbf{1}_{A}(x),
  \end{equation*}
  with $X$ bounded and $\mathcal{F}_{a}$-measurable, as
  in~\eqref{eq:elementary-process}.  By the
  definition~\eqref{eq:walsh-elementary} of the integral,
  \begin{equation}\label{eq:elem-integral}
    \int_{0}^{T}\!\int_{\R^{d}}g\,W(ds,dx)
    \,=\,X\bigl(M_{T\wedge b}(A)-M_{T\wedge a}(A)\bigr).
  \end{equation}
  The process $t\mapsto M_{t}(A)$ is a square-integrable
  $(\mathcal{F}_{t})$-martingale (Notation~\ref{not:mg-measure}), and its
  increments after time $a$ are independent of $\mathcal{F}_{a}$. Hence the
  covariation formula~\eqref{eq:mg-covariation} gives the conditional second
  moment
  \begin{equation}\label{eq:incr-cov}
    \E\!\bigl[\bigl(M_{T\wedge b}(A)-M_{T\wedge a}(A)\bigr)^{2}\,\big|\,
      \mathcal{F}_{a}\bigr]
    \,=\,(T\wedge b-T\wedge a)\int_{A}\!\int_{A}\Lambda(x-y)\,dx\,dy.
  \end{equation}
  Since $X$ is $\mathcal{F}_{a}$-measurable, multiplying~\eqref{eq:incr-cov} by
  $X^{2}$, taking expectations, and using~\eqref{eq:elem-integral} yields
  \begin{equation}\label{eq:elem-isometry}
    \E\!\left[\left|\int_{0}^{T}\!\int_{\R^{d}}g\,W(ds,dx)\right|^{2}\right]
    \,=\,\E[X^{2}]\,(T\wedge b-T\wedge a)\int_{A}\!\int_{A}\Lambda(x-y)\,dx\,dy.
  \end{equation}
  Feeding the same $g$ into the bilinear form on the right-hand side
  of~\eqref{eq:walsh-isometry} and using
  $\int_{0}^{T}\mathbf{1}_{(a,b]}(s)\,ds=T\wedge b-T\wedge a$ produces exactly
  the right-hand side of~\eqref{eq:elem-isometry}.  This proves  equality~\eqref{eq:walsh-isometry}
  for a single elementary process.

  \smallskip

  \noindent
  {\it Step 2: Polarization to finite sums.}
  For a general $g=\sum_{i}g_{i}\in\mathcal{E}$, define
  $I_{i}$ as 
  \begin{equation*}
I_{i}
=
\int_{0}^{T}\!\int_{\R^{d}}g_{i}\,W(ds,dx) .
\end{equation*}
Expanding the square on
  the left of~\eqref{eq:walsh-isometry} and the bilinear form on its right gives
  the two double sums
  \begin{equation}\label{eq:polarization}
    \E\!\left[\Bigl|\sum_{i}I_{i}\Bigr|^{2}\right]
    \,=\,\sum_{i,j}\E\bigl[I_{i}I_{j}\bigr],
    \qquad
    \sum_{i,j}\E\!\left[\int_{0}^{T}\!\int_{\R^{d}}\!\int_{\R^{d}}
      g_{i}(s,x)\,\Lambda(x-y)\,g_{j}(s,y)\,dx\,dy\,ds\right].
  \end{equation}
  Each cross term $\E[I_{i}I_{j}]$ is evaluated exactly as
  in~\eqref{eq:incr-cov}-\eqref{eq:elem-isometry}, conditioning on
  $\mathcal{F}_{a}$ for the later of the two starting times. Martingale
  increments over disjoint time intervals are orthogonal. On overlapping
  intervals, the covariation~\eqref{eq:mg-covariation} produces again the
  spatial factor $\int_{A_{i}}\!\int_{A_{j}}\Lambda(x-y)\,dx\,dy$. Term by term
  the two sums in~\eqref{eq:polarization} agree, which is the equality
  in~\eqref{eq:walsh-isometry} for every $g\in\mathcal{E}$.

  \smallskip

  \noindent
  {\it Step 3: Inequality and extension.}
  The subsequent inequality in~\eqref{eq:walsh-isometry} is the pointwise bound
  $g(s,x)\,g(s,y)\le|g(s,x)|\,|g(s,y)|$ together with $\Lambda\ge 0$.  The
  integral then extends to all of $P_{+}$ by $L^{2}(\Omega)$-closure.  See
  Walsh~\cite[Chapter~2]{walsh:86:introduction} for the worthy
  martingale-measure framework on which~\eqref{eq:incr-cov} relies.
\end{proof}

We close this subsection with a result which allows to exchange a Walsh
integral with an integral with respect to an auxiliary parameter. It will be
used several times in the sequel, either with a Lebesgue integral in the space
variable or with the expectation with respect to a Brownian motion independent
of $W$. Below, $\mathcal{P}$ denotes the predictable $\sigma$-field on
$[0,T]\times\R^{d}\times\Omega$, that is the $\sigma$-field generated by the
elementary processes of Definition~\ref{def:walsh-integral}.

\begin{assumption}[Stochastic Fubini setting]\label{hyp:stoch-fubini}
  Let $(G,\mathcal{G},\nu)$ be a $\sigma$-finite measure space and $T>0$.
  Consider a function $f:G\times[0,T]\times\R^{d}\times\Omega\to\R$, and write
  $f(z)$ for the process $(s,y,\omega)\mapsto f(z,s,y,\omega)$. We assume the
  following.
  \begin{itemize}
  \item The map $f$ is $\mathcal{G}\otimes\mathcal{P}$-measurable.
  \item The norms~\eqref{eq:walsh-norm-plus} of the processes $f(z)$ are
  $\nu$-integrable, namely
  \begin{equation}\label{eq:fubini-cond}
    \mathcal{N}(f)\,:=\,\int_{G}\|f(z)\|_{+}\,\nu(dz)\,<\,\infty.
  \end{equation}
  \end{itemize}
\end{assumption}

\begin{proposition}[Stochastic Fubini theorem]\label{prop:stoch-fubini}
  Let $f$ satisfy Assumption~\ref{hyp:stoch-fubini}. Then the following
  assertions hold true.
  \begin{enumerate}
  \item For $\nu$-almost every $z\in G$ we have $f(z)\in P_{+}$. The stochastic
  integrals
  \begin{equation*}
    I(z)\,:=\,\int_{0}^{T}\!\int_{\R^{d}}f(z,s,y)\,W(ds,dy)
  \end{equation*}
  admit a $\mathcal{G}\otimes\mathcal{F}_{T}$-measurable version, which
  satisfies $\E\int_{G}|I(z)|\,\nu(dz)\le\mathcal{N}(f)$.
  \item The process $\bar{f}:=\int_{G}|f(z)|\,\nu(dz)$ satisfies
  $\|\bar{f}\|_{+}\le\mathcal{N}(f)$. In particular the process
  $F:=\int_{G}f(z)\,\nu(dz)$, set to $0$ where $\bar{f}=\infty$, belongs to
  $P_{+}$ and $\|F\|_{+}\le\mathcal{N}(f)$.
  \item Almost surely, we have
  \begin{equation}\label{eq:stoch-fubini}
    \int_{G}\Bigl(\int_{0}^{T}\!\int_{\R^{d}}f(z,s,y)\,W(ds,dy)\Bigr)\nu(dz)
    \,=\,
    \int_{0}^{T}\!\int_{\R^{d}}\Bigl(\int_{G}f(z,s,y)\,\nu(dz)\Bigr)W(ds,dy).
  \end{equation}
  \end{enumerate}
\end{proposition}

\begin{proof}
  We only sketch the argument, which follows the usual pattern. Namely we check
  the identity for simple integrands, and both sides of~\eqref{eq:stoch-fubini}
  are then shown to be continuous with respect to the quantity $\mathcal{N}(f)$
  of~\eqref{eq:fubini-cond}.

  \smallskip

  \noindent
  {\it Step 1: Bounds on both sides.} For $\nu$-almost every $z$ we have
  $\|f(z)\|_{+}<\infty$ owing to~\eqref{eq:fubini-cond}, so that $f(z)\in P_{+}$.
  The measurability of $(z,\omega)\mapsto I(z)$ is a technical point which we
  take for granted (see~\cite[Theorem~2.6]{walsh:86:introduction}
  and~\cite[Section~2.6]{dalang:26:stochastic}). Then the Cauchy-Schwarz
  inequality and the isometry~\eqref{eq:walsh-isometry} give
  $\E|I(z)|\le\|f(z)\|_{+}$, and integrating in $z$ yields the bound in
  item~(i). As far as item~(ii) is concerned, let us introduce the bilinear form
  \begin{equation*}
    \langle g,h\rangle_{+}
    \,:=\,\E\!\left[\int_{0}^{T}\!\int_{\R^{d}}\!\int_{\R^{d}}
      g(s,x)\,\Lambda(x-y)\,h(s,y)\,dx\,dy\,ds\right],
  \end{equation*}
  so that $\|g\|_{+}^{2}=\langle|g|,|g|\rangle_{+}$. Since $\Lambda$ is
  non-negative definite, $\langle\cdot,\cdot\rangle_{+}$ satisfies the
  Cauchy-Schwarz inequality. Since moreover $\Lambda\ge0$, Tonelli's theorem
  applies to non-negative integrands, and we get
  \begin{equation*}
    \|\bar{f}\|_{+}^{2}
    \,=\,\int_{G}\bigl\langle|f(z)|,\bar{f}\bigr\rangle_{+}\,\nu(dz)
    \,\le\,\|\bar{f}\|_{+}\int_{G}\|f(z)\|_{+}\,\nu(dz).
  \end{equation*}
  This is $\|\bar{f}\|_{+}\le\mathcal{N}(f)$, once $\|\bar{f}\|_{+}$ is known
  to be finite. The latter restriction is removed by first truncating $f$ and
  then letting the truncation level go to infinity, by monotone convergence.
  Finally $|F|\le\bar{f}$, so that $\|F\|_{+}\le\|\bar{f}\|_{+}$.

  \smallskip

  \noindent
  {\it Step 2: Simple integrands.} Assume that $f$ is of the form
  \begin{equation*}
    f(z,s,y)\,=\,\sum_{i=1}^{n}\mathbf{1}_{A_{i}}(z)\,g_{i}(s,y),
  \end{equation*}
  where $A_{i}\in\mathcal{G}$ satisfies $\nu(A_{i})<\infty$ and
  $g_{i}\in\mathcal{E}$. Then by linearity of the Walsh integral, both sides
  of~\eqref{eq:stoch-fubini} are equal to
  $\sum_{i=1}^{n}\nu(A_{i})\int_{0}^{T}\!\int_{\R^{d}}g_{i}\,W(ds,dy)$.

  \smallskip

  \noindent
  {\it Step 3: Conclusion.} Both sides of~\eqref{eq:stoch-fubini} are linear in
  $f$. Owing to Step~1, their $L^{1}(\Omega)$ norms are bounded by
  $\mathcal{N}(f)$. In addition, a standard density argument shows that the
  simple integrands of Step~2 are dense in the space of functions $f$
  satisfying~\eqref{eq:fubini-cond}, for the norm $\mathcal{N}$. Hence the
  identity of Step~2 extends to all such $f$, which proves~\eqref{eq:stoch-fubini}.
\end{proof}

\begin{remark}\label{rmk:stoch-fubini}
  Proposition~\ref{prop:stoch-fubini} is due to
  Walsh~\cite[Theorem~2.6]{walsh:86:introduction}, who states it for a finite
  measure $\nu$ and worthy martingale measures. For space-time white noise, a
  sharper version is given by Dalang and
  Sanz-Sol\'{e}~\cite[Theorem~2.4.1]{dalang:26:stochastic}, where
  condition~\eqref{eq:fubini-cond} is replaced by an almost sure condition (see
  also~\cite[Corollary~2.4.3]{dalang:26:stochastic} for the comparison with
  Walsh's condition). We shall apply Proposition~\ref{prop:stoch-fubini} in two
  situations. Either $G=\R^{d}$ and $\nu$ is absolutely continuous with respect
  to the Lebesgue measure, which allows to move a deterministic convolution
  through the stochastic integral. Or $G$ is the path space of a Brownian motion
  $B$ independent of $W$ and $\nu$ is the law of $B$, in which case
  relation~\eqref{eq:stoch-fubini} exchanges the expectation $\E_{B}$ with the
  Walsh integral.
\end{remark}

\subsection{Dalang's extension}
\label{sec:dalang-ext}
The space $P_{+}$ is not sufficient for our purposes. Indeed, the norm
$\|\cdot\|_{+}$ involves the absolute values $|g(s,x)|$, which are meaningless
when $g(s,\cdot)$ is a distribution. The decisive observation of
Dalang~\cite{dalang:99:extending} is that the isometry
in~\eqref{eq:walsh-isometry} uses the inner product $\int
g(s,x)\Lambda(x-y)g(s,y)\,dx\,dy$ \emph{without} absolute values. The
left-hand side is thus controlled by a strictly weaker norm. The completion
$\mathcal{H}_{\omega}$ for this norm is much larger than $P_{+}$, and it does
accommodate distributional integrands.

\begin{definition}[Extended norm and Hilbert space $\mathcal{H}_{\omega}$]\label{def:extended-norm}
  For predictable integrands $g,h$, average the Hilbert-space inner
  product~\eqref{eq:H-innprod} over the randomness:
  \begin{equation}\label{eq:innprod-extended}
    \langle g,h\rangle_{\mathcal{H}_{\omega}}
    \,:=\,\E\bigl[\langle g,h\rangle_{\mathcal{H}}\bigr]
    \,=\,\E\!\left[\int_{0}^{T}\!\int_{\R^{d}}\!\int_{\R^{d}}
        g(s,x)\,\Lambda(x-y)\,h(s,y)\,dx\,dy\,ds\right],
  \end{equation}
  where $\langle\cdot,\cdot\rangle_{\mathcal{H}}$ is the inner
  product~\eqref{eq:H-innprod} applied to the sample paths
  $(s,x)\mapsto g(s,x;\omega)$ (supported in $[0,T]$).  In other words,
  $\langle\cdot,\cdot\rangle_{\mathcal{H}_{\omega}}$ is the inner product of
  $L^{2}(\Omega;\mathcal{H})$, the random counterpart of $\mathcal{H}$.  The
  associated norm
  \begin{equation}\label{eq:norm-0}
    \|g\|_{\mathcal{H}_{\omega}}\,:=\,\langle g,g\rangle_{\mathcal{H}_{\omega}}^{1/2}
    \,=\,\E\bigl[\|g\|_{\mathcal{H}}^{2}\bigr]^{1/2}
  \end{equation}
  is the norm under which Dalang~\cite{dalang:99:extending} performs the
  extension, where it is denoted $\|\cdot\|_{0}$.  We write $\mathcal{H}_{\omega}$
  for the completion of $(\mathcal{E},\|\cdot\|_{\mathcal{H}_{\omega}})$; it is a
  Hilbert space, canonically identified with $L^{2}(\Omega;\mathcal{H})$, and the
  subscript $\mathcal{H}_{\omega}$ in~\eqref{eq:innprod-extended}
  and~\eqref{eq:norm-0} records both the inner product and the space at once.
  The same construction will be performed for the modified martingale measure
  $M^{Z}$ in Notation~\ref{not:norm-Z}, yielding the norm
  $\|\cdot\|_{\mathcal{H}_{Z}}$ of~\eqref{eq:norm-Z} and its associated Hilbert
  space $\mathcal{H}_{Z}$.
\end{definition}

\begin{remark}[Fourier form of the extended inner product]\label{rmk:innprod-extended-fourier}
  Applying the Fourier identity~\eqref{eq:H-innprod-fourier} pathwise to the
  sample paths $(s,x)\mapsto g(s,x;\omega)$ and $(s,x)\mapsto h(s,x;\omega)$ and
  then averaging, the right-hand side of~\eqref{eq:innprod-extended} can be
  written in Fourier mode as
  \begin{equation}\label{eq:innprod-extended-fourier}
    \langle g,h\rangle_{\mathcal{H}_{\omega}}
    \,=\,\E\!\left[\int_{0}^{T}\!\int_{\R^{d}}
      \mathcal{F}_{x}g(s,\cdot)(\xi)\,\overline{\mathcal{F}_{x}h(s,\cdot)(\xi)}\,
      \mu(d\xi)\,ds\right],
  \end{equation}
  where $\mathcal{F}_{x}$ is the spatial Fourier transform and
  $\mu=\mathcal{F}\Lambda$ is the spectral measure of the
  noise~\eqref{eq:noise-cov}. Dropping the expectation
  in~\eqref{eq:innprod-extended-fourier} recovers the deterministic
  identity~\eqref{eq:innprod-fourier} (Remark~\ref{rmk:deterministic-fourier}).
\end{remark}

The point of this weaker norm is that the Walsh integral remains an isometry for
it, so that the construction of Proposition~\ref{prop:walsh-isometry} carries
over verbatim to the larger space $\mathcal{H}_{\omega}$.

\begin{theorem}[Extension of the Walsh integral to $\mathcal{H}_{\omega}$]\label{thm:extended-norm}
  One has $P_{+}\subset \mathcal{H}_{\omega}$, and the Walsh integral of
  Proposition~\ref{prop:walsh-isometry} extends by continuity to an isometry
  \begin{eqnarray}\label{eq:P0-isometry}
    \mathcal{H}_{\omega}&\longrightarrow&
    \bigl\{\text{continuous square-integrable }(\mathcal{F}_{t})\text{-martingales}\bigr\} \notag\\
     g&\longmapsto&\int_{0}^{\cdot}\!\int_{\R^{d}}g\,W(ds,dx).
  \end{eqnarray}
  In particular, for every $g\in\mathcal{H}_{\omega}$ the second moment of the
  integral is given by the isometry identity
  \begin{multline}\label{eq:extended-isometry}
    \E\!\left[\left|\int_{0}^{T}\!\int_{\R^{d}}g\,W(ds,dx)\right|^{2}\right]
    \,=\,\|g\|_{\mathcal{H}_{\omega}}^{2}  \\
    \,=\,\E\!\left[\int_{0}^{T}\!\int_{\R^{d}}\!\int_{\R^{d}}
        g(s,x)\,\Lambda(x-y)\,g(s,y)\,dx\,dy\,ds\right],
  \end{multline}
  the norm $\|\cdot\|_{\mathcal{H}_{\omega}}$ being the one
  of~\eqref{eq:norm-0}.
\end{theorem}

\begin{proof}[Proof of Theorem~\ref{thm:extended-norm}]
  Since $\Lambda\ge 0$, the integrand in~\eqref{eq:innprod-extended} is
  dominated, for $g=h$, by the one in~\eqref{eq:walsh-norm-plus}, so that
  $\|g\|_{\mathcal{H}_{\omega}}\le\|g\|_{+}$ for every $g\in\mathcal{E}$.  Hence the identity map on
  $\mathcal{E}$ is $\|\cdot\|_{+}$-to-$\|\cdot\|_{\mathcal{H}_{\omega}}$ continuous and extends to
  the inclusion $P_{+}\subset \mathcal{H}_{\omega}$.  By the
  isometry~\eqref{eq:walsh-isometry} of
  Proposition~\ref{prop:walsh-isometry}, the middle term of which equals
  $\|g\|_{\mathcal{H}_{\omega}}^{2}$, the map $g\mapsto\int g\,W(ds,dx)$ is an isometry from
  $(\mathcal{E},\|\cdot\|_{\mathcal{H}_{\omega}})$ into $L^{2}(\Omega)$; since $\mathcal{E}$ is
  dense in $\mathcal{H}_{\omega}$ by construction, it extends to all of $\mathcal{H}_{\omega}$,
  yielding both the isometry~\eqref{eq:P0-isometry} and the
  identity~\eqref{eq:extended-isometry}.
\end{proof}

\begin{remark}[The deterministic case]\label{rmk:deterministic-fourier}
  For deterministic $g,h$ the expectation in the Fourier
  identity~\eqref{eq:innprod-extended-fourier}
  (Remark~\ref{rmk:innprod-extended-fourier}) is vacuous, so
  $\langle g,h\rangle_{\mathcal{H}_{\omega}}=\langle g,h\rangle_{\mathcal{H}}$ and it
  reduces to
  \begin{equation}\label{eq:innprod-fourier}
    \langle g,h\rangle_{\mathcal{H}}
    \,=\,\int_{0}^{T}\!\int_{\R^{d}}
    \mathcal{F}_{x}g(s,\cdot)(\xi)\,\overline{\mathcal{F}_{x}h(s,\cdot)(\xi)}\,
    \mu(d\xi)\,ds.
  \end{equation}
  This Fourier representation extends to distributional integrands: it is
  well-defined as soon as $\mathcal{F}_{x}S(t,\cdot)$ is a function, even when
  $S(t,\cdot)\in\mathcal{S}'(\R^{d})$ (Notation~\ref{not:schwartz}).
\end{remark}

The following proposition makes the claim of Remark~\ref{rmk:filtration}
rigorous for processes in $\mathcal{H}_{\omega}$. Once the Walsh integral has
been extended to $\mathcal{H}_{\omega}$ by Theorem~\ref{thm:extended-norm},
every adapted process in $\mathcal{H}_{\omega}$ belongs to
$\mathrm{Dom}(\delta)$, and the Skorohod integral reduces to the
It\^{o}--Walsh integral. 
We adapt this result from~\cite[Proposition~1.3.3]{nualart:06:malliavin}, 
which handled the case of integrals with respect to a classical Brownian motion.

\begin{proposition}[It\^{o}--Walsh $=$ Skorohod for adapted integrands]
  \label{prop:ito-sko}
  Let $v\in\mathcal{H}_{\omega}$ be adapted to the filtration
  $(\mathcal{F}_{s})_{s\ge 0}$ of~\eqref{eq:filtration}. Then
  $v\in\mathrm{Dom}(\delta)$ and
  \begin{equation}\label{eq:ito-sko-H}
    \delta(v)
    \;=\;\int_{0}^{T}\!\int_{\R^{d}}v(s,y)\,W(ds,dy),
  \end{equation}
  where the right-hand side is the It\^{o}--Walsh integral of
  Theorem~\ref{thm:extended-norm}.
\end{proposition}

\begin{proof}
  \smallskip
  \noindent
  {\it Step 1: Locality of $D$ for adapted functionals.}
  Let $0\le a\le T$, let $F\in\mathbb{D}^{1,2}$ be $\mathcal{F}_{a}$-measurable,
  and let $g\in\mathcal{H}$ be a deterministic element satisfying $g(r,\cdot)=0$
  for $r\le a$. We claim that
  \begin{equation}\label{eq:D-locality}
    \langle DF,g\rangle_{\mathcal{H}} = 0 \qquad\text{a.s.}
  \end{equation}
  Indeed, a typical smooth cylindrical functional which is
  $\mathcal{F}_{a}$-measurable is of the form
  \begin{equation*}
F=f(W(h_{1}),\dots,W(h_{n})) \, ,
\end{equation*}
with each $h_{j}\in\mathcal{H}$
  satisfying $h_{j}(r,\cdot)=0$ for $r>a$. For this type of functional, the derivative
  formula~\eqref{eq:malliavin-D} gives
  \begin{equation*}
    DF = \sum_{j=1}^{n}
    \frac{\partial f}{\partial x_{j}}(W(h_{1}),\dots,W(h_{n}))\,h_{j}.
  \end{equation*}
  Since $h_{j}(r,\cdot)=0$ for $r>a$ and $g(r,\cdot)=0$ for $r\le a$, the
  inner product~\eqref{eq:H-innprod} yields $\langle h_{j},g\rangle_{\mathcal{H}}=0$
  for every $j$, hence $\langle DF,g\rangle_{\mathcal{H}}=0$. The general case
  $F\in\mathbb{D}^{1,2}\cap L^{2}(\mathcal{F}_{a})$ follows by approximating $F$
  by $\mathcal{F}_{a}$-measurable smooth cylindricals (which are dense in
  $L^{2}(\mathcal{F}_{a})$) and invoking the closability of $D$
  (Proposition~\ref{prop:D-closable}).

  \smallskip
  \noindent
  {\it Step 2: Elementary adapted integrands.}
  Let $0\le a<b\le T$, let $g\in\mathcal{H}$ be deterministic with $g(r,\cdot)=0$
  for $r\notin(a,b]$, and let $F\in\mathbb{D}^{1,2}\cap L^{2}(\mathcal{F}_{a})$.
  The product rule~\eqref{eq:product-rule} applied to $F$ and the deterministic
  element $g\in\mathcal{H}$ gives
  \begin{equation*}
    F\,W(g) = \delta(Fg) + \langle DF,g\rangle_{\mathcal{H}}.
  \end{equation*}
  By Step~1, $\langle DF,g\rangle_{\mathcal{H}}=0$ a.s. On the other hand,
  the It\^{o}--Walsh integral of the elementary adapted process $Fg$ equals
  $\int_{0}^{T}\!\int_{\R^{d}}F\,g\,W(ds,dy) = F\,W(g)$ by the isometry
  of Theorem~\ref{thm:extended-norm}. Hence $Fg\in\mathrm{Dom}(\delta)$ and
  \begin{equation*}
    \delta(Fg) = F\,W(g) = \int_{0}^{T}\!\int_{\R^{d}}Fg\,W(ds,dy).
  \end{equation*}
  Equation~\eqref{eq:ito-sko-H} therefore holds for every such elementary
  adapted integrand $v=Fg$.

  \smallskip
  \noindent
  {\it Step 3: Extension by density.}
  Finite linear combinations of elementary adapted processes $Fg$ as in Step~2
  are dense in the subspace of adapted processes in $\mathcal{H}_{\omega}$. To
  see this, note that for each $a\ge 0$ the smooth cylindricals
  $f(W(h_{1}),\dots,W(h_{n}))$ with $h_{j}$ supported in $[0,a]\times\R^{d}$
  form a subspace of $\mathbb{D}^{1,2}\cap L^{2}(\mathcal{F}_{a})$ which is
  dense in $L^{2}(\mathcal{F}_{a})$. Combining this with step-function
  approximations in the time variable, one gets the claimed density in
  $\mathcal{H}_{\omega}$. Let
  $v\in\mathcal{H}_{\omega}$ be adapted and let $(v_{n})$ be such an approximating
  sequence. By Step~2 and linearity, $\delta(v_{n})=\int v_{n}\,W$ for each $n$.
  The Walsh isometry~\eqref{eq:extended-isometry} gives
  \begin{equation*}
    \left\|\int v_{n}\,W - \int v\,W\right\|_{L^{2}(\Omega)}
    = \|v_{n}-v\|_{\mathcal{H}_{\omega}} \to 0.
  \end{equation*}
  Since $\delta$ is a closed operator (it is the adjoint of $D$;
  Definition~\ref{def:skorohod}), the convergence $v_{n}\to v$ in
  $\mathcal{H}_{\omega}$ and $\delta(v_{n})\to\int v\,W$ in $L^{2}(\Omega)$
  together imply $v\in\mathrm{Dom}(\delta)$ and
  $\delta(v)=\int_{0}^{T}\!\int_{\R^{d}}v\,W(ds,dy)$.
\end{proof}

\subsection{Integration against a modified martingale measure}\label{sec:intg-modif}
The second difficulty announced at the beginning of this section is that the
integrand in~\eqref{eq:mild-pam} is the product of the (deterministic) heat
kernel and the (random) solution $u$.  In Picard's iteration scheme, this means
one needs to integrate a deterministic function $S(s,\cdot)$ against the
modified noise measure $u_{n}\cdot W(ds,dx)$, where $u_{n}$ is the $n$-th
iterate.  Dalang~\cite{dalang:99:extending} handles this by replacing the
worthy martingale measure $M$ of Notation~\ref{not:mg-measure} with a
\emph{modified} martingale measure $M^{Z}$, in which $Z$ plays the role of
the current iterate.

\begin{definition}[Modified martingale measure $M^{Z}$]\label{def:mg-Z}
  Let $Z=(Z(s,x),\;0\le s\le T,\;x\in\R^{d})$ be an $(\mathcal{F}_{s})$-adapted
  process satisfying
  \begin{equation}\label{eq:Z-L2-bound}
    \sup_{0\le s\le T}\;\sup_{x\in\R^{d}}\;\E\bigl[Z(s,x)^{2}\bigr]<+\infty.
  \end{equation}
  The \emph{modified martingale measure} $M^{Z}$ is defined by
  \begin{equation}\label{eq:mg-Z-def}
    M^{Z}_{t}(B)
    \,:=\,\int_{0}^{t}\!\int_{B}Z(s,y)\,W(ds,dy),
    \qquad B\in\mathcal{B}_{b}(\R^{d}),
  \end{equation}
  where the integral is a Walsh stochastic integral in the sense of
  Definition~\ref{def:walsh-integral} and Proposition~\ref{prop:walsh-isometry}.
\end{definition}

In order to extend the integration against $M^{Z}$ to distributional
integrands, the natural analogue of~\eqref{eq:innprod-fourier} involves the
effective spectral measure of the product $Z\cdot W$. This measure depends on
the two-point statistics of $Z$. The following spatial homogeneity condition
(see~\cite[Section~2]{dalang:99:extending}) ensures that those statistics are
well structured.

\begin{assumption}[Spatial homogeneity of $Z$]\label{hyp:spatial-hom}
  The process $Z$ satisfying~\eqref{eq:Z-L2-bound} is
  \emph{spatially homogeneous in mean square}: for all $s\ge 0$ and
  $x,y\in\R^{d}$,
  \begin{equation}\label{eq:spatial-hom}
    \E\bigl[Z(s,x)\,Z(s,y)\bigr]
    \,=\,\E\bigl[Z(s,0)\,Z(s,x-y)\bigr].
  \end{equation}
\end{assumption}

Under this homogeneity, the two-point statistics of $Z$ depend on the spatial
increment only. One can thus attach a spectral measure to $M^{Z}$, exactly as
the noise itself carries the spectral measure $\mu=\mathcal{F}\Lambda$
of~\eqref{eq:noise-cov}.

\begin{proposition}[Effective spectral measure of $M^{Z}$]\label{prop:effective-spectral}
  Under Assumption~\ref{hyp:spatial-hom}, the function
  \begin{equation}\label{eq:gs-def}
    g_{s}(z)\,:=\,\E\bigl[Z(s,0)\,Z(s,z)\bigr],
    \qquad z\in\R^{d},
  \end{equation}
  is well-defined, non-negative definite, and independent of the reference
  point.  Consequently the pointwise product $\Lambda(\cdot)\,g_{s}(\cdot)$ is
  non-negative definite, and its Fourier transform is a non-negative tempered
  measure
  \begin{equation}\label{eq:muZ-def}
    \mu^{Z}_{s}(d\xi)
    \,:=\,\mathcal{F}\bigl[\Lambda(\cdot)\,g_{s}(\cdot)\bigr](d\xi),
  \end{equation}
  called the \emph{effective spectral measure} of $M^{Z}$ at time $s$.
\end{proposition}

\begin{proof}
  By~\eqref{eq:spatial-hom} the right-hand side of~\eqref{eq:gs-def} depends
  only on the increment $z=x-y$, so $g_{s}$ is well-defined and independent of
  the reference point; it is non-negative definite, being the covariance
  function of the mean-square stationary field $x\mapsto Z(s,x)$.  Since
  $\Lambda$ is non-negative definite as well (the standing hypothesis following
  the covariance structure~\eqref{eq:noise-cov}), the pointwise product
  $\Lambda\,g_{s}$ is again non-negative definite.  Indeed, fix points
  $x_{1},\dots,x_{n}\in\R^{d}$; saying that a function $\phi$ is non-negative
  definite means precisely that the matrix $[\phi(x_{i}-x_{j})]_{1\le i,j\le n}$
  is positive semidefinite.  The two matrices
  \begin{equation*}
    \bigl[\Lambda(x_{i}-x_{j})\bigr]_{i,j}
    \quad\text{and}\quad
    \bigl[g_{s}(x_{i}-x_{j})\bigr]_{i,j}
  \end{equation*}
  are thus positive semidefinite, and Schur's product theorem (the Hadamard,
  i.e.\ entrywise, product of two positive semidefinite matrices is positive
  semidefinite) shows that their entrywise product
  $[\Lambda(x_{i}-x_{j})\,g_{s}(x_{i}-x_{j})]_{i,j}$ is positive
  semidefinite as well. Since $x_{1},\dots,x_{n}$ were arbitrary, $\Lambda\,g_{s}$
  is non-negative definite.  By Bochner's theorem its Fourier
  transform~\eqref{eq:muZ-def} is therefore a non-negative measure.
\end{proof}

\begin{notation}[Norm $\|\cdot\|_{\mathcal{H}_{Z}}$ and space $\mathcal{H}_{Z}$]\label{not:norm-Z}
  For a deterministic function $S:[0,T]\to\mathcal{S}'(\R^{d})$ such that
  $\mathcal{F}_{x}S(t,\cdot)$ is a measurable function for each~$t$, set
  \begin{equation}\label{eq:norm-Z}
    \|S\|_{\mathcal{H}_{Z}}^{2}
    \,:=\,\int_{0}^{T}\!\int_{\R^{d}}
       \bigl|\mathcal{F}_{x}S(t,\cdot)(\xi)\bigr|^{2}\,\mu^{Z}_{t}(d\xi)\,dt.
  \end{equation}
  The completion of the space of deterministic Schwartz functions under
  $\|\cdot\|_{\mathcal{H}_{Z}}$ is denoted $\mathcal{H}_{Z}$.  For $S\in \mathcal{H}_{Z}$, the stochastic
  integral $\int S\,dM^{Z}$ is defined via the isometry
  \begin{multline}\label{eq:isometry-MZ}
    \E\!\left[\left|\int_{0}^{t}\!\int_{\R^{d}}S(s,x)\,M^{Z}(ds,dx)\right|^{2}
    \right] \\
    \,=\,\bigl\|S\,\mathbf{1}_{[0,t]}\bigr\|_{\mathcal{H}_{Z}}^{2} 
    \,=\,\int_{0}^{t}\!\int_{\R^{d}}
    \bigl|\mathcal{F}_{x}S(s,\cdot)(\xi)\bigr|^{2}\,\mu^{Z}_{s}(d\xi)\,ds.
  \end{multline}
\end{notation}

The main result of this section gives a sufficient condition for a
deterministic integrand $S$ to belong to $\mathcal{H}_{Z}$. This condition is
expressed in terms of the spectral measure $\mu$ of the noise, and not in terms
of the effective measure $\mu^{Z}$. It is thus easily checked in practice,
thanks to Assumption~\ref{hyp:spectral} only. Let us also mention that the
hypothesis of a \emph{non-negative} distribution $S(t)$ ensures that the
absolute values can be removed from the Walsh norm. This is the key step in the
proof.

\begin{theorem}[Dalang's stochastic integral]\label{thm:dalang-integral}
  Let $Z$ satisfy~\eqref{eq:Z-L2-bound} and Assumption~\ref{hyp:spatial-hom}.
  Let $t\mapsto S(t)$ be a deterministic function with values in the space of
  \emph{non-negative tempered distributions with rapid decrease}, i.e.\
  $S(t)\in\mathcal{S}'(\R^{d})$ (Notation~\ref{not:schwartz}),
  $\langle S(t),\varphi\rangle\ge 0$ for all $\varphi\ge 0$, and
  $\mathcal{F}_{x}S(t,\cdot)$ is a locally bounded function for each~$t$.
  Suppose that
  \begin{equation}\label{eq:dalang-cond-S}
    \int_{0}^{T}\!\int_{\R^{d}}
    \bigl|\mathcal{F}_{x}S(t,\cdot)(\xi)\bigr|^{2}\,\mu(d\xi)\,dt
    \;<\;+\infty.
  \end{equation}
  Then $S\in \mathcal{H}_{Z}$ (see Notation~\ref{not:norm-Z}) and the
  stochastic integral $\int S\,dM^{Z}$ is well-defined.  Moreover, for every
  $0\le t\le T$,
  \begin{equation}\label{eq:dalang-isometry-thm}
    \E\!\left[\left|\int_{0}^{t}\!\int_{\R^{d}}S(s,x)\,M^{Z}(ds,dx)
    \right|^{2}\right]
    \,=\,\int_{0}^{t}\!\int_{\R^{d}}
    \bigl|\mathcal{F}_{x}S(s,\cdot)(\xi)\bigr|^{2}\,\mu^{Z}_{s}(d\xi)\,ds,
  \end{equation}
  and the following $L^{2}$ bound holds:
  \begin{multline}\label{eq:dalang-bound}
    \E\!\left[\left|\int_{0}^{t}\!\int_{\R^{d}}S(s,x)\,M^{Z}(ds,dx)
    \right|^{2}\right]  \\
    \,\le\,\left(\sup_{0\le s\le t}\sup_{x\in\R^{d}}\E\bigl[Z(s,x)^{2}\bigr]
    \right)
    \int_{0}^{t}\!\int_{\R^{d}}
    \bigl|\mathcal{F}_{x}S(s,\cdot)(\xi)\bigr|^{2}\,\mu(d\xi)\,ds.
  \end{multline}
\end{theorem}

It is convenient to abbreviate the prefactor in~\eqref{eq:dalang-bound}. This
quantity plays a central role in the contraction argument for the stochastic
heat equation (Section~\ref{sec:contraction}). There the weight $Z$ is the
unknown solution itself, and the factor $\mathcal{N}_{t,2}(Z)$ is exactly what
has to be controlled in order to close the fixed point estimate.

\begin{notation}\label{not:N-t-2}
  For a weight $Z$ satisfying~\eqref{eq:Z-L2-bound} and $t\ge 0$, we set
  \begin{equation}\label{eq:N-t-2-def}
    \mathcal{N}_{t,2}(Z)
    \,:=\, \sup_{0\le s\le t}\sup_{x\in\R^{d}}\E\bigl[Z(s,x)^{2}\bigr].
  \end{equation}
\end{notation}

\begin{remark}\label{rmk:dalang-bound-compact}
  With the notation~\eqref{eq:N-t-2-def},  writing $\|\cdot\|_{\mathcal{H}}$
  for the deterministic (unweighted) norm defined by~\eqref{eq:H-innprod-fourier}  and $\|\cdot\|_{\mathcal{H}_{Z}}$ for its weighted
  counterpart in Notation~\ref{not:norm-Z}, the $L^{2}$
  bound~\eqref{eq:dalang-bound} can be written compactly as
  \begin{equation}\label{eq:dalang-bound-compact}
    \bigl\|S\,\mathbf{1}_{[0,t]}\bigr\|_{\mathcal{H}_{Z}}^{2}
    \,\le\, \mathcal{N}_{t,2}(Z)\,
    \bigl\|S\,\mathbf{1}_{[0,t]}\bigr\|_{\mathcal{H}}^{2}.
  \end{equation}
  Indeed, the left-hand side of~\eqref{eq:dalang-bound-compact} is the
  $L^{2}(\Omega)$ norm of the stochastic integral by the
  isometry~\eqref{eq:isometry-MZ}, while the integral against $\mu$ on the
  right-hand side of~\eqref{eq:dalang-bound} is exactly
  $\|S\,\mathbf{1}_{[0,t]}\|_{\mathcal{H}}^{2}$
  by~\eqref{eq:H-innprod-fourier}.
\end{remark}

\begin{proof}[Proof of Theorem \ref{thm:dalang-integral}]
  Fix a standard mollifier $\psi\in C^{\infty}_{0}(\R^{d})$ with $\psi\ge 0$,
  $\int_{\R^{d}}\psi(x)\,dx=1$, $\mathrm{supp}(\psi)\subset\{|x|\le 1\}$, and
  set $\psi_{n}(x):=n^{d}\psi(nx)$.  Since $\psi_{n}\to\delta_{0}$ in
  $\mathcal{S}'(\R^{d})$, one has $\mathcal{F}\psi_{n}(\xi)\to 1$ pointwise and
  $|\mathcal{F}\psi_{n}(\xi)|\le 1$ for all $\xi$.  Define the mollified
  integrand
  \begin{equation}\label{eq:S-mollified}
    S_{n}(t) \,:=\, \psi_{n}*S(t), \qquad n\ge 1.
  \end{equation}
  Since $S(t)$ is a distribution with rapid decrease, $S_{n}(t)=\psi_{n}*S(t)$
  belongs to $\mathcal{S}(\R^{d})$ for each $t$.  Moreover, $S_{n}(t)\ge 0$
  because both $\psi_{n}\ge 0$ and $S(t)\ge 0$.  We prove in three steps that
  $\|S_{n}-S\|_{\mathcal{H}_{Z}}\to 0$, which yields $S\in \mathcal{H}_{Z}$ and the claimed
  identities.

  \smallskip
  \noindent
  {\it Step 1: $S_{n}\in P_{+,Z}$ with a uniform bound.}
  Let $\|\cdot\|_{+,Z}$ be the Walsh norm adapted to the modified martingale
  measure $M^{Z}$ of Definition~\ref{def:mg-Z}, i.e.\ the analogue
  of~\eqref{eq:walsh-norm-plus} in which the integrand is weighted by $Z$:
  \begin{equation}\label{eq:walsh-norm-plus-Z}
    \|S\|_{+,Z}^{2}
    \,:=\,\E\!\left[\int_{0}^{T}\!\int_{\R^{d}}\!\int_{\R^{d}}
      |S(s,x)|\,|Z(s,x)|\,\Lambda(x-y)\,|Z(s,y)|\,|S(s,y)|\,
      dx\,dy\,ds\right],
  \end{equation}
  and let $P_{+,Z}$ be the completion of the deterministic Schwartz functions
  under $\|\cdot\|_{+,Z}$.  Since $S_{n}(t)\in\mathcal{S}(\R^{d})$ and
  $S_{n}(t)\ge 0$, the absolute values in~\eqref{eq:walsh-norm-plus-Z} can be
  removed.  Applying the Cauchy--Schwarz inequality to bound
  $\E[|Z(s,x)||Z(s,y)|]\le\sup_{x'}\E[Z(s,x')^{2}]$, and then the
  Fourier identity~\eqref{eq:innprod-fourier}, we obtain
  \begin{equation}\label{eq:Sn-norm-bound}
    \|S_{n}\|_{+,Z}^{2}
    \,\le\,\Bigl(\sup_{0\le s\le T}\sup_{x\in\R^{d}}\E\bigl[Z(s,x)^{2}\bigr]
    \Bigr)
    \int_{0}^{T}\!\int_{\R^{d}}
    \bigl|\mathcal{F}S_{n}(t)(\xi)\bigr|^{2}\,\mu(d\xi)\,dt.
  \end{equation}
  Since $|\mathcal{F}S_{n}(t)(\xi)|=|\mathcal{F}\psi_{n}(\xi)|\,|\mathcal{F}S(t)(\xi)|\le
  |\mathcal{F}S(t)(\xi)|$, the right-hand side of~\eqref{eq:Sn-norm-bound} is
  bounded uniformly in $n$ by conditions~\eqref{eq:Z-L2-bound}
  and~\eqref{eq:dalang-cond-S}.  Hence $S_{n}\in P_{+,Z}\subset \mathcal{H}_{Z}$
  for every $n\ge 1$.

  \smallskip
  \noindent
  {\it Step 2: $\|S\|_{\mathcal{H}_{Z}}<\infty$.}
  Using the definition~\eqref{eq:norm-Z} of $\|\cdot\|_{\mathcal{H}_{Z}}$ and
  $|\mathcal{F}S_{n}(t)(\xi)|^{2}\nearrow|\mathcal{F}S(t)(\xi)|^{2}$
  as $n\to\infty$, Fatou's lemma gives
  \begin{equation*}
    \|S\|_{\mathcal{H}_{Z}}^{2}
    \,\le\,\liminf_{n\to\infty}\|S_{n}\|_{\mathcal{H}_{Z}}^{2}
    \,\le\,\liminf_{n\to\infty}\|S_{n}\|_{+,Z}^{2}
    \,<\,\infty,
  \end{equation*}
  where the second inequality uses $\|\cdot\|_{\mathcal{H}_{Z}}\le\|\cdot\|_{+,Z}$ (valid
  since $S_{n}\ge 0$) and the last one is the uniform bound from Step~1.

  \smallskip
  \noindent
  {\it Step 3: $\|S_{n}-S\|_{\mathcal{H}_{Z}}\to 0$.}
  From the definition~\eqref{eq:norm-Z} and
  $\mathcal{F}(S_{n}(t)-S(t))(\xi)=(\mathcal{F}\psi_{n}(\xi)-1)\,\mathcal{F}S(t)(\xi)$,
  \begin{equation}\label{eq:Sn-minus-S}
    \|S_{n}-S\|_{\mathcal{H}_{Z}}^{2}
    \,=\,\int_{0}^{T}\!\int_{\R^{d}}
    |\mathcal{F}\psi_{n}(\xi)-1|^{2}\,
    |\mathcal{F}S(t)(\xi)|^{2}\,\mu^{Z}_{t}(d\xi)\,dt.
  \end{equation}
  The integrand in~\eqref{eq:Sn-minus-S} converges pointwise to zero since
  $\mathcal{F}\psi_{n}(\xi)\to 1$, and is dominated by
  $4\,|\mathcal{F}S(t)(\xi)|^{2}$, whose integral is $4\|S\|_{\mathcal{H}_{Z}}^{2}<\infty$
  by Step~2.  The dominated convergence theorem therefore gives
  $\|S_{n}-S\|_{\mathcal{H}_{Z}}\to 0$.

  \smallskip
  \noindent
  Since $S_{n}\in \mathcal{H}_{Z}$ and $\|S_{n}-S\|_{\mathcal{H}_{Z}}\to 0$, the completeness of
  $\mathcal{H}_{Z}$ implies $S\in \mathcal{H}_{Z}$, and the
  isometry~\eqref{eq:dalang-isometry-thm} is~\eqref{eq:isometry-MZ} applied to
  $S\,\mathbf{1}_{[0,t]}$.  Finally, the bound~\eqref{eq:dalang-bound} is
  obtained by passing to the limit in the Step~1 estimate.  Since $S_{n}\ge 0$
  one has $\|\cdot\|_{\mathcal{H}_{Z}}\le\|\cdot\|_{+,Z}$, so
  restricting~\eqref{eq:Sn-norm-bound} to the interval $[0,t]$ gives
  \begin{equation}\label{eq:Sn-bound-t}
    \|S_{n}\,\mathbf{1}_{[0,t]}\|_{\mathcal{H}_{Z}}^{2}
    \,\le\,\Bigl(\sup_{0\le s\le t}\sup_{x\in\R^{d}}\E\bigl[Z(s,x)^{2}\bigr]\Bigr)
    \int_{0}^{t}\!\int_{\R^{d}}
    \bigl|\mathcal{F}_{x}S_{n}(s,\cdot)(\xi)\bigr|^{2}\,\mu(d\xi)\,ds.
  \end{equation}
  Here the supremum in the prefactor may be restricted to $[0,t]$ because, after
  restriction to $[0,t]$, only the times $s\in[0,t]$ enter the integral, and the
  Cauchy--Schwarz step bounding $\E[Z(s,x)Z(s,y)]\le\sup_{x'}\E[Z(s,x')^{2}]$ is
  applied for such $s$ only.
  As $n\to\infty$, the left-hand side of~\eqref{eq:Sn-bound-t} converges to
  $\|S\,\mathbf{1}_{[0,t]}\|_{\mathcal{H}_{Z}}^{2}
  =\E[\,|\int_{0}^{t}\!\int_{\R^{d}}S\,dM^{Z}|^{2}\,]$ (by Step~3
  and~\eqref{eq:isometry-MZ}), while
  $|\mathcal{F}_{x}S_{n}(s,\cdot)(\xi)|^{2}\nearrow|\mathcal{F}_{x}S(s,\cdot)(\xi)|^{2}$
  allows the monotone convergence theorem on the right-hand side.  Passing to
  the limit yields exactly~\eqref{eq:dalang-bound}.
\end{proof}

\subsection{A Burkholder--Davis--Gundy inequality for \texorpdfstring{$L^{p}$}{Lp} moments}
\label{sec:bdg-Lp}

The $L^{2}$ bound~\eqref{eq:dalang-bound} of
Theorem~\ref{thm:dalang-integral} controls the \emph{second} moment of the
modified stochastic integral $\int S\,dM^{Z}$ in terms of the deterministic
spectral integral~\eqref{eq:dalang-cond-S}. For the fixed point argument of
Section~\ref{sec:contraction}, however, the contraction has to be set up in a
space of processes with finite $p$-th moments, for a fixed $p\ge 2$. This is
also the natural setting in which one reads off the moment Lyapunov exponents
responsible for intermittency. We thus need the analogue
of~\eqref{eq:dalang-bound} for every $L^{p}(\Omega)$ norm. Such an estimate is
the colored-noise counterpart of Khoshnevisan's Burkholder--Davis--Gundy (BDG)
inequality (\cite[Proposition~4.4]{khoshnevisan:14:analysis}), which is stated
there for the space-time white noise. We adapt it to the spatially homogeneous
noise of Subsection~\ref{ssec:noise}, and to the modified martingale measure
$M^{Z}$ of Definition~\ref{def:mg-Z}.

To begin with, we introduce a new piece of notation.
Namely the role played by the prefactor $\mathcal{N}_{t,2}(Z)$
of~\eqref{eq:N-t-2-def} in the $L^{2}$ bound is now played by its $L^{p}$
upgrade.

\begin{notation}[$L^{p}$ moment functional]\label{not:N-t-p}
  For $p\ge 2$, an $(\mathcal{F}_{s})$-adapted process $Z$ with
  $Z(s,x)\in L^{p}(\Omega)$ for every $(s,x)$, and $t\ge 0$, we set
  \begin{equation}\label{eq:N-t-p-def}
    \mathcal{N}_{t,p}(Z)
    \,:=\,\sup_{0\le s\le t}\sup_{x\in\R^{d}}
      \bigl\|Z(s,x)\bigr\|_{L^{p}(\Omega)}^{2}
    \,=\,\sup_{0\le s\le t}\sup_{x\in\R^{d}}
      \bigl(\E\bigl[\,|Z(s,x)|^{p}\,\bigr]\bigr)^{2/p}.
  \end{equation}
  The requirement that $\mathcal{N}_{t,p}(Z)$ be \emph{finite} is the
  $L^{p}$ analogue of the second-moment bound~\eqref{eq:Z-L2-bound}. It will be assumed
   in Theorem~\ref{thm:bdg-Lp} below.
\end{notation}

\begin{remark}\label{rmk:N-t-p-consistency}
  We have kept the \emph{square} of the $L^{p}(\Omega)$ norm
  in~\eqref{eq:N-t-p-def} on purpose: for $p=2$ the functional
  $\mathcal{N}_{t,p}(Z)$ reduces exactly to
  $\mathcal{N}_{t,2}(Z)=\sup_{s\le t}\sup_{x}\E[Z(s,x)^{2}]$
  of~\eqref{eq:N-t-2-def}, so that~\eqref{eq:N-t-p-def} extends the earlier
  Notation~\ref{not:N-t-2}, and the BDG bound~\eqref{eq:bdg-Lp} below has the
  same algebraic shape as the $L^{2}$ bound~\eqref{eq:dalang-bound}. Moreover
  $p\mapsto\mathcal{N}_{t,p}(Z)$ is non-decreasing in $p$ by Jensen's
  inequality, and non-decreasing in $t$.
\end{remark}

The functional $\mathcal{N}_{t,p}$ singles out a natural space of admissible
weights, the one in which the fixed-point argument of
Section~\ref{sec:contraction} will be carried out.

\begin{definition}[The space $\mathcal{A}_{t,p}$]\label{def:A-t-p}
  For $p\ge 2$ and $t\ge 0$, let $\mathcal{A}_{t,p}$ be the set of all
  jointly measurable, $(\mathcal{F}_{s})$-adapted processes
  $Z=(Z(s,x))_{0\le s\le t,\,x\in\R^{d}}$ for which $\mathcal{N}_{t,p}(Z)<+\infty$
  (Notation~\ref{not:N-t-p}). In this context, observe that two processes coinciding
  $ds\,dx\,d\Prob$-almost everywhere are being identified. Equipped with the norm
  \begin{equation}\label{eq:A-norm}
    \|Z\|_{\mathcal{A}_{t,p}}
    \,:=\,\mathcal{N}_{t,p}(Z)^{1/2}
    \,=\,\sup_{0\le s\le t}\sup_{x\in\R^{d}}\|Z(s,x)\|_{L^{p}(\Omega)},
  \end{equation}
  $\mathcal{A}_{t,p}$ is a Banach space. Indeed~\eqref{eq:A-norm} is the norm of
  $L^{\infty}\bigl([0,t]\times\R^{d};L^{p}(\Omega)\bigr)$, which is easily seen to be complete.
\end{definition}

\begin{remark}\label{rmk:A-inclusions}
  By the monotonicity recorded in Remark~\ref{rmk:N-t-p-consistency}, the
  norm~\eqref{eq:A-norm} increases with $p$ and with $t$, so that
  $\mathcal{A}_{t,q}\subseteq\mathcal{A}_{t,p}$ for $q\ge p\ge 2$ and the
  restriction of any $Z\in\mathcal{A}_{t',p}$ to $[0,t]$ lies in
  $\mathcal{A}_{t,p}$ for $t'\ge t$. In
  particular $\mathcal{A}_{t,p}\subseteq\mathcal{A}_{t,2}$, so that every
  $Z\in\mathcal{A}_{t,p}$ satisfies the second-moment
  bound~\eqref{eq:Z-L2-bound} on $[0,t]$.
\end{remark}

\begin{theorem}[BDG inequality for $M^{Z}$]\label{thm:bdg-Lp}
  Let $p\ge 2$, $0\le t\le T$, let $Z\in\mathcal{A}_{t,p}$
  (Definition~\ref{def:A-t-p}), and let $t\mapsto S(t)$ be a deterministic
  non-negative integrand as in Theorem~\ref{thm:dalang-integral}, satisfying the
  integrability condition~\eqref{eq:dalang-cond-S}. Then the modified stochastic
  integral $\int S\,dM^{Z}$, which by~\eqref{eq:mg-Z-def} coincides with the
  Walsh integral $\int S\,Z\,dW$ of Definition~\ref{def:walsh-integral}, is
  well-defined, belongs to $L^{p}(\Omega)$, and satisfies
  \begin{equation}\label{eq:bdg-Lp}
    \left\|\int_{0}^{t}\!\int_{\R^{d}}S(s,y)\,M^{Z}(ds,dy)
    \right\|_{L^{p}(\Omega)}^{2}
    \,\le\,4p\,\mathcal{N}_{t,p}(Z)\,
    \int_{0}^{t}\!\int_{\R^{d}}
      \bigl|\mathcal{F}_{x}S(s,\cdot)(\xi)\bigr|^{2}\,\mu(d\xi)\,ds.
  \end{equation}
  Equivalently, with the deterministic norm $\|\cdot\|_{\mathcal{H}}$
  of~\eqref{eq:H-innprod-fourier}, the bound~\eqref{eq:bdg-Lp} reads
  \begin{equation}\label{eq:bdg-Lp-compact}
    \left\|\int_{0}^{t}\!\int_{\R^{d}}S\,dM^{Z}
    \right\|_{L^{p}(\Omega)}^{2}
    \,\le\,4p\,\mathcal{N}_{t,p}(Z)\,
    \bigl\|S\,\mathbf{1}_{[0,t]}\bigr\|_{\mathcal{H}}^{2}.
  \end{equation}
\end{theorem}

\begin{remark}\label{rmk:no-homogeneity}
  In contrast with Theorem~\ref{thm:dalang-integral}, the moment
  bound~\eqref{eq:bdg-Lp} does \emph{not} require the spatial homogeneity
  Assumption~\ref{hyp:spatial-hom}. It is read off directly from the quadratic
  variation of the Walsh integral $\int S\,Z\,dW$, and the effective spectral
  measure $\mu^{Z}$ of~\eqref{eq:muZ-def} is bypassed. Homogeneity is needed
  for the exact isometry~\eqref{eq:dalang-isometry-thm} only, and not for the
  inequality driving the contraction. This is why $\mathcal{A}_{t,p}$, whose
  elements are merely measurable and adapted, is a legitimate domain for the
  stochastic integral.
\end{remark}

\begin{remark}\label{rmk:bdg-vs-L2}
  For $p=2$, the right-hand side of~\eqref{eq:bdg-Lp} coincides with that
  of the $L^{2}$ bound~\eqref{eq:dalang-bound} up to the universal constant
  $4p=8$; the sharp isometric constant $1$ of~\eqref{eq:dalang-bound} is of
  course lost, since the BDG inequality replaces the exact second-moment
  identity~\eqref{eq:dalang-isometry-thm} by the martingale inequality of
  Step~1 below. What matters for Section~\ref{sec:contraction} is the
  structure of~\eqref{eq:bdg-Lp-compact}, rather than the constant. Namely the
  random weight $Z$ enters through the scalar $\mathcal{N}_{t,p}(Z)$ only,
  while the deterministic kernel $S$ enters through
  $\|S\,\mathbf{1}_{[0,t]}\|_{\mathcal{H}}$ only.
\end{remark}

\begin{proof}[Proof of Theorem~\ref{thm:bdg-Lp}]

  \noindent
  {\it Well-definedness.}
  Since $Z\in\mathcal{A}_{t,p}\subseteq\mathcal{A}_{t,2}$
  (Remark~\ref{rmk:A-inclusions}) and $S\ge 0$
  satisfies~\eqref{eq:dalang-cond-S}, the product $S\,Z$ is predictable and, by
  the Cauchy--Schwarz step bounding $\E[|Z(s,x)||Z(s,y)|]\le\mathcal{N}_{t,2}(Z)$
  in the Walsh norm~\eqref{eq:walsh-norm-plus},
  \begin{equation}\label{eq:SZ-walsh-finite}
    \|S\,Z\|_{+}^{2}
    \,\le\,\mathcal{N}_{t,2}(Z)\,\bigl\|S\,\mathbf{1}_{[0,t]}\bigr\|_{\mathcal{H}}^{2}
    \,<\,+\infty,
  \end{equation}
  so $S\,Z\in P_{+}$ and $\int S\,Z\,dW$ is a well-defined Walsh integral by
  Proposition~\ref{prop:walsh-isometry}.

  \smallskip
  \noindent
  {\it Step~0: the martingale and its bracket.}
  Set
  \begin{equation}\label{eq:Mt-mart-def}
    M_{t}
    \,:=\,\int_{0}^{t}\!\int_{\R^{d}}S(s,y)\,M^{Z}(ds,dy)
    \,=\,\int_{0}^{t}\!\int_{\R^{d}}S(s,y)\,Z(s,y)\,W(ds,dy),
  \end{equation}
  the second equality being the definition~\eqref{eq:mg-Z-def} of the modified
  martingale measure. Then $t\mapsto M_{t}$ is a continuous, square-integrable
  $(\mathcal{F}_{t})$-martingale. Thanks to the covariation
  formula~\eqref{eq:mg-covariation} for the worthy martingale measure $M$ (or,
  equivalently, by the Walsh isometry~\eqref{eq:walsh-isometry} read at the
  level of brackets), its quadratic variation is the random measure
  \begin{equation}\label{eq:bracket-Mz}
    \langle M\rangle_{t}
    \,=\,\int_{0}^{t}\!\!\int_{\R^{d}}\!\int_{\R^{d}}
      S(s,x)\,Z(s,x)\,\Lambda(x-y)\,Z(s,y)\,S(s,y)\,
      dx\,dy\,ds.
  \end{equation}

  \smallskip
  \noindent
  {\it Step 1: BDG for the continuous martingale $M$.}
  We apply the Burkholder--Davis--Gundy inequality for continuous $L^{2}$
  martingales in the sharp form of~\cite[Theorem~B.1]{khoshnevisan:14:analysis},
  namely $\|M_{t}\|_{L^{p}(\Omega)}^{2}\le 4p\,\|\langle M\rangle_{t}\|_{L^{p/2}(\Omega)}$.
  With the notation $\|Y\|_{q}:=\|Y\|_{L^{q}(\Omega)}$, this is
  \begin{equation}\label{eq:bdg-step1}
    \bigl\|M_{t}\bigr\|_{p}^{2}
    \,\le\,4p\,\bigl\|\langle M\rangle_{t}\bigr\|_{p/2}.
  \end{equation}

  \smallskip
  \noindent
  {\it Step 2: Minkowski's inequality on the bracket.}
  Since $S(s,\cdot)\ge 0$ and the covariance $\Lambda\ge 0$ (recall that the
  noise model of Subsection~\ref{ssec:noise} takes
  $\Lambda:\R^{d}\to\R_{+}$ non-negative in~\eqref{eq:noise-cov}, the same
  property used for the worthiness of $M$ in
  Remark~\ref{rmk:covariation-measure}), the deterministic
  weight $S(s,x)\,\Lambda(x-y)\,S(s,y)$ in~\eqref{eq:bracket-Mz} is non-negative.
  Minkowski's integral inequality in $L^{p/2}(\Omega)$ (valid since $p/2\ge 1$)
  therefore lets us pull the norm inside the $dx\,dy\,ds$ integral,
  \begin{equation}\label{eq:bdg-step2}
    \bigl\|\langle M\rangle_{t}\bigr\|_{p/2}
    \,\le\,\int_{0}^{t}\!\!\int_{\R^{d}}\!\int_{\R^{d}}
      S(s,x)\,\Lambda(x-y)\,S(s,y)\,
      \bigl\|Z(s,x)\,Z(s,y)\bigr\|_{p/2}\,
      dx\,dy\,ds.
  \end{equation}
  Moreover, owing to the Cauchy--Schwarz inequality on $\Omega$, we have
  \begin{equation*}
\|Z(s,x)\,Z(s,y)\|_{p/2}\le\|Z(s,x)\|_{p}\,\|Z(s,y)\|_{p} ,
\end{equation*}
and by the
  definition~\eqref{eq:N-t-p-def} of $\mathcal{N}_{t,p}(Z)$ each factor is at
  most $\mathcal{N}_{t,p}(Z)^{1/2}$ for $s\le t$, so that
  \begin{equation}\label{eq:Z-prod-bound}
    \bigl\|Z(s,x)\,Z(s,y)\bigr\|_{p/2}
    \,\le\,\sup_{x'\in\R^{d}}\|Z(s,x')\|_{p}^{2}
    \,\le\,\mathcal{N}_{t,p}(Z),
    \qquad 0\le s\le t.
  \end{equation}

  \smallskip
  \noindent
  {\it Step 3: Recognizing the deterministic norm.}
  Substituting~\eqref{eq:Z-prod-bound} into~\eqref{eq:bdg-step2} pulls the
  scalar $\mathcal{N}_{t,p}(Z)$ out of the integral and leaves a purely
  deterministic double integral against $\Lambda$, which is exactly the
  $\mathcal{H}$-norm of $S\,\mathbf{1}_{[0,t]}$ in the
  form~\eqref{eq:H-innprod} (equivalently~\eqref{eq:H-innprod-fourier}):
  \begin{align}\label{eq:bdg-step3}
    \bigl\|\langle M\rangle_{t}\bigr\|_{p/2}
    &\,\le\,\mathcal{N}_{t,p}(Z)
      \int_{0}^{t}\!\!\int_{\R^{d}}\!\int_{\R^{d}}
      S(s,x)\,\Lambda(x-y)\,S(s,y)\,dx\,dy\,ds \nonumber\\
    &\,=\,\mathcal{N}_{t,p}(Z)\,
      \bigl\|S\,\mathbf{1}_{[0,t]}\bigr\|_{\mathcal{H}}^{2}
    \,=\,\mathcal{N}_{t,p}(Z)
      \int_{0}^{t}\!\int_{\R^{d}}
      \bigl|\mathcal{F}_{x}S(s,\cdot)(\xi)\bigr|^{2}\,\mu(d\xi)\,ds.
  \end{align}
  The last quantity is finite by the integrability
  condition~\eqref{eq:dalang-cond-S} and $Z\in\mathcal{A}_{t,p}$, so
  $M_{t}\in L^{p}(\Omega)$. Inserting the
  chain~\eqref{eq:bdg-step3} into~\eqref{eq:bdg-step1} yields both the compact
  form~\eqref{eq:bdg-Lp-compact} and, after expanding
  $\|S\,\mathbf{1}_{[0,t]}\|_{\mathcal{H}}^{2}$ through the Fourier
  identity~\eqref{eq:H-innprod-fourier}, the bound~\eqref{eq:bdg-Lp}.
\end{proof}

In the contraction argument of Section~\ref{sec:contraction}, the integrand is
not a single function $S$. It is the time-space shifted heat kernel
$S^{(\tau,z)}(s,y)=p_{\tau-s}(z-y)$, with one copy for each output point
$(\tau,z)$, and the integral~\eqref{eq:bdg-Lp} is evaluated at that point.
Applying Theorem~\ref{thm:bdg-Lp} uniformly over $(\tau,z)$, the moment bound
becomes a genuine \emph{stability} statement. Namely the map sending a weight
$Z$ to its stochastic convolution preserves the space $\mathcal{A}_{T,p}$ of
Definition~\ref{def:A-t-p}.

\begin{theorem}[Stability of $\mathcal{A}_{T,p}$]\label{thm:stability}
  Let $p\ge 2$, $T>0$, and let
  $S=\bigl(S^{(\tau,z)}\bigr)_{(\tau,z)\in[0,T]\times\R^{d}}$ be a family of
  deterministic non-negative integrands, each as in
  Theorem~\ref{thm:dalang-integral}, jointly measurable in $(\tau,z,s,y)$, and
  such that
  \begin{equation}\label{eq:S-star-norm}
    \|S\|_{*,T}^{2}
    \,:=\,\sup_{0\le\tau\le T}\sup_{z\in\R^{d}}
      \bigl\|S^{(\tau,z)}\,\mathbf{1}_{[0,\tau]}\bigr\|_{\mathcal{H}}^{2}
    \,<\,+\infty.
  \end{equation}
  Then for every $Z\in\mathcal{A}_{T,p}$ the random field
  \begin{equation}\label{eq:I-S-def}
    \mathcal{I}_{S}Z(\tau,z)
    \,:=\,\int_{0}^{\tau}\!\int_{\R^{d}}S^{(\tau,z)}(s,y)\,Z(s,y)\,W(ds,dy),
    \qquad (\tau,z)\in[0,T]\times\R^{d},
  \end{equation}
  which by~\eqref{eq:mg-Z-def} is the modified stochastic integral
  $\int_{0}^{\tau}\!\int_{\R^{d}}S^{(\tau,z)}\,dM^{Z}$ to which
  Theorem~\ref{thm:bdg-Lp} applies,
  admits a jointly measurable, adapted version belonging to
  $\mathcal{A}_{T,p}$, with
  \begin{equation}\label{eq:I-S-bound}
    \mathcal{N}_{T,p}\bigl(\mathcal{I}_{S}Z\bigr)
    \,\le\,4p\,\|S\|_{*,T}^{2}\,\mathcal{N}_{T,p}(Z).
  \end{equation}
  In particular $Z\mapsto\mathcal{I}_{S}Z$ is a bounded linear operator on
  $\mathcal{A}_{T,p}$ with operator norm at most $(4p)^{1/2}\,\|S\|_{*,T}$.
\end{theorem}

\begin{proof}[Proof of Theorem~\ref{thm:stability}]
  Fix $Z\in\mathcal{A}_{T,p}$. For each $(\tau,z)\in[0,T]\times\R^{d}$, the
  integrand $S^{(\tau,z)}$ satisfies the integrability
  condition~\eqref{eq:dalang-cond-S} (since
  $\|S^{(\tau,z)}\mathbf{1}_{[0,\tau]}\|_{\mathcal{H}}\le\|S\|_{*,T}<\infty$), and
  $Z\in\mathcal{A}_{\tau,p}$ by restriction (Remark~\ref{rmk:A-inclusions}). So
  Theorem~\ref{thm:bdg-Lp}, applied on the interval $[0,\tau]$ in the compact
  form~\eqref{eq:bdg-Lp-compact}, gives
  \begin{equation}\label{eq:I-S-pointwise}
    \bigl\|\mathcal{I}_{S}Z(\tau,z)\bigr\|_{p}^{2}
    \,\le\,4p\,\mathcal{N}_{\tau,p}(Z)\,
      \bigl\|S^{(\tau,z)}\,\mathbf{1}_{[0,\tau]}\bigr\|_{\mathcal{H}}^{2}
    \,\le\,4p\,\|S\|_{*,T}^{2}\,\mathcal{N}_{T,p}(Z),
  \end{equation}
  where the second inequality uses $\mathcal{N}_{\tau,p}(Z)\le\mathcal{N}_{T,p}(Z)$
  for $\tau\le T$ (Remark~\ref{rmk:N-t-p-consistency}) and the
  definition~\eqref{eq:S-star-norm} of $\|S\|_{*,T}$. The right-hand side
  of~\eqref{eq:I-S-pointwise} is independent of $(\tau,z)$; taking the supremum
  over $(\tau,z)\in[0,T]\times\R^{d}$ on the left yields
  precisely~\eqref{eq:I-S-bound}. Since the kernel family is jointly measurable,
  the random field~\eqref{eq:I-S-def} admits a jointly measurable, adapted
  version (obtained, e.g., as the limit of the Walsh--Riemann sums defining the
  integral, each of which is jointly measurable and adapted); by~\eqref{eq:I-S-bound}
  this version belongs to $\mathcal{A}_{T,p}$. Finally, $Z\mapsto\mathcal{I}_{S}Z$
  is linear, since the integral~\eqref{eq:I-S-def} is linear in the
  integrand $S^{(\tau,z)}Z$ and hence in $Z$, and~\eqref{eq:I-S-bound} is exactly
  the operator-norm bound
  $\|\mathcal{I}_{S}\|_{\mathcal{A}_{T,p}\to\mathcal{A}_{T,p}}
  \le(4p)^{1/2}\|S\|_{*,T}$.
\end{proof}

\begin{remark}\label{rmk:stability-contraction}
  With the heat-kernel choice $S^{(\tau,z)}(s,y)=p_{\tau-s}(z-y)$, the operator
  $\mathcal{I}_{S}$ of~\eqref{eq:I-S-def} is exactly the stochastic convolution,
  and Theorem~\ref{thm:stability} says it maps the solution space
  $\mathcal{A}_{T,p}$ into itself, with a norm governed by the deterministic
  quantity $\|S\|_{*,T}$ alone. Tuning the horizon $T$ (or inserting an
  exponential weight in time) so that $4p\,\|S\|_{*,T}^{2}<1$ is what closes
  the fixed point in Section~\ref{sec:contraction}.
\end{remark}

\section[Additive noise case]{Additive noise case}\label{sec:additive}

We begin with the simplest instance of the stochastic heat equation, namely
the one where the noise enters \emph{additively}. This is the tame model already
announced in item~\textbf{(i)} below the PAM~\eqref{eq:pam}. It lacks the
multiplicative feedback responsible for intermittency. However, its definition
already requires to check that the heat kernel is an admissible integrand for
the noise, and to give a sense to the resulting stochastic convolution in
$L^{2}(\Omega)$. The purpose of this section is to isolate those difficulties
from the fixed point argument of Section~\ref{sec:contraction}. The solution
constructed here will also reappear there, as an intermediate object.

\subsection{The additive stochastic heat equation}

Replacing the multiplicative term $\beta\,u\,\dot{W}$ in the PAM~\eqref{eq:pam}
by the forcing $\beta\,\dot{W}$ alone yields the \emph{additive} stochastic heat
equation
\begin{equation}\label{eq:additive-she}
  \partial_{t} u(t,x)
  \,=\, \tfrac{1}{2}\Delta\, u(t,x) \,+\, \beta\,\dot{W}(t,x),
  \qquad u(0,x)=u_{0}(x),
\end{equation}
where $\dot{W}$ is the spatially homogeneous, white-in-time Gaussian noise of
Subsection~\ref{ssec:noise}, with spatial covariance $\Lambda$ and spectral
measure $\mu=\mathcal{F}\Lambda$ determined by~\eqref{eq:noise-cov}, and
$u_{0}:\R^{d}\to\R$ is a bounded measurable initial datum. Throughout this
section the noise is assumed to satisfy Dalang's
condition~\eqref{eq:dalang-cond}, that is Assumption~\ref{hyp:spectral}.

\subsection{Mild formulation}

Since the product $u\,\dot{W}$ has disappeared, equation~\eqref{eq:additive-she}
is \emph{linear} with a prescribed random forcing. The variation of constants
(Duhamel) principle thus applies. The heat semigroup with
kernel~$p_{t}$~\eqref{eq:heat-kernel-explicit} propagates the initial datum,
while the noise is accumulated and transported over $[0,t]$. Specializing the
mild PAM~\eqref{eq:mild-pam} to the unit multiplier $u(s,y)\equiv 1$ inside the
stochastic integral, we call $u$ a \emph{mild solution}
of~\eqref{eq:additive-she} if
\begin{equation}\label{eq:additive-mild}
  u(t,x) \,=\, P_{t}u_{0}(x)
  \,+\, \beta\int_{0}^{t}\!\int_{\R^{d}}
    p_{t-s}(x-y)\,W(ds,dy),
\end{equation}
the stochastic integral being understood in the Walsh--Dalang sense of
Section~\ref{sec:integration}. Here $P_{t}=e^{t\frac{1}{2}\Delta}$ is again the
heat semigroup with Gaussian kernel~$p_{t}$ introduced in~\eqref{eq:heat-kernel-explicit}, so
that $P_{t}u_{0}=p_{t}*u_{0}$. The first term $P_{t}u_{0}(x)$ is an ordinary
convolution, the classical solution of the homogeneous heat equation, and is
bounded whenever $u_{0}$ is. All the probabilistic content
of~\eqref{eq:additive-mild} is therefore carried by the second term, the
\emph{stochastic convolution}
\begin{equation}\label{eq:stoch-conv}
  v(t,x) \,:=\, \beta\int_{0}^{t}\!\int_{\R^{d}}
    p_{t-s}(x-y)\,W(ds,dy),
\end{equation}
so that existence and uniqueness for~\eqref{eq:additive-she} reduce to showing
that $v$ is a well-defined random field.

\begin{remark}\label{rmk:mild-weak}
We have phrased existence and uniqueness in the \emph{mild} sense
of~\eqref{eq:additive-mild}. This is no loss of generality: under the standing
assumptions, the mild solution is also a \emph{weak} (distributional) solution
of~\eqref{eq:additive-she}, and the two notions coincide. A proof of this
equivalence can be found in~\cite[Section~5.3]{khoshnevisan:14:analysis}.
\end{remark}

\subsection{The stochastic convolution is well defined}

The integrand in~\eqref{eq:stoch-conv} is the \emph{deterministic} function
$S(s,y)=p_{t-s}(x-y)$, so~\eqref{eq:stoch-conv} is an integral against the
martingale measure $M$ associated with the noise $W$
(Section~\ref{sec:integration}), equivalently the modified measure $M^{Z}$ of
Definition~\ref{def:mg-Z} with the trivial weight $Z\equiv 1$. For this weight
the correlation~\eqref{eq:gs-def} is $g_{s}\equiv 1$, so that the effective
spectral measure~\eqref{eq:muZ-def} reduces to $\mu^{Z}_{s}=\mu$ and the
admissibility of $S$ is governed by Dalang's condition for $\mu$ itself. We
record this computation, in the form reused for general weights in
Section~\ref{sec:contraction}, in the following proposition.

\begin{proposition}[The heat kernel as integrand]\label{ex:heat-integrand}
  Fix $(t,x)\in(0,T]\times\R^{d}$. Under Dalang's condition
  (Assumption~\ref{hyp:spectral}), the deterministic integrand
  $S(s,y)=p_{t-s}(x-y)$ is admissible for the trivial weight $Z\equiv 1$, so
  that the stochastic integral
  \begin{equation}\label{eq:stoch-int-heat}
    \int_{0}^{t}\!\int_{\R^{d}}p_{t-s}(x-y)\,W(ds,dy)
  \end{equation}
  is well-defined.
\end{proposition}

\begin{proof}
  Consider the deterministic integrand
  \begin{equation}\label{eq:S-heat}
    S(s,y) \,:=\, p_{t-s}(x-y),
    \qquad 0\le s< t,
  \end{equation}
  where $p_{t}$ is the heat kernel~\eqref{eq:heat-kernel-explicit}.  Since $S$ is
  deterministic, we may dispense with the weighted theory and work directly with
  the extended norm $\|\cdot\|_{\mathcal{H}_{\omega}}$ of Definition~\ref{def:extended-norm}:
  by the deterministic Fourier identity~\eqref{eq:innprod-fourier}
  (Remark~\ref{rmk:deterministic-fourier}, with $g=h=S\,\mathbf{1}_{[0,t]}$),
  \begin{equation}\label{eq:S-Homega-norm}
    \bigl\|S\,\mathbf{1}_{[0,t]}\bigr\|_{\mathcal{H}_{\omega}}^{2}
    \,=\,\int_{0}^{t}\!\int_{\R^{d}}
    \bigl|\mathcal{F}_{x}S(s,\cdot)(\xi)\bigr|^{2}\,\mu(d\xi)\,ds.
  \end{equation}
  The spatial Fourier transform of the integrand is
  \begin{equation}\label{eq:FT-heat-integrand}
    \mathcal{F}_{x}S(s,\cdot)(\xi)
    \,=\,e^{-\imath\xi\cdot x}\,e^{-\frac12(t-s)|\xi|^{2}}.
  \end{equation}
  Inserting~\eqref{eq:FT-heat-integrand} into~\eqref{eq:S-Homega-norm} and
  integrating in $s\in[0,t]$,
  \begin{equation}\label{eq:dalang-heat-check}
    \bigl\|S\,\mathbf{1}_{[0,t]}\bigr\|_{\mathcal{H}_{\omega}}^{2}
    \,=\,\int_{\R^{d}}\frac{1-e^{-t|\xi|^{2}}}{|\xi|^{2}}\,\mu(d\xi)
    \,\lesssim\,\int_{\R^{d}}\frac{\mu(d\xi)}{1+|\xi|^{2}},
  \end{equation}
  where the last inequality uses $\frac{1-e^{-t|\xi|^{2}}}{|\xi|^{2}}\le
  \frac{1}{|\xi|^{2}}\lesssim\frac{1}{1+|\xi|^{2}}$ for $|\xi|\ge 1$, and
  $\frac{1-e^{-t|\xi|^{2}}}{|\xi|^{2}}\le t$ for $|\xi|\le 1$.  The
  right-hand side of~\eqref{eq:dalang-heat-check} is finite by
  Assumption~\ref{hyp:spectral} (equation~\eqref{eq:dalang-cond}).  Hence
  $S\,\mathbf{1}_{[0,t]}\in\mathcal{H}_{\omega}$ under Dalang's condition alone, and
  Theorem~\ref{thm:extended-norm} shows that the stochastic
  integral~\eqref{eq:stoch-int-heat} is well-defined, with second moment given
  by the isometry~\eqref{eq:extended-isometry}.
\end{proof}

Up to the factor $\beta$, the integral~\eqref{eq:stoch-int-heat} is exactly the
stochastic convolution~\eqref{eq:stoch-conv} driving the additive
equation~\eqref{eq:additive-she}. The general weighted case $Z=u$, needed when
the noise is multiplied by the unknown solution, is treated in
Section~\ref{sec:contraction}.

\begin{proposition}[Existence and uniqueness for the additive equation]\label{prop:additive-existence}
  Assume Dalang's condition (Assumption~\ref{hyp:spectral}) and let $u_{0}$ be a
  bounded measurable initial datum. Then:
  \begin{enumerate}
  \item The stochastic convolution~\eqref{eq:stoch-conv} is, for every
  $(t,x)\in(0,T]\times\R^{d}$, a well-defined centered Gaussian random variable
  in $L^{2}(\Omega)$, with variance
  \begin{equation}\label{eq:stoch-conv-var}
    \E\bigl[v(t,x)^{2}\bigr]
    \,=\,\beta^{2}\int_{\R^{d}}\frac{1-e^{-t|\xi|^{2}}}{|\xi|^{2}}\,\mu(d\xi)
    \,<\,\infty,
  \end{equation}
  finite and independent of~$x$;
  \item The additive equation~\eqref{eq:additive-she} admits a unique mild
  solution, given explicitly by~\eqref{eq:additive-mild}. It is a Gaussian field,
  bounded in $L^{2}(\Omega)$ uniformly on $[0,T]\times\R^{d}$.
  \end{enumerate}
\end{proposition}

\begin{proof}
  By Proposition~\ref{ex:heat-integrand}, the deterministic integrand
  $S(s,y)=p_{t-s}(x-y)$ satisfies $S\,\mathbf{1}_{[0,t]}\in\mathcal{H}_{\omega}$,
  so the integral~\eqref{eq:stoch-conv} is well-defined by
  Theorem~\ref{thm:extended-norm}. Its second moment is the
  isometry~\eqref{eq:extended-isometry} evaluated at $S$:
  \begin{equation*}
    \E\bigl[v(t,x)^{2}\bigr]
    \,=\,\beta^{2}\bigl\|S\,\mathbf{1}_{[0,t]}\bigr\|_{\mathcal{H}_{\omega}}^{2}
    \,=\,\beta^{2}\int_{\R^{d}}\frac{1-e^{-t|\xi|^{2}}}{|\xi|^{2}}\,\mu(d\xi),
  \end{equation*}
  the last equality being the computation~\eqref{eq:dalang-heat-check}; it is
  finite by Dalang's condition~\eqref{eq:dalang-cond} and, depending on $t$
  alone, independent of~$x$. As an $L^{2}(\Omega)$-limit of Wiener integrals of
  deterministic integrands, $v(t,x)$ is centered Gaussian. This proves~(i).

  \smallskip
  \noindent
  For~(ii), define $u$ by the right-hand side of~\eqref{eq:additive-mild}, that
  is, $u(t,x)=P_{t}u_{0}(x)+v(t,x)$. The first term is the bounded deterministic
  function $P_{t}u_{0}$, and the second is the well-defined random field $v$ just
  constructed, so $u$ exists as a random field and solves~\eqref{eq:additive-she}
  in the mild sense; this gives existence. Conversely, the forcing
  in~\eqref{eq:additive-she} being prescribed (it does not involve the unknown),
  the mild formulation~\eqref{eq:additive-mild} determines $u(t,x)$ explicitly
  for every $(t,x)$, so the mild solution is unique. Finally, $u$ is Gaussian as
  the sum of the deterministic term $P_{t}u_{0}$ and the Gaussian field $v$, and
  the uniform $L^{2}$ bound follows from~\eqref{eq:stoch-conv-var} together with
  the boundedness of $P_{t}u_{0}$.
\end{proof}

\subsection{H\"older continuity of the stochastic convolution}\label{ssec:holder-convolution}

Proposition~\ref{prop:additive-existence} constructs the stochastic
convolution $v$ of~\eqref{eq:stoch-conv} as an $L^{2}(\Omega)$-valued field. It
says nothing about the regularity of the trajectories $(t,x)\mapsto v(t,x)$.
Under Dalang's condition~\eqref{eq:dalang-cond} alone, one can only show that
$v$ is continuous in $L^{2}(\Omega)$. A genuine \emph{pathwise} regularity
requires the small amount of extra integrability already alluded to for the
noise itself in Assumption~\ref{hyp:spectral-holder}. We recall that assumption
here for convenience.

\smallskip

\noindent
{\it Strong Dalang condition.} Throughout this subsection we assume that the
spectral measure $\mu$ satisfies the strong Dalang condition of
Assumption~\ref{hyp:spectral-holder}: there exists $\alpha\in(0,1)$ such that
\begin{equation}\label{eq:strong-dalang-recall}
  \int_{\R^{d}}\frac{\mu(d\xi)}{1+|\xi|^{2(1-\alpha)}}\,<\,\infty,
\end{equation}
which is~\eqref{eq:dalang-holder}. As noted after
Assumption~\ref{hyp:spectral-holder}, relation~\eqref{eq:strong-dalang-recall} is stronger
than Dalang's condition~\eqref{eq:dalang-cond} (to which it reduces as
$\alpha\to 0$), and the gap between the two is precisely the price for pathwise
regularity. The parameter $\alpha$ will turn out to be the spatial H\"older
exponent. Owing to the parabolic scaling of the heat equation, it is halved in
the time variable.

In order to quantify pathwise H\"older regularity, we will rely on the
following space-time version of Kolmogorov's continuity criterion, which we state
in the parabolic form adapted to the heat equation
(see~\cite[Chapter~3]{dalang:09:minicourse}).

\begin{theorem}[Kolmogorov's criterion, parabolic form]\label{thm:kolmogorov}
  Equip $[0,T]\times\R^{d}$ with the parabolic distance
  \begin{equation}\label{eq:parabolic-metric}
    \rho\bigl((t,x),(s,y)\bigr)
    \,:=\,\bigl(|t-s|+|x-y|^{2}\bigr)^{1/2}.
  \end{equation}
  Let $\{X(t,x):(t,x)\in[0,T]\times\R^{d}\}$ be a real-valued random field and
  suppose there exist $p\in[2,\infty)$, an exponent $\gamma\in(0,1]$ and a
  constant $C$ such that
  \begin{equation}\label{eq:kolmo-hyp}
    \E\bigl[\,|X(t,x)-X(s,y)|^{p}\,\bigr]
    \,\le\, C\,\rho\bigl((t,x),(s,y)\bigr)^{\gamma p},
    \qquad (t,x),(s,y)\in[0,T]\times\R^{d}.
  \end{equation}
  If $\gamma p>d+2$, then $X$ admits a modification whose trajectories are, on
  every compact subset of $[0,T]\times\R^{d}$, almost surely H\"older continuous
  of every order $\theta<\gamma-\frac{d+2}{p}$ with respect to $\rho$;
  equivalently, of order $\theta$ in the space variable and $\theta/2$ in the
  time variable. The corresponding H\"older seminorms have finite moments of all
  orders.
\end{theorem}

\noindent
Observe that the exponent $d+2$ in the threshold $\gamma p>d+2$ is the volume-growth
dimension of the parabolic metric~\eqref{eq:parabolic-metric}: a ball of radius
$r$ for $\rho$ has Lebesgue measure of order $r^{2}\cdot r^{d}=r^{d+2}$, the
factor $r^{2}$ accounting for the time direction. We can now state the main
result of this subsection.

\begin{theorem}[H\"older continuity of the stochastic convolution]
  \label{thm:holder-convolution}
  Assume the strong Dalang condition~\eqref{eq:strong-dalang-recall}
  (Assumption~\ref{hyp:spectral-holder}) for some $\alpha\in(0,1)$, and let
  $u_{0}$ be a bounded measurable initial datum. Then the stochastic
  convolution $v$ of~\eqref{eq:stoch-conv} admits a modification whose
  trajectories are, on every compact subset of $[0,T]\times\R^{d}$, almost surely
  jointly H\"older continuous of every order $\theta_{x}<\alpha$ in the space
  variable and every order $\theta_{t}<\tfrac{\alpha}{2}$ in the time variable.
\end{theorem}

\noindent
The proof follows a single principle: because the integrand
in~\eqref{eq:stoch-conv} is deterministic, every increment of $v$ is a centered
Gaussian variable, so its $p$-th moment is a power of its variance, and that
variance is computed exactly through the Fourier identity~\eqref{eq:innprod-fourier}
as in the proof of Proposition~\ref{ex:heat-integrand}. The space and time
increments are estimated separately, each through a decomposition of the
heat-kernel integrand of the type used in the existence theory for the
multiplicative equation.

\begin{proof}[Proof of Theorem~\ref{thm:holder-convolution}]
  Since the integrand $\beta\,p_{t-s}(x-y)\,\mathbf{1}_{[0,t]}(s)$
  in~\eqref{eq:stoch-conv} is deterministic, for each $(t,x)$ the random
  variable $v(t,x)$ is a centered Gaussian; so is every increment
  $v(t,x)-v(s,y)$, being a finite linear combination of Wiener integrals of
  deterministic integrands. Hence for every $p\ge 2$,
  \begin{equation}\label{eq:gauss-increment-moment}
    \E\bigl[\,|v(t,x)-v(s,y)|^{p}\,\bigr]
    \,=\, c_{p}\,\E\bigl[\,|v(t,x)-v(s,y)|^{2}\,\bigr]^{p/2} ,
  \end{equation}
  where the constant $c_{p}$ can be expressed as $c_{p}=\E[|N|^{p}]$ for a standard Gaussian random variable 
  $N\sim\mathcal{N}(0,1)$.
  By Kolmogorov's criterion (Theorem~\ref{thm:kolmogorov}) it therefore suffices
  to control the second moment of the increments by the parabolic
  distance~\eqref{eq:parabolic-metric}. Otherwise stated we have to prove
  \begin{equation}\label{eq:increment-target}
    \E\bigl[\,|v(t,x)-v(s,y)|^{2}\,\bigr]
    \,\le\, C\,\bigl(|t-s|+|x-y|^{2}\bigr)^{\alpha}.
  \end{equation}
  Indeed, inserting~\eqref{eq:increment-target} into~\eqref{eq:gauss-increment-moment}
  yields~\eqref{eq:kolmo-hyp} with $\gamma=\alpha$ for every $p$, and since the
  Gaussian field $v$ has moments of all orders, letting $p\to\infty$ in
  Theorem~\ref{thm:kolmogorov} produces every space H\"older order
  $\theta_{x}<\alpha$ and every time H\"older order $\theta_{t}<\tfrac{\alpha}{2}$,
  which is the assertion of the theorem.

  Kolmogorov's criterion is a local statement: only the
  bound~\eqref{eq:kolmo-hyp} for points at parabolic distance at most $1$ enters,
  so we may and do assume $|x-y|\le 1$ throughout (no restriction on
  $|t-s|\in[0,T]$ is needed). In addition, a simple application of the triangle inequality in $L^{2}(\Omega)$ yields
  \begin{equation}\label{eq:increment-triangle}
    \E\bigl[|v(t,x)-v(s,y)|^{2}\bigr]^{1/2}
    \,\le\,\E\bigl[|v(t,x)-v(t,y)|^{2}\bigr]^{1/2}
          +\E\bigl[|v(t,y)-v(s,y)|^{2}\bigr]^{1/2}.
  \end{equation}
  Therefore~\eqref{eq:increment-target} follows once we establish the spatial bound
  \begin{equation}\label{eq:space-increment-target}
    \E\bigl[\,|v(t,x)-v(t,y)|^{2}\,\bigr]\,\le\, C\,|x-y|^{2\alpha}
  \end{equation}
  as well as the temporal bound
  \begin{equation}\label{eq:time-increment-target}
    \E\bigl[\,|v(t,x)-v(s,x)|^{2}\,\bigr]\,\le\, C\,|t-s|^{\alpha} .
  \end{equation}
  We will prove \eqref{eq:space-increment-target} and \eqref{eq:time-increment-target} in the two steps
  below.

  \smallskip

  \noindent
  {\it Step 1: Spatial increment.} The increment $v(t,x)-v(t,y)$ is the
  stochastic integral of the deterministic integrand
  $\beta[p_{t-\cdot}(x-\cdot)-p_{t-\cdot}(y-\cdot)]\mathbf{1}_{[0,t]}$.
  By the isometry~\eqref{eq:extended-isometry} in its deterministic form
  (Remark~\ref{rmk:deterministic-fourier}) and the Fourier
  identity~\eqref{eq:innprod-fourier}, exactly as in the proof of
  Proposition~\ref{ex:heat-integrand}, we get
  \begin{equation}\label{eq:space-incr-fourier}
    \E\bigl[\,|v(t,x)-v(t,y)|^{2}\,\bigr]
    \,=\,\beta^{2}\int_{0}^{t}\!\int_{\R^{d}}
      \bigl|\mathcal{F}_{x}\bigl[p_{t-s}(x-\cdot)-p_{t-s}(y-\cdot)\bigr](\xi)\bigr|^{2}
      \,\mu(d\xi)\,ds.
  \end{equation}
  Owing to the form of the spatial Fourier transform~\eqref{eq:FT-heat-integrand}
  we have
  \begin{equation*}
\mathcal{F}_{x}[p_{t-s}(x-\cdot)](\xi)=e^{-\imath\xi\cdot x}\,e^{-\frac12(t-s)|\xi|^{2}} ,
\end{equation*}
  so the integrand in~\eqref{eq:space-incr-fourier} equals
  \begin{equation}\label{eq:space-incr-mod}
    \bigl|e^{-\imath\xi\cdot x}-e^{-\imath\xi\cdot y}\bigr|^{2}\,e^{-(t-s)|\xi|^{2}}
    \,=\,\bigl|1-e^{\imath\xi\cdot(x-y)}\bigr|^{2}\,e^{-(t-s)|\xi|^{2}}.
  \end{equation}
  Integrating in $s\in[0,t]$ as in~\eqref{eq:dalang-heat-check}, with
  $\int_{0}^{t}e^{-(t-s)|\xi|^{2}}\,ds=\frac{1-e^{-t|\xi|^{2}}}{|\xi|^{2}}
  \le\frac{1}{|\xi|^{2}}$, we obtain
  \begin{equation}\label{eq:space-incr-bound}
    \E\bigl[\,|v(t,x)-v(t,y)|^{2}\,\bigr]
    \,\le\,\beta^{2}\int_{\R^{d}}
      \frac{|1-e^{\imath\xi\cdot(x-y)}|^{2}}{|\xi|^{2}}\,\mu(d\xi).
  \end{equation}
  We now invoke the elementary inequality: for all $\xi\in\R^{d}$ and all
  $h\in\R^{d}$ with $|h|\le 1$,
  \begin{equation}\label{eq:elem-space}
    \frac{|1-e^{\imath\xi\cdot h}|^{2}}{|\xi|^{2}}
    \,\le\, C_{\alpha}\,\frac{|h|^{2\alpha}}{1+|\xi|^{2(1-\alpha)}}.
  \end{equation}
  To prove~\eqref{eq:elem-space}, we start from the two elementary bounds
  $|1-e^{\imath\xi\cdot h}|\le|\xi|\,|h|$ and $|1-e^{\imath\xi\cdot h}|\le 2$, and we
  split the ratio of interest into two factors,
  \begin{equation}\label{eq:elem-space-split}
    \frac{|1-e^{\imath\xi\cdot h}|^{2}}{|\xi|^{2}}
    \,=\,\left(\frac{|1-e^{\imath\xi\cdot h}|^{2}}{|\xi|^{2}}\right)^{1-\alpha}
        \left(\frac{|1-e^{\imath\xi\cdot h}|^{2}}{|\xi|^{2}}\right)^{\alpha}.
  \end{equation}
  In the first factor of~\eqref{eq:elem-space-split} we use $|1-e^{\imath\xi\cdot h}|\le 2$,
  while in the second factor we use $|1-e^{\imath\xi\cdot h}|\le|\xi|\,|h|$. This yields
  \begin{equation}\label{eq:elem-space-factors}
    \left(\frac{|1-e^{\imath\xi\cdot h}|^{2}}{|\xi|^{2}}\right)^{1-\alpha}
    \,\le\,\frac{4^{1-\alpha}}{|\xi|^{2(1-\alpha)}},
    \qquad
    \left(\frac{|1-e^{\imath\xi\cdot h}|^{2}}{|\xi|^{2}}\right)^{\alpha}
    \,\le\,|h|^{2\alpha}.
  \end{equation}
  Multiplying the two bounds in~\eqref{eq:elem-space-factors} gives 
  \begin{equation}\label{eq:elem-space-product}
    \frac{|1-e^{\imath\xi\cdot h}|^{2}}{|\xi|^{2}}
    \,\le\, C_{\alpha}\,\frac{|h|^{2\alpha}}{|\xi|^{2(1-\alpha)}}.
  \end{equation}
  We now compare the right-hand side of~\eqref{eq:elem-space-product} with that
  of~\eqref{eq:elem-space} according to the size of $|\xi|$. For $|\xi|\ge 1$ we
  have $|\xi|^{2(1-\alpha)}\ge\tfrac12\bigl(1+|\xi|^{2(1-\alpha)}\bigr)$, so
  \eqref{eq:elem-space-product} gives
  \eqref{eq:elem-space}. For $|\xi|\le 1$ we
  instead bound the ratio on the left hand side of ~\eqref{eq:elem-space-product} directly, using
  $|1-e^{\imath\xi\cdot h}|\le|\xi|\,|h|$ and $|h|\le 1$,
  \begin{equation}\label{eq:elem-space-small}
    \frac{|1-e^{\imath\xi\cdot h}|^{2}}{|\xi|^{2}}\,\le\,|h|^{2}\,\le\,|h|^{2\alpha},
  \end{equation}
  which together with $1+|\xi|^{2(1-\alpha)}\le 2$
  again yields~\eqref{eq:elem-space}. This proves the elementary
  inequality~\eqref{eq:elem-space}. Now inserting~\eqref{eq:elem-space} with $h=x-y$
  into~\eqref{eq:space-incr-bound}, we end up with
  \begin{equation*}
    \E\bigl[\,|v(t,x)-v(t,y)|^{2}\,\bigr]
    \,\le\, C_{\alpha}\,\beta^{2}\,|x-y|^{2\alpha}
      \int_{\R^{d}}\frac{\mu(d\xi)}{1+|\xi|^{2(1-\alpha)}}
    \,=\, C\,|x-y|^{2\alpha},
  \end{equation*}
  finite by the strong Dalang condition~\eqref{eq:strong-dalang-recall}. This
  proves~\eqref{eq:space-increment-target}.

  \smallskip

  \noindent
  {\it Step 2: Temporal increment.} Fix $0\le s\le t\le T$. We split the time
  increment into a near-diagonal and an off-diagonal part,
  $v(t,x)-v(s,x)=\beta\,(I_{1}+I_{2})$, where
  \begin{eqnarray}
    I_{1} & := & \int_{s}^{t}\!\int_{\R^{d}}p_{t-r}(x-y)\,W(dr,dy),
      \label{eq:time-I1}\\
    I_{2} & := & \int_{0}^{s}\!\int_{\R^{d}}
      \bigl[p_{t-r}(x-y)-p_{s-r}(x-y)\bigr]\,W(dr,dy).
      \label{eq:time-I2}
  \end{eqnarray}
  Since $I_{1}$ and $I_{2}$ integrate the noise over the disjoint time slabs
  $[s,t]$ and $[0,s]$, and $W$ has independent increments in time, they are
  independent centered Gaussian variables. Hence we obtain
  \begin{equation}\label{eq:time-split}
    \E\bigl[\,|v(t,x)-v(s,x)|^{2}\,\bigr]
    \,=\,\beta^{2}\bigl(\E[I_{1}^{2}]+\E[I_{2}^{2}]\bigr).
  \end{equation}

  \noindent
  {\it The near-diagonal term $I_{1}$.} Similarly to what we did for~\eqref{eq:space-incr-fourier}, we apply the isometry~\eqref{eq:extended-isometry}
  and the Fourier identity~\eqref{eq:innprod-fourier}, using
  $|\mathcal{F}_{x}p_{t-r}(x-\cdot)(\xi)|^{2}=e^{-(t-r)|\xi|^{2}}$
  from~\eqref{eq:FT-heat-integrand}. We discover that
  \begin{equation}\label{eq:time-I1-bound}
    \E[I_{1}^{2}]
    \,=\,\int_{s}^{t}\!\int_{\R^{d}}e^{-(t-r)|\xi|^{2}}\,\mu(d\xi)\,dr
    \,=\,\int_{\R^{d}}\frac{1-e^{-(t-s)|\xi|^{2}}}{|\xi|^{2}}\,\mu(d\xi),
  \end{equation}
  the last equality being obtained by an elementary integration in the time variable.
  We will show how to bound the right hand side of~\eqref{eq:time-I1-bound} below.

  \smallskip

  \noindent
  {\it The off-diagonal term $I_{2}$.} As for $I_{1}$, invoke~\eqref{eq:extended-isometry}
  and~\eqref{eq:innprod-fourier} to get
  \begin{equation}\label{eq:time-I2-fourier}
    \E[I_{2}^{2}]
    \,=\,\int_{0}^{s}\!\int_{\R^{d}}
      \bigl|\mathcal{F}_{x}\bigl[p_{t-r}(x-\cdot)-p_{s-r}(x-\cdot)\bigr](\xi)\bigr|^{2}
      \,\mu(d\xi)\,dr.
  \end{equation}
  Using $\mathcal{F}_{x}p_{\tau}(x-\cdot)(\xi)=e^{-\imath\xi\cdot x}e^{-\frac12\tau|\xi|^{2}}$,
  the squared modulus of the bracket factors as
  \begin{equation}\label{eq:time-I2-mod}
    \bigl(e^{-\frac12(s-r)|\xi|^{2}}-e^{-\frac12(t-r)|\xi|^{2}}\bigr)^{2}
    \,=\, e^{-(s-r)|\xi|^{2}}\bigl(1-e^{-\frac12(t-s)|\xi|^{2}}\bigr)^{2}.
  \end{equation}
  Integrating in $r\in[0,s]$, with
  $\int_{0}^{s}e^{-(s-r)|\xi|^{2}}\,dr=(1-e^{-s|\xi|^{2}})/|\xi|^{2}\le 1/|\xi|^{2}$,
  and bounding the square by itself,
  $(1-e^{-\frac12(t-s)|\xi|^{2}})^{2}\le 1-e^{-\frac12(t-s)|\xi|^{2}}$, we get
  \begin{equation}\label{eq:time-I2-bound}
    \E[I_{2}^{2}]
    \,\le\,\int_{\R^{d}}\frac{1-e^{-\frac12(t-s)|\xi|^{2}}}{|\xi|^{2}}\,\mu(d\xi).
  \end{equation}
  
  The right hand sides of~\eqref{eq:time-I1-bound} and~\eqref{eq:time-I2-bound} are of the same kind. Both are controlled by a
  second elementary inequality: for every $c>0$, every $0\le t-s\le T$, and every
  $\xi\in\R^{d}$,
  \begin{equation}\label{eq:elem-time}
    \frac{1-e^{-c(t-s)|\xi|^{2}}}{|\xi|^{2}}
    \,\le\, C_{\alpha,c,T}\,\frac{(t-s)^{\alpha}}{1+|\xi|^{2(1-\alpha)}}.
  \end{equation}
  To prove~\eqref{eq:elem-time}, we start from the elementary bound
  \begin{equation}\label{eq:elem-time-exp}
    1-e^{-a}\,\le\,\min(1,a)\,\le\, a^{\alpha},
    \qquad a\ge 0,\ \alpha\in(0,1).
  \end{equation}
  We now apply~\eqref{eq:elem-time-exp} according to the size of $|\xi|$, with
  $a=c(t-s)|\xi|^{2}$. For $|\xi|\ge 1$, the third bound in~\eqref{eq:elem-time-exp}
  gives $1-e^{-c(t-s)|\xi|^{2}}\le c^{\alpha}(t-s)^{\alpha}|\xi|^{2\alpha}$, hence
  \begin{equation}\label{eq:elem-time-large}
    \frac{1-e^{-c(t-s)|\xi|^{2}}}{|\xi|^{2}}
    \,\le\, c^{\alpha}\,\frac{(t-s)^{\alpha}}{|\xi|^{2(1-\alpha)}}
    \,\le\, 2\,c^{\alpha}\,\frac{(t-s)^{\alpha}}{1+|\xi|^{2(1-\alpha)}},
  \end{equation}
  where in the last step we used $|\xi|^{2(1-\alpha)}\ge\tfrac12(1+|\xi|^{2(1-\alpha)})$
  for $|\xi|\ge 1$. For $|\xi|\le 1$, the first bound in~\eqref{eq:elem-time-exp}
  together with $t-s\le T$ gives instead
  \begin{equation}\label{eq:elem-time-small}
    \frac{1-e^{-c(t-s)|\xi|^{2}}}{|\xi|^{2}}
    \,\le\, c\,(t-s)
    \,\le\, c\,T^{1-\alpha}\,(t-s)^{\alpha}
    \,\le\, 2\,c\,T^{1-\alpha}\,\frac{(t-s)^{\alpha}}{1+|\xi|^{2(1-\alpha)}},
  \end{equation}
  using $1+|\xi|^{2(1-\alpha)}\le 2$ in the last step. Bounds
  \eqref{eq:elem-time-large} and~\eqref{eq:elem-time-small} together
  yield~\eqref{eq:elem-time} with $C_{\alpha,c,T}=2\max(c^{\alpha},c\,T^{1-\alpha})$.
  
  We can now conclude:
plugging~\eqref{eq:elem-time} (with $c=1$ and $c=\tfrac12$) into~\eqref{eq:time-I1-bound} and~\eqref{eq:time-I2-bound}, then
  inserting into~\eqref{eq:time-split}, we end up with
  \begin{equation*}
    \E\bigl[\,|v(t,x)-v(s,x)|^{2}\,\bigr]
    \,\le\, C\,(t-s)^{\alpha}\int_{\R^{d}}\frac{\mu(d\xi)}{1+|\xi|^{2(1-\alpha)}}
    \,=\, C\,(t-s)^{\alpha} ,
  \end{equation*}
 where the constant $C$ is finite by~\eqref{eq:strong-dalang-recall}. This
  proves~\eqref{eq:time-increment-target} and completes the proof.
\end{proof}

\noindent
The H\"older orders $\alpha$ in space and $\tfrac{\alpha}{2}$ in time exhibit the
parabolic scaling of the heat equation: time counts as two space derivatives. In
the model case of space-time white noise in dimension $d=1$, where
$\mu(d\xi)=d\xi$, the strong Dalang condition~\eqref{eq:strong-dalang-recall}
holds for every $\alpha<\tfrac12$, and Theorem~\ref{thm:holder-convolution}
recovers the classical regularity of the stochastic heat equation: H\"older
continuous of every order below $\tfrac12$ in space and below $\tfrac14$ in time.

\begin{remark}\label{rmk:holder-solution}
  The deterministic part $P_{t}u_{0}$ of the mild solution~\eqref{eq:additive-mild}
  is, for $t>0$, of class $\mathcal{C}^{\infty}$ in $(t,x)$ by the smoothing
  property of the heat semigroup. Consequently the mild solution
  $u=P_{t}u_{0}+v$ of~\eqref{eq:additive-she} inherits, on every compact subset
  of $(0,T]\times\R^{d}$, the H\"older regularity of $v$ stated in
  Theorem~\ref{thm:holder-convolution}: order $\theta_{x}<\alpha$ in space and
  $\theta_{t}<\tfrac{\alpha}{2}$ in time.
\end{remark}

\section[PAM case]{Existence and uniqueness of PAM}\label{sec:contraction}

Throughout this section we fix $p\ge 2$, a coupling constant $\beta\in\R$, and a
bounded measurable initial datum $u_{0}:\R^{d}\to\R$, and we assume Dalang's
condition (Assumption~\ref{hyp:spectral}). A \emph{mild solution} of the PAM on
$[0,T]$ is a process $u\in\mathcal{A}_{T,p}$ (Definition~\ref{def:A-t-p})
satisfying the mild equation~\eqref{eq:mild-pam} for every
$(t,x)\in[0,T]\times\R^{d}$.

Most existence and uniqueness results in this context (see e.g.
Dalang~\cite[Theorem~13]{dalang:99:extending}) construct such a solution by a
Picard iteration. We give here a slight variation, our construction being based
on a single application of the Banach fixed point theorem. The whole argument
rests on the one quantitative input of Subsection~\ref{sec:bdg-Lp}, namely that
$\mathcal{I}_{S}$ contracts $\mathcal{A}_{\tau,p}$ once the horizon $\tau$ is
small.

\subsection{The solution map and its contraction constant}

Let us specialize the kernel family of Theorem~\ref{thm:stability} to the heat
kernel. Namely, for $(t,x)\in[0,T]\times\R^{d}$ we put
\begin{equation}\label{eq:heat-family}
  S^{(t,x)}(s,y)\,:=\,p_{t-s}(x-y),
  \qquad 0\le s\le t,\ y\in\R^{d},
\end{equation}
with $p$ the heat kernel~\eqref{eq:heat-kernel-explicit}. Each $S^{(t,x)}$ is a
non-negative integrand of the type required by
Theorem~\ref{thm:dalang-integral}, and
$(t,x,s,y)\mapsto S^{(t,x)}(s,y)\,\mathbf{1}_{[0,t]}(s)$ is jointly measurable.
For this family the operator $\mathcal{I}_{S}$ of~\eqref{eq:I-S-def} is the
stochastic convolution
\begin{equation}\label{eq:conv-operator}
  \mathcal{I}_{S}Z(t,x)
  \,=\,\int_{0}^{t}\!\int_{\R^{d}}p_{t-s}(x-y)\,Z(s,y)\,W(ds,dy),
\end{equation}
so that the mild equation~\eqref{eq:mild-pam} reads compactly as
$u=P_{\cdot}u_{0}+\beta\,\mathcal{I}_{S}u$.

We first prove that the size of $\mathcal{I}_{S}$ over the horizon $\tau$ is
measured by the constant $\|S\|_{*,\tau}^{2}$ of~\eqref{eq:S-star-norm}.

\begin{lemma}\label{lem:Phi}
  For $\tau\ge 0$ set
  \begin{equation}\label{eq:Phi-def}
    \Phi(\tau)\,:=\,\int_{\R^{d}}\frac{1-e^{-\tau|\xi|^{2}}}{|\xi|^{2}}\,
      \mu(d\xi).
  \end{equation}
  Under Assumption~\ref{hyp:spectral} the function $\Phi$ is finite,
  non-decreasing, and $\Phi(\tau)\to 0$ as $\tau\downarrow 0$. Moreover, for the
  heat family~\eqref{eq:heat-family} we have
  \begin{equation}\label{eq:S-star-Phi}
    \|S\|_{*,\tau}^{2}
    \,=\,\sup_{0\le t\le\tau}\sup_{x\in\R^{d}}
      \bigl\|S^{(t,x)}\,\mathbf{1}_{[0,t]}\bigr\|_{\mathcal{H}}^{2}
    \,=\,\Phi(\tau).
  \end{equation}
\end{lemma}

\begin{proof}
  By the Fourier identity~\eqref{eq:H-innprod-fourier} and the
  transform~\eqref{eq:FT-heat-integrand}, we get
  \begin{equation}\label{eq:S-tx-Phi}
    \bigl\|S^{(t,x)}\,\mathbf{1}_{[0,t]}\bigr\|_{\mathcal{H}}^{2}
    \,=\,\int_{0}^{t}\!\int_{\R^{d}}e^{-(t-s)|\xi|^{2}}\,\mu(d\xi)\,ds
    \,=\,\Phi(t),
  \end{equation}
  where the last equality is easily obtained by a Fubini type argument. In
  addition, the right hand side of~\eqref{eq:S-tx-Phi} is independent of $x$ and
  non-decreasing in $t$, which gives~\eqref{eq:S-star-Phi}. Finiteness of
  $\Phi(\tau)$ then follows from~\eqref{eq:dalang-heat-check}. Namely, one has
  $$
  \Phi(\tau)\lesssim\int_{\R^{d}}\frac{\mu(d\xi)}{1+|\xi|^{2}}<\infty ,
  $$ 
  by Dalang's condition~\eqref{eq:dalang-cond}.  The integrand in~\eqref{eq:Phi-def} is
  non-negative and non-decreasing in $\tau$, so $\Phi$ is non-decreasing; for
  $\tau\le 1$ it is dominated by $\min(1,|\xi|^{-2})\lesssim(1+|\xi|^{2})^{-1}\in
  L^{1}(\mu)$ and tends to $0$ pointwise as $\tau\downarrow 0$, so
  $\Phi(\tau)\to 0$ by dominated convergence.
\end{proof}

We can now introduce the solution map.

\begin{definition}[Solution map]\label{def:Gamma}
  For $\tau\in(0,T]$ and $Z\in\mathcal{A}_{\tau,p}$, let $\Gamma_{\tau}Z$ be the
  process
  \begin{equation}\label{eq:Gamma-def}
    \Gamma_{\tau}Z(t,x)
    \,:=\,P_{t}u_{0}(x)\,+\,\beta\,\mathcal{I}_{S}Z(t,x),
    \qquad (t,x)\in[0,\tau]\times\R^{d}.
  \end{equation}
\end{definition}

\begin{proposition}[$\Gamma_{\tau}$ is an affine contraction]\label{prop:Gamma-contraction}
  Let $\tau\in(0,T]$. Then $\Gamma_{\tau}$ maps $\mathcal{A}_{\tau,p}$ into
  itself, and for all $Z_{1},Z_{2}\in\mathcal{A}_{\tau,p}$,
  \begin{equation}\label{eq:Gamma-lip}
    \bigl\|\Gamma_{\tau}Z_{1}-\Gamma_{\tau}Z_{2}\bigr\|_{\mathcal{A}_{\tau,p}}
    \,\le\,2\,|\beta|\,\sqrt{p\,\Phi(\tau)}\;
      \bigl\|Z_{1}-Z_{2}\bigr\|_{\mathcal{A}_{\tau,p}}.
  \end{equation}
  In particular $\Gamma_{\tau}$ is a contraction of $\mathcal{A}_{\tau,p}$ as
  soon as
  \begin{equation}\label{eq:tau-smallness}
    4\,p\,\beta^{2}\,\Phi(\tau)\,<\,1.
  \end{equation}
\end{proposition}

\begin{proof}
  \noindent
  {\it Stability.} The free term $f(t,x):=P_{t}u_{0}(x)$ is deterministic,
  jointly measurable, and bounded by $\|u_{0}\|_{\infty}$ (the heat kernel being a
  probability density). Therefore we get that $f\in\mathcal{A}_{\tau,p}$ and
  \begin{equation*}
\mathcal{N}_{\tau,p}(f)=\sup_{t\le\tau}\sup_{x}|P_{t}u_{0}(x)|^{2}
  \le\|u_{0}\|_{\infty}^{2}<\infty .
\end{equation*}
In addition, owing to  Lemma~\ref{lem:Phi}, $\|S\|_{*,\tau}=\Phi(\tau)^{1/2}<\infty$.
  Hence Theorem~\ref{thm:stability} gives $\mathcal{I}_{S}Z\in\mathcal{A}_{\tau,p}$.
  Since $\mathcal{A}_{\tau,p}$ is a vector space, we obtain
  \begin{equation*}
    \Gamma_{\tau}Z\,=\,f+\beta\,\mathcal{I}_{S}Z\,\in\,\mathcal{A}_{\tau,p}.
  \end{equation*}

  \smallskip
  \noindent
  {\it Contraction.} The free term $f$ cancels in the difference, and $\mathcal{I}_{S}$
  is linear (Theorem~\ref{thm:stability}), so
  $\Gamma_{\tau}Z_{1}-\Gamma_{\tau}Z_{2}=\beta\,\mathcal{I}_{S}(Z_{1}-Z_{2})$.
  By the operator bound~\eqref{eq:I-S-bound} and the
  identity~\eqref{eq:S-star-Phi}, we thus get
  \begin{align*}
    \bigl\|\Gamma_{\tau}Z_{1}-\Gamma_{\tau}Z_{2}\bigr\|_{\mathcal{A}_{\tau,p}}^{2}
    &\,=\,\beta^{2}\,\mathcal{N}_{\tau,p}\bigl(\mathcal{I}_{S}(Z_{1}-Z_{2})\bigr)
    \,\le\,4\,p\,\beta^{2}\,\|S\|_{*,\tau}^{2}\,
      \mathcal{N}_{\tau,p}(Z_{1}-Z_{2}) \\
    &\,=\,4\,p\,\beta^{2}\,\Phi(\tau)\,
      \bigl\|Z_{1}-Z_{2}\bigr\|_{\mathcal{A}_{\tau,p}}^{2},
  \end{align*}
  which is~\eqref{eq:Gamma-lip}. Moreover, the Lipschitz constant is $<1$ exactly
  under~\eqref{eq:tau-smallness}.
\end{proof}

\subsection{Local existence and uniqueness}

The contraction constant in~\eqref{eq:Gamma-lip} tends to $0$ with $\tau$, by
Lemma~\ref{lem:Phi}. This is all the contraction argument needs.

\begin{theorem}[Local well-posedness]\label{thm:local-wp}
  Assume Dalang's condition (Assumption~\ref{hyp:spectral}) and let $u_{0}$ be
  bounded measurable. There is a horizon $\tau_{*}>0$, depending only on $p$,
  $\beta$ and the spectral measure $\mu$ (and \emph{not} on $u_{0}$), such that
  for every $\tau\le\tau_{*}$ the PAM~\eqref{eq:mild-pam} has a unique mild
  solution $u\in\mathcal{A}_{\tau,p}$ on $[0,\tau]$, where we recall that $\mathcal{A}_{\tau,p}$ is the space of
  Definition~\ref{def:A-t-p}.
\end{theorem}

\begin{proof}
  Since $\Phi(\tau)\to 0$ as $\tau\downarrow 0$ (Lemma~\ref{lem:Phi}), there is a
  largest $\tau_{*}>0$ with $4p\beta^{2}\Phi(\tau_{*})\le\tfrac12$, and $\tau_{*}$
  depends only on $p,\beta,\mu$. For every $\tau\le\tau_{*}$ the
  smallness~\eqref{eq:tau-smallness} holds, so by
  Proposition~\ref{prop:Gamma-contraction} the map $\Gamma_{\tau}$ is a
  contraction (of Lipschitz constant at most $1/\sqrt{2}$) of the Banach space
  $\mathcal{A}_{\tau,p}$ (Definition~\ref{def:A-t-p}). By the Banach fixed point
  theorem it has a unique fixed point $u\in\mathcal{A}_{\tau,p}$, and the fixed
  point equation $u=\Gamma_{\tau}u$ is precisely the mild
  equation~\eqref{eq:mild-pam} on $[0,\tau]$. Uniqueness of the fixed point is
  uniqueness of the mild solution in $\mathcal{A}_{\tau,p}$.
\end{proof}

\subsection{Globalization by patching}

The horizon $\tau_{*}$ of Theorem~\ref{thm:local-wp} depends only on
$p,\beta,\mu$, not on the initial datum. This uniformity is what allows the local
construction to be iterated over consecutive windows of equal length and the
pieces patched into a global solution.

\begin{theorem}[Global existence and uniqueness for the PAM]\label{thm:global-wp}
  Assume Dalang's condition (Assumption~\ref{hyp:spectral}) and let $u_{0}$ be
  bounded measurable. Then for every $T>0$ and every $p\ge 2$ the
  PAM~\eqref{eq:mild-pam} has a unique mild solution $u\in\mathcal{A}_{T,p}$.
\end{theorem}

\begin{proof}
  Fix $\tau=\tau_{*}$ from Theorem~\ref{thm:local-wp} and let
  $N=\lceil T/\tau\rceil$, so that $[0,T]\subseteq\bigcup_{k=0}^{N-1}
  [k\tau,(k+1)\tau]$. We build the solution one window at a time.

  \smallskip
  \noindent
  {\it Step 1: The restarted map.}
  On a window $[a,a+\tau]$, given an $\mathcal{F}_{a}$-measurable random initial
  datum $w$ with $\|w\|_{\infty,p}:=\sup_{x}\|w(x)\|_{p}<\infty$, define for a
  process $Z$ on that window
  \begin{equation}\label{eq:restart-map}
    \Gamma^{a,w}Z(t,x)
    \,:=\,P_{t-a}w(x)
      \,+\,\beta\int_{a}^{t}\!\int_{\R^{d}}p_{t-s}(x-y)\,Z(s,y)\,W(ds,dy),
    \quad t\in[a,a+\tau].
  \end{equation}
  Two facts make $\Gamma^{a,w}$ a contraction with the \emph{same} constant as on
  $[0,\tau]$. First, the free term $P_{t-a}w$ is $\mathcal{F}_{a}$-measurable,
  hence adapted on $[a,a+\tau]$, and Minkowski's inequality gives
  $\sup_{t}\sup_{x}\|P_{t-a}w(x)\|_{p}\le\|w\|_{\infty,p}$ (again because the heat
  kernel is a probability density). Second, the convolution over the window has
  the same size as on $[0,\tau]$, namely
  \begin{equation*}
    \sup_{a\le t\le a+\tau}\sup_{x\in\R^{d}}
      \bigl\|p_{t-\cdot}(x-\cdot)\,\mathbf{1}_{[a,t]}\bigr\|_{\mathcal{H}}^{2}
    \,=\,\sup_{a\le t\le a+\tau}\Phi(t-a)\,=\,\Phi(\tau).
  \end{equation*}
  Note that the identity above is obtained by the very same
  computation~\eqref{eq:dalang-heat-check}, since it depends only
  on the elapsed time $t-a$ and the noise is homogeneous in time. Hence the proof
  of Proposition~\ref{prop:Gamma-contraction} applies verbatim (with
  $\mathcal{A}_{\tau,p}$ replaced by its time-translate to $[a,a+\tau]$), and for
  $\tau=\tau_{*}$ the map $\Gamma^{a,w}$ has a unique fixed point on
  $[a,a+\tau]$.

  \smallskip
  \noindent
  {\it Step 2: Patching.}
  Apply Step~1 successively. On $[0,\tau]$ with $w=u_{0}$ we recover the unique
  solution $u$ of Theorem~\ref{thm:local-wp}. Suppose $u$ has been built up to
  time $k\tau$ with $\sup_{t\le k\tau}\sup_{x}\|u(t,x)\|_{p}<\infty$. Set
  $w:=u(k\tau,\cdot)$, which is $\mathcal{F}_{k\tau}$-measurable with
  $\|w\|_{\infty,p}\le\mathcal{N}_{k\tau,p}(u)^{1/2}<\infty$. Then applying Step~1 on
  $[k\tau,(k+1)\tau]$, we produce the unique continuation. The concatenation solves
  the mild equation~\eqref{eq:mild-pam} on $[0,(k+1)\tau]$. Indeed, since $u$
  already solves~\eqref{eq:mild-pam} up to time $k\tau$, the restart datum is
  \begin{equation}\label{eq:restart-datum}
    u(k\tau,z)\,=\,P_{k\tau}u_{0}(z)\,+\,\beta\int_{0}^{k\tau}\!\int_{\R^{d}}
      p_{k\tau-s}(z-y)\,u(s,y)\,W(ds,dy).
  \end{equation}
  Applying the operator $P_{t-k\tau}$, i.e. integrating~\eqref{eq:restart-datum}
  against $p_{t-k\tau}(x-\cdot)$, the restarted free term expands as
  \begin{align*}
  &  P_{t-k\tau}u(k\tau,\cdot)(x)\\
    &\,=\,P_{t-k\tau}P_{k\tau}u_{0}(x)\,+\,\beta\int_{0}^{k\tau}\!\int_{\R^{d}}
      \Bigl(\int_{\R^{d}}p_{t-k\tau}(x-z)\,p_{k\tau-s}(z-y)\,dz\Bigr)
      u(s,y)\,W(ds,dy)\\
    &\,=\,P_{t}u_{0}(x)\,+\,\beta\int_{0}^{k\tau}\!\int_{\R^{d}}
      p_{t-s}(x-y)\,u(s,y)\,W(ds,dy),
  \end{align*}
  where the first equality moves $P_{t-k\tau}$ through the stochastic integral by
  the stochastic Fubini theorem (see Proposition~\ref{prop:stoch-fubini},
  applied with $G=\R^{d}$ and $\nu(dz)=p_{t-k\tau}(x-z)\,dz$), and the second uses the semigroup property
  $P_{t-k\tau}P_{k\tau}=P_{t}$ on the free term together with the
  Chapman--Kolmogorov identity
  \begin{equation}\label{eq:chapman-kolmogorov}
    \int_{\R^{d}}p_{t-k\tau}(x-z)\,p_{k\tau-s}(z-y)\,dz\,=\,p_{t-s}(x-y)
  \end{equation}
  on the stochastic term. Hence the two stochastic integrals
  in~\eqref{eq:restart-map} over $[0,k\tau]$ and $[k\tau,t]$ recombine into the
  single integral of~\eqref{eq:mild-pam} over $[0,t]$. After $N$ steps we obtain $u\in\mathcal{A}_{T,p}$, its
  $\mathcal{A}_{T,p}$-norm being the maximum of the $N$ finite window norms.

  \smallskip
  \noindent
  {\it Step 3: Uniqueness.}
  Let $u,\tilde{u}\in\mathcal{A}_{T,p}$ both solve~\eqref{eq:mild-pam}. On
  $[0,\tau]$ both are fixed points of $\Gamma_{\tau}$, hence equal by
  Theorem~\ref{thm:local-wp}; in particular $u(\tau,\cdot)=\tilde{u}(\tau,\cdot)$.
  Assuming $u=\tilde{u}$ up to time $k\tau$, both restrict on
  $[k\tau,(k+1)\tau]$ to a fixed point of the same restarted map
  $\Gamma^{k\tau,\,u(k\tau,\cdot)}$, so they agree there too. After $N$ steps, we get
  $u=\tilde{u}$ on $[0,T]$.
\end{proof}

\begin{remark}\label{rmk:contraction-vs-picard}
  The only place where the noise enters is through the scalar $\Phi(\tau)$
  of~\eqref{eq:Phi-def}, via the stability estimate~\eqref{eq:I-S-bound}. Compared
  with Dalang's Picard scheme~\cite[Theorem~13]{dalang:99:extending}, the
  contraction argument replaces the explicit iterates $u_{n}$ by a single
  fixed point statement. The price to pay is that one works on a small horizon
  $\tau_{*}$ and then patches the solutions. The gain is that convergence is
  immediate from the completeness of $\mathcal{A}_{\tau,p}$, with no need to
  track the series $\sum_{n}\|u_{n+1}-u_{n}\|$.
  The same scheme will be reused, essentially unchanged, in the geometric settings
  of Chapter~\ref{chap:geom-pam}, where only the kernel $p$ and the
  constant~$\Phi$ change.
\end{remark}

\subsection{The nonlinear case: multiplicative Lipschitz noise}
\label{sec:nonlinear}

The contraction scheme above never used the linearity of the map
$r\mapsto\beta r$ multiplying the noise. It only used that this map is globally
Lipschitz, and that $\mathcal{I}_{S}$ is a bounded operator on
$\mathcal{A}_{\tau,p}$ (Theorem~\ref{thm:stability}). Replacing $\beta\,u$ by a
nonlinear coefficient $\sigma(u)$ therefore costs nothing beyond bookkeeping,
as we now record.

\begin{assumption}\label{hyp:sigma-lip}
  The function $\sigma\colon\R\to\R$ is globally Lipschitz: there is a constant
  $L_{\sigma}<\infty$ such that
  \begin{equation}\label{eq:sigma-lip}
    |\sigma(r)-\sigma(r')|\,\le\,L_{\sigma}\,|r-r'|,
    \qquad r,r'\in\R.
  \end{equation}
  In particular $\sigma$ has linear growth,
  $|\sigma(r)|\le|\sigma(0)|+L_{\sigma}|r|$.
\end{assumption}

\begin{corollary}[Global well-posedness for multiplicative Lipschitz noise]
\label{cor:nonlinear-wp}
  Assume Dalang's condition (Assumption~\ref{hyp:spectral}), let $u_{0}$ be
  bounded measurable, and let $\sigma$ satisfy Assumption~\ref{hyp:sigma-lip}.
  Then for every $T>0$ and every $p\ge 2$ the nonlinear stochastic heat equation
  \begin{equation}\label{eq:mild-she-nonlinear}
    u(t,x)\,=\,P_{t}u_{0}(x)\,+\,\int_{0}^{t}\!\int_{\R^{d}}
      p_{t-s}(x-y)\,\sigma\bigl(u(s,y)\bigr)\,W(ds,dy)
  \end{equation}
  has a unique mild solution $u\in\mathcal{A}_{T,p}$ (Definition~\ref{def:A-t-p}).
  The linear PAM~\eqref{eq:mild-pam} is the special case $\sigma(r)=\beta r$, for
  which $L_{\sigma}=|\beta|$.
\end{corollary}

\begin{proof}
  For a process $Z$ write $Z^{\sigma}$ for the process
  $(s,y)\mapsto\sigma(Z(s,y))$, and define the solution map
  \begin{equation}\label{eq:Gamma-sigma-def}
    \Gamma_{\tau}^{\sigma}Z(t,x)
    \,:=\,P_{t}u_{0}(x)\,+\,\mathcal{I}_{S}Z^{\sigma}(t,x),
    \qquad (t,x)\in[0,\tau]\times\R^{d},
  \end{equation}
  with $\mathcal{I}_{S}$ the stochastic convolution of~\eqref{eq:conv-operator}.
  The proofs of Proposition~\ref{prop:Gamma-contraction},
  Theorem~\ref{thm:local-wp} and Theorem~\ref{thm:global-wp} carry over verbatim
  once the following three points, the only places where $\sigma$ enters, are
  checked.

  \begin{enumerate}
  \item {\it Stability.} Since $\sigma$ is continuous, $Z^{\sigma}$ is jointly
  measurable and adapted whenever $Z$ is. The linear growth in
  Assumption~\ref{hyp:sigma-lip} gives, pointwise in $\omega$,
  $|Z^{\sigma}(s,y)|\le|\sigma(0)|+L_{\sigma}|Z(s,y)|$, so taking the
  $L^{p}(\Omega)$ norm and the supremum over $(s,y)\in[0,\tau]\times\R^{d}$,
  \begin{equation*}
    \|Z^{\sigma}\|_{\mathcal{A}_{\tau,p}}
    \,\le\,|\sigma(0)|\,+\,L_{\sigma}\,\|Z\|_{\mathcal{A}_{\tau,p}}
    \,<\,\infty,
    \qquad Z\in\mathcal{A}_{\tau,p}.
  \end{equation*}
  Hence $Z^{\sigma}\in\mathcal{A}_{\tau,p}$, and Theorem~\ref{thm:stability} gives
  $\mathcal{I}_{S}Z^{\sigma}\in\mathcal{A}_{\tau,p}$; with the deterministic free
  term, $\Gamma_{\tau}^{\sigma}Z\in\mathcal{A}_{\tau,p}$. (In the linear case
  $Z^{\sigma}=\beta Z$ and this step is immediate.)

  \item {\it Contraction.} As in the proof of
  Proposition~\ref{prop:Gamma-contraction}, the free term cancels in the
  difference and $\mathcal{I}_{S}$ is linear, so
  $\Gamma_{\tau}^{\sigma}Z_{1}-\Gamma_{\tau}^{\sigma}Z_{2}
  =\mathcal{I}_{S}\bigl(Z_{1}^{\sigma}-Z_{2}^{\sigma}\bigr)$. The operator
  bound~\eqref{eq:I-S-bound} and the identity~\eqref{eq:S-star-Phi} give
  \begin{equation*}
    \bigl\|\Gamma_{\tau}^{\sigma}Z_{1}-\Gamma_{\tau}^{\sigma}Z_{2}
      \bigr\|_{\mathcal{A}_{\tau,p}}^{2}
    \,\le\,4\,p\,\Phi(\tau)\,
      \mathcal{N}_{\tau,p}\bigl(Z_{1}^{\sigma}-Z_{2}^{\sigma}\bigr).
  \end{equation*}
  The new ingredient is the pointwise Lipschitz bound~\eqref{eq:sigma-lip},
  applied $\omega$ by $\omega$, which after taking $L^{p}(\Omega)$ norms reads
  $\|Z_{1}^{\sigma}(s,y)-Z_{2}^{\sigma}(s,y)\|_{p}
  \le L_{\sigma}\,\|Z_{1}(s,y)-Z_{2}(s,y)\|_{p}$, whence
  $\mathcal{N}_{\tau,p}(Z_{1}^{\sigma}-Z_{2}^{\sigma})
  \le L_{\sigma}^{2}\,\mathcal{N}_{\tau,p}(Z_{1}-Z_{2})$. Therefore
  \begin{equation}\label{eq:Gamma-sigma-lip}
    \bigl\|\Gamma_{\tau}^{\sigma}Z_{1}-\Gamma_{\tau}^{\sigma}Z_{2}
      \bigr\|_{\mathcal{A}_{\tau,p}}
    \,\le\,2\,L_{\sigma}\,\sqrt{p\,\Phi(\tau)}\;
      \bigl\|Z_{1}-Z_{2}\bigr\|_{\mathcal{A}_{\tau,p}},
  \end{equation}
  which is exactly~\eqref{eq:Gamma-lip} with $|\beta|$ replaced by $L_{\sigma}$.
  Thus $\Gamma_{\tau}^{\sigma}$ is a contraction as soon as
  $4\,p\,L_{\sigma}^{2}\,\Phi(\tau)<1$, the analogue
  of~\eqref{eq:tau-smallness}; since $\Phi(\tau)\to 0$
  (Lemma~\ref{lem:Phi}), this fixes a horizon $\tau_{*}>0$ depending only on
  $p$, $L_{\sigma}$ and $\mu$.

  \item {\it Patching.} Since $u^{\sigma}$ is adapted whenever $u$ is, the
  restarted map and the Chapman--Kolmogorov
  recombination~\eqref{eq:restart-datum}--\eqref{eq:chapman-kolmogorov} hold
  verbatim, with $u(s,y)$ replaced by $u^{\sigma}(s,y)$ in the stochastic
  integrals. The $N$-step patch then produces the global solution
  $u\in\mathcal{A}_{T,p}$, and uniqueness follows window by window exactly as in
  Theorem~\ref{thm:global-wp}.
  \end{enumerate}

  In short, the only change with respect to the linear PAM is that the role of
  $|\beta|$ is now played by the Lipschitz constant $L_{\sigma}$. The noise
  still enters through the scalar $\Phi(\tau)$ only.
\end{proof}

\section*{Bibliographic notes}
\addcontentsline{toc}{section}{Bibliographic notes}
\markright{BIBLIOGRAPHIC NOTES}

\smallskip

\noindent
{\it Integration theory.} The worthy martingale measure of
Notation~\ref{not:mg-measure} and the integral of
Definition~\ref{def:walsh-integral} are Walsh's, and our presentation of them
is a condensed version of Chapter~2 of his Saint-Flour
notes~\cite{walsh:86:introduction}. Walsh works with space-time white noise,
for which the integrands are genuine functions. The extension to the
distributional integrands forced on us by a spatially homogeneous noise, that
is to say Theorem~\ref{thm:dalang-integral} and the modified martingale
measure $M^{Z}$ of Subsection~\ref{sec:intg-modif}, is due to
Dalang~\cite{dalang:99:extending}. The $L^{p}$ theory of
Subsection~\ref{sec:bdg-Lp} follows Khoshnevisan's
lectures~\cite{khoshnevisan:14:analysis}. A parallel and equally complete
theory is available in the Da Prato and Zabczyk framework of stochastic
evolution equations in Hilbert spaces~\cite{daprato:14:stochastic}, where the
noise is a Hilbert-space-valued Wiener process rather than a martingale
measure. The two theories are compared, and shown to give the same solutions,
by Dalang and Quer-Sardanyons~\cite{dalang:11:stochastic}. For the
one-dimensional equations driven by space-time white noise, where the whole
apparatus of Section~\ref{sec:integration} simplifies considerably, the recent
book of Dalang and Sanz-Sol\'{e}~\cite{dalang:26:stochastic} is the natural
reference. See also the minicourse~\cite{dalang:09:minicourse}.

\smallskip

\noindent
{\it The additive case.} Section~\ref{sec:additive} is standard material, and
is treated at length both in Da Prato and
Zabczyk~\cite{daprato:14:stochastic} and in Dalang and
Sanz-Sol\'{e}~\cite{dalang:26:stochastic}. It is worth stressing that the
additive equation is not merely a warm-up for the multiplicative one. Several
questions are genuinely easier, or only accessible, in that case, and we list
four of them.

\begin{enumerate}
  \item \emph{Ergodic theory.} The invariant measure is Gaussian and can be
    identified explicitly, which makes the long-time behavior of the equation
    comparatively simple to analyze (see Da Prato and
    Zabczyk~\cite{daprato:96:ergodicity}, and~\cite{chen:26:ergodic} for the
    difficulties met in the multiplicative case).
  \item \emph{Regularization by noise.} The equation can be solved with a drift
    that would not be valid for a deterministic PDE. When $\sigma$ is a nonzero
    constant, existence, uniqueness and a comparison theorem hold for the
    one-dimensional equation with a drift $f$ that is merely measurable, as
    shown by Gy\"{o}ngy and Pardoux~\cite{gyongy:93:on}, and in fact for an $f$
    that need not even be locally bounded~\cite{gyongy:93:regularization}.
    Bally, Gy\"{o}ngy and Pardoux~\cite{bally:94:white} extend this to a
    nondegenerate, non-constant $\sigma$, the argument resting on a
    Malliavin-calculus estimate for the density of the driftless equation. The
    noise is thus restoring a well-posedness that the deterministic equation
    does not have, and this is by now a subject of its own. Butkovsky and
    Mytnik~\cite{butkovsky:19:regularization} construct the flow of solutions,
    and prove path-by-path uniqueness, for a bounded measurable drift. Athreya,
    Butkovsky, L\^{e} and Mytnik~\cite{athreya:24:well-posedness} reach drifts
    which are distributions.
  \item \emph{Numerical schemes.} They converge at better rates in the additive
    case. For a scheme built from the increments of the Wiener process the
    temporal rate cannot exceed $\tfrac14$, a barrier going back to the lower
    bounds of Davie and Gaines~\cite{davie:01:convergence}. The additive
    structure allows one instead to sample the stochastic convolution itself,
    which is the idea behind the accelerated exponential Euler scheme of
    Jentzen and Kloeden~\cite{jentzen:09:overcoming}. That this genuinely
    breaks the barrier, and gives a temporal rate arbitrarily close to $1$ for
    a large class of nonlinearities, was proved by Djurdjevac, Gerencs\'{e}r
    and Kremp~\cite{djurdjevac:24:higher}.
  \item \emph{Singular equations.} Those which currently lie at the frontier of
    the field are, with few exceptions, equations with additive noise and a
    nonlinear drift. This is the case of the KPZ
    equation~\cite{hairer:13:solving,gubinelli:17:kpz} and of the dynamical
    $\Phi^{4}_{3}$
    model~\cite{hairer:14:theory,catellier:18:paracontrolled,mourrat:17:dynamic}.
\end{enumerate}

\smallskip

\noindent
{\it PAM and the nonlinear case.} The existence and uniqueness results of
Section~\ref{sec:contraction}, namely Theorem~\ref{thm:global-wp} and its
nonlinear counterpart Corollary~\ref{cor:nonlinear-wp}, go back to the same
sources as the integration theory itself, that is
to~\cite{walsh:86:introduction},
to~\cite[Theorem~13]{dalang:99:extending}, to~\cite{daprato:14:stochastic} and
to~\cite{dalang:26:stochastic}. Our own route through a single application of
the Banach fixed point theorem, rather than a Picard iteration, is a matter of
presentation only.

Beyond the globally Lipschitz coefficient $\sigma$ of
Assumption~\ref{hyp:sigma-lip}, the picture becomes considerably more
demanding, and we only indicate a few directions. Superlinear or non-Lipschitz
reaction terms, where solutions may blow up in finite time, are studied by
Salins~\cite{salins:22:global} and by Salins and
Tindel~\cite{salins:24:regularity}. Rough initial conditions, where $u_{0}$ is
merely a measure, are handled by Chen and Dalang~\cite{chen:15:moments}, and
are the natural setting for the delta initial datum of the polymer
interpretation. Noises with jumps are covered by the random field approach of
Balan~\cite{balan:14:spdes,balan:15:integration} and by the evolution equation
approach of Peszat and Zabczyk~\cite{peszat:07:stochastic}. Finally, when
Dalang's condition fails the product $u\,\dot{W}$ has no meaning and the
equation has to be renormalized. The two systematic theories available for this
are the regularity structures of Hairer~\cite{hairer:14:theory} and the
paracontrolled calculus of Gubinelli, Imkeller and
Perkowski~\cite{gubinelli:15:paracontrolled}. The two-dimensional PAM is solved
in that framework in~\cite{hairer:15:simple}, and moment estimates for
renormalized models are obtained in~\cite{chen:21:moment}.

\smallskip

\noindent
{\it The critical dimension two.} Dalang's condition
(Assumption~\ref{hyp:spectral}) is what keeps the whole theory alive, and for a
space-time white noise it holds only in dimension $d=1$. Dimension two is
exactly critical, and the failure there is marginal, which is what makes it
interesting. Mollify the noise at spatial scale $\varepsilon$. The mollified
equation is perfectly well posed, but its solution $u^{\varepsilon}$ degenerates
as $\varepsilon\to 0$ unless the coupling constant is tuned to zero at the rate
\begin{equation*}
  \beta_{\varepsilon}\,\approx\,\hat{\beta}\,
  \bigl(\log(1/\varepsilon)\bigr)^{-1/2}.
\end{equation*}
This renormalization was identified by Bertini and
Cancrini~\cite{bertini:98:two-dimensional}. For $\hat{\beta}<1$ the fluctuations
of $u^{\varepsilon}$ are Gaussian in the limit, of Edwards--Wilkinson type, as
shown by Caravenna, Sun and Zygouras~\cite{caravenna:17:universality}. The
critical window sits at $\hat{\beta}=1$, where a second-order tuning of
$\beta_{\varepsilon}$ is required, and the limit there is a genuinely new
object. It is neither a function nor the exponential of a Gaussian field. It is rather described by a
process of random measures on $\R^{2}\times\R^{2}$ indexed by pairs of times called
the \emph{critical 2d stochastic heat flow}. Its construction is the content
of~\cite{caravenna:23:critical-a}. The moments in the critical window had been
computed beforehand by Gu, Quastel and Tsai~\cite{gu:21:moments} and
in~\cite{caravenna:19:on}. In addition, sharper lower bounds on those moments show that the
flow not to be a Gaussian multiplicative
chaos~\cite{caravenna:23:critical-b}, that is to say not the exponential of a
log-correlated Gaussian field. Tsai~\cite{tsai:24:stochastic} gives an
independent construction, characterizing the flow by a short list of axioms of
which the essential one is the matching of the first four moments. The
survey~\cite{caravenna:24:critical} is a good entry point to this circle of
ideas.

\smallskip

\noindent
{\it The Anderson Hamiltonian.} A second object attached to a rough noise, and
this one not tied to $d=2$, is the operator $\mathcal{H}=\tfrac12\Delta+\dot{W}$
itself, the \emph{Anderson Hamiltonian}. When the noise is rough enough,
$\mathcal{H}$ is not an operator in the usual sense and has to be renormalized
before it can be diagonalized. Labb\'{e}~\cite{labbe:19:continuous} constructs
it on a box as a self-adjoint operator with pure point spectrum, in dimension
$d\le 3$. Gubinelli, Ugurcan and Zachhuber~\cite{gubinelli:20:semilinear} study
the semigroup it generates and the semilinear equations that semigroup drives,
in dimensions two and three. In two dimensions the asymptotics of the
eigenvalues of $\mathcal{H}$, which govern the moment asymptotics in the same
way as the Lyapunov exponents~\eqref{eq:lyapunov} do here, are obtained by
Chouk and van Zuijlen~\cite{chouk:21:asymptotics}.

\chapter{Moment estimates}
\label{chap:moments}

\section{Chaos expansions}\label{sec:chaos}

Chapter~\ref{chap:existence} established existence and uniqueness of the
mild solution~\eqref{eq:mild-pam} via a Banach fixed-point argument in the
space $\mathcal{A}_{T,p}$. In the present section we revisit existence and
uniqueness with a different tool, namely the \emph{Wiener chaos expansion} of
Section~\ref{ssec:isonormal}. Instead of iterating the solution map, one looks
for an explicit term-by-term formula for $u(t,x)$, as a series of multiple
stochastic integrals. This approach yields an explicit $L^{2}$ moment
formula,
\begin{equation}\label{eq:L2-moment-chaos}
  \E\bigl[u(t,x)^{2}\bigr]
  = \sum_{n=0}^{\infty} n!\,\|f_{n}(\cdot,t,x)\|^{2}_{\mathcal{H}^{\otimes n}},
\end{equation}
where $f_{n}$ are the chaos kernels derived below. This representation is at
the heart of our moment estimates. Moreover, the convergence
of~\eqref{eq:L2-moment-chaos} under Dalang's
condition~\eqref{eq:dalang-cond} provides an independent proof of the
$L^{2}$-finiteness of $u(t,x)$.

\subsection{Derivation of the chaos kernels}\label{ssec:chaos-kernels}

Since the solution $u(s,y)$ is adapted to the
filtration~\eqref{eq:filtration}, Proposition~\ref{prop:ito-sko} shows that
the It\^{o}--Walsh integral in~\eqref{eq:mild-pam} coincides with the
Skorohod integral $\delta$ (Definition~\ref{def:skorohod}): for any adapted
$v\in\mathcal{H}_{\omega}$, one can write
\begin{equation}\label{eq:ito-sko}
  \int_{0}^{t}\!\int_{\R^{d}} v(s,y)\,W(ds,dy)
  \;=\;\delta\!\left(v\,\mathbf{1}_{[0,t]}\right).
\end{equation}
We therefore work with the Skorohod-sense mild equation
\begin{equation}\label{eq:mild-sko}
  u(t,x)
  \;=\;P_{t}u_{0}(x)
  \;+\;\beta\int_{0}^{t}\!\int_{\R^{d}}
    p_{t-s}(x-y)\,u(s,y)\;\delta W_{s,y},
\end{equation}
which is just another notation for~\eqref{eq:mild-pam} in case of an adapted solution.

By Theorem~\ref{thm:chaos-decomp}, for each fixed $(t,x)$ the random variable
$u(t,x)$ has a Wiener chaos decomposition
\begin{equation}\label{eq:u-chaos}
  u(t,x) = \sum_{n=0}^{\infty} I_{n}\!\left(f_{n}(\cdot,t,x)\right),
\end{equation}
where $f_{n}(\cdot,t,x)\in\mathcal{H}^{\odot n}$ is a symmetric
kernel depending on $(t,x)$ as parameters, and $f_{0}(t,x)=\E[u(t,x)]$.
The kernels $f_{n}$ are not specified by~\eqref{eq:u-chaos} yet. They will be
identified by substituting the expansion back into the mild
equation~\eqref{eq:mild-sko}.

\begin{lemma}[Chaos kernels]\label{lem:chaos-kernels}
  The kernels $f_{n}(\,\cdot\,;\,t,x)\in\mathcal{H}^{\odot n}$
  in~\eqref{eq:u-chaos} are given by the initial value
  \begin{equation}\label{eq:f0-val}
    f_{0}(t,x) = P_{t}u_{0}(x),
  \end{equation}
  and, for $n\ge 1$, by the explicit symmetric kernel
  \begin{multline}\label{eq:fn-def}
    f_{n}(s_{1},y_{1},\ldots,s_{n},y_{n};\,t,x) 
    = \frac{\beta^{n}}{n!}\,
      p_{t-s_{\sigma(n)}}(x-y_{\sigma(n)})\,
      p_{s_{\sigma(n)}-s_{\sigma(n-1)}}(y_{\sigma(n)}-y_{\sigma(n-1)}) \\
      \times\cdots\times
      p_{s_{\sigma(2)}-s_{\sigma(1)}}(y_{\sigma(2)}-y_{\sigma(1)})\,
      P_{s_{\sigma(1)}}u_{0}(y_{\sigma(1)}),
  \end{multline}
  where $\sigma$ is the unique permutation of $\{1,\ldots,n\}$ such that
  \begin{equation}\label{eq:sigma-order}
    0 < s_{\sigma(1)} < s_{\sigma(2)} < \cdots < s_{\sigma(n)} < t.
  \end{equation}
\end{lemma}

\begin{proof}
The strategy is to substitute the chaos expansion~\eqref{eq:u-chaos} into the
mild equation~\eqref{eq:mild-sko} and match chaos levels using
identity~\eqref{eq:delta-chaos}. This yields a recursion on the kernels
(Step~1), which is then iterated explicitly for $n=1$ and $n=2$
(Steps~2--3) before being resolved in general by induction (Step~4).

\smallskip

\noindent
{\it Step 1: Derivation of the recursion.}
Substitute the chaos expansion~\eqref{eq:u-chaos} of $u(s,y)$ into the
right-hand side of~\eqref{eq:mild-sko}.
The integrand $v(s,y):=p_{t-s}(x-y)\,u(s,y)$ then admits the chaos decomposition
\begin{equation}\label{eq:v-chaos}
  v(s,y) = \sum_{n=0}^{\infty} I_{n}\!\left(g_{n}(\,\cdot\,;\,s,y)\right),
 \end{equation}
where the kernel $g_{n}$ is defined by 
\begin{equation*}
g_{n}(s_{1},y_{1},\ldots,s_{n},y_{n};\,s,y)
  = p_{t-s}(x-y)\,f_{n}(s_{1},y_{1},\ldots,s_{n},y_{n};\,s,y),
\end{equation*}
and where $(s,y)$ plays the role of the free $\mathcal{H}$-variable $\star$
in~\eqref{eq:delta-chaos}.
Applying~\eqref{eq:delta-chaos} to $v$ we end up with
\begin{equation}\label{eq:delta-v-chaos}
  \delta(v) = \sum_{n=0}^{\infty} I_{n+1}\!\left(\widetilde{g}_{n}\right),
\end{equation}
where, writing $(s_{n+1},y_{n+1})$ for the free variable $(s,y)$,
$\widetilde{g}_{n}$ is the full symmetrization of
$g_{n}(s_{1},y_{1},\ldots,s_{n},y_{n};\,s_{n+1},y_{n+1})$
in all $n+1$ pairs $(s_{k},y_{k})$:
\begin{multline}\label{eq:gn-sym}
  \widetilde{g}_{n}(s_{1},y_{1},\ldots,s_{n+1},y_{n+1}) \\
  = \frac{1}{n+1}\sum_{k=1}^{n+1}
    p_{t-s_{k}}(x-y_{k})\,
    f_{n}\!\left(s_{1},y_{1},\ldots,\widehat{(s_{k},y_{k})},\ldots,s_{n+1},y_{n+1};\,s_{k},y_{k}\right).
\end{multline}
The right-hand side of~\eqref{eq:mild-sko} therefore has chaos decomposition
\begin{equation}\label{eq:rhs-chaos}
  P_{t}u_{0}(x) + \beta\,\delta(v)
  = I_{0}\!\left(P_{t}u_{0}(x)\right)
    + \sum_{n=1}^{\infty} I_{n}\!\left(\beta\,\widetilde{g}_{n-1}\right).
\end{equation}
By uniqueness of the chaos decomposition in Theorem~\ref{thm:chaos-decomp},
matching level $n=0$ gives~\eqref{eq:f0-val}, and matching level $n\ge 1$ gives
$f_{n} = \beta\,\widetilde{g}_{n-1}$.
Substituting~\eqref{eq:gn-sym} with $n$ replaced by $n-1$ yields the recursion,
for $n\ge 1$:
\begin{multline}\label{eq:fn-recursion}
  f_{n}(s_{1},y_{1},\ldots,s_{n},y_{n};\,t,x) \\
  = \frac{\beta}{n}\sum_{k=1}^{n}
    p_{t-s_{k}}(x-y_{k})\,
    f_{n-1}\!\left(s_{1},y_{1},\ldots,\widehat{(s_{k},y_{k})},\ldots,s_{n},y_{n};\,s_{k},y_{k}\right),
\end{multline}
where $\widehat{\ }$ indicates the omitted pair. For sake of clarity, we will now detail the cases $n=1$ and $n=2$.

\smallskip

\noindent
{\it Step 2: Case $n=1$.}
In~\eqref{eq:fn-recursion} with $n=1$, the sum reduces to a single term
(there is one variable, the hat omits nothing, and no averaging is needed):
\begin{equation}\label{eq:f1-chaos}
  f_{1}(s_{1},y_{1};\,t,x)
  = \beta\,p_{t-s_{1}}(x-y_{1})\,f_{0}(s_{1},y_{1})
  = \beta\,p_{t-s_{1}}(x-y_{1})\,P_{s_{1}}u_{0}(y_{1}).
\end{equation}
This is~\eqref{eq:fn-def} for $n=1$ (with the trivial permutation $\sigma=\mathrm{id}$).

\smallskip

\noindent
{\it Step 3: Case $n=2$.}
Applying~\eqref{eq:fn-recursion} with $n=2$ and substituting~\eqref{eq:f1-chaos} we obtain
\begin{align}\label{eq:f2-chaos}
  &f_{2}(s_{1},y_{1},s_{2},y_{2};\,t,x) \notag\\
  &= \frac{\beta}{2}\Bigl[
       p_{t-s_{1}}(x-y_{1})\,f_{1}(s_{2},y_{2};\,s_{1},y_{1})
     + p_{t-s_{2}}(x-y_{2})\,f_{1}(s_{1},y_{1};\,s_{2},y_{2})
     \Bigr] \\
  &= \frac{\beta^{2}}{2}\Bigl[
       p_{t-s_{1}}(x-y_{1})\,p_{s_{1}-s_{2}}(y_{1}-y_{2})\,P_{s_{2}}u_{0}(y_{2})
     + p_{t-s_{2}}(x-y_{2})\,p_{s_{2}-s_{1}}(y_{2}-y_{1})\,P_{s_{1}}u_{0}(y_{1})
     \Bigr]. \notag
\end{align}
The two summands correspond to the two orderings of $(s_{1},y_{1})$ and $(s_{2},y_{2})$:
the first is nonzero only when $s_{2}<s_{1}<t$ (permutation $\sigma(1)=2,\,\sigma(2)=1$),
and the second only when $s_{1}<s_{2}<t$ (identity permutation).
For Lebesgue-almost every $(s_{1},y_{1},s_{2},y_{2})$, exactly one summand is nonzero,
and the factor $\tfrac{1}{2}=\tfrac{1}{2!}$ is the reciprocal of the number of orderings.
This is~\eqref{eq:fn-def} for $n=2$.

\smallskip

\noindent
{\it Step 4: General formula by induction.}
An induction on $n$, using~\eqref{eq:fn-recursion} and the inductive hypothesis for $f_{n-1}$,
shows that the double sum over $k\in\{1,\ldots,n\}$ and permutations $\tau\in\mathfrak{S}_{n-1}$
(ordering the remaining $n-1$ variables) is in bijection with $\mathfrak{S}_{n}$,
yielding
\begin{multline}\label{eq:fn-sum}
  f_{n}(s_{1},y_{1},\ldots,s_{n},y_{n};\,t,x) \\
  = \frac{\beta^{n}}{n!}
    \sum_{\sigma\in\mathfrak{S}_{n}}
    p_{t-s_{\sigma(n)}}(x-y_{\sigma(n)})\,
    \prod_{j=1}^{n-1}p_{s_{\sigma(j+1)}-s_{\sigma(j)}}(y_{\sigma(j+1)}-y_{\sigma(j)})\,
    P_{s_{\sigma(1)}}u_{0}(y_{\sigma(1)}),
\end{multline}
where $\mathfrak{S}_{n}$ denotes the symmetric group on $\{1,\ldots,n\}$.
Since $p_{t'}(\cdot)=0$ for $t'\le 0$, the summand for $\sigma$ is nonzero only on
$\{s_{\sigma(1)}<\cdots<s_{\sigma(n)}<t\}$, and for Lebesgue-almost every
$(s_{1},\ldots,s_{n})$ at most one $\sigma$ satisfies this.
Hence~\eqref{eq:fn-sum} reduces to the single-term formula~\eqref{eq:fn-def}.
\end{proof}

\begin{remark}[Relation with directed polymer models]\label{rem:polymer}
Formula~\eqref{eq:fn-def} admits a transparent probabilistic reading. Namely
one places $n$ spacetime marks $(s_{i},y_{i})$ in the slab
$[0,t]\times\R^{d}$, orders them by time, and weights each consecutive pair by
the heat kernel. This is exactly the $n$-th term in the partition function of
a directed polymer in a random environment. The marks play the role of the
polymer's intermediate vertices, and the heat kernels are the transition
weights of the underlying random walk. As for the sum over chaos levels, it
corresponds to a sum over all possible numbers of contacts with the
environment.
We refer to~\cite{comets:17:directed} for a detailed account of directed polymers
and their connection with the PAM.
\end{remark}

\subsection{Convergence and the \texorpdfstring{$L^{2}$}{L2} moment formula}
\label{ssec:chaos-convergence}

Lemma~\ref{lem:chaos-kernels} identifies the kernels $f_{n}$. It says nothing
about the convergence of the series~\eqref{eq:u-chaos} in $L^{2}(\Omega)$.
This convergence requires a control on the growth of the norms
$n!\,\|f_{n}(\cdot,t,x)\|^{2}_{\mathcal{H}^{\otimes n}}$, and this is where
Dalang's condition~\eqref{eq:dalang-cond} enters. Once the convergence is
established, the $L^{2}$ moment formula~\eqref{eq:L2-moment-formula} follows
from the orthogonality of distinct chaos levels~\eqref{eq:In-isometry}.

\begin{theorem}[Existence-uniqueness and $L^{2}$ moments via chaos]
  \label{thm:chaos-eu}
  Assume Dalang's condition~\eqref{eq:dalang-cond} and let $u_{0}\in L^{\infty}(\R^{d})$.
  Then, for every $(t,x)\in(0,\infty)\times\R^{d}$,
  \begin{equation}\label{eq:chaos-L2-converge}
    \sum_{n=0}^{\infty}\,n!\,\|f_{n}(\cdot,t,x)\|^{2}_{\mathcal{H}^{\otimes n}}
    \;<\;\infty.
  \end{equation}
  Consequently, the series~\eqref{eq:u-chaos} converges in $L^{2}(\Omega)$,
  defining a unique mild solution to~\eqref{eq:mild-pam} with $L^{2}$-moment
  \begin{equation}\label{eq:L2-moment-formula}
    \E\bigl[u(t,x)^{2}\bigr]
    = \sum_{n=0}^{\infty}\,n!\,\|f_{n}(\cdot,t,x)\|^{2}_{\mathcal{H}^{\otimes n}}.
  \end{equation}
\end{theorem}

\begin{proof}
The $L^{2}$-moment formula~\eqref{eq:L2-moment-formula} follows from the
orthogonality of different chaos levels~\eqref{eq:In-isometry}, so it
suffices to establish~\eqref{eq:chaos-L2-converge}. Throughout the proof we
write $\mu^{\otimes n}(d\xi):=\prod_{i=1}^{n}\mu(d\xi_{i})$ for the
$n$-fold product spectral measure.

\smallskip

\noindent
{\it Step 1: Preliminary considerations.}

\smallskip

\noindent
Before passing to Fourier variables, we record a trivial bound on the kernel
$f_{n}$ itself in direct variables. The heat semigroup $(P_{s})_{s\ge 0}$
is a contraction on $L^{\infty}(\R^{d})$, since $p_{s}$~\eqref{eq:heat-kernel-explicit}
is a probability density for every $s>0$. Namely we have
\begin{equation}\label{eq:Ps-Linfty}
  |P_{s}u_{0}(y)| \;\le\; \|u_{0}\|_{\infty}, \qquad s\ge 0,\ y\in\R^{d}.
\end{equation}
Let us isolate the pure chain of heat-kernel transitions appearing
in~\eqref{eq:fn-def} by setting, on the ordered simplex~\eqref{eq:sigma-order}
(with $\sigma$ as in~\eqref{eq:sigma-order}),
\begin{equation}\label{eq:Kn-def}
  K_{n}(s_{1},y_{1},\ldots,s_{n},y_{n};\,t,x)
  :=
  p_{t-s_{\sigma(n)}}(x-y_{\sigma(n)})\,
  \prod_{j=1}^{n-1}p_{s_{\sigma(j+1)}-s_{\sigma(j)}}(y_{\sigma(j+1)}-y_{\sigma(j)}),
\end{equation}
and $K_{n}:=0$ off the simplex~\eqref{eq:sigma-order}. Substituting~\eqref{eq:Ps-Linfty}
into~\eqref{eq:fn-def}, and using that the heat kernel $p_{t}$ is
non-negative, gives the direct (pre-Fourier) bound
\begin{equation}\label{eq:fn-direct-bound}
  \bigl|f_{n}(s_{1},y_{1},\ldots,s_{n},y_{n};\,t,x)\bigr|
  \;\le\;
  \frac{\beta^{n}}{n!}\,\|u_{0}\|_{\infty}\,
  K_{n}(s_{1},y_{1},\ldots,s_{n},y_{n};\,t,x).
\end{equation}
This trivial bound isolates the role of the initial datum $u_{0}$ from the
chain of transitions $K_{n}$. This is exactly the bound on $u_{0}$ which we
shall pull out of the Fourier transform of $f_{n}$, in order to control
$\|f_{n}(\cdot,t,x)\|^{2}_{\mathcal{H}^{\otimes n}}$ below.

\smallskip

\noindent
As a second preliminary step, we reduce the integral defining
$\|K_{n}(\cdot,t,x)\|^{2}_{\mathcal{H}^{\otimes n}}$ (and hence, through the
bound~\eqref{eq:fn-direct-bound}, the one defining
$\|f_{n}(\cdot,t,x)\|^{2}_{\mathcal{H}^{\otimes n}}$) to a single ordered
region in time. Both $K_{n}$~\eqref{eq:Kn-def} and $f_{n}$ (which lies in
$\mathcal{H}^{\odot n}$ by Lemma~\ref{lem:chaos-kernels}) are symmetric
functions of the $n$ pairs $(s_{1},y_{1}),\ldots,(s_{n},y_{n})$. Indeed,
relabeling the pairs by a permutation $\pi\in\mathfrak{S}_{n}$ only relabels
which index plays the role of $\sigma(1),\ldots,\sigma(n)$
in~\eqref{eq:sigma-order}. It leaves the sorted chain of heat-kernel
transitions unchanged, and thus the values of $K_{n}$ and $f_{n}$.
Consequently, the function
\begin{equation*}
  \Psi_{n}(s_{1},\ldots,s_{n})
  :=
  \int_{\R^{nd}}\bigl|\mathcal{F}_{y}\,K_{n}(s,\cdot;\,t,x)(\xi)\bigr|^{2}\,
  \mu^{\otimes n}(d\xi),
\end{equation*}
which is the $s$-density of $\|K_{n}(\cdot,t,x)\|^{2}_{\mathcal{H}^{\otimes n}}$
under the inner product~\eqref{eq:H-innprod-fourier}, is itself symmetric in
$s=(s_{1},\ldots,s_{n})$. Indeed, permuting the $s_{i}$ amounts to permuting
the dummy variables $y_{i}$, equivalently the $\xi_{i}$, inside the integral.
This leaves its value unchanged, since $\mu^{\otimes n}(d\xi)$ is itself a
symmetric product measure. Now decompose the cube $[0,t]^{n}$ into the $n!$ disjoint images
$\pi\cdot T_{n}(t)$ (for a generic $\pi\in\mathfrak{S}_{n}$) of the ordered simplex
\begin{equation}\label{eq:simplex}
  T_{n}(t) = \bigl\{(s_{1},\ldots,s_{n})\in[0,\infty)^{n}:
  0 < s_{1} < \cdots < s_{n} < t\bigr\},
\end{equation}
so that symmetry of $\Psi_{n}$ gives
\begin{equation}\label{eq:Kn-simplex-reduction}
  \|K_{n}(\cdot,t,x)\|^{2}_{\mathcal{H}^{\otimes n}}
  \;=\;
  \int_{[0,t]^{n}}\Psi_{n}(s)\,ds
  \;=\;
  n!\int_{T_{n}(t)}\Psi_{n}(s)\,ds.
\end{equation}
We can perform the same manipulation for $f_{n}$: writing
\begin{equation*}
  \Phi_{n}(s_{1},\ldots,s_{n})
  :=
  \int_{\R^{nd}}\bigl|\mathcal{F}_{y}\,f_{n}(s,\cdot;\,t,x)(\xi)\bigr|^{2}\,
  \mu^{\otimes n}(d\xi)
\end{equation*}
for the analogous $s$-density of $\|f_{n}(\cdot,t,x)\|^{2}_{\mathcal{H}^{\otimes n}}$,
an identical symmetry argument gives
\begin{equation*}
  \|f_{n}(\cdot,t,x)\|^{2}_{\mathcal{H}^{\otimes n}}
  =
  n!\int_{T_{n}(t)}\Phi_{n}(s)\,ds.
\end{equation*}
We may therefore restrict our attention, in what follows, to the single ordered
region $T_{n}(t)$. On that set we have $\sigma=\mathrm{id}$ in~\eqref{eq:sigma-order} and
the chain~\eqref{eq:Kn-def} takes the simple, unsymmetrized form
\begin{equation}\label{eq:Kn-simplex-value}
  K_{n}(s_{1},y_{1},\ldots,s_{n},y_{n};\,t,x)
  =
  p_{t-s_{n}}(x-y_{n})\,
  \prod_{i=1}^{n-1}p_{s_{i+1}-s_{i}}(y_{i+1}-y_{i}),
  \qquad (s_{1},\ldots,s_{n})\in T_{n}(t).
\end{equation}

\smallskip

\noindent
Let us now record, in terms of $\Psi_{n}$, the bound towards~\eqref{eq:chaos-L2-converge}
that the two preliminary considerations above already yield. Since
$\Lambda\ge 0$ pointwise (by definition of the covariance~\eqref{eq:noise-cov}),
the $\mathcal{H}^{\otimes n}$ inner product is monotone under pointwise
domination of non-negative kernels. Indeed, let $|h|\le g$ pointwise with
$g\ge 0$. Then $h(y)h(y')\le|h(y)||h(y')|\le g(y)g(y')$ for all $y,y'$.
Multiplying by the non-negative weight $\prod_{i}\Lambda(y_{i}-y_{i}')$ and
integrating preserves the inequality, so that
$\|h\|^{2}_{\mathcal{H}^{\otimes n}}\le\|g\|^{2}_{\mathcal{H}^{\otimes n}}$.
Applying this for each fixed $s$ to the direct
bound~\eqref{eq:fn-direct-bound}, we get
\begin{equation*}
  \Phi_{n}(s) \;\le\; \frac{\beta^{2n}}{(n!)^{2}}\,\|u_{0}\|_{\infty}^{2}\,\Psi_{n}(s).
\end{equation*}
Integrating over $T_{n}(t)$, using $\|f_{n}(\cdot,t,x)\|^{2}_{\mathcal{H}^{\otimes n}}=n!\int_{T_{n}(t)}\Phi_{n}(s)\,ds$
above, and multiplying through by $n!$, we arrive at
\begin{equation}\label{eq:fn-Psi-bound}
  n!\,\|f_{n}(\cdot,t,x)\|^{2}_{\mathcal{H}^{\otimes n}}
  \;\le\;
  \beta^{2n}\,\|u_{0}\|_{\infty}^{2}\int_{T_{n}(t)}\Psi_{n}(s)\,ds ,
\end{equation}
where $\Psi_{n}(s)$ will be evaluated thanks to~\eqref{eq:Kn-simplex-value}.
This is exactly the $n$-th term on the left-hand side
of~\eqref{eq:chaos-L2-converge}. To summarize, the two preliminary steps have
reduced the proof of~\eqref{eq:chaos-L2-converge} to an estimate, for each
$n$, of the single integral $\int_{T_{n}(t)}\Psi_{n}(s)\,ds$. This integral is
built from the explicit and $u_{0}$-free heat-kernel chain $K_{n}$ on the
ordered simplex~\eqref{eq:Kn-simplex-value}. This is the task of Step~2.

\smallskip

\noindent
{\it Step 2: Fourier bound on each term.}

\smallskip

\noindent
By~\eqref{eq:fn-Psi-bound}, it remains to compute $\Psi_{n}(s)$ for
$s\in T_{n}(t)$, where, by definition,
\begin{equation*}
  \Psi_{n}(s)
  =
  \int_{\R^{nd}}\bigl|\mathcal{F}_{y}\,K_{n}(s,\cdot;\,t,x)(\xi)\bigr|^{2}\,
  \mu^{\otimes n}(d\xi).
\end{equation*}
We compute the spatial Fourier transform of $K_{n}$ explicitly, from its
simplified form~\eqref{eq:Kn-simplex-value} on $T_{n}(t)$. That is, for $s\in T_{n}(t)$ and
$\xi=(\xi_{1},\ldots,\xi_{n})\in\R^{nd}$, by definition of the Fourier transform we have
\begin{equation}\label{eq:Kn-FT-plain}
  \mathcal{F}_{y}\,K_{n}(s,\cdot;\,t,x)(\xi)
  =
  \int_{\R^{nd}}
  p_{t-s_{n}}(x-y_{n})\,
  \prod_{i=1}^{n-1}p_{s_{i+1}-s_{i}}(y_{i+1}-y_{i})\,
  e^{-\imath\sum_{j=1}^{n}\xi_{j}\cdot y_{j}}\,
  dy_{1}\cdots dy_{n}.
\end{equation}
To evaluate this integral, set $u_{i}:=y_{i+1}-y_{i}$ for $1\le i\le n-1$ and
$u_{n}:=x-y_{n}$, so that
$K_{n}(s,y;\,t,x)=\prod_{i=1}^{n}p_{s_{i+1}-s_{i}}(u_{i})$ (where we set
$s_{n+1}:=t$). Note that this change of variables has unit Jacobian. Furthermore, since $y_{i}=x-\sum_{k=i}^{n}u_{k}$, we have
\begin{equation*}
  \sum_{i=1}^{n}\xi_{i}\cdot y_{i}
  =
  x\cdot\Bigl(\sum_{i=1}^{n}\xi_{i}\Bigr)
  -
  \sum_{k=1}^{n}u_{k}\cdot(\xi_{1}+\cdots+\xi_{k}).
\end{equation*}
Substituting into~\eqref{eq:Kn-FT-plain} (the Jacobian is $1$, and
$K_{n}(s,u;\,t,x)=\prod_{i=1}^{n}p_{s_{i+1}-s_{i}}(u_{i})$), we obtain
\begin{equation}\label{eq:Kn-FT-substituted}
  \mathcal{F}_{y}\,K_{n}(s,\cdot;\,t,x)(\xi)
  =
  e^{-\imath x\cdot\sum_{i=1}^{n}\xi_{i}}
  \int_{\R^{nd}}
  \prod_{i=1}^{n}
  p_{s_{i+1}-s_{i}}(u_{i})\,
  e^{\imath\,u_{i}\cdot(\xi_{1}+\cdots+\xi_{i})}\,
  du_{1}\cdots du_{n}.
\end{equation}
Since the integrand in~\eqref{eq:Kn-FT-substituted} factorizes over $i$,
we thus get
\begin{equation*}
  \mathcal{F}_{y}\,K_{n}(s,\cdot;\,t,x)(\xi)
  =
  e^{-\imath x\cdot\sum_{i=1}^{n}\xi_{i}}
  \prod_{i=1}^{n}
  \Bigl[\int_{\R^{d}}
  p_{s_{i+1}-s_{i}}(u_{i})\,
  e^{\imath\,u_{i}\cdot(\xi_{1}+\cdots+\xi_{i})}\,
  du_{i}\Bigr].
\end{equation*}
Using the Fourier representation of the heat
kernel~\eqref{eq:heat-kernel-explicit}, that is
$\mathcal{F}p_{w}(\zeta)=e^{-w|\zeta|^{2}/2}$, we end up with
\begin{equation}\label{eq:Kn-FT}
  \mathcal{F}_{y}\,K_{n}(s,\cdot;\,t,x)(\xi)
  =
  e^{-\imath x\cdot\sum_{i}\xi_{i}}
  \prod_{i=1}^{n}e^{-\frac12(s_{i+1}-s_{i})|\xi_{1}+\cdots+\xi_{i}|^{2}},
  \qquad s\in T_{n}(t).
\end{equation}
Taking the modulus square (the phase $e^{-\imath x\cdot\sum_{i}\xi_{i}}$ has
modulus $1$) and integrating against $\mu^{\otimes n}(d\xi)$ gives the
\emph{exact} value of $\Psi_{n}$ on $T_{n}(t)$:
\begin{equation}\label{eq:Psi-n-value}
  \Psi_{n}(s)
  =
  \int_{\R^{nd}}
  \prod_{i=1}^{n}e^{-(s_{i+1}-s_{i})|\xi_{1}+\cdots+\xi_{i}|^{2}}\,
  \mu^{\otimes n}(d\xi),
  \qquad s\in T_{n}(t).
\end{equation}
Combining~\eqref{eq:Psi-n-value} with the bound~\eqref{eq:fn-Psi-bound} from
Step~1, we conclude this step with the upper bound
\begin{equation}\label{eq:fn-simplex-bound}
  n!\,\|f_{n}(\cdot,t,x)\|^{2}_{\mathcal{H}^{\otimes n}}
  \;\le\;
  \|u_{0}\|_{\infty}^{2}\,\beta^{2n}
  \int_{T_{n}(t)}\!\int_{\R^{nd}}
  \prod_{i=1}^{n}e^{-(s_{i+1}-s_{i})|\xi_{1}+\cdots+\xi_{i}|^{2}}
  \,\mu^{\otimes n}(d\xi)\,ds.
\end{equation}

\smallskip

\noindent
{\it Step 3: Reduction to gap variables.}

\smallskip

\noindent
The integral in~\eqref{eq:fn-simplex-bound} is over the ordered simplex $T_{n}(t)$,
and the exponent contains the nested spectral sums $|\xi_{1}+\cdots+\xi_{i}|^{2}$
that couple all $n$ variables. This step consists in two operations. We first
perform a change of time variables, which maps $T_{n}(t)$ to a product-like
domain. We then decouple the spectral variables.

\smallskip

\noindent
\textit{Change of time variables.}
Set $w_{i}:=s_{i+1}-s_{i}$ for $1\le i\le n-1$ and $w_{n}:=t-s_{n}$.
The map $s\mapsto w$ is a bijection from $T_{n}(t)$ onto
\begin{equation}\label{eq:Stn-def}
  S_{t,n}:=\{w\in[0,\infty)^{n}:\,w_{1}+\cdots+w_{n}\le t\}
\end{equation}
with unit Jacobian (it is a lower-triangular shear: $s_{i}=t-w_{n}-\cdots-w_{i}$).
The bound~\eqref{eq:fn-simplex-bound} becomes
\begin{equation}\label{eq:fn-gap-intermediate}
  n!\,\|f_{n}(\cdot,t,x)\|^{2}_{\mathcal{H}^{\otimes n}}
  \;\le\;
  \|u_{0}\|_{\infty}^{2}\,\beta^{2n}
  \int_{S_{t,n}}\int_{\R^{nd}}
  \prod_{i=1}^{n}e^{-w_{i}|\xi_{1}+\cdots+\xi_{i}|^{2}}\,
  \mu^{\otimes n}(d\xi)\,dw.
\end{equation}

\smallskip

\noindent
\textit{Decoupling of spectral variables.}
We claim that for any $w>0$ and $\eta\in\R^{d}$,
\begin{equation}\label{eq:spectral-decouple}
  \int_{\R^{d}} e^{-w|\xi+\eta|^{2}}\,\mu(d\xi)
  \;\le\;
  \int_{\R^{d}} e^{-w|\xi|^{2}}\,\mu(d\xi).
\end{equation}
In order to prove \eqref{eq:spectral-decouple}, we use the Fourier representation~\eqref{eq:heat-kernel-explicit}, namely
$e^{-w|\zeta|^{2}}=\mathcal{F}p_{2w}(\zeta)=\int_{\R^{d}}p_{2w}(x)\,e^{-\imath\zeta\cdot x}\,dx$,
as well as Fubini. This enables to write
\begin{equation}\label{d1}
  \int_{\R^{d}} e^{-w|\xi+\eta|^{2}}\,\mu(d\xi)
  \;=\;
  \int_{\R^{d}} p_{2w}(x)\,e^{-\imath\eta\cdot x}
  \Bigl[\int_{\R^{d}} e^{-\imath\xi\cdot x}\,\mu(d\xi)\Bigr]dx.
\end{equation}
Set $\Phi(x):=\int_{\R^{d}}e^{\imath\xi\cdot x}\,\mu(d\xi)$.
By Fourier inversion of $\mu=\mathcal{F}\Lambda$ we have $\Phi(x)=(2\pi)^{d}\Lambda(x)$,
which is non-negative since $\Lambda:\R^{d}\to\R_{+}$ in~\eqref{eq:noise-cov}.
By the symmetry of $\mu$ (substitute $\xi\to-\xi$), the bracket satisfies
\begin{equation*}
\int_{\R^{d}}e^{-\imath\xi\cdot x}\,\mu(d\xi)=\int_{\R^{d}}e^{\imath\xi\cdot x}\,\mu(d\xi)=\Phi(x) \ge 0.
\end{equation*}
Plugging this information in \eqref{d1}, we obtain
\begin{equation}\label{d2}
  \int_{\R^{d}} e^{-w|\xi+\eta|^{2}}\,\mu(d\xi)
  \;=\;
  \Bigl|\int_{\R^{d}} p_{2w}(x)\,e^{-\imath\eta\cdot x}\,\Phi(x)\,dx\Bigr|
  \;\le\;
  \int_{\R^{d}} p_{2w}(x)\,\Phi(x)\,dx.
\end{equation}
To identify the right-hand side, substitute $\Phi(x)=(2\pi)^{d}\Lambda(x)$
and apply Parseval's identity with $\widehat{p_{2w}}(\xi)=e^{-w|\xi|^{2}}$
(from~\eqref{eq:heat-kernel-explicit}), $\hat{\Lambda}(\xi)\,d\xi=\mu(d\xi)$,
and the evenness of $\Lambda$:
\begin{equation}\label{d3}
  \int_{\R^{d}} p_{2w}(x)\,\Phi(x)\,dx
  \;=\;
  (2\pi)^{d}\int_{\R^{d}} p_{2w}(x)\,\Lambda(x)\,dx
  \;=\;
  \int_{\R^{d}} e^{-w|\xi|^{2}}\,\mu(d\xi).
\end{equation}
Hence putting together \eqref{d2} and \eqref{d3} establishes~\eqref{eq:spectral-decouple}.

\smallskip

\noindent
Returning to~\eqref{eq:fn-gap-intermediate}, apply~\eqref{eq:spectral-decouple}
iteratively for $i=n,n-1,\ldots,1$. At step $i$, all the variables
$\xi_{i+1},\ldots,\xi_{n}$ have already been integrated out. For fixed
$\xi_{1},\ldots,\xi_{i-1}$ and $w$, we set $\eta:=\xi_{1}+\cdots+\xi_{i-1}$
and integrate over $\xi_{i}$, which yields
\begin{equation*}
  \int_{\R^{d}} e^{-w_{i}|\xi_{1}+\cdots+\xi_{i}|^{2}}\,\mu(d\xi_{i})
  \;\le\;
  \int_{\R^{d}} e^{-w_{i}|\xi_{i}|^{2}}\,\mu(d\xi_{i}).
\end{equation*}
After $n$ such steps the spectral variables are fully decoupled, and~\eqref{eq:fn-gap-intermediate}
gives
\begin{equation}\label{eq:fn-gap-bound}
  n!\,\|f_{n}(\cdot,t,x)\|^{2}_{\mathcal{H}^{\otimes n}}
  \;\le\;
  \|u_{0}\|_{\infty}^{2}\,\beta^{2n}
  \int_{\R^{nd}}\!\int_{S_{t,n}}
  e^{-\sum_{i=1}^{n}w_{i}|\xi_{i}|^{2}}
  \,dw\,\mu^{\otimes n}(d\xi),
\end{equation}
where we recall that $S_{t,n}$ is the simplex~\eqref{eq:Stn-def}, i.e.\ $S_{t,n}=\{w\in[0,\infty)^{n}:\,w_{1}+\cdots+w_{n}\le t\}$.
This is the key estimate from which Step~4 will produce a convergent series.

\smallskip

\noindent
{\it Step 4: Summation via Lemma~\ref{lem:simplex-integral}.}

\smallskip

\noindent
To evaluate the right-hand side of~\eqref{eq:fn-gap-bound}, fix $N\ge 1$ and
decompose the spectral measure into a low-frequency and a high-frequency
part by setting
\begin{equation}\label{eq:DN-CN}
  D_{N}:=\mu\bigl(\{\xi:\,|\xi|\le N\}\bigr)
  \qquad\text{and}\qquad
  C_{N}:=\int_{|\xi|>N}\frac{\mu(d\xi)}{|\xi|^{2}}.
\end{equation}
Dalang's condition~\eqref{eq:dalang-cond} ensures that $C_{N}<\infty$ for all $N$, and
$C_{N}\to 0$ as $N\to\infty$ by dominated convergence.
With $D_{N}$ and $C_{N}$ at hand, Lemma~\ref{lem:simplex-integral} below applied to
the right-hand side of~\eqref{eq:fn-gap-bound} gives
\begin{equation}\label{eq:fn-binom-bound}
  n!\,\|f_{n}(\cdot,t,x)\|^{2}_{\mathcal{H}^{\otimes n}}
  \;\le\;
  \|u_{0}\|_{\infty}^{2}\,\beta^{2n}
  \sum_{k=0}^{n}\binom{n}{k}\frac{t^{k}}{k!}\,D_{N}^{k}\,C_{N}^{n-k}.
\end{equation}
Since $C_{N}\to 0$ as $N\to\infty$, choose $N$ large enough so that
$2\beta^{2}C_{N}<1$. Using $\binom{n}{k}\le 2^{n}$, we sum over $n\ge 0$:
\begin{align}\label{eq:chaos-series-sum}
  \sum_{n=0}^{\infty}n!\,\|f_{n}(\cdot,t,x)\|^{2}_{\mathcal{H}^{\otimes n}}
  &\;\le\;
  \|u_{0}\|_{\infty}^{2}
  \sum_{k=0}^{\infty}\sum_{n=k}^{\infty}
  \beta^{2n}\cdot 2^{n}\,\frac{t^{k}}{k!}\,D_{N}^{k}\,C_{N}^{n-k} \nonumber\\
  &\;=\;
  \|u_{0}\|_{\infty}^{2}
  \sum_{k=0}^{\infty}\frac{(2\beta^{2}t\,D_{N})^{k}}{k!}
  \sum_{n=k}^{\infty}(2\beta^{2}C_{N})^{n-k} \nonumber\\
  &\;=\;
  \frac{\|u_{0}\|_{\infty}^{2}}{1-2\beta^{2}C_{N}}
  \exp\bigl(2\beta^{2}t\,D_{N}\bigr)
  \;<\;\infty.
\end{align}
This proves~\eqref{eq:chaos-L2-converge} and concludes our proof.
\end{proof}

It remains to establish the technical estimate invoked in Step~4 above. The
lemma below shows how Dalang's condition~\eqref{eq:dalang-cond} allows to split
the integration over $S_{t,n}$ and $\R^{nd}$ into a tractable product. The
high-frequency modes are controlled by $C_{N}$, that is the tail of $\mu$
weighted by $|\xi|^{-2}$. As for the low-frequency modes, they only contribute
a finite mass $D_{N}$.

\begin{lemma}[Simplex integral estimate]\label{lem:simplex-integral}
  Let $\mu$ satisfy Dalang's condition~\eqref{eq:dalang-cond}, and let $D_{N}$
  and $C_{N}$ be as in~\eqref{eq:DN-CN}. For any $n\ge 1$ and $t>0$,
  with $S_{t,n}$ the simplex defined in~\eqref{eq:Stn-def}, it holds
  \begin{equation}\label{eq:simplex-integral-bound}
    \int_{\R^{nd}}\!\int_{S_{t,n}}
    e^{-\sum_{i=1}^{n}w_{i}|\xi_{i}|^{2}}\,dw\,\mu^{\otimes n}(d\xi)
    \;\le\;
    \sum_{k=0}^{n}\binom{n}{k}\frac{t^{k}}{k!}\,D_{N}^{k}\,C_{N}^{n-k}.
  \end{equation}
\end{lemma}

\begin{proof}
  Write $1=\mathbf{1}_{\{|\xi_{i}|\le N\}}+\mathbf{1}_{\{|\xi_{i}|>N\}}$
  inside each factor $e^{-w_{i}|\xi_{i}|^{2}}$, expand the resulting product,
  and collect terms according to the subset $I\subset\{1,\ldots,n\}$
  of low-frequency indices (those with $|\xi_{i}|\le N$). This yields
  \begin{align}\label{eq:simplex-expand}
    &\int_{\R^{nd}}\!\int_{S_{t,n}}
    e^{-\sum_{i=1}^{n}w_{i}|\xi_{i}|^{2}}\,dw\,\mu^{\otimes n}(d\xi) \nonumber\\
    &\quad= \sum_{I\subset\{1,\ldots,n\}}
    \int_{S_{t,n}}\!\int_{\R^{nd}}
    \prod_{i\in I}e^{-w_{i}|\xi_{i}|^{2}}\mathbf{1}_{\{|\xi_{i}|\le N\}}
    \prod_{j\in I^{c}}e^{-w_{j}|\xi_{j}|^{2}}\mathbf{1}_{\{|\xi_{j}|>N\}}
    \,\mu^{\otimes n}(d\xi)\,dw.
  \end{align}
  For each term in~\eqref{eq:simplex-expand}, bound $e^{-w_{i}|\xi_{i}|^{2}}\le 1$
  for $i\in I$, and for $j\in I^{c}$ integrate in $w_{j}$ over $[0,t]$
  (bounding the simplex $S_{t,n}$~\eqref{eq:Stn-def} by the product simplex
  $S^{I}_{t}\times S^{I^{c}}_{t}$, where
  $S^{I}_{t}=\{(w_{i})_{i\in I}:w_{i}\ge 0,\,\sum_{i\in I}w_{i}\le t\}$). We get
  \begin{equation}\label{eq:simplex-split}
    \int_{S^{I^{c}}_{t}}
    \prod_{j\in I^{c}}e^{-w_{j}|\xi_{j}|^{2}}\,dw
    \;\le\;
    \prod_{j\in I^{c}}\int_{0}^{t}e^{-w_{j}|\xi_{j}|^{2}}\,dw_{j}
    \;=\;
    \prod_{j\in I^{c}}\frac{1-e^{-t|\xi_{j}|^{2}}}{|\xi_{j}|^{2}}
    \;\le\;
    \prod_{j\in I^{c}}\frac{1}{|\xi_{j}|^{2}}.
  \end{equation}
  For the low-frequency part, we use the trivial bound
  \begin{equation}\label{eq:low-freq-bound}
    \int_{S^{I}_{t}}\prod_{i\in I} e^{-w_{i}|\xi_{i}|^{2}} \,dw
    \;\le\;
    \int_{S^{I}_{t}}\prod_{i\in I}1\,dw
    \;\le\;
    \frac{t^{|I|}}{|I|!}.
  \end{equation}
  Combining~\eqref{eq:simplex-split} and~\eqref{eq:low-freq-bound} and integrating
  over $\mu^{\otimes n}(d\xi)$, the term in~\eqref{eq:simplex-expand} indexed by
  $I$ satisfies
  \begin{equation}\label{eq:simplex-per-I}
    \int_{S_{t,n}}\!\int_{\R^{nd}}
    \prod_{i\in I}e^{-w_{i}|\xi_{i}|^{2}}\mathbf{1}_{\{|\xi_{i}|\le N\}}
    \prod_{j\in I^{c}}e^{-w_{j}|\xi_{j}|^{2}}\mathbf{1}_{\{|\xi_{j}|>N\}}
    \,\mu^{\otimes n}(d\xi)\,dw
    \;\le\;
    \frac{t^{|I|}}{|I|!}\,D_{N}^{|I|}\,C_{N}^{n-|I|}.
  \end{equation}
  Summing~\eqref{eq:simplex-per-I} over all $I\subset\{1,\ldots,n\}$
  and collecting by $|I|=k$ (there are $\binom{n}{k}$ subsets of size $k$) yields
  \begin{eqnarray}\label{eq:simplex-final}
    \int_{\R^{nd}}\!\int_{S_{t,n}}
    e^{-\sum_{i=1}^{n}w_{i}|\xi_{i}|^{2}}\,dw\,\mu^{\otimes n}(d\xi)
    &\le&
    \sum_{k=0}^{n}\binom{n}{k}
    \frac{t^{k}}{k!}\,D_{N}^{k}\,
    \Bigl(\int_{|\xi|>N}\frac{\mu(d\xi)}{|\xi|^{2}}\Bigr)^{n-k} \notag \\
    &=&
    \sum_{k=0}^{n}\binom{n}{k}\frac{t^{k}}{k!}\,D_{N}^{k}\,C_{N}^{n-k},
  \end{eqnarray}
  which is the claimed bound~\eqref{eq:simplex-integral-bound}.
\end{proof}

\begin{remark}\label{rmk:chaos-vs-contraction}
  Theorem~\ref{thm:chaos-eu} and Theorem~\ref{thm:global-wp} (contraction)
  assert the same $L^{2}$-existence result by different methods. The chaos
  approach requires only Dalang's condition~\eqref{eq:dalang-cond} and
  $u_{0}\in L^{\infty}$, and delivers an explicit series representation with a
  concrete $L^{2}$-moment formula. The contraction approach works in the
  larger space $A_{T,p}$ for $p\ge 2$, gives $L^{p}$-estimates, and adapts to
  nonlinear multiplicative noise (Corollary~\ref{cor:nonlinear-wp}).  The two
  routes are thus complementary. The chaos expansion is essential for the
  moment analysis of Chapter~\ref{chap:moments}, while the contraction
  argument extends to a wider class of equations.
\end{remark}

\section[Feynman-Kac]{Feynman-Kac representation of moments}\label{sec:fk-moments}

The informal Feynman-Kac formula~\eqref{eq:feynman-kac} of the introduction
described the solution $u(t,x)$ as a Brownian expectation of an exponential
functional of the noise. For the \emph{moments} of the solution, a rigorous
and explicit version of this formula is available. It is one of our key tools
for the quantitative analysis of intermittency. The approach goes back to
Carmona--Molchanov~\cite{carmona:94:parabolic}, who were dealing with discrete
noises or smooth noises in space. As we shall see below, our setting is
handled by smooth approximations in space. Since the moment formula involves
$k$ independent Brownian paths, one for each factor of $u$, we first fix a
notation for this double randomness.

\begin{notation}\label{not:BM-copies}
Let $B^{1},\ldots,B^{k}$ be independent standard $d$-dimensional Brownian
motions starting at the origin, independent of the noise $W$.  We write
$\E_{B}$ for the expectation with respect to the Brownian motions only
(holding $W$ fixed), and $\E_{W}$ for the expectation with respect to
$W$ only (holding the paths fixed).
\end{notation}

\subsection{The spatially regularized process $u^{\varepsilon}$}\label{sec:spatial-reg}

For a Brownian path $B$ starting at the origin and $\varepsilon>0$, define
the $\mathcal{H}$-valued functional
\begin{equation}\label{eq:Aeps-def}
  A^{\varepsilon}_{t,x}(s,y)
  \,:=\,
  p_{\varepsilon}(B_{t-s}+x-y)\,\mathbf{1}_{[0,t]}(s_{{}})
  \qquad (s,y)\in\R_{+}\times\R^{d},
\end{equation}
where $p_{\varepsilon}$~\eqref{eq:heat-kernel-explicit} is the Gaussian heat
kernel at time $\varepsilon$.
Since $\widehat{p_{\varepsilon}}(\xi)=e^{-\varepsilon|\xi|^{2}/2}$~\eqref{eq:heat-kernel-explicit},
the Fourier formula~\eqref{eq:H-innprod-fourier} gives
\begin{equation*}
  \|A^{\varepsilon}_{t,x}\|^{2}_{\mathcal{H}}
  \,=\,
  t\int_{\R^{d}}e^{-\varepsilon|\xi|^{2}}\,\mu(d\xi),
\end{equation*}
which is finite for every $\varepsilon>0$ by Dalang's condition~\eqref{eq:dalang-cond},
so $A^{\varepsilon}_{t,x}\in\mathcal{H}$.
Next recalling that $B$ is independent of the noise $W$, and $\E_{B}$ denotes expectation
with respect to $B$ alone (see Notation~\ref{not:BM-copies}),
we define the spatially regularized process $u^{\varepsilon}$ by
\begin{equation}\label{eq:FK-eps}
  u^{\varepsilon}(t,x)
  \,:=\,
  \E_{B}\!\left[
    u_{0}(B_{t}+x)\,
    \exp\left(
      \beta\,W(A^{\varepsilon}_{t,x})
      \,-\,\frac{\beta^{2}}{2}\,\|A^{\varepsilon}_{t,x}\|^{2}_{\mathcal{H}}
    \right)
  \right].
\end{equation}
Observe that the Gaussian exponential above is well-defined path-by-path since
$A^{\varepsilon}_{t,x}\in\mathcal{H}$ almost surely. In addition, the It\^o correction
$-\tfrac{\beta^{2}}{2}\|A^{\varepsilon}_{t,x}\|^{2}_{\mathcal{H}}$
normalizes the conditional expectation. Namely, as we will see in~\eqref{e1}, we have
\begin{equation*}
\E_{W}\bigl[e^{\beta W(A^{\varepsilon}_{t,x})-\frac{\beta^{2}}{2}\|A^{\varepsilon}_{t,x}\|^{2}_{\mathcal{H}}}\bigr]=1
\end{equation*}
for every fixed Brownian path $B$.
We will show in Section~\ref{sec:FK} that $u^{\varepsilon}$ converges to the
mild solution $u$ and that~\eqref{eq:FK-eps} leads to the Feynman-Kac moment
formula (see~\eqref{eq:fk-moment-formula} below).  As a first step, we establish
uniform-in-$\varepsilon$ moment bounds.

\begin{proposition}[Uniform moment bounds]\label{prop:ueps-moments}
Let $u_{0}\in L^{\infty}(\R^{d})$ and suppose Dalang's
condition~\eqref{eq:dalang-cond} holds.  Then for every integer $k\ge 1$,
$t>0$, and $x\in\R^{d}$, we have
\begin{equation}\label{eq:ueps-moment-bound}
  \sup_{\varepsilon>0}\,\E\bigl[\left|u^{\varepsilon}(t,x)\right|^{k}\bigr]<\infty.
\end{equation}
\end{proposition}

\begin{proof}
Let $B^{1},\ldots,B^{k}$ be the independent standard Brownian motions of
Notation~\ref{not:BM-copies}, and write $A^{\varepsilon,j}_{t,x}$ for the
functional~\eqref{eq:Aeps-def} built from $B^{j}$, $j=1,\ldots,k$.
We proceed in two steps.

\smallskip

\noindent
{\it Step 1: Reduction to inner products.}

\smallskip

\noindent
Since $B^{1},\ldots,B^{k}$ are independent copies independent of $W$,
substituting~\eqref{eq:FK-eps} for each factor and invoking Fubini's theorem gives
\begin{equation*}
  \E\bigl[(u^{\varepsilon}(t,x))^{k}\bigr]
  \;=\;
  \E_{B}\!\left[
    \prod_{j=1}^{k}u_{0}(B^{j}_{t}+x)\;
    \E_{W}\!\left[
      \exp\left(
        \beta\sum_{j=1}^{k}W(A^{\varepsilon,j}_{t,x})
        -\frac{\beta^{2}}{2}\sum_{j=1}^{k}
        \|A^{\varepsilon,j}_{t,x}\|^{2}_{\mathcal{H}}
      \right)
    \right]
  \right].
\end{equation*}
For fixed paths, $\beta\sum_{j}W(A^{\varepsilon,j}_{t,x})
=W\!\bigl(\beta\sum_{j}A^{\varepsilon,j}_{t,x}\bigr)$
is Gaussian with variance
$\beta^{2}\|\sum_{j}A^{\varepsilon,j}_{t,x}\|^{2}_{\mathcal{H}}$
by the isometry~\eqref{eq:W-isometry}, so
\begin{equation}\label{e1}
  \E_{W}\!\left[
    \exp\left(\beta\sum_{j}W(A^{\varepsilon,j}_{t,x})\right)
  \right]
  \;=\;
  \exp\left(
    \frac{\beta^{2}}{2}
    \Bigl\|\sum_{j}A^{\varepsilon,j}_{t,x}\Bigr\|^{2}_{\mathcal{H}}
  \right).
\end{equation}
Expanding $\|\sum_{j}A^{\varepsilon,j}_{t,x}\|^{2}_{\mathcal{H}}
=\sum_{j}\|A^{\varepsilon,j}_{t,x}\|^{2}_{\mathcal{H}}
+2\sum_{i<j}\langle A^{\varepsilon,i}_{t,x},A^{\varepsilon,j}_{t,x}\rangle_{\mathcal{H}}$
and canceling the diagonal terms gives
\begin{equation}\label{eq:mom-eps-formula}
  \E\bigl[(u^{\varepsilon}(t,x))^{k}\bigr]
  \,=\,
  \E_{B}\!\left[
    \prod_{i=1}^{k}u_{0}(B^{i}_{t}+x)\,
    \exp\left(
      \beta^{2}\!\!\sum_{1\le i<j\le k}
      \langle A^{\varepsilon,i}_{t,x},A^{\varepsilon,j}_{t,x}\rangle_{\mathcal{H}}
    \right)
  \right].
\end{equation}
Since $u_{0}\in L^{\infty}$, it suffices to bound
$\E_{B}[\exp(\beta^{2}\sum_{i<j}\langle A^{\varepsilon,i},A^{\varepsilon,j}\rangle_{\mathcal{H}})]$
uniformly in $\varepsilon$.

\smallskip

\noindent
{\it Step 2: Uniform exponential bound.}

\smallskip

\noindent
By the Fourier formula~\eqref{eq:H-innprod-fourier} and
$\widehat{p_{\varepsilon}}(\xi)=e^{-\varepsilon|\xi|^{2}/2}$~\eqref{eq:heat-kernel-explicit},
\begin{equation}\label{eq:Aeps-inner}
  \langle A^{\varepsilon,i}_{t,x},A^{\varepsilon,j}_{t,x}\rangle_{\mathcal{H}}
  \;=\;
  \int_{0}^{t}\!\int_{\R^{d}}
  e^{-\varepsilon|\xi|^{2}}\,e^{i\xi\cdot(B^{j}_{r}-B^{i}_{r})}\,\mu(d\xi)\,dr
  \;=:\;\int_{0}^{t}\Lambda_{\varepsilon}(B^{j}_{r}-B^{i}_{r})\,dr,
\end{equation}
where the change of variable $r=t-s$ was used and
$\Lambda_{\varepsilon}(z):=\int_{\R^{d}}e^{-\varepsilon|\xi|^{2}}e^{i\xi\cdot z}\mu(d\xi)$
is real and non-negative (as a spatial average of $\Lambda\ge0$~\eqref{eq:noise-cov}
against a Gaussian weight).  In particular each inner product
$\langle A^{\varepsilon,i},A^{\varepsilon,j}\rangle_{\mathcal{H}}$ is positive.

Set $M:=\binom{k}{2}$ and $\lambda_{0}:=M\beta^{2}$.
Since every inner product $\langle A^{\varepsilon,i},A^{\varepsilon,j}\rangle_{\mathcal{H}}\ge0$,
the trivial bound $\sum_{i<j}a_{ij}\le M\max_{i<j}a_{ij}$ and the
inequality $e^{M\max_{i<j}a_{ij}}\le\sum_{i<j}e^{Ma_{ij}}$ together give
\begin{equation*}
  \E_{B}\!\left[
    \exp\left(\beta^{2}\!\sum_{i<j}
    \langle A^{\varepsilon,i},A^{\varepsilon,j}\rangle_{\mathcal{H}}\right)
  \right]
  \;\le\;
  \sum_{1\le i<j\le k}
  \E_{B^{i},B^{j}}\!\left[
    \exp\left(\lambda_{0}
    \langle A^{\varepsilon,i},A^{\varepsilon,j}\rangle_{\mathcal{H}}\right)
  \right].
\end{equation*}
Therefore by symmetry it suffices to bound
$\E[\exp(\lambda_{0}\langle A^{\varepsilon,1},A^{\varepsilon,2}\rangle_{\mathcal{H}})]$
uniformly in $\varepsilon$.
To this aim, expanding the exponential and using~\eqref{eq:Aeps-inner} we can write
\begin{equation}\label{eq:An-n-power}
  \E\bigl[\langle A^{\varepsilon,1}_{t,x},A^{\varepsilon,2}_{t,x}
  \rangle_{\mathcal{H}}^{n}\bigr]
  \;=\;
  \int_{[0,t]^{n}}\!\int_{\R^{nd}}
  e^{-\varepsilon\sum_{l=1}^{n}|\xi_{l}|^{2}}\,
  \E\!\left[e^{i\sum_{l=1}^{n}\xi_{l}\cdot(B^{2}_{r_{l}}-B^{1}_{r_{l}})}\right]
  \mu^{\otimes n}(d\xi)\,dr.
\end{equation}
Furthermore, since $B^{1}$ and $B^{2}$ are independent standard Brownian motions, their
characteristic function can be expressed as
\begin{equation*}
  \E\!\left[e^{i\sum_{l}\xi_{l}\cdot(B^{2}_{r_{l}}-B^{1}_{r_{l}})}\right]
  \;=\;
  e^{-\sum_{l,m}(r_{l}\wedge r_{m})\,\xi_{l}\cdot\xi_{m}}.
\end{equation*}
Dropping $e^{-\varepsilon\sum|\xi_{l}|^{2}}\le 1$, restricting the integral
in~\eqref{eq:An-n-power} by symmetry to the simplex 
$T_{n}(t)$ defined by~\eqref{eq:simplex} (with factor $n!$), and writing $w_{l}=r_{l}-r_{l-1}$,
the standard identity
\begin{equation}\label{eq:quadform}
  \sum_{l,m}(r_{l}\wedge r_{m})\,\xi_{l}\cdot\xi_{m}
  \;=\;
  \sum_{l=1}^{n}w_{l}\,|\xi_{l}+\cdots+\xi_{n}|^{2}
  \quad\text{on }T_{n}(t)
\end{equation}
yields
\begin{equation*}
  \E\bigl[\langle A^{\varepsilon,1},A^{\varepsilon,2}\rangle_{\mathcal{H}}^{n}\bigr]
  \;\le\;
  n!\int_{S_{t,n}}\!\int_{\R^{nd}}
  \exp\left(-\sum_{l=1}^{n}w_{l}|\xi_{l}+\cdots+\xi_{n}|^{2}\right)
  \mu^{\otimes n}(d\xi)\,dw.
\end{equation*}
From this point the argument is the same as in the proof of
Theorem~\ref{thm:chaos-eu} (Steps~3--4).
Applying~\eqref{eq:spectral-decouple} iteratively (integrate first over $\xi_{1}$
with $\xi_{2},\ldots,\xi_{n}$ fixed, then over $\xi_{2}$, etc.) decouples the
spectral variables:
\begin{equation*}
  \int_{\R^{nd}}
  \exp\left(-\sum_{l=1}^{n}w_{l}|\xi_{l}+\cdots+\xi_{n}|^{2}\right)
  \mu^{\otimes n}(d\xi)
  \;\le\;
  \prod_{l=1}^{n}\int_{\R^{d}}e^{-w_{l}|\xi|^{2}}\,\mu(d\xi).
\end{equation*}
Next by Lemma~\ref{lem:simplex-integral}~\eqref{eq:simplex-integral-bound} we obtain:
\begin{equation*}
  \E\bigl[\langle A^{\varepsilon,1},A^{\varepsilon,2}\rangle_{\mathcal{H}}^{n}\bigr]
  \;\le\;
  n!\sum_{j=0}^{n}\binom{n}{j}\frac{t^{j}}{j!}\,D_{N}^{j}\,C_{N}^{n-j}.
\end{equation*}
Choose $N$ large enough so that $2\lambda_{0}C_{N}<1$.  By the
same summation as~\eqref{eq:chaos-series-sum}, we end up with
\begin{equation*}
  \E\!\left[\exp\left(
  \lambda_{0}\langle A^{\varepsilon,1},A^{\varepsilon,2}\rangle_{\mathcal{H}}\right)\right]
  \;=\;
  \sum_{n=0}^{\infty}\frac{\lambda_{0}^{n}}{n!}
  \E\bigl[\langle A^{\varepsilon,1},A^{\varepsilon,2}\rangle_{\mathcal{H}}^{n}\bigr]
  \;\le\;
  \frac{\exp(2\lambda_{0}t\,D_{N})}{1-2\lambda_{0}C_{N}}
  \;<\;\infty,
\end{equation*}
uniformly in $\varepsilon>0$ (since we dropped $e^{-\varepsilon\sum|\xi_{l}|^{2}}$
before applying the lemma).  Combined with~\eqref{eq:mom-eps-formula}
and $u_{0}\in L^{\infty}$, this proves~\eqref{eq:ueps-moment-bound}.
\end{proof}

With the uniform moment bounds of Proposition~\ref{prop:ueps-moments} in hand,
we can now identify the mild SHE that $u^{\varepsilon}$ satisfies.
It turns out that the Feynman-Kac representation~\eqref{eq:FK-eps}
forces $u^{\varepsilon}$ to solve the mild SHE driven by the
spatially regularized noise $W^{\varepsilon}:=p_{\varepsilon}*W$.

\begin{lemma}[Mild SHE with regularized noise]\label{lem:fk-duality}
Let $u^{\varepsilon}$ be as in~\eqref{eq:FK-eps}.
Define the spatially regularized noise $W^{\varepsilon}$ by setting,
for any adapted process $\phi\in\mathcal{H}$,
\begin{equation}\label{eq:Weps-def}
  \int_{0}^{T}\!\int_{\R^{d}}\phi(s,y)\,W^{\varepsilon}(ds,dy)
  \;:=\;
  \int_{0}^{T}\!\int_{\R^{d}}
  (p_{\varepsilon}*_{y}\phi_{s})(y)\,W(ds,dy),
\end{equation}
where $(p_{\varepsilon}*_{y}\phi_{s})(y):=\int_{\R^{d}}p_{\varepsilon}(y-z)\,\phi(s,z)\,dz$
is the spatial convolution of $\phi_{s}$ with $p_{\varepsilon}$.
Then for all $(t,x)\in\R_{+}\times\R^{d}$,
\begin{equation}\label{eq:mild-ueps-reg}
  u^{\varepsilon}(t,x)
  \;=\; p_{t}u_{0}(x)
  \;+\; \beta\int_{0}^{t}\!\int_{\R^{d}}
        p_{t-s}(x-y)\,u^{\varepsilon}(s,y)\,W^{\varepsilon}(ds,dy).
\end{equation}
\end{lemma}

\begin{proof}
The argument proceeds in two steps.
In Step~1, we derive a mild It\^{o}--Walsh equation for $u^{\varepsilon}$
by expanding the Feynman-Kac exponential via It\^{o}'s formula,
averaging over the Brownian motion $B$, and identifying the resulting
kernel via a conditional independence argument.
In Step~2, we rewrite that equation as the regularized mild
SHE~\eqref{eq:mild-ueps-reg} by recognizing the Walsh kernel
as an integral against the regularized noise $W^{\varepsilon}$,
using the stochastic Fubini theorem (Proposition~\ref{prop:stoch-fubini}).

\smallskip

\noindent
{\it Step 1: Mild It\^{o}--Walsh equation for $u^{\varepsilon}$.}
In order to obtain a mild equation, introduce for $(s,y)\in[0,t]\times\R^{d}$ the kernel
\begin{equation}\label{eq:Zeps-def}
  Z^{\varepsilon,t,x}(s,y)
  \,:=\,
  \int_{\R^{d}}p_{t-s}(x-z)\,p_{\varepsilon}(z-y)\,u^{\varepsilon}(s,z)\,dz.
\end{equation}
We claim that $u^{\varepsilon}$ satisfies the mild equation
\begin{equation}\label{eq:mild-ueps}
  u^{\varepsilon}(t,x)
  \;=\; p_{t}u_{0}(x)
  \;+\; \beta\int_{0}^{t}\!\int_{\R^{d}}
        Z^{\varepsilon,t,x}(s,y)\,W(ds,dy).
\end{equation}

\noindent
{\it Step 1-a: Stochastic exponential.}
For a fixed Brownian path $B$, set
\begin{equation}\label{eq:NsB}
  N_{s}^{B}
  \;:=\; \beta\int_{0}^{s}\!\int_{\R^{d}}
         p_{\varepsilon}(B_{t-r}+x-y)\,W(dr,dy),
  \qquad s\in[0,t].
\end{equation}
For fixed $B$, the space-time function $(r,y)\mapsto p_{\varepsilon}(B_{t-r}+x-y)\,\mathbf{1}_{[0,t]}(r)$
is a deterministic element of $\mathcal{H}$ (by the same norm computation
as for $A^{\varepsilon}_{t,x}$ in~\eqref{eq:Aeps-def}), so $N^{B}$ is a
$W$-martingale.  By the Walsh isometry~\eqref{eq:walsh-isometry}, its quadratic
variation is
\begin{equation*}
  \langle N^{B}\rangle_{s}
  \;=\; \beta^{2}\!\int_{0}^{s}\!\int_{\R^{2d}}
        p_{\varepsilon}(B_{t-r}+x-y_{1})\,p_{\varepsilon}(B_{t-r}+x-y_{2})\,
        \Lambda(y_{1}-y_{2})\,dy_{1}\,dy_{2}\,dr.
\end{equation*}
The change of variables $z_{i}:=B_{t-r}+x-y_{i}$ and the symmetry
$\Lambda(y_{1}-y_{2})=\Lambda(z_{1}-z_{2})$~\eqref{eq:noise-cov} give
\begin{multline}\label{e3}
  \int_{\R^{2d}}
  p_{\varepsilon}(B_{t-r}+x-y_{1})\,p_{\varepsilon}(B_{t-r}+x-y_{2})\,
  \Lambda(y_{1}-y_{2})\,dy_{1}\,dy_{2}  \\
  \;=\;
  \int_{\R^{2d}}
  p_{\varepsilon}(z_{1})\,p_{\varepsilon}(z_{2})\,
  \Lambda(z_{1}-z_{2})\,dz_{1}\,dz_{2}.
\end{multline}
By the Fourier identity~\eqref{eq:H-innprod-fourier} (restricted to a single
time slice) and $\widehat{p_{\varepsilon}}(\xi)=e^{-\varepsilon|\xi|^{2}/2}$~\eqref{eq:heat-kernel-explicit},
the right-hand side of~\eqref{e3} equals
$\int_{\R^{d}}e^{-\varepsilon|\xi|^{2}}\mu(d\xi)=\Lambda_{\varepsilon}(0)$
(the same computation as in~\eqref{eq:Aeps-inner} with $i=j$
and $B^{i}_{r}=B^{j}_{r}$).  Hence
\begin{equation}\label{e4}
  \langle N^{B}\rangle_{s}
  \;=\; \beta^{2}\!\int_{0}^{s}\Lambda_{\varepsilon}(0)\,dr
  \;=\; \beta^{2}s\,\Lambda_{\varepsilon}(0).
\end{equation}
The Dol\'eans--Dade stochastic exponential of $N^{B}$ is thus
\begin{equation}\label{eq:MsB}
  M_{s}^{B}
  \;:=\; \exp\Bigl(N_{s}^{B}-\tfrac{1}{2}\langle N^{B}\rangle_{s}\Bigr).
\end{equation}
At $s=t$ this coincides with the Feynman-Kac exponential of~\eqref{eq:FK-eps}.
Indeed, evaluating~\eqref{eq:NsB} at $s=t$ and recalling from~\eqref{eq:Aeps-def} that $A^{\varepsilon}_{t,x}(r,y)=p_{\varepsilon}(B_{t-r}+x-y)\,\mathbf{1}_{[0,t]}(r)$ we obtain
\begin{equation*}
  N^{B}_{t}
  \;=\; \beta\int_{0}^{t}\!\int_{\R^{d}}p_{\varepsilon}(B_{t-r}+x-y)\,W(dr,dy)
  \;=\; \beta\,W(A^{\varepsilon}_{t,x}).
\end{equation*}
For the correction term, setting $i=j$ in~\eqref{eq:Aeps-inner}
 gives
 \begin{equation*}
\|A^{\varepsilon}_{t,x}\|^{2}_{\mathcal{H}}=\int_{0}^{t}\Lambda_{\varepsilon}(0)\,dr=t\,\Lambda_{\varepsilon}(0) 
=
\frac{\langle N^{B}\rangle_{t}}{\beta^{2}},
\end{equation*}
where the last identity stems from~\eqref{e3}.
Therefore
\begin{equation}\label{eq:MtB-fk}
  M_{t}^{B}
  \;=\;\exp\Bigl(N^{B}_{t}-\tfrac{1}{2}\langle N^{B}\rangle_{t}\Bigr)
  \;=\;\exp\Bigl(\beta\,W(A^{\varepsilon}_{t,x})-\tfrac{\beta^{2}}{2}\|A^{\varepsilon}_{t,x}\|^{2}_{\mathcal{H}}\Bigr).
\end{equation}
Setting $X_{s}:=N^{B}_{s}-\tfrac{1}{2}\langle N^{B}\rangle_{s}$,
so that $M^{B}_{s}=e^{X_{s}}$~\eqref{eq:MsB}, 
It\^{o}'s formula applied to $f(x)=e^{x}$ and the semimartingale $X$ gives
\begin{eqnarray*}
  dM^{B}_{s}
  &=& M^{B}_{s}\,dX_{s}+\tfrac{1}{2}M^{B}_{s}\,d\langle X\rangle_{s} \\
  &=& M^{B}_{s}\!\left(dN^{B}_{s}-\tfrac{1}{2}d\langle N^{B}\rangle_{s}\right)
        +\tfrac{1}{2}M^{B}_{s}\,d\langle N^{B}\rangle_{s}
  \;=\; M^{B}_{s}\,dN^{B}_{s}.
\end{eqnarray*}
Substituting the differential $dN^{B}_{s}=\beta\int_{\R^{d}}p_{\varepsilon}(B_{t-s}+x-y)\,W(ds,dy)$
from~\eqref{eq:NsB} then yields
\begin{equation}\label{eq:MsB-sde}
  dM_{s}^{B}
  \;=\; \beta\,M_{s}^{B}\!\int_{\R^{d}}p_{\varepsilon}(B_{t-s}+x-y)\,W(ds,dy).
\end{equation}
so integrating from $0$ to $t$ and using $M_{0}^{B}=1$ we end up with:
\begin{equation}\label{eq:MtB-expand}
  M_{t}^{B}
  = 1 + \beta\int_{0}^{t}\!\int_{\R^{d}}
    M_{s}^{B}\,p_{\varepsilon}(B_{t-s}+x-y)\,W(ds,dy).
\end{equation}

\noindent
{\it Step 1-b: Averaging over $B$.}
By~\eqref{eq:FK-eps} and~\eqref{eq:MtB-fk}, observe that
\begin{equation*}
  u^{\varepsilon}(t,x)
  \;=\; \E_{B}\!\bigl[u_{0}(B_{t}+x)\,M^{B}_{t}\bigr].
\end{equation*}
Substituting~\eqref{eq:MtB-expand} for $M^{B}_{t}$, using that $u_{0}(B_{t}+x)$
is $W$-independent for fixed $B$, and applying Proposition~\ref{prop:stoch-fubini}
(with $\nu$ the law of $B$) to exchange $\E_{B}$ with the Walsh integral, we thus get
\begin{equation}\label{eq:ueps-prefk}
  u^{\varepsilon}(t,x)
  = p_{t}u_{0}(x)
  + \beta\int_{0}^{t}\!\int_{\R^{d}}
    \E_{B}\!\bigl[u_{0}(B_{t}+x)\,M_{s}^{B}\,p_{\varepsilon}(B_{t-s}+x-y)\bigr]\,
    W(ds,dy).
\end{equation}

\noindent
{\it Step 1-c: Identification of the kernel.}
Our next objective is to show that
\begin{equation}\label{e5}
\E_{B}[u_{0}(B_{t}+x)\,M_{s}^{B}\,p_{\varepsilon}(B_{t-s}+x-y)]
=Z^{\varepsilon,t,x}(s,y) ,
\end{equation}
with $Z^{\varepsilon,t,x}$ as defined in~\eqref{eq:Zeps-def}.
To this aim,
condition on $B_{t-s}$ and set $b:=B_{t-s}+x$.
Write $\tilde{B}_{v}:=B_{t-s+v}-B_{t-s}$ for $v\in[0,s]$. This is a
Brownian motion starting at $0$, independent of $B_{t-s}$.  Then for $r\in[0,s]$
one can write
\begin{equation}\label{e6}
  B_{t-r}+x = b+\tilde{B}_{s-r} ,
  \qquad
  B_{t}+x = b+\tilde{B}_{s},
  \qquad
  p_{\varepsilon}(B_{t-s}+x-y)=p_{\varepsilon}(b-y),
\end{equation}
Substituting the first identity of~\eqref{e6} 
into the integrand of $N_{s}^{B}$ in~\eqref{eq:NsB},
and using~\eqref{e4} for the quadratic variation,
the definition~\eqref{eq:MsB} of $M_{s}^{B}$ gives
\begin{equation}\label{eq:MsB-conditioned}
  M_{s}^{B}\big|_{B_{t-s}=b-x}
  =\exp\left(\beta\int_{0}^{s}\!\int_{\R^{d}}
   p_{\varepsilon}(b+\tilde{B}_{s-r}-y)\,W(dr,dy)
   -\tfrac{\beta^{2}}{2}s\,\Lambda_{\varepsilon}(0)\right) ,
\end{equation}
where we recall that $(\tilde{B}_{v})_{v\in[0,s]}$ is now a standard Brownian motion starting at $0$.
Therefore comparing~\eqref{eq:MsB-conditioned} with
the Feynman-Kac formula~\eqref{eq:FK-eps} applied at $(s,b)$, we obtain
\begin{equation}\label{eq:kern-id}
  \E_{\tilde{B}}\!\left[u_{0}(b+\tilde{B}_{s})\,M_{s}^{B}\big|_{B_{t-s}=b-x}\right]
  \;=\; u^{\varepsilon}(s,b).
\end{equation}
Integrating over the law $\Prob_{B}(B_{t-s}+x\in db)=p_{t-s}(x-b)\,db$ and
using~\eqref{eq:kern-id}, as well as the definition~\eqref{eq:Zeps-def} of $Z^{\varepsilon,t,x}$, we obtain
\begin{equation*}
  \E_{B}\!\bigl[u_{0}(B_{t}+x)\,M_{s}^{B}\,p_{\varepsilon}(B_{t-s}+x-y)\bigr]
  = \int_{\R^{d}}p_{t-s}(x-b)\,p_{\varepsilon}(b-y)\,u^{\varepsilon}(s,b)\,db
  = Z^{\varepsilon,t,x}(s,y),
\end{equation*}
which proves~\eqref{e5}.
Substituting~\eqref{e5} back into~\eqref{eq:ueps-prefk}
replaces the kernel $\E_{B}\!\bigl[u_{0}(B_{t}+x)\,M_{s}^{B}\,p_{\varepsilon}(B_{t-s}+x-y)\bigr]$
by $Z^{\varepsilon,t,x}(s,y)$, yielding~\eqref{eq:mild-ueps}.

\smallskip

\noindent
{\it Step 2: Identification as the regularized mild SHE~\eqref{eq:mild-ueps-reg}.}
Substituting the definition~\eqref{eq:Zeps-def} of $Z^{\varepsilon,t,x}$
into the Walsh integral in~\eqref{eq:mild-ueps}, and using the symmetry
$p_{\varepsilon}(z-y)=p_{\varepsilon}(y-z)$, gives
\begin{equation*}
  \int_{0}^{t}\!\int_{\R^{d}}Z^{\varepsilon,t,x}(s,y)\,W(ds,dy)
  \;=\;
  \int_{0}^{t}\!\int_{\R^{d}}
  \Bigl(\int_{\R^{d}}p_{\varepsilon}(y-z)\,p_{t-s}(x-z)\,u^{\varepsilon}(s,z)\,dz\Bigr)
  \,W(ds,dy).
\end{equation*}
By Proposition~\ref{prop:stoch-fubini} (applied with $G=\R^{d}$ and $\nu(dz)=dz$) and the
definition~\eqref{eq:Weps-def} of $W^{\varepsilon}$
(applied with $\phi(s,y):=p_{t-s}(x-y)\,u^{\varepsilon}(s,y)$), it is readily checked that
\begin{equation*}
  \int_{0}^{t}\!\int_{\R^{d}}Z^{\varepsilon,t,x}(s,y)\,W(ds,dy)
  \;=\;
  \int_{0}^{t}\!\int_{\R^{d}}
  p_{t-s}(x-y)\,u^{\varepsilon}(s,y)\,W^{\varepsilon}(ds,dy).
\end{equation*}
Substituting back into~\eqref{eq:mild-ueps} yields~\eqref{eq:mild-ueps-reg}.
\end{proof}

The mild equation~\eqref{eq:mild-ueps} also admits a reformulation in terms
of Malliavin calculus. Namely, test the equation against a functional
$F\in\mathbb{D}^{1,2}$ and apply the Skorohod duality relation. This replaces
the Walsh integral by an $\mathcal{H}$-inner product involving the Malliavin
derivative $DF$.

\begin{lemma}[Mild equation in duality form]\label{lem:mild-duality}
Let $u^{\varepsilon}$ be as in~\eqref{eq:FK-eps} and let
$Z^{\varepsilon,t,x}$ be the kernel of~\eqref{eq:Zeps-def}.
Then $Z^{\varepsilon,t,x}\in\mathrm{Dom}(\delta)$, and for every
$F\in\mathbb{D}^{1,2}$,
\begin{equation}\label{eq:fk-duality}
  \E\bigl[F\,u^{\varepsilon}(t,x)\bigr]
  \;=\;
  \E[F]\,p_{t}u_{0}(x)
  \;+\;
  \beta\,\E\bigl[\langle DF,Z^{\varepsilon,t,x}\rangle_{\mathcal{H}}\bigr].
\end{equation}
\end{lemma}

\begin{proof}
Multiply~\eqref{eq:mild-ueps} by $F$ and take expectations:
\begin{equation*}
  \E\bigl[F\,u^{\varepsilon}(t,x)\bigr]
  \;=\;
  \E[F]\,p_{t}u_{0}(x)
  \;+\;
  \beta\,\E\!\left[F\int_{0}^{t}\!\int_{\R^{d}}
  Z^{\varepsilon,t,x}(s,y)\,W(ds,dy)\right].
\end{equation*}
Since $u^{\varepsilon}(s,z)$ is $\mathcal{F}_{s}$-measurable
by~\eqref{eq:FK-eps}, the kernel $Z^{\varepsilon,t,x}(s,y)$
of~\eqref{eq:Zeps-def} is also $\mathcal{F}_{s}$-measurable
for each $(s,y)$, hence adapted to $(\mathcal{F}_{s})_{s\ge 0}$
of~\eqref{eq:filtration}.
Proposition~\ref{prop:ito-sko} (formula~\eqref{eq:ito-sko-H}) therefore gives
\begin{equation}\label{eq:Zeps-skorohod}
  \int_{0}^{t}\!\int_{\R^{d}}
  Z^{\varepsilon,t,x}(s,y)\,W(ds,dy)
  \;=\; \delta\bigl(Z^{\varepsilon,t,x}\bigr),
\end{equation}
so in particular $Z^{\varepsilon,t,x}\in\mathrm{Dom}(\delta)$.
Applying the Malliavin duality relation~\eqref{eq:duality} gives
$\E[F\,\delta(Z^{\varepsilon,t,x})]
=\E[\langle DF,Z^{\varepsilon,t,x}\rangle_{\mathcal{H}}]$,
which yields~\eqref{eq:fk-duality}.
\end{proof}

\begin{remark}[Chaos expansion of $u^{\varepsilon}$]\label{rem:ueps-chaos}
Iterating the mild equation~\eqref{eq:mild-ueps} exactly as in the identification of the chaos kernels of $u$ (Lemma~\ref{lem:chaos-kernels}) shows that $u^{\varepsilon}$ has its own chaos expansion. Specifically one can write
\begin{equation}\label{eq:ueps-chaos}
  u^{\varepsilon}(t,x)
  =
  p_{t}u_{0}(x)+\sum_{n\ge 1}I_{n}\bigl(f^{\varepsilon}_{n}(\cdot,t,x)\bigr),
\end{equation}
where $f^{\varepsilon}_{n}$ is obtained from the kernel $f_{n}$ of~\eqref{eq:fn-def} by a convolution with $p_{\varepsilon}$ in each of its $n$ space variables. Namely, suppressing the time variables $s_{i}$ and the parameters $t,x$ in the kernels and writing the spatial arguments only for notational sake, it holds
\begin{equation}\label{eq:fneps-def}
  f^{\varepsilon}_{n}(y_{1},\ldots,y_{n})
  :=
  \int_{\R^{nd}}f_{n}(z_{1},\ldots,z_{n})\prod_{i=1}^{n}p_{\varepsilon}(y_{i}-z_{i})\,dz_{i}.
\end{equation}
Indeed, each iteration of~\eqref{eq:mild-ueps} inserts the kernel $Z^{\varepsilon,t,x}$ of~\eqref{eq:Zeps-def}, which contains one convolution with $p_{\varepsilon}$ in the noise variable $y$, while the zeroth-order term $p_{t}u_{0}(x)$ is the same as for $u$. Taking Fourier transforms in the space variables and recalling that $\widehat{p_{\varepsilon}}(\xi)=e^{-\varepsilon|\xi|^{2}/2}$~\eqref{eq:heat-kernel-explicit}, relation~\eqref{eq:fneps-def} can be recast as
\begin{equation}\label{eq:fneps-fourier}
  \mathcal{F}_{x}^{\otimes n}f^{\varepsilon}_{n}(s,\xi)
  =
  \prod_{i=1}^{n}e^{-\varepsilon|\xi_{i}|^{2}/2}\,
  \mathcal{F}_{x}^{\otimes n}f_{n}(s,\xi),
  \qquad (s,\xi)\in\R_{+}^{n}\times\R^{nd},
\end{equation}
where $\mathcal{F}_{x}^{\otimes n}$ denotes the Fourier transform in the $n$ space variables, in coherence with the notation of~\eqref{eq:H-innprod-fourier}.
\end{remark}

The expansions~\eqref{eq:u-chaos} and~\eqref{eq:ueps-chaos} live on the same chaos levels, so that $u^{\varepsilon}(t,x)-u(t,x)$ can be compared kernel by kernel. Combined with Nelson's hypercontractivity inequality~\eqref{eq:hypercontractivity-sum}, this observation yields the following convergence result.

\begin{proposition}[$L^{k}$ convergence of $u^{\varepsilon}$]\label{prop:ueps-Lk-conv}
Let $u_{0}\in L^{\infty}(\R^{d})$ and suppose that Dalang's condition~\eqref{eq:dalang-cond} holds. Let $u$ be the mild solution of~\eqref{eq:mild-pam} from Theorem~\ref{thm:chaos-eu}, with chaos kernels $f_{n}$ of~\eqref{eq:fn-def}. Let also $u^{\varepsilon}$ be the regularized process~\eqref{eq:FK-eps}. Then for every integer $k\ge 2$, every $t>0$ and every $x\in\R^{d}$, it holds:
\begin{equation*}
  u^{\varepsilon}(t,x)\;\longrightarrow\;u(t,x)
  \qquad\text{in } L^{k}(\Omega),\ \text{as } \varepsilon\to 0.
\end{equation*}
\end{proposition}

\begin{proof}
Owing to~\eqref{eq:u-chaos} and~\eqref{eq:ueps-chaos}, the difference $u^{\varepsilon}(t,x)-u(t,x)$ admits the chaos decomposition $\sum_{n\ge1}I_{n}(f^{\varepsilon}_{n}-f_{n})$. We bound its $L^{k}(\Omega)$ norm in three steps. We first record a contraction property of the kernels $f^{\varepsilon}_{n}$. We then prove a $k$-dependent summability estimate for the series of chaos norms. We finally conclude by hypercontractivity.

\smallskip

\noindent
{\it Step 1: Contraction and convergence at fixed chaos level.} Since $0<e^{-\varepsilon|\xi|^{2}/2}\le 1$, the Fourier expression~\eqref{eq:fneps-fourier} together with the Fourier formula~\eqref{eq:H-innprod-fourier} for the norm in $\mathcal{H}^{\otimes n}$ yields
\begin{equation}\label{eq:fneps-contract}
  \|f^{\varepsilon}_{n}(\cdot,t,x)\|_{\mathcal{H}^{\otimes n}}
  \;\le\;
  \|f_{n}(\cdot,t,x)\|_{\mathcal{H}^{\otimes n}},
  \qquad\text{for all }\varepsilon>0,\ n\ge 1.
\end{equation}
Resorting to~\eqref{eq:fneps-fourier} and~\eqref{eq:H-innprod-fourier} again, and dropping the arguments $(\cdot,t,x)$ for notational sake, we also get
\begin{equation}\label{eq:fneps-diff}
  \|f^{\varepsilon}_{n}-f_{n}\|^{2}_{\mathcal{H}^{\otimes n}}
  =
  \int_{\R_{+}^{n}}\int_{\R^{nd}}
  \Bigl(1-\prod_{i=1}^{n}e^{-\varepsilon|\xi_{i}|^{2}/2}\Bigr)^{2}
  \bigl|\mathcal{F}_{x}^{\otimes n}f_{n}(s,\xi)\bigr|^{2}
  \prod_{i=1}^{n}\mu(d\xi_{i})\,ds.
\end{equation}
The factor $(1-\prod_{i=1}^{n}e^{-\varepsilon|\xi_{i}|^{2}/2})^{2}$ in~\eqref{eq:fneps-diff} is bounded by $1$ and converges pointwise to $0$ as $\varepsilon\to0$, while $|\mathcal{F}_{x}^{\otimes n}f_{n}|^{2}$ is integrable (we have $\|f_{n}(\cdot,t,x)\|_{\mathcal{H}^{\otimes n}}<\infty$ by~\eqref{eq:chaos-series-sum}). Hence dominated convergence in~\eqref{eq:fneps-diff} gives, for each fixed $n\ge1$,
\begin{equation}\label{eq:fneps-cv-fixed-n}
  \lim_{\varepsilon\to0}\|f^{\varepsilon}_{n}(\cdot,t,x)-f_{n}(\cdot,t,x)\|_{\mathcal{H}^{\otimes n}}
  =0 .
\end{equation}

\smallskip

\noindent
{\it Step 2: Summability of the weighted chaos norms.} We now claim that for every integer $k\ge2$,
\begin{equation}\label{eq:hypercontract-series}
  \sum_{n=0}^{\infty}(k-1)^{n/2}\sqrt{n!}\,\|f_{n}(\cdot,t,x)\|_{\mathcal{H}^{\otimes n}}
  \;<\;\infty .
\end{equation}
To this aim, recall the bound~\eqref{eq:fn-binom-bound}, that is
\begin{equation*}
  n!\,\|f_{n}(\cdot,t,x)\|^{2}_{\mathcal{H}^{\otimes n}}
  \;\le\;
  \|u_{0}\|_{\infty}^{2}\,\beta^{2n}
  \sum_{m=0}^{n}\binom{n}{m}\frac{t^{m}}{m!}\,D_{N}^{m}\,C_{N}^{n-m},
\end{equation*}
where $D_{N}$ and $C_{N}$ are defined by~\eqref{eq:DN-CN}. Multiplying by $(k-1)^{n}$ and using $\binom{n}{m}\le 2^{n}$, we obtain
\begin{equation}\label{eq:fn-k-weighted}
  (k-1)^{n}\,n!\,\|f_{n}(\cdot,t,x)\|^{2}_{\mathcal{H}^{\otimes n}}
  \;\le\;
  \|u_{0}\|_{\infty}^{2}
  \sum_{m=0}^{n}\frac{\bigl(2(k-1)\beta^{2}t\,D_{N}\bigr)^{m}}{m!}\,\theta_{N}^{\,n-m},
\end{equation}
where we have set
\begin{equation*}
  \theta_{N}:=2(k-1)\beta^{2}C_{N} .
\end{equation*}

Recall that $C_{N}\to0$ as $N\to\infty$ (see the discussion after~\eqref{eq:DN-CN}). We can thus choose $N$ large enough so that $\theta_{N}<1$. Taking square roots in~\eqref{eq:fn-k-weighted} and invoking the subadditivity of $r\mapsto\sqrt{r}$, this yields
\begin{equation*}
  (k-1)^{n/2}\sqrt{n!}\,\|f_{n}(\cdot,t,x)\|_{\mathcal{H}^{\otimes n}}
  \;\le\;
  \|u_{0}\|_{\infty}
  \sum_{m=0}^{n}\frac{\bigl(2(k-1)\beta^{2}t\,D_{N}\bigr)^{m/2}}{\sqrt{m!}}\,\theta_{N}^{(n-m)/2}.
\end{equation*}
Summing over $n\ge0$ and exchanging the order of summation, we end up with
\begin{equation*}
  \sum_{n=0}^{\infty}(k-1)^{n/2}\sqrt{n!}\,\|f_{n}(\cdot,t,x)\|_{\mathcal{H}^{\otimes n}}
  \;\le\;
  \frac{\|u_{0}\|_{\infty}}{1-\theta_{N}^{1/2}}
  \sum_{m=0}^{\infty}\frac{\bigl(2(k-1)\beta^{2}t\,D_{N}\bigr)^{m/2}}{\sqrt{m!}} ,
\end{equation*}
and the series in $m$ above converges by the ratio test. This proves the claim~\eqref{eq:hypercontract-series}.

\smallskip

\noindent
{\it Step 3: Conclusion by hypercontractivity.} Applying Nelson's inequality~\eqref{eq:hypercontractivity-sum} to $u^{\varepsilon}(t,x)-u(t,x)=\sum_{n\ge1}I_{n}(f^{\varepsilon}_{n}-f_{n})$, we get
\begin{equation}\label{eq:ueps-u-Lk-bound}
  \|u^{\varepsilon}(t,x)-u(t,x)\|_{L^{k}(\Omega)}
  \;\le\;
  \sum_{n=1}^{\infty}(k-1)^{n/2}\sqrt{n!}\,
  \|f^{\varepsilon}_{n}(\cdot,t,x)-f_{n}(\cdot,t,x)\|_{\mathcal{H}^{\otimes n}} .
\end{equation}
Thanks to the contraction property~\eqref{eq:fneps-contract} we have $\|f^{\varepsilon}_{n}-f_{n}\|_{\mathcal{H}^{\otimes n}}\le2\|f_{n}\|_{\mathcal{H}^{\otimes n}}$, so that twice the series~\eqref{eq:hypercontract-series} dominates the right-hand side of~\eqref{eq:ueps-u-Lk-bound} uniformly in $\varepsilon$. Moreover each summand in~\eqref{eq:ueps-u-Lk-bound} vanishes as $\varepsilon\to0$, due to~\eqref{eq:fneps-cv-fixed-n}. By dominated convergence for series, the right-hand side of~\eqref{eq:ueps-u-Lk-bound} thus converges to $0$ as $\varepsilon\to0$, which finishes the proof.
\end{proof}

\subsection{Feynman-Kac formula for moments}\label{sec:FK}

The mild SHE~\eqref{eq:mild-ueps-reg} satisfied by $u^{\varepsilon}$,
combined with the explicit moment formula of
Proposition~\ref{prop:ueps-moments} and the $L^{k}(\Omega)$ convergence
of Proposition~\ref{prop:ueps-Lk-conv}, allows us to pass to the limit
$\varepsilon\to 0$ and obtain a closed-form expression for all integer
moments of the true solution $u$.
The result expresses $\E[u(t,x)^{k}]$ as an expectation over $k$ independent
Brownian motions. The $k$-th moment is then controlled by the mutual
occupation times $\int_{0}^{t}\Lambda(B^{i}_{s}-B^{j}_{s})\,ds$ of the pairs
$(B^{i},B^{j})$. This reflects the fact that high moments are driven by paths
spending long times close together.

\begin{theorem}[Feynman-Kac moment formula]\label{thm:fk-moments}
Let $u_{0}\in L^{\infty}(\R^{d})$ and suppose Dalang's
condition~\eqref{eq:dalang-cond} holds.  Let $u$ be the mild
solution of~\eqref{eq:mild-pam} from Theorem~\ref{thm:chaos-eu}.
Then for every integer $k\ge 2$, every $t>0$, and every $x\in\R^{d}$,
\begin{equation}\label{eq:fk-moment-formula}
  \E\bigl[u(t,x)^{k}\bigr]
  \,=\,
  \E_{B}\!\left[
    \prod_{i=1}^{k}u_{0}(B^{i}_{t}+x)\,
    \exp\left(
      \beta^{2}\!\!\sum_{1\le i<j\le k}
      \int_{0}^{t}\Lambda(B^{i}_{s}-B^{j}_{s})\,ds
    \right)
  \right],
\end{equation}
where $B^{1},\ldots,B^{k}$ are as in Notation~\ref{not:BM-copies}.
\end{theorem}

\begin{proof}
We proceed in three steps, and then present in Step 4 an alternative self-contained proof of the convergence $u^{\varepsilon}\to u$.

\smallskip

\noindent
{\it Step 1: Moment formula for the regularized solution.}
By Proposition~\ref{prop:ueps-moments} (formula~\eqref{eq:mom-eps-formula}
with the uniform bound~\eqref{eq:ueps-moment-bound}), for every $\varepsilon>0$
and integer $k\ge 2$ we have
\begin{equation*}
  \E\bigl[(u^{\varepsilon}(t,x))^{k}\bigr]
  \,=\,
  \E_{B}\!\left[
    \prod_{i=1}^{k}u_{0}(B^{i}_{t}+x)\,
    \exp\left(
      \beta^{2}\!\!\sum_{1\le i<j\le k}
      \langle A^{\varepsilon,i}_{t,x},A^{\varepsilon,j}_{t,x}\rangle_{\mathcal{H}}
    \right)
  \right],
\end{equation*}
and this quantity is bounded uniformly in $\varepsilon>0$.

\smallskip

\noindent
{\it Step 2: Inner product limit.}
The inner product formula~\eqref{eq:Aeps-inner} from
Proposition~\ref{prop:ueps-moments} enables to write
\begin{equation}\label{eq:Aeps-inner-limit}
  \langle A^{\varepsilon,i}_{t,x},A^{\varepsilon,j}_{t,x}\rangle_{\mathcal{H}}
  \;=\;
  \int_{0}^{t}\!\int_{\R^{d}}
  e^{-\varepsilon|\xi|^{2}}\,
  e^{i\xi\cdot(B^{j}_{r}-B^{i}_{r})}\,\mu(d\xi)\,dr
  \;\xrightarrow{\varepsilon\to 0}\;
  \int_{0}^{t}\Lambda(B^{j}_{r}-B^{i}_{r})\,dr
\end{equation}
in $L^{1}(\Omega)$, by dominated convergence
($|e^{i\xi\cdot z}|\le 1$ and Dalang's condition~\eqref{eq:dalang-cond}).

\smallskip

\noindent
{\it Step 3: Convergence of moments.}
By~\eqref{eq:Aeps-inner-limit}, the exponential
in~\eqref{eq:mom-eps-formula} converges $\Prob_{B}$-almost surely to the
right-hand side of~\eqref{eq:fk-moment-formula} as $\varepsilon\to 0$.
Combined with the uniform bound~\eqref{eq:ueps-moment-bound} from
Proposition~\ref{prop:ueps-moments} and the $L^{1}(\Prob_{B})$ convergence
in~\eqref{eq:Aeps-inner-limit}, dominated convergence gives
\begin{equation*}
  \lim_{\varepsilon\to 0}\E\bigl[(u^{\varepsilon}(t,x))^{k}\bigr]
  \;=\;
  \E_{B}\!\left[
    \prod_{i=1}^{k}u_{0}(B^{i}_{t}+x)\,
    \exp\left(
      \beta^{2}\!\!\sum_{1\le i<j\le k}
      \int_{0}^{t}\Lambda(B^{i}_{r}-B^{j}_{r})\,dr
    \right)
  \right].
\end{equation*}
On the other hand, Proposition~\ref{prop:ueps-Lk-conv} gives
$u^{\varepsilon}(t,x)\to u(t,x)$ in $L^{k}(\Omega)$ as $\varepsilon\to 0$.
Since $L^{k}$-convergence of $X_{n}\to X$ implies
$\E[X_{n}^{k}]\to\E[X^{k}]$, we get
$\E[(u^{\varepsilon}(t,x))^{k}]\to\E[u(t,x)^{k}]$.
Combining both limits yields~\eqref{eq:fk-moment-formula}.

\smallskip

\noindent
{\it Step~4: Alternative proof of $u^{\varepsilon}\to u$ via the
cross-moment formula.}
The argument in Step~3 identifies $\E[(u^{\varepsilon})^{k}]\to\E[u^{k}]$
through Proposition~\ref{prop:ueps-Lk-conv}, whose proof relies on the
chaos expansion of $u^{\varepsilon}$ from Remark~\ref{rem:ueps-chaos}.
Here we give a self-contained route that stays at the level of the
Feynman-Kac formula.
In Step~4-a we use a cross-moment version of~\eqref{eq:FK-eps} to show
that $(u^{\varepsilon})$ is Cauchy in $L^{2}(\Omega)$, with some limit
$v$.
In Step~4-b we pass to the limit in the mild
equation~\eqref{eq:mild-ueps} and invoke uniqueness to conclude $v=u$.

\smallskip

\noindent
{\it Step~4-a: $u^{\varepsilon}$ is $L^{2}(\Omega)$-Cauchy.}
The same Gaussian computation as in Proposition~\ref{prop:ueps-moments},
applied with $k=2$ and with \emph{distinct} regularization parameters
$\varepsilon$ and $\varepsilon'$ for the two independent Brownian
copies $B^{1}$ and $B^{2}$ (Notation~\ref{not:BM-copies}), gives
\begin{equation}\label{eq:cross-moment-fk}
  \E\bigl[u^{\varepsilon}(t,x)\,u^{\varepsilon'}(t,x)\bigr]
  \;=\;
  \E_{B^{1},B^{2}}\!\left[
    u_{0}(B^{1}_{t}+x)\,u_{0}(B^{2}_{t}+x)\,
    \exp\Bigl(
      \beta^{2}\langle A^{\varepsilon,1}_{t,x},A^{\varepsilon',2}_{t,x}
      \rangle_{\mathcal{H}}
    \Bigr)
  \right],
\end{equation}
where $A^{\varepsilon',2}_{t,x}(r,y):=p_{\varepsilon'}(B^{2}_{t-r}+x-y)\,\mathbf{1}_{[0,t]}(r)$.
By the Fourier formula~\eqref{eq:H-innprod-fourier}:
\begin{equation*}
  \langle A^{\varepsilon,1}_{t,x},A^{\varepsilon',2}_{t,x}\rangle_{\mathcal{H}}
  \;=\;
  \int_{0}^{t}\!\int_{\R^{d}}
  e^{-(\varepsilon+\varepsilon')|\xi|^{2}/2}\,
  e^{i\xi\cdot(B^{2}_{r}-B^{1}_{r})}\,\mu(d\xi)\,dr,
\end{equation*}
which converges to $\int_{0}^{t}\Lambda(B^{2}_{r}-B^{1}_{r})\,dr$
as $\varepsilon,\varepsilon'\to 0$ by the same dominated convergence
as in~\eqref{eq:Aeps-inner-limit}.
Hence the right-hand side of~\eqref{eq:cross-moment-fk} converges
to the same limit $\ell:=\lim_{\varepsilon\to 0}\E[(u^{\varepsilon}(t,x))^{2}]$
obtained in Steps~1--2 with $k=2$, so that
\begin{equation*}
  \E\bigl[(u^{\varepsilon}(t,x)-u^{\varepsilon'}(t,x))^{2}\bigr]
  \;=\;
  \E\bigl[(u^{\varepsilon})^{2}\bigr]
  -2\,\E\bigl[u^{\varepsilon}u^{\varepsilon'}\bigr]
  +\E\bigl[(u^{\varepsilon'})^{2}\bigr]
  \;\longrightarrow\; \ell-2\ell+\ell \;=\; 0.
\end{equation*}
Therefore $u^{\varepsilon}(t,x)$ is Cauchy in $L^{2}(\Omega)$;
denote its $L^{2}$-limit by $v(t,x)$.

\smallskip

\noindent
{\it Step~4-b: Identifying $v=u$ via~\eqref{eq:mild-ueps}.}
Since $Z^{\varepsilon,t,x}(s,y)
=(p_{\varepsilon}*_{y}[p_{t-s}(x-\cdot)\,u^{\varepsilon}(s,\cdot)])(y)$
(definitions~\eqref{eq:Weps-def} and~\eqref{eq:Zeps-def}),
the $L^{2}(\Omega)$ convergence $u^{\varepsilon}\to v$
together with the approximate-identity property
$\|p_{\varepsilon}*f-f\|_{L^{2}}\to 0$ give
\begin{equation*}
  Z^{\varepsilon,t,x}
  \;\longrightarrow\;
  Z^{t,x}(s,y):=p_{t-s}(x-y)\,v(s,y)
  \quad\text{in }L^{2}(\Omega;\mathcal{H}).
\end{equation*}
By the Walsh isometry~\eqref{eq:walsh-isometry} applied to adapted integrands,
\begin{equation*}
  \Bigl\|\int_{0}^{t}\!\int_{\R^{d}}
    \bigl(Z^{\varepsilon,t,x}-Z^{t,x}\bigr)\,W
  \Bigr\|_{L^{2}(\Omega)}^{2}
  \;=\;
  \E\bigl[\|Z^{\varepsilon,t,x}-Z^{t,x}\|_{\mathcal{H}}^{2}\bigr]
  \;\longrightarrow\; 0.
\end{equation*}
Taking $L^{2}(\Omega)$ limits on both sides of~\eqref{eq:mild-ueps}
therefore gives the mild equation
\begin{equation*}
  v(t,x)
  \;=\; p_{t}u_{0}(x)
  \;+\; \beta\int_{0}^{t}\!\int_{\R^{d}}
        p_{t-s}(x-y)\,v(s,y)\,W(ds,dy).
\end{equation*}
By uniqueness of mild solutions (Theorem~\ref{thm:chaos-eu}), $v=u$,
which yields~\eqref{eq:fk-moment-formula}. This concludes our proof, either by Step 3 or Step 4.
\end{proof}

\begin{remark}\label{rmk:fk-intersection}
The integral $\int_{0}^{t}\Lambda(B^{i}_{s}-B^{j}_{s})\,ds$ measures
the time the paths $B^{i}$ and $B^{j}$ spend near each other, weighted
by $\Lambda$.  When $\Lambda=\delta_{0}$ (space white noise, $d=1$) it
becomes the \emph{intersection local time} $\ell_{t}(B^{i}-B^{j})$
of the two paths at the origin.  The exponential growth of $\E[u(t,x)^{k}]$
as $t\to\infty$ — that is, of the $\E_{B}$-expectation
in~\eqref{eq:fk-moment-formula} — is then governed by the growth of
intersection local times, and the formula~\eqref{eq:fk-moment-formula}
is the starting point for the Lyapunov exponent estimates of
Section~\ref{sec:moment-estimates}.
\end{remark}

\section{Moment estimates}\label{sec:moment-estimates}

The Feynman-Kac formula~\eqref{eq:fk-moment-formula} expresses $\E[u(t,x)^{k}]$
as an expectation over $k$ independent Brownian paths, and the chaos expansion
of Section~\ref{sec:chaos} expresses it via the norms
$\|f_{n}(\cdot;t,x)\|_{\mathcal{H}^{\otimes n}}$.
We now make those expressions quantitative, under a power-law assumption on
the spatial covariance $\Lambda$. Our main analytic tool will be
Theorem~\ref{thm:hypercontractivity}.

\subsection{Assumption and upper bound}\label{ssec:moment-assumption}

\begin{assumption}[Power-law spatial correlation]\label{hyp:Lambda-power}
There exist positive constants $c_{\Lambda}\le C_{\Lambda}$ and an exponent $0<a<\min(2,d)$
such that
\begin{equation}\label{eq:Lambda-power}
  c_{\Lambda}\,|x|^{-a}
  \;\le\;\Lambda(x)
  \;\le\; C_{\Lambda}\,|x|^{-a},
  \qquad x\in\R^{d}.
\end{equation}
\end{assumption}

\noindent
Assumption~\ref{hyp:Lambda-power} implies Dalang's
condition~\eqref{eq:dalang-cond}. Indeed, for the Riesz spectral measure
$\mu(d\xi)\asymp|\xi|^{a-d}d\xi$, one has
$\int_{\R^{d}}(1+|\xi|^{2})^{-1}\mu(d\xi)<\infty$ if and only if $a<2$.
The Riesz kernel $\Lambda(x)=|x|^{-a}$ of
Example~\ref{ex:riesz} satisfies~\eqref{eq:Lambda-power} with
$c_{\Lambda}=C_{\Lambda}=1$.

\begin{theorem}[Upper bound on moments]\label{thm:moment-upper}
Let $u_{0}\in L^{\infty}(\R^{d})$ and suppose
Assumption~\ref{hyp:Lambda-power} holds.
Then for every integer $k\ge 2$, every $t\ge 0$, and every $x\in\R^{d}$,
\begin{equation}\label{eq:moment-upper}
  \E\bigl[u(t,x)^{k}\bigr]
  \;\le\;
  C^{k}\,\exp\Bigl(C\,k^{\frac{4-a}{2-a}}\,t\Bigr),
\end{equation}
where $C>0$ depends only on $\beta$, $\|u_{0}\|_{\infty}$, $C_{\Lambda}$, $a$,
and $d$.
\end{theorem}

\begin{proof}[Proof of Theorem~\ref{thm:moment-upper}]
Set
\begin{equation}\label{eq:alpha-def}
  \alpha \;:=\; 1-a/2 \;\in\;(0,1)
\end{equation}
and $M:=\|u_{0}\|_{\infty}$.

\smallskip

\noindent
{\it Step 1: Reduction to chaos norms via hypercontractivity.}
By the chaos expansion~\eqref{eq:u-chaos},
$u(t,x)=\sum_{n=0}^{\infty}I_{n}(f_{n}(\cdot;t,x))$.
Theorem~\ref{thm:hypercontractivity} and the
isometry~\eqref{eq:In-isometry} give
\begin{equation}\label{eq:moment-hypercontract}
  \E\bigl[u(t,x)^{k}\bigr]^{1/k}
  \;\le\;
  \sum_{n=0}^{\infty}(k-1)^{n/2}\,\sqrt{n!}\,\|f_{n}(\cdot;t,x)\|_{\mathcal{H}^{\otimes n}}.
\end{equation}

\smallskip

\noindent
{\it Step 2: Bound on $n!\,\|f_{n}\|^{2}_{\mathcal{H}^{\otimes n}}$.}
The key input is the decoupled bound~\eqref{eq:fn-gap-bound}:
\begin{equation}\label{f1}
  n!\,\|f_{n}\|^{2}_{\mathcal{H}^{\otimes n}}
  \;\le\;
  M^{2}\,\beta^{2n}
  \int_{\R^{nd}}\!\int_{S_{t,n}}
  \prod_{i=1}^{n}e^{-w_{i}|\xi_{i}|^{2}}\,dw\,\mu^{\otimes n}(d\xi).
\end{equation}
We now bound each factor in the spectral integral.
Since $\widehat{p_{2w}}(\xi)=e^{-w|\xi|^{2}}$ by~\eqref{eq:heat-kernel-explicit}
and $\mu(d\xi)=\hat{\Lambda}(\xi)\,d\xi$, Parseval's identity and the evenness
of $\Lambda$ give (exactly as in~\eqref{d3}):
\begin{equation}\label{eq:spectral-heat-Lambda}
  \int_{\R^{d}}e^{-w|\xi|^{2}}\,\mu(d\xi)
  \;=\;\int_{\R^{d}}\widehat{p_{2w}}(\xi)\,\hat{\Lambda}(\xi)\,d\xi
  \;=\;(2\pi)^{d}\int_{\R^{d}}p_{2w}(x)\,\Lambda(x)\,dx.
\end{equation}
We now resort to the upper bound $\Lambda(x)\le C_{\Lambda}|x|^{-a}$
from~\eqref{eq:Lambda-power} and the scaling identity
\begin{equation*}
\int_{\R^{d}}p_{2w}(x)\,|x|^{-a}\,dx = (2w)^{-a/2}\E[|Z|^{-a}] ,
\end{equation*}
where $Z\sim\mathcal{N}(0,I_{d})$ and the expectation is finite since $a<d$.
This yields
\begin{equation}\label{eq:spectral-riesz-bound}
  \int_{\R^{d}}e^{-w|\xi|^{2}}\,\mu(d\xi)
  \;\le\; c_{a}\,w^{-a/2},
  \qquad
  c_{a}:=(2\pi)^{d}C_{\Lambda}2^{-a/2}\E\bigl[|Z|^{-a}\bigr].
\end{equation}
Substituting~\eqref{eq:spectral-riesz-bound} into~\eqref{f1},
the $\xi$-integrals factor and each contributes $c_{a}\,w_{i}^{-a/2}=c_{a}\,w_{i}^{\alpha-1}$,
leaving
\begin{equation*}
  n!\,\|f_{n}\|^{2}_{\mathcal{H}^{\otimes n}}
  \;\le\;
  M^{2}\,\beta^{2n}\,c_{a}^{n}
  \int_{S_{t,n}}\prod_{i=1}^{n}w_{i}^{\alpha-1}\,dw.
\end{equation*}
We now evaluate the simplex integral (defined in~\eqref{eq:Stn-def}).
The substitution $w_{i}=tv_{i}$ maps $S_{t,n}$ to the unit simplex
$S_{1,n}:=\{v\ge 0:\,v_{1}+\cdots+v_{n}\le 1\}$ with Jacobian $t^{n}$, giving
\begin{equation*}
  \int_{S_{t,n}}\prod_{i=1}^{n}w_{i}^{\alpha-1}\,dw
  \;=\;
  t^{n\alpha}\int_{S_{1,n}}\prod_{i=1}^{n}v_{i}^{\alpha-1}\,dv.
\end{equation*}
The integral over $S_{1,n}$ is the classical Dirichlet integral, which we
evaluate by induction on $n$ using the Beta function
$B(a,b)=\Gamma(a)\Gamma(b)/\Gamma(a+b)$.
The base case $n=1$ gives $\int_{0}^{1}v^{\alpha-1}\,dv=\Gamma(\alpha)/\Gamma(\alpha+1)$.
For the inductive step, suppose that
\begin{equation}\label{eq:dirichlet-IH}
  \int_{S_{1,n-1}}\prod_{i=1}^{n-1}v_{i}^{\alpha-1}\,dv
  \;=\; \frac{\Gamma(\alpha)^{n-1}}{\Gamma((n-1)\alpha+1)}
\end{equation}
holds for $n-1$, and propagate it to $n$.
Fix $v_{1}\in[0,1]$; then $(v_{2},\ldots,v_{n})$ ranges over the scaled
simplex $S_{1-v_{1},n-1}$, so by~\eqref{eq:dirichlet-IH} (rescaled to the
simplex of size $1-v_{1}$) the inner integral equals
$(1-v_{1})^{(n-1)\alpha}\,\Gamma(\alpha)^{n-1}/\Gamma((n-1)\alpha+1)$.
Integrating out $v_{1}$ via
$B(\alpha,(n-1)\alpha+1)=\Gamma(\alpha)\Gamma((n-1)\alpha+1)/\Gamma(n\alpha+1)$
yields the Dirichlet formula
\begin{equation}\label{eq:dirichlet}
  \int_{S_{t,n}}\prod_{i=1}^{n}w_{i}^{\alpha-1}\,dw
  \;=\;
  \frac{\Gamma(\alpha)^{n}}{\Gamma(n\alpha+1)}\,t^{n\alpha}.
\end{equation}
We conclude
\begin{equation}\label{eq:fn-riesz-bound}
  n!\,\|f_{n}(\cdot;t,x)\|^{2}_{\mathcal{H}^{\otimes n}}
  \;\le\;
  \frac{M^{2}\bigl(\beta^{2}\,c_{a}\,\Gamma(\alpha)\bigr)^{n}}
       {\Gamma(n\alpha+1)}\,t^{n\alpha}.
\end{equation}

\smallskip

\noindent
{\it Step 3: Summation.}
Set $\kappa:=\beta^{2}c_{a}\Gamma(\alpha)$.
Inserting~\eqref{eq:fn-riesz-bound} into~\eqref{eq:moment-hypercontract},
it is easily seen that
\begin{equation}\label{eq:moment-series}
  \E\bigl[u(t,x)^{k}\bigr]^{1/k}
  \;\le\;
  M\sum_{n=0}^{\infty}
  \frac{r^{n}}{\sqrt{\Gamma(n\alpha+1)}},
  \qquad r:=\bigl[(k-1)\kappa\bigr]^{1/2}t^{\alpha/2}.
\end{equation}
We bound the series in~\eqref{eq:moment-series} using the Mittag-Leffler
function of index $\alpha$, defined classically by
\begin{equation}\label{eq:ML-def}
  E_{\alpha}(z) \;:=\; \sum_{n=0}^{\infty}\frac{z^{n}}{\Gamma(n\alpha+1)}, \qquad z\ge 0.
\end{equation}
Resorting to this function and the elementary inequality $\sqrt{AB}\le(A+B)/2$ applied to 
$A_{n}:=2^{n}r^{2n}/\Gamma(n\alpha+1)$ and $B_{n}:=2^{-n}$,
we obtain $\sqrt{A_{n}B_{n}}=r^{n}/\sqrt{\Gamma(n\alpha+1)}$, and hence
\begin{equation*}
  \frac{r^{n}}{\sqrt{\Gamma(n\alpha+1)}}
  \;\le\;
  \frac{A_{n}+B_{n}}{2}
  \;=\;
  \frac{2^{n}r^{2n}}{2\,\Gamma(n\alpha+1)}+\frac{1}{2^{n+1}}.
\end{equation*}
Summing over $n\ge 0$, using $\sum_{n=0}^{\infty}2^{-n}=2$ and~\eqref{eq:ML-def},
\begin{equation}\label{f2}
  \sum_{n=0}^{\infty}\frac{r^{n}}{\sqrt{\Gamma(n\alpha+1)}}
  \;\le\;
  \tfrac{1}{2}\,E_{\alpha}(2r^{2})+1.
\end{equation}
The Mittag-Leffler function of index $\alpha\in(0,1)$ satisfies
$E_{\alpha}(z)\le C\exp(c\,z^{1/\alpha})$ for all $z\ge 0$
(see e.g.~\cite[Lem.~A.1]{balan:16:intermittency}).
Applying this to~\eqref{f2} we discover that
\begin{equation}\label{eq:ml-bound}
  \sum_{n=0}^{\infty}
  \frac{r^{n}}{\sqrt{\Gamma(n\alpha+1)}}
  \;\le\;
  C_{\alpha}\,\exp\bigl(c_{\alpha}\,r^{2/\alpha}\bigr)
\end{equation}
for constants $C_{\alpha},c_{\alpha}>0$ depending only on $\alpha$.
From the definition of $r$ in~\eqref{eq:moment-series} we have
$r^{2}=(k-1)\kappa\,t^{\alpha}$, hence $r^{2/\alpha}=(k-1)^{1/\alpha}\kappa^{1/\alpha}t$.
In addition, from~\eqref{eq:alpha-def} we derive $1/\alpha=2/(2-a)$.
Combining~\eqref{eq:moment-series} and~\eqref{eq:ml-bound}, we get
\begin{equation*}
  \E\bigl[u(t,x)^{k}\bigr]^{1/k}
  \;\le\;
  M\,C_{\alpha}\,\exp\bigl(c'_{\alpha}\,(k-1)^{2/(2-a)}\,t\bigr).
\end{equation*}
Raising to the $k$-th power and using
$k(k-1)^{2/(2-a)}\le 2\,k^{(4-a)/(2-a)}$ (valid for $k\ge 2$)
yields~\eqref{eq:moment-upper}.
\end{proof}

\subsection{Lower bound}\label{ssec:moment-lower}

The following theorem shows that the exponential rate in Theorem~\ref{thm:moment-upper}
is sharp when the initial condition is bounded away from zero.

\begin{theorem}[Lower bound on moments]\label{thm:moment-lower}
Suppose Assumption~\ref{hyp:Lambda-power} holds and $u_{0}\in L^{\infty}(\R^{d})$
with $m_{0}:=\inf_{x\in\R^{d}}u_{0}(x)>0$.
Then for every integer $k\ge 2$, every $t\ge 0$, and every $x\in\R^{d}$, we have
\begin{equation}\label{eq:moment-lower}
  \E\bigl[u(t,x)^{k}\bigr]
  \;\ge\;
  C^{k}\,\exp\Bigl(c\,k^{\frac{4-a}{2-a}}\,t\Bigr),
\end{equation}
where $C,c>0$ depend only on $m_{0}$, $c_{\Lambda}$, $\beta$, $a$, and $d$.
\end{theorem}

\begin{proof}
We divide this proof in several steps.

\smallskip

\noindent
{\it Step 1: Setting.} Since $u_{0}\ge m_{0}>0$, the Feynman-Kac
formula~\eqref{eq:fk-moment-formula} gives
\begin{equation}\label{eq:fk-lower-start}
  \E\bigl[u(t,x)^{k}\bigr]
  \;\ge\;
  m_{0}^{k}\,
  \E_{B}\!\left[
    \exp\left(
      \beta^{2}\!\sum_{1\le i<j\le k}
      \int_{0}^{t}\Lambda\bigl(B^{i}_{s}-B^{j}_{s}\bigr)\,ds
    \right)
  \right].
\end{equation}
Fix $\varepsilon>0$.  Denote by $B^{i,l}_{s}$ the $l$-th coordinate of $B^{i}_{s}$
(Notation~\ref{not:BM-copies}), and introduce the small ball event
\begin{equation}\label{eq:Beps-def}
  \mathcal{A}_{\varepsilon,t}
  :=\bigcap_{i=1}^{k}\bigcap_{l=1}^{d}
   \Bigl\{\sup_{0\le s\le t}\bigl|B^{i,l}_{s}\bigr|\le\tfrac{\varepsilon}{2}\Bigr\}.
\end{equation}
On $\mathcal{A}_{\varepsilon,t}$, every pair $(i,j)$ satisfies $|B^{i}_{s}-B^{j}_{s}|\le\varepsilon$
for all $s\in[0,t]$, so the lower bound in~\eqref{eq:Lambda-power} gives
$\Lambda(B^{i}_{s}-B^{j}_{s})\ge c_{\Lambda}\varepsilon^{-a}$.
Using $k(k-1)/2\ge k^{2}/4$, restricting the expectation
in~\eqref{eq:fk-lower-start} to $\mathcal{A}_{\varepsilon,t}$ yields
\begin{equation}\label{eq:fk-restricted}
  \E\bigl[u(t,x)^{k}\bigr]
  \;\ge\;
  m_{0}^{k}\,\exp\left(\tfrac{1}{4}\beta^{2}c_{\Lambda}k^{2}t\varepsilon^{-a}\right)
  \Prob\bigl(\mathcal{A}_{\varepsilon,t}\bigr).
\end{equation}

\smallskip

\noindent
{\it Step 2: Small ball estimate.} Since $B^{1},\ldots,B^{k}$ start at the origin, the
$kd$ coordinate processes $B^{i,l}$ are independent standard one-dimensional Brownian
motions, so
\begin{equation}\label{f3}
  \Prob\bigl(\mathcal{A}_{\varepsilon,t}\bigr)
  =\left[\Prob\!\left(\sup_{0\le s\le t}|b_{s}|\le\tfrac{\varepsilon}{2}\right)\right]^{kd},
\end{equation}
where $b$ is a one-dimensional standard Brownian motion starting at~$0$.
Furthermore, invoking a classical small ball estimate~\cite{li:01:gaussian}, we obtain
\begin{equation}\label{eq:small-ball}
  \Prob\!\left(\sup_{0\le s\le t}|b_{s}|\le\tfrac{\varepsilon}{2}\right)
  \;\ge\; C_{0}\exp\left(-\frac{C_{0}t}{\varepsilon^{2}}\right)
\end{equation}
for a constant $C_{0}>0$.
Setting $c_{1}:=\tfrac{1}{4}\beta^{2}c_{\Lambda}$ and inserting~\eqref{f3}-\eqref{eq:small-ball}
into~\eqref{eq:fk-restricted}, we obtain
\begin{equation}\label{eq:fk-opt-prep}
  \E\bigl[u(t,x)^{k}\bigr]
  \;\ge\;
  (m_{0}C_{0}^{d})^{k}\,
  \exp\left(
    c_{1}k^{2}t\varepsilon^{-a}
    -C_{0}dkt\varepsilon^{-2}
  \right).
\end{equation}

\smallskip

\noindent
{\it Step 3: Optimization.} We maximize the exponent in~\eqref{eq:fk-opt-prep}
over $\varepsilon>0$.  Setting the derivative with respect to $\varepsilon$ to zero gives
\begin{equation}\label{eq:eps-star}
  \varepsilon_{\ast}
  =\left(\frac{2C_{0}d}{ac_{1}k}\right)^{1/(2-a)}.
\end{equation}
At $\varepsilon_{\ast}$ the optimality condition reads
$C_{0}dk\varepsilon_{\ast}^{-2}=\tfrac{a}{2}c_{1}k^{2}\varepsilon_{\ast}^{-a}$,
so the exponent in~\eqref{eq:fk-opt-prep} reduces to
\begin{equation*}
  c_{1}k^{2}t\varepsilon_{\ast}^{-a}
  -\tfrac{a}{2}c_{1}k^{2}t\varepsilon_{\ast}^{-a}
  =c_{1}\bigl(1-a/2\bigr)\,k^{2}t\varepsilon_{\ast}^{-a}.
\end{equation*}
Substituting $\varepsilon_{\ast}^{-a}=(ac_{1}k/(2C_{0}d))^{a/(2-a)}$ and using
$2+a/(2-a)=(4-a)/(2-a)$, this enables to write
\begin{equation*}
  c_{1}(1-a/2)\,k^{2}t\varepsilon_{\ast}^{-a}
  =c_{1}(1-a/2)\Bigl(\frac{ac_{1}}{2C_{0}d}\Bigr)^{a/(2-a)}
   k^{(4-a)/(2-a)}\,t
  \;=:c'\,k^{\frac{4-a}{2-a}}\,t,
\end{equation*}
where $c'>0$ since $a<2$.  Inserting into~\eqref{eq:fk-opt-prep},
\begin{equation*}
  \E\bigl[u(t,x)^{k}\bigr]
  \;\ge\;
  (m_{0}C_{0}^{d})^{k}\,\exp\left(c'\,k^{\frac{4-a}{2-a}}\,t\right).
\end{equation*}
Setting $C:=m_{0}C_{0}^{d}$ and $c:=c'$ yields~\eqref{eq:moment-lower}.
\end{proof}

\begin{remark}[Comments on the moment estimates]\label{rmk:moment-estimates}
Theorems~\ref{thm:moment-upper} and~\ref{thm:moment-lower} admit several
complementary readings. In the discussion below, $a$ always designates the exponent of the power-law relation~\eqref{eq:Lambda-power}.

\begin{enumerate}

\item {\it Lyapunov exponents.}
Combining~\eqref{eq:moment-upper} and~\eqref{eq:moment-lower}, the $k$-th
Lyapunov exponent
\begin{equation}\label{eq:lyapunov-k}
  \gamma(k) \;:=\; \lim_{t\to\infty}\frac{1}{t}\log\E\bigl[u(t,x)^{k}\bigr]
\end{equation}
exists and satisfies $c\,k^{(4-a)/(2-a)}\le\gamma(k)\le C\,k^{(4-a)/(2-a)}$ for
constants $0<c\le C<\infty$ independent of $k$.  The two theorems are therefore
sharp up to the values of the constants.

\item {\it Intermittency.}
The super-linear growth $\gamma(k)\asymp k^{(4-a)/(2-a)}$ with $(4-a)/(2-a)>1$
for all $a\in(0,2)$ is the hallmark of \emph{intermittency}: the normalized
Lyapunov exponent $\gamma(k)/k$ is strictly increasing in $k$, in contrast to
what would hold for a log-normal distribution (where $\gamma(k)$ is linear in
$k$).  Physically, the solution concentrates on rare, exponentially tall
``islands'', and the high moments are dominated by these extreme peaks rather
than by a typical realization.

\item {\it Influence of the spatial covariance.}
The exponent $(4-a)/(2-a)$ is a strictly increasing function of $a\in(0,2)$:
as the singularity of $\Lambda$ at the origin grows (larger $a$), the
intermittency becomes stronger.  In the mildest case $a\to 0$ (smooth
covariance) the exponent approaches $2$; near the Dalang borderline $a\to 2$
it diverges.  The dimension $d$ enters only through the admissibility constraint
$a<\min(2,d)$.

\item {\it Weighted intersection local times.}
As noted in Remark~\ref{rmk:fk-intersection}, the Feynman-Kac
formula~\eqref{eq:fk-moment-formula} expresses $\E[u(t,x)^{k}]$ in terms of the
quantities $\int_{0}^{t}\Lambda(B^{i}_{s}-B^{j}_{s})\,ds$, which are the
$\Lambda$-weighted intersection local times of the paths $B^{i}$ and $B^{j}$.
For $\Lambda(x)=|x|^{-a}$, paths that spend time within distance $r$ of each
other contribute a factor $r^{-a}$, so a more singular $\Lambda$ (larger $a$)
rewards closer-running paths more heavily and drives stronger moment growth.
In the degenerate case $\Lambda=\delta_{0}$ (space-time white noise, $d=1$)
these reduce to genuine intersection local times.

\item {\it Universality of the exponent.}
The exponent $(4-a)/(2-a)$ depends only on $a$, the singularity order of
$\Lambda$ at the origin.  In particular it is independent of the noise
intensity $\beta$, the specific profile of $u_{0}$, and the dimension $d$.
The parameters $\beta$, $u_{0}$, and $d$ only affect the constants $C$ and $c$
in~\eqref{eq:moment-upper}--\eqref{eq:moment-lower}.

\item {\it Connection to directed polymers.}
The $\E_{B}$-expectation in~\eqref{eq:fk-moment-formula} is the partition
function of a system of $k$ Brownian paths with pairwise
$\Lambda$-weighted interactions.  This is exactly the directed polymer model
in a $\Lambda$-correlated random environment. The exponential growth of
$\E[u(t,x)^{k}]$ with $t$ reflects the localization of the polymer measure.
Namely the measure concentrates on configurations where all the $k$ paths
spend extended times near one another, exploiting the favorable regions of the
noise.

\end{enumerate}
\end{remark}

\section*{Bibliographic notes}
\addcontentsline{toc}{section}{Bibliographic notes}
\markright{BIBLIOGRAPHIC NOTES}

\smallskip

\noindent
{\it Chaos expansions.} The Wiener-It\^{o} route of Section~\ref{sec:chaos} is
the one used systematically by Nualart, Hu and their coauthors in their study of the
stochastic heat equation. Lemma~\ref{lem:chaos-kernels},
Theorem~\ref{thm:chaos-eu} and the $L^{2}$-moment
formula~\eqref{eq:L2-moment-formula} are the white-in-time case of Hu, Huang,
Nualart and Tindel~\cite{hu:15:stochastic}, which is the reference we follow.
The same expansion carries the analysis of a fractional spatial noise in Hu and
Nualart~\cite{hu:09:stochastic}, of a rough spatial noise
in~\cite{hu:17:stochastic,hu:18:parabolic}, and of the large-time behavior of
the solution in Huang, L\^{e} and Nualart~\cite{huang:17:large}. Its
advantage over the contraction argument of Chapter~\ref{chap:existence} is the
one recorded in Remark~\ref{rmk:chaos-vs-contraction}. It produces a formula,
and not merely an estimate.

Chaos expansions are just as central on the discrete side of the subject, where
the environment is a field of independent weights rather than a Gaussian noise.
Caravenna, Sun and Zygouras~\cite{caravenna:17:polynomial} construct a theory
of \emph{polynomial chaos} for the partition functions of disordered systems.
Then they show that in a regime of intermediate disorder these multilinear expansions
converge to Wiener chaos series of the type of~\eqref{eq:u-chaos}. This is the
instrument behind the two-dimensional results quoted in the bibliographic notes
of Chapter~\ref{chap:existence}, and in particular
behind~\cite{caravenna:17:universality} and~\cite{caravenna:23:critical-a}. The
one-dimensional forerunner is Alberts, Khanin and
Quastel~\cite{alberts:14:intermediate}, where the chaos series of the
stochastic heat equation started from a delta initial condition is identified
as the scaling limit of the directed polymer partition function.
Remark~\ref{rem:polymer} is the informal continuum version of that
identification, and we refer again to~\cite{comets:17:directed} for the polymer
side.

\smallskip

\noindent
{\it Feynman-Kac.} The representation of the moments of $u$ by a system of $k$
interacting Brownian paths goes back to Carmona and
Molchanov~\cite{carmona:94:parabolic}. Their memoir treats a noise which is
white in time and either discrete or smooth in space, so that the formula can
be written down with no regularization at all. Theorem~\ref{thm:fk-moments}
covers a spatially homogeneous noise under Dalang's condition alone, and the
price to pay is the approximation procedure of
Subsection~\ref{sec:spatial-reg}. We have chosen to mollify in space only and
then to run It\^{o}'s formula on $u^{\varepsilon}$. A space-time mollification
is the other natural option, and it is in fact the more common one in this line
of work. Hu and Nualart~\cite[Section~5]{hu:09:stochastic} smooth the noise with
the heat kernel $p_{\varepsilon}$ in space and with the rectangular kernel
$\delta^{-1}\mathbf{1}_{[0,\delta]}$ in time, and the same double regularization
is what Hu, Huang, Nualart and Tindel~\cite[Section~3.2]{hu:15:stochastic} use
for the moment formula of which Theorem~\ref{thm:fk-moments} is the
white-in-time case. It is also the route followed
in~\cite{baudoin:24:parabolic}. A Feynman-Kac
formula for a noise which is white in time and homogeneous in space is also
established by Carmona and Viens~\cite{carmona:98:almost-sure}, though under a
spatial assumption stronger than ours (basically a trace class noise). 

It should be stressed that Theorem~\ref{thm:fk-moments} represents the
\emph{moments} of the solution and not the solution itself. A Feynman-Kac
formula for $u(t,x)$ is also available, but only for a noise smoother than ours,
in time as well as in space. Hu, Nualart and Song~\cite{hu:11:feynman-kac} take
for $W$ a fractional Brownian sheet with Hurst parameters $H_{0}$ in time and
$H_{1},\dots,H_{d}$ in space, and work under the condition
$2H_{0}+\sum_{i=1}^{d}H_{i}>d+1$, which forces $H_{0}>\frac{1}{2}$ and therefore
rules out a noise white in time. Read in the Stratonovich sense their formula is
the expected one, and read in the Skorokhod sense the exponential has to be
corrected by a variance term. Chen~\cite{chen:24:feynman-kac} allows a general
spatial covariance, with a temporal covariance $|t-s|^{-\alpha_{0}}$ for some
$\alpha_{0}\in(0,1)$, which is again smoother than white in time. The
representation is then obtained under the condition
\begin{equation*}
  \int_{\R^{d}}\left(\frac{1}{1+|\xi|^{2}}\right)^{1-\alpha_{0}}\mu(d\xi)<\infty,
\end{equation*}
strictly stronger than Dalang's condition~\eqref{eq:dalang-cond}. The singularity
of the noise in time is thus compensated by extra integrability in space.
Note that a true Feynman-Kac representation is very useful if one wishes to
get  the \emph{almost sure} behavior of
$u(t,x)$. For the lattice model with a small
diffusivity $\kappa$, a lower bound of order $1/\log(1/\kappa)$ on the almost
sure Lyapunov exponent is obtained in~\cite{carmona:94:parabolic}, and a
matching upper bound by Carmona, Viens and Molchanov~\cite{carmona:96:sharp}.
The continuous-space counterpart is~\cite{carmona:98:almost-sure}, and the
sharp constant, with its expression in terms of the correlation function of the
potential, is due to Carmona, Koralov and
Molchanov~\cite{carmona:01:asymptotics}.

\smallskip

\noindent
{\it Moment estimates.} Two rather different strategies lead to bounds of the
kind of Theorems~\ref{thm:moment-upper} and~\ref{thm:moment-lower}. Let us summarize both.

\begin{enumerate}
  \item \emph{By It\^{o} calculus and martingale estimates.} One stays inside
    the framework of Section~\ref{sec:integration}, applies the
    Burkholder-Davis-Gundy inequality to the mild formulation, and closes a
    nonlinear renewal equation on the quantity
    $t\mapsto\sup_{x}\E[|u(t,x)|^{k}]$. This is the route of Foondun and
    Khoshnevisan~\cite{foondun:09:intermittence} and of Khoshnevisan's
    lectures~\cite{khoshnevisan:14:analysis}, where the argument is presented
    at length. Its merit is that it never uses the linearity of the equation,
    so it applies verbatim to a Lipschitz $\sigma$, and it yields in that
    generality the \emph{weak intermittency} $\gamma(2)>0$, which by the
    convexity of $k\mapsto\gamma(k)$ already implies the strict monotonicity of
    $\gamma(k)/k$. What it gives less easily is the sharp dependence on $k$.
    For space-time white noise it produces $\gamma(k)=O(k^{3})$, which is the
    value of the exponent $(4-a)/(2-a)$ of~\eqref{eq:moment-upper} at the
    formal endpoint $a=d=1$. The dependence of the moments on the noise intensity
    $\beta$ is studied in the same framework by Khoshnevisan and
    Kim~\cite{khoshnevisan:15:non-linear-a,khoshnevisan:15:non-linear-b}.
  \item \emph{By chaos for the upper bound and Feynman-Kac for the lower one.}
    This is the route of Section~\ref{sec:moment-estimates}. It was developed by
    Hu, Nualart and their coauthors, the hypercontractivity argument behind
    Theorem~\ref{thm:moment-upper} going back to Balan and
    Conus~\cite{balan:16:intermittency}, and the pair of bounds
    in~\eqref{eq:moment-upper}--\eqref{eq:moment-lower} being the white-in-time
    case of~\cite{hu:15:stochastic}. Its interest is that neither half needs a
    martingale structure, which opens the door to noises that are fractional,
    or even rough, in space or in time. See~\cite{hu:09:stochastic}
    and~\cite{hu:17:stochastic,hu:18:parabolic} for a rough spatial noise,
    and Chen, Hu, Kalbasi and Nualart~\cite{chen:18:intermittency} for a rough
    fractional noise in time.
\end{enumerate}

A third body of work, due to Xia Chen and his coauthors, extracts the sharp
constants rather than the sharp exponents. The starting point is again the
Feynman-Kac formula. This reduces the problem to a large-time
behavior estimate of the exponential moments of the weighted intersection local times 
(see Remark~\ref{rmk:fk-intersection}). Chen, Hu, Song and
Xing~\cite{chen:15:exponential} treat those functionals for their own sake, and
the large-deviation machinery for intersection local times is collected in
Chen's monograph~\cite{chen:10:random}. The results this route produces are
precise. For the $(1+1)$-dimensional space-time white noise the moment
Lyapunov exponent $\gamma(k)$ is computed exactly in~\cite{chen:15:precise},
where it is identified with the ground state energy of a $k$-body problem with
pairwise Dirac interaction, which is the $\Lambda=\delta_{0}$ case of the
functional in Remark~\ref{rmk:fk-intersection}. The spatial asymptotics of the
solution are obtained in~\cite{chen:16:spatial}, the moment asymptotics for a
fractional noise in the Skorokhod regime in~\cite{chen:17:moment}, and the
rough and critical cases in~\cite{chen:19:parabolic}. A noise which is
constant in time, as opposed to white, is a genuinely different problem and
calls for other techniques (see~\cite{chen:14:quenched}). The whole of this
chapter then has a counterpart, on $\R^{d}$ as well as on $\Z^{d}$, in which the
moment asymptotics are governed by a variational problem for the principal
eigenvalue rather than by a Feynman-Kac functional, and the phenomenology is
correspondingly richer. G\"{a}rtner and K\"{o}nig~\cite{gartner:00:moment} is
the entry point in the continuous case, with the
survey~\cite{gartner:05:parabolic} and K\"{o}nig's
monograph~\cite{konig:16:parabolic} for the systematic account. That eigenvalue
route is carried out on $\R$ for a \emph{rough} spatial noise, of Hurst
parameter $H<\frac{1}{2}$, in~\cite{chakraborty:20:quenched}. Therein the quenched asymptotics of the
solution are read off from those of the principal eigenvalue of
$\frac{1}{2}\Delta+\dot{W}$.

\smallskip

\noindent
{\it Beyond the moments.} The estimates of
Section~\ref{sec:moment-estimates} are only the first chapter of the
intermittency story. They say how large $\E[u(t,x)^{k}]$ is, and hence that
the solution concentrates on rare tall peaks. However, they say nothing about how
tall a given peak is, where it sits, or how many of them there are. These are
the questions taken up by Khoshnevisan and his coauthors. Conus and
Khoshnevisan~\cite{conus:12:on} locate the farthest high peak at a fixed time,
Conus, Joseph and Khoshnevisan~\cite{conus:13:on} measure the height of the
tall peaks over large spatial regions, and Khoshnevisan, Kim and
Xiao~\cite{khoshnevisan:17:intermittency} carry out a multifractal analysis of
the set of high peaks. The tall islands and the intermittency fronts are
treated systematically in~\cite{khoshnevisan:14:analysis}. For a noise which is
rough in space, where none of the martingale arguments is available, the
corresponding spatial asymptotics are obtained by Chen, Hu, Nualart and
Tindel~\cite{chen:17:spatial}.

\chapter[Geometric setting]{A general geometric setting}
\label{chap:geom-setting}

This chapter and the next one carry the program developed in Chapters~\ref{chap:noise}--\ref{chap:moments} over from $\R^{d}$ to a wide class of geometric structures. The current chapter focuses on the underlying geometric and analytic framework. Namely Section~\ref{sec:geom-settings} introduces the metric measure spaces we shall work with, together with the related heat semigroups and fractional Laplace operators on domains. Then Section~\ref{sec:geom-examples} illustrates this framework with four model examples. The parabolic Anderson model on those structures will be analyzed in Chapter~\ref{chap:geom-pam}.

\section[General setting]{General geometric and analytic settings}\label{sec:geom-settings}

In this section we introduce the geometric framework in which the PAM will be studied.  We work on a metric measure space $(X,d,\mu)$ equipped with a Dirichlet form whose heat kernel satisfies sub-Gaussian estimates.  This general setting covers a wide range of situations of interest: the classical Euclidean case $\R^{d}$ (recovered when $d_{h}=d$ and $d_{w}=2$), compact Riemannian manifolds, metric graphs, and fractal spaces such as the Sierpi\'nski gasket for which the walk dimension $d_{w}$ is strictly greater than $2$.  We first describe the standing assumptions on the ambient space, and then specialize to bounded domains with Dirichlet and Neumann boundary conditions.

\subsection{Recurrent metric measure spaces}\label{ssec:mmspace}

The ambient space throughout this chapter and Chapter~\ref{chap:geom-pam} is a metric measure space satisfying several structural properties, which we collect into a single definition.

\begin{definition}[Metric measure space]\label{def:mmspace}
A \emph{locally compact complete geodesic metric measure space} is a triple $(X,d,\mu)$ where:
\begin{enumerate}
\item \emph{Metric structure.}  $(X,d)$ is a complete metric space that is locally compact (every point has a compact neighborhood) and geodesic (any two points $x,y\in X$ are joined by a curve of length $d(x,y)$).

\item \emph{Measure.}  $\mu$ is a Radon measure on $X$ with full support, i.e., $\mu(U)>0$ for every non-empty open set $U\subset X$.
\end{enumerate}
We write $\mathrm{diam}(X)\in(0,+\infty]$ for the diameter of $X$, and $B(x,r):=\{y\in X:d(x,y)<r\}$ for the open ball of radius $r$ centered at $x$.
\end{definition}

The analytic structure on $X$ is encoded by a Dirichlet form, which replaces the classical gradient and plays the role of the Laplacian in our setting.  Most of the notions introduced below are standard and can be found in~\cite{fukushima:94:dirichlet}.

\begin{definition}[Dirichlet form]\label{def:dirichlet-form}
A \emph{Dirichlet form} on $L^{2}(X,\mu)$ is a pair $(\mathcal{E},\mathcal{F})$, where $\mathcal{F}:=\mathrm{dom}\,\mathcal{E}$ is a dense subspace of $L^{2}(X,\mu)$ and $\mathcal{E}:\mathcal{F}\times\mathcal{F}\to\R$ is a bilinear form satisfying the following properties.
\begin{enumerate}
\item \emph{Closed and symmetric.}  $\mathcal{E}$ is symmetric, and $\mathcal{F}$ is complete for the $\mathcal{E}_{1}$-norm
\begin{equation}\label{eq:E1-norm}
  \|f\|_{\mathcal{E}_{1}} := \bigl(\|f\|_{L^{2}(X,\mu)}^{2}+\mathcal{E}(f,f)\bigr)^{1/2}.
\end{equation}

\item \emph{Markovian.}  If $u\in\mathcal{F}$ and $v$ is a normal contraction of $u$ (i.e., $|v(x)-v(y)|\le|u(x)-u(y)|$ and $|v(x)|\le|u(x)|$ $\mu$-a.e.), then $v\in\mathcal{F}$ and $\mathcal{E}(v,v)\le\mathcal{E}(u,u)$.

\item \emph{Regular.}  There exists a core $\mathcal{C}\subseteq C_{c}(X)\cap\mathcal{F}$ that is dense in $C_{c}(X)$ in the sup-norm and dense in $\mathcal{F}$ in the $\mathcal{E}_{1}$-norm.

\item \emph{Strongly local.}  $\mathcal{E}(u,v)=0$ whenever $u,v\in\mathcal{F}$ have compact support and $u$ is constant on a neighborhood of $\mathrm{supp}(v)$.
\end{enumerate}
\end{definition}

A regular strongly local Dirichlet form carries two fundamental objects: an energy measure and a generator.

\begin{notation}[Energy measure and generator]\label{not:energy-gen}
Let $(\mathcal{E},\mathcal{F})$ be a regular strongly local Dirichlet form on $L^{2}(X,\mu)$.
\begin{enumerate}
\item \emph{Energy measure.}  For every $u,v\in\mathcal{F}\cap L^{\infty}(X,\mu)$, the \emph{energy measure} $\Gamma(u,v)$ is the unique signed Radon measure satisfying
\begin{equation}\label{eq:energy-measure}
  \int_{X}\phi\,d\Gamma(u,v)
  = \frac{1}{2}\bigl[\mathcal{E}(\phi u,v)+\mathcal{E}(\phi v,u)-\mathcal{E}(\phi,uv)\bigr],
  \qquad \phi\in\mathcal{F}\cap C_{c}(X).
\end{equation}

\item \emph{Generator and semigroup.}  The \emph{generator} $\Delta$ is the self-adjoint operator on $L^{2}(X,\mu)$ defined on the dense domain $\mathrm{dom}(\Delta)\subset\mathcal{F}$ by the integration-by-parts formula
\begin{equation}\label{eq:IBP}
  \mathcal{E}(f,g) = -\int_{X}f\,\Delta g\,d\mu, \qquad f\in\mathcal{F},\ g\in\mathrm{dom}(\Delta).
\end{equation}
The associated strongly continuous semigroup is $P_{t}=e^{t\Delta}$~\cite{fukushima:94:dirichlet}.
\end{enumerate}
\end{notation}

A key identity links the energy measure back to the form itself. Its proof can be found in~\cite[Thm.~3.2.1]{fukushima:94:dirichlet}. 

\begin{proposition}[Beurling--Deny]\label{prop:BD}
Under the conditions of Notation~\ref{not:energy-gen},
\begin{equation}\label{eq:BD}
  \mathcal{E}(u,v) = \int_{X}d\Gamma(u,v), \qquad u,v\in\mathcal{F}.
\end{equation}
\end{proposition}

\begin{remark}\label{rmk:energy-gen}
We wish to stress the following facts about the measure $\Gamma$ and the operator $\Delta$.
\begin{enumerate}
\item \emph{Interpretation of $\Gamma$.}  On $\R^{d}$ with $\mathcal{E}(f,g)=\int_{\R^{d}}\nabla f\cdot\nabla g\,dx$, one has $d\Gamma(f,g)=\nabla f\cdot\nabla g\,dx$, so $\Gamma$ is the analogue of the squared gradient.  In particular, $\Gamma(f,f)(A)$ measures the energy of $f$ concentrated on the set $A$.  On a fractal, no pointwise gradient exists, but $\Gamma$ is still well-defined and captures the local energy content of a function.

\item \emph{Definition of $e^{t\Delta}$.}  Since $\Delta$ is a non-positive self-adjoint operator on $L^{2}(X,\mu)$, the spectral theorem provides a projection-valued measure $E(\cdot)$ on $(-\infty,0]$ such that $\Delta=\int_{-\infty}^{0}\lambda\,dE(\lambda)$.  The semigroup is then defined spectrally by $P_{t}=e^{t\Delta}=\int_{-\infty}^{0}e^{t\lambda}\,dE(\lambda)$, which is a contraction on $L^{2}(X,\mu)$ for each $t\ge 0$.  Under Assumption~\ref{hyp:mmspace} below, $P_{t}$ extends to a bounded operator on $L^{\infty}(X,\mu)$ and admits a continuous integral kernel $p_{t}(x,y)$, the \emph{heat kernel}.
\end{enumerate}
\end{remark}

We can now state the precise conditions on $(X,d,\mu,\mathcal{E},\mathcal{F})$ that will be in force throughout this chapter and the next one.  They amount to requiring that balls have a polynomial volume growth (controlled by a Hausdorff dimension $d_{h}$) and that the heat kernel satisfies sub-Gaussian off-diagonal estimates governed by a walk dimension $d_{w}\ge 2$.  The condition $d_{h}<d_{w}$ is absent in the Euclidean setting, where $d_{h}=d$, $d_{w}=2$ and $d\ge 2$.  It forces the recurrence of the associated Brownian motion, and it is characteristic of fractal-like spaces.

\begin{assumption}[Sub-Gaussian framework]\label{hyp:mmspace}
Let $(X,d,\mu)$ be the locally compact complete geodesic metric measure space of Definition~\ref{def:mmspace}, equipped with the regular strongly local Dirichlet form $(\mathcal{E},\mathcal{F})$ of Definition~\ref{def:dirichlet-form}.  We assume:
\begin{enumerate}
\item \emph{Ahlfors regularity.}  There exist $c_{1},c_{2},d_{h}\in(0,\infty)$ such that for all $x\in X$ and $r\in(0,\mathrm{diam}(X))$,
\begin{equation}\label{eq:ahlfors}
  c_{1}\,r^{d_{h}} \;\le\; \mu\bigl(B(x,r)\bigr) \;\le\; c_{2}\,r^{d_{h}}.
\end{equation}
The constant $d_{h}$ is the \emph{Hausdorff dimension} of $(X,d)$.

\item \emph{Regularity and strong locality.}  $(\mathcal{E},\mathcal{F})$ is a regular and strongly local Dirichlet form in the sense of Definition~\ref{def:dirichlet-form}.

\item \emph{Stochastic completeness.}  $P_{t}1=1$ for every $t\ge 0$.

\item \emph{Sub-Gaussian heat kernel.}  The semigroup $\{P_{t}\}_{t\ge 0}$ has a continuous heat kernel $p_{t}(x,y)$ satisfying, for some $c_{3},c_{4},c_{5},c_{6}\in(0,\infty)$ and $d_{w}\in[2,+\infty)$,
\begin{multline}\label{eq:subGauss}
  c_{5}\,t^{-d_{h}/d_{w}}\exp\left(-c_{6}\Bigl(\frac{d(x,y)^{d_{w}}}{t}\Bigr)^{\!\frac{1}{d_{w}-1}}\right)
  \le p_{t}(x,y) \\
  \le c_{3}\,t^{-d_{h}/d_{w}}\exp\left(-c_{4}\Bigl(\frac{d(x,y)^{d_{w}}}{t}\Bigr)^{\!\frac{1}{d_{w}-1}}\right)
\end{multline}
for $\mu\!\times\!\mu$-a.e.\ $(x,y)\in X\times X$ and each $t\in(0,\mathrm{diam}(X)^{d_{w}})$.  The constant $d_{w}$ is called the \emph{walk dimension}.

\item \emph{Recurrence condition.}  $d_{h}<d_{w}$.
\end{enumerate}
\end{assumption}

\begin{remark}
When $d_{w}=2$, the sub-Gaussian estimates~\eqref{eq:subGauss} reduce to the classical Gaussian estimates.  The Euclidean space $\R^{n}$ satisfies Assumption~\ref{hyp:mmspace} with $d_{h}=n$ and $d_{w}=2$ (and the recurrence condition $d_{h}<d_{w}$ fails for $n\ge 2$, consistently with the transience of Brownian motion in $\R^{n}$, $n\ge 2$).  Fractal spaces such as the Sierpi\'nski gasket satisfy $d_{w}>2$; see Section~\ref{ssec:sierpinski}.
\end{remark}

By a fundamental theorem on Dirichlet forms~\cite[Thm.~7.2.1]{fukushima:94:dirichlet}, Assumption~\ref{hyp:mmspace} guarantees the existence of a canonical Markov process associated with $(\mathcal{E},\mathcal{F})$, which we call the Brownian motion on $X$.

\begin{definition}[Brownian motion on $X$]\label{def:BM}
A \emph{Brownian motion} on $(X,d,\mu)$ is a continuous $X$-valued process $(B_{t})_{t\ge 0}$ defined on a filtered probability space $(\Omega^{B},\mathcal{F}^{B},(\mathcal{F}^{B}_{t})_{t\ge 0})$, together with a family of probability measures $(\Prob_{x})_{x\in X}$ on $(\Omega^{B},\mathcal{F}^{B})$, such that:
\begin{enumerate}
\item \emph{Initial condition.}  $\Prob_{x}(B_{0}=x)=1$ for every $x\in X$.
\item \emph{Markov property.}  For every bounded measurable $f:X\to\R$, every $s,t\ge 0$, and every $x\in X$,
\begin{equation}\label{eq:BM-Markov}
  \E_{x}\bigl[f(B_{s+t})\,\big|\,\mathcal{F}^{B}_{t}\bigr]
  \;=\; \E_{B_{t}}\bigl[f(B_{s})\bigr]
  \qquad \Prob_{x}\text{-a.s.}
\end{equation}
\item \emph{Semigroup.}  The semigroup of $(B_{t})_{t\ge 0}$ coincides with $P_{t}$ of Notation~\ref{not:energy-gen}: for every $f\in L^{2}(X,\mu)$ bounded,
\begin{equation}\label{eq:BM-semigroup}
  P_{t}f(x) \;=\; \E_{x}\bigl[f(B_{t})\bigr], \qquad x\in X,\ t\ge 0.
\end{equation}
\end{enumerate}
We write $\E_{x}$ for expectation under $\Prob_{x}$, and $\E_{B}$ when the starting point is clear from context or held fixed.  In joint computations with the noise $W$ defined on $(\Omega,\mathcal{F},\Prob)$, we treat $B$ and $W$ as independent by working on the product space $(\Omega\times\Omega^{B},\mathcal{F}\otimes\mathcal{F}^{B},\Prob\otimes\Prob_{x})$.
\end{definition}

\begin{remark}\label{rmk:g1}
The exponents $d_{h}$ and $d_{w}$ in Assumption~\ref{hyp:mmspace} have the following geometric and probabilistic interpretations.
\begin{enumerate}
\item The Hausdorff dimension $d_{h}$ describes the volume growth of metric balls. Namely we have $\mu(B(x,r))\asymp r^{d_{h}}$.
\item The walk dimension $d_{w}$ describes the spreading of the Brownian motion $(B_{t})_{t\ge 0}$. That is it holds  $\E_{x}[d(x,B_{t})]\asymp t^{1/d_{w}}$.
\end{enumerate}
The condition $d_{h}<d_{w}$ therefore means that the process spreads more slowly than the volume grows, which forces recurrence.
\end{remark}

As mentioned in Remark~\ref{rmk:g1}, the recurrence of the Brownian motion $B$ is a consequence of the condition $d_{h}<d_{w}$ in Assumption~\ref{hyp:mmspace}. We state this property as a lemma, and refer to~\cite[Cor.~3.29]{barlow:98:diffusions} for a complete proof.

\begin{lemma}[Recurrence]\label{lem:recurrence}
Under Assumption~\ref{hyp:mmspace}, the Brownian motion $(B_{t})_{t\ge 0}$ is strongly recurrent: for every $x,y\in X$,
\begin{equation*}
  \Prob_{x}\bigl(\exists\,t\ge 0 : B_{t}=y\bigr) = 1.
\end{equation*}
\end{lemma}

\subsection{Dirichlet boundary conditions}\label{ssec:dirichlet-bc}

We now specialize to a bounded open set $U\subset X$ with compact closure $\overline{U}$.  This subsection introduces the heat semigroup on $U$ with Dirichlet boundary conditions, its heat kernel, and the associated fractional powers of the Laplacian.

\subsubsection{Heat semigroup}\label{sssec:dirichlet-semi}

\begin{definition}[Localized domains]\label{def:FU}
Let $U\subset X$ be a non-empty bounded open set with compact closure.  Set
\begin{equation}\label{eq:F-U}
  \mathcal{F}(U) := \bigl\{u\in L^{2}(U,\mu) : \exists\,f\in\mathcal{F},\ f|_{U}=u\ \mu\text{-a.e.}\bigr\},
\end{equation}
and
\begin{equation}\label{eq:Fc-U}
  \mathcal{F}_{c}(U) := \bigl\{u\in\mathcal{F}(U) : \text{ess-supp}(u)\text{ is compact in }U\bigr\}.
\end{equation}
The \emph{Dirichlet domain} $\mathcal{F}^{D}(U)$ is the closure of $\mathcal{F}_{c}(U)$ in the $\mathcal{E}_{1}$-norm~\eqref{eq:E1-norm}.  Functions in $\mathcal{F}^{D}(U)$ satisfy the Dirichlet boundary condition.
\end{definition}

Restricting the global form $(\mathcal{E},\mathcal{F})$ to $\mathcal{F}^{D}(U)$ yields a well-defined Dirichlet form on $U$, whose generator provides the Laplacian with Dirichlet boundary conditions. We label this as a proposition (borrowed from~\cite{fukushima:94:dirichlet} for further use).

\begin{proposition}[Dirichlet form on $U$]\label{prop:E-D}
The restricted form
\begin{equation}\label{eq:E-D}
  \mathcal{E}^{D}_{U}(u,v) := \mathcal{E}(u,v), \qquad u,v\in\mathcal{F}^{D}(U),
\end{equation}
is a strongly local regular Dirichlet form on $L^{2}(U,\mu)$.  Its generator is denoted $\Delta_{D,U}$ and called the \emph{Dirichlet Laplacian on $U$}.
\end{proposition}

The associated semigroup $P^{D,U}_{t}$ admits a probabilistic representation via the Brownian motion $(B_{t})_{t\ge 0}$ of Lemma~\ref{lem:recurrence}.  Let $T_{U}:=\inf\{t\ge 0:B_{t}\notin U\}$ be the first exit time from $U$.  Then, for $f\in L^{2}(U,\mu)$ bounded and $x\in U$, the following is shown in~\cite[Thm.~4.4.2]{fukushima:94:dirichlet}:
\begin{equation}\label{eq:PDU-prob}
  P^{D,U}_{t}f(x) = \E_{x}\bigl[f(B_{t})\,\mathbf{1}_{t<T_{U}}\bigr].
\end{equation}
Now to obtain sharp heat kernel bounds on $U$ we restrict to \emph{inner uniform domains}, defined as follows.  Recall that a curve $\gamma:[0,1]\to X$ is \emph{rectifiable} if its length is finite, where the length is defined as
\begin{equation*}
\mathrm{length}(\gamma):=\sup\sum_{i}d(\gamma(t_{i}),\gamma(t_{i+1})) ,
\end{equation*}
the supremum above being over all finite subdivisions of $[0,1]$.

\begin{definition}[Inner metric and inner uniform domain]\label{def:inner-uniform}
For a connected open set $U\subset X$, the \emph{inner metric} is
\begin{equation}\label{eq:inner-metric}
  d_{U}(x,y) := \inf\bigl\{\mathrm{length}(\gamma) : \gamma:[0,1]\to U \text{ rectifiable}, \gamma(0)=x, \gamma(1)=y\bigr\}.
\end{equation}
We denote by $\widetilde{U}$ the completion of $U$ for $d_{U}$.

Fix $c\in(0,1)$ and $C>1$.  A rectifiable curve $\gamma:[\alpha,\beta]\to U$ is a \emph{$(c,C)$ inner uniform curve} if:
\begin{enumerate}
\item for all $t\in[\alpha,\beta]$, $d_{U}(\gamma(t),\widetilde{U}\setminus U)\ge c\min\{d_{U}(\gamma(\alpha),\gamma(t)),d_{U}(\gamma(t),\gamma(\beta))\}$;
\item $\mathrm{length}(\gamma)\le C\,d_{U}(\gamma(\alpha),\gamma(\beta))$.
\end{enumerate}
The domain $U$ is called \emph{inner uniform} if any two points in $U$ can be joined by a $(c,C)$ inner uniform curve for some fixed $c,C$.
\end{definition}

\begin{figure}[h]
\centering
\begin{tikzpicture}[line width=1pt]

\begin{scope}[xshift=-4.2cm]
  \draw[fill=gray!12, line width=1.2pt] (0,0) ellipse (2.5 and 1.8);
  \coordinate (A) at (-1.3,-0.2);
  \coordinate (B) at ( 1.3, 0.2);
  \fill (A) circle (2.5pt) node[left=3pt] {$A$};
  \fill (B) circle (2.5pt) node[right=3pt] {$B$};
  \draw[blue!80!black, line width=1.5pt]
    (A) .. controls (-0.3,0.9) and (0.3,0.9) .. (B)
    node[above, midway, font=\small] {$\gamma$};
  \draw[red!70!black, dashed, line width=0.8pt] (0,0.78) circle (0.9);
  \node[red!70!black, font=\scriptsize] at (0.65,1.55) {$d_{\partial U}$};
  \node[font=\small] at (0,-2.2) {(a) inner uniform};
\end{scope}

\begin{scope}[xshift=4.2cm]
  \fill[gray!12] (0,0) circle (2.2);
  \fill[white] (0,-0.09) rectangle (2.25,0.09);
  \draw[line width=1.2pt] (0,0) circle (2.2);
  \draw[black, line width=2pt] (0,0) -- (2.2,0);
  \fill[red!80!black] (0,0) circle (3pt);
  \node[above left, font=\scriptsize] at (0.05,0.05) {tip};
  \coordinate (A) at (0.8, 0.8);
  \coordinate (B) at (0.8,-0.8);
  \fill (A) circle (2.5pt) node[above right=1pt] {$A$};
  \fill (B) circle (2.5pt) node[below right=1pt] {$B$};
  \draw[blue!80!black, line width=1.5pt]
    (A) .. controls (-0.4,0.5) and (-0.4,-0.5) .. (B);
  \node[blue!80!black, left=4pt, font=\small] at (-0.4,0) {$\gamma$};
  \draw[red!70!black, dashed, line width=0.8pt] (0,0) circle (0.32);
  \node[red!70!black, font=\scriptsize] at (0.55,-0.25) {$d_{\partial U}\!\approx 0$};
  \node[font=\small] at (0,-2.6) {(b) non-inner-uniform};
\end{scope}

\end{tikzpicture}
\caption{(a)~An inner uniform curve $\gamma$ in a smooth convex domain: it stays at distance $d_{\partial U}$ from $\partial U$, proportional to the shorter arc length to each endpoint (condition~(i) of Definition~\ref{def:inner-uniform}), and its length is comparable to $d_{U}(A,B)$ (condition~(ii)).
(b)~A slit disk (disk with a radial cut shown as a thick segment).  Any curve from $A$ to $B$ in $U$ must cross $\{y=0\}$ to the left of the slit tip (red dot), since the slit blocks all crossings for $x\ge 0$.  A curve that passes close to the tip (as drawn) violates condition~(i) of Definition~\ref{def:inner-uniform}, while a curve that swings far to the left to avoid the tip incurs excess length, violating condition~(ii).  No fixed $(c,C)$ accommodates both constraints simultaneously, so the domain is not inner uniform.}
\label{fig:inner-uniform}
\end{figure}

\begin{remark}
In Euclidean space, every bounded convex domain is inner uniform, as is every bounded domain with piecewise smooth boundary and non-zero interior angles.  The class of inner uniform domains also includes sets with fractal boundaries such as the Koch snowflake~\cite{lierl:22:dirichlet}.  On the other hand, domains with inward cusps or slits are not inner uniform.  The slit disk of Figure~\ref{fig:inner-uniform} is one example.  Another is the cusp domain $U = \{(x,y)\in\mathbb{R}^{2} : 0 < x < 1,\ |y| < x^{2}\}$: two points $(x_{0}, \pm x_{0}^{2}/2)$ symmetric about the $x$-axis lie at Euclidean distance $x_{0}^{2}$, yet any path in $U$ joining them must travel around the cusp tip at the origin, incurring length of order $x_{0}$; as $x_{0}\to 0$ no fixed constant $C$ can bound the ratio $\mathrm{length}(\gamma)/d_{U}(A,B)$, so condition~(ii) of Definition~\ref{def:inner-uniform} fails.
\end{remark}

The following theorem gives sharp heat kernel estimates for $P^{D,U}_{t}$ on inner uniform domains.  Let $\Phi_{1}$ be the first (positive) eigenfunction of $-\Delta_{D,U}$ with eigenvalue $\lambda_{1}>0$ (see Notation~\ref{not:spectrum-D} below).

\begin{theorem}[\protect{\cite[Thm.~3.3]{lierl:22:dirichlet}}]\label{thm:hk-diri}
Let $U\subset X$ be an inner uniform domain.  Then $(\mathcal{E}^{D}_{U},\mathcal{F}^{D}(U))$ is a strongly local regular Dirichlet form admitting a continuous symmetric heat kernel $p^{D,U}_{t}(x,y)$.  For all $(x,y)\in U\times U$ and $t>0$,
\begin{equation}\label{eq:subGauss-D-upper}
  p^{D,U}_{t}(x,y)
  \;\le\; \frac{c_{1}\,e^{-\lambda_{1}t}}{\min(1,t^{d_{h}/d_{w}})}
  \exp\left(-c_{2}\Bigl(\frac{d(x,y)^{d_{w}}}{t}\Bigr)^{\!\frac{1}{d_{w}-1}}\right),
\end{equation}
and
\begin{equation}\label{eq:subGauss-D-lower}
  p^{D,U}_{t}(x,y)
  \;\ge\; \frac{c_{3}\,\Phi_{1}(x)\Phi_{1}(y)\,e^{-\lambda_{1}t}}{\min(1,t^{d_{h}/d_{w}})}
  \exp\left(-c_{4}\Bigl(\frac{d_{U}(x,y)^{d_{w}}}{t}\Bigr)^{\!\frac{1}{d_{w}-1}}\right),
\end{equation}
for constants $c_{1},c_{2},c_{3},c_{4}>0$.
\end{theorem}

\begin{remark}[Comments on the bounds of Theorem~\ref{thm:hk-diri}]\label{rem:Phi1-positive}
The upper bound~\eqref{eq:subGauss-D-upper} for small $t$ follows from the global estimate~\eqref{eq:subGauss} and $p^{D,U}_{t}(x,y)\le p_{t}(x,y)$, which is immediate from~\eqref{eq:PDU-prob}.  The lower bound involves the first eigenfunction $\Phi_{1}$, which is positive on $U$ by~\cite{lierl:22:dirichlet}. In particular, for every compact $A\subset U$ one has $\inf_{y\in A}\Phi_{1}(y)>0$.
\end{remark}

The positivity of $\Phi_{1}$ invoked in Remark~\ref{rem:Phi1-positive} relies on the inner uniformity of $U$. The next remark records a positivity property of the ground state which holds on general connected open sets.

\begin{remark}[Positivity of the ground state on connected sets]\label{rem:ground-state-positive}
Let $V\subset X$ be a connected bounded open set with compact closure, and assume that $-\Delta_{D,V}$ admits a discrete spectrum together with an orthonormal basis of eigenfunctions, as in Notation~\ref{not:spectrum-D} below. Then the smallest eigenvalue $\lambda_{1}(V)$ is simple, and the associated ground state $\Phi_{1}$ can be chosen to be strictly positive on $V$. This is the classical Perron-Frobenius argument for positivity improving semigroups, for which we refer to~\cite[Thm.~XIII.43 and~XIII.44]{reed:78:methods}. Let us briefly recall the mechanism of the proof. First, owing to the Markov property of the Dirichlet form $\mathcal{E}^{D}_{V}$, replacing $\Phi_{1}$ by $|\Phi_{1}|$ does not increase the  quotient
\begin{equation*}
\frac{\mathcal{E}^{D}_{V}(f,f)}{\|f\|_{L^{2}(V,\mu)}^{2}} .
\end{equation*}
 Hence, up to a change of sign, one can assume $\Phi_{1}\ge0$ on $V$. Second, since $V$ is connected the killed semigroup $P^{D,V}_{t}$ is positivity improving. That is, its kernel $p^{D,V}_{t}(y,z)$ is strictly positive for all $t>0$ and $y,z\in V$. Therefore the eigenfunction relation
 \begin{equation}\label{f4}
\Phi_{1}=e^{\lambda_{1}(V)t}\,P^{D,V}_{t}\Phi_{1}
\end{equation}
 yields $\Phi_{1}(y)>0$ for every $y\in V$. Indeed, since $\Phi_{1}\ge0$ and $\|\Phi_{1}\|_{L^{2}(V,\mu)}=1$, the set $A:=\{z\in V:\,\Phi_{1}(z)>0\}$ has positive measure. Thus writing~\eqref{f4} at a point $y\in V$ we get
\begin{equation*}
\Phi_{1}(y)
=
e^{\lambda_{1}(V)t}\int_{V}p^{D,V}_{t}(y,z)\,\Phi_{1}(z)\,d\mu(z)
\;\ge\;
e^{\lambda_{1}(V)t}\int_{A}p^{D,V}_{t}(y,z)\,\Phi_{1}(z)\,d\mu(z)
\;>\;
0,
\end{equation*}
where the strict inequality is due to the fact that the integrand is strictly positive on $A$ and $\mu(A)>0$. Summarizing, for a connected set $V$ the ground state $\Phi_{1}$ can always be chosen strictly positive on $V$. This property will be applied to small metric balls in the proof of Lemma~\ref{lem:small-ball-geom} in Section~\ref{ssec:geom-moment-lower}.
\end{remark}

The following corollaries of Theorem~\ref{thm:hk-diri} will be used repeatedly.

\begin{corollary}[$L^{2}$ bound]\label{cor:L2-bound-D}
Under the conditions of Theorem~\ref{thm:hk-diri}, there exists $C>0$ such that for all $x\in U$ and $t>0$,
\begin{equation}\label{eq:L2-bound-D}
  \int_{U}p^{D,U}_{t}(x,y)^{2}\,d\mu(y)
  \;\le\;
  C\,\frac{e^{-2\lambda_{1}t}}{\min(1,t^{d_{h}/d_{w}})}.
\end{equation}
\end{corollary}

\begin{proof}
By symmetry of $p^{D,U}_{t}$ (Theorem~\ref{thm:hk-diri}) and the semigroup property, we have
\begin{equation*}
\int_{U}p^{D,U}_{t}(x,y)^{2}\,d\mu(y)=p^{D,U}_{2t}(x,x) .
\end{equation*}
The bound then follows from~\eqref{eq:subGauss-D-upper}.
\end{proof}

\begin{corollary}[$L^{\infty}$ bound]\label{cor:Linf-bound-D}
Under the same conditions, there exists $C>0$ such that for every $t\ge 0$ and $u_{0}\in L^{\infty}(U,\mu)$,
\begin{equation}\label{eq:Linf-bound-D}
  \|P^{D,U}_{t}u_{0}\|_{\infty} \;\le\; C\,e^{-\lambda_{1}t}\|u_{0}\|_{\infty}.
\end{equation}
\end{corollary}

\begin{proof}
For $0\le t\le 1$ use~\eqref{eq:PDU-prob} to write $\|P^{D,U}_{t}u_{0}\|_{\infty}\le\|u_{0}\|_{\infty}$.  For $t>1$, write $P^{D,U}_{t}u_{0}(x)=\int_{U}p^{D,U}_{t}(x,y)\,u_{0}(y)\,d\mu(y)$ and apply Cauchy--Schwarz:
\begin{equation}\label{eq:Linf-CS}
  |P^{D,U}_{t}u_{0}(x)|
  \;\le\; \Bigl(\int_{U}p^{D,U}_{t}(x,y)^{2}\,d\mu(y)\Bigr)^{1/2}
           \Bigl(\int_{U}u_{0}(y)^{2}\,d\mu(y)\Bigr)^{1/2}.
\end{equation}
Since $\mu(U)<\infty$ we have $\|u_{0}\|_{L^{2}(U)}\le\mu(U)^{1/2}\|u_{0}\|_{\infty}$.  For $t>1$ we have $\min(1,t^{d_{h}/d_{w}})=1$, so~\eqref{eq:L2-bound-D} gives $\int_{U}p^{D,U}_{t}(x,y)^{2}\,d\mu(y)\le C\,e^{-2\lambda_{1}t}$.  Substituting into~\eqref{eq:Linf-CS} and taking the supremum over $x\in U$ yields~\eqref{eq:Linf-bound-D} with constant $C^{1/2}\mu(U)^{1/2}$.
\end{proof}

\subsubsection{Fractional Laplacian}\label{sssec:dirichlet-frac}

With the heat-kernel bounds of the previous subsection in hand, we turn to the \emph{fractional powers} $(-\Delta_{D,U})^{-\alpha}$ of the Dirichlet Laplacian and the associated function spaces.  These operators are defined by spectral calculus. We will see that this type of operator admits a Riesz kernel $G^{D,U}_{\alpha}$, obtained as a Laplace-type transform of the heat kernel $p^{D,U}_{t}$ (Lemma~\ref{lem:kernel-D}), and we derive sharp two-sided bounds on $G^{D,U}_{\alpha}$.  These bounds yield the $L^{2}$-estimate of Lemma~\ref{lem:DeltaHK-D} and motivate the Dirichlet Sobolev space $\mathcal{W}^{-\alpha}_{D}(U)$. Those objects are of prime importance for us, since they provide the natural function-space setting for the noise in the PAM on~$U$.

\begin{notation}[Dirichlet spectrum]\label{not:spectrum-D}
The operator $-\Delta_{D,U}$ admits a discrete spectrum $0<\lambda_{1}\le\lambda_{2}\le\cdots$ together with an orthonormal basis $\{\Phi_{j}\}_{j\ge 1}$ of $L^{2}(U,\mu)$ such that $\Delta_{D,U}\Phi_{j}=-\lambda_{j}\Phi_{j}$ for each $j\ge 1$.
\end{notation}

With the spectrum fixed by Notation~\ref{not:spectrum-D}, we record the spectral expansion of the heat kernel and a summability condition on the eigenvalues. Both will be needed in order to show that $(-\Delta_{D,U})^{-\alpha}$ is a bounded operator on $L^{2}(U,\mu)$.

\begin{lemma}[Spectral expansion and eigenvalue series]\label{lem:spectral-D}
Under Notation~\ref{not:spectrum-D}, the Dirichlet heat kernel admits the uniformly convergent expansion
\begin{equation}\label{eq:HK-spectral-D}
  p^{D,U}_{t}(x,y) = \sum_{j=1}^{\infty}e^{-\lambda_{j}t}\Phi_{j}(x)\Phi_{j}(y), \qquad t>0,\ x,y\in U.
\end{equation}
Moreover, for every $\alpha>d_{h}/d_{w}$,
\begin{equation}\label{eq:eigen-series-D}
  \sum_{j=1}^{\infty}\frac{1}{\lambda_{j}^{\alpha}} < +\infty.
\end{equation}
\end{lemma}

\begin{proof}
The expansion~\eqref{eq:HK-spectral-D} is classical, and we refer to~\cite[Thm.~3.44]{barlow:98:diffusions} for its proof.  As for the series~\eqref{eq:eigen-series-D}, we integrate~\eqref{eq:HK-spectral-D} along the diagonal and resort to the on-diagonal part of the upper bound~\eqref{eq:subGauss-D-upper}.  Since $\mu(U)<\infty$, this yields, for $t\in(0,1]$,
\begin{equation}\label{eq:heat-trace-D}
  \sum_{j\ge 1}e^{-\lambda_{j}t}
  = \int_{U}p^{D,U}_{t}(x,x)\,d\mu(x)
  \le c_{1}\,\mu(U)\,t^{-d_{h}/d_{w}}.
\end{equation}
Let us write $N(\lambda):=\#\{j\ge 1:\lambda_{j}\le\lambda\}$ for the eigenvalue counting function.  Choosing $t=1/\lambda$ in~\eqref{eq:heat-trace-D} and bounding the $N(\lambda)$ first terms of the sum from below by $e^{-1}$, we get for $\lambda\ge 1$,
\begin{equation}\label{eq:weyl-upper-D}
  N(\lambda) \le e\sum_{j\ge 1}e^{-\lambda_{j}/\lambda}
  \le c_{1}\,e\,\mu(U)\,\lambda^{d_{h}/d_{w}}.
\end{equation}
In particular only finitely many eigenvalues sit in $(0,1]$, and those contribute a finite amount to~\eqref{eq:eigen-series-D}.  For the remaining ones we split the sum over the dyadic blocks $\{2^{k}<\lambda_{j}\le 2^{k+1}\}$, $k\ge 0$, and apply~\eqref{eq:weyl-upper-D} on each of them.  This gives
\begin{equation*}
  \sum_{\lambda_{j}>1}\frac{1}{\lambda_{j}^{\alpha}}
  \le \sum_{k\ge 0}\frac{N(2^{k+1})}{2^{k\alpha}}
  \le c_{1}\,e\,\mu(U)\,2^{d_{h}/d_{w}}\sum_{k\ge 0}2^{-k(\alpha-d_{h}/d_{w})},
\end{equation*}
which is finite whenever $\alpha>d_{h}/d_{w}$.
\end{proof}

The spectral expansion~\eqref{eq:HK-spectral-D} suggests a natural way to define fractional powers of $-\Delta_{D,U}$: replace the heat-semigroup weights $e^{-\lambda_{j}t}$ by $\lambda_{j}^{-\alpha}$ in the eigenfunction series.

\begin{definition}[Dirichlet fractional Laplacian]\label{def:frac-lap-D}
Working under Notation~\ref{not:spectrum-D}, for $\alpha\ge 0$ and $f\in L^{2}(U,\mu)$ define
\begin{equation}\label{eq:frac-lap-D}
  (-\Delta_{D,U})^{-\alpha}f(x)
  := \sum_{j=1}^{\infty}\frac{1}{\lambda_{j}^{\alpha}}\Bigl(\int_{U}\Phi_{j}(y)f(y)\,d\mu(y)\Bigr)\Phi_{j}(x).
\end{equation}
\end{definition}

The formula~\eqref{eq:frac-lap-D} defines an $L^{2}$ series. The following proposition confirms that this series converges and yields a bounded operator.

\begin{proposition}[$L^{2}$ boundedness of $(-\Delta_{D,U})^{-\alpha}$]\label{prop:frac-lap-D-bdd}
Under Notation~\ref{not:spectrum-D}, for every $\alpha\ge 0$ the operator $(-\Delta_{D,U})^{-\alpha}$ defined in~\eqref{eq:frac-lap-D} is bounded on $L^{2}(U,\mu)$, with
\begin{equation}\label{eq:frac-lap-D-bdd}
  \|(-\Delta_{D,U})^{-\alpha}\|_{L^{2}(U,\mu)\to L^{2}(U,\mu)}\le\lambda_{1}^{-\alpha}.
\end{equation}
\end{proposition}

\begin{proof}
Set $\hat{f}_{j}:=\int_{U}\Phi_{j}(y)f(y)\,d\mu(y)$ for each $j\ge 1$.  Since $\{\Phi_{j}\}$ is an orthonormal basis of $L^{2}(U,\mu)$, Parseval's identity gives $\sum_{j=1}^{\infty}|\hat{f}_{j}|^{2}=\|f\|^{2}_{L^{2}(U,\mu)}$ and 
\begin{equation*}
\|(-\Delta_{D,U})^{-\alpha}f\|^{2}_{L^{2}(U,\mu)}=\sum_{j=1}^{\infty}\lambda_{j}^{-2\alpha}|\hat{f}_{j}|^{2}.
\end{equation*}
Since $\lambda_{j}\ge\lambda_{1}>0$ for all $j\ge 1$, we have $\lambda_{j}^{-\alpha}\le\lambda_{1}^{-\alpha}$. Therefore we get
\begin{equation*}
  \|(-\Delta_{D,U})^{-\alpha}f\|^{2}_{L^{2}(U,\mu)}
  = \sum_{j=1}^{\infty}\lambda_{j}^{-2\alpha}|\hat{f}_{j}|^{2}
  \;\le\; \lambda_{1}^{-2\alpha}\sum_{j=1}^{\infty}|\hat{f}_{j}|^{2}
  = \lambda_{1}^{-2\alpha}\|f\|^{2}_{L^{2}(U,\mu)}.
\end{equation*}
Taking square roots gives~\eqref{eq:frac-lap-D-bdd}, which finishes the proof.
\end{proof}

To identify a kernel for $(-\Delta_{D,U})^{-\alpha}$, we introduce a class of smooth test functions.

\begin{definition}[Dirichlet--Schwartz space]\label{def:schwartz-D}
\begin{equation*}
  \mathcal{S}_{D}(U) := \Bigl\{f\in L^{2}(U,\mu) : \forall\,k\in\N,\ \lim_{n\to\infty}n^{k}\Bigl|\int_{U}\Phi_{n}(y)f(y)\,d\mu(y)\Bigr|=0\Bigr\}.
\end{equation*}
\end{definition}

On $\mathcal{S}_{D}(U)$, the spectral definition~\eqref{eq:frac-lap-D} can be rewritten as a convolution against an explicit kernel. Indeed, the identity $\lambda_{j}^{-\alpha}=\frac{1}{\Gamma(\alpha)}\int_{0}^{\infty}t^{\alpha-1}e^{-\lambda_{j}t}\,dt$ allows to exchange the sum over $j$ with the time integral, and to recover the heat kernel via~\eqref{eq:HK-spectral-D}.

\begin{lemma}[Riesz kernel]\label{lem:kernel-D}
For $\alpha>0$ the operator $(-\Delta_{D,U})^{-\alpha}$ admits a kernel $G^{D,U}_{\alpha}$ given by
\begin{equation}\label{eq:RK-D}
  G^{D,U}_{\alpha}(x,y)
  := \frac{1}{\Gamma(\alpha)}\int_{0}^{+\infty}t^{\alpha-1}p^{D,U}_{t}(x,y)\,dt,
  \qquad x,y\in U,\ x\ne y ,
\end{equation}
such that for $f\in\mathcal{S}_{D}(U)$ and $y\in U$ we have
\begin{equation}\label{eq:RK-D-action}
  (-\Delta_{D,U})^{-\alpha}f(y) = \int_{U}f(x)\,G^{D,U}_{\alpha}(x,y)\,d\mu(x),
\end{equation}
and for $f,g\in\mathcal{S}_{D}(U)$ it holds
\begin{equation}\label{eq:RK-D-bilinear}
  \int_{U}(-\Delta_{D,U})^{-\alpha}f\cdot(-\Delta_{D,U})^{-\alpha}g\,d\mu
  = \int_{U}\!\int_{U}f(x)g(y)\,G^{D,U}_{2\alpha}(x,y)\,d\mu(x)\,d\mu(y).
\end{equation}
\end{lemma}

\begin{proof}
See~\cite[Lem.~2.9 and Cor.~2.11]{baudoin:23:dirichlet} for full details. Here we will just spell out the main steps.

\smallskip

\noindent
{\it Step 1: Convergence of~\eqref{eq:RK-D}.}  For $x\ne y$, the sub-Gaussian upper bound~\eqref{eq:subGauss-D-upper} implies that $p^{D,U}_{t}(x,y)$ decays faster than any power of $t$ as $t\to 0^{+}$, while the spectral expansion~\eqref{eq:HK-spectral-D} gives $p^{D,U}_{t}(x,y)\le Ce^{-\lambda_{1}t}$ for $t\ge 1$.  Hence the integrand $t^{\alpha-1}p^{D,U}_{t}(x,y)$ is integrable on $(0,+\infty)$ for every $\alpha>0$, and~\eqref{eq:RK-D} is well defined.

\smallskip

\noindent
{\it Step 2: Proof of~\eqref{eq:RK-D-action}.}  Set $\hat{f}_{j}:=\int_{U}\Phi_{j}(y)f(y)\,d\mu(y)$.  Since $f\in\mathcal{S}_{D}(U)$ is bounded, the sub-probability bound~\eqref{eq:PDU-prob} for $t\in(0,1]$ and the exponential decay of~\eqref{eq:subGauss-D-upper} for $t>1$ give
\begin{equation*}
  \frac{1}{\Gamma(\alpha)}\int_{0}^{+\infty}t^{\alpha-1}\int_{U}|f(x)|\,p^{D,U}_{t}(x,y)\,d\mu(x)\,dt
  < +\infty,
\end{equation*}
Therefore applying Fubini's theorem to $G^{D,U}_{\alpha}$ defined by~\eqref{eq:RK-D}, we get
\begin{equation*}
  \int_{U}f(x)\,G^{D,U}_{\alpha}(x,y)\,d\mu(x)
  = \frac{1}{\Gamma(\alpha)}\int_{0}^{+\infty}t^{\alpha-1}P^{D,U}_{t}f(y)\,dt.
\end{equation*}
Moreover substituting the spectral expansion~\eqref{eq:HK-spectral-D}, one can write
\begin{equation*}
P^{D,U}_{t}f(y)=\sum_{j=1}^{\infty}e^{-\lambda_{j}t}\hat{f}_{j}\Phi_{j}(y) .
\end{equation*}
Hence interchanging sum and time integral (the rapid decay $|\hat{f}_{j}|=o(j^{-k})$ for every $k\ge 1$, guaranteed by Definition~\ref{def:schwartz-D}, makes this interchange valid) then yield
\begin{equation*}
  \frac{1}{\Gamma(\alpha)}\int_{0}^{+\infty}t^{\alpha-1}P^{D,U}_{t}f(y)\,dt
  = \sum_{j=1}^{\infty}\hat{f}_{j}\Phi_{j}(y)\cdot\frac{1}{\Gamma(\alpha)}\int_{0}^{+\infty}t^{\alpha-1}e^{-\lambda_{j}t}\,dt
  = \sum_{j=1}^{\infty}\frac{\hat{f}_{j}}{\lambda_{j}^{\alpha}}\,\Phi_{j}(y),
\end{equation*}
where the second equality uses $\frac{1}{\Gamma(\alpha)}\int_{0}^{\infty}t^{\alpha-1}e^{-\lambda t}\,dt=\lambda^{-\alpha}$.  By~\eqref{eq:frac-lap-D} the right-hand side equals $(-\Delta_{D,U})^{-\alpha}f(y)$, which gives~\eqref{eq:RK-D-action}.

\smallskip

\noindent
{\it Step 3: Proof of~\eqref{eq:RK-D-bilinear}.}  Let $\hat{g}_{j}:=\int_{U}\Phi_{j}(y)g(y)\,d\mu(y)$.  By Parseval's identity applied twice,
\begin{equation*}
  \int_{U}(-\Delta_{D,U})^{-\alpha}f\cdot(-\Delta_{D,U})^{-\alpha}g\,d\mu
  = \sum_{j=1}^{\infty}\frac{\hat{f}_{j}\,\hat{g}_{j}}{\lambda_{j}^{2\alpha}}
  = \int_{U}f\cdot(-\Delta_{D,U})^{-2\alpha}g\,d\mu.
\end{equation*}
Applying~\eqref{eq:RK-D-action} with $2\alpha$ in place of $\alpha$ to $(-\Delta_{D,U})^{-2\alpha}g$ then gives~\eqref{eq:RK-D-bilinear}.
\end{proof}

The following proposition gives sharp two-sided bounds on $G^{D,U}_{\alpha}$ in three regimes.

\begin{proposition}[Bounds on $G^{D,U}_{\alpha}$]\label{prop:RK-D-bounds}
Let $U\subset X$ be an inner uniform domain and $G^{D,U}_{\alpha}$ be as in~\eqref{eq:RK-D}.  Then:
\begin{enumerate}
\item If $0<\alpha<d_{h}/d_{w}$: there exist $c,C>0$ such that for all $x,y\in U$, $x\ne y$,
\begin{equation}\label{eq:RK-D-small}
  c\,\frac{\Phi_{1}(x)\Phi_{1}(y)}{d_{U}(x,y)^{d_{h}-\alpha d_{w}}}
  \;\le\; G^{D,U}_{\alpha}(x,y)
  \;\le\; \frac{C}{d(x,y)^{d_{h}-\alpha d_{w}}}.
\end{equation}

\item If $\alpha=d_{h}/d_{w}$: there exist $c,C>0$ such that for all $x,y\in U$, $x\ne y$,
\begin{equation}\label{eq:RK-D-crit}
  c\,\Phi_{1}(x)\Phi_{1}(y)\max\bigl(1,|\ln d_{U}(x,y)|\bigr)
  \;\le\; G^{D,U}_{\alpha}(x,y)
  \;\le\; C\max\bigl(1,|\ln d(x,y)|\bigr).
\end{equation}

\item If $\alpha>d_{h}/d_{w}$: there exist $c,C>0$ such that for all $x,y\in U$,
\begin{equation}\label{eq:RK-D-large}
  c\,\Phi_{1}(x)\Phi_{1}(y)
  \;\le\; G^{D,U}_{\alpha}(x,y)
  \;\le\; C.
\end{equation}
\end{enumerate}
\end{proposition}

\begin{proof}
We give the main computations; see~\cite[Prop.~2.6]{baudoin:22:fractional} for complete details.

\smallskip

\noindent
{\it Step 1: Reduction of the upper bounds to $t\in(0,1]$.}  Split~\eqref{eq:RK-D} at $t=1$.  For $t\ge 1$, the bound~\eqref{eq:subGauss-D-upper} gives $p^{D,U}_{t}(x,y)\le c_{1}e^{-\lambda_{1}t}$. Therefore we have
\begin{equation*}
  \frac{1}{\Gamma(\alpha)}\int_{1}^{+\infty}t^{\alpha-1}p^{D,U}_{t}(x,y)\,dt
  \;\le\; \frac{c_{1}}{\Gamma(\alpha)}\int_{1}^{+\infty}t^{\alpha-1}e^{-\lambda_{1}t}\,dt
  \;=:\; C_{\alpha} \;<\; +\infty.
\end{equation*}
For $t\le 1$, we have $\min(1,t^{d_{h}/d_{w}})=t^{d_{h}/d_{w}}$, so~\eqref{eq:subGauss-D-upper} gives
\begin{equation*}
  \frac{1}{\Gamma(\alpha)}\int_{0}^{1}t^{\alpha-1}p^{D,U}_{t}(x,y)\,dt
  \;\le\; \frac{c_{1}}{\Gamma(\alpha)}\int_{0}^{1}t^{\alpha-1-d_{h}/d_{w}}
  \exp\left(-c_{2}\Bigl(\frac{d(x,y)^{d_{w}}}{t}\Bigr)^{\frac{1}{d_{w}-1}}\right)dt.
\end{equation*}
The substitution $t=u\,d(x,y)^{d_{w}}$ transforms the right-hand side into
\begin{equation}\label{eq:RK-upper-subst}
  \frac{c_{1}}{\Gamma(\alpha)}\,d(x,y)^{\alpha d_{w}-d_{h}}
  \int_{0}^{d(x,y)^{-d_{w}}}u^{\alpha-1-d_{h}/d_{w}}
  \exp\left(-c_{2}\Bigl(\frac{1}{u}\Bigr)^{\frac{1}{d_{w}-1}}\right)du.
\end{equation}

\smallskip

\noindent
{\it Step 2: Upper bounds in the three regimes.}

In case (i), $\alpha<d_{h}/d_{w}$: extend the upper limit in~\eqref{eq:RK-upper-subst} to $+\infty$.  The resulting integral is finite because the exponential factor controls $u\to 0$ and the power $u^{\alpha-1-d_{h}/d_{w}}$ is integrable at $+\infty$ (its exponent equals $-1-(d_{h}/d_{w}-\alpha)<-1$).  This gives
\begin{equation}\label{eq:RK-D-up-i}
  G^{D,U}_{\alpha}(x,y) \;\le\; C\,d(x,y)^{\alpha d_{w}-d_{h}}
  \;=\; \frac{C}{d(x,y)^{d_{h}-\alpha d_{w}}}.
\end{equation}

In case (iii), $\alpha>d_{h}/d_{w}$: drop the exponential in the $t$-integral (it is $\le 1$) to get
\begin{equation*}
  \int_{0}^{1}t^{\alpha-1-d_{h}/d_{w}}\,dt \;=\; \frac{d_{w}}{\alpha d_{w}-d_{h}} \;<\; +\infty,
\end{equation*}
which yields
\begin{equation}\label{eq:RK-D-up-iii}
  G^{D,U}_{\alpha}(x,y) \;\le\; C.
\end{equation}

In case (ii), $\alpha=d_{h}/d_{w}$: the integrand in~\eqref{eq:RK-upper-subst} reduces to $u^{-1}\exp(-c_{2}(1/u)^{1/(d_{w}-1)})$.  For $u\in(0,R^{-d_{w}})$ (where $R=\mathrm{diam}(U)$) the super-exponential decay contributes a finite constant.  For $u\in(R^{-d_{w}},d(x,y)^{-d_{w}})$ bound the exponential by $1$ and integrate:
\begin{equation*}
  \int_{R^{-d_{w}}}^{d(x,y)^{-d_{w}}}u^{-1}\,du \;=\; d_{w}|\ln d(x,y)| + C,
\end{equation*}
This gives
\begin{equation}\label{eq:RK-D-up-ii}
  G^{D,U}_{\alpha}(x,y) \;\le\; C\max\bigl(1,|\ln d(x,y)|\bigr).
\end{equation}

\smallskip

\noindent
{\it Step 3: Lower bounds.}  Inserting~\eqref{eq:subGauss-D-lower} into~\eqref{eq:RK-D} gives
\begin{equation}\label{eq:RK-lower-insert}
  G^{D,U}_{\alpha}(x,y)
  \;\ge\; \frac{c_{3}\,\Phi_{1}(x)\Phi_{1}(y)}{\Gamma(\alpha)}
  \int_{0}^{+\infty}
  \frac{t^{\alpha-1}\,e^{-\lambda_{1}t}}{\min(1,t^{d_{h}/d_{w}})}
  \exp\left(-c_{4}\Bigl(\frac{d_{U}(x,y)^{d_{w}}}{t}\Bigr)^{\!\frac{1}{d_{w}-1}}\right)dt.
\end{equation}
For case (ii) and for small $d_{U}(x,y)$ in case (i) we apply the substitution $r=t/d_{U}(x,y)^{d_{w}}$ to~\eqref{eq:RK-lower-insert}; for case (iii) and for large $d_{U}(x,y)$ in case (i) we restrict directly to $t\ge 1$, where $\min(1,t^{d_{h}/d_{w}})=1$. More precisely, one can argue as follows:

In case (i): restrict to $t\in(0,1]$ and substitute $r=t/d_{U}(x,y)^{d_{w}}$ to get
\begin{multline}\label{g2}
  G^{D,U}_{\alpha}(x,y) \;\ge\; 
  \frac{c_{3}\Phi_{1}(x)\Phi_{1}(y)}{\Gamma(\alpha)}\,
  d_{U}(x,y)^{\alpha d_{w}-d_{h}} \\ \times
  \int_{0}^{d_{U}(x,y)^{-d_{w}}}r^{\alpha-1-d_{h}/d_{w}}
  e^{-\lambda_{1}r\,d_{U}(x,y)^{d_{w}}}
  \exp\left(-c_{4}\Bigl(\frac{1}{r}\Bigr)^{\frac{1}{d_{w}-1}}\right)dr.
\end{multline}
Note that for all $r\ge 1$,
\begin{equation}\label{eq:exp-bdd-r1}
  \exp\left(-c_{4}\Bigl(\frac{1}{r}\Bigr)^{\frac{1}{d_{w}-1}}\right) \ge e^{-c_{4}} > 0.
\end{equation}
Since the integrand is non-negative, we may restrict the integration to any subinterval of $(0,d_{U}(x,y)^{-d_{w}}]$.  For $d_{U}(x,y)\le 2^{-1/d_{w}}$ the interval $[1,2]$ lies in $(0,d_{U}(x,y)^{-d_{w}}]$, and on $[1,2]$ each factor in~\eqref{g2} is bounded below by a positive constant: $r^{\alpha-1-d_{h}/d_{w}}$ is bounded away from $0$ on a compact interval; $\exp(-c_{4}(1/r)^{1/(d_{w}-1)})\ge e^{-c_{4}}>0$ by~\eqref{eq:exp-bdd-r1}; and $e^{-\lambda_{1}r\,d_{U}(x,y)^{d_{w}}}\ge e^{-2\lambda_{1}d_{U}(x,y)^{d_{w}}}$, which tends to $1$ as $d_{U}(x,y)\to 0$.  Hence the integral over $[1,2]$ is bounded below by a positive constant, so that for $d_{U}(x,y)\le 2^{-1/d_{w}}$,
\begin{equation}\label{eq:RK-D-low-i-small}
  G^{D,U}_{\alpha}(x,y)
  \;\ge\; c\,\frac{\Phi_{1}(x)\Phi_{1}(y)}{d_{U}(x,y)^{d_{h}-\alpha d_{w}}}.
\end{equation}  As mentioned above, for $d_{U}(x,y)\ge 2^{-1/d_{w}}$ we don't use the above substitution and we argue directly from~\eqref{eq:RK-lower-insert}. Namely, since $d_{U}(x,y)\ge 2^{-1/d_{w}}$, the lower-bound side of~\eqref{eq:RK-D-small} satisfies
\begin{equation*}
\frac{c\,\Phi_{1}(x)\Phi_{1}(y)}{d_{U}(x,y)^{d_{h}-\alpha d_{w}}}\le c\,2^{(d_{h}-\alpha d_{w})/d_{w}}\Phi_{1}(x)\Phi_{1}(y) ,
\end{equation*}
so it suffices to show $G^{D,U}_{\alpha}(x,y)\ge c'\Phi_{1}(x)\Phi_{1}(y)$ for some $c'>0$.  This follows from~\eqref{eq:RK-lower-insert} restricted to $t\ge 1$: for $t\ge 1$ we have $\min(1,t^{d_{h}/d_{w}})=1$, and since $d_{U}(x,y)\le\mathrm{diam}(U)$,
\begin{equation}\label{eq:exp-bdd-t1}
  \Bigl(\frac{d_{U}(x,y)^{d_{w}}}{t}\Bigr)^{\!\frac{1}{d_{w}-1}}
  \;\le\;\mathrm{diam}(U)^{\frac{d_{w}}{d_{w}-1}},
  \qquad t\ge 1,
\end{equation}
so the exponential in~\eqref{eq:RK-lower-insert} satisfies $\exp(-c_{4}(\cdots))\ge c_{5}>0$ uniformly. Gathering this into~\eqref{eq:RK-lower-insert} gives
\begin{equation}\label{eq:RK-D-low-i-large}
G^{D,U}_{\alpha}(x,y)
\;\ge\; c\,\Phi_{1}(x)\Phi_{1}(y)\int_{1}^{+\infty}t^{\alpha-1}e^{-\lambda_{1}t}\,dt
\;=\; c'\,\Phi_{1}(x)\Phi_{1}(y),
\end{equation}
which is the desired inequality. Together with~\eqref{eq:RK-D-low-i-small} this concludes the lower bound for case (i).

In case (iii): restrict the integration in~\eqref{eq:RK-lower-insert} to $t\ge 1$.  For $t\ge 1$ we have $\min(1,t^{d_{h}/d_{w}})=1$, and~\eqref{eq:exp-bdd-t1} gives $\exp(-c_{4}(\cdots))\ge c_{5}>0$ uniformly.  Hence
\begin{equation}\label{eq:RK-D-low-iii}
  G^{D,U}_{\alpha}(x,y) \;\ge\; c\,\Phi_{1}(x)\Phi_{1}(y)\int_{1}^{+\infty}t^{\alpha-1}e^{-\lambda_{1}t}\,dt
  \;=\; c'\,\Phi_{1}(x)\Phi_{1}(y).
\end{equation}

In case (ii): since $\alpha=d_{h}/d_{w}$, for $t\in(0,1]$ we have $\min(1,t^{d_{h}/d_{w}})=t^{d_{h}/d_{w}}$, so the integrand in~\eqref{eq:RK-lower-insert} reduces to $t^{-1}e^{-\lambda_{1}t}\exp(-c_{4}(d_{U}(x,y)^{d_{w}}/t)^{1/(d_{w}-1)})$.  Restricting~\eqref{eq:RK-lower-insert} to $t\in(0,1]$ and applying the substitution $r=t/d_{U}(x,y)^{d_{w}}$ gives
\begin{equation}\label{eq:RK-lower-crit-subst}
  G^{D,U}_{\alpha}(x,y) \;\ge\; \frac{c_{3}\,\Phi_{1}(x)\Phi_{1}(y)}{\Gamma(d_{h}/d_{w})}
  \int_{0}^{d_{U}(x,y)^{-d_{w}}} r^{-1}\,e^{-\lambda_{1}r\,d_{U}(x,y)^{d_{w}}}\,
  \exp\left(-c_{4}\Bigl(\frac{1}{r}\Bigr)^{\frac{1}{d_{w}-1}}\right)dr.
\end{equation}
Restricting to $r\in[1,d_{U}(x,y)^{-d_{w}}]$ in~\eqref{eq:RK-lower-crit-subst} and using~\eqref{eq:exp-bdd-r1} together with $e^{-\lambda_{1}rd_{U}(x,y)^{d_{w}}}\ge e^{-\lambda_{1}}>0$ (since $rd_{U}(x,y)^{d_{w}}\le 1$ on this range), one obtains
\begin{equation}\label{eq:RK-D-low-ii}
  G^{D,U}_{\alpha}(x,y) \;\ge\; c\,\Phi_{1}(x)\Phi_{1}(y)\int_{1}^{d_{U}(x,y)^{-d_{w}}}r^{-1}\,dr
  \;=\; c\,d_{w}\,\Phi_{1}(x)\Phi_{1}(y)\,|\ln d_{U}(x,y)|.
\end{equation}

It remains to gather the six bounds obtained in Steps~2 and~3, case by case.  Assertion~(i) of the proposition, that is~\eqref{eq:RK-D-small}, follows from the upper bound~\eqref{eq:RK-D-up-i} together with the lower bounds~\eqref{eq:RK-D-low-i-small} and~\eqref{eq:RK-D-low-i-large}.  Assertion~(ii), that is~\eqref{eq:RK-D-crit}, follows from~\eqref{eq:RK-D-up-ii} and~\eqref{eq:RK-D-low-ii}.  Assertion~(iii), that is~\eqref{eq:RK-D-large}, follows from~\eqref{eq:RK-D-up-iii} and~\eqref{eq:RK-D-low-iii}.  This completes the proof of Proposition~\ref{prop:RK-D-bounds}.
\end{proof}

The next lemma controls the $L^{2}$-norm of $(-\Delta_{D,U})^{-\alpha}p^{D,U}_{t}(x,\cdot)$, which is one of the key estimates for the existence theory of the PAM on $U$.

\begin{lemma}\label{lem:DeltaHK-D}
Let $\alpha\ge 0$.  There exists $C>0$ such that for every $t>0$ and $x\in U$,
\begin{equation}\label{eq:DeltaHK-D}
  \sup_{x\in U}\|(-\Delta_{D,U})^{-\alpha}p^{D,U}_{t}(x,\cdot)\|_{L^{2}(U,\mu)}
  \;\le\;
  \begin{cases}
    C\,e^{-\lambda_{1}t}(1+t^{\alpha-d_{h}/(2d_{w})}), & 0\le\alpha<\tfrac{d_{h}}{2d_{w}},\\[4pt]
    C\,e^{-\lambda_{1}t}|\ln\min(\tfrac{1}{2},t)|, & \alpha=\tfrac{d_{h}}{2d_{w}},\\[4pt]
    C\,e^{-\lambda_{1}t}, & \alpha>\tfrac{d_{h}}{2d_{w}}.
  \end{cases}
\end{equation}
\end{lemma}

\begin{proof}
Using~\eqref{eq:RK-D} and the semigroup property one writes
\begin{equation}\label{eq:RK-D-semigroup}
(-\Delta_{D,U})^{-\alpha}p^{D,U}_{t}(x,y)
=
\frac{1}{\Gamma(\alpha)}\int_{0}^{\infty}s^{\alpha-1}p^{D,U}_{s+t}(x,y)\,ds ,
\qquad t>0, \; x,y\in U.
\end{equation}
Minkowski's inequality then reduces the estimate to $\int_{0}^{\infty}s^{\alpha-1}\|p^{D,U}_{s+t}(x,\cdot)\|_{L^{2}}\,ds$, which is bounded using~\eqref{eq:L2-bound-D} together with elementary calculus.
\end{proof}

The kernel bounds of Proposition~\ref{prop:RK-D-bounds} and the $L^{2}$-estimate of Lemma~\ref{lem:DeltaHK-D} together suggest a natural Hilbert space of distributions on $U$: its inner product is given by the Riesz kernel $G^{D,U}_{2\alpha}$, or equivalently, by pairing through $(-\Delta_{D,U})^{-\alpha}$ in $L^{2}(U,\mu)$.  This space provides the natural function-space setting for the noise in the PAM on $U$.

\begin{definition}[Dirichlet Sobolev space]\label{def:sobolev-D}
For $\alpha>0$, the \emph{Dirichlet Sobolev space} $\mathcal{W}^{-\alpha}_{D}(U)$ is the completion of $\mathcal{S}_{D}(U)$ with respect to the inner product
\begin{eqnarray}\label{eq:sobolev-D}
  \langle f,g\rangle_{\mathcal{W}^{-\alpha}_{D}(U)}
  &:=& \int_{U}\!\int_{U}f(x)g(y)\,G^{D,U}_{2\alpha}(x,y)\,d\mu(x)\,d\mu(y) \notag\\
  &=& \int_{U}(-\Delta_{D,U})^{-\alpha}f\cdot(-\Delta_{D,U})^{-\alpha}g\,d\mu,
\end{eqnarray}
where the second equality uses~\eqref{eq:RK-D-bilinear}.
\end{definition}

\subsection{Neumann boundary conditions}\label{ssec:neumann-bc}

The Neumann setting mirrors the Dirichlet one, with two key differences. First, the domain of the form is the full $\mathcal{F}(U)$, with no vanishing condition at the boundary. Second, one requires the slightly stronger notion of a \emph{uniform domain}, in order to ensure that the resulting Dirichlet form is regular.

\subsubsection{Heat semigroup}\label{sssec:neumann-semi}

The Neumann form is defined by the same integral expression as the Dirichlet form~\eqref{eq:E-D}, but on the larger function space $\mathcal{F}(U)$ (Definition~\ref{def:FU}) rather than $\mathcal{F}^{D}(U)$: no vanishing condition is imposed at the boundary.  At the probabilistic level, this corresponds to \emph{reflected} rather than killed Brownian motion in $U$ (see Theorem~\ref{thm:hk-neum} below).

\begin{definition}[Neumann Dirichlet form]\label{def:neumann-form}
For a non-empty bounded open set $U\subset X$ with compact closure, the \emph{Neumann Dirichlet form} is
\begin{equation}\label{eq:E-N}
  \mathcal{E}^{N}_{U}(u,v) := \int_{U}d\Gamma(u,v), \qquad u,v\in\mathcal{F}(U),
\end{equation}
where $\Gamma$ is the energy measure~\eqref{eq:energy-measure}.
\end{definition}

For the Dirichlet form~\eqref{eq:E-D}, the inner uniform condition of Definition~\ref{def:inner-uniform} was sufficient to guarantee sharp heat kernel bounds, the two conditions on the curve $\gamma$ being stated in terms of the inner metric $d_{U}$.  For the Neumann form we assume a stronger geometric condition on $U$, namely that the same two conditions hold with the inner metric $d_{U}$ replaced throughout by the ambient metric $d$ of $X$.  This is what the available heat kernel bounds require (see Theorem~\ref{thm:hk-neum}), and the strengthening is a feature of their proof rather than of the reflected process itself.  The inner uniform condition is expected to be enough here as well, with the bounds then written in $d_{U}$, as discussed in the bibliographic notes at the end of this chapter.

\begin{definition}[Uniform domain]\label{def:uniform}
A connected bounded open set $U\subset X$ (with compact closure) is called a \emph{$(c,C)$ uniform domain} if any two points in $U$ can be joined by a rectifiable curve $\gamma\subset U$ satisfying:
\begin{enumerate}
\item $d(\gamma(t),\partial U)\ge c\min\{d(\gamma(\alpha),\gamma(t)),d(\gamma(t),\gamma(\beta))\}$ for all $t$;
\item $\mathrm{length}(\gamma)\le C\,d(\gamma(\alpha),\gamma(\beta))$.
\end{enumerate}
Every uniform domain is inner uniform (Definition~\ref{def:inner-uniform})~\cite{gyrya:11:neumann}.
\end{definition}

\begin{figure}[h]
\centering
\begin{tikzpicture}[line width=1pt]
  \begin{scope}
    \clip (0,0) circle (2.05);
    \fill[gray!15] (0,0) circle (2.05);
    \fill[white] (0,0) circle (1.0);
  \end{scope}
  \fill[white] (0.98,-0.12) rectangle (2.07, 0.12);
  \draw[line width=1.2pt] (0,0) circle (2.0);
  \draw[line width=1.2pt] (0,0) circle (1.0);
  \draw[black, line width=2pt] (1.0,0) -- (2.0,0);
  \coordinate (A) at ({1.5*cos(8)}, {1.5*sin(8)});
  \coordinate (B) at ({1.5*cos(-8)}, {1.5*sin(-8)});
  \fill (A) circle (2.5pt) node[above right=2pt] {$A$};
  \fill (B) circle (2.5pt) node[below right=2pt] {$B$};
  \draw[red!80!black, <->, line width=0.9pt]
    ({1.5*cos(8)},{1.5*sin(8)}) -- ({1.5*cos(-8)},{1.5*sin(-8)})
    node[right=3pt, midway, font=\scriptsize, red!80!black] {$d(A,B)$};
  \draw[blue!80!black, line width=1.6pt, ->]
    (8:1.5) arc[start angle=8, end angle=352, radius=1.5]
    node[above, near start, font=\small, xshift=-12pt] {};
  \node[blue!80!black, font=\small, left=6pt] at (-1.5,0.25) {$\gamma$};
  \draw[SeaGreen!80!black, dashed, line width=0.8pt] (-1.5,0) circle (0.5);
  \node[SeaGreen!70!black, font=\scriptsize] at (-1.5,-0.65) {$d_{U}(\gamma(t),\partial U)$};
  \node[blue!80!black, font=\scriptsize] at (0.1,-1.7)
    {$d_{U}(A,B)\approx 2\pi\bar{r}\gg d(A,B)$};
\end{tikzpicture}
\caption{The slit annulus $U$.  Points $A$ and $B$ lie at radius $\bar{r}$ just above and below the slit; their Euclidean distance $d(A,B)=2\bar{r}\sin\varepsilon$ is small, but every path $\gamma\subset U$ must travel almost all the way around (inner geodesic in blue), so $d_{U}(A,B)\approx 2\pi\bar{r}\gg d(A,B)$.  The dashed circle illustrates that the midpoint of $\gamma$ stays at inner distance $\approx (R_{2}-R_{1})/2$ from the boundary, which is proportional to the shorter arc-length to each endpoint, confirming condition~(i) of Definition~\ref{def:inner-uniform}.  Condition~(ii) of Definition~\ref{def:uniform} fails because no fixed $C$ satisfies $\mathrm{length}(\gamma)\le C\,d(A,B)$ for all $\varepsilon>0$.}
\label{fig:slit-annulus}
\end{figure}
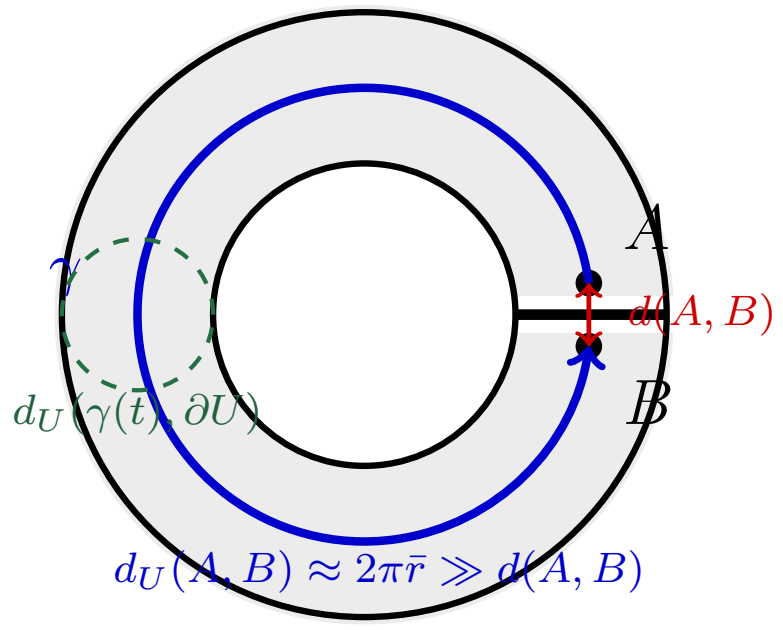

\begin{example}[Slit annulus: inner uniform but not uniform]\label{ex:slit-annulus}
Let $X=\mathbb{R}^{2}$ with the Euclidean metric, fix $0<R_{1}<R_{2}$, and set $\bar{r}:=(R_{1}+R_{2})/2$.  In polar coordinates $(r,\theta)$, define
\begin{equation*}
  U \;:=\; \bigl\{(r,\theta) : R_{1}<r<R_{2},\; 0<\theta<2\pi\bigr\},
\end{equation*}
i.e., the open annulus $\{R_{1}<|x|<R_{2}\}$ with the positive real axis removed (Figure~\ref{fig:slit-annulus}).  Then $U$ is inner uniform but \emph{not} uniform.
\end{example}

\begin{proof}

\smallskip\noindent
{\it $U$ is not uniform.}  Let $A=(\bar{r}\cos\varepsilon,\,\bar{r}\sin\varepsilon)$ and $B=(\bar{r}\cos\varepsilon,\,-\bar{r}\sin\varepsilon)$ for small $\varepsilon>0$.  Their Euclidean distance is $d(A,B)=2\bar{r}\sin\varepsilon\to 0$ as $\varepsilon\to 0$.  Since the slit $\{(x,0):x>0\}$ blocks any direct crossing, every path $\gamma\subset U$ from $A$ to $B$ must travel almost all the way around the annulus, so $\mathrm{length}(\gamma)\ge 2\pi\bar{r}(1-\varepsilon/\pi)$.  No fixed constant $C$ can satisfy $\mathrm{length}(\gamma)\le C\,d(A,B)$ for all $\varepsilon>0$, and condition~(ii) of Definition~\ref{def:uniform} fails.

\smallskip\noindent
{\it $U$ is inner uniform.}  In the completion $\widetilde{U}$ under $d_{U}$, the two sides of the slit are distinct boundary components.  For any $x=(r,\theta)\in U$ the inner distance to $\widetilde{U}\setminus U$ satisfies
\begin{equation}\label{eq:slit-annulus-bdry}
  d_{U}(x,\widetilde{U}\setminus U) \;\ge\; \min\!\bigl(r-R_{1},\;R_{2}-r,\;r\theta,\;r(2\pi-\theta)\bigr),
\end{equation}
since the nearest boundary component is either one of the two circles or one of the two sides of the slit (reached by an arc of length $r\theta$ or $r(2\pi-\theta)$ respectively).  We verify condition~(i) of Definition~\ref{def:inner-uniform} for the representative pair $A=(\bar{r},\varepsilon)$, $B=(\bar{r},2\pi-\varepsilon)$ above; the general case is analogous.  The inner geodesic is the arc $\gamma(t)=(\bar{r},\theta(t))$ with $\theta(t)=\varepsilon+(2\pi-2\varepsilon)t$, $t\in[0,1]$.  Since $\bar{r}-R_{1}=R_{2}-\bar{r}=(R_{2}-R_{1})/2$, \eqref{eq:slit-annulus-bdry} gives
\begin{equation*}
  d_{U}(\gamma(t),\widetilde{U}\setminus U)
  \;\ge\; \min\!\Bigl(\tfrac{R_{2}-R_{1}}{2},\;\bar{r}\,\theta(t),\;\bar{r}(2\pi-\theta(t))\Bigr).
\end{equation*}
The inner arc-lengths from the endpoints are $d_{U}(A,\gamma(t))=\bar{r}(\theta(t)-\varepsilon)$ and $d_{U}(\gamma(t),B)=\bar{r}(2\pi-\varepsilon-\theta(t))$.  We consider three cases.
\begin{itemize}
\item $\bar{r}\,\theta(t)\le(R_{2}-R_{1})/2$ (near $A$, close to the slit): then $d_{U}(\gamma(t),\widetilde{U}\setminus U)\ge\bar{r}\,\theta(t)\ge\bar{r}(\theta(t)-\varepsilon)=d_{U}(A,\gamma(t))$, so condition~(i) holds with $c=1$.
\item $\bar{r}(2\pi-\theta(t))\le(R_{2}-R_{1})/2$ (near $B$, close to the slit): symmetric, condition~(i) holds with $c=1$.
\item Both $\bar{r}\,\theta(t)>(R_{2}-R_{1})/2$ and $\bar{r}(2\pi-\theta(t))>(R_{2}-R_{1})/2$ (middle of the arc): then $d_{U}(\gamma(t),\widetilde{U}\setminus U)\ge(R_{2}-R_{1})/2$, while $\min\{d_{U}(A,\gamma(t)),d_{U}(\gamma(t),B)\}\le\bar{r}\pi$ (worst at the midpoint $\theta=\pi$).  Condition~(i) requires $(R_{2}-R_{1})/2\ge c\bar{r}\pi$, i.e., $c\le(R_{2}-R_{1})/(2\pi\bar{r})$.  Since $\bar{r}\le R_{2}$, taking $c=(R_{2}-R_{1})/(2\pi R_{2})$ suffices.
\end{itemize}
This proves condition (i) for our case of study. Condition~(ii) is trivially satisfied since $\gamma$ is a geodesic.
\end{proof}

The Neumann analogue of Theorem~\ref{thm:hk-diri} is~\cite[Thm.~2.8]{murugan:24:heat}.  It requires the stronger uniform condition of Definition~\ref{def:uniform}, and the resulting bounds differ from the Dirichlet case in three respects.  First, there is no factor $e^{-\lambda_{1}t}$: the Neumann Laplacian has $0$ as an eigenvalue (the constant function), so the semigroup does not decay exponentially.  Second, the lower bound carries no $\Phi_{1}(x)\Phi_{1}(y)$ factor, since the reflected Brownian motion does not vanish at the boundary.  Third, the ambient metric $d(x,y)$ appears throughout (in both upper and lower bounds), whereas the Dirichlet lower bound~\eqref{eq:subGauss-D-lower} involves the inner metric $d_{U}$.

\begin{theorem}\label{thm:hk-neum}
Let $U\subset X$ be a uniform domain.  Then $(\mathcal{E}^{N}_{U},\mathcal{F}(U))$ is a strongly local regular Dirichlet form with a continuous heat kernel $p^{N,U}_{t}(x,y)$ satisfying, for all $(x,y)\in U\times U$ and $t>0$,
\begin{equation}\label{eq:subGauss-N-upper}
  p^{N,U}_{t}(x,y)
  \;\le\; \frac{c_{1}}{\min(1,t^{d_{h}/d_{w}})}
  \exp\left(-c_{2}\Bigl(\frac{d(x,y)^{d_{w}}}{t}\Bigr)^{\!\frac{1}{d_{w}-1}}\right),
\end{equation}
\begin{equation}\label{eq:subGauss-N-lower}
  p^{N,U}_{t}(x,y)
  \;\ge\; \frac{c_{3}}{\min(1,t^{d_{h}/d_{w}})}
  \exp\left(-c_{4}\Bigl(\frac{d(x,y)^{d_{w}}}{t}\Bigr)^{\!\frac{1}{d_{w}-1}}\right).
\end{equation}
The associated diffusion process is the \emph{reflected Brownian motion} in $U$.
\end{theorem}

The following $L^{2}$ bound is the Neumann analogue of Corollary~\ref{cor:L2-bound-D}.  The proof is the same. Namely, the symmetry of $p^{N,U}_{t}$ together with the semigroup property give $\int_{U}p^{N,U}_{t}(x,y)^{2}\,d\mu(y)=p^{N,U}_{2t}(x,x)$, and the bound follows from~\eqref{eq:subGauss-N-upper} with $d(x,x)=0$.  Compared with~\eqref{eq:L2-bound-D}, the only difference is the absence of the factor $e^{-2\lambda_{1}t}$. Indeed, the Neumann spectrum contains $0$, so that there is no exponential decay.

\begin{corollary}\label{cor:L2-bound-N}
Under the conditions of Theorem~\ref{thm:hk-neum}, there exists $C>0$ such that
\begin{equation}\label{eq:L2-bound-N}
  \int_{U}p^{N,U}_{t}(x,y)^{2}\,d\mu(y)
  \;\le\;
  C\min\!\bigl(t^{-d_{h}/d_{w}},1\bigr), \qquad t>0,\ x\in U.
\end{equation}
\end{corollary}

\subsubsection{Fractional Laplacian}\label{sssec:neumann-frac}

The spectral theory of the Neumann Laplacian $-\Delta_{N,U}$ (the generator of $(\mathcal{E}^{N}_{U},\mathcal{F}(U))$) parallels Notation~\ref{not:spectrum-D} and Lemma~\ref{lem:spectral-D}, with one key difference. Namely constant functions belong to $\mathcal{F}(U)$ and have zero energy ($\mathcal{E}^{N}_{U}(c,v)=\int_{U}d\Gamma(c,v)=0$ for every $v\in\mathcal{F}(U)$), so that $0$ is an eigenvalue.  As a result, the heat kernel expansion acquires an additional constant term $1/\mu(U)$. This reflects the fact that $p^{N,U}_{t}(x,\cdot)\to 1/\mu(U)$ in $L^{2}$ as $t\to\infty$, that is the reflected Brownian motion equidistributes over $U$.

\begin{proposition}[Neumann spectrum]\label{prop:spectrum-N}
Consider a uniform domain $U$.  The operator $-\Delta_{N,U}$ has a discrete spectrum $0=\nu_{0}<\nu_{1}\le\nu_{2}\le\cdots$ with corresponding eigenfunctions: $\Psi_{0}=\mu(U)^{-1/2}$ (constant), and $\{\Psi_{j}\}_{j\ge 1}$ an orthonormal basis of $L^{2}_{0}(U,\mu):=\{f\in L^{2}(U,\mu):\int_{U}f\,d\mu=0\}$.  The Neumann heat kernel admits the expansion
\begin{equation}\label{eq:HK-spectral-N}
  p^{N,U}_{t}(x,y) = \frac{1}{\mu(U)}+\sum_{j=1}^{\infty}e^{-\nu_{j}t}\Psi_{j}(x)\Psi_{j}(y),
  \qquad t>0,\ x,y\in U.
\end{equation}
Moreover, $\sum_{j=1}^{\infty}\nu_{j}^{-\alpha}<+\infty$ for every $\alpha>d_{h}/d_{w}$.
\end{proposition}

The fractional Neumann Laplacian is defined by the same spectral formula as Definition~\ref{def:frac-lap-D}, with two modifications forced by $\nu_{0}=0$.  First, the sum must start at $j=1$.  Second, since $\Psi_{0}$ is a constant, $\langle f,\Psi_{0}\rangle_{L^{2}}\Psi_{0}$ is the mean of $f$, The sum over $j\ge 1$ thus projects $f$ onto its mean-zero part, and $(-\Delta_{N,U})^{-\alpha}$ maps $L^{2}(U,\mu)$ into $L^{2}_{0}(U,\mu)$.

\begin{definition}[Neumann fractional Laplacian]\label{def:frac-lap-N}
For $\alpha>0$ and $f\in L^{2}(U,\mu)$,
\begin{equation}\label{eq:frac-lap-N}
  (-\Delta_{N,U})^{-\alpha}f
  := \sum_{j=1}^{\infty}\frac{1}{\nu_{j}^{\alpha}}\Psi_{j}\int_{U}\Psi_{j}(y)f(y)\,d\mu(y).
\end{equation}
This is a bounded operator $L^{2}(U,\mu)\to L^{2}_{0}(U,\mu)$ with norm at most $\nu_{1}^{-\alpha}$.
\end{definition}

Defining the Neumann--Schwartz space $\mathcal{S}_{N}(U)$ analogously to Definition~\ref{def:schwartz-D} (replacing $\Phi_{n}$ by $\Psi_{n}$), the operator $(-\Delta_{N,U})^{-\alpha}$ has a kernel.

\begin{lemma}[Neumann Riesz kernel]\label{lem:kernel-N}
For $\alpha>0$, $(-\Delta_{N,U})^{-\alpha}$ admits the kernel
\begin{equation}\label{eq:RK-N}
  G^{N,U}_{\alpha}(x,y)
  := \frac{1}{\Gamma(\alpha)}\int_{0}^{+\infty}t^{\alpha-1}\!\left(p^{N,U}_{t}(x,y)-\frac{1}{\mu(U)}\right)dt,
  \qquad x,y\in U,\ x\ne y,
\end{equation}
with the representation $(-\Delta_{N,U})^{-\alpha}f(y)=\int_{U}f(x)G^{N,U}_{\alpha}(x,y)\,d\mu(x)$ for $f\in\mathcal{S}_{N}(U)$, and the bilinear formula analogous to~\eqref{eq:RK-D-bilinear}.
\end{lemma}

\begin{proof}
The argument follows the same three steps as the proof of Lemma~\ref{lem:kernel-D}; see~\cite{baudoin:22:fractional} for full details.  We highlight the differences due to the zero eigenvalue.

\smallskip

\noindent
{\it Step 1: Convergence of~\eqref{eq:RK-N}.}  For small $t$, the sub-Gaussian upper bound~\eqref{eq:subGauss-N-upper} controls $p^{N,U}_{t}(x,y)-1/\mu(U)$ exactly as~\eqref{eq:subGauss-D-upper} controls $p^{D,U}_{t}(x,y)$ in the Dirichlet case.  For large $t$ the key difference appears. Unlike $p^{D,U}_{t}$, the Neumann heat kernel does \emph{not} decay to $0$ as $t\to\infty$, since it converges to $1/\mu(U)$.  The quantity $t^{\alpha-1}p^{N,U}_{t}(x,y)$ would thus not be integrable at $+\infty$.  The subtraction of $1/\mu(U)$ restores integrability. Indeed, owing to~\eqref{eq:HK-spectral-N},
\begin{equation*}
  p^{N,U}_{t}(x,y)-\frac{1}{\mu(U)}
  = \sum_{j=1}^{\infty}e^{-\nu_{j}t}\Psi_{j}(x)\Psi_{j}(y),
\end{equation*}
which decays like $e^{-\nu_{1}t}$ as $t\to\infty$, making $t^{\alpha-1}(p^{N,U}_{t}(x,y)-1/\mu(U))$ integrable on $(1,+\infty)$ for any $\alpha>0$.

\smallskip

\noindent
{\it Step 2: Kernel representation.}  Set $\hat{f}_{j}:=\int_{U}\Psi_{j}(y)f(y)\,d\mu(y)$.  For $f\in\mathcal{S}_{N}(U)$ the same Fubini argument as in Step~2 of Lemma~\ref{lem:kernel-D} gives
\begin{equation*}
  \int_{U}f(x)\,G^{N,U}_{\alpha}(x,y)\,d\mu(x)
  = \frac{1}{\Gamma(\alpha)}\int_{0}^{+\infty}t^{\alpha-1}
  \Bigl(P^{N,U}_{t}f(y)-\frac{1}{\mu(U)}\int_{U}f\,d\mu\Bigr)\,dt.
\end{equation*}
By~\eqref{eq:HK-spectral-N}, since the $j=0$ terms cancel, it is readily checked that
\begin{equation*}
P^{N,U}_{t}f(y)-\frac{1}{\mu(U)}\int_{U}f\,d\mu = \sum_{j=1}^{\infty}e^{-\nu_{j}t}\hat{f}_{j}\Psi_{j}(y)
\end{equation*}
Interchanging sum and integral via the rapid decay of $\hat{f}_{j}$ (Definition of $\mathcal{S}_{N}(U)$) and applying the gamma identity $\frac{1}{\Gamma(\alpha)}\int_{0}^{\infty}t^{\alpha-1}e^{-\nu t}\,dt=\nu^{-\alpha}$ to each $j\ge 1$ yields
\begin{equation*}
  \int_{U}f(x)\,G^{N,U}_{\alpha}(x,y)\,d\mu(x)
  = \sum_{j=1}^{\infty}\frac{\hat{f}_{j}}{\nu_{j}^{\alpha}}\Psi_{j}(y)
  = (-\Delta_{N,U})^{-\alpha}f(y).
\end{equation*}

\smallskip

\noindent
{\it Step 3: Bilinear formula.}  Identical to Step~3 of Lemma~\ref{lem:kernel-D}, with $\nu_{j}$ and $\Psi_{j}$ in place of $\lambda_{j}$ and $\Phi_{j}$, and Parseval's identity applied to the sum over $j\ge 1$.
\end{proof}

The following proposition gives sharp two-sided bounds on $G^{N,U}_{\alpha}$ in the same three regimes as Proposition~\ref{prop:RK-D-bounds}.  The structure is parallel. The leading singularity near $x=y$ is again governed by $d(x,y)^{\alpha d_{w}-d_{h}}$, and the three regimes are still determined by the comparison of $\alpha$ with $d_{h}/d_{w}$.  Two differences appear.  First, the ambient metric $d(x,y)$ replaces the inner metric $d_{U}(x,y)$ throughout.  This is inherited from the Neumann heat kernel bounds of Theorem~\ref{thm:hk-neum}.  Second, the lower bounds carry a negative additive constant $-C$. Indeed, $G^{N,U}_{\alpha}$ is defined through the mean-subtracted kernel $p^{N,U}_{t}-1/\mu(U)$, which changes sign. The kernel is therefore not pointwise non-negative, and no $\Phi_{1}(x)\Phi_{1}(y)$ factor appears in the lower bounds.

\begin{proposition}[Bounds on $G^{N,U}_{\alpha}$]\label{prop:RK-N-bounds}
Let $U$ be a uniform domain and $G^{N,U}_{\alpha}$ as in~\eqref{eq:RK-N}.  Then:
\begin{enumerate}
\item If $0<\alpha<d_{h}/d_{w}$: there exist $C_{1},C_{2},C_{3}>0$ such that for all $x,y\in U$, $x\ne y$,
\begin{equation}\label{eq:RK-N-small}
  C_{1}\,d(x,y)^{\alpha d_{w}-d_{h}}-C_{2}
  \;\le\; G^{N,U}_{\alpha}(x,y)
  \;\le\; C_{3}\,d(x,y)^{\alpha d_{w}-d_{h}}.
\end{equation}
\item If $\alpha=d_{h}/d_{w}$: there exist $C_{1},C_{2},C_{3}>0$ such that for all $x,y\in U$, $x\ne y$,
\begin{equation}\label{eq:RK-N-crit}
  -C_{1}\ln d(x,y)-C_{2}
  \;\le\; G^{N,U}_{\alpha}(x,y)
  \;\le\; C_{3}\max\bigl(|\ln d(x,y)|,1\bigr).
\end{equation}
\item If $\alpha>d_{h}/d_{w}$: there exist $C_{1},C_{2}>0$ such that for all $x,y\in U$,
\begin{equation}\label{eq:RK-N-large}
  -C_{1} \;\le\; G^{N,U}_{\alpha}(x,y) \;\le\; C_{2}.
\end{equation}
\end{enumerate}
\end{proposition}

\begin{proof}
The argument parallels the proof of Proposition~\ref{prop:RK-D-bounds}; see~\cite[Prop.~2.6]{baudoin:22:fractional} for full details.  We describe the main differences.

\smallskip

\noindent
{\it Step 1: Upper bounds.}  Split~\eqref{eq:RK-N} at $t=1$:
\begin{multline}\label{eq:RK-N-split}
  G^{N,U}_{\alpha}(x,y)
  = \frac{1}{\Gamma(\alpha)}\int_{0}^{1}t^{\alpha-1}\!\left(p^{N,U}_{t}(x,y)-\frac{1}{\mu(U)}\right)dt \\
  + \frac{1}{\Gamma(\alpha)}\int_{1}^{+\infty}t^{\alpha-1}\!\left(p^{N,U}_{t}(x,y)-\frac{1}{\mu(U)}\right)dt ,
\end{multline}
where we recall that $p^{N,U}_{t}$ is defined by~\eqref{eq:HK-spectral-N}.
Hence for $t\ge 1$, the spectral expansion~\eqref{eq:HK-spectral-N} gives $|p^{N,U}_{t}(x,y)-1/\mu(U)|\le Ce^{-\nu_{1}t}$, so the second integral in~\eqref{eq:RK-N-split} is bounded in absolute value by a constant $C$ independent of $x$ and $y$.  Note that in the Dirichlet case, it is $p^{D,U}_{t}(x,y)$ itself that decays like $e^{-\lambda_{1}t}$, via~\eqref{eq:subGauss-D-upper}. Here the decay comes only after subtracting $1/\mu(U)$, via~\eqref{eq:HK-spectral-N}.  Therefore, the large-$t$ integral contributes only a bounded additive constant to $G^{N,U}_{\alpha}(x,y)$, uniformly in $(x,y)\in U\times U$; all singular dependence on the pair $(x,y)$ comes from the integral over $[0,1]$.

Let us now handle the case $t\le 1$. In this case the term $-1/\mu(U)$ contributes $-c_{\alpha}/\mu(U)$, a bounded constant, and the main contribution is the term 
\begin{equation*}
\frac{1}{\Gamma(\alpha)}\int_{0}^{1}t^{\alpha-1}p^{N,U}_{t}(x,y)\,dt .
\end{equation*}
For this term we resort to~\eqref{eq:subGauss-N-upper} and the substitution $t=ud(x,y)^{d_{w}}$, with the \emph{ambient} metric $d(x,y)$ in place of the inner metric $d_{U}(x,y)$ used in Proposition~\ref{prop:RK-D-bounds}. This yields the upper bounds~\eqref{eq:RK-N-small}--\eqref{eq:RK-N-large} by the same three-case analysis as Steps~1--2 of Proposition~\ref{prop:RK-D-bounds}.

\smallskip

\noindent
{\it Step 2: Lower bounds.}  The large-$t$ integral in~\eqref{eq:RK-N-split} satisfies
\begin{equation*}
\frac{1}{\Gamma(\alpha)}\int_{1}^{\infty}t^{\alpha-1}(p^{N,U}_{t}-1/\mu(U))\,dt\ge -C ,
\end{equation*}
by the same spectral estimate $|p^{N,U}_{t}(x,y)-1/\mu(U)|\le Ce^{-\nu_{1}t}$ used in Step~1, and the small-$t$ subtracted term contributes the bounded constant $-c_{\alpha}/\mu(U)$.  Combining gives
\begin{equation}\label{eq:RK-N-lower-reduce}
  G^{N,U}_{\alpha}(x,y) \;\ge\; \frac{1}{\Gamma(\alpha)}\int_{0}^{1}t^{\alpha-1}p^{N,U}_{t}(x,y)\,dt \;-\; C.
\end{equation}
Inserting the lower bound~\eqref{eq:subGauss-N-lower} into~\eqref{eq:RK-N-lower-reduce} and applying the substitution $t=ud(x,y)^{d_{w}}$ produces the lower bounds in~\eqref{eq:RK-N-small}--\eqref{eq:RK-N-large}, with three differences from Step~3 of Proposition~\ref{prop:RK-D-bounds}:
\begin{enumerate}
\item The ambient metric $d(x,y)$ appears throughout, because~\eqref{eq:subGauss-N-lower} is stated in terms of $d$.
\item There is no factor $\Phi_{1}(x)\Phi_{1}(y)$, because the Neumann lower bound~\eqref{eq:subGauss-N-lower} requires no eigenfunction weight.
\item The additive $-C$ in~\eqref{eq:RK-N-lower-reduce} shifts all three lower bounds downward, producing the negative additive constants in~\eqref{eq:RK-N-small}--\eqref{eq:RK-N-large}.  In particular, for case~(iii) the contribution from $\int_{0}^{1}t^{\alpha-1}p^{N,U}_{t}\,dt$ is bounded (not growing in $d(x,y)$), so the lower bound reduces to $G^{N,U}_{\alpha}(x,y)\ge -C_{1}$.
\end{enumerate}
Combining Steps~1 and~2, the upper bounds in~\eqref{eq:RK-N-small}--\eqref{eq:RK-N-large} follow from Step~1, and the lower bounds follow from Step~2 together with the three items above.  This completes the proof of Proposition~\ref{prop:RK-N-bounds}.
\end{proof}

The kernel bounds of Proposition~\ref{prop:RK-N-bounds} likewise suggest a natural Hilbert space of distributions on $U$, paralleling Definition~\ref{def:sobolev-D}.  There is one essential difference. The lower bounds in~\eqref{eq:RK-N-small}--\eqref{eq:RK-N-large} all carry a negative additive constant, so that $G^{N,U}_{2\alpha}$ is not pointwise non-negative. The kernel alone does not define a positive inner product.  Adding a sufficiently large constant $C_{U}$ restores positivity, and the resulting topology on $\mathcal{W}^{-\alpha}_{N}(U)$ does not depend on the choice of $C_{U}$.

\begin{definition}[Neumann Sobolev space]\label{def:sobolev-N}
For $\alpha>0$ and a constant $C_{U}>0$ chosen so that $G^{N,U}_{2\alpha}(x,y)+C_{U}\ge 0$, the \emph{Neumann Sobolev space} $\mathcal{W}^{-\alpha}_{N}(U)$ is the completion of $\mathcal{S}_{N}(U)$ with respect to
\begin{equation}\label{eq:sobolev-N}
  \langle f,g\rangle_{\mathcal{W}^{-\alpha}_{N}(U)}
  := \int_{U}\!\int_{U}f(x)g(y)\bigl(G^{N,U}_{2\alpha}(x,y)+C_{U}\bigr)\,d\mu(x)\,d\mu(y).
\end{equation}
The value of $C_{U}$ does not affect the topology of $\mathcal{W}^{-\alpha}_{N}(U)$.
\end{definition}

\section{Examples}\label{sec:geom-examples}

The abstract framework of Assumption~\ref{hyp:mmspace} is broad enough to cover a wide variety of geometric structures, ranging from classical Euclidean domains to fractals with anomalous diffusion.  We illustrate this with four families of examples.  The first two are the compact intervals and the metric graphs. They are classical in the sense that $d_{h}=1$ and $d_{w}=2$, so that the heat kernel satisfies standard Gaussian estimates. They differ in their topology.  The third family consists of the Riemannian manifolds, for which $d_{h}=n$ is the dimension and $d_{w}=2$. The curvature only enters through the constants in the heat kernel bounds, and the recurrence condition $d_{h}<d_{w}$ confines Assumption~\ref{hyp:mmspace} to the case $n=1$.  The fourth example is the Sierpi\'nski gasket, which is genuinely fractal. Its walk dimension $d_{w}=\ln 5/\ln 2>2$ reflects the anomalously slow diffusion caused by the fractal geometry, and the sub-Gaussian heat kernel estimates~\eqref{eq:subGauss} are sharp.  In each case we examine Assumption~\ref{hyp:mmspace}, and we identify the dimensional parameters $d_{h}$ and $d_{w}$ governing the kernel bounds of Propositions~\ref{prop:RK-D-bounds} and~\ref{prop:RK-N-bounds}.

\subsection{Compact intervals}\label{ssec:interval}

The simplest example is $X=\R$ with Euclidean distance $d(x,y)=|x-y|$, Lebesgue measure $\mu=dx$, and the standard Dirichlet form $\mathcal{E}(f,g)=\int_{\R}f'g'\,dx$ on $\mathcal{F}=W^{1,2}(\R)$.  The heat kernel is the Gaussian kernel 
\begin{equation*}
p_{t}(x,y)=(4\pi t)^{-1/2}e^{-(x-y)^{2}/(4t)} 
\end{equation*}
giving $d_{h}=1$ and $d_{w}=2$.  The recurrence condition $d_{h}<d_{w}$ holds.

Let us now focus on compact sets, the prototypical example being $U=(0,1)$ (or any interval $(a,b)$). In this case one can separate two cases:

\begin{enumerate}
\item {\it Dirichlet boundary conditions.}  The form $(\mathcal{E}^{D}_{U},\mathcal{F}^{D}(U))$ becomes
\begin{equation*}
(\mathcal{E}^{D},W^{1,2}_{0}([0,1])),
\quad\text{where}\quad
\mathcal{E}^{D}(f,g)=\int_{0}^{1}f'g'\,dx .
\end{equation*}
Here the Sobolev space $W^{1,2}_{0}([0,1])$ is defined by
\begin{equation*}
W^{1,2}_{0}([0,1]):=\{f\in L^{2}([0,1]):f'\in L^{2}([0,1]),\,f(0)=f(1)=0\} .
\end{equation*}
It is the concrete instance of the Dirichlet domain $\mathcal{F}^{D}(U)$ (see Proposition~\ref{prop:E-D}).  By a classical Littlewood--Paley characterization of Sobolev spaces, one can identify $W^{1,2}_{0}([0,1])$ isometrically with the Besov space $B^{1}_{2,2}(\R)$. That is in the notation of Definition~\ref{def:besov-Lp}, $B^{1}_{2,2}(\R)$ coincides with $\mathcal{B}^{1}_{2,w_{0}}(\R)$ with trivial weight $w_{0}\equiv 1$.
Moreover, the Dirichlet Laplacian is given by
\begin{equation*}
\Delta_{D}=\frac{d^{2}}{dx^{2}},
\quad\text{with domain}\quad
\{f\in W^{2,2}([0,1]):f(0)=f(1)=0\}.
\end{equation*}
Eventually, the heat kernel and eigenvalues are explicit:
\begin{equation*}
  p^{D}_{t}(x,y) = \frac{2}{\sqrt{\pi t}}\sum_{k\in\Z}\!\left(e^{-(x-y+2k)^{2}/(4t)}-e^{-(x+y+2k)^{2}/(4t)}\right), 
\end{equation*}
and
\begin{equation*}
\lambda_{k}=k^{2}\pi^{2},
\quad\text{with}\quad
\Phi_{k}(x)=\sqrt{2}\sin(k\pi x).
\end{equation*}

\item {\it Neumann boundary conditions.}  The form $(\mathcal{E}^{N}_{U},\mathcal{F}(U))$ becomes
\begin{equation*}
(\mathcal{E}^{N},W^{1,2}([0,1])),
\quad\text{where}\quad
\mathcal{E}^{N}(f,g)=\int_{0}^{1}f'g'\,dx .
\end{equation*}
Here the Sobolev space $W^{1,2}([0,1])$ is defined by
\begin{equation*}
W^{1,2}([0,1]):=\{f\in L^{2}([0,1]):f'\in L^{2}([0,1])\} .
\end{equation*}
It is the concrete instance of the Neumann domain $\mathcal{F}(U)$ (see Definition~\ref{def:neumann-form}); unlike the Dirichlet case, no vanishing condition is imposed at the boundary.  By restriction of $B^{1}_{2,2}(\R)$ to $[0,1]$, one has $W^{1,2}([0,1])=\mathcal{B}^{1}_{2,w_{0}}([0,1])$ with trivial weight $w_{0}\equiv 1$ (Definition~\ref{def:besov-Lp}).
Moreover, the Neumann Laplacian is given by
\begin{equation*}
\Delta_{N}=\frac{d^{2}}{dx^{2}},
\quad\text{with domain}\quad
\{f\in W^{2,2}([0,1]):f'(0)=f'(1)=0\}.
\end{equation*}
The heat kernel and eigenvalues are explicit:
\begin{equation*}
  p^{N}_{t}(x,y) = \frac{1}{\sqrt{\pi t}}\sum_{k\in\Z}\!\left(e^{-(x-y+2k)^{2}/(4t)}+e^{-(x+y+2k)^{2}/(4t)}\right),
\end{equation*}
and
\begin{equation*}
\nu_{k}=k^{2}\pi^{2},
\quad\text{with}\quad
\Psi_{k}(x)=\sqrt{2}\cos(k\pi x).
\end{equation*}
The PAM on $(0,1)$ with these boundary conditions has been studied in~\cite{foondun:09:intermittence,foondun:15:behaviour}.
\end{enumerate}

\subsection{Metric graphs}\label{ssec:metric-graphs}

Metric graphs generalize the interval example of Section~\ref{ssec:interval} by gluing finitely many one-dimensional edges and half-lines at a finite set of vertices.  Each edge is a copy of a compact interval, so that the local geometry remains Euclidean. The heat kernel exponents are $d_{h}=1$ and $d_{w}=2$, and the recurrence condition $d_{h}<d_{w}$ holds.  The new ingredient is a \emph{Kirchhoff vertex condition} at each junction. It encodes flux conservation, and it singles out the correct self-adjoint extension of the edge-wise Laplacian $-f''$.

\begin{definition}[Metric graph]\label{def:metric-graph}
A \emph{metric graph} $\mathbf{G}=(V,\mathbf{E},\mathbf{R})$ consists of:
\begin{enumerate}
  \item a finite set of vertices $V$;
  \item a finite set of edges $\mathbf{E}$, each identified with a closed interval $I_{e}=[0,r(e)]$ of finite length $r(e)>0$, with two endpoints $e^{-},e^{+}\in V$;
  \item a finite set of rays $\mathbf{R}$, each identified with a half-line $[0,\infty)$, with one finite endpoint $e^{-}\in V$.
\end{enumerate}
The graph $\mathbf{G}$ is assumed connected.  It is equipped with the geodesic distance $d$ induced by the Euclidean lengths of edges and rays, and with the measure $\mu$ that restricts to Lebesgue measure on each edge and ray.
\end{definition}

\begin{figure}[htb]
\centering
\begin{tikzpicture}[line width=2pt, >=stealth]

  \colorlet{cOne}{blue!70!black}
  \colorlet{cTwo}{red!65!black}
  \colorlet{cThree}{ForestGreen}
  \colorlet{cFour}{violet!80!black}
  \colorlet{cRay}{gray!65!black}

  \fill (0,0)     circle (4pt) node[below=6pt, font=\small] {$v_{1}$};
  \fill (4.5,0)   circle (4pt) node[below=6pt, font=\small] {$v_{2}$};
  \fill (5.5,2.3) circle (4pt) node[right=5pt,  font=\small] {$v_{3}$};
  \fill (0.8,2.6) circle (4pt) node[above=5pt,  font=\small] {$v_{4}$};

  \draw[cOne] (0,0) -- (4.5,0);
  \node[above, font=\small, cOne] at (2.25,0.08) {$e_{1}$};

  \draw[<->, thin, black] (0,-0.75) -- node[below, font=\small]{$r(e_{1})$} (4.5,-0.75);
  \draw[thin, black] (0,-0.35) -- (0,-0.9);
  \draw[thin, black] (4.5,-0.35) -- (4.5,-0.9);

  \draw[cTwo] (4.5,0) -- (5.5,2.3);
  \node[right=2pt, font=\small, cTwo] at (5.2,1.15) {$e_{2}$};

  \draw[cThree] (0,0) -- (0.8,2.6);
  \node[left=2pt, font=\small, cThree] at (0.25,1.3) {$e_{3}$};

  \draw[cFour] (0.8,2.6) -- (5.5,2.3);
  \node[above=3pt, font=\small, cFour] at (3.15,2.55) {$e_{4}$};

  \draw[cRay, ->] (4.5,0) -- (6.8,-1.8) node[right, font=\small, cRay] {$+\infty$};
  \node[right=2pt, font=\small, cRay] at (5.85,-0.7) {$r_{1}$};

\end{tikzpicture}
\caption{A metric graph $\mathbf{G}=(V,\mathbf{E},\mathbf{R})$ with vertices $V=\{v_{1},v_{2},v_{3},v_{4}\}$, edges $\mathbf{E}=\{e_{1},e_{2},e_{3},e_{4}\}$, and one ray $r_{1}\in\mathbf{R}$ based at $v_{2}$.  The length $r(e_{1})$ of $e_{1}$ is indicated; the ray $r_{1}$ extends to $+\infty$.}
\label{fig:metric-graph}
\end{figure}
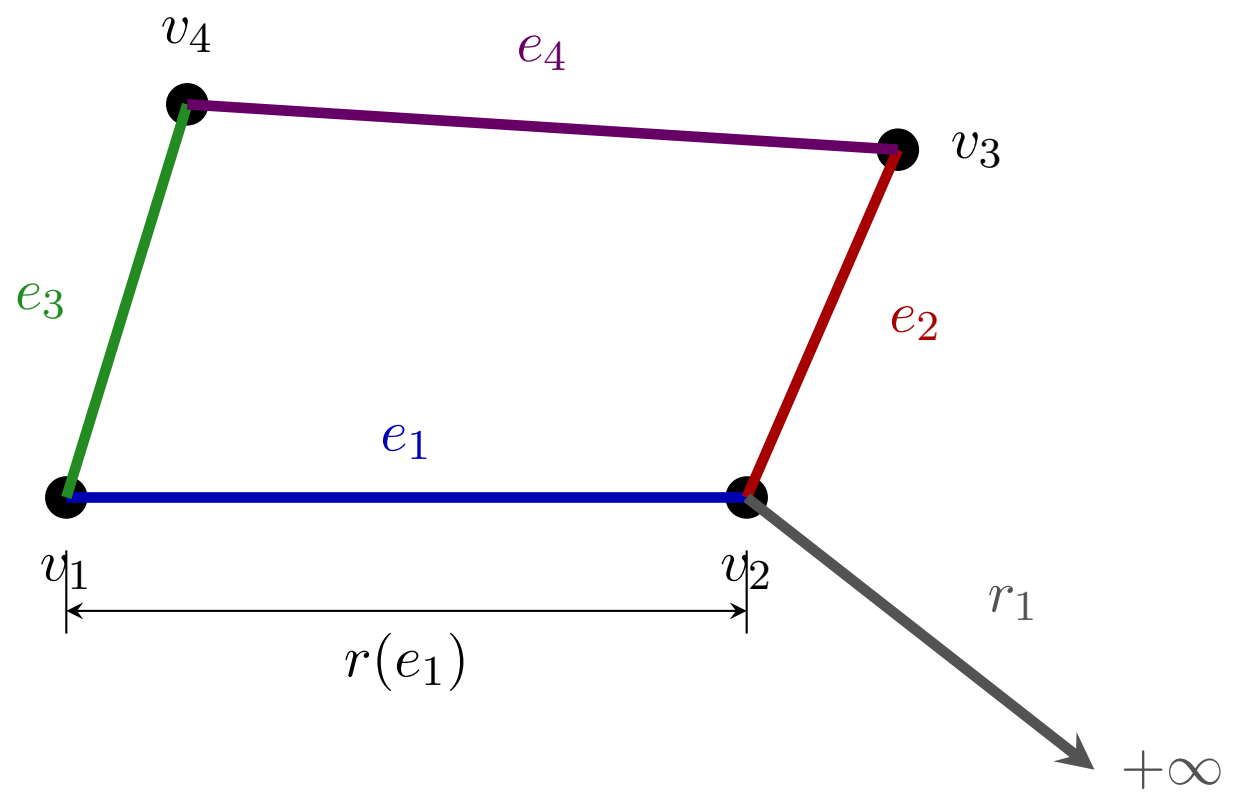

Before defining function spaces, we fix notation for the edges incident at each vertex and their orientation, which will appear systematically in the vertex conditions.

\begin{notation}[Adjacent edges and orientation]\label{not:graph-orient}
Let $\mathbf{G}=(V,\mathbf{E},\mathbf{R})$  be a metric graph. Then
\begin{enumerate}
  \item For each vertex $v\in V$, the set of \emph{adjacent edges} is
\begin{equation*}
\mathbf{E}_{v} = \bigl\{e\in\mathbf{E}\cup\mathbf{R} : v=e^{-}\text{ or }v=e^{+}\bigr\}.
\end{equation*}
  \item The \emph{orientation function} $U_{v}:\mathbf{E}_{v}\to\{-1,1\}$ is defined by $U_{v}(e)=1$ if $v=e^{-}$ and $U_{v}(e)=-1$ if $v=e^{+}$.
\end{enumerate}
\end{notation}

With this orientation in place, we define the natural function spaces on $\mathbf{G}$, assembled edge-wise from the standard $L^{2}$ and $H^{1}$ spaces on each interval $I_{e}$, with continuity enforced at vertices.

\begin{definition}[Function spaces on $\mathbf{G}$]\label{def:fspaces-graph}
Let $\mathbf{G}=(V,\mathbf{E},\mathbf{R})$ be a metric graph (Definition~\ref{def:metric-graph}).  Functions on $\mathbf{G}$ are families $f=(f_{e})_{e\in\mathbf{E}\cup\mathbf{R}}$ with $f_{e}:I_{e}\to\R$.
\begin{enumerate}
  \item The space $L^{2}(\mathbf{G},\mu)$ decomposes as
\begin{equation}\label{eq:L2-graph}
L^{2}(\mathbf{G},\mu) = \bigoplus_{e\in\mathbf{E}\cup\mathbf{R}} L^{2}(I_{e}).
\end{equation}
  \item A function $f$ is \emph{continuous on $\mathbf{G}$} if $f_{e}\in C(I_{e})$ for all $e$ and $f_{e_{1}}(v)=f_{e_{2}}(v)$ whenever $e_{1},e_{2}\in\mathbf{E}_{v}$ (see Notation~\ref{not:graph-orient}\,(i)).
  \item The \emph{Sobolev space} $H^{1}_{0}(\mathbf{G})$ is the natural domain for Dirichlet forms on $\mathbf{G}$:
\begin{multline}\label{eq:H1-graph}
H^{1}_{0}(\mathbf{G}) = \bigl\{f=(f_{e})_{e} : f_{e}\in H^{1}(I_{e})\text{ for all }e\in\mathbf{E}\cup\mathbf{R},\; \\
\text{ and } f\text{ continuous at every vertex}\bigr\}.
\end{multline}
  \item For $f=(f_{e})_{e}$ with $f_{e}\in H^{1}(I_{e})$, the \emph{edge derivative at a vertex} $f'_{e}(v)$ is the trace of $f'_{e}$ at the endpoint of $I_{e}$ corresponding to $v$:
\begin{equation}\label{eq:trace-graph}
f'_{e}(v) = \begin{cases} f'_{e}(0) & \text{if }v=e^{-}, \\ f'_{e}(r(e)) & \text{if }v=e^{+}. \end{cases}
\end{equation}
  \item The \emph{inward-pointing normal derivative} of $f$ along edge $e$ at vertex $v$ is $U_{v}(e)f'_{e}(v)$, with $U_{v}$ as in Notation~\ref{not:graph-orient}\,(ii) and $f'_{e}(v)$ as in~(iv).
\end{enumerate}
\end{definition}

The Dirichlet form on $\mathbf{G}$ is obtained by summing the edge-wise $H^{1}$ energies. The vertex continuity condition in~\eqref{eq:H1-graph} ensures that $H^{1}_{0}(\mathbf{G})$ is a Hilbert space, and it makes the sum below well defined.

\begin{definition}[Kirchhoff Dirichlet form]\label{def:kirchhoff-form}
Let $\mathbf{G}=(V,\mathbf{E},\mathbf{R})$ be a metric graph.  The \emph{Kirchhoff Dirichlet form} on $L^{2}(\mathbf{G},\mu)$ (see Definition~\ref{def:fspaces-graph}\,(i)) is
\begin{equation}\label{eq:E-graph}
\mathcal{E}(f,g) = \sum_{e\in\mathbf{E}\cup\mathbf{R}}\int_{I_{e}}f'_{e}(x)g'_{e}(x)\,dx,
\qquad f,g\in H^{1}_{0}(\mathbf{G}),
\end{equation}
with domain $H^{1}_{0}(\mathbf{G})$ as in Definition~\ref{def:fspaces-graph}\,(iii).
\end{definition}

The main analytical properties of $\mathcal{E}$ (its generator, the precise form of the vertex condition, and the heat kernel bounds that place $\mathbf{G}$ in the framework of Section~\ref{sec:geom-examples}) are recorded in the following theorem.

\begin{theorem}\label{thm:kirchhoff-props}
The Kirchhoff Dirichlet form $(\mathcal{E},H^{1}_{0}(\mathbf{G}))$ of Definition~\ref{def:kirchhoff-form} is a strongly local regular Dirichlet form on $L^{2}(\mathbf{G},\mu)$.  Moreover:
\begin{enumerate}
  \item {\it Generator.}  The generator of $\mathcal{E}$ is $\Delta f=-f''$, acting edge-wise, with domain
\begin{multline}\label{eq:dom-Delta-graph}
\mathrm{Dom}(\Delta) = \bigl\{f\in H^{1}_{0}(\mathbf{G}) : f''_{e}\in L^{2}(I_{e})\text{ for all }e\in\mathbf{E}\cup\mathbf{R},\;\; \\
\text{ and }
\textstyle\sum_{e\in\mathbf{E}_{v}}U_{v}(e)f'_{e}(v)=0\text{ for all }v\in V\bigr\},
\end{multline}
where $\mathbf{E}_{v}$, $U_{v}$ are as in Notation~\ref{not:graph-orient} and $f'_{e}(v)$ is as in Definition~\ref{def:fspaces-graph}\,(iv).  The condition $f''_{e}\in L^{2}(I_{e})$ ensures $f'_{e}\in H^{1}(I_{e})$, so the traces in~\eqref{eq:trace-graph} are well defined.
  \item {\it Heat kernel bounds.}  Since $V$, $\mathbf{E}$ and $\mathbf{R}$ are finite and every edge has finite positive length, $\mu(B(x,r))\asymp r$ uniformly in $x$ and $r$, and $\mathbf{G}$ carries a Poincar\'e inequality.  The criterion of~\cite{grigoryan:12:two-sided}, in the case $V(x,r)\asymp r^{d_{h}}$ and $R(t)\asymp t^{1/d_{w}}$, then shows that the heat kernel $p_{t}(x,y)$ associated with $\mathcal{E}$ satisfies the Gaussian estimates~\eqref{eq:subGauss} with $d_{h}=1$ and $d_{w}=2$.
\end{enumerate}
\end{theorem}

\begin{remark}\label{rem:kirchhoff}
A few comments on Theorem \ref{thm:kirchhoff-props} are in order:
\begin{enumerate}
  \item {\it Conservation of flux.}  The Kirchhoff condition in~\eqref{eq:dom-Delta-graph} requires the sum of inward-pointing normal derivatives (Definition~\ref{def:fspaces-graph}\,(v)) to vanish at every vertex $v\in V$.  This is the graph analogue of the Neumann condition $\partial_{n}f=0$ on a manifold with boundary.
  \item {\it Independence of orientation.}  Neither $\mathcal{E}$ nor $\Delta$ depends on the choice of orientation.  If the sign convention for $U_{v}(e)$ is reversed for some edge $e$, then $f'_{e}(v)$ changes sign accordingly by~\eqref{eq:trace-graph}, leaving both the bilinear form~\eqref{eq:E-graph} and the vertex condition in~\eqref{eq:dom-Delta-graph} unchanged.
\end{enumerate}
\end{remark}

Connected bounded open subsets of $\mathbf{G}$ are uniform domains, so that both the Dirichlet and Neumann setups of Sections~\ref{sssec:dirichlet-semi}--\ref{sssec:neumann-frac} apply. We refer to~\cite{post:12:spectral} for details on function spaces on metric graphs.

\subsection{Riemannian manifolds}\label{ssec:manifold}

Riemannian manifolds generalize $\R^{n}$ by allowing non-trivial curvature while retaining the local Euclidean structure that makes classical analysis available.  The heat kernel exponents remain $d_{h}=n$ (the manifold dimension) and $d_{w}=2$ (standard Gaussian diffusion), exactly as for flat $\R^{n}$; the curvature governs the volume growth of geodesic balls and enters through the constants in the heat kernel bounds.  Non-negative Ricci curvature is the standard assumption for these estimates, and includes the sphere $\mathbb{S}^{n}$ and flat tori $\mathbb{T}^{n}$ as prototypical examples.

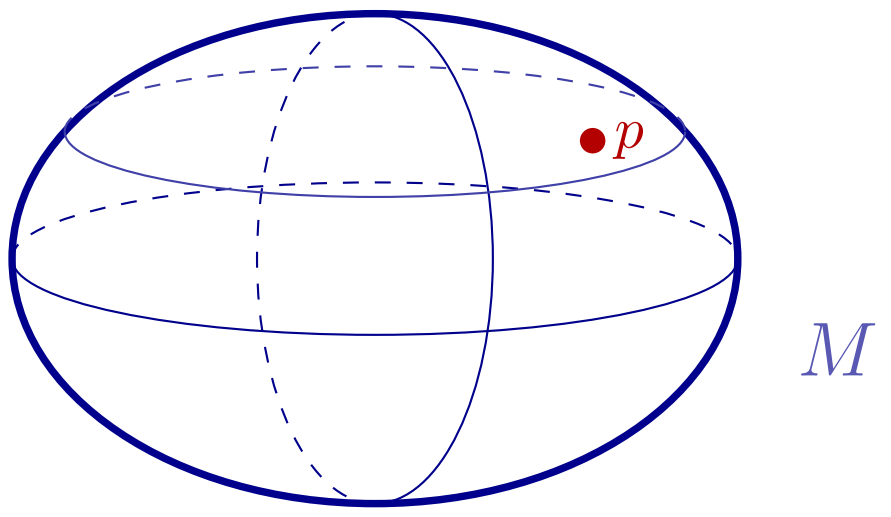
\begin{figure}[h]
\centering
\begin{tikzpicture}[scale=1.2]
\colorlet{surfcol}{blue!55!black}
\colorlet{ptcol}{red!70!black}


\draw[surfcol, line width=1.4pt] (0,0) ellipse (2.0 and 1.35);

\draw[surfcol] (-2,0) arc (180:360:2.0 and 0.42);
\draw[surfcol, dashed] (-2,0) arc (180:0:2.0 and 0.42);

\draw[surfcol] (0,-1.35) arc (-90:90:0.65 and 1.35);
\draw[surfcol, dashed] (0,1.35) arc (90:270:0.65 and 1.35);

\draw[surfcol!75] (-1.71,0.7) arc (180:360:1.71 and 0.36);
\draw[surfcol!75, dashed] (-1.71,0.7) arc (180:0:1.71 and 0.36);

\fill[ptcol] (1.2, 0.65) circle (2pt);
\node[right, ptcol, font=\small] at (1.2, 0.65) {$p$};

\node[surfcol!65, font=\large\itshape] at (2.55, -0.5) {$M$};

\end{tikzpicture}
\caption{A compact Riemannian manifold $(M,g)$ without boundary, with positive
sectional curvature (an ellipsoid with three distinct semi-axes), where $p\in M$ is a generic point.}
\label{fig:riem-surface}
\end{figure}

\begin{definition}[Riemannian manifold]\label{def:riem-manifold}
A \emph{Riemannian manifold} $(M,g)$ consists of a connected smooth manifold $M$ of dimension $n\ge 1$ equipped with a Riemannian metric $g$.  We equip $M$ with:
\begin{enumerate}
  \item The \emph{Riemannian measure} $\mu=\mathrm{vol}_{g}$, given in local coordinates by 
  \begin{equation*}
d\mu=\sqrt{\det(g_{ij})}\,dx^{1}\cdots dx^{n} .
\end{equation*}

  \item The \emph{Riemannian distance}, defined by 
  \begin{equation*}
d(x,y)=\inf_{\gamma}\int_{0}^{1}|\dot{\gamma}(t)|_{g}\,dt \, ,
\end{equation*}
where the infimum runs over smooth curves $\gamma:[0,1]\to M$ with $\gamma(0)=x$ and $\gamma(1)=y$;
\end{enumerate}
\end{definition}

The derivative structure on $M$ is encoded by the Riemannian gradient, which extends the Euclidean gradient $\partial_{i}f$ to curved geometry and gives rise to the Sobolev space serving as the natural domain for the Dirichlet form.

\begin{definition}[Riemannian gradient and Sobolev space]\label{def:riem-gradient}
Let $(M,g)$ be a Riemannian manifold.
\begin{enumerate}
  \item The \emph{Riemannian gradient} $\nabla_{g}f$ of $f\in C^{\infty}(M)$ is the vector field characterized by $\langle\nabla_{g}f,X\rangle_{g}=df(X)$ for every smooth vector field $X$ on $M$.  In local coordinates $(\nabla_{g}f)^{i}=g^{ij}\partial_{j}f$, so that
\begin{equation}\label{eq:grad-sq}
|\nabla_{g}f|_{g}^{2} = g^{ij}\partial_{i}f\,\partial_{j}f.
\end{equation}
  \item The \emph{Sobolev space} $W^{1,2}(M)$ is the completion of $C^{\infty}_{c}(M)$ under the norm
\begin{equation}\label{eq:W12-manifold}
\|f\|_{W^{1,2}(M)}^{2} = \int_{M}\bigl(f^{2}+|\nabla_{g}f|_{g}^{2}\bigr)\,d\mu,
\end{equation}
with $|\nabla_{g}f|_{g}^{2}$ as in~\eqref{eq:grad-sq}.
\end{enumerate}
\end{definition}

The Riemannian Dirichlet form is the natural generalization of $\int_{\R^{n}}|\nabla f|^{2}\,dx$: it integrates the pointwise squared gradient length against the Riemannian measure.

\begin{definition}[Riemannian Dirichlet form]\label{def:riem-form}
Let $(M,g)$ be a complete Riemannian manifold.  The \emph{Riemannian Dirichlet form} on $L^{2}(M,\mu)$ is
\begin{equation}\label{eq:E-riem}
\mathcal{E}(f,h) = \int_{M}\langle\nabla_{g}f,\nabla_{g}h\rangle_{g}\,d\mu,
\qquad f,h\in W^{1,2}(M),
\end{equation}
with domain $W^{1,2}(M)$ as in Definition~\ref{def:riem-gradient}\,(ii).
\end{definition}

The Riemannian Dirichlet form is strongly local and regular, and its generator is the Laplace--Beltrami operator.  This is a classical result of functional analysis on manifolds, and we refer to~\cite{grigoryan:09:heat} for a detailed treatment.

\begin{proposition}\label{prop:riem-form}
Let $(M,g)$ be a complete Riemannian manifold of dimension $n$.  The form $(\mathcal{E},W^{1,2}(M))$ of Definition~\ref{def:riem-form} is a strongly local regular Dirichlet form on $L^{2}(M,\mu)$.  Its generator is the \emph{Laplace--Beltrami operator}
\begin{equation}\label{eq:LB-def}
\Delta_{M}f = -\mathrm{div}_{g}(\nabla_{g}f)
= -\frac{1}{\sqrt{|g|}}\,\partial_{i}\!\left(\sqrt{|g|}\,g^{ij}\partial_{j}f\right),
\end{equation}
with domain
\begin{equation}\label{eq:dom-LB}
\mathrm{Dom}(\Delta_{M}) = \bigl\{f\in W^{1,2}(M) : \mathrm{div}_{g}(\nabla_{g}f)\in L^{2}(M,\mu)\bigr\}.
\end{equation}
The restriction of $\Delta_{M}$ to $C^{\infty}_{c}(M)$ is essentially self-adjoint on $L^{2}(M,\mu)$.
\end{proposition}

The generator is now identified.  The question of how the associated heat semigroup $e^{t\Delta_{M}}$ spreads mass is governed by the geometry of $(M,g)$, and specifically by its curvature.  Under non-negative Ricci curvature, Li and Yau~\cite{li:86:parabolic} established sharp two-sided Gaussian estimates for the heat kernel, which translate directly into the sub-Gaussian framework of Assumption~\ref{hyp:mmspace} with $d_{h}=n$ and $d_{w}=2$.

\begin{theorem}[Li--Yau]\label{thm:riem-hk}
Let $(M,g)$ be a complete Riemannian manifold of dimension $n$ with $\mathrm{Ric}_{g}\ge 0$.  The heat kernel $p_{t}(x,y)$ of the semigroup generated by $\Delta_{M}$ (Proposition~\ref{prop:riem-form}) satisfies the Gaussian estimates~\eqref{eq:subGauss} with $d_{h}=n$ and $d_{w}=2$~\cite{li:86:parabolic,grigoryan:09:heat}.  In particular, $(M,d_{g},\mu,\mathcal{E})$ fits into the framework of Assumption~\ref{hyp:mmspace} whenever $n<2$, i.e.\ $n=1$.
\end{theorem}

\begin{remark}\label{rem:riem-remarks}
Let us comment on two aspects of the Riemannian setting, namely the status of the recurrence assumption and the prototypical compact examples to which the theorem applies.
\begin{enumerate}
  \item {\it Recurrence.}  The recurrence condition $d_{h}<d_{w}$ of Assumption~\ref{hyp:mmspace} requires $n<2$, hence $n=1$ for a manifold.  For compact Riemannian manifolds of dimension $n\ge 2$ the condition is not satisfied, so the general framework of this book does not apply directly.  The PAM on such manifolds can nonetheless be studied by techniques similar to ours; see~\cite{chen:25:global,chen:26:sharp}.
  \item {\it Compact examples.}  The sphere $\mathbb{S}^{n}$ (with $\mathrm{Ric}_{g}=(n-1)g>0$) and the flat torus $\mathbb{T}^{n}=\R^{n}/\Z^{n}$ (with $\mathrm{Ric}_{g}=0$) are compact manifolds satisfying $\mathrm{Ric}_{g}\ge 0$; their heat kernels satisfy the bounds in Theorem~\ref{thm:riem-hk} for small times $t\le\mathrm{diam}(M)^{2}$.
\end{enumerate}
\end{remark}

\subsection{Sierpi\'nski gasket}\label{ssec:sierpinski}

The Sierpi\'nski gasket is one of the most studied examples of a post-critically finite (p.c.f.) fractal~\cite{kigami:01:analysis,barlow:98:diffusions}.  Unlike the interval and the metric graph, it exhibits a genuinely sub-Gaussian diffusion ($d_{w}=\ln 5/\ln 2>2$), providing a non-trivial verification of Assumption~\ref{hyp:mmspace}.

\begin{definition}[Sierpi\'nski gasket]\label{def:SG}
Let $V_{0}=\{p_{1},p_{2},p_{3}\}$ be the vertices of an equilateral triangle of side~$1$ in~$\C$.  For $i=1,2,3$, define the contraction maps
\begin{equation}\label{eq:SK-maps}
  \mathfrak{f}_{i}(z) = \frac{z-p_{i}}{2}+p_{i}.
\end{equation}
The \emph{Sierpi\'nski gasket} $K$ is the unique non-empty compact set satisfying
\begin{equation}\label{eq:SK-attractor}
  K = \bigcup_{i=1}^{3}\mathfrak{f}_{i}(K).
\end{equation}
The set $\partial K:=V_{0}$ is called the \emph{boundary} of~$K$.  We equip $K$ with the Euclidean distance $d(x,y)=|x-y|$.
\end{definition}

\newcommand{\SGfill}[5]{%
  \begingroup
  \ifnum#4=0%
    \fill[#5]
      (#1,#2)
      -- ({(#1)+(#3)},{#2})
      -- ({(#1)+0.5*(#3)},{(#2)+0.8660254*(#3)})
      -- cycle;%
  \else%
    \pgfmathsetmacro\sgfH{(#3)*0.5}%
    \pgfmathtruncatemacro\sgfD{(#4)-1}%
    \pgfmathsetmacro\sgfRX{(#1)+\sgfH}%
    \pgfmathsetmacro\sgfTX{(#1)+(#3)*0.25}%
    \pgfmathsetmacro\sgfTY{(#2)+(#3)*0.4330127}%
    \edef\sgfA{\noexpand\SGfill{#1}{#2}{\sgfH}{\sgfD}{#5}}%
    \edef\sgfB{\noexpand\SGfill{\sgfRX}{#2}{\sgfH}{\sgfD}{#5}}%
    \edef\sgfC{\noexpand\SGfill{\sgfTX}{\sgfTY}{\sgfH}{\sgfD}{#5}}%
    \sgfA\sgfB\sgfC%
  \fi%
  \endgroup
}
\colorlet{SGcolorBL}{teal!65}
\colorlet{SGcolorBR}{orange!80}
\colorlet{SGcolorTop}{violet!55}

\begin{center}
\begin{tikzpicture}[scale=1]
  \SGfill{0}{0}{3}{3}{SGcolorBL}              
  \SGfill{3}{0}{3}{3}{SGcolorBR}              
  \SGfill{1.5}{2.598076}{3}{3}{SGcolorTop}    
  \draw[black!30, line width=0.4pt] (0,0) -- (6,0) -- (3,5.196152) -- cycle;
  \fill (0,0)        circle (2pt) node[below left,  font=\small] {$p_{1}$};
  \fill (6,0)        circle (2pt) node[below right, font=\small] {$p_{2}$};
  \fill (3,5.196152) circle (2pt) node[above,       font=\small] {$p_{3}$};
\end{tikzpicture}
\end{center}

The self-similar structure of $K$ determines a canonical measure and a pair of dimensional exponents.

\begin{definition}[Hausdorff measure and dimensions of $K$]\label{def:SG-measure}
The (normalized) Hausdorff measure $\mu$ on $K$ is the unique Borel measure satisfying
\begin{equation}\label{eq:SG-mu}
  \mu\bigl(\mathfrak{f}_{i_{1}}\circ\cdots\circ\mathfrak{f}_{i_{n}}(K)\bigr)=3^{-n}
\end{equation}
for every finite word $(i_{1},\ldots,i_{n})\in\{1,2,3\}^{n}$.  The measure $\mu$ is Ahlfors $d_{h}$-regular in the sense of~\eqref{eq:ahlfors}, with dimensional exponents
\begin{equation}\label{eq:gasket-dims}
  d_{h} = \frac{\ln 3}{\ln 2}, \qquad d_{w} = \frac{\ln 5}{\ln 2}.
\end{equation}
Since $d_{h}\approx 1.585<d_{w}\approx 2.322$, the recurrence condition of Assumption~\ref{hyp:mmspace} holds.
\end{definition}

\begin{remark}[Intuition for $d_{h}=\ln 3/\ln 2$]\label{rem:gasket-dh}
The gasket consists of exactly $3$ copies of itself, each contracted by factor $1/2$ (the maps $\mathfrak{f}_{i}$ in~\eqref{eq:SK-maps}).  For a self-similar set made of $N$ copies at scaling ratio $r$, the Hausdorff dimension satisfies $N\cdot r^{d_{h}}=1$.  With $N=3$ and $r=1/2$:
\begin{equation}\label{eq:dh-heuristic}
  3\cdot\left(\frac{1}{2}\right)^{d_{h}}=1,\qquad\text{hence}\qquad d_{h}=\frac{\ln 3}{\ln 2}.
\end{equation}
Concretely, at level $m$ the gasket is covered by $3^{m}$ cells of diameter $2^{-m}$ and mass $3^{-m}$, so a ball of radius $r\asymp 2^{-m}$ satisfies $\mu(B(x,r))\asymp 3^{-m}\asymp r^{d_{h}}$, consistently with~\eqref{eq:ahlfors}.
\end{remark}

The construction of the Dirichlet form on $K$ proceeds by passing to the limit along a sequence of finite graphs.  We first set up the necessary notation.

\begin{notation}[Discrete approximation of $K$]\label{not:SG-approx}
Starting from $V_{0}=\{p_{1},p_{2},p_{3}\}$, define inductively
\begin{equation}\label{eq:SG-Vm}
  V_{m+1}=\bigcup_{i=1}^{3}\mathfrak{f}_{i}(V_{m}), \qquad m\ge 0.
\end{equation}
One has $\#V_{m}=\frac{3(3^{m}+1)}{2}$.  Set $V_{*}=\bigcup_{m\ge 0}V_{m}$ and $V_{*}^{0}=V_{*}\setminus V_{0}$.  For each $m\ge 0$, the \emph{pre-gasket graph} $G_{m}$ has vertex set $V_{m}$ and edges between all pairs at Euclidean distance $2^{-m}$; write $p\sim_{m} q$ for adjacency in $G_{m}$.
\end{notation}

\newcommand{\drawGm}[7]{%
  \begingroup
  \ifnum#7=0%
    \draw[gmEdge, line width=\gmEW]
      (#1,#2) -- (#3,#4) -- (#5,#6) -- cycle;%
    \fill[gmVtx] (#1,#2) circle (\gmVR);%
    \fill[gmVtx] (#3,#4) circle (\gmVR);%
    \fill[gmVtx] (#5,#6) circle (\gmVR);%
  \else%
    \pgfmathsetmacro\gmMabX{0.5*(#1)+0.5*(#3)}%
    \pgfmathsetmacro\gmMabY{0.5*(#2)+0.5*(#4)}%
    \pgfmathsetmacro\gmMacX{0.5*(#1)+0.5*(#5)}%
    \pgfmathsetmacro\gmMacY{0.5*(#2)+0.5*(#6)}%
    \pgfmathsetmacro\gmMbcX{0.5*(#3)+0.5*(#5)}%
    \pgfmathsetmacro\gmMbcY{0.5*(#4)+0.5*(#6)}%
    \pgfmathtruncatemacro\gmND{(#7)-1}%
    \edef\gmCA{\noexpand\drawGm{#1}{#2}{\gmMabX}{\gmMabY}{\gmMacX}{\gmMacY}{\gmND}}%
    \edef\gmCB{\noexpand\drawGm{\gmMabX}{\gmMabY}{#3}{#4}{\gmMbcX}{\gmMbcY}{\gmND}}%
    \edef\gmCC{\noexpand\drawGm{\gmMacX}{\gmMacY}{\gmMbcX}{\gmMbcY}{#5}{#6}{\gmND}}%
    \gmCA\gmCB\gmCC%
  \fi%
  \endgroup
}
\colorlet{gmEdge}{olive!90!black}
\colorlet{gmVtx}{magenta!80!black}
\def\gmVR{1.8pt}
\def\gmEW{0.7pt}

\begin{center}
\begin{tikzpicture}[scale=1]
  %
  \drawGm{0}{3.2}{2.2}{3.2}{1.1}{5.1053}{0}
  \node[below=2pt, font=\small] at (1.1,3.2) {$G_{0}$};
  %
  \drawGm{3.8}{3.2}{6.0}{3.2}{4.9}{5.1053}{1}
  \node[below=2pt, font=\small] at (4.9,3.2) {$G_{1}$};
  %
  \drawGm{0}{0}{2.2}{0}{1.1}{1.9053}{2}
  \node[below=2pt, font=\small] at (1.1,0) {$G_{2}$};
  %
  \drawGm{3.8}{0}{6.0}{0}{4.9}{1.9053}{3}
  \node[below=2pt, font=\small] at (4.9,0) {$G_{3}$};
\end{tikzpicture}
\end{center}

The Dirichlet form on $K$ (Definition~\ref{def:SG}) is built as a limit of discrete energies defined on the approximating graphs $G_{m}$ of Notation~\ref{not:SG-approx}.  The key point is to choose the correct renormalization factor so that these discrete energies are compatible as $m$ increases.

\begin{definition}[Discrete energy on $G_{m}$]\label{def:E-gasket-m}
Here $K$ is the Sierpi\'nski gasket of Definition~\ref{def:SG} and $G_{m}=(V_{m},\sim_{m})$ are the pre-gasket graphs of Notation~\ref{not:SG-approx}.  For $m\ge 0$ and $f:V_{m}\to\R$, the \emph{discrete energy} of $f$ on $G_{m}$ is
\begin{equation}\label{eq:E-gasket-m}
  \mathcal{E}_{m}(f,f) = \left(\frac{5}{3}\right)^{m}\sum_{\substack{p,q\in V_{m}\\p\sim_{m} q}}(f(p)-f(q))^{2}.
\end{equation}
\end{definition}

The factor $(5/3)^{m}$ in~\eqref{eq:E-gasket-m} is the \emph{self-similar resistance scaling} of the gasket.  It is the unique choice making the sequence $m\mapsto\mathcal{E}_{m}(f|_{V_{m}},f|_{V_{m}})$ non-decreasing for every continuous $f:K\to\R$: since $(5/3)^{m}$ is exactly the ratio that preserves energy under harmonic extension from $V_{m}$ to $V_{m+1}$, the restriction $f|_{V_{m+1}}$ (which is generally not harmonic) satisfies $\mathcal{E}_{m+1}(f|_{V_{m+1}},f|_{V_{m+1}})\ge\mathcal{E}_{m}(f|_{V_{m}},f|_{V_{m}})$; see~\cite[Proposition~3.1.3 and Corollary~3.4.7]{kigami:01:analysis}.  The two following remarks explain where this factor comes from: Remark~\ref{rem:dirichlet-electric} sets up the resistor-network language needed for the computation, and Remark~\ref{rem:gasket-renorm} carries out the explicit resistance calculation on $G_{1}$ that identifies $5/3$ as the correct ratio.

\begin{remark}[Discrete Dirichlet forms and resistor networks]\label{rem:dirichlet-electric}
There is a classical correspondence between discrete Dirichlet forms and electrical resistor networks~\cite{kigami:01:analysis,barlow:98:diffusions}, which underlies Remarks~\ref{rem:gasket-renorm} and~\ref{rem:gasket-dw}.  On a finite graph $G=(V,E)$ with unit edge resistances, let $f:V\to\R$ be any function, interpreted as a voltage distribution on the nodes.  By Ohm's law, the current on edge $\{p,q\}$ is $I_{pq}=f(p)-f(q)$, so the total power dissipated is
\begin{equation}\label{eq:dirichlet-electric}
  P(f)=\sum_{\{p,q\}\in E}I_{pq}^{2}=\sum_{\{p,q\}\in E}(f(p)-f(q))^{2}.
\end{equation}
The quantity $P(f)$ is the discrete Dirichlet form on $G$ with unit conductances. On the pre-gasket graph $G_{m}$, it is exactly the unrenormalized energy $\sum_{p\sim_{m}q}(f(p)-f(q))^{2}$, that is the summand in~\eqref{eq:E-gasket-m} before multiplication by $(5/3)^{m}$.

The unique voltage distribution $h$ minimizing $P(h)$ with prescribed boundary values $h(a)=1$ and $h(b)=0$ is the solution to the discrete Laplace equation (Kirchhoff's current law: total current into each node is zero)
\begin{equation}\label{eq:kirchhoff-law}
  \Delta_{G}h(p):=\sum_{q\sim p}(h(q)-h(p))=0, \qquad p\in V\setminus\{a,b\},
\end{equation}
supplemented by the boundary conditions $h(a)=1$, $h(b)=0$.  This is a linear system in the values $\{h(p)\}_{p\in V\setminus\{a,b\}}$ with a unique solution (the matrix $-\Delta_{G}$ restricted to interior nodes is positive definite).
The \emph{effective resistance} between $a$ and $b$ is
\begin{equation*}
R_{\mathrm{eff}}(a,b)=\frac{h(a)-h(b)}{I_{a}} ,
\quad\text{where}\quad
I_{a}=\sum_{q\sim a}(h(a)-h(q))
\end{equation*}
is the current injected at~$a$.  In order to give another expression for $R_{\mathrm{eff}}(a,b)$, let us compute $P(h)$. To this aim, for every edge write
\begin{equation*}
(h(p)-h(q))^{2}=h(p)(h(p)-h(q))+h(q)(h(q)-h(p)) , 
\end{equation*}
and sum over all directed pairs $(p,q)$ with $p\sim q$. Regrouping by base node gives
\begin{equation}\label{h1}
P(h)=\sum_{p\in V}h(p)\sum_{q\sim p}(h(p)-h(q)) .
\end{equation}
Using~\eqref{eq:kirchhoff-law} at all nodes $p\in V\setminus\{a,b\}$, the only surviving terms in~\eqref{h1} are $h(a)I_{a}+h(b)(-I_{a})=I_{a}$, giving
\begin{equation}\label{eq:eff-resist}
  R_{\mathrm{eff}}(a,b)=\frac{1}{P(h)}.
\end{equation}
The function $h$ also admits a probabilistic interpretation. Namely, for the simple random walk on $G$ we have $h(p)=P_{p}(\tau_{a}<\tau_{b})$, that is the probability of hitting $a$ before $b$ when started at $p$. This follows from the fact that the hitting probability satisfies~\eqref{eq:kirchhoff-law} at every interior node, which is the same equation characterizing $h$.  The effective resistance is then linked to the random walk by the \emph{commute-time formula}~\cite{barlow:98:diffusions}. Namely, denoting by $|E|$ the number of edges of $G$,
\begin{equation}\label{eq:commute-time}
  E_{a}[\tau_{b}]+E_{b}[\tau_{a}] = 2\,|E|\,R_{\mathrm{eff}}(a,b).
\end{equation}
Larger resistance thus means a longer expected round-trip time between $a$ and $b$.
\end{remark}

\begin{remark}[Origin of the renormalization factor $5/3$]\label{rem:gasket-renorm}
By Remark~\ref{rem:dirichlet-electric}, resistor networks obey the standard rules of circuit theory: resistances in series add ($R_{\mathrm{series}}=R_{1}+R_{2}$), resistances in parallel combine as $R_{\mathrm{par}}^{-1}=R_{1}^{-1}+R_{2}^{-1}$, and Ohm's law $I=V/R$ holds on each edge.  We apply these rules to $G_{0}$ and $G_{1}$, each with unit resistance on every edge.  In $G_{0}$ (a triangle with three unit-resistance edges), the effective resistance between any two corner vertices is $R_{0}=2/3$ (two parallel paths of resistance $1$ and $2$).  To find the effective resistance in $G_{1}$, assign voltages $V(p_{1})=1$, $V(p_{2})=0$.

\begin{center}
\begin{tikzpicture}[line width=1.5pt, font=\small]
  \colorlet{cEdge}{blue!70!black}
  \colorlet{cVtx}{red!65!black}

  \coordinate (p1) at (0, 0);
  \coordinate (p2) at (4, 0);
  \coordinate (p3) at (2, {2*sqrt(3)});

  \draw[cEdge] (p1) -- (p2);
  \draw[cEdge] (p1) -- (p3);
  \draw[cEdge] (p2) -- (p3);

  \node[cEdge, below=4pt] at (2, 0)         {$1$};
  \node[cEdge, left=4pt]  at (1, {sqrt(3)}) {$1$};
  \node[cEdge, right=4pt] at (3, {sqrt(3)}) {$1$};

  \fill[cVtx] (p1) circle (3.5pt);
  \fill[cVtx] (p2) circle (3.5pt);
  \fill[cVtx] (p3) circle (3.5pt);

  \node[cVtx, below left=2pt]  at (p1) {$p_{1}$};
  \node[cVtx, below right=2pt] at (p2) {$p_{2}$};
  \node[cVtx, above=3pt]       at (p3) {$p_{3}$};

  \node[font=\small\bfseries] at (0, {2*sqrt(3)+0.3}) {$G_{0}$};
\end{tikzpicture}
\hspace{1.5cm}
\begin{tikzpicture}[line width=1.5pt, font=\small]
  \colorlet{cEdge}{blue!70!black}
  \colorlet{cVtx}{red!65!black}

  \coordinate (p1)  at (0, 0);
  \coordinate (p2)  at (4, 0);
  \coordinate (p3)  at (2, {2*sqrt(3)});
  \coordinate (q12) at (2, 0);
  \coordinate (q13) at (1, {sqrt(3)});
  \coordinate (q23) at (3, {sqrt(3)});

  \draw[cEdge] (p1)  -- (q12);
  \draw[cEdge] (p1)  -- (q13);
  \draw[cEdge] (q12) -- (q13);
  \draw[cEdge] (q12) -- (p2);
  \draw[cEdge] (q12) -- (q23);
  \draw[cEdge] (p2)  -- (q23);
  \draw[cEdge] (q13) -- (q23);
  \draw[cEdge] (q13) -- (p3);
  \draw[cEdge] (q23) -- (p3);

  \node[cEdge, below=4pt]  at (1, 0)               {$1$};
  \node[cEdge, below=4pt]  at (3, 0)               {$1$};
  \node[cEdge, left=4pt]   at (0.5, {sqrt(3)/2})   {$1$};
  \node[cEdge, left=4pt]   at (1.5, {3*sqrt(3)/2}) {$1$};
  \node[cEdge, right=4pt]  at (3.5, {sqrt(3)/2})   {$1$};
  \node[cEdge, right=4pt]  at (2.5, {3*sqrt(3)/2}) {$1$};
  \node[cEdge, below left=2pt]  at (1.5, {sqrt(3)/2}) {$1$};
  \node[cEdge, below right=2pt] at (2.5, {sqrt(3)/2}) {$1$};
  \node[cEdge, above=3pt]       at (2,   {sqrt(3)})   {$1$};

  \fill[cVtx] (p1)  circle (3.5pt);
  \fill[cVtx] (p2)  circle (3.5pt);
  \fill[cVtx] (p3)  circle (3.5pt);
  \fill[cVtx] (q12) circle (3.5pt);
  \fill[cVtx] (q13) circle (3.5pt);
  \fill[cVtx] (q23) circle (3.5pt);

  \node[cVtx, below left=2pt]  at (p1)  {$p_{1}$};
  \node[cVtx, below right=2pt] at (p2)  {$p_{2}$};
  \node[cVtx, above=3pt]       at (p3)  {$p_{3}$};
  \node[cVtx, below=4pt]       at (q12) {$q_{12}$};
  \node[cVtx, left=4pt]        at (q13) {$q_{13}$};
  \node[cVtx, right=4pt]       at (q23) {$q_{23}$};

  \node[font=\small\bfseries] at (0, {2*sqrt(3)+0.3}) {$G_{1}$};
\end{tikzpicture}
\end{center}

\noindent
Here $q_{ij}$ denotes the midpoint of the segment $p_{i}p_{j}$.  Since the graph $G_{1}$ has a reflection symmetry swapping $p_{1}$ and $p_{2}$ (reflection across the perpendicular bisector of $p_{1}p_{2}$), and the boundary values also swap ($V(p_{1})=1\leftrightarrow V(p_{2})=0$), the unique Kirchhoff solution satisfies $V\mapsto 1-V$ under this reflection.  Both $q_{12}$ and $p_{3}$ lie on the symmetry axis, so $V(q_{12})=1-V(q_{12})$ and $V(p_{3})=1-V(p_{3})$, giving $V(q_{12})=V(p_{3})=\frac{1}{2}$.  Setting $a=V(q_{13})$ (so $V(q_{23})=1-a$ by the same symmetry), Kirchhoff's equation~\eqref{eq:kirchhoff-law} at $q_{13}$ (which has four neighbors $p_{1}$, $q_{12}$, $q_{23}$, $p_{3}$ in $G_{1}$) reads
\begin{equation}\label{eq:gasket-kirchhoff}
  (1-a)+\Bigl(\frac{1}{2}-a\Bigr)+(1-2a)+\Bigl(\frac{1}{2}-a\Bigr)=0,
  \qquad\text{giving}\qquad a=\frac{3}{5}.
\end{equation}
By Remark~\ref{rem:dirichlet-electric}, $R_{1}=(V(p_{1})-V(p_{2}))/I_{p_{1}}$, where $I_{p_{1}}=\sum_{q\sim_{1}p_{1}}(V(p_{1})-V(q))$ is the current injected at $p_{1}$.  Since $p_{1}$ has exactly two neighbors $q_{12}$ and $q_{13}$ in $G_{1}$,
\begin{equation}\label{eq:gasket-current}
  I_{p_{1}}
  = \bigl(V(p_{1})-V(q_{12})\bigr)+\bigl(V(p_{1})-V(q_{13})\bigr)
  = \Bigl(1-\frac{1}{2}\Bigr)+\Bigl(1-\frac{3}{5}\Bigr)
  = \frac{9}{10}.
\end{equation}
Therefore we obtain $R_{1}=(V(p_{1})-V(p_{2}))/I_{p_{1}}=\frac{10}{9}$, and
\begin{equation}\label{eq:gasket-renorm-ratio}
  \frac{R_{1}}{R_{0}}=\frac{10/9}{2/3}=\frac{5}{3}.
\end{equation}
Iterating this procedure, we obtain that each step of self-similar refinement multiplies the effective resistance by $5/3$.  Hence without the factor $(5/3)^{m}$ in~\eqref{eq:E-gasket-m}, the raw sum $\sum_{p\sim_{m}q}(f(p)-f(q))^{2}$ would decay as $(3/5)^{m}\to 0$ for any fixed continuous $f$; the factor $(5/3)^{m}$ is the unique compensation making the limit energy finite and non-degenerate.
\end{remark}

\begin{remark}[Intuition for $d_{w}=\ln 5/\ln 2$]\label{rem:gasket-dw}
As established in Remark~\ref{rem:gasket-renorm}, the effective resistance between two points at Euclidean distance $r\asymp 2^{-m}$ scales as $(3/5)^{m}$.  Setting $r\asymp 2^{-m}$ identifies the resistance exponent $\zeta$:
\begin{equation}\label{eq:resist-scaling}
  R(x,y)\asymp\left(\frac{3}{5}\right)^{m}\asymp d(x,y)^{\zeta},\qquad
  \zeta=\frac{\ln(5/3)}{\ln 2}.
\end{equation}
The walk dimension is then given by the \emph{Einstein relation} $d_{w}=d_{h}+\zeta$, which links the exit-time exponent to the volume-growth and resistance exponents.  It is stated as equation~(1.2) in~\cite{lee:23:relations}, and it is how $\zeta$ is introduced in~\cite[Sec.~6.1]{barlow:17:random}, whose Theorem~2.72 is the estimate~\eqref{eq:resist-scaling} for the graphical gasket.  Kigami~\cite[Thm.~4.2.1]{kigami:01:analysis} proves the counterpart in the effective resistance metric, where $\zeta=1$ by construction and the relation reads $d_{s}=2d_{H}/(d_{H}+1)$, with $d_{H}$ the Hausdorff dimension in that metric rather than in $d$.  In the Euclidean metric we obtain
\begin{equation}\label{eq:einstein-gasket}
  d_{w}=d_{h}+\zeta=\frac{\ln 3}{\ln 2}+\frac{\ln(5/3)}{\ln 2}=\frac{\ln 5}{\ln 2}.
\end{equation}
The fact that $\zeta>0$ (resistance grows with distance faster than linearly) is responsible for $d_{w}>2$: the Brownian motion on $K$ spreads anomalously slowly compared to flat space.
\end{remark}

Remarks~\ref{rem:gasket-renorm} and~\ref{rem:gasket-dw} have established two facts. First, the factor $(5/3)^{m}$ in~\eqref{eq:E-gasket-m} is the unique compensation which keeps the sequence $m\mapsto\mathcal{E}_{m}(f|_{V_{m}},f|_{V_{m}})$ non-decreasing. Second, this same factor encodes the anomalous geometry of $K$ through the exponents $d_{h}=\ln 3/\ln 2$ and $d_{w}=\ln 5/\ln 2$.  Since the sequence is non-decreasing, the limit $\lim_{m\to\infty}\mathcal{E}_{m}(f|_{V_{m}},f|_{V_{m}})$ exists in $[0,+\infty]$ for every $f:V_{*}\to\R$. It is finite exactly for those $f$ which carry a bounded energy across all approximating graphs.  This class of functions, equipped with the limiting energy, is the Dirichlet form on $K$.

\begin{definition}[Dirichlet form on $K$]\label{def:SG-form}
Recall that $K$ is the Sierpi\'nski gasket of Definition~\ref{def:SG}.  Set
\begin{equation}\label{eq:SG-F-star}
  \mathcal{F}_{*}=\Bigl\{f:V_{*}\to\R \;:\; \sup_{m\ge 0}\,\mathcal{E}_{m}(f|_{V_{m}},f|_{V_{m}})<\infty\Bigr\},
  \qquad
  \mathcal{E}(f,f)=\lim_{m\to\infty}\mathcal{E}_{m}(f|_{V_{m}},f|_{V_{m}}).
\end{equation}
Every $f\in\mathcal{F}_{*}$ extends uniquely to a continuous function on $K$; let $\mathcal{F}$ denote the space of such extensions.  Then $(\mathcal{E},\mathcal{F})$ is a strongly local regular Dirichlet form on $L^{2}(K,\mu)$~\cite{kigami:01:analysis,barlow:98:diffusions}.
\end{definition}

The pair $(\mathcal{E},\mathcal{F})$ is the ambient Dirichlet form in Assumption~\ref{hyp:mmspace} with $X=K$.  Boundary conditions at $\partial K=V_{0}$ are imposed by restricting the domain of the form.

\begin{definition}[Neumann and Dirichlet Laplacians on $K$]\label{def:SG-laplacians}
\begin{enumerate}
\item {\it Neumann Laplacian.}  The generator of $(\mathcal{E},\mathcal{F})$ on $L^{2}(K,\mu)$ is the \emph{Neumann Laplacian} $\Delta_{N}$, with domain
\begin{equation}\label{eq:dom-N-gasket}
  \mathrm{Dom}(\Delta_{N})=\Bigl\{u\in\mathcal{F} \;:\; \exists\, f\in L^{2}(K,\mu),\quad \mathcal{E}(u,v)=\int_{K}fv\,d\mu\quad \forall\, v\in\mathcal{F}\Bigr\},
\end{equation}
and $\Delta_{N}u=-f$.

\item {\it Dirichlet Laplacian.}  Set $\mathcal{F}_{0}=\{f\in\mathcal{F}:f|_{V_{0}}=0\}$.  Then $(\mathcal{E},\mathcal{F}_{0})$ is itself a strongly local regular Dirichlet form on $L^{2}(K,\mu)$~\cite[Corollary~3.4.7]{kigami:01:analysis}, and its generator is the \emph{Dirichlet Laplacian} $\Delta_{D}$, with domain
\begin{equation}\label{eq:dom-D-gasket}
  \mathrm{Dom}(\Delta_{D})=\Bigl\{u\in\mathcal{F}_{0} \;:\; \exists\, f\in L^{2}(K,\mu),\quad \mathcal{E}(u,v)=\int_{K}fv\,d\mu\quad \forall\, v\in\mathcal{F}_{0}\Bigr\},
\end{equation}
and $\Delta_{D}u=-f$.
\end{enumerate}
\end{definition}

Both Laplacians admit a concrete pointwise description, introduced by Kigami~\cite{kigami:93:harmonic,kigami:01:analysis}, which parallels the classical notion of a Laplacian in smooth settings.  For $p\in V_{m}$, let $V_{m,p}$ denote the set of neighbors of $p$ in $G_{m}$; one has $\#V_{m,p}=4$ for $p\in V_{m}\setminus V_{0}$ and $\#V_{m,p}=2$ for $p\in V_{0}$.

\begin{definition}[Kigami Laplacian]\label{def:kigami-lap}
For $m\ge 1$ and $p\in V_{m}\setminus V_{0}$, the discrete Laplacian on $G_{m}$ is
\begin{equation}\label{eq:discrete-lap-gasket}
  \Delta_{m}f(p)=\frac{3}{2}\cdot 5^{m}\sum_{q\in V_{m,p}}(f(q)-f(p)).
\end{equation}
The domain of the \emph{Kigami Laplacian} is
\begin{equation}\label{eq:dom-kigami}
  \mathcal{D}=\Bigl\{f\in C(K) \;:\; \exists\, g\in C(K),\quad \lim_{m\to\infty}\max_{p\in V_{m}\setminus V_{0}}|\Delta_{m}f(p)-g(p)|=0\Bigr\},
\end{equation}
and one sets $\Delta f:=g$ for $f\in\mathcal{D}$.  For $p\in V_{0}$, the \emph{Neumann derivative} of $f$ at~$p$ is
\begin{equation}\label{eq:SG-neumann-deriv}
  (df)_{p}=-\lim_{m\to\infty}\left(\frac{5}{3}\right)^{m}\sum_{q\in V_{m,p}}(f(q)-f(p)),
\end{equation}
whenever the limit exists.
\end{definition}

\begin{remark}[Origin of the factor $\frac{3}{2}\cdot 5^{m}$ in~\eqref{eq:discrete-lap-gasket}]\label{rem:kigami-renorm}
The factor $\frac{3}{2}\cdot 5^{m}$ in the Kigami Laplacian~\eqref{eq:discrete-lap-gasket} arises from two independent contributions that we identify in turn.

\smallskip

\noindent
{\it Factor $5^{m}$: balancing energy against measure.}  The Kigami Laplacian is defined so that the discrete Green's identity
\begin{equation}\label{eq:discrete-green}
  \mathcal{E}_{m}(f,g) = -\sum_{p\in V_{m}\setminus V_{0}} g(p)\,\Delta_{m}f(p)\,\mu_{m}(p) + \text{boundary terms}
\end{equation}
holds, in analogy with $\mathcal{E}(f,g)=-\int_{K}g\,\Delta f\,d\mu$.  Here $\mu_{m}(p)$ is the natural discrete approximation to $\mu$ at level $m$.  Since $\mu$ assigns mass $3^{-m}$,
comparing the right hand side of~\eqref{eq:discrete-green} with~\eqref{eq:E-gasket-m} yields the factor $5^{m}$ in~\eqref{eq:discrete-lap-gasket}. 

\smallskip

\noindent
{\it Factor $3/2$: cell-sharing at interior vertices.}  The mass $\mu_{m}(p)$ is computed by distributing the $\mu$-mass of each level-$m$ cell equally among its three corner vertices.  One can check by induction that every interior vertex $p\in V_{m}\setminus V_{0}$ is a corner of exactly two level-$m$ cells (as is already visible in $G_{1}$ from Remark~\ref{rem:gasket-renorm}), so
\begin{equation}\label{eq:vertex-mass}
  \mu_{m}(p) = 2\cdot\frac{3^{-m}}{3} = \frac{2}{3^{m+1}}.
\end{equation}
The full prefactor is then $(5/3)^{m}/\mu_{m}(p)=5^{m}\cdot 3^{m+1}/(2\cdot 3^{m})=\frac{3}{2}\cdot 5^{m}$.
\end{remark}

\begin{remark}[Another explanation for the factor $5^{m}$ in~\eqref{eq:discrete-lap-gasket}]\label{rem:kigami-renorm-2}
In $\R^{n}$ on a grid of mesh size $h=2^{-m}$ (namely $2^{-m}\mathbb{Z}$), the discrete Laplacian carries a factor $h^{-2}=4^{m}$ because a random walk exits a ball of radius $h$ in time $\asymp h^{2}$.  On the Sierpi\'nski gasket, the walk dimension $d_{w}$ is $\ln 5/\ln 2$ (Remark~\ref{rem:gasket-dw}). This means that the exit time from a ball of radius $2^{-m}$ scales as $(2^{-m})^{d_{w}}=5^{-m}$, so the generator is renormalized by $5^{m}$ rather than $4^{m}$.
\end{remark}

The Kigami Laplacian $\Delta$ is related to $\Delta_{N}$ and $\Delta_{D}$ through boundary conditions at $V_{0}$, in full analogy with the classical smooth setting.

\begin{proposition}[\protect{\cite[Theorem~3.7.9]{kigami:01:analysis}}]\label{prop:kigami-bc}
Let $\mathcal{D}_{N}=\{f\in\mathcal{D}:(df)_{p}=0\ \forall\, p\in V_{0}\}$ and $\mathcal{D}_{D}=\{f\in\mathcal{D}:f|_{V_{0}}=0\}$.
\begin{enumerate}
\item One has $\mathcal{D}_{N}=\mathrm{Dom}(\Delta_{N})\cap\mathcal{D}$ and $\Delta_{N}|_{\mathcal{D}_{N}}=\Delta|_{\mathcal{D}_{N}}$; moreover, $\Delta_{N}$ is the Friedrichs extension of $\Delta|_{\mathcal{D}_{N}}$.
\item One has $\mathcal{D}_{D}=\mathrm{Dom}(\Delta_{D})\cap\mathcal{D}$ and $\Delta_{D}|_{\mathcal{D}_{D}}=\Delta|_{\mathcal{D}_{D}}$; moreover, $\Delta_{D}$ is the Friedrichs extension of $\Delta|_{\mathcal{D}_{D}}$.
\end{enumerate}
\end{proposition}

With the Neumann and Dirichlet Laplacians $\Delta_{N}$ and $\Delta_{D}$ of Definition~\ref{def:SG-laplacians} in place, it remains to verify that the Sierpi\'nski gasket fits into the sub-Gaussian framework of Assumption~\ref{hyp:mmspace}.  The following proposition confirms this and closes the verification of the abstract setup for $X=K$.

\begin{theorem}[Sub-Gaussian heat kernel on $K$]\label{prop:SG-hk}
Here $K$ is the Sierpi\'nski gasket of Definition~\ref{def:SG}, $\mu$ is the self-similar measure of Definition~\ref{def:SG-measure}, $\Delta_{N}$ and $\Delta_{D}$ are the Laplacians of Definition~\ref{def:SG-laplacians}, and $d_{h}=\ln 3/\ln 2$, $d_{w}=\ln 5/\ln 2$ are the exponents of~\eqref{eq:gasket-dims}.  The semigroups $e^{t\Delta_{N}}$ and $e^{t\Delta_{D}}$ each admit a continuous heat kernel $p_{t}(x,y)$ satisfying the sub-Gaussian estimates~\eqref{eq:subGauss}.  In particular, $(K,d,\mu,\mathcal{E})$ satisfies Assumption~\ref{hyp:mmspace} with $X=K$~\cite{barlow:98:diffusions}.
\end{theorem}

We close this section by mentioning that open sets of the form $U=K\setminus V_{0}$ (the gasket with the boundary vertices removed) or $U=\mathfrak{f}_{i_{1}}\circ\cdots\circ\mathfrak{f}_{i_{n}}(K\setminus V_{0})$ (scaled copies thereof) are inner uniform domains~\cite[Sec.~2.4.3]{baudoin:24:parabolic}, so all the results of Sections~\ref{sssec:dirichlet-semi}--\ref{sssec:neumann-frac} apply.  For a systematic account of inner uniform and uniform domains within the Sierpi\'nski gasket, we refer to~\cite{lierl:22:dirichlet,barlow:98:diffusions}.

\section*{Bibliographic notes}
\addcontentsline{toc}{section}{Bibliographic notes}
\markright{BIBLIOGRAPHIC NOTES}

\smallskip

\noindent
{\it Recurrent metric measure spaces.} The Dirichlet form material of
Subsection~\ref{ssec:mmspace} is classical, and the reference we follow is
Fukushima, \=Oshima and Takeda~\cite{fukushima:94:dirichlet}. We borrow from it
the Beurling-Deny decomposition (as stated in Proposition~\ref{prop:BD}).
The construction of the Markov process attached to a regular
Dirichlet form, which is what turns Assumption~\ref{hyp:mmspace} into a
statement about a Brownian motion on $X$, is also taken from there. Chen and
Fukushima~\cite{chen:12:symmetric} develop the same theory from the side of the
process, where our Dirichlet form~\eqref{eq:E-D} becomes the form of the
\emph{part process} of $B$ killed at the exit time $T_{U}$
(see~\cite[Sec.~3.3]{chen:12:symmetric}). As far as our Neumann
form~\eqref{eq:E-N} is concerned, a \emph{reflected Dirichlet space} is introduced
(see~\cite[Chap.~6]{chen:12:symmetric}). 

The study of diffusions with a walk dimension $d_{w}>2$ begins with Barlow and
Perkins~\cite{barlow:88:brownian}. They construct Brownian motion on the
Sierpi\'nski gasket as a scaling limit of random walks on the approximating
graphs, and they prove for it the two-sided sub-Gaussian
bounds~\eqref{eq:subGauss} with $d_{w}=\ln5/\ln2$. That paper is where the
exponent $\frac{1}{d_{w}-1}$ of~\eqref{eq:subGauss} first appears. It also
becomes clear that the Gaussian shape of the heat kernel is related to
$d_{w}=2$ rather than a general fact. The spectral
dimension $d_{s}$ of~\eqref{eq:spectral-dimension}, which is an important
parameter for them, plays no role in the geometric chapters, where everything
is written in $d_{h}$ and $d_{w}$. For instance their on-diagonal decay
$t^{-d_{s}/2}$ is our $t^{-d_{h}/d_{w}}$.
Lindstr{\o}m~\cite{lindstrom:90:brownian} carries the construction over to the
larger class of nested fractals. Barlow's Saint-Flour
lectures~\cite{barlow:98:diffusions} remain the best introduction to this
circle of ideas, and they are the source we quote for the recurrence
Lemma~\ref{lem:recurrence} and for the spectral
expansion~\eqref{eq:HK-spectral-D}. Their counterpart for graphs, where the
same exponents govern the random walk, is Barlow's
book~\cite{barlow:17:random}.

On the analytic side the Dirichlet form itself can be built without any
probability, as a limit of discrete energies along the approximating graphs.
This is Kigami's construction, carried out for the gasket and then for
post-critically finite self-similar sets in~\cite{kigami:93:harmonic}.
The book~\cite{kigami:01:analysis} is a systematic account on the topic. Note that we use both
the probabilistic and the analytic routes in Section~\ref{sec:geom-examples}. Indeed, the harmonic structure of
Subsection~\ref{ssec:sierpinski} is Kigami's while the heat kernel bounds are
Barlow and Perkins'.

We have taken Assumption~\ref{hyp:mmspace} as a hypothesis, and said nothing
about how one checks it (although we have checked it on examples). Deciding which spaces satisfy a sub-Gaussian estimate
is a subject in itself, and the sharpest results are due to Grigor'yan and his
coauthors. On a metric measure space carrying a strongly local regular
Dirichlet form, Grigor'yan and Telcs~\cite{grigoryan:12:two-sided} find
conditions equivalent to a two-sided bound of the general shape
\begin{equation*}
  p_{t}(x,y) \asymp \frac{C}{V(x,R(t))}
  \exp\left(-c\,t\,\Psi\Bigl(c\,\frac{d(x,y)}{t}\Bigr)\right),
\end{equation*}
of which~\eqref{eq:subGauss} is the case $V(x,r)\asymp r^{d_{h}}$ and
$R(t)\asymp t^{1/d_{w}}$. The graph
version is cleaner still. In~\cite{grigoryan:01:sub-gaussian} the same authors
show that a two-sided sub-Gaussian estimate on an infinite weighted graph holds
if and only if the volume growth is uniformly polynomial and the Green kernel
decays uniformly polynomially. The Gaussian case $d_{w}=2$ on graphs, which is
the discrete counterpart of Li and Yau's theorem, is treated by Coulhon and
Grigor'yan~\cite{coulhon:98:random}, and the manifold case by Grigor'yan's
book~\cite{grigoryan:09:heat}. A reader who wants the probabilistic
counterpart of all this, on disordered media and with the scaling limits
included, will find it in Kumagai's Saint-Flour
lectures~\cite{kumagai:14:random}.

\smallskip

\noindent
{\it Inner uniform domains.} The notion of a uniform domain
(Definition~\ref{def:uniform}) comes from the quasiconformal literature, where
it was introduced by Martio and Sarvas~\cite{martio:79:injectivity}. The inner
uniform variant of Definition~\ref{def:inner-uniform}, in which the ambient
metric is replaced by the inner metric $d_{U}$, is what makes a slit domain
admissible. This is  illustrated by Example~\ref{ex:slit-annulus}. The heat kernel
theory on such domains is due to Gyrya and
Saloff-Coste~\cite{gyrya:11:neumann}, who treat the Dirichlet and the Neumann
problem together. Their setting is the Gaussian one, governed by a parabolic Harnack
inequality, so $d_{w}=2$ throughout. Their domains are moreover unbounded (the
harmonic profile is defined for an unbounded $U$), where ours is bounded with
compact closure. Bounded domains, non-symmetric forms and local
versions of the theory are covered by Lierl and
Saloff-Coste~\cite{lierl:14:dirichlet}. Local is meant there in two senses. The
estimates hold on a fixed time interval $(0,T)$, with constants depending on
$T$, and the hypotheses (the parabolic Harnack inequality, and the inner
uniformity of $U$) are required only near a subdomain rather than on all of
$X$. 

What we actually need is the extension of that theory past $d_{w}=2$, and this
is Lierl's~\cite{lierl:22:dirichlet}. Theorem~\ref{thm:hk-diri} is her result.
It is worth noticing where the features of that statement come from. The
factor $\Phi_{1}(x)\Phi_{1}(y)$ of the lower
bound~\eqref{eq:subGauss-D-lower}, which is the analytic trace of the killing
at $\partial U$, and the inner metric $d_{U}$ in place of $d$, are not fractal
effects. Both are already present in the Gaussian bounded-domain estimate
of~\cite{lierl:14:dirichlet}, where $\Phi_{1}$ takes over the role played on an
unbounded domain by the harmonic profile of~\cite{gyrya:11:neumann}. What the
fractal setting adds to~\eqref{eq:subGauss-D-lower} is the exponent
$\frac{1}{d_{w}-1}$ and the normalization by $\min(1,t^{d_{h}/d_{w}})$. The Neumann counterpart, Theorem~\ref{thm:hk-neum},
is due to Murugan~\cite{murugan:24:heat} and requires the strictly stronger
uniform condition of Definition~\ref{def:uniform}. That hypothesis comes from
the method rather than from the phenomenon. Indeed, Murugan's route to the heat kernel bounds
relies on extension operators.  He conjectures in his concluding remarks
that the result should hold on inner uniform domains as well, at the price of an intrinsic argument bypassing
 the extension operator machinery. In the Gaussian case that version is already
known (see the Neumann half of~\cite{gyrya:11:neumann}).

\smallskip

\noindent
{\it Examples.} Each of the examples of Section~\ref{sec:geom-examples} has a classical
literature of its own, of which we quote only the entry points.

For metric graphs, the function spaces and the Kirchhoff vertex conditions of
Definition~\ref{def:metric-graph} are treated in Post's
book~\cite{post:12:spectral}. The spectral theory of the resulting
operators, under the name of quantum graphs, is handled in Berkolaiko and
Kuchment~\cite{berkolaiko:13:introduction}. A reader checking this reference
should be aware that they call our Kirchhoff condition  Neumann (or
Neumann-Kirchhoff) condition. The Gaussian bounds of
Theorem~\ref{thm:kirchhoff-props} are not quoted from a metric graph
reference. They are the case $d_{w}=2$ of the general criterion
of~\cite{grigoryan:12:two-sided} already discussed above, applied to a space
of volume growth $V(x,r)\asymp r$.

For Riemannian manifolds the two-sided Gaussian estimate under non-negative
Ricci curvature is Li and Yau's~\cite{li:86:parabolic}, and
Grigor'yan's book~\cite{grigoryan:09:heat} is the systematic treatment.
Davies~\cite{davies:89:heat} gives an account on the semigroup and spectral side. Recall
however that this example is largely formal for us, since the recurrence
condition $d_{h}<d_{w}$ forces $n=1$ (see item~(i) of
Remark~\ref{rem:riem-remarks}). The parabolic Anderson model on a compact
manifold of dimension $n\ge2$ is studied by other means
in~\cite{chen:25:global,chen:26:sharp}.

For the Sierpi\'nski gasket the two sources already quoted,
\cite{kigami:01:analysis} for the harmonic structure and
\cite{barlow:98:diffusions} for the diffusion, cover everything we use.
Strichartz's tutorial~\cite{strichartz:06:differential} is a gentler entry
point, written for a reader who wants to compute on the gasket rather than to
build the theory. The spectrum of the gasket Laplacian is described completely by Fukushima and
Shima~\cite{fukushima:92:spectral}. Their tool is the decimation method of
Rammal and Toulouse, which computes the eigenvalues at one level of the
approximation from those at the previous one. In addition, 
the eigenvalue counting function
$\rho(\lambda)$ is of order $\lambda^{d_{s}/2}$ (recall that $d_{s}$ is the spectral
dimension~\eqref{eq:spectral-dimension}), but the ratio
$\rho(\lambda)/\lambda^{d_{s}/2}$ oscillates and does not converge. This is in contrast with the Euclidean case. Finally, the Einstein relation
behind~\eqref{eq:gasket-dims}, which ties the walk dimension to the volume and
resistance exponents, holds well beyond the self-similar setting.
Lee~\cite[eq.~(1.2)]{lee:23:relations} proves it for unimodular random
networks, a setting far from ours in that it is random and discrete where the
gasket is deterministic and continuous. His statement is conditional, the
volume and resistance exponents being assumed to exist, and the equality
further requires the resistance exponent to be non-negative. Without that sign
condition what remains is the inequality $d_{w}\ge d_{h}+\zeta$ (recall that
$\zeta$ is the resistance exponent of~\eqref{eq:resist-scaling}). That paper is
also a good entry point for these scaling exponents and the relations among
them.

\chapter{PAM on metric spaces}
\label{chap:geom-pam}

With the geometric framework of Chapter~\ref{chap:geom-setting} in hand, we can now analyze the parabolic Anderson model on a metric measure space satisfying Assumption~\ref{hyp:mmspace}. Our route parallels the Euclidean one. Namely we first construct a natural family of noises on a bounded domain $U\subset X$ (Section~\ref{sec:geom-noise}). We then solve the equation by means of chaos expansions (Section~\ref{sec:geom-existence}). Finally we quantify the intermittent behavior of the solution through upper and lower bounds on its moments (Section~\ref{sec:geom-moments}).

\section[Noises for PAM]{Noises for PAM on metric measure spaces}\label{sec:geom-noise}

In Chapter~\ref{chap:noise}, the driving noise for the PAM on $\R^{d}$ was a centered
Gaussian family white in time and correlated in space with covariance
$\Lambda(x-y)$, and its Hilbert space $\mathcal{H}$ was built by completing the test
functions under the corresponding inner product (Section~\ref{ssec:noise}).  In the
geometric setting of Chapter~\ref{chap:geom-setting}, the translation-invariant kernel $\Lambda(x-y)$
is no longer available. The role it played is taken over by the Riesz kernels
$G^{D,U}_{2\alpha}$ or $G^{N,U}_{2\alpha}$ introduced in
Sections~\ref{sssec:dirichlet-frac} and~\ref{sssec:neumann-frac}.  The parameter
$\alpha\ge 0$ controls the spatial regularity of the noise. Namely $\alpha=0$
corresponds to the space-time white noise on $U$, while $\alpha>0$ introduces a
spatially colored covariance, produced by the fractional operator
$(-\Delta_{U})^{-2\alpha}$.

\subsection{The fractional noise family}\label{ssec:geom-noise-def}

We work throughout this section under the standing framework of the previous sections:
$(X,d,\mu)$ is a metric measure space satisfying Assumption~\ref{hyp:mmspace} and $U\subset X$
is a bounded open set with compact closure.  We define the noise for the Dirichlet
setting first, and then for the Neumann setting. The two coincide when $\alpha=0$.
More specifically, the covariance of the Dirichlet noise is governed by the Riesz kernel $G^{D,U}_{2\alpha}$
of~\eqref{eq:RK-D}, and the natural Hilbert space of test functions is the Dirichlet
Sobolev space $\mathcal{W}^{-\alpha}_{D}(U)$ of Definition~\ref{def:sobolev-D}, extended
to space-time functions.

\begin{definition}[Dirichlet fractional noise]\label{def:geom-noise-D}
Let $U\subset X$ be an inner uniform domain and $\alpha>0$.  The space-time Hilbert space corresponding to the Dirichlet fractional noise is
\begin{equation}\label{eq:H-geom-D}
  \mathcal{H}^{D}_{\alpha}
  \;:=\;
  L^{2}\!\bigl(\R_{+},\;\mathcal{W}^{-\alpha}_{D}(U)\bigr),
\end{equation}
with inner product
$\langle\varphi,\psi\rangle_{\mathcal{H}^{D}_{\alpha}}
:=\int_{\R_{+}}\langle\varphi(t,\cdot),\psi(t,\cdot)\rangle_{\mathcal{W}^{-\alpha}_{D}(U)}\,dt$.
The \emph{Dirichlet fractional noise} is a centered Gaussian family
$\{W^{D}_{\alpha}(\varphi);\,\varphi\in\mathcal{H}^{D}_{\alpha}\}$ on
$(\Omega,\mathcal{F},\Prob)$ whose covariance is
\begin{equation}\label{eq:geom-cov-D}
  \E\bigl[W^{D}_{\alpha}(\varphi)\,W^{D}_{\alpha}(\psi)\bigr]
  \;=\;
  \langle\varphi,\psi\rangle_{\mathcal{H}^{D}_{\alpha}},
  \qquad\varphi,\psi\in\mathcal{H}^{D}_{\alpha}.
\end{equation}
Unfolding the inner product of Definition~\ref{def:sobolev-D}, the covariance structure can be spelled out as
\begin{equation}\label{eq:geom-cov-D-kernel}
  \E\bigl[W^{D}_{\alpha}(\varphi)\,W^{D}_{\alpha}(\psi)\bigr]
  =
  \int_{\R_{+}}\!\int_{U\times U}
    \varphi(t,x)\,\psi(t,y)\,G^{D,U}_{2\alpha}(x,y)\,d\mu(x)\,d\mu(y)\,dt.
\end{equation}
\end{definition}

Since $G^{D,U}_{2\alpha}(x,y)\ge 0$ for all $x,y\in U$ (see~\eqref{eq:RK-D}), the right-hand side
of~\eqref{eq:geom-cov-D-kernel} defines a genuine inner product on
$\mathcal{W}^{-\alpha}_{D}(U)$.  The covariance identity~\eqref{eq:geom-cov-D} means that
the map $\varphi\mapsto W^{D}_{\alpha}(\varphi)$ is a linear isometry from
$\mathcal{H}^{D}_{\alpha}$ to $L^{2}(\Omega)$, of the form
\begin{equation}\label{eq:geom-W-isometry-D}
  \E\bigl[W^{D}_{\alpha}(h)\,W^{D}_{\alpha}(k)\bigr]
  \;=\;
  \langle h,k\rangle_{\mathcal{H}^{D}_{\alpha}}
  \qquad\text{for all }h,k\in\mathcal{H}^{D}_{\alpha}.
\end{equation}
In particular, $W^{D}_{\alpha}$ is an isonormal Gaussian process on
$\mathcal{H}^{D}_{\alpha}$ in the sense of Definition~\ref{def:isonormal}.  Following
the convention of Section~\ref{ssec:noise}, we write Wiener integrals as
\begin{equation*}
  W^{D}_{\alpha}(h)
  \;=\;
  \int_{0}^{\infty}\!\int_{U}h(s,x)\,W^{D}_{\alpha}(ds,dx),
  \qquad h\in\mathcal{H}^{D}_{\alpha}.
\end{equation*}

The case $\alpha=0$ plays a special role. Indeed, the covariance kernel
$G^{D,U}_{0}$ reduces to the Dirac mass on the diagonal, and the noise becomes
the space-time white noise on $U$.

\begin{definition}[Space-time white noise on $U$]\label{def:geom-noise-white}
The Hilbert space for the white noise is
\begin{equation}\label{eq:H-geom-0}
  \mathcal{H}_{0}
  \;:=\;
  L^{2}(\R_{+},\,L^{2}(U,\mu)),
\end{equation}
with inner product
$\langle\varphi,\psi\rangle_{\mathcal{H}_{0}}
=\int_{\R_{+}}\int_{U}\varphi(t,x)\,\psi(t,x)\,d\mu(x)\,dt$.
The \emph{space-time white noise} on $U$ is a centered Gaussian family
$\{W_{0}(\varphi);\,\varphi\in\mathcal{H}_{0}\}$ with covariance
\begin{equation}\label{eq:geom-cov-0}
  \E\bigl[W_{0}(\varphi)\,W_{0}(\psi)\bigr]
  \;=\;
  \int_{\R_{+}}\!\int_{U}\varphi(t,x)\,\psi(t,x)\,d\mu(x)\,dt.
\end{equation}
As for $\alpha>0$, the map $W_{0}:\mathcal{H}_{0}\to L^{2}(\Omega)$ is a linear
isometry, and $W_{0}$ is an isonormal Gaussian process on $\mathcal{H}_{0}$.
\end{definition}

The three regimes identified in Proposition~\ref{prop:RK-D-bounds} for the kernel
$G^{D,U}_{2\alpha}$ have a direct probabilistic interpretation in terms of the regularity
of the field $W^{D}_{\alpha}$.

\begin{remark}[Regularity regimes]\label{rem:geom-noise-reg}
The critical exponent $d_{h}/(2d_{w})$ separates three re\-gi\-mes.
The precise Sobolev and H\"{o}lder regularity of $W^{D}_{\alpha}$ in each regime
is worked out in Section~\ref{ssec:sobolev-inclusion}, but we can already summarize what we get as follows:
\begin{enumerate}
\item {\it Distributional regime ($0\le\alpha<d_{h}/(2d_{w})$).}
The covariance kernel $G^{D,U}_{2\alpha}(x,y)$ diverges like
$d_{U}(x,y)^{2\alpha d_{w}-d_{h}}\to+\infty$ near the diagonal
(see~\eqref{eq:RK-D-small}), so the field $W^{D}_{\alpha}$ is
a distribution in the space variable.  This includes the white noise case $\alpha=0$.

\item {\it Log-correlated regime ($\alpha=d_{h}/(2d_{w})$).}
The kernel grows logarithmically: $G^{D,U}_{2\alpha}(x,y)\asymp|\ln d_{U}(x,y)|$
near the diagonal (see~\eqref{eq:RK-D-crit}).
\item {\it Pointwise regime ($\alpha>d_{h}/(2d_{w})$).}
The kernel $G^{D,U}_{2\alpha}$ is bounded on $U\times U$
(see~\eqref{eq:RK-D-large}), and the field $W^{D}_{\alpha}$ can be
realized as a pointwise continuous Gaussian field on $\R_{+}\times U$.
\end{enumerate}
\end{remark}

We now turn to the Neumann setting.  The construction mirrors Definition~\ref{def:geom-noise-D},
with $\mathcal{W}^{-\alpha}_{N}(U)$ replacing $\mathcal{W}^{-\alpha}_{D}(U)$.  Recall
from Definition~\ref{def:sobolev-N} that the inner product on
$\mathcal{W}^{-\alpha}_{N}(U)$ involves the shifted kernel
$G^{N,U}_{2\alpha}(x,y)+C_{U}$. Indeed, since $G^{N,U}_{2\alpha}$ is not pointwise non-negative
(Proposition~\ref{prop:RK-N-bounds}), the constant $C_{U}$ is chosen large enough so that
\begin{equation*}
G^{N,U}_{2\alpha}+C_{U}\ge 0 \quad\text{pointwise,}
\end{equation*}
restoring positivity.

\begin{definition}[Neumann fractional noise]\label{def:geom-noise-N}
Let $U\subset X$ be a uniform domain and $\alpha>0$.  The space-time Hilbert space is
\begin{equation}\label{eq:H-geom-N}
  \mathcal{H}^{N}_{\alpha}
  \;:=\;
  L^{2}\!\bigl(\R_{+},\;\mathcal{W}^{-\alpha}_{N}(U)\bigr),
\end{equation}
with the inner product inherited from Definition~\ref{def:sobolev-N}.  The \emph{Neumann
fractional noise} is a centered Gaussian family
$\{W^{N}_{\alpha}(\varphi);\,\varphi\in\mathcal{H}^{N}_{\alpha}\}$ on
$(\Omega,\mathcal{F},\Prob)$ with covariance
\begin{equation}\label{eq:geom-cov-N}
  \E\bigl[W^{N}_{\alpha}(\varphi)\,W^{N}_{\alpha}(\psi)\bigr]
  \;=\;
  \langle\varphi,\psi\rangle_{\mathcal{H}^{N}_{\alpha}},
  \qquad\varphi,\psi\in\mathcal{H}^{N}_{\alpha}.
\end{equation}
Unfolding the inner product of Definition~\ref{def:sobolev-N}, the expression for the covariance function reads
\begin{equation}\label{eq:geom-cov-N-kernel}
  \E\bigl[W^{N}_{\alpha}(\varphi)\,W^{N}_{\alpha}(\psi)\bigr]
  =
  \int_{\R_{+}}\!\int_{U\times U}
    \varphi(t,x)\,\psi(t,y)\,\bigl(G^{N,U}_{2\alpha}(x,y)+C_{U}\bigr)\,
    d\mu(x)\,d\mu(y)\,dt.
\end{equation}
The map $W^{N}_{\alpha}:\mathcal{H}^{N}_{\alpha}\to L^{2}(\Omega)$ is a linear
isometry and $W^{N}_{\alpha}$ is an isonormal Gaussian process on
$\mathcal{H}^{N}_{\alpha}$.  We write $W^{N}_{\alpha}(h)
=\int_{0}^{\infty}\int_{U}h(s,x)\,W^{N}_{\alpha}(ds,dx)$ for
$h\in\mathcal{H}^{N}_{\alpha}$.
\end{definition}

\begin{remark}\label{rem:geom-noise-alpha0}
For $\alpha=0$ both boundary conditions yield the same noise.  The formal limit of
$\mathcal{W}^{-\alpha}_{D}(U)$ and $\mathcal{W}^{-\alpha}_{N}(U)$ as $\alpha\to 0$
is $L^{2}(U,\mu)$, so both $W^{D}_{\alpha}$ and $W^{N}_{\alpha}$ degenerate to the
space-time white noise $W_{0}$ of Definition~\ref{def:geom-noise-white}, with
$\mathcal{H}^{D}_{0}=\mathcal{H}^{N}_{0}=\mathcal{H}_{0}$.  The choice of boundary
condition is therefore relevant only for $\alpha>0$.
\end{remark}

\subsection{Wiener chaos and Malliavin calculus in the geometric setting}
\label{ssec:geom-noise-chaos}

Since each of $W^{D}_{\alpha}$ and $W^{N}_{\alpha}$ is an isonormal Gaussian process on
its respective Hilbert space $\mathcal{H}_{\alpha}$ (either $\mathcal{H}^{D}_{\alpha}$ or
$\mathcal{H}^{N}_{\alpha}$), the full Wiener chaos framework of
Section~\ref{ssec:isonormal} carries over verbatim with the abstract $\mathcal{H}$ of
that section replaced by the present $\mathcal{H}_{\alpha}$.  We briefly record the main
consequences in this subsection.

First, every $F\in L^{2}(\Omega,\sigma(W_{\alpha}),\Prob)$ admits a unique Wiener chaos decomposition
\begin{equation}\label{eq:geom-chaos-decomp}
  F \;=\; \E[F] + \sum_{n=1}^{\infty}I_{n}(f_{n}),
\end{equation}
where $f_{n}\in(\mathcal{H}_{\alpha})^{\odot n}$ for each $n\ge 1$ and the series converges
in $L^{2}(\Omega)$ (Theorem~\ref{thm:chaos-decomp}).  The $n$-th multiple stochastic
integral $I_{n}:\mathcal{H}_{\alpha}^{\odot n}\to\mathbf{H}_{n}$ is an isometry satisfying (see~\eqref{eq:In-isometry})
\begin{equation*}
\E[I_{n}(f)\,I_{m}(g)]=\mathbf{1}_{n=m}\cdot n!\,\langle f,g\rangle_{\mathcal{H}_{\alpha}^{\otimes n}}.
\end{equation*}

The Malliavin derivative $D$ and Skorohod integral $\delta$ of
Section~\ref{ssec:malliavin} are equally available in this setting.  For
$F=\sum_{n\ge 0}I_{n}(f_{n})\in\mathbb{D}^{1,2}$ the Malliavin derivative acts chaos-level
by chaos-level according to
\begin{equation}\label{eq:geom-D-chaos}
  DF \;=\; \sum_{n=1}^{\infty}n\,I_{n-1}\bigl(f_{n}(\,\cdot\,,\star)\bigr),
\end{equation}
where $\star$ denotes the free $\mathcal{H}_{\alpha}$-variable
(see~\eqref{eq:D-chaos}).  The duality relation~\eqref{eq:duality} between $D$ and
$\delta$, the product rule~\eqref{eq:product-rule}, and Nelson's hypercontractivity
estimate~\eqref{eq:hypercontractivity} all carry over without modification.  For adapted
integrands the Skorohod integral coincides with the It\^o--Walsh integral
(Proposition~\ref{prop:ito-sko}), which is the stochastic integral used in the mild
formulation of the PAM in the following section.

\subsection{Inclusion in Sobolev spaces}\label{ssec:sobolev-inclusion}

For this regularity study, we will focus on the Dirichlet setting for conciseness.  The Neumann case can be treated
along the same lines, replacing $-\Delta_{D,U}$ by $-\Delta_{N,U}$ throughout. We
refer to~\cite{baudoin:22:fractional} for the details.

The Dirichlet Sobolev space $\mathcal{W}^{-\alpha}_{D}(U)$ (Definition~\ref{def:sobolev-D}) and the noise $W^{D}_{\alpha}$ (Definition~\ref{def:geom-noise-D}) are defined
for any metric measure space satisfying Assumption~\ref{hyp:mmspace}.  On the
Sierpi\'nski gasket $K$, the precise spectral asymptotics of $-\Delta_{D,U}$ allow a
sharper analysis. Namely one can characterize, thanks to a classical scale of Sobolev
spaces on $K$, exactly how singular $W^{D}_{\alpha}$ is in the spatial variable. We can also identify a pointwise
H\"older-continuous realization of the noise when $\alpha>d_{h}/(2d_{w})$.  Throughout this
subsection we work in the setting of Subsection~\ref{ssec:sierpinski}: $K$ is the Sierpi\'nski
gasket of Definition~\ref{def:SG}, $V_{0}=\{p_{1},p_{2},p_{3}\}$ is the boundary,
$U=K\setminus V_{0}$, $\mu$ is the self-similar measure of Definition~\ref{def:SG-measure}, and
$(\lambda_{j},\Phi_{j})_{j\ge 1}$ are the eigenpairs of $-\Delta_{D,U}$ from
Notation~\ref{not:spectrum-D}.

\smallskip

\noindent
Let us start by discussing spectral asymptotics.
Lemma~\ref{lem:spectral-D} recorded the sufficient condition $\alpha>d_{h}/d_{w}$ for the
convergence of $\sum_{j}\lambda_{j}^{-\alpha}$.  For the Sierpi\'nski gasket, Fukushima and
Shima~\cite{fukushima:92:spectral} proved a sharp two-sided bound on the eigenvalue counting
function $N(t):=\#\{j:\lambda_{j}\le t\}$, turning that condition into an exact equivalence.

\begin{lemma}[Weyl law on $K$]\label{lem:weyl}
There exist constants $0<c\le C<\infty$ such that
\begin{equation}\label{eq:weyl}
  c\,t^{d_{h}/d_{w}} \;\le\; N(t) \;\le\; C\,t^{d_{h}/d_{w}},
  \qquad t>0.
\end{equation}
As a consequence, for every $\beta>0$,
\begin{equation}\label{eq:weyl-series}
  \sum_{j=1}^{\infty}\lambda_{j}^{-\beta}\;<\;+\infty
  \quad\text{if and only if}\quad \beta>\frac{d_{h}}{d_{w}}.
\end{equation}
\end{lemma}

\begin{proof}
The bound~\eqref{eq:weyl} is due to~\cite{fukushima:92:spectral}.  For~\eqref{eq:weyl-series},
a Stieltjes integration gives $\sum_{j}\lambda_{j}^{-\beta}=\int_{0}^{\infty}t^{-\beta}\,dN(t)$.
An integration by parts, combined with the two-sided bound~\eqref{eq:weyl}, shows that this
integral is finite if and only if
$\int_{1}^{\infty}t^{-\beta+d_{h}/d_{w}-1}\,dt<\infty$, i.e., $\beta>d_{h}/d_{w}$.
\end{proof}

As a direct application of Lemma~\ref{lem:weyl}, we have that $\sum_{j}\lambda_{j}^{-\beta-2\alpha}<\infty$ if and only if $\beta>d_{h}/d_{w}-2\alpha$.
This sharpens the one-directional statement~\eqref{eq:eigen-series-D} and is the key input for the
Sobolev regularity analysis below. Namely we introduce the Hilbert scale on $U$ naturally associated with the Dirichlet eigenvalues,
following~\cite[Sec.~3.1]{baudoin:23:dirichlet}.

\begin{definition}[Sobolev spaces on $K$]\label{def:sobolev-K}
For $f\in\mathcal{S}_{D}(U)$ set $\hat{f}_{j}:=\int_{U}f\,\Phi_{j}\,d\mu$, where $(\lambda_{j},\Phi_{j})_{j\ge 1}$ are the eigenpairs of $-\Delta_{D,U}$ from Notation~\ref{not:spectrum-D}.  For $\alpha\ge 0$,
the \emph{Sobolev space} $H^{\alpha}(K)$ is the completion of $\mathcal{S}_{D}(U)$ under the norm
\begin{equation}\label{eq:H-alpha-norm}
  \|f\|^{2}_{H^{\alpha}(K)}
  := \sum_{j=1}^{\infty}\lambda_{j}^{\alpha}\,\hat{f}_{j}^{2}.
\end{equation}
Its dual $H^{-\alpha}(K):=(H^{\alpha}(K))^{*}$, viewed as a subspace of $\mathcal{S}'_{D}(U)$,
carries the canonical norm
\begin{equation}\label{eq:H-neg-alpha-norm}
  \|\psi\|^{2}_{H^{-\alpha}(K)}
  := \sum_{j=1}^{\infty}\lambda_{j}^{-\alpha}\,\psi(\Phi_{j})^{2}.
\end{equation}
Here $\psi(\Phi_{j}) := \langle \psi, \Phi_{j} \rangle_{\mathcal{S}'_{D}(U),\,\mathcal{S}_{D}(U)}$
denotes the action of the distribution $\psi$ on the test function $\Phi_{j}$.
\end{definition}

That~\eqref{eq:H-neg-alpha-norm} is indeed the dual norm is the content of the next lemma.

\begin{lemma}[Dual norm on $H^{-\alpha}(K)$]\label{lem:sobolev-K-dual}
For every $\psi\in H^{-\alpha}(K)$, the right-hand side of~\eqref{eq:H-neg-alpha-norm} is such that
\begin{equation}\label{i1}
\|\psi\|^{2}_{H^{-\alpha}(K)}
  = \sum_{j=1}^{\infty}\lambda_{j}^{-\alpha}\,\psi(\Phi_{j})^{2}
  =
  \sup\{|\psi(f)|^{2}:\|f\|_{H^{\alpha}(K)}\le 1\} .
\end{equation}
\end{lemma}

\begin{proof}
By the Riesz representation theorem there exists $f_{\psi}\in H^{\alpha}(K)$ satisfying
$\psi(g)=(f_{\psi},g)_{H^{\alpha}(K)}$ for all $g\in\mathcal{S}_{D}(U)$, and
$\|\psi\|_{H^{-\alpha}}=\|f_{\psi}\|_{H^{\alpha}}$.  Setting
$\hat{f}^{\psi}_{j}:=\int_{U}f_{\psi}\Phi_{j}\,d\mu$ and evaluating at $g=\Phi_{j}$ gives
\begin{equation*}
  \psi(\Phi_{j})
  = (f_{\psi},\Phi_{j})_{H^{\alpha}(K)}
  = \sum_{k=1}^{\infty}\lambda_{k}^{\alpha}\,\hat{f}^{\psi}_{k}\int_{U}\Phi_{j}\Phi_{k}\,d\mu
  = \lambda_{j}^{\alpha}\hat{f}^{\psi}_{j},
\end{equation*}
where the last step uses the $L^{2}(U,\mu)$-orthonormality of the functions $\Phi_{j}$. Therefore we get
\begin{equation*}
  \|\psi\|^{2}_{H^{-\alpha}}
  = \|f_{\psi}\|^{2}_{H^{\alpha}}
  = \sum_{j=1}^{\infty}\lambda_{j}^{\alpha}\bigl(\hat{f}^{\psi}_{j}\bigr)^{2}
  = \sum_{j=1}^{\infty}\lambda_{j}^{-\alpha}\psi(\Phi_{j})^{2}.
\end{equation*}
This proves the first identity in \eqref{i1}.
For the second, fix $f\in H^{\alpha}(K)$ with $\|f\|_{H^{\alpha}}\le 1$.
By the Cauchy-Schwarz inequality in $H^{\alpha}(K)$,
\begin{equation*}
  |\psi(f)|^{2}
  = |(f_{\psi},f)_{H^{\alpha}}|^{2}
  \le \|f_{\psi}\|^{2}_{H^{\alpha}}\,\|f\|^{2}_{H^{\alpha}}
  \le \|f_{\psi}\|^{2}_{H^{\alpha}}
  = \|\psi\|^{2}_{H^{-\alpha}},
\end{equation*}
with equality when $f=f_{\psi}/\|f_{\psi}\|_{H^{\alpha}}$.
Hence $\sup\{|\psi(f)|^{2}:\|f\|_{H^{\alpha}}\le 1\}=\|f_{\psi}\|^{2}_{H^{\alpha}}=\|\psi\|^{2}_{H^{-\alpha}}$,
which is the second identity in \eqref{i1}.
\end{proof}

\begin{remark}[Identification with $H^{-2\alpha}(K)$]\label{rem:sobolev-comparison}
The Sobolev spaces of Definition~\ref{def:sobolev-K} are directly related to the Dirichlet Sobolev
space of Definition~\ref{def:sobolev-D}:
\begin{equation}\label{eq:sobolev-comparison}
  \mathcal{W}^{-\alpha}_{D}(U)
  = H^{-2\alpha}(K)
  \quad\text{(isometrically).}
\end{equation}
In particular, $\mathcal{H}^{D}_{\alpha}=L^{2}(\R_{+},H^{-2\alpha}(K))$, and the noise
$W^{D}_{\alpha}$ can equivalently be viewed as an isonormal Gaussian process on
$L^{2}(\R_{+},H^{-2\alpha}(K))$.
\end{remark}

\begin{proof}
For $f\in\mathcal{S}_{D}(U)$, the second line of~\eqref{eq:sobolev-D} gives
\begin{equation*}
  \|f\|^{2}_{\mathcal{W}^{-\alpha}_{D}(U)}
  = \int_{U}\bigl[(-\Delta_{D,U})^{-\alpha}f\bigr]^{2}\,d\mu
  = \sum_{j=1}^{\infty}\lambda_{j}^{-2\alpha}\,\hat{f}_{j}^{2},
\end{equation*}
where the second step uses $(-\Delta_{D,U})^{-\alpha}f = \sum_{j}\lambda_{j}^{-\alpha}\hat{f}_{j}\Phi_{j}$
and the $L^{2}(U,\mu)$-orthonormality of $(\Phi_{j})$ (Notation~\ref{not:spectrum-D}).
Applying~\eqref{eq:H-neg-alpha-norm} with $\alpha$ replaced by $2\alpha$
and noting that $f(\Phi_{j}) = \hat{f}_{j}$ for $f\in\mathcal{S}_{D}(U)$,
\begin{equation*}
  \|f\|^{2}_{H^{-2\alpha}(K)}
  = \sum_{j=1}^{\infty}\lambda_{j}^{-2\alpha}\,\hat{f}_{j}^{2}.
\end{equation*}
Since both norms coincide on the dense subspace $\mathcal{S}_{D}(U)$, the two completions are isometrically equal.
\end{proof}

\smallskip

We now turn to the Sobolev regularity of $W^{D}_{\alpha}$, which will be more easily expressed thanks to another piece of notation.

\begin{notation}[Time-integrated noise]\label{not:time-int-noise}
For $t>0$, the \emph{time-integrated noise} $W^{D}_{\alpha}[t]$ is the centred Gaussian family
\begin{equation}\label{eq:time-int-noise}
  W^{D}_{\alpha}[t]
  := \bigl\{W^{D}_{\alpha}(\mathbf{1}_{[0,t]}\otimes f)
  : f\in\mathcal{S}_{D}(U)\bigr\}.
\end{equation}
We also set $\xi_{j}(t):=W^{D}_{\alpha}(\mathbf{1}_{[0,t]}\otimes\Phi_{j})$.
\end{notation}

The key consequence of Notation~\ref{not:time-int-noise} is that the variables $\xi_{j}(t)$ are
independent centred Gaussians, which makes the Sobolev regularity of $W^{D}_{\alpha}[t]$
tractable.

\begin{proposition}[Sobolev regularity of $W^{D}_{\alpha}$]\label{prop:sobolev-reg}
Let $0\le\alpha\le d_{h}/(2d_{w})$ and $t>0$.  The time-integrated noise $W^{D}_{\alpha}[t]$ (Notation~\ref{not:time-int-noise})
belongs almost surely to $H^{-\beta}(K)$ for every $\beta>d_{h}/d_{w}-2\alpha$.
\end{proposition}

\begin{proof}
Recall that $(\lambda_{j},\Phi_{j})_{j\ge 1}$ are the eigenpairs of $-\Delta_{D,U}$
(Notation~\ref{not:spectrum-D}).
The isometry~\eqref{eq:geom-W-isometry-D} and the inner product formula~\eqref{eq:sobolev-D} give
\begin{equation}\label{eq:xi-j-variance}
  \E\bigl[\xi_{j}(t)\,\xi_{k}(t)\bigr]
  = t\,\langle\Phi_{j},\Phi_{k}\rangle_{\mathcal{W}^{-\alpha}_{D}}
  = t\,\lambda_{j}^{-2\alpha}\delta_{jk},
\end{equation}
so the $\xi_{j}(t)$ are independent centred Gaussians with variance $t\lambda_{j}^{-2\alpha}$.
Lemma~\ref{lem:sobolev-K-dual} then identifies $\|W^{D}_{\alpha}[t]\|^{2}_{H^{-\beta}}
=\sum_{j}\lambda_{j}^{-\beta}\xi_{j}(t)^{2}$.  Taking expectations and using~\eqref{eq:xi-j-variance} enables to write
\begin{equation*}
  \E\bigl[\|W^{D}_{\alpha}[t]\|^{2}_{H^{-\beta}}\bigr]
  = \sum_{j=1}^{\infty}\lambda_{j}^{-\beta}\,\E\bigl[\xi_{j}(t)^{2}\bigr]
  = t\,\sum_{j=1}^{\infty}\lambda_{j}^{-\beta-2\alpha}.
\end{equation*}
By Lemma~\ref{lem:weyl}, this series converges if and only if $\beta+2\alpha>d_{h}/d_{w}$, i.e.,
$\beta>d_{h}/d_{w}-2\alpha$.  For $\beta$ strictly above this threshold,
$\E[\|W^{D}_{\alpha}[t]\|^{2}_{H^{-\beta}}]<\infty$ and hence $W^{D}_{\alpha}[t]\in H^{-\beta}(K)$
almost surely.
\end{proof}

As $\alpha$ increases toward $d_{h}/(2d_{w})$, the threshold $d_{h}/d_{w}-2\alpha$ decreases to
$0$.  At the log-correlated value $\alpha=d_{h}/(2d_{w})$
(Remark~\ref{rem:geom-noise-reg}),
$W^{D}_{\alpha}[t]$ belongs to $H^{-\beta}(K)$ for every $\beta>0$, but not to $L^{2}(U,\mu)
=H^{0}(K)$.
The regularity at this critical value is more delicate to handle, and we skip this topic for the sake of conciseness.
When $\alpha>d_{h}/(2d_{w})$, Lemma~\ref{lem:weyl} gives $\sum_{j}\lambda_{j}^{-2\alpha}<\infty$,
which implies that $G^{D,U}_{2\alpha}(x,x)=\sum_{j}\lambda_{j}^{-2\alpha}\Phi_{j}(x)^{2}$ is
finite for each $x\in U$, and hence $G^{D,U}_{\alpha}(x,\cdot)\in L^{2}(U,\mu)$ for each $x$.
This allows one to evaluate $W^{D}_{\alpha}$ pointwise in space, leading to the following definition.

\begin{definition}[Density field of $W^{D}_{\alpha}$]\label{def:density-field}
For $\alpha>d_{h}/(2d_{w})$, the \emph{density field} of $W^{D}_{\alpha}$ is the process
\begin{equation}\label{eq:density-field}
  \tilde{W}^{D}_{\alpha}(t,x)
  := W_{0}\bigl(\mathbf{1}_{[0,t]}\otimes G^{D,U}_{\alpha}(x,\cdot)\bigr),
  \qquad t>0,\ x\in U,
\end{equation}
where $W_{0}$ is the white noise of Definition~\ref{def:geom-noise-white}.
By the bilinear formula~\eqref{eq:RK-D-bilinear}, its covariance is
\begin{equation}\label{eq:density-cov}
  \E\bigl[\tilde{W}^{D}_{\alpha}(t,x)\,\tilde{W}^{D}_{\alpha}(s,y)\bigr]
  = \min(t,s)\,G^{D,U}_{2\alpha}(x,y),
\end{equation}
which coincides with the covariance of $W^{D}_{\alpha}$ applied to point masses. The
field $\tilde{W}^{D}_{\alpha}$ is thus the pointwise spatial density of the noise.
\end{definition}

With the density field in hand, we turn to its spatial regularity.
Since $\tilde{W}^{D}_{\alpha}(t,x)=W_{0}(\mathbf{1}_{[0,t]}\otimes G^{D,U}_{\alpha}(x,\cdot))$,
the regularity of $x\mapsto\tilde{W}^{D}_{\alpha}(t,x)$ is governed by the regularity of
$x\mapsto G^{D,U}_{\alpha}(x,\cdot)$ in $L^{2}(U,\mu)$, which in turn reflects the H\"older
continuity of the operator $(-\Delta_{D,U})^{-\alpha}:L^{2}(U,\mu)\to C(U)$.  The natural
H\"older exponent for this operator is
\begin{equation}\label{eq:Holder-exp}
  H_{\alpha} := \min\Bigl(\alpha\,d_{w}-\frac{d_{h}}{2},\;d_{w}-d_{h}\Bigr),
\end{equation}
which is positive for $\alpha>d_{h}/(2d_{w})$, grows linearly in $\alpha$ in the subcritical
regime, and saturates at the value $d_{w}-d_{h}$ once $\alpha\ge 1-d_{h}/(2d_{w})$.
The following theorem makes this precise.

\begin{theorem}[H\"older regularity of $(-\Delta_{D,U})^{-\alpha}$]\label{thm:holder-frac-lap}
Let $\alpha>d_{h}/(2d_{w})$, and let $H_{\alpha}$ be as in~\eqref{eq:Holder-exp}.  There exists
a constant $C>0$ such that for every $f\in L^{2}(U,\mu)$ and $x,y\in U$ (where $d$ is the metric of Definition~\ref{def:mmspace}, which can be reduced to the Euclidean distance in the case of the Sierpi\'nski gasket; see Section~\ref{ssec:sierpinski}):
\begin{enumerate}
\item If $\alpha\ne 1-d_{h}/(2d_{w})$:
\begin{equation}\label{eq:holder-frac-lap-noncrit}
  \bigl|(-\Delta_{D,U})^{-\alpha}f(x)-(-\Delta_{D,U})^{-\alpha}f(y)\bigr|
  \;\le\; C\,d(x,y)^{H_{\alpha}}\,\|f\|_{L^{2}(U,\mu)}.
\end{equation}

\item If $\alpha=1-d_{h}/(2d_{w})$, for $d(x,y)\le 1/2$:
\begin{equation}\label{eq:holder-frac-lap-crit}
  \bigl|(-\Delta_{D,U})^{-\alpha}f(x)-(-\Delta_{D,U})^{-\alpha}f(y)\bigr|
  \;\le\; C\,d(x,y)^{H_{\alpha}}\,|\ln d(x,y)|\,\|f\|_{L^{2}(U,\mu)}.
\end{equation}
\end{enumerate}
\end{theorem}

\begin{proof}
See~\cite[Theorem~3.4]{baudoin:23:dirichlet}.  Throughout, the integral
representation~\eqref{eq:RK-D} gives
\begin{equation*}
  \bigl|(-\Delta_{D,U})^{-\alpha}f(x)-(-\Delta_{D,U})^{-\alpha}f(y)\bigr|
  \;\le\; \frac{1}{\Gamma(\alpha)}\int_{0}^{\infty}t^{\alpha-1}
  \bigl|P^{D,U}_{t}f(x)-P^{D,U}_{t}f(y)\bigr|\,dt.
\end{equation*}
We will use two bounds.  The first one is a pointwise bound. Namely, owing to
Cauchy--Schwarz and Corollary~\ref{cor:L2-bound-D}, it is readily checked that
\begin{equation}\label{eq:ptwise-semi}
  |P^{D,U}_{t}f(x)|
  \;\le\; \Bigl(\int_{U}p^{D,U}_{t}(x,z)^{2}\,d\mu(z)\Bigr)^{1/2}\|f\|_{L^{2}}
  \;\le\; C\,\frac{e^{-\lambda_{1}t}}{\min(1,t^{d_{h}/d_{w}})^{1/2}}\,\|f\|_{L^{2}}.
\end{equation}
In particular, for $t\le 1$ this simplifies to
$|P^{D,U}_{t}f(x)|\le C\,t^{-d_{h}/(2d_{w})}\|f\|_{L^{2}}$.
The second is the semigroup H\"older bound (see~\cite{baudoin:23:dirichlet,barlow:98:diffusions}):
\begin{equation}\label{eq:semi-holder}
  \bigl|P^{D,U}_{t}f(x)-P^{D,U}_{t}f(y)\bigr|
  \;\le\; C\,e^{-\lambda_{1}t/2}\,\frac{d(x,y)^{d_{w}-d_{h}}}{t^{1-d_{h}/(2d_{w})}}\,\|f\|_{L^{2}}.
\end{equation}
Set $\delta:=d(x,y)^{d_{w}}$, so $\delta\le 1$ when $d(x,y)\le 1$. We now divide the proof in several steps.

\smallskip

\noindent
{\it Step 1: $\alpha<1-d_{h}/(2d_{w})$.}  Split the time integral at $\delta$:
\begin{equation}\label{eq:split-I1-I2}
  \int_{0}^{\infty}t^{\alpha-1}\bigl|P^{D,U}_{t}f(x)-P^{D,U}_{t}f(y)\bigr|\,dt
  \;=\; I_{1}+I_{2},
\end{equation}
where $I_{1}:=\int_{0}^{\delta}t^{\alpha-1}|P^{D,U}_{t}f(x)-P^{D,U}_{t}f(y)|\,dt$ and $I_{2}:=\int_{\delta}^{\infty}t^{\alpha-1}|P^{D,U}_{t}f(x)-P^{D,U}_{t}f(y)|\,dt$.

For $I_{1}$, bound the integrand by $|P^{D,U}_{t}f(x)|+|P^{D,U}_{t}f(y)|$ and
apply~\eqref{eq:ptwise-semi} with $t\le\delta\le 1$. We obtain
\begin{equation*}
  I_{1} \;\le\; C\int_{0}^{\delta}t^{\alpha-1-d_{h}/(2d_{w})}\,dt\,\|f\|_{L^{2}}
  \;=\; \frac{C}{\alpha-d_{h}/(2d_{w})}\,\delta^{\alpha-d_{h}/(2d_{w})}\,\|f\|_{L^{2}},
\end{equation*}
where convergence at $0$ uses $\alpha>d_{h}/(2d_{w})$.  Substituting $\delta=d(x,y)^{d_{w}}$, we discover that
\begin{equation*}
  I_{1} \;\le\; C\,d(x,y)^{d_{w}(\alpha-d_{h}/(2d_{w}))}\,\|f\|_{L^{2}}
  \;=\; C\,d(x,y)^{\alpha d_{w}-d_{h}/2}\,\|f\|_{L^{2}}
  \;=\; C\,d(x,y)^{H_{\alpha}}\,\|f\|_{L^{2}}.
\end{equation*}

For $I_{2}$ (the second term in~\eqref{eq:split-I1-I2}), apply~\eqref{eq:semi-holder} and drop $e^{-\lambda_{1}t/2}\le 1$. This yields
\begin{equation*}
  I_{2} \;\le\; C\,d(x,y)^{d_{w}-d_{h}}\int_{\delta}^{\infty}
  t^{\alpha-2+d_{h}/(2d_{w})}\,dt\,\|f\|_{L^{2}}.
\end{equation*}
Since $\alpha<1-d_{h}/(2d_{w})$, the exponent satisfies $\alpha-2+d_{h}/(2d_{w})<-1$, so
the integral converges at $\infty$:
\begin{equation*}
  \int_{\delta}^{\infty}t^{\alpha-2+d_{h}/(2d_{w})}\,dt
  \;=\; \frac{\delta^{\alpha-1+d_{h}/(2d_{w})}}{1-\alpha-d_{h}/(2d_{w})}.
\end{equation*}
Substituting $\delta=d(x,y)^{d_{w}}$ and computing the exponent, we end up with
\begin{eqnarray*}
  I_{2} &\le& C\,d(x,y)^{d_{w}-d_{h}}\cdot d(x,y)^{d_{w}(\alpha-1+d_{h}/(2d_{w}))}\,\|f\|_{L^{2}} \\
  &=& C\,d(x,y)^{\alpha d_{w}-d_{h}/2}\,\|f\|_{L^{2}}
  \;=\; C\,d(x,y)^{H_{\alpha}}\,\|f\|_{L^{2}}.
\end{eqnarray*}
Plugging the bounds on $I_{1}$ and $I_{2}$ into~\eqref{eq:split-I1-I2} yields~\eqref{eq:holder-frac-lap-noncrit} in this sub-case.

\smallskip

\noindent
{\it Step 2: $\alpha>1-d_{h}/(2d_{w})$.}  No splitting is needed for this case.
Specifically, since $\alpha-2+d_{h}/(2d_{w})>-1$, applying~\eqref{eq:semi-holder} over all of
$[0,\infty)$ gives
\begin{equation}\label{i2}
  \bigl|(-\Delta_{D,U})^{-\alpha}f(x)-(-\Delta_{D,U})^{-\alpha}f(y)\bigr|
  \;\le\; C\,d(x,y)^{d_{w}-d_{h}}\int_{0}^{\infty}t^{\alpha-2+d_{h}/(2d_{w})}
  e^{-\lambda_{1}t/2}\,dt\,\|f\|_{L^{2}}.
\end{equation}
The remaining integral is a Gamma integral, namely
\begin{equation*}
\int_{0}^{\infty}t^{\alpha-2+d_{h}/(2d_{w})}e^{-\lambda_{1}t/2}\,dt
=(\lambda_{1}/2)^{-(\alpha-1+d_{h}/(2d_{w}))}\Gamma(\alpha-1+d_{h}/(2d_{w}))<\infty .
\end{equation*}
Since $H_{\alpha}=d_{w}-d_{h}$ in this regime, relation~\eqref{i2} trivially yields~\eqref{eq:holder-frac-lap-noncrit}.

\smallskip

\noindent
{\it Step 3: $\alpha=1-d_{h}/(2d_{w})$.}  Use the same split~\eqref{eq:split-I1-I2} at $\delta$.
The bound for $I_{1}$ is identical to Step~1 and gives $I_{1}\le C\,d(x,y)^{d_{w}-d_{h}}\,\|f\|_{L^{2}}$.
For $I_{2}$, observe that $\alpha-2+d_{h}/(2d_{w})=-1$, so~\eqref{eq:semi-holder} gives
\begin{equation*}
  I_{2} \;\le\; C\,d(x,y)^{d_{w}-d_{h}}\int_{\delta}^{\infty}t^{-1}\,e^{-\lambda_{1}t/2}\,dt
  \,\|f\|_{L^{2}}.
\end{equation*}
Splitting the remaining integral at $t=1$ and bounding $e^{-\lambda_{1}t/2}\le 1$ on
$[\delta,1]$, we get
\begin{equation*}
  \int_{\delta}^{1}t^{-1}\,dt+\int_{1}^{\infty}t^{-1}e^{-\lambda_{1}t/2}\,dt
  \;=\; |\ln\delta|+C_{1}
  \;\le\; C\,|\ln\delta|,
\end{equation*}
for $\delta\le 1/2$.  Substituting $\delta=d(x,y)^{d_{w}}$ gives $|\ln\delta|=d_{w}\,|\ln d(x,y)|$,
so 
\begin{equation*}
I_{2}\le C\,d(x,y)^{d_{w}-d_{h}}\,|\ln d(x,y)|\,\|f\|_{L^{2}}.
\end{equation*}
Since $I_{2}$ dominates $I_{1}$, this yields~\eqref{eq:holder-frac-lap-crit}.

\smallskip

Combining Steps~1-3, in all three regimes we have obtained
\begin{equation*}
\bigl|(-\Delta_{D,U})^{-\alpha}f(x)-(-\Delta_{D,U})^{-\alpha}f(y)\bigr|
\le C\,d(x,y)^{H_{\alpha}}\,\|f\|_{L^{2}(U,\mu)} ,
\end{equation*}
with an additional factor $|\ln d(x,y)|$ at criticality. This is
exactly the content of~\eqref{eq:holder-frac-lap-noncrit}--\eqref{eq:holder-frac-lap-crit}.
\end{proof}

By $L^{2}(U,\mu)$ duality, the H\"older bounds~\eqref{eq:holder-frac-lap-noncrit}--\eqref{eq:holder-frac-lap-crit} on $(-\Delta_{D,U})^{-\alpha}$ transfer to the $L^{2}$-increments of the Riesz kernel $G^{D,U}_{\alpha}$.

\begin{corollary}[H\"older regularity of $G^{D,U}_{\alpha}$]\label{cor:holder-riesz}
Under the hypotheses of Theorem~\ref{thm:holder-frac-lap}, with $H_{\alpha}$ as
in~\eqref{eq:Holder-exp}, there exists $C>0$ such that for all $x,y\in U$:
\begin{enumerate}
\item If $\alpha\ne 1-d_{h}/(2d_{w})$:
\begin{equation}\label{eq:holder-riesz-noncrit}
  \int_{U}\bigl(G^{D,U}_{\alpha}(x,z)-G^{D,U}_{\alpha}(y,z)\bigr)^{2}\,d\mu(z)
  \;\le\; C\,d(x,y)^{2H_{\alpha}}.
\end{equation}

\item If $\alpha=1-d_{h}/(2d_{w})$, for $d(x,y)\le 1/2$:
\begin{equation}\label{eq:holder-riesz-crit}
  \int_{U}\bigl(G^{D,U}_{\alpha}(x,z)-G^{D,U}_{\alpha}(y,z)\bigr)^{2}\,d\mu(z)
  \;\le\; C\,d(x,y)^{2H_{\alpha}}\,|\ln d(x,y)|^{2}.
\end{equation}
\end{enumerate}
\end{corollary}

\begin{proof}
The action formula~\eqref{eq:RK-D-action} gives, for any $f\in L^{2}(U,\mu)$,
\begin{equation*}
  \int_{U}\bigl(G^{D,U}_{\alpha}(x,z)-G^{D,U}_{\alpha}(y,z)\bigr)f(z)\,d\mu(z)
  = (-\Delta_{D,U})^{-\alpha}f(x)-(-\Delta_{D,U})^{-\alpha}f(y).
\end{equation*}
Taking the supremum over $\|f\|_{L^{2}}=1$ identifies the $L^{2}(U,\mu)$-norm of
$G^{D,U}_{\alpha}(x,\cdot)-G^{D,U}_{\alpha}(y,\cdot)$ with the operator modulus of continuity of
$(-\Delta_{D,U})^{-\alpha}$, and the three bounds follow directly from
Theorem~\ref{thm:holder-frac-lap}.
\end{proof}

The corollary controls the second moment of the increment of the density field:
by~\eqref{eq:density-cov} and Corollary~\ref{cor:holder-riesz},
\begin{equation}\label{eq:field-increment-var}
  \E\bigl[(\tilde{W}^{D}_{\alpha}(t,x)-\tilde{W}^{D}_{\alpha}(t,y))^{2}\bigr]
  = t\int_{U}\bigl(G^{D,U}_{\alpha}(x,z)-G^{D,U}_{\alpha}(y,z)\bigr)^{2}\,d\mu(z).
\end{equation}
Applying the entropy method of~\cite{adler:07:random} to this Gaussian field yields the
following result.

\begin{theorem}[H\"older continuity of $\tilde{W}^{D}_{\alpha}$]\label{thm:holder-field}
Let $\alpha>d_{h}/(2d_{w})$ and $t>0$.  Define the modulus
\begin{equation}\label{eq:modulus}
  \omega_{\alpha}(r)
  :=
  \begin{cases}
    r^{H_{\alpha}}\,\sqrt{|\ln r|}, & \alpha\ne 1-d_{h}/(2d_{w}), \\[4pt]
    r^{d_{w}-d_{h}}\,|\ln r|^{3/2}, & \alpha=1-d_{h}/(2d_{w}),
  \end{cases}
\end{equation}
with $H_{\alpha}$ as in~\eqref{eq:Holder-exp}.  There exists an almost surely continuous
version of $x\mapsto\tilde{W}^{D}_{\alpha}(t,x)$, still denoted $\tilde{W}^{D}_{\alpha}(t,\cdot)$, such that
\begin{equation}\label{eq:holder-field}
  \lim_{\delta\to 0}\sup_{\substack{x,y\in U,\\0<d(x,y)\le\delta}}
  \frac{|\tilde{W}^{D}_{\alpha}(t,x)-\tilde{W}^{D}_{\alpha}(t,y)|}{\omega_{\alpha}(d(x,y))}
  < +\infty
  \quad\text{almost surely.}
\end{equation}
\end{theorem}

\begin{proof}
See~\cite[Thm.~3.8]{baudoin:23:dirichlet}.  Set
$\rho(x,y):=\E[(\tilde{W}^{D}_{\alpha}(t,x)-\tilde{W}^{D}_{\alpha}(t,y))^{2}]^{1/2}$.
By~\eqref{eq:field-increment-var} and Corollary~\ref{cor:holder-riesz}, one has
$\rho(x,y)\le C\sqrt{t}\,d(x,y)^{H_{\alpha}}$ (with an additional logarithmic factor at
$\alpha=1-d_{h}/(2d_{w})$).  Let $\mathcal{N}_{\rho}(\varepsilon)$ denote the smallest number of
$\rho$-balls of radius $\varepsilon$ needed to cover $U$.  Since $\rho(x,y)\le
C\sqrt{t}\,d(x,y)^{H_{\alpha}}$, the metric $\rho$ and the intrinsic metric $d$ are
H\"older-equivalent, giving $\mathcal{N}_{\rho}(\varepsilon)=O(\varepsilon^{-d_{h}/H_{\alpha}})$.
The Dudley entropy integral bound~\cite[Theorem~1.3.5]{adler:07:random} then yields
\begin{equation*}
  \sup_{\rho(x,y)\le\tau}|\tilde{W}^{D}_{\alpha}(t,x)-\tilde{W}^{D}_{\alpha}(t,y)|
  \;\le\; C'\int_{0}^{\tau}\sqrt{\ln\mathcal{N}_{\rho}(\varepsilon)}\,d\varepsilon,
\end{equation*}
and an integration-by-parts argument applied to the right-hand side gives the modulus
$\omega_{\alpha}$; see the proof of~\cite[Thm.~3.8]{baudoin:23:dirichlet} for the complete computation.
\end{proof}

\begin{remark}
The Hölder regularity result of Theorem~\ref{thm:holder-field} has a natural counterpart
in the Euclidean setting that we have not developed explicitly.  If $W$ is a Gaussian
noise on $\R^{d}$ with covariance $\E[W(A)W(B)]=\int_{A\times B}f(x-y)\,dx\,dy$
for a smooth non-negative definite function $f$, then the associated density field
$x\mapsto W(f_{x})$ (where $f_{x}(\cdot)=f(x-\cdot)$) is a stationary Gaussian process
whose regularity is governed by the decay of $f$ near the origin.  This type of analysis
is standard for Gaussian processes on $\R^{d}$ (see Section~\ref{sec:noise-def} for the
framework used in this book) and can be carried out by the same entropy method as
in the proof of Theorem~\ref{thm:holder-field}.  We have chosen not to include it here
for the sake of conciseness.
\end{remark}

\section{Existence and uniqueness}\label{sec:geom-existence}

This section establishes existence and uniqueness for the parabolic Anderson
model on a bounded domain $U\subset X$.  In Chapter~\ref{chap:existence} we proved this result for the PAM on $\R^{d}$
via a Banach fixed-point argument (Theorem~\ref{thm:global-wp}), and in
Section~\ref{sec:chaos} we revisited it via the Wiener chaos expansion
(Theorem~\ref{thm:chaos-eu}).  We follow the chaos
route here, since it simultaneously yields the $L^{2}$-moment formula that drives
the moment estimates of Section~\ref{sec:geom-moments}.

All the necessary ingredients are now available: the Dirichlet noise
$W^{D}_{\alpha}$ and its Hilbert space $\mathcal{H}^{D}_{\alpha}$
(Section~\ref{sec:geom-noise}), the Dirichlet heat kernel $p^{D,U}_{t}$ with
its $L^{\infty}$ bound (Corollary~\ref{cor:Linf-bound-D}), the Riesz kernel
$G^{D,U}_{2\alpha}$ with the $L^{2}$-operator estimate of
Lemma~\ref{lem:DeltaHK-D}, and the chaos framework adapted to the geometric
setting (Section~\ref{ssec:geom-noise-chaos}).  The key challenge which
distinguishes the geometric proof from the Euclidean one is the absence of
Fourier analysis. Indeed, the norms
$\|f_{k}\|^{2}_{(\mathcal{H}^{D}_{\alpha})^{\otimes k}}$ cannot be evaluated
through a spectral measure as in Theorem~\ref{thm:chaos-eu}.  Instead,
Lemma~\ref{lem:DeltaHK-D} serves as a substitute, leading to a renewal-type
recursion in physical space that mimics the role of Dalang's
condition~\eqref{eq:dalang-cond}.

Section~\ref{ssec:geom-mild} sets up the mild equation, recalls the solution
concept, and reduces existence and uniqueness to the convergence of a chaos
series.  Section~\ref{ssec:geom-chaos-est} establishes that convergence via the
renewal argument and deduces the $L^{2}$-moment bound.

\subsection{Mild formulation and chaos characterisation}
\label{ssec:geom-mild}

Recall that we now turn to the parabolic
Anderson model on the bounded domain $U\subset X$.  The equation of interest is
\begin{equation}\label{eq:pam-D-geom}
  \partial_{t}u(t,x)
  \;=\;
  \Delta_{D,U}\,u(t,x)
  \;+\;
  \beta\,u(t,x)\,\dot{W}^{D}_{\alpha}(t,x),
  \qquad (t,x)\in\R_{+}\times U,
\end{equation}
where $\beta\ge 0$, $\alpha\ge 0$, and the initial condition is
$u(0,\cdot)=u_{0}\in L^{\infty}(U,\mu)$.  For $\alpha=0$ the noise is the
space-time white noise $W_{0}$ of Definition~\ref{def:geom-noise-white}; for
$\alpha>0$ it is the colored noise of Definition~\ref{def:geom-noise-D}.

Equation~\eqref{eq:pam-D-geom} is solved in its \emph{mild form}.  That is,
using the Dirichlet semigroup $P^{D,U}_{t}$ and heat kernel $p^{D,U}_{t}$
(Theorem~\ref{thm:hk-diri}), we rewrite~\eqref{eq:pam-D-geom} as the integral
equation
\begin{equation}\label{eq:mild-pam-D-geom}
  u(t,x)
  \;=\;
  P^{D,U}_{t}u_{0}(x)
  \;+\;
  \beta\int_{[0,t]\times U}p^{D,U}_{t-s}(x,y)\,u(s,y)\,W^{D}_{\alpha}(ds,dy).
\end{equation}
The stochastic integral in~\eqref{eq:mild-pam-D-geom} is understood in the
It\^{o}--Walsh sense. Indeed, $W^{D}_{\alpha}$ is a worthy martingale measure, with
covariation measure $G^{D,U}_{2\alpha}(x,y)\,d\mu(x)\,d\mu(y)\,ds$ (see
relation~\eqref{eq:geom-cov-D-kernel}), which is non-negative. Hence Walsh's
construction of Section~\ref{sec:integration}, as well as the stochastic Fubini
theorem of Proposition~\ref{prop:stoch-fubini}, apply verbatim, with
$\Lambda(x-y)\,dx\,dy$ replaced by $G^{D,U}_{2\alpha}(x,y)\,d\mu(x)\,d\mu(y)$. Since the integrand $p^{D,U}_{t-s}(x,y)u(s,y)$ is
adapted, Proposition~\ref{prop:ito-sko} identifies it with the Skorohod integral
$\delta$ (see also Section~\ref{ssec:geom-noise-chaos}).
Observe that
equation~\eqref{eq:mild-pam-D-geom} is the geometric counterpart of the
Euclidean mild equation~\eqref{eq:mild-pam}:
\begin{equation}\label{eq:mild-pam-recall}
  u(t,x) = P_{t}u_{0}(x)
  + \beta\int_{0}^{t}\!\int_{\R^{d}}
    p_{t-s}(x-y)\,u(s,y)\,W(ds,dy).
\end{equation}
Both equations have the same Duhamel (variation-of-constants) structure, but 
one has to handle two major changes when passing from $\R^{d}$ to~$(U,d_{U},\mu)$.

\begin{enumerate}
\item {\it Heat kernel.}
In~\eqref{eq:mild-pam-recall} the free kernel $p_{t-s}(x-y)$ depends only on
the difference $x-y$ (translation invariance) and is a probability density for
all times.  In~\eqref{eq:mild-pam-D-geom} the kernel $p^{D,U}_{t-s}(x,y)$ is
not translation-invariant (it depends on the geometry of $U$), and it satisfies
the Dirichlet boundary condition, which forces an exponential decay
$p^{D,U}_{t}(x,y)\le Ce^{-\lambda_{1}t}$ for large $t$ (Theorem~\ref{thm:hk-diri},
equation~\eqref{eq:subGauss-D-upper}).

\item {\it Noise and spatial covariance.}
In~\eqref{eq:mild-pam-recall} the noise $W$ has spatial covariance given by the
kernel $\Lambda$ (see~\eqref{eq:noise-cov}), and the natural space-time Hilbert
space $\mathcal{H}$ is built from the spectral measure $\mu$ via Fourier
analysis.  In~\eqref{eq:mild-pam-D-geom} the noise $W^{D}_{\alpha}$ has spatial
covariance $G^{D,U}_{2\alpha}$ (Definition~\ref{def:geom-noise-D},
equation~\eqref{eq:geom-cov-D-kernel}), and $\mathcal{H}^{D}_{\alpha}$ is the
space-time Sobolev space $L^{2}(\R_{+},\mathcal{W}^{-\alpha}_{D}(U))$.  In the
Euclidean setting the integrability of the free kernel against $\mathcal{H}$ is
guaranteed by Dalang's condition~\eqref{eq:dalang-cond}; in the geometric
setting it is guaranteed by the $L^{2}(U,\mu)$-norm bound of
Lemma~\ref{lem:DeltaHK-D}, which plays the same role.
\end{enumerate}

Despite these differences, the notion of solution carries over without change,
with $W^{D}_{\alpha}$ and $\mathcal{H}^{D}_{\alpha}$ in place of $W$ and
$\mathcal{H}$. We spell it out below.

\begin{definition}[It\^{o}--Skorohod solution]\label{def:ito-sol-geom}
A process $u=\{u(t,x);\,(t,x)\in\R_{+}\times U\}$ is a
\emph{random field solution} of~\eqref{eq:pam-D-geom} in the It\^{o}--Skorohod
sense if the following five conditions are satisfied.
\begin{enumerate}
\item The process $u$ is adapted to the filtration generated by
  $W^{D}_{\alpha}$.
\item The map $(t,x)\mapsto u(t,x)$ is jointly measurable.
\item We have $\E[I(t,x)^{2}]<\infty$ for all $(t,x)$, where $I(t,x)$ denotes
  the stochastic integral in~\eqref{eq:mild-pam-D-geom}.
\item The map $(t,x)\mapsto I(t,x)$ is continuous in $L^{2}(\Omega)$.
\item The process $u$ satisfies~\eqref{eq:mild-pam-D-geom} almost surely for
  all $(t,x)$.
\end{enumerate}
\end{definition}

\smallskip

\noindent
\textbf{Chaos kernels.}
By a standard iteration of~\eqref{eq:mild-pam-D-geom} and the chaos
decomposition of Section~\ref{ssec:geom-noise-chaos} (see also
Section~\ref{ssec:chaos-kernels} for the analogous computation in the Euclidean
setting), any solution $u$ has the Wiener chaos expansion
\begin{equation}\label{eq:u-chaos-geom}
  u(t,x)
  \;=\;
  P^{D,U}_{t}u_{0}(x)
  \;+\;
  \sum_{k=1}^{\infty}I_{k}(f_{k}(\,\cdot\,,t,x)),
\end{equation}
where the chaos kernels $f_{k}\in(\mathcal{H}^{D}_{\alpha})^{\odot k}$ are given, on
the ordered simplex $\{0<s_{\sigma(1)}<\cdots<s_{\sigma(k)}<t\}$, by
\begin{equation}\label{eq:fk-geom}
  f_{k}(s_{1},y_{1},\ldots,s_{k},y_{k};\,t,x)
  =
  \frac{\beta^{k}}{k!}\,
  p^{D,U}_{t-s_{\sigma(k)}}(x,y_{\sigma(k)})\cdots
  p^{D,U}_{s_{\sigma(2)}-s_{\sigma(1)}}(y_{\sigma(2)},y_{\sigma(1)})\,
  P^{D,U}_{s_{\sigma(1)}}u_{0}(y_{\sigma(1)}),
\end{equation}
where $\sigma$ is the permutation of $\{1,\ldots,k\}$ such that
$0<s_{\sigma(1)}<\cdots<s_{\sigma(k)}<t$, and $f_{k}:=0$ off the simplex.  This is the
geometric analogue of~\eqref{eq:fn-def}.

The following proposition characterizes solutions via the square-summability of
the chaos series.  Its proof is identical to the Euclidean case
(Theorem~\ref{thm:chaos-eu}): the stochastic integral in~\eqref{eq:mild-pam-D-geom}
admits a Wiener chaos expansion by the same algebraic iteration, and orthogonality
of distinct chaos levels~\eqref{eq:geom-chaos-decomp} yields the $L^{2}$-moment
formula.

\begin{proposition}[Chaos characterisation]\label{prop:chaos-char-geom}
A process $u$ with $\E[u(t,x)^{2}]<\infty$ for all $(t,x)$ solves~\eqref{eq:mild-pam-D-geom}
in the sense of Definition~\ref{def:ito-sol-geom} if and only if, for every $(t,x)$, the
random variable $u(t,x)$ has the chaos expansion~\eqref{eq:u-chaos-geom}--\eqref{eq:fk-geom}
with
\begin{equation}\label{eq:chaos-finite-geom}
  \sum_{k=1}^{\infty}k!\,\|f_{k}(\,\cdot\,,t,x)\|^{2}_{(\mathcal{H}^{D}_{\alpha})^{\otimes k}}
  \;<\;+\infty.
\end{equation}
Here $\mathcal{H}^{D}_{\alpha}$ is the space-time Hilbert space introduced in
Definition~\ref{def:geom-noise-D}.
Moreover, the unique solution satisfies the $L^{2}$-moment formula
\begin{equation}\label{eq:L2-moment-geom}
  \E\bigl[u(t,x)^{2}\bigr]
  \;=\;
  (P^{D,U}_{t}u_{0}(x))^{2}
  \;+\;
  \sum_{k=1}^{\infty}k!\,\|f_{k}(\,\cdot\,,t,x)\|^{2}_{(\mathcal{H}^{D}_{\alpha})^{\otimes k}}.
\end{equation}
\end{proposition}

\subsection{Chaos estimates and existence}
\label{ssec:geom-chaos-est}

According to Proposition \eqref{prop:chaos-char-geom}, the key step to get existence and uniqueness results for~\eqref{eq:mild-pam-D-geom} is reduced to prove~\eqref{eq:chaos-finite-geom}.  Now recall that in the
Euclidean setting the norm $\|f_{n}\|^{2}_{\mathcal{H}^{\otimes n}}$ is
computed via Fourier transforms and the spectral measure $\mu$ (see
Theorem~\ref{thm:chaos-eu}). In our general geometric setting there is no natural
Fourier calculus.  Instead we use the $L^{2}$-estimate of
Lemma~\ref{lem:DeltaHK-D} as a substitute, establishing a recursive bound
directly in physical space.

\begin{lemma}[Kernel norm estimate]\label{lem:fk-bound-geom}
Let $\alpha\ge 0$ and $\rho>0$.  There exists a constant $C>0$ depending only
on $\alpha$, $d_{h}$, $d_{w}$, and $U$ such that for all $k\ge 1$, $t\ge 0$,
and $x\in U$,
\begin{align}
  \label{eq:fk-bound-small}
  \|f_{k}(\,\cdot\,,t,x)\|^{2}_{(\mathcal{H}^{D}_{\alpha})^{\otimes k}}
  &\;\le\;
  \frac{C^{k}\,\beta^{2k}\,\|u_{0}\|_{\infty}^{2}}{k!}\,
  \hat{\Psi}(\rho)^{k}\,e^{(\rho-2\lambda_{1})t},
\end{align}
where $\lambda_{1}=\lambda_{1}(U)$ is the first eigenvalue of $-\Delta_{D,U}$
(Notation~\ref{not:spectrum-D}) and $\hat{\Psi}(\rho)$ is the Laplace transform
\begin{equation}\label{eq:PsiHat-geom}
  \hat{\Psi}(\rho)
  :=
  \int_{0}^{\infty}e^{-\rho t}\,\Psi(t)\,dt,
  \quad\text{where}\quad
  \Psi(t)
  :=
  e^{2\lambda_{1}t}\,
  \sup_{x\in U}\|(-\Delta_{D,U})^{-\alpha}p^{D,U}_{t}(x,\cdot)\|^{2}_{L^{2}(U,\mu)}.
\end{equation}
\end{lemma}

\begin{proof}
We drop the superscripts $D,U$ for legibility and write $p_{t}$, $P_{t}$,
$\Delta$, $G_{2\alpha}$ for $p^{D,U}_{t}$, $P^{D,U}_{t}$, $\Delta_{D,U}$,
$G^{D,U}_{2\alpha}$ throughout.

\smallskip

\noindent
{\it Step 1: Reduction to $M_{k}(t,x)$.}
This step is the geometric counterpart of the simplex reduction performed in the
Euclidean setting (see~\eqref{eq:Kn-simplex-reduction}--\eqref{eq:Kn-simplex-value}).
Since $f_{k}\in(\mathcal{H}^{D}_{\alpha})^{\odot k}$ is symmetric in its $k$ pairs
$(s_{i},y_{i})$ by~\eqref{eq:fk-geom}, the full cube $[0,t]^{k}$ decomposes into
$k!$ congruent simplices, just as $[0,t]^{n}$ decomposed into $n!$ copies of
$T_{n}(t)$ (the notation comes from~\eqref{eq:simplex} in the Euclidean case).  Restricting to the ordered simplex
\begin{equation*}
  [0,t]^{k}_{<}
  \;:=\;
  \{(s_{1},\ldots,s_{k})\in[0,t]^{k} : 0<s_{1}<\cdots<s_{k}<t\},
\end{equation*}
and multiplying by $k!$ recovers the full norm, exactly as in~\eqref{eq:Kn-simplex-reduction}.
Moreover, on $[0,t]^{k}_{<}$ the permutation $\sigma$ in~\eqref{eq:fk-geom} is the identity, so
the heat-kernel chain simplifies to
\begin{equation}\label{eq:gk-chain}
  g_{k}(s,y;\,t,x)
  :=
  p_{t-s_{k}}(x,y_{k})\cdots p_{s_{2}-s_{1}}(y_{2},y_{1}),
\end{equation}
which is the geometric analogue of $K_{n}$ on the simplex~\eqref{eq:Kn-simplex-value}.
Thus, on $[0,t]^{k}_{<}$, the kernel~\eqref{eq:fk-geom} takes the explicit form
\begin{equation}\label{eq:fk-simplex}
  f_{k}(s_{1},y_{1},\ldots,s_{k},y_{k};\,t,x)
  \;=\;
  \frac{\beta^{k}}{k!}\,
  g_{k}(s,y;\,t,x)\,P_{s_{1}}u_{0}(y_{1}).
\end{equation}
We bound each of the two factors in the right-hand side of~\eqref{eq:fk-simplex} separately.
\begin{enumerate}
\item {\it The kernel factor $g_{k}(s,y;\,t,x)$.}
Since $f_{k}$ is symmetric, the norm
$\|f_{k}\|^{2}_{(\mathcal{H}^{D}_{\alpha})^{\otimes k}}$ equals $k!$ times the
integral of $f_{k}(s,y)\,\prod_{i}G_{2\alpha}(y_{i},y'_{i})\,f_{k}(s,y')$ over
$[0,t]^{k}_{<}\times U^{2k}$, exactly as in~\eqref{eq:Kn-simplex-reduction}.
After substituting~\eqref{eq:fk-simplex} for both occurrences of $f_{k}$, the factor of $k!$
from symmetry cancels one power of $k!$ in the denominator, and the remaining integral over the
chain of heat kernels $g_{k}$ defines $M_{k}(t,x)$ below in~\eqref{eq:Mk-geom}.

\item {\it The initial condition factor $P_{s_{1}}u_{0}(y_{1})$.}
In the Euclidean case this factor was bounded directly by
$\|u_{0}\|_{\infty}$~\eqref{eq:Ps-Linfty}, since the flat heat semigroup is an $L^{\infty}$
contraction.  Here we instead apply Corollary~\ref{cor:Linf-bound-D}, which gives
$|P_{s_{1}}u_{0}(y_{1})|\le Ce^{-\lambda_{1}s_{1}}\|u_{0}\|_{\infty}$. The
exponential decay $e^{-\lambda_{1}s_{1}}$ reflects the spectral gap of
$-\Delta_{D,U}$. It is also what eventually shifts the exponent from $\rho$ to
$\rho-2\lambda_{1}$ in~\eqref{eq:fk-bound-small}.
\end{enumerate}

\noindent
Gathering those two arguments into~\eqref{eq:fk-simplex}, we end up with the following upper bound:
\begin{equation}\label{eq:fk-Mk}
  \|f_{k}(\,\cdot\,,t,x)\|^{2}_{(\mathcal{H}^{D}_{\alpha})^{\otimes k}}
  \;\le\;
  C\,\frac{\beta^{2k}}{k!}\,\|u_{0}\|_{\infty}^{2}\,M_{k}(t,x),
\end{equation}
where
\begin{equation}\label{eq:Mk-geom}
  M_{k}(t,x)
  :=
  \int_{[0,t]^{k}_{<}}
  \!\int_{U^{2k}}
  g_{k}(s,y;\,t,x)\,\prod_{i=1}^{k}G_{2\alpha}(y_{i},y'_{i})\,
  g_{k}(s,y';\,t,x)\,e^{-2\lambda_{1}s_{1}}
  \,d\mu(y)\,d\mu(y')\,ds.
\end{equation}

\smallskip

\noindent
{\it Step 2: Recursion.}
The aim of this step is to prove that
\begin{equation}\label{eq:Mk-recursion}
  e^{2\lambda_{1}t}M_{k+1}(t,x)
  \;\le\;
  \int_{0}^{t}
    \Psi(t-s)
    \cdot
    e^{2\lambda_{1}s}\sup_{y\in U}M_{k}(s,y)
  \,ds ,
\end{equation}
where $\Psi$ is defined in~\eqref{eq:PsiHat-geom} and we initiate with
$M_{0}(t,x):=e^{-2\lambda_{1}t}$. 
To this aim, we will use the chain structure of $g_{k}$ in~\eqref{eq:gk-chain} and the covariance formula~\eqref{eq:geom-cov-D-kernel}.
 Indeed, starting from \eqref{eq:Mk-geom}, the following recursion is easily seen:
\begin{multline}\label{i3}
  e^{2\lambda_{1}t}M_{k+1}(t,x)
  =
  \int_{0}^{t}e^{2\lambda_{1}(t-s_{k+1})}
  \int_{U^{2}}
  p_{t-s_{k+1}}(x,y_{k+1}) G_{2\alpha}(y_{k+1},y'_{k+1})\,p_{t-s_{k+1}}(x,y'_{k+1})\\
  \cdot e^{2\lambda_{1}s_{k+1}}M_{k}(s_{k+1},y_{k+1})
  \,d\mu(y_{k+1})\,d\mu(y'_{k+1})\,ds_{k+1}. 
\end{multline}
Let us explain the exponential factors in~\eqref{i3}.
The factor $e^{2\lambda_{1}t}$ on the left-hand side is introduced to normalize $M_{k+1}$,
and on the right-hand side it is split as
$e^{2\lambda_{1}t}=e^{2\lambda_{1}(t-s_{k+1})}\cdot e^{2\lambda_{1}s_{k+1}}$:
the first part accompanies the outer heat kernels and will be absorbed into
$\Psi(t-s_{k+1})$ below, while the second part is grouped with $M_{k}(s_{k+1},y_{k+1})$
to form the normalized inner quantity $e^{2\lambda_{1}s_{k+1}}M_{k}(s_{k+1},y_{k+1})$
appearing in~\eqref{eq:Mk-recursion}.
The factor $e^{-2\lambda_{1}s_{1}}$ from definition~\eqref{eq:Mk-geom}, where $s_{1}$
is the smallest time in the inner $k$-step chain, is still present inside
$M_{k}(s_{k+1},y_{k+1})$ and is a distinct quantity from the outer prefactor
$e^{2\lambda_{1}s_{k+1}}$.

Since $M_{k}(s_{k+1},y_{k+1})$ depends on $y_{k+1}$, we cannot apply the bilinear
formula~\eqref{eq:RK-D-bilinear} to the double spatial integral directly.
We first bound $M_{k}(s_{k+1},y_{k+1})\le\sup_{y\in U}M_{k}(s_{k+1},y)$,
which decouples the two spatial integration variables $y_{k+1}$ and $y'_{k+1}$.
Once this substitution is made, the bilinear formula~\eqref{eq:RK-D-bilinear} identifies
the remaining double spatial integral as
\begin{equation}\label{eq:bilinear-p}
  \int_{U^{2}}p_{t-s}(x,y)\,G_{2\alpha}(y,y')\,p_{t-s}(x,y')\,d\mu(y)\,d\mu(y')
  =
  \|(-\Delta)^{-\alpha}p_{t-s}(x,\cdot)\|^{2}_{L^{2}(U,\mu)}.
\end{equation}
Multiplying by the prefactor $e^{2\lambda_{1}(t-s_{k+1})}$ and taking the supremum over
$x\in U$, definition~\eqref{eq:PsiHat-geom} gives
\begin{equation}\label{eq:Psi-id}
  e^{2\lambda_{1}(t-s_{k+1})}\sup_{x\in U}
  \|(-\Delta_{D,U})^{-\alpha}p^{D,U}_{t-s_{k+1}}(x,\cdot)\|^{2}_{L^{2}(U,\mu)}
  =
  \Psi(t-s_{k+1}).
\end{equation}
Plugging~\eqref{eq:bilinear-p} and~\eqref{eq:Psi-id} into~\eqref{i3},
and using $M_{k}(s_{k+1},y_{k+1})\le\sup_{y\in U}M_{k}(s_{k+1},y)$, gives
\begin{align*}
  e^{2\lambda_{1}t}M_{k+1}(t,x)
  &\le
  \int_{0}^{t}
    e^{2\lambda_{1}(t-s_{k+1})}
    \|(-\Delta)^{-\alpha}p_{t-s_{k+1}}(x,\cdot)\|^{2}_{L^{2}(U,\mu)}
    \, e^{2\lambda_{1}s_{k+1}}\sup_{y\in U}M_{k}(s_{k+1},y)
  \,ds_{k+1}
  \\
  &\le
  \int_{0}^{t}
    \Psi(t-s_{k+1})
    \, e^{2\lambda_{1}s_{k+1}}\sup_{y\in U}M_{k}(s_{k+1},y)
  \,ds_{k+1},
\end{align*}
which is~\eqref{eq:Mk-recursion}.

\smallskip

\noindent
{\it Step 3: Renewal argument.}
Set $N_{k}(t):=e^{2\lambda_{1}t}\sup_{x}M_{k}(t,x)$, so $N_{0}\equiv 1$ and
$N_{k+1}(t)\le\int_{0}^{t}\Psi(t-s)N_{k}(s)\,ds$ by~\eqref{eq:Mk-recursion}.
For the additional parameter $\rho>0$, define
$\mathcal{N}_{k}(\rho):=\sup_{t\ge 0}e^{-\rho t}N_{k}(t)$.
Then
\begin{equation}\label{i4}
  e^{-\rho t}N_{k+1}(t)
  \;\le\;
  \int_{0}^{t}e^{-\rho(t-s)}\Psi(t-s)\cdot e^{-\rho s}N_{k}(s)\,ds
  \;\le\;
  \hat{\Psi}(\rho)\,\mathcal{N}_{k}(\rho).
\end{equation}
Since $\mathcal{N}_{0}(\rho)=1$, a simple induction procedure gives $\mathcal{N}_{k}(\rho)\le\hat{\Psi}(\rho)^{k}$. Hence
\begin{equation}\label{eq:Nk-bound-geom}
  N_{k}(t) \;\le\; \hat{\Psi}(\rho)^{k}\,e^{\rho t}.
\end{equation}
We can now assemble the three ingredients.
From~\eqref{eq:fk-Mk},
\begin{equation*}
  \|f_{k}(\,\cdot\,,t,x)\|^{2}_{(\mathcal{H}^{D}_{\alpha})^{\otimes k}}
  \;\le\;
  C\,\frac{\beta^{2k}}{k!}\,\|u_{0}\|_{\infty}^{2}\,M_{k}(t,x).
\end{equation*}
The definition of $N_{k}$ gives $M_{k}(t,x)\le\sup_{y\in U}M_{k}(t,y)=e^{-2\lambda_{1}t}N_{k}(t)$,
so substituting and then applying~\eqref{eq:Nk-bound-geom} yields
\begin{equation*}
  \|f_{k}(\,\cdot\,,t,x)\|^{2}_{(\mathcal{H}^{D}_{\alpha})^{\otimes k}}
  \;\le\;
  C\,\frac{\beta^{2k}}{k!}\,\|u_{0}\|_{\infty}^{2}\,e^{-2\lambda_{1}t}\,N_{k}(t)
  \;\le\;
  C\,\frac{\beta^{2k}}{k!}\,\|u_{0}\|_{\infty}^{2}\,\hat{\Psi}(\rho)^{k}\,e^{(\rho-2\lambda_{1})t},
\end{equation*}
which is exactly~\eqref{eq:fk-bound-small}.
\end{proof}

\begin{remark}[On the role of the renewal argument]\label{rem:renewal-vs-fourier}
In this remark we are trying to answer the following question: 
why Step~2--3 above needed a renewal-type recursion,
while the analogous estimate~\eqref{eq:chaos-L2-converge} on $\R^{d}$
(Theorem~\ref{thm:chaos-eu}) did not. The comparison shows that the renewal
argument is not an artifact of the proof. Instead, it should be seen as a substitute for a tool
that is simply unavailable outside $\R^{d}$.

In the Euclidean setting, translation invariance lets us diagonalize the
whole chain of $n$ heat kernels at once. The Fourier transform~\eqref{eq:Kn-FT}
turns the $n$-fold convolution $K_{n}$ into the explicit product
$\prod_{i=1}^{n}e^{-\frac12(s_{i+1}-s_{i})|\xi_{1}+\cdots+\xi_{i}|^{2}}$, and
the decoupling inequality~\eqref{eq:spectral-decouple} can then be applied
successively, for $i=n,n-1,\ldots,1$, to separate all $n$ spectral variables
$\xi_{1},\ldots,\xi_{n}$ in a single pass, see~\eqref{eq:fn-gap-bound}. Once
decoupled, the remaining $n$-dimensional integral over the simplex $S_{t,n}$
is evaluated directly, for each fixed $n$, by the combinatorial argument of
Lemma~\ref{lem:simplex-integral}. Then the resulting bounds are summed over
$n$ in a single step~\eqref{eq:chaos-series-sum}. At no point is
$n!\,\|f_{n}\|^{2}_{\mathcal{H}^{\otimes n}}$ related back to
$(n-1)!\,\|f_{n-1}\|^{2}_{\mathcal{H}^{\otimes(n-1)}}$. Each term of the
series is controlled on its own.

On the metric measure space $(X,d_{U},\mu)$ this shortcut is unavailable,
because $U$ carries no group structure and $p^{D,U}_{t}$ is not translation
invariant. The factorization~\eqref{eq:Kn-FT} relies entirely on the
character identity $e^{i\xi\cdot(a+b)}=e^{i\xi\cdot a}\,e^{i\xi\cdot b}$,
which has no counterpart on a general $U$. The only substitute we have,
namely the bilinear identity~\eqref{eq:RK-D-bilinear}, diagonalizes a
\emph{single} pairing
$\int_{U^{2}}p_{t-s}(x,y)\,G_{2\alpha}(y,y')\,p_{t-s}(x,y')\,d\mu(y)\,d\mu(y')$
against the covariance kernel $G_{2\alpha}$, and gives no information on how
to integrate out several nested spatial variables $y_{1},\ldots,y_{k}$ coming
from a chain of $k$ heat kernels. Consequently, at each stage
the only variable that can be integrated out is the last one, $y_{k+1}$,
paired with the outer kernel $p_{t-s_{k+1}}$. This is exactly what produces
the recursion~\eqref{eq:Mk-recursion} relating $M_{k+1}$ to
$\sup_{y}M_{k}$, in place of a single joint bound on all $k+1$ variables at
once.

Given this recursion, the renewal step (Step~3 above) is the standard way to sum a convolution-type inequality
$N_{k+1}\le\Psi*N_{k}$ in the absence of a closed form for $N_{k}$.
Introducing the exponential weight $e^{-\rho t}$ turns the convolution into
the product bound $\mathcal{N}_{k+1}(\rho)\le\hat{\Psi}(\rho)\,\mathcal{N}_{k}(\rho)$~\eqref{i4},
which iterates directly into the geometric bound~\eqref{eq:Nk-bound-geom}.
In this sense the Laplace-transform trick used here plays exactly the role
that the spectral decoupling~\eqref{eq:spectral-decouple} (together with
Lemma~\ref{lem:simplex-integral}) plays in the Euclidean proof.\end{remark}

The bound~\eqref{eq:fk-bound-small} depends on $\hat{\Psi}(\rho)$ through the quantity defined
in~\eqref{eq:PsiHat-geom}.  Its size (and hence the range of $\beta$ for which the chaos
series converges) is governed by the position of $\alpha$ relative to the critical threshold
$d_{h}/(2d_{w})$.  The next lemma makes this precise.

\begin{lemma}[Three-regime estimate for $\hat{\Psi}$]\label{lem:PsiHat-bounds}
Under the hypotheses of Lemma~\ref{lem:fk-bound-geom}, the function $\hat{\Psi}(\rho)$
defined in~\eqref{eq:PsiHat-geom} satisfies the following three-regime bound, with
a possibly different constant $C$ in each case:
\begin{enumerate}
\item If $0\le\alpha<\dfrac{d_{h}}{2d_{w}}$, then
  \begin{equation}\label{eq:PsiHat-small}
    \hat{\Psi}(\rho)
    \;\le\;
    C\!\left(\frac{1}{\rho}+\frac{1}{\rho^{1+2\alpha-d_{h}/d_{w}}}\right).
  \end{equation}
\item If $\alpha=\dfrac{d_{h}}{2d_{w}}$, then
  \begin{equation}\label{eq:PsiHat-crit}
    \hat{\Psi}(\rho)
    \;\le\;
    C\,\frac{(\ln\max(3/2,\rho))^{2}}{\rho}.
  \end{equation}
\item If $\alpha>\dfrac{d_{h}}{2d_{w}}$, then
  \begin{equation}\label{eq:PsiHat-large}
    \hat{\Psi}(\rho)
    \;\le\;
    \frac{C}{\rho}.
  \end{equation}
\end{enumerate}
In all three cases, $\hat{\Psi}(\rho)\to 0$ as $\rho\to+\infty$.
\end{lemma}

\begin{proof}
From Lemma~\ref{lem:DeltaHK-D}, the function $\Psi$ defined by~\eqref{eq:PsiHat-geom}
satisfies
\begin{equation}\label{eq:Psi-bound-geom}
  \Psi(t)
  \;\le\;
  \begin{cases}
    C\,(1+t^{2\alpha-d_{h}/d_{w}}),
    & 0\le\alpha<\tfrac{d_{h}}{2d_{w}},\\[4pt]
    C\,|\ln\min(\tfrac{1}{2},t)|^{2},
    & \alpha=\tfrac{d_{h}}{2d_{w}},\\[4pt]
    C,
    & \alpha>\tfrac{d_{h}}{2d_{w}}.
  \end{cases}
\end{equation}
We integrate~\eqref{eq:Psi-bound-geom} against $e^{-\rho t}$ and treat each case separately.

\begin{enumerate}
\item {\it Subcritical case $0\le\alpha<d_{h}/(2d_{w})$.} Set $\gamma:=2\alpha-d_{h}/d_{w}<0$.
Integrating the two terms of~\eqref{eq:Psi-bound-geom} separately against $e^{-\rho t}$ gives
\begin{equation}\label{eq:crit-comp-sub1}
  \int_{0}^{\infty}e^{-\rho t}\,dt=\frac{1}{\rho},
\end{equation}
Since $\alpha\ge 0$ in the subcritical range, $\gamma=2\alpha-d_{h}/d_{w}\ge -d_{h}/d_{w}$, and
the recurrence condition $d_{h}<d_{w}$ of Assumption~\ref{hyp:mmspace} gives $d_{h}/d_{w}<1$, hence
$\gamma>-1$. Thus, by the Gamma function identity $\int_{0}^{\infty}e^{-\rho t}\,t^{\gamma}\,dt=\Gamma(1+\gamma)/\rho^{1+\gamma}$
(valid for $\gamma>-1$),
\begin{equation}\label{eq:crit-comp-sub2}
  \int_{0}^{\infty}e^{-\rho t}\,t^{\gamma}\,dt=\frac{\Gamma(1+\gamma)}{\rho^{1+\gamma}}.
\end{equation}
Adding~\eqref{eq:crit-comp-sub1} and~\eqref{eq:crit-comp-sub2}, multiplied by $C$, gives~\eqref{eq:PsiHat-small}.

\item {\it Critical case $\alpha=d_{h}/(2d_{w})$.} Here $\Psi(t)\le C|\ln\min(1/2,t)|^{2}$,
which is bounded by the constant $C|\ln(1/2)|^{2}$ for $t\ge 1/2$ and blows up like
$|\ln t|^{2}$ only as $t\to0$. Splitting the integral at $t=1/2$,
\begin{equation}\label{eq:crit-split}
  \hat{\Psi}(\rho)
  \;\le\;
  C\int_{0}^{1/2}e^{-\rho t}\,|\ln t|^{2}\,dt
  \;+\;
  C|\ln(1/2)|^{2}\int_{1/2}^{\infty}e^{-\rho t}\,dt,
\end{equation}
the second term of~\eqref{eq:crit-split} is $O(e^{-\rho/2}/\rho)$, hence negligible for large
$\rho$. In the first term, the change of variables $u=\rho t$ gives
\begin{equation}\label{eq:crit-subst}
  \int_{0}^{1/2}e^{-\rho t}\,|\ln t|^{2}\,dt
  \;=\;
  \frac{1}{\rho}\int_{0}^{\rho/2}e^{-u}\,(\ln\rho-\ln u)^{2}\,du.
\end{equation}
Expanding the square in~\eqref{eq:crit-subst}, the three resulting integrals
$\int_{0}^{\rho/2}e^{-u}\,du$, $\int_{0}^{\rho/2}e^{-u}\ln u\,du$, and
$\int_{0}^{\rho/2}e^{-u}(\ln u)^{2}\,du$ each converge, uniformly in $\rho$, to a finite
constant as $\rho\to\infty$. Hence~\eqref{eq:crit-subst} is of order $(\ln\rho)^{2}/\rho$, and
combining with the negligible term in~\eqref{eq:crit-split} yields~\eqref{eq:PsiHat-crit}. This
computation follows the same lines as Step~4 in the proof of \cite[Lem.~3.9]{baudoin:24:parabolic}.

\item {\it Supercritical case $\alpha>d_{h}/(2d_{w})$.} The bound $\Psi(t)\le C$ integrates
directly:
\begin{equation}\label{eq:crit-comp-super}
  \hat{\Psi}(\rho)\;\le\;C\int_{0}^{\infty}e^{-\rho t}\,dt\;=\;\frac{C}{\rho},
\end{equation}
which is~\eqref{eq:PsiHat-large}.
\end{enumerate}

\noindent
It is readily checked that in all three cases above, $\hat{\Psi}(\rho)$ tends to $0$ as $\rho\to+\infty$. This concludes our proof.
\end{proof}

We can now state the main result of this section.

\begin{theorem}[Existence, uniqueness, and $L^{2}$ moments]\label{thm:existence-geom}
Let $(X,d,\mu)$ satisfy Assumption~\ref{hyp:mmspace}, let $U\subset X$ be an
inner uniform domain (Definition~\ref{def:inner-uniform}), let $\alpha\ge 0$ and
$\beta\ge 0$, and let $u_{0}\in L^{\infty}(U,\mu)$.  Then equation~\eqref{eq:pam-D-geom}
has a unique solution in the sense of Definition~\ref{def:ito-sol-geom}, and
this solution satisfies the $L^{2}$ moment bound
\begin{equation}\label{eq:L2-asmp-geom}
  \limsup_{t\to+\infty}
  \frac{\ln\E(u(t,x)^{2})}{t}
  \;\le\;
  \Theta_{\alpha}(\beta) - 2\lambda_{1},
\end{equation}
where $\lambda_{1}=\lambda_{1}(U)$ is the first eigenvalue of $-\Delta_{D,U}$ and
$\Theta_{\alpha}:(0,\infty)\to(0,\infty)$ is a continuous increasing function whose asymptotics
are governed by $\hat{\Psi}$ via~\eqref{eq:PsiHat-small}--\eqref{eq:PsiHat-large}:
\begin{enumerate}
\item If $0\le\alpha<\dfrac{d_{h}}{2d_{w}}$, then $\Theta_{\alpha}(\beta)\sim C\beta^{2}$ as
  $\beta\to 0$ and $\Theta_{\alpha}(\beta)\sim C\beta^{2d_{w}/((1+2\alpha)d_{w}-d_{h})}$ as
  $\beta\to+\infty$.
\item If $\alpha=\dfrac{d_{h}}{2d_{w}}$, then $\Theta_{\alpha}(\beta)\sim C\beta^{2}$ as
  $\beta\to 0$ and $\Theta_{\alpha}(\beta)\sim C\beta^{2}(\ln\beta)^{2}$ as $\beta\to+\infty$.
\item If $\alpha>\dfrac{d_{h}}{2d_{w}}$, then $\Theta_{\alpha}(\beta)=C\beta^{2}$ for all $\beta>0$.
\end{enumerate}
\end{theorem}

\begin{proof}
By Proposition~\ref{prop:chaos-char-geom}, it suffices to verify~\eqref{eq:chaos-finite-geom}.
We treat case~(i) in detail, and indicate at the end of the proof the (minor) changes needed
for cases~(ii) and~(iii).

\smallskip

\noindent
{\it Step 1: choice of $\rho$.} Introduce the function
\begin{equation}\label{eq:F-def-geom}
  F(\rho) := \frac{1}{\rho}+\frac{1}{\rho^{1+\gamma}}, \qquad \rho>0, \qquad\text{where}\qquad
  \gamma:=2\alpha-\frac{d_{h}}{d_{w}},
\end{equation}
so that the bound~\eqref{eq:PsiHat-small} of Lemma~\ref{lem:PsiHat-bounds} reads
$\hat{\Psi}(\rho)\le C\,F(\rho)$. As noted in the proof of Lemma~\ref{lem:PsiHat-bounds}
(subcritical case), $\gamma\in(-1,0)$, so $F$ is continuous and strictly decreasing on
$(0,\infty)$, with $F(\rho)\to+\infty$ as $\rho\to0$ and $F(\rho)\to0$ as $\rho\to+\infty$.
Hence, for every $\beta>0$, there is a unique $\rho_{c}=\rho_{c}(\beta)>0$ such that
\begin{equation}\label{eq:rho-c-geom}
  C\beta^{2}F(\rho_{c}) = 1.
\end{equation}
Fix any $\rho>\rho_{c}(\beta)$. Since $F$ is decreasing, $F(\rho)<F(\rho_{c})$, and since
$\hat{\Psi}(\rho)\le F(\rho)$,
\begin{equation}\label{eq:delta-rho-geom}
  C\beta^{2}\hat{\Psi}(\rho) \;\le\; \delta_{\rho} \;<\; 1,
  \qquad\text{where}\qquad
  \delta_{\rho} := C\beta^{2}F(\rho).
\end{equation}

\smallskip

\noindent
{\it Step 2: $L^{2}$ bound.} Recall from the $L^{2}$-moment formula~\eqref{eq:L2-moment-geom}
of Proposition~\ref{prop:chaos-char-geom} that
\begin{equation}\label{eq:L2-decomp-geom}
  \E(u(t,x)^{2})
  \;=\;
  (P^{D,U}_{t}u_{0}(x))^{2}
  \;+\;
  \sum_{k=1}^{\infty}k!\,\|f_{k}(\,\cdot\,,t,x)\|^{2}_{(\mathcal{H}^{D}_{\alpha})^{\otimes k}}.
\end{equation}
Bound~\eqref{eq:fk-bound-small} of Lemma~\ref{lem:fk-bound-geom} controls each term of the
series in~\eqref{eq:L2-decomp-geom} through $\hat{\Psi}(\rho)$:
\begin{equation}\label{eq:fk-term-bound}
  k!\,\|f_{k}(\,\cdot\,,t,x)\|^{2}_{(\mathcal{H}^{D}_{\alpha})^{\otimes k}}
  \;\le\;
  C^{k}\,\beta^{2k}\,\|u_{0}\|_{\infty}^{2}\,\hat{\Psi}(\rho)^{k}\,
  e^{(\rho-2\lambda_{1})t}.
\end{equation}
By~\eqref{eq:delta-rho-geom}, $C\beta^{2}\hat{\Psi}(\rho)\le\delta_{\rho}$, so raising this to
the $k$-th power turns~\eqref{eq:fk-term-bound} into the $\delta_{\rho}$-bound
\begin{equation}\label{eq:fk-term-bound2}
  k!\,\|f_{k}(\,\cdot\,,t,x)\|^{2}_{(\mathcal{H}^{D}_{\alpha})^{\otimes k}}
  \;\le\;
  \|u_{0}\|_{\infty}^{2}\,\delta_{\rho}^{k}\,e^{(\rho-2\lambda_{1})t}.
\end{equation}
Since $\delta_{\rho}<1$ by~\eqref{eq:delta-rho-geom}, the series~\eqref{eq:chaos-finite-geom}
converges by comparison with~\eqref{eq:fk-term-bound2}, and Proposition~\ref{prop:chaos-char-geom}
gives existence and uniqueness of the solution of~\eqref{eq:pam-D-geom}.

We now proceed to the $L^{2}$ bound. Plugging~\eqref{eq:fk-term-bound2} into the series
in~\eqref{eq:L2-decomp-geom}, summing the resulting geometric series
$\sum_{k=1}^{\infty}\delta_{\rho}^{k}=\delta_{\rho}/(1-\delta_{\rho})$, and bounding the first
term of~\eqref{eq:L2-decomp-geom} by Corollary~\ref{cor:Linf-bound-D}
($(P^{D,U}_{t}u_{0}(x))^{2}\le C\|u_{0}\|_{\infty}^{2}e^{-2\lambda_{1}t}$), we get
\begin{equation}\label{eq:L2-est-geom}
  \E(u(t,x)^{2})
  \;\le\;
  C\|u_{0}\|_{\infty}^{2}e^{-2\lambda_{1}t}
  \;+\;
  \|u_{0}\|_{\infty}^{2}\,e^{(\rho-2\lambda_{1})t}
  \sum_{k=1}^{\infty}\delta_{\rho}^{k}
  \;\le\;
  C_{\rho}\,\|u_{0}\|_{\infty}^{2}\,e^{(\rho-2\lambda_{1})t},
\end{equation}
for a constant $C_{\rho}$ depending on $\rho$ (and $\beta,\alpha,d_{h},d_{w},U$) but not on
$t,x$. Consequently
\begin{equation}\label{eq:limsup-rho-geom}
  \limsup_{t\to+\infty}\frac{\ln\E(u(t,x)^{2})}{t} \;\le\; \rho-2\lambda_{1},
  \qquad\text{for every } \rho>\rho_{c}(\beta).
\end{equation}
Letting $\rho\downarrow\rho_{c}(\beta)$ in~\eqref{eq:limsup-rho-geom} yields~\eqref{eq:L2-asmp-geom}
with $\Theta_{\alpha}(\beta):=\rho_{c}(\beta)$.

\smallskip

\noindent
{\it Step 3: asymptotics of $\Theta_{\alpha}$.} From the definition~\eqref{eq:F-def-geom} of
$F$ and $\gamma\in(-1,0)$, it is readily checked that
\begin{equation}\label{eq:F-limits}
  \rho\,F(\rho)=1+\rho^{-\gamma}\xrightarrow[\rho\to0]{}1,
  \quad\text{and}\quad
  \rho^{1+\gamma}F(\rho)=\rho^{\gamma}+1\xrightarrow[\rho\to+\infty]{}1,
\end{equation}
By~\eqref{eq:rho-c-geom},
$F(\rho_{c}(\beta))=1/(C\beta^{2})$. As $\beta\to0$, $1/(C\beta^{2})\to+\infty$, so
$\rho_{c}(\beta)\to0$, and the first limit in~\eqref{eq:F-limits} gives
$\rho_{c}(\beta)\sim1/F(\rho_{c}(\beta))=C\beta^{2}$, i.e., $\Theta_{\alpha}(\beta)\sim C\beta^{2}$
as $\beta\to0$. As $\beta\to+\infty$, $1/(C\beta^{2})\to0$, so $\rho_{c}(\beta)\to+\infty$, and
the second limit in~\eqref{eq:F-limits} gives $\rho_{c}(\beta)^{1+\gamma}\sim1/F(\rho_{c}(\beta))=C\beta^{2}$,
i.e.,
\begin{equation}\label{eq:theta-large-beta}
  \Theta_{\alpha}(\beta) \;\sim\; (C\beta^{2})^{1/(1+\gamma)}
  \;=\; C'\beta^{2d_{w}/((1+2\alpha)d_{w}-d_{h})}, \qquad \beta\to+\infty,
\end{equation}
since $1+\gamma=((1+2\alpha)d_{w}-d_{h})/d_{w}$. This is the asymptotics claimed in case~(i).

\smallskip

\noindent
{\it Step 4: cases~(ii) and~(iii).} The critical case is treated by the same three steps,
with $F$ in~\eqref{eq:F-def-geom} replaced by $F_{2}(\rho):=(\ln\max(3/2,\rho))^{2}/\rho$ from
the bound~\eqref{eq:PsiHat-crit}. As above, $F_{2}$ is continuous and strictly decreasing, with
$F_{2}(\rho)\to+\infty$ as $\rho\to0$ and $F_{2}(\rho)\to0$ as $\rho\to+\infty$, and
\begin{equation*}
  \rho\,F_{2}(\rho)\xrightarrow[\rho\to0]{}(\ln(3/2))^{2},
  \qquad
  \frac{\rho\,F_{2}(\rho)}{(\ln\rho)^{2}}\xrightarrow[\rho\to+\infty]{}1.
\end{equation*}
Repeating Step~1--3 with $F_{2}$ in place of $F$, the definition~\eqref{eq:rho-c-geom} of
$\rho_{c}(\beta)$ becomes
\begin{equation*}
C\beta^{2}F_{2}(\rho_{c}(\beta))=1 .
\end{equation*}
 As $\beta\to0$, the first
limit above gives, exactly as in Step~3,
\begin{equation}\label{eq:theta-crit-small-beta}
  \Theta_{\alpha}(\beta) \;=\; \rho_{c}(\beta) \;\sim\; \frac{1}{F_{2}(\rho_{c}(\beta))}
  \;=\; C\beta^{2}, \qquad \beta\to0,
\end{equation}
which is the first asymptotic
in case~(ii). As $\beta\to+\infty$, the second limit above gives
\begin{equation}\label{eq:theta-crit-large-beta-1}
  \rho_{c}(\beta) \;\sim\; C\beta^{2}\,(\ln\rho_{c}(\beta))^{2}, \qquad \beta\to+\infty.
\end{equation}
Taking logarithms in~\eqref{eq:theta-crit-large-beta-1} gives $\ln\rho_{c}(\beta)\sim2\ln\beta$
as $\beta\to+\infty$. Hence substituting this back into~\eqref{eq:theta-crit-large-beta-1} yields
\begin{equation}\label{eq:theta-crit-large-beta-2}
  \Theta_{\alpha}(\beta) \;=\; \rho_{c}(\beta) \;\sim\; C\beta^{2}(\ln\beta)^{2},
  \qquad \beta\to+\infty,
\end{equation}
which is the second asymptotic in case~(ii).

\smallskip

\noindent
In the supercritical case~(iii), the bound~\eqref{eq:PsiHat-large} reads
$\hat{\Psi}(\rho)\le C/\rho$, so the definition~\eqref{eq:rho-c-geom} of $\rho_{c}$ becomes
\begin{equation}\label{eq:rho-c-super}
  \frac{C\beta^{2}}{\rho_{c}(\beta)} \;=\; 1,
\end{equation}
which gives, exactly rather than only asymptotically,
\begin{equation}\label{eq:theta-super}
  \Theta_{\alpha}(\beta) \;=\; \rho_{c}(\beta) \;=\; C\beta^{2}, \qquad \text{for every } \beta>0.
\end{equation}
Steps~2--3 above apply verbatim with this value of $\rho_{c}(\beta)$, and~\eqref{eq:theta-super}
is exactly case~(iii). This finishes the proof.
\end{proof}

\begin{remark}[Neumann boundary conditions]\label{rem:existence-neumann}
The same strategy applies to the Neumann PAM, replacing $\Delta_{D,U}$ by
$\Delta_{N,U}$, $p^{D,U}_{t}$ by $p^{N,U}_{t}$, $W^{D}_{\alpha}$ by
$W^{N}_{\alpha}$, and $G^{D,U}_{2\alpha}$ by $G^{N,U}_{2\alpha}+C_{U}$ (see
Definition~\ref{def:geom-noise-N} and Proposition~\ref{prop:RK-N-bounds}).  The
only difference is that $G^{N,U}_{2\alpha}$ is not pointwise non-negative. The
shifted kernel $G^{N,U}_{2\alpha}+C_{U}$ is non-negative, and the recursion
in~\eqref{eq:Mk-recursion} goes through with the Neumann heat kernel. The proof
of existence and uniqueness, as well as the $L^{2}$ bound, hold verbatim. We
refer to~\cite[Sec.~4.2]{baudoin:24:parabolic} for the details.
\end{remark}

\begin{remark}[Feynman--Kac formula for $\alpha>d_{h}/(2d_{w})$]\label{rem:fk-geom}
When $\alpha>d_{h}/(2d_{w})$, Theorem~\ref{thm:holder-field} shows that the density
field $x\mapsto\tilde{W}^{D}_{\alpha}(t,x)$ is H\"older continuous.  In this regime
the solution admits the Feynman--Kac representation
\begin{equation}\label{eq:fk-geom-formula}
  u(t,x)
  =
  \E^{B}\!\left(u_{0}(B^{x}_{t})
  \exp\!\left(\beta\int_{0}^{t}\dot{W}^{D}_{\alpha}(s,B^{x}_{t-s})\,ds
  -\frac{\beta^{2}}{2}\int_{0}^{t}G^{D,U}_{2\alpha}(B^{x}_{s},B^{x}_{s})\,ds\right)
  \mathbf{1}_{(t<T)}\right),
\end{equation}
where $B^{x}$ is a Brownian motion on $(X,d,\mu)$ started at $x$ and $T$ is its
hitting time of $U^{c}$.  The proof follows an approximation procedure
(mollify $W^{D}_{\alpha}$ in space, write the regularized Feynman--Kac formula,
pass to the limit). See~\cite[Section~3.2]{baudoin:24:parabolic}, as well
as~\cite[Section~3.2]{hu:15:stochastic} for the approximation argument.
\end{remark}

\section{Moment estimates}\label{sec:geom-moments}

In the Euclidean setting, Section~\ref{sec:moment-estimates} took the Feynman-Kac moment formula of Theorem~\ref{thm:fk-moments} as its starting point for the Lyapunov exponent estimates of the parabolic Anderson model. We follow the same strategy on the bounded domain $U\subset X$. In this section we develop moment estimates for the solution $u$ of~\eqref{eq:pam-D-geom} furnished by Theorem~\ref{thm:existence-geom}, under the standing framework of Assumption~\ref{hyp:mmspace} and Definition~\ref{def:inner-uniform}. Remark~\ref{rem:fk-geom} already recorded a heuristic Feynman-Kac formula for a single realization of $u(t,x)$, valid in the pointwise regime $\alpha>d_{h}/(2d_{w})$ and obtained by a space-time mollification of the noise. We begin by making this formula precise, extending it to a genuine \emph{moment} formula for every integer $k\ge2$ and every $\alpha>0$, via the same spatial mollification and It\^{o} calculus approach used in Section~\ref{sec:fk-moments}.

\subsection{Feynman-Kac representation of moments}\label{ssec:geom-fk-moments}

Recall from Theorem~\ref{thm:fk-moments} that on $\R^{d}$ the $k$-th moment of the solution is expressed as a Brownian expectation
\begin{equation*}
\E\bigl[u(t,x)^{k}\bigr]
=
\E_{B}\!\left[
\prod_{i=1}^{k}u_{0}(B^{i}_{t}+x)\,
\exp\left(
\beta^{2}\!\!\sum_{1\le i<j\le k}
\int_{0}^{t}\Lambda(B^{i}_{s}-B^{j}_{s})\,ds
\right)
\right],
\end{equation*}
obtained by a spatial mollification $u^{\varepsilon}$ of $u$ built from $k$ independent Brownian paths, and by letting $\varepsilon\to0$ (Section~\ref{sec:spatial-reg}). We now adapt this strategy to the geometric setting. To begin with, let us highlight three structural differences with the Euclidean case:

\begin{enumerate}
\item \emph{Bounded domain and exit times.} The noise $W^{D}_{\alpha}$ and the heat kernel $p^{D,U}_{t}$ are only defined on the bounded domain $U$ (Section~\ref{sec:geom-noise}), while the underlying Brownian motion $B^{x}$ lives on the full space $X$ (Definition~\ref{def:BM}). Each Brownian path entering the formula must therefore be paired with its exit time from $U$, exactly as in Remark~\ref{rem:fk-geom}.

\item \emph{Substitute for Dalang's condition.} In the Euclidean case, finiteness of $\|A^{\varepsilon}_{t,x}\|^{2}_{\mathcal{H}}$ relied on Dalang's condition~\eqref{eq:dalang-cond}. Here this role is played by the $L^{2}$-estimate of Lemma~\ref{lem:DeltaHK-D}, which controls $(-\Delta_{D,U})^{-\alpha}p^{D,U}_{t}(x,\cdot)$ for every $\alpha\ge0$.

\item \emph{No Fourier analysis.} The translation-invariance and Fourier computations of Section~\ref{sec:spatial-reg} (e.g.\ the change of variables leading to~\eqref{e3}) are unavailable on $U$. As for the existence and uniqueness results, they are replaced throughout by the Markov property of $B^{x}$, the Chapman-Kolmogorov identity for $p^{D,U}_{t}$, and the bilinear Riesz-kernel identity~\eqref{eq:RK-D-bilinear}.
\end{enumerate}

We record a notation for the $k$ independent Brownian paths needed for the $k$-th moment, in the spirit of Notation~\ref{not:BM-copies}.

\begin{notation}[Brownian copies and exit times]\label{not:BM-copies-geom}
Let $k\ge1$. We let $B^{x,1},\ldots,B^{x,k}$ be $k$ independent copies of the Brownian motion on $(X,d,\mu)$ of Definition~\ref{def:BM}, all started at the same point $x\in U$, and independent of the noise $W^{D}_{\alpha}$. For $1\le j\le k$ we write
\begin{equation*}
T^{j} := \inf\{t\ge0 : B^{x,j}_{t}\notin U\}
\end{equation*}
for the exit time of $B^{x,j}$ from $U$. As in Notation~\ref{not:BM-copies}, $\E_{B}$ denotes expectation with respect to the Brownian paths only (holding $W^{D}_{\alpha}$ fixed), and $\E_{W}$ denotes expectation with respect to $W^{D}_{\alpha}$ only (holding the paths fixed).
\end{notation}

\subsubsection{The spatially regularized process $u^{\varepsilon}$}\label{sssec:geom-ueps}

For $\varepsilon>0$ we use the Dirichlet heat kernel $p^{D,U}_{\varepsilon}$ itself as spatial mollifier, in place of the Gaussian kernel $p_{\varepsilon}$ of~\eqref{eq:Aeps-def}. For a Brownian path $B=B^{x}$ started at $x\in U$ and $t>0$, define the $\mathcal{H}^{D}_{\alpha}$-valued functional
\begin{equation}\label{eq:Aeps-def-geom}
A^{\varepsilon}_{t,x}(s,y)
:=
p^{D,U}_{\varepsilon}(B_{t-s},y)\,\mathbf{1}_{[0,t]}(s)\,\mathbf{1}_{\{t<T\}},
\qquad (s,y)\in\R_{+}\times U,
\end{equation}
where $T$ is the exit time of $B$ from $U$ (Notation~\ref{not:BM-copies-geom}). On the event $\{t<T\}$ the whole path $(B_{r})_{0\le r\le t}$ stays in $U$, so $B_{t-s}\in U$ for every $s\in[0,t]$ and~\eqref{eq:Aeps-def-geom} is well defined. On the complementary event we set $A^{\varepsilon}_{t,x}:=0$.

To compute inner products of functionals of the form~\eqref{eq:Aeps-def-geom}, introduce the auxiliary kernel
\begin{eqnarray}\label{eq:Xi-eps-def}
\Xi_{\varepsilon}(y,y')
&:=&
\bigl\langle
p^{D,U}_{\varepsilon}(y,\cdot),\,p^{D,U}_{\varepsilon}(y',\cdot)
\bigr\rangle_{\mathcal{W}^{-\alpha}_{D}(U)} \\
&=&
\bigl\langle
(-\Delta_{D,U})^{-\alpha}p^{D,U}_{\varepsilon}(y,\cdot),\,
(-\Delta_{D,U})^{-\alpha}p^{D,U}_{\varepsilon}(y',\cdot)
\bigr\rangle_{L^{2}(U,\mu)},
\qquad y,y'\in U, \notag
\end{eqnarray}
where the second equality is the bilinear Riesz-kernel identity~\eqref{eq:RK-D-bilinear}. The kernel $\Xi_{\varepsilon}$ expresses those inner products as time integrals along the Brownian paths. Namely, let $B^{x,i}$ and $B^{x,j}$ be two of the independent copies of Notation~\ref{not:BM-copies-geom}, with exit times $T^{i}$ and $T^{j}$, and write $A^{\varepsilon,i}_{t,x}$ and $A^{\varepsilon,j}_{t,x}$ for the functionals~\eqref{eq:Aeps-def-geom} built from $B^{x,i}$ and $B^{x,j}$ respectively. Unfolding the definition~\eqref{eq:H-geom-D} of the inner product in $\mathcal{H}^{D}_{\alpha}$ and the definition~\eqref{eq:Xi-eps-def} of $\Xi_{\varepsilon}$, on the event $\{t<T^{i}\}\cap\{t<T^{j}\}$ we get
\begin{equation}\label{eq:Aeps-inner-geom}
\langle A^{\varepsilon,i}_{t,x},A^{\varepsilon,j}_{t,x}\rangle_{\mathcal{H}^{D}_{\alpha}}
=
\int_{0}^{t}\Xi_{\varepsilon}(B^{x,i}_{t-s},B^{x,j}_{t-s})\,ds
=
\int_{0}^{t}\Xi_{\varepsilon}(B^{x,i}_{r},B^{x,j}_{r})\,dr,
\end{equation}
where the change of variable $r=t-s$ was used for the second equality (off the event $\{t<T^{i}\}\cap\{t<T^{j}\}$ at least one of the two functionals vanishes, and so does the inner product). Relation~\eqref{eq:Aeps-inner-geom} is the geometric analogue of the Fourier computation~\eqref{eq:Aeps-inner}. In particular, computing both functionals along the same path $B=B^{x}$ with exit time $T$ (that is taking $i=j$ in~\eqref{eq:Aeps-inner-geom}), on $\{t<T\}$ we get
\begin{equation}\label{eq:ueps-norm-geom}
\|A^{\varepsilon}_{t,x}\|^{2}_{\mathcal{H}^{D}_{\alpha}}
=
\int_{0}^{t}\Xi_{\varepsilon}(B_{t-r},B_{t-r})\,dr
=
\int_{0}^{t}
\bigl\|(-\Delta_{D,U})^{-\alpha}p^{D,U}_{\varepsilon}(B_{t-r},\cdot)\bigr\|^{2}_{L^{2}(U,\mu)}
\,dr.
\end{equation}
By Lemma~\ref{lem:DeltaHK-D}, applied at the fixed time $\varepsilon$ and valid for every $\alpha\ge0$, the integrand in~\eqref{eq:ueps-norm-geom} is bounded uniformly over $B_{t-r}\in U$ by a finite constant $\kappa(\varepsilon)^{2}$ (with $\kappa(\varepsilon)$ given by the right-hand side of~\eqref{eq:DeltaHK-D} at time $\varepsilon$), so
\begin{equation}\label{j1}
\|A^{\varepsilon}_{t,x}\|^{2}_{\mathcal{H}^{D}_{\alpha}} \le t\,\kappa(\varepsilon)^{2} < \infty
\qquad \text{on } \{t<T\}.
\end{equation}
In particular $A^{\varepsilon}_{t,x}\in\mathcal{H}^{D}_{\alpha}$ a.s.\ for \emph{every} $\alpha\ge0$, including the white noise case $\alpha=0$. Unlike the unmollified formula~\eqref{eq:fk-geom-formula} of Remark~\ref{rem:fk-geom}, which requires $\alpha>d_{h}/(2d_{w})$ for the exponent to make sense pointwise, the spatial regularization tames the noise for every $\alpha\ge0$, exactly as in the Euclidean case (compare~\eqref{eq:ueps-moment-bound}). Only the limit $\varepsilon\to0$ taken in Section~\ref{sssec:geom-fk-formula} below will force $\alpha>0$.

We may now define the regularized process itself.

\begin{definition}[Spatially regularized process]\label{def:ueps-geom}
For $\varepsilon>0$, $\beta\ge0$ and $u_{0}\in L^{\infty}(U,\mu)$, define
\begin{equation}\label{eq:FK-eps-geom}
u^{\varepsilon}(t,x)
:=
\E_{B}\!\left[
u_{0}(B_{t})\,\mathbf{1}_{\{t<T\}}\,
\exp\left(
\beta\,W^{D}_{\alpha}(A^{\varepsilon}_{t,x})
-\frac{\beta^{2}}{2}\|A^{\varepsilon}_{t,x}\|^{2}_{\mathcal{H}^{D}_{\alpha}}
\right)
\right],
\qquad (t,x)\in\R_{+}\times U,
\end{equation}
where $B=B^{x}$ is the Brownian motion of Notation~\ref{not:BM-copies-geom} and $A^{\varepsilon}_{t,x}$ is as in~\eqref{eq:Aeps-def-geom}.
\end{definition}

\begin{remark}\label{rem:ueps-martingale-geom}
The Gaussian exponential in~\eqref{eq:FK-eps-geom} is well defined path-by-path on $\{t<T\}$, since $A^{\varepsilon}_{t,x}\in\mathcal{H}^{D}_{\alpha}$ there. The stochastic exponential built from $A^{\varepsilon}_{t,x}$ is moreover a \emph{true} (as opposed to merely local) martingale in $s\in[0,t]$: its quadratic variation is deterministically bounded by $\beta^{2}t\,\kappa(\varepsilon)^{2}$, by the computation above. Novikov's condition therefore holds trivially, and $\E_{W}[\exp(\beta W^{D}_{\alpha}(A^{\varepsilon}_{t,x})-\tfrac{\beta^{2}}{2}\|A^{\varepsilon}_{t,x}\|^{2}_{\mathcal{H}^{D}_{\alpha}})]=1$ for every path $B$ with $t<T$, exactly as in the Euclidean case (below~\eqref{eq:FK-eps}).
\end{remark}

\subsubsection{The mild equation for $u^{\varepsilon}$}\label{sssec:geom-mild-ueps}

As in Lemma~\ref{lem:fk-duality}, the Feynman-Kac exponential defining $u^{\varepsilon}$ can be expanded via It\^{o}'s formula into a Dol\'eans-Dade stochastic exponential, and averaging over $B$ identifies $u^{\varepsilon}$ as the solution of a mild stochastic heat equation driven by a spatially mollified noise. Unlike~\eqref{eq:Weps-def}, the mollification here uses the Dirichlet semigroup $P^{D,U}_{\varepsilon}$ rather than a spatial convolution.

\begin{proposition}[Mild equation for $u^{\varepsilon}$]\label{prop:mild-ueps-geom}
Let $u^{\varepsilon}$ be as in Definition~\ref{def:ueps-geom}. Define the spatially mollified noise $W^{D,\varepsilon}_{\alpha}$ by setting, for every adapted process $\phi\in\mathcal{H}^{D}_{\alpha}$,
\begin{equation}\label{eq:Weps-def-geom}
\int_{0}^{\infty}\!\int_{U}\phi(s,y)\,W^{D,\varepsilon}_{\alpha}(ds,dy)
:=
\int_{0}^{\infty}\!\int_{U}(P^{D,U}_{\varepsilon}\phi_{s})(y)\,W^{D}_{\alpha}(ds,dy),
\end{equation}
where $(P^{D,U}_{\varepsilon}\phi_{s})(y):=\int_{U}p^{D,U}_{\varepsilon}(y,z)\,\phi(s,z)\,d\mu(z)$. Then for every $(t,x)\in\R_{+}\times U$,
\begin{equation}\label{eq:mild-ueps-geom}
u^{\varepsilon}(t,x)
=
P^{D,U}_{t}u_{0}(x)
+
\beta\int_{0}^{t}\!\int_{U}p^{D,U}_{t-s}(x,y)\,u^{\varepsilon}(s,y)\,W^{D,\varepsilon}_{\alpha}(ds,dy).
\end{equation}
\end{proposition}

\begin{proof}
We mimic Steps 1-a to 1-c of the proof of Lemma~\ref{lem:fk-duality}, replacing the translation-invariance argument of that proof by the Markov property of $B$ at the fixed time $t-s$.

\smallskip

\noindent
{\it Step 1: Stochastic exponential.} Mimicking Step 1-a of the proof of Lemma~\ref{lem:fk-duality}, for a fixed path $B$ with $t<T$ we set
\begin{equation}\label{eq:NsB-geom}
N^{B}_{s}
:=
\beta\int_{0}^{s}\!\int_{U}p^{D,U}_{\varepsilon}(B_{t-r},y)\,W^{D}_{\alpha}(dr,dy),
\qquad s\in[0,t].
\end{equation}
Since $B$ stays in $U$ on $[0,t]$, the space-time function $(r,y)\mapsto p^{D,U}_{\varepsilon}(B_{t-r},y)\mathbf{1}_{[0,t]}(r)$ is, path by path, an element of $\mathcal{H}^{D}_{\alpha}$, by the same computation as in~\eqref{eq:ueps-norm-geom}, so $N^{B}$ is a $W^{D}_{\alpha}$-martingale. By the isometry~\eqref{eq:geom-W-isometry-D} and the definition~\eqref{eq:Xi-eps-def} of $\Xi_{\varepsilon}$, its quadratic variation is
\begin{equation}\label{eq:NB-qv-geom}
\langle N^{B}\rangle_{s}
=
\beta^{2}\int_{0}^{s}\Xi_{\varepsilon}(B_{t-r},B_{t-r})\,dr,
\qquad s\in[0,t].
\end{equation}
Unlike the Euclidean computation~\eqref{e4}, this quadratic variation is not a deterministic function of $s$ alone, since there is no translation invariance available to eliminate the dependence on $B$. It is nonetheless deterministically bounded by $\beta^{2}s\,\kappa(\varepsilon)^{2}$, by the same estimate as in~\eqref{j1}. At $s=t$, comparing~\eqref{eq:NsB-geom}--\eqref{eq:NB-qv-geom} with~\eqref{eq:Aeps-def-geom} and~\eqref{eq:ueps-norm-geom} gives $N^{B}_{t}=\beta\,W^{D}_{\alpha}(A^{\varepsilon}_{t,x})$ and $\langle N^{B}\rangle_{t}=\beta^{2}\|A^{\varepsilon}_{t,x}\|^{2}_{\mathcal{H}^{D}_{\alpha}}$, so the Dol\'eans-Dade exponential $M^{B}_{s}:=\exp(N^{B}_{s}-\tfrac{1}{2}\langle N^{B}\rangle_{s})$ satisfies, by the same It\^{o}-formula computation as in~\eqref{eq:MsB}--\eqref{eq:MtB-expand},
\begin{equation}\label{eq:MtB-expand-geom}
M^{B}_{t}
=
1+\beta\int_{0}^{t}\!\int_{U}M^{B}_{s}\,p^{D,U}_{\varepsilon}(B_{t-s},y)\,W^{D}_{\alpha}(ds,dy).
\end{equation}

\smallskip

\noindent
{\it Step 2: Averaging over $B$ and identifying the kernel.} Mimicking Step 1-b of the proof of Lemma~\ref{lem:fk-duality}, we now average over $B$. By Definition~\ref{def:ueps-geom}, $u^{\varepsilon}(t,x)=\E_{B}[u_{0}(B_{t})\mathbf{1}_{\{t<T\}}M^{B}_{t}]$. Substituting~\eqref{eq:MtB-expand-geom} and applying Proposition~\ref{prop:stoch-fubini} (in its version for $W^{D}_{\alpha}$, with $\nu$ the law of $B$) to exchange $\E_{B}$ with the Walsh integral gives
\begin{equation}\label{eq:ueps-prefk-geom}
u^{\varepsilon}(t,x)
=
P^{D,U}_{t}u_{0}(x)
+
\beta\int_{0}^{t}\!\int_{U}
Z^{\varepsilon,t,x}(s,y)
\,W^{D}_{\alpha}(ds,dy),
\end{equation}
where the first term uses $\E_{B}[u_{0}(B_{t})\mathbf{1}_{\{t<T\}}]=P^{D,U}_{t}u_{0}(x)$~\eqref{eq:PDU-prob}, and where the stochastic integrand $Z^{\varepsilon,t,x}$ is defined by
\begin{equation}\label{eq:Zeps-integrand-geom}
Z^{\varepsilon,t,x}(s,y)
:=
\E_{B}\bigl[u_{0}(B_{t})\,\mathbf{1}_{\{t<T\}}\,M^{B}_{s}\,p^{D,U}_{\varepsilon}(B_{t-s},y)\bigr].
\end{equation}
We now wish to prove the identity
\begin{equation}\label{eq:Zeps-def-geom}
Z^{\varepsilon,t,x}(s,y)
=
\int_{U}p^{D,U}_{t-s}(x,z)\,p^{D,U}_{\varepsilon}(z,y)\,u^{\varepsilon}(s,z)\,d\mu(z).
\end{equation}
Once~\eqref{eq:Zeps-def-geom} is established, it turns~\eqref{eq:ueps-prefk-geom} into~\eqref{eq:mild-ueps-geom} by~\eqref{eq:Weps-def-geom} and Proposition~\ref{prop:stoch-fubini}. Note that this step mimmicks Step~2 of  Lemma~\ref{lem:fk-duality}. This will conclude the proof.

\smallskip

\noindent
{\it Step 3: Proof of~\eqref{eq:Zeps-def-geom}.}
In order to prove~\eqref{eq:Zeps-def-geom} we mimic Step 1-c of the proof of Lemma~\ref{lem:fk-duality}. To this aim, fix $s\in(0,t)$, condition on $B_{t-s}$ on the event $\{t-s<T\}$, and set $b:=B_{t-s}\in U$.
Write $\tilde{B}_{v}:=B_{t-s+v}$ for $v\in[0,s]$. By the Markov property~\eqref{eq:BM-Markov} of Definition~\ref{def:BM}, conditionally on $B_{t-s}=b$ this is a Brownian motion starting at $b$, independent of $\mathcal{F}^{B}_{t-s}$. Then for $r\in[0,s]$ one can write
\begin{equation}\label{eq:e6-geom}
B_{t-r}=\tilde{B}_{s-r},
\qquad
B_{t}=\tilde{B}_{s},
\qquad
p^{D,U}_{\varepsilon}(B_{t-s},y)=p^{D,U}_{\varepsilon}(b,y).
\end{equation}
In addition (a feature with no counterpart in the whole space), the exit time has to be split at $t-s$. Namely, writing $\tilde{T}$ for the exit time of $\tilde{B}$ from $U$, we have
\begin{equation}\label{eq:exit-split-geom}
\mathbf{1}_{\{t<T\}}=\mathbf{1}_{\{t-s<T\}}\,\mathbf{1}_{\{s<\tilde{T}\}}.
\end{equation}
Substituting the first identity of~\eqref{eq:e6-geom} into the integrand of $N^{B}_{s}$ in~\eqref{eq:NsB-geom} and into the quadratic variation~\eqref{eq:NB-qv-geom}, the definition $M^{B}_{s}=\exp(N^{B}_{s}-\tfrac{1}{2}\langle N^{B}\rangle_{s})$ of Step~1 gives
\begin{eqnarray}\label{eq:MsB-conditioned-geom}
M^{B}_{s}\big|_{B_{t-s}=b}
&=&
\exp\left(
\beta\int_{0}^{s}\!\int_{U}p^{D,U}_{\varepsilon}(\tilde{B}_{s-r},y)\,W^{D}_{\alpha}(dr,dy)
-\frac{\beta^{2}}{2}\int_{0}^{s}\Xi_{\varepsilon}(\tilde{B}_{s-r},\tilde{B}_{s-r})\,dr
\right) \notag \\ 
&=&
M^{\tilde{B}}_{s},
\end{eqnarray}
where $M^{\tilde{B}}$ is the stochastic exponential of Step~1. Therefore comparing~\eqref{eq:MsB-conditioned-geom} with the Feynman-Kac formula~\eqref{eq:FK-eps-geom} applied at $(s,b)$, and taking the second identity of~\eqref{eq:e6-geom} and the splitting~\eqref{eq:exit-split-geom} into account, we obtain
\begin{equation}\label{eq:kern-id-geom}
\E_{\tilde{B}}\bigl[u_{0}(\tilde{B}_{s})\,\mathbf{1}_{\{s<\tilde{T}\}}\,M^{\tilde{B}}_{s}\bigr]
=
u^{\varepsilon}(s,b).
\end{equation}
Integrating over the law $\Prob_{x}(B_{t-s}\in db,\,t-s<T)=p^{D,U}_{t-s}(x,b)\,d\mu(b)$~\eqref{eq:PDU-prob} of $B_{t-s}$ on $\{t-s<T\}$ and using~\eqref{eq:kern-id-geom} together with the third identity of~\eqref{eq:e6-geom}, as well as the definition~\eqref{eq:Zeps-integrand-geom} of $Z^{\varepsilon,t,x}$, we obtain
\begin{equation*}
Z^{\varepsilon,t,x}(s,y)
=
\int_{U}p^{D,U}_{t-s}(x,b)\,p^{D,U}_{\varepsilon}(b,y)\,u^{\varepsilon}(s,b)\,d\mu(b),
\end{equation*}
which proves~\eqref{eq:Zeps-def-geom} and hence~\eqref{eq:mild-ueps-geom}.
\end{proof}

\begin{remark}[Chaos expansion of $u^{\varepsilon}$]\label{rem:ueps-chaos-geom}
Iterating the mild equation~\eqref{eq:mild-ueps-geom} exactly as in the derivation of the chaos expansion~\eqref{eq:u-chaos-geom}-\eqref{eq:fk-geom} shows that $u^{\varepsilon}$ has its own chaos expansion.
Specifically one can write
\begin{equation}\label{eq:ueps-chaos-geom}
u^{\varepsilon}(t,x)
=
P^{D,U}_{t}u_{0}(x)+\sum_{k\ge1}I_{k}(f^{\varepsilon}_{k}(\cdot,t,x)).
\end{equation}
Suppressing $s_{i},t,x$ in the kernels $f^{\varepsilon}_{k}$ above and writing the spatial arguments only for notational sake, it holds
\begin{equation}\label{eq:fkeps-def}
f^{\varepsilon}_{k}(y_{1},\ldots,y_{k})
:=
\int_{U^{k}}f_{k}(z_{1},\ldots,z_{k})\prod_{i=1}^{k}p^{D,U}_{\varepsilon}(y_{i},z_{i})\,d\mu(z_{i}) .
\end{equation}
Otherwise stated, we have
$f^{\varepsilon}_{k}=(P^{D,U}_{\varepsilon})^{\otimes k}f_{k}$, similarly to relation~\eqref{eq:fneps-def} in Remark~\ref{rem:ueps-chaos} for the Euclidean case.
\end{remark}

The expansions~\eqref{eq:u-chaos-geom} and~\eqref{eq:ueps-chaos-geom} live on the same chaos levels, so $u^{\varepsilon}(t,x)-u(t,x)$ can be compared kernel by kernel. Combined with Nelson's hypercontractivity inequality~\eqref{eq:hypercontractivity-sum}, this observation yields the following convergence result, which is the geometric counterpart of Proposition~\ref{prop:ueps-Lk-conv}.

\begin{proposition}[$L^{k}$ convergence of $u^{\varepsilon}$]\label{prop:ueps-Lk-conv-geom}
Let $\alpha\ge0$, $\beta\ge0$, $u_{0}\in L^{\infty}(U,\mu)$, and let $u$ be the solution of~\eqref{eq:pam-D-geom} with chaos kernels $f_{k}$ of~\eqref{eq:fk-geom}. Then for every integer $k\ge2$,
\begin{equation*}
u^{\varepsilon}(t,x)\;\longrightarrow\;u(t,x)
\qquad\text{in } L^{k}(\Omega),\ \text{as } \varepsilon\to0.
\end{equation*}
\end{proposition}

\begin{proof}
The proof follows the same pattern as that of its Euclidean counterpart, Proposition~\ref{prop:ueps-Lk-conv}, and we point out the correspondence between the two sets of inequalities as we go.

Since $P^{D,U}_{\varepsilon}$ is self-adjoint with spectrum $\{e^{-\lambda_{j}\varepsilon}\}_{j\ge1}\subset(0,1]$ (Notation~\ref{not:spectrum-D}), it is a contraction on $\mathcal{W}^{-\alpha}_{D}(U)$ (as well as on $L^{2}(U,\mu)$), and $P^{D,U}_{\varepsilon}\to\mathrm{Id}$ strongly as $\varepsilon\to0$ on both spaces. Consequently the kernels $f^{\varepsilon}_{n}=(P^{D,U}_{\varepsilon})^{\otimes n}f_{n}$ of~\eqref{eq:fkeps-def} satisfy $\|f^{\varepsilon}_{n}(\cdot,t,x)\|_{(\mathcal{H}^{D}_{\alpha})^{\otimes n}}\le\|f_{n}(\cdot,t,x)\|_{(\mathcal{H}^{D}_{\alpha})^{\otimes n}}$ for every $\varepsilon>0$ and every $n\ge1$, and $f^{\varepsilon}_{n}\to f_{n}$ in $(\mathcal{H}^{D}_{\alpha})^{\otimes n}$ as $\varepsilon\to0$, for each fixed $n$. These two properties are the respective analogues of the contraction bound~\eqref{eq:fneps-contract} and of the fixed-level convergence~\eqref{eq:fneps-cv-fixed-n} in Step~1 of the proof of Proposition~\ref{prop:ueps-Lk-conv}, with the spectral bound $e^{-\lambda_{j}\varepsilon}\le1$ playing the role of the Fourier multiplier bound $e^{-\varepsilon|\xi|^{2}/2}\le1$ appearing in~\eqref{eq:fneps-fourier}.

Fix an integer $k\ge2$. By Lemma~\ref{lem:PsiHat-bounds}, $\hat{\Psi}(\rho)\to0$ as $\rho\to\infty$, so we may choose $\rho$ large enough that $\delta:=(k-1)\,C\beta^{2}\hat{\Psi}(\rho)<1$, where $C$ is the constant of Lemma~\ref{lem:fk-bound-geom}. (This choice of $\rho$ is the geometric substitute for the choice of the frequency threshold $N$ achieving $\theta_{N}=2(k-1)\beta^{2}C_{N}<1$ in~\eqref{eq:fn-k-weighted}. In both cases one tunes a free parameter, $\rho$ here and $N$ there, so that the small quantity it governs, namely $\hat{\Psi}(\rho)$ here and the high-frequency tail $C_{N}$ of~\eqref{eq:DN-CN} there, beats the factor $(k-1)\beta^{2}$ stemming from hypercontractivity at order $k$.) Bound~\eqref{eq:fk-bound-small} then gives
\begin{equation}\label{eq:hypercontract-term-geom}
(k-1)^{n/2}\sqrt{n!}\,\|f_{n}(\cdot,t,x)\|_{(\mathcal{H}^{D}_{\alpha})^{\otimes n}}
\;\le\;
\|u_{0}\|_{\infty}\,\delta^{n/2}\,e^{(\rho-2\lambda_{1})t/2},
\qquad n\ge1,
\end{equation}
which is summable in $n$ since $\delta<1$. Observe that~\eqref{eq:hypercontract-term-geom} is the exact analogue of the summability estimate~\eqref{eq:hypercontract-series} obtained in Step~2 of the proof of Proposition~\ref{prop:ueps-Lk-conv}, the geometric series in $\delta^{n/2}$ replacing the series produced there by the binomial bound~\eqref{eq:fn-k-weighted}. From here the conclusion is reached exactly as in Step~3 of the proof of Proposition~\ref{prop:ueps-Lk-conv}. Namely, since $\|f^{\varepsilon}_{n}-f_{n}\|\le2\|f_{n}\|$, twice the bound~\eqref{eq:hypercontract-term-geom} dominates $(k-1)^{n/2}\sqrt{n!}\|f^{\varepsilon}_{n}-f_{n}\|$ uniformly in $\varepsilon$, and each term of the latter converges to $0$ as $\varepsilon\to0$ (for fixed $n$). By dominated convergence for series,
\begin{equation*}
\sum_{n=1}^{\infty}(k-1)^{n/2}\sqrt{n!}\,\|f^{\varepsilon}_{n}(\cdot,t,x)-f_{n}(\cdot,t,x)\|_{(\mathcal{H}^{D}_{\alpha})^{\otimes n}}
\;\longrightarrow\;0
\qquad\text{as }\varepsilon\to0.
\end{equation*}
Nelson's hypercontractivity inequality~\eqref{eq:hypercontractivity-sum}, applied to $u^{\varepsilon}(t,x)-u(t,x)=\sum_{n\ge1}I_{n}(f^{\varepsilon}_{n}-f_{n})$ as in~\eqref{eq:ueps-u-Lk-bound}, bounds $\|u^{\varepsilon}(t,x)-u(t,x)\|_{L^{k}(\Omega)}$ by the left-hand side above, which therefore tends to $0$.
\end{proof}

\subsubsection{The Feynman-Kac moment formula}\label{sssec:geom-fk-formula}

We can now state the main result of this section.

\begin{theorem}[Feynman-Kac moment formula]\label{thm:fk-moments-geom}
Let $\alpha>0$, $\beta\ge0$ and $u_{0}\in L^{\infty}(U,\mu)$, and let $u$ be the solution of~\eqref{eq:pam-D-geom} given by Theorem~\ref{thm:existence-geom}. Then for every integer $k\ge2$, every $t\ge0$ and every $x\in U$,
\begin{equation}\label{eq:fk-moment-formula-geom}
\E\bigl[u(t,x)^{k}\bigr]
=
\E_{B}\!\left[
\prod_{j=1}^{k}u_{0}(B^{x,j}_{t})\,\mathbf{1}_{\{t<T^{j}\}}\,
\exp\left(
\beta^{2}\!\!\sum_{1\le i<j\le k}
\int_{0}^{t}G^{D,U}_{2\alpha}(B^{x,i}_{s},B^{x,j}_{s})\,ds
\right)
\right],
\end{equation}
where $B^{x,1},\ldots,B^{x,k}$ and $T^{1},\ldots,T^{k}$ are as in Notation~\ref{not:BM-copies-geom}.
\end{theorem}

\begin{proof}
We proceed in three steps, exactly as in the proof of Theorem~\ref{thm:fk-moments}.

\smallskip

\noindent
{\it Step 1: Moment formula for $u^{\varepsilon}$.} Let $B^{x,1},\ldots,B^{x,k}$ be as in Notation~\ref{not:BM-copies-geom}, and write $A^{\varepsilon,j}_{t,x}$ for the functional~\eqref{eq:Aeps-def-geom} built from $B^{x,j}$. Since the $B^{x,j}$ are independent copies independent of $W^{D}_{\alpha}$, we substitute~\eqref{eq:FK-eps-geom} for each factor. Then using Fubini, the isometry~\eqref{eq:geom-W-isometry-D}, and canceling the diagonal terms exactly as in Step~1 of the proof of Proposition~\ref{prop:ueps-moments}, we obtain
\begin{equation}\label{eq:mom-eps-formula-geom}
\E\bigl[(u^{\varepsilon}(t,x))^{k}\bigr]
=
\E_{B}\!\left[
\prod_{j=1}^{k}u_{0}(B^{x,j}_{t})\,\mathbf{1}_{\{t<T^{j}\}}\,
\exp\left(
\beta^{2}\!\!\sum_{1\le i<j\le k}
\langle A^{\varepsilon,i}_{t,x},A^{\varepsilon,j}_{t,x}\rangle_{\mathcal{H}^{D}_{\alpha}}
\right)
\right].
\end{equation}

\smallskip

\noindent
{\it Step 2: Limit of the cross terms.} By~\eqref{eq:Xi-eps-def} and the semigroup property of $P^{D,U}$, we have
\begin{equation}\label{eq:Xi-eps-semigroup}
\Xi_{\varepsilon}(y,y')
=
\frac{1}{\Gamma(2\alpha)}\int_{0}^{\infty}r^{2\alpha-1}p^{D,U}_{2\varepsilon+r}(y,y')\,dr
=
\frac{1}{\Gamma(2\alpha)}\int_{2\varepsilon}^{\infty}(s-2\varepsilon)^{2\alpha-1}p^{D,U}_{s}(y,y')\,ds,
\end{equation}
where the first equality follows from~\eqref{eq:RK-D} (with $\alpha$ replaced by $2\alpha$) together with the Chapman-Kolmogorov identity 
\begin{equation}\label{j2}
\int_{U}p^{D,U}_{\varepsilon}(y,z)p^{D,U}_{r}(z,z')\,d\mu(z)=p^{D,U}_{\varepsilon+r}(y,z') ,
\end{equation}
applied at both endpoints $y,y'$. Note that~\eqref{j2} is the geometric substitute for the Fourier identity~\eqref{eq:Aeps-inner} in the $\mathbb{R}^{d}$ case. 

Let us now take limits in~\eqref{eq:Xi-eps-semigroup}, for which we should elaborate slightly on~\eqref{eq:Aeps-inner-limit}. Namely
for $y\ne y'$ we claim that
\begin{equation}\label{eq:Xi-eps-limit}
\lim_{\varepsilon\to0}\Xi_{\varepsilon}(y,y')
=
G^{D,U}_{2\alpha}(y,y').
\end{equation}
In order to prove~\eqref{eq:Xi-eps-limit}, we apply dominated convergence to the middle expression of~\eqref{eq:Xi-eps-semigroup}. The pointwise convergence of the dominated argument is clear. That is, for every $r>0$ we have
\begin{equation*}
\lim_{\varepsilon\to0} p^{D,U}_{2\varepsilon+r}(y,y') = p^{D,U}_{r}(y,y'),
\end{equation*}
which stems easily from the spectral expansion~\eqref{eq:HK-spectral-D}. The domination can be argued as follows: for $y\ne y'$, the sub-Gaussian upper bound~\eqref{eq:subGauss-D-upper} for $t\le1$ and the spectral expansion~\eqref{eq:HK-spectral-D} for $t\ge1$ furnish a constant $c=c(y,y')$ such that $p^{D,U}_{t}(y,y')\le c\,e^{-\lambda_{1}t}$ for every $t>0$. Hence for $\varepsilon\in(0,1)$ the integrand $r^{2\alpha-1}p^{D,U}_{2\varepsilon+r}(y,y')$ in~\eqref{eq:Xi-eps-semigroup} is bounded by $c\,r^{2\alpha-1}e^{-\lambda_{1}r}$, which is integrable on $(0,+\infty)$. This finishes the proof of our claim~\eqref{eq:Xi-eps-limit}.

Next we upgrade the pointwise convergence~\eqref{eq:Xi-eps-limit} to a convergence result for the cross inner products appearing in~\eqref{eq:mom-eps-formula-geom}. Namely we claim that for $i\ne j$ we have
\begin{equation}\label{eq:Aeps-inner-limit-geom}
\lim_{\varepsilon\to0}
\E_{B}\left[
\mathbf{1}_{\{t<T^{i}\}}\,\mathbf{1}_{\{t<T^{j}\}}
\left|
\langle A^{\varepsilon,i}_{t,x},A^{\varepsilon,j}_{t,x}\rangle_{\mathcal{H}^{D}_{\alpha}}
-
\int_{0}^{t}G^{D,U}_{2\alpha}(B^{x,i}_{r},B^{x,j}_{r})\,dr
\right|
\right]
=0 .
\end{equation}
Relation~\eqref{eq:Aeps-inner-limit-geom} is the geometric analogue of the $L^{1}(\Prob_{B})$ convergence~\eqref{eq:Aeps-inner-limit}, and the indicators in~\eqref{eq:Aeps-inner-limit-geom} are harmless.  
Let us also highlight the main obstacle in the proof of~\eqref{eq:Aeps-inner-limit-geom}. In the Euclidean case, the convergence~\eqref{eq:Aeps-inner-limit} stemmed from an elementary dominated convergence argument, the domination being provided by the uniform bound $|e^{-\varepsilon|\xi|^{2}}e^{i\xi\cdot z}|\le1$ together with Dalang's condition~\eqref{eq:dalang-cond}. In our geometric setting no such uniform bound is available. We circumvent this problem thanks to Scheff\'e's lemma, spelled out under the following form: if a family of nonnegative functions $F_{\varepsilon}$ converges almost everywhere to $F$ and $\int F_{\varepsilon}\to\int F<\infty$, then $F_{\varepsilon}\to F$ in $L^{1}$.

In order to use the above lemma, let us call $\Omega_{B}$ the probability space supporting the Brownian copies of Notation~\ref{not:BM-copies-geom}. Then we apply Scheff\'e's lemma on the product space $\Omega_{B}\times[0,t]$ equipped with the measure $\Prob_{B}\otimes dr$, to the following functions:
\begin{eqnarray}\label{eq:Scheffe-Feps-def}
F_{\varepsilon}(\omega,r) &=& 
\mathbf{1}_{\{r<T^{i}\}}\,\mathbf{1}_{\{r<T^{j}\}}\,\Xi_{\varepsilon}\bigl(B^{x,i}_{r},B^{x,j}_{r}\bigr),
\notag\\
F(\omega,r) &=& 
\mathbf{1}_{\{r<T^{i}\}}\,\mathbf{1}_{\{r<T^{j}\}}\,G^{D,U}_{2\alpha}\bigl(B^{x,i}_{r},B^{x,j}_{r}\bigr).
\end{eqnarray}
Those functions are nonnegative, thanks to the positivity of the heat kernel $p^{D,U}$ in the representations~\eqref{eq:Xi-eps-semigroup} and~\eqref{eq:RK-D}. We now check the two remaining hypotheses of Scheff\'e's lemma.

\smallskip

\noindent
{\it (a) Almost everywhere convergence.} The Ahlfors regularity~\eqref{eq:ahlfors} implies that $\mu$ has no atoms. Since $B^{x,i}_{r}$ and $B^{x,j}_{r}$ are independent, each with a law admitting a density with respect to $\mu$, we get $\Prob(B^{x,i}_{r}=B^{x,j}_{r})=0$ for every fixed $r>0$. By Fubini's theorem, the set of couples $(\omega,r)$ such that $B^{x,i}_{r}(\omega)=B^{x,j}_{r}(\omega)$ is thus $\Prob_{B}\otimes dr$-negligible. Owing to the convergence~\eqref{eq:Xi-eps-limit}, valid off the diagonal, we conclude that $F_{\varepsilon}\to F$, $\Prob_{B}\otimes dr$-almost everywhere.

\smallskip

\noindent
{\it (b) Convergence of the integrals.} We first derive an alternative expression for $\Xi_{\varepsilon}$. Starting from the middle expression of~\eqref{eq:Xi-eps-semigroup}, apply the Chapman-Kolmogorov identity~\eqref{j2} at both endpoints as before, exchange the time and space integrals, and identify the time integral by means of~\eqref{eq:RK-D}. This yields
\begin{equation}\label{eq:Xi-eps-G-form}
\Xi_{\varepsilon}(z,z')
=
\int_{U}\int_{U}p^{D,U}_{\varepsilon}(z,w)\,p^{D,U}_{\varepsilon}(z',w')\,G^{D,U}_{2\alpha}(w,w')\,d\mu(w)\,d\mu(w'),
\qquad z,z'\in U.
\end{equation}
Next recall from~\eqref{eq:PDU-prob} that $\Prob(B^{x,i}_{r}\in dz,\,r<T^{i})=p^{D,U}_{r}(x,z)\,d\mu(z)$. Hence, invoking the independence of $B^{x,i}$ and $B^{x,j}$, then~\eqref{eq:Xi-eps-G-form}, the Chapman-Kolmogorov identity~\eqref{j2}, and reading~\eqref{eq:Xi-eps-G-form} backwards at scale $r+\varepsilon$ (note that the definition~\eqref{eq:Xi-eps-def} of $\Xi_{\varepsilon}$ makes sense for any positive value of the parameter $\varepsilon$, which gives a meaning to $\Xi_{r+\varepsilon}$ and $\Xi_{s}$ below), we obtain
\begin{eqnarray}\label{eq:Scheffe-mass-geom}
\int F_{\varepsilon}\,d(\Prob_{B}\otimes dr)
&=&
\int_{0}^{t}\int_{U}\int_{U}p^{D,U}_{r}(x,z)\,p^{D,U}_{r}(x,z')\,\Xi_{\varepsilon}(z,z')\,d\mu(z)\,d\mu(z')\,dr \notag\\
&=&
\int_{0}^{t}\int_{U}\int_{U}p^{D,U}_{r+\varepsilon}(x,w)\,p^{D,U}_{r+\varepsilon}(x,w')\,G^{D,U}_{2\alpha}(w,w')\,d\mu(w)\,d\mu(w')\,dr \notag\\
&=&
\int_{0}^{t}\Xi_{r+\varepsilon}(x,x)\,dr
\;=\;
\int_{\varepsilon}^{t+\varepsilon}\Xi_{s}(x,x)\,ds .
\end{eqnarray}
The same computation, performed without the smoothing at scale $\varepsilon$, gives
\begin{equation}\label{eq:Scheffe-mass-geom-limit}
\int F\,d(\Prob_{B}\otimes dr)
=
\int_{0}^{t}\Xi_{s}(x,x)\,ds .
\end{equation}
In addition, the function $s\mapsto\Xi_{s}(x,x)$ is integrable on $(0,t+1]$. Indeed, comparing the diagonal value $\Xi_{s}(x,x)$ given by the second line of~\eqref{eq:Xi-eps-def} with the definition~\eqref{eq:PsiHat-geom} of $\Psi$, we have $\Xi_{s}(x,x)\le e^{-2\lambda_{1}s}\Psi(s)\le\Psi(s)$, and thus
\begin{equation}\label{eq:Xi-diag-integrable}
\int_{0}^{t+1}\Xi_{s}(x,x)\,ds
\;\le\;
e^{t+1}\int_{0}^{t+1}e^{-s}\Psi(s)\,ds
\;\le\;
e^{t+1}\,\hat{\Psi}(1).
\end{equation}
Now the right hand side of \eqref{eq:Xi-diag-integrable} is finite by Lemma~\ref{lem:PsiHat-bounds} (this is one of the places where the recurrence condition $d_{h}<d_{w}$ of Assumption~\ref{hyp:mmspace} is used), which yields the integrability of the function $\Xi_{\cdot}(x,x)$ on $(0,t+1]$. Since the integrals~\eqref{eq:Scheffe-mass-geom} and~\eqref{eq:Scheffe-mass-geom-limit} differ by the integral of $\Xi_{s}(x,x)$ over the two small intervals $(0,\varepsilon]$ and $(t,t+\varepsilon]$, the integrability of $\Xi_{\cdot}(x,x)$ shows that $\int F_{\varepsilon}\,d(\Prob_{B}\otimes dr)\to\int F\,d(\Prob_{B}\otimes dr)<\infty$ as $\varepsilon\to0$.

\smallskip

Having checked items (a) and (b), Scheff\'e's lemma asserts that $F_{\varepsilon}\to F$ in $L^{1}(\Prob_{B}\otimes dr)$. In order to deduce~\eqref{eq:Aeps-inner-limit-geom}, observe that $\mathbf{1}_{\{t<T^{i}\}}\le\mathbf{1}_{\{r<T^{i}\}}$ for $r\le t$. Together with the expression~\eqref{eq:Aeps-inner-geom} for the cross inner products and the triangle inequality, this entails
\begin{equation*}
\E_{B}\left[
\mathbf{1}_{\{t<T^{i}\}}\,\mathbf{1}_{\{t<T^{j}\}}
\left|
\langle A^{\varepsilon,i}_{t,x},A^{\varepsilon,j}_{t,x}\rangle_{\mathcal{H}^{D}_{\alpha}}
-
\int_{0}^{t}G^{D,U}_{2\alpha}(B^{x,i}_{r},B^{x,j}_{r})\,dr
\right|
\right] 
\le
\int\bigl|F_{\varepsilon}-F\bigr|\,d(\Prob_{B}\otimes dr),
\end{equation*}
and the right-hand side vanishes as $\varepsilon\to0$. The claim~\eqref{eq:Aeps-inner-limit-geom} is thus proved.

\smallskip

\noindent
{\it Step 3: Convergence of moments.} This step is entirely similar to Step 3 in the proof of the Euclidean Theorem~\ref{thm:fk-moments}, and we mostly have to translate the notation. The only genuine difference lies in the justification of the uniform integrability, for which we invoke the moment formula~\eqref{eq:mom-eps-formula-geom} itself in place of the Euclidean uniform bound~\eqref{eq:ueps-moment-bound}.

We first take limits on the left-hand side of~\eqref{eq:mom-eps-formula-geom}. Proposition~\ref{prop:ueps-Lk-conv-geom} gives $u^{\varepsilon}(t,x)\to u(t,x)$ in $L^{k}(\Omega)$. Since $L^{k}$-convergence of $X_{n}\to X$ implies $\E[X_{n}^{k}]\to\E[X^{k}]$, we obtain
\begin{equation}\label{eq:lhs-moment-limit-geom}
\lim_{\varepsilon\to0}\E\bigl[(u^{\varepsilon}(t,x))^{k}\bigr]
=
\E\bigl[u(t,x)^{k}\bigr].
\end{equation}

We now handle the right-hand side of~\eqref{eq:mom-eps-formula-geom}. On the probability space $\Omega_{B}$ supporting the Brownian copies (introduced before~\eqref{eq:Scheffe-Feps-def}), define the random variables
\begin{eqnarray}\label{eq:Veps-def-geom}
V_{\varepsilon}
&:=&
\prod_{j=1}^{k}u_{0}(B^{x,j}_{t})\,\mathbf{1}_{\{t<T^{j}\}}\,
\exp\left(
\beta^{2}\sum_{1\le i<j\le k}
\langle A^{\varepsilon,i}_{t,x},A^{\varepsilon,j}_{t,x}\rangle_{\mathcal{H}^{D}_{\alpha}}
\right),
\notag\\
V
&:=&
\prod_{j=1}^{k}u_{0}(B^{x,j}_{t})\,\mathbf{1}_{\{t<T^{j}\}}\,
\exp\left(
\beta^{2}\sum_{1\le i<j\le k}
\int_{0}^{t}G^{D,U}_{2\alpha}(B^{x,i}_{s},B^{x,j}_{s})\,ds
\right),
\end{eqnarray}
so that the moment formula~\eqref{eq:mom-eps-formula-geom} reads $\E[(u^{\varepsilon}(t,x))^{k}]=\E_{B}[V_{\varepsilon}]$, while the desired identity~\eqref{eq:fk-moment-formula-geom} reads $\E[u(t,x)^{k}]=\E_{B}[V]$. Taking~\eqref{eq:lhs-moment-limit-geom} into account, we are thus left with proving that
\begin{equation}\label{eq:Veps-L1-limit-geom}
\lim_{\varepsilon\to0}\E_{B}[V_{\varepsilon}]
=
\E_{B}[V].
\end{equation}
We will get~\eqref{eq:Veps-L1-limit-geom} from Vitali's convergence theorem, which is ensured thanks to two ingredients: convergence in probability of $V_{\varepsilon}$ to $V$, and boundedness in $L^{p}(\mathbb{P}_{B})$ for $p>1$ of the family $\{V_{\varepsilon}\}$. Note that below we will choose $p=2$.

Let us first check that $V_{\varepsilon}\to V$ in $\Prob_{B}$-probability. Off the event $\cap_{j=1}^{k}\{t<T^{j}\}$, both $V_{\varepsilon}$ and $V$ vanish thanks to the indicators in~\eqref{eq:Veps-def-geom}, and there is nothing to prove. On this event, the $L^{1}(\Prob_{B})$ convergence~\eqref{eq:Aeps-inner-limit-geom}, summed over the pairs $1\le i<j\le k$, shows that the exponent in $V_{\varepsilon}$ converges in $\Prob_{B}$-probability to the exponent in $V$. It then suffices to compose with the continuous function $\exp$ and to multiply by the bounded random variable $\prod_{j=1}^{k}u_{0}(B^{x,j}_{t})$.

We now turn to boundedness in $L^{2}(\Prob_{B})$. Namely, bounding $|u_{0}|$ by $\|u_{0}\|_{\infty}$ in~\eqref{eq:Veps-def-geom} and using $\mathbf{1}_{\{t<T^{j}\}}^{2}=\mathbf{1}_{\{t<T^{j}\}}$, we get
\begin{equation}\label{eq:Veps-L2-bound-geom}
\E_{B}\bigl[V_{\varepsilon}^{2}\bigr]
\le
\|u_{0}\|_{\infty}^{2k}\,
\E_{B}\left[
\prod_{j=1}^{k}\mathbf{1}_{\{t<T^{j}\}}\,
\exp\left(
2\beta^{2}\sum_{1\le i<j\le k}
\langle A^{\varepsilon,i}_{t,x},A^{\varepsilon,j}_{t,x}\rangle_{\mathcal{H}^{D}_{\alpha}}
\right)
\right].
\end{equation}
Now the expectation on the right-hand side of~\eqref{eq:Veps-L2-bound-geom} is exactly the right-hand side of the moment formula~\eqref{eq:mom-eps-formula-geom}, applied with the initial condition $u_{0}\equiv1$ and with $\beta$ replaced by $\sqrt{2}\,\beta$. With a slight abuse of notation, we denote again by $u^{\varepsilon}$ the regularized process of Definition~\ref{def:ueps-geom} corresponding to those two parameters. Then the moment formula~\eqref{eq:mom-eps-formula-geom} recasts~\eqref{eq:Veps-L2-bound-geom} as
\begin{equation}\label{eq:Veps-L2-ustar-geom}
\E_{B}\bigl[V_{\varepsilon}^{2}\bigr]
\le
\|u_{0}\|_{\infty}^{2k}\,
\E\bigl[(u^{\varepsilon}(t,x))^{k}\bigr].
\end{equation}
Now observe that Proposition~\ref{prop:ueps-Lk-conv-geom} asserts that $u^{\varepsilon}(t,x)$ converges in $L^{k}(\Omega)$ as $\varepsilon\to0$. In particular the right-hand side of~\eqref{eq:Veps-L2-ustar-geom} is bounded for $\varepsilon$ in an interval $(0,\varepsilon_{0}]$. The family $\{V_{\varepsilon}\,;\,\varepsilon\in(0,\varepsilon_{0}]\}$ is therefore bounded in $L^{2}(\Prob_{B})$, and hence uniformly integrable.

Having checked its two hypotheses, Vitali's convergence theorem yields the claim~\eqref{eq:Veps-L1-limit-geom}. Combining~\eqref{eq:lhs-moment-limit-geom}, the moment formula~\eqref{eq:mom-eps-formula-geom} and~\eqref{eq:Veps-L1-limit-geom}, we thus obtain
\begin{equation*}
\E[u(t,x)^{k}]=\E_{B}[V],
\end{equation*}
which is exactly the Feynman-Kac formula~\eqref{eq:fk-moment-formula-geom}.
\end{proof}

\begin{remark}[The white noise case $\alpha=0$]\label{rem:fk-intersection-geom}
As in Remark~\ref{rmk:fk-intersection}, the kernel $G^{D,U}_{2\alpha}$ measures the time the paths $B^{x,i}$ and $B^{x,j}$ spend near each other. When $\alpha=0$, $G^{D,U}_{0}$ degenerates to a Dirac mass on the diagonal (Section~\ref{sec:geom-noise}), and the exponent in~\eqref{eq:fk-moment-formula-geom} would need to be interpreted through an intersection local time for two independent Brownian motions on $(X,d,\mu)$, an object that lies beyond the scope of this book. We therefore restrict Theorem~\ref{thm:fk-moments-geom} to $\alpha>0$, and refer to~\cite{hu:09:stochastic} for the analogous Euclidean local-time theory.
\end{remark}

\begin{remark}[Neumann boundary conditions]\label{rem:existence-neumann-fk-geom}
The same strategy applies to the Neumann PAM (Remark~\ref{rem:existence-neumann}): replacing $\Delta_{D,U}$, $p^{D,U}_{t}$, $W^{D}_{\alpha}$ and $G^{D,U}_{2\alpha}$ by $\Delta_{N,U}$, $p^{N,U}_{t}$, $W^{N}_{\alpha}$ and $G^{N,U}_{2\alpha}+C_{U}$ throughout Sections~\ref{sssec:geom-ueps}--\ref{sssec:geom-fk-formula} yields the Neumann counterpart of the Feynman-Kac moment formula~\eqref{eq:fk-moment-formula-geom}.
\end{remark}

The formula~\eqref{eq:fk-moment-formula-geom} will be our main tool for deriving \emph{lower} bounds on the moments of the solution. Before that, we establish the corresponding upper bounds, which follow a different route. As in the Euclidean Section~\ref{ssec:moment-assumption}, they rest on the chaos estimates of Section~\ref{ssec:geom-chaos-est}, combined with hypercontractivity.

\subsection{Upper bounds on the moments}\label{ssec:geom-moment-upper}

An analysis of the asymptotic behavior of the moments $\E[u(t,x)^{k}]$ as $t\to\infty$ brings useful information about the parabolic Anderson model, in particular about its intermittency. Recall that Theorem~\ref{thm:existence-geom} already contains an upper bound for the second moment, namely the estimate~\eqref{eq:L2-asmp-geom} involving the function $\Theta_{\alpha}$ and the spectral gap $\lambda_{1}$. In this subsection we extend~\eqref{eq:L2-asmp-geom} to moments of every order $k\ge2$. The strategy is the same as for the Euclidean Theorem~\ref{thm:moment-upper}: we combine the chaos expansion of $u$ with Nelson's hypercontractivity inequality (Theorem~\ref{thm:hypercontractivity}), which upgrades $L^{2}$ information on each chaos to $L^{k}$ information at the price of a factor $(k-1)^{n/2}$ on the $n$-th chaos. At the level of the kernel estimates, this price amounts to replacing the noise intensity $\beta$ by $\sqrt{k-1}\,\beta$, and the upper bound below is accordingly obtained from~\eqref{eq:L2-asmp-geom} through the substitution $\beta\mapsto\sqrt{k-1}\,\beta$. Note that the Feynman-Kac formula~\eqref{eq:fk-moment-formula-geom} is not needed at this point. It will come into play for the lower bounds.

\begin{theorem}[Upper bound on moments]\label{thm:moment-upper-geom}
Let $(X,d,\mu)$ satisfy Assumption~\ref{hyp:mmspace}, let $U\subset X$ be an inner uniform domain (Definition~\ref{def:inner-uniform}), let $\alpha\ge0$, $\beta>0$, and $u_{0}\in L^{\infty}(U,\mu)$. Let $u$ be the unique solution of~\eqref{eq:pam-D-geom} given by Theorem~\ref{thm:existence-geom}. Then for every integer $k\ge2$ and every $x\in U$ we have
\begin{equation}\label{eq:moment-upper-geom}
\limsup_{t\to+\infty}\frac{\ln\E\bigl[u(t,x)^{k}\bigr]}{t}
\;\le\;
\frac{k}{2}\,\Theta_{\alpha}\bigl(\sqrt{k-1}\,\beta\bigr)-k\lambda_{1},
\end{equation}
where $\lambda_{1}=\lambda_{1}(U)$ is the first eigenvalue of $-\Delta_{D,U}$ (Notation~\ref{not:spectrum-D}) and $\Theta_{\alpha}$ is the function of Theorem~\ref{thm:existence-geom}. For $k=2$ the right-hand side of~\eqref{eq:moment-upper-geom} is $\Theta_{\alpha}(\beta)-2\lambda_{1}$, so that~\eqref{eq:moment-upper-geom} reduces to the $L^{2}$ bound~\eqref{eq:L2-asmp-geom}.
\end{theorem}

\begin{proof}
We proceed in three steps, which parallel the three steps in the proof of the Euclidean Theorem~\ref{thm:moment-upper}. As we will see, Step 1 is identical in both cases, while Steps 2 and 3 are simpler here. Indeed, the Laplace-transform estimate~\eqref{eq:fk-bound-small} replaces the Mittag-Leffler analysis of the Euclidean proof. Throughout the proof we fix an integer $k\ge2$ and we set $M:=\|u_{0}\|_{\infty}$.

\smallskip

\noindent
{\it Step 1: Reduction to chaos norms via hypercontractivity.} Recall from~\eqref{eq:u-chaos-geom} that the solution can be decomposed as
\begin{equation*}
u(t,x)
=
P^{D,U}_{t}u_{0}(x)
+
\sum_{n=1}^{\infty}I_{n}(f_{n}(\,\cdot\,,t,x)),
\end{equation*}
with the chaos kernels $f_{n}$ of~\eqref{eq:fk-geom}. Nelson's inequality~\eqref{eq:hypercontractivity-sum} only depends on the underlying isonormal Gaussian structure, and therefore applies to the noise $W^{D}_{\alpha}$ exactly as it did in the proof of Proposition~\ref{prop:ueps-Lk-conv-geom}. Bounding the deterministic term $P^{D,U}_{t}u_{0}(x)$ by Corollary~\ref{cor:Linf-bound-D} it is readily checked that
\begin{equation}\label{eq:moment-hypercontract-geom}
\|u(t,x)\|_{L^{k}(\Omega)}
\;\le\;
C\,M\,e^{-\lambda_{1}t}
+
\sum_{n=1}^{\infty}(k-1)^{n/2}\,\sqrt{n!}\,\|f_{n}(\,\cdot\,,t,x)\|_{(\mathcal{H}^{D}_{\alpha})^{\otimes n}}.
\end{equation}
Relation~\eqref{eq:moment-hypercontract-geom} is the exact analogue of~\eqref{eq:moment-hypercontract} in the Euclidean case.

\smallskip

\noindent
{\it Step 2: Choice of $\rho$.} We now invoke the bound~\eqref{eq:fk-bound-small} of Lemma~\ref{lem:fk-bound-geom}, applied at chaos level $n$ (that is with $k$ replaced by $n$ in~\eqref{eq:fk-bound-small}). For every $\rho>0$ it yields
\begin{equation}\label{eq:fkn-bound-recall}
\sqrt{n!}\,\|f_{n}(\,\cdot\,,t,x)\|_{(\mathcal{H}^{D}_{\alpha})^{\otimes n}}
\;\le\;
M\,\bigl(C\beta^{2}\hat{\Psi}(\rho)\bigr)^{n/2}\,e^{(\rho-2\lambda_{1})t/2},
\end{equation}
where $\hat{\Psi}$ is the Laplace transform defined in~\eqref{eq:PsiHat-geom}. Multiplying~\eqref{eq:fkn-bound-recall} by $(k-1)^{n/2}$ and absorbing this factor into the noise intensity, the general term of the series in~\eqref{eq:moment-hypercontract-geom} is bounded by
\begin{equation}\label{eq:moment-term-geom}
(k-1)^{n/2}\,\sqrt{n!}\,\|f_{n}(\,\cdot\,,t,x)\|_{(\mathcal{H}^{D}_{\alpha})^{\otimes n}}
\;\le\;
M\,\Bigl(C\bigl(\sqrt{k-1}\,\beta\bigr)^{2}\,\hat{\Psi}(\rho)\Bigr)^{n/2}\,e^{(\rho-2\lambda_{1})t/2}.
\end{equation}
We already encountered the right-hand side of~\eqref{eq:moment-term-geom} in the proof of Proposition~\ref{prop:ueps-Lk-conv-geom} (see~\eqref{eq:hypercontract-term-geom}). There, any $\rho$ achieving $(k-1)C\beta^{2}\hat{\Psi}(\rho)<1$ was sufficient. Here we additionally track the smallest admissible $\rho$, since it dictates the exponential rate in~\eqref{eq:moment-upper-geom}.

To this aim, recall from Step 1 of the proof of Theorem~\ref{thm:existence-geom} (case (i), that is $0\le\alpha<d_{h}/(2d_{w})$) that the critical parameter $\rho_{c}(\beta)$ is defined by relation~\eqref{eq:rho-c-geom} as the unique solution of $C\beta^{2}F(\rho_{c})=1$, with $F$ given by~\eqref{eq:F-def-geom}, and that $\Theta_{\alpha}(\beta)=\rho_{c}(\beta)$. The key observation is that the noise intensity enters this definition only through the product $C\beta^{2}$. Hence the whole of Step 1 there can be applied verbatim with noise intensity $\sqrt{k-1}\,\beta$, producing the critical parameter
\begin{equation}\label{eq:rho-c-k-geom}
\rho_{c,k}
:=
\rho_{c}\bigl(\sqrt{k-1}\,\beta\bigr)
=
\Theta_{\alpha}\bigl(\sqrt{k-1}\,\beta\bigr).
\end{equation}
Exactly as in~\eqref{eq:delta-rho-geom}, for every $\rho>\rho_{c,k}$ we then get
\begin{equation}\label{eq:delta-rho-k-geom}
C(k-1)\beta^{2}\,\hat{\Psi}(\rho)
\;\le\;
\delta_{\rho,k}
\;<\;
1,
\qquad\text{where}\qquad
\delta_{\rho,k}:=C(k-1)\beta^{2}F(\rho).
\end{equation}
The same conclusion holds in the critical case (ii), replacing $F$ by the function $F_{2}$ of Step~4 of the proof of Theorem~\ref{thm:existence-geom}. It also holds in the supercritical case (iii), replacing $F$ by $F_{3}(\rho):=1/\rho$ (see~\eqref{eq:rho-c-super}). Since the three regimes are treated identically from this point on, we do not distinguish them in the sequel.

\smallskip

\noindent
{\it Step 3: Summation.} Fix $\rho>\Theta_{\alpha}(\sqrt{k-1}\,\beta)$. Plugging~\eqref{eq:moment-term-geom} and~\eqref{eq:delta-rho-k-geom} into~\eqref{eq:moment-hypercontract-geom}, the chaos series is dominated by the geometric series $\sum_{n\ge1}\delta_{\rho,k}^{n/2}=\delta_{\rho,k}^{1/2}/(1-\delta_{\rho,k}^{1/2})$. Bounding in addition the first term of~\eqref{eq:moment-hypercontract-geom} according to $e^{-\lambda_{1}t}\le e^{(\rho-2\lambda_{1})t/2}$ (recall that $\rho>0$), we end up with
\begin{equation}\label{eq:Lk-est-geom}
\|u(t,x)\|_{L^{k}(\Omega)}
\;\le\;
C_{\rho,k}\,M\,e^{(\rho-2\lambda_{1})t/2},
\end{equation}
for a constant $C_{\rho,k}$ depending on $\rho,k,\beta,\alpha,d_{h},d_{w},U$ but not on $t,x$. Estimate~\eqref{eq:Lk-est-geom} is the $L^{k}$ analogue of the $L^{2}$ estimate~\eqref{eq:L2-est-geom}. Observe that, in contrast with Step 3 of the Euclidean proof, no Mittag-Leffler asymptotics (see~\eqref{eq:ML-def} and~\eqref{eq:ml-bound}) is needed here. Indeed, the Laplace-transform bound~\eqref{eq:fk-bound-small} has already traded the factorial gain of the Dirichlet formula~\eqref{eq:dirichlet} for the exponential factor $e^{(\rho-2\lambda_{1})t}$, and the series is summed by merely tuning $\rho$.

We can now conclude. Since $\E[u(t,x)^{k}]\le\E[|u(t,x)|^{k}]=\|u(t,x)\|_{L^{k}(\Omega)}^{k}$, raising~\eqref{eq:Lk-est-geom} to the power $k$ gives
\begin{equation*}
\E\bigl[u(t,x)^{k}\bigr]
\;\le\;
C_{\rho,k}^{k}\,M^{k}\,e^{k(\frac{\rho}{2}-\lambda_{1})t},
\end{equation*}
and therefore
\begin{equation}\label{eq:limsup-rho-k-geom}
\limsup_{t\to+\infty}\frac{\ln\E\bigl[u(t,x)^{k}\bigr]}{t}
\;\le\;
k\left(\frac{\rho}{2}-\lambda_{1}\right),
\qquad\text{for every }\rho>\Theta_{\alpha}\bigl(\sqrt{k-1}\,\beta\bigr).
\end{equation}
Letting $\rho\downarrow\Theta_{\alpha}(\sqrt{k-1}\,\beta)$ in~\eqref{eq:limsup-rho-k-geom}, exactly as in~\eqref{eq:limsup-rho-geom}, yields~\eqref{eq:moment-upper-geom} and finishes the proof.
\end{proof}

\begin{remark}\label{rem:moment-growth-geom}
Let us compare~\eqref{eq:moment-upper-geom} with its Euclidean counterpart~\eqref{eq:moment-upper}, whose exponential rate is of order $k^{(4-a)/(2-a)}$, where $a$ is the exponent of the power-law relation~\eqref{eq:Lambda-power} in Assumption~\ref{hyp:Lambda-power}. Fix $\beta>0$, consider the subcritical regime $0<\alpha<d_{h}/(2d_{w})$, and let $k\to\infty$. The large-$\beta$ asymptotics of $\Theta_{\alpha}$ in case (i) of Theorem~\ref{thm:existence-geom} yield
\begin{equation}\label{eq:moment-growth-geom}
\frac{k}{2}\,\Theta_{\alpha}\bigl(\sqrt{k-1}\,\beta\bigr)
\;\asymp\;
k^{\frac{(2+2\alpha)d_{w}-d_{h}}{(1+2\alpha)d_{w}-d_{h}}},
\qquad k\to\infty.
\end{equation}
On the other hand, by Proposition~\ref{prop:RK-D-bounds} (i) (applied with $\alpha$ replaced by $2\alpha$), the covariance kernel $G^{D,U}_{2\alpha}$ of the noise behaves like $d(x,y)^{-a}$, where the parameter $a$ is defined by
\begin{equation}\label{eq:a-def}
a
:=
d_{h}-2\alpha d_{w}.
\end{equation}
This parameter will feature prominently in our computations below. The subcriticality of $\alpha$ amounts to $a>0$, while the recurrence condition $d_{h}<d_{w}$ of Assumption~\ref{hyp:mmspace} gives $a<d_{w}$. In the Euclidean normalization $d_{w}=2$ and $d_{h}=d$, we thus get $a\in(0,\min(2,d))$, which is the range prescribed for the exponent in~\eqref{eq:Lambda-power}, and a direct computation shows that the exponent in~\eqref{eq:moment-growth-geom} becomes
\begin{equation*}
\frac{(2+2\alpha)d_{w}-d_{h}}{(1+2\alpha)d_{w}-d_{h}}
=
\frac{4+4\alpha-d}{2+4\alpha-d}
=
\frac{4-a}{2-a}.
\end{equation*}
Hence it is possible to match the growth in $k$ of the geometric bound~\eqref{eq:moment-upper-geom}
with  the Euclidean exponent in~\eqref{eq:moment-upper}. Observe however that playing with $d_{w}$ and $d_{h}$, one gets a wide variety of growth exponents in~\eqref{eq:moment-upper-geom}.
\end{remark}

\subsection{Lower bounds on the moments}\label{ssec:geom-moment-lower}

We now derive lower bounds matching (at least partially) the upper bounds of Theorem~\ref{thm:moment-upper-geom}, following~\cite{baudoin:24:parabolic}. This is where the Feynman-Kac formula of Theorem~\ref{thm:fk-moments-geom} comes into play, as announced at the beginning of Section~\ref{ssec:geom-moment-upper}. Indeed,  the upper bounds only required the chaos expansion of $u$. However, the lower bounds hinge on the fact that the exponent in~\eqref{eq:fk-moment-formula-geom} is nonnegative and can be bounded below on suitable confinement events. Our strategy is modeled on the Euclidean Theorem~\ref{thm:moment-lower}, whose proof can be briefly summarized as follows:
\begin{enumerate}[label=\textbf{(\alph*)}]
\item
We restricted the expectation in the Feynman-Kac formula~\eqref{eq:fk-moment-formula} to the small ball event~\eqref{eq:Beps-def}, as mentioned in Remark~\ref{rmk:moment-estimates}.

\item
On this event the interaction $\Lambda(B^{i}_{s}-B^{j}_{s})$ is at least $c_{\Lambda}\varepsilon^{-a}$ owing to the power lower bound in~\eqref{eq:Lambda-power}, while the probability of the event itself is controlled by the Gaussian small ball estimate~\eqref{eq:small-ball}.

\item
An optimization over the radius $\varepsilon$ then produced the rate $k^{(4-a)/(2-a)}$ in~\eqref{eq:moment-lower}.
\end{enumerate}
 
\noindent 
The geometric proof will follow the same pattern, with three substitutions:

\begin{enumerate}

\item The small ball event~\eqref{eq:Beps-def} was defined through the coordinates $B^{i,l}$ of the Brownian motions, which are not available on $X$. It is replaced by the event that each path $B^{x,j}$ stays in a small ball $B(x,\varepsilon)$ up to time $t$, expressed by means of exit times as in Notation~\ref{not:BM-copies-geom}.

\item The power lower bound in~\eqref{eq:Lambda-power} is replaced by the lower bounds on the Riesz kernel $G^{D,U}_{2\alpha}$ gathered in Proposition~\ref{prop:RK-D-bounds}. Note that the three regimes of Proposition~\ref{prop:RK-D-bounds} will produce three different rates, mirroring the three regimes of the function $\Theta_{\alpha}$ in Theorem~\ref{thm:existence-geom}. In addition, the eigenfunction factor $\Phi_{1}(y)\Phi_{1}(z)$ in those lower bounds vanishes at the boundary, which forces us to confine the paths in a compact subset of $U$. 

\item The Gaussian small ball estimate~\eqref{eq:small-ball} is replaced by a spectral estimate for the exit time of a small ball (see Lemma~\ref{lem:small-ball-geom} below). The confinement cost over a time interval of length $t$ becomes $\exp(-c\,t/\varepsilon^{d_{w}})$ instead of $\exp(-C_{0}t/\varepsilon^{2})$, in accordance with the interpretation of the walk dimension $d_{w}$ given in Remark~\ref{rmk:g1}.

\end{enumerate}

We start with the small ball estimate mentioned in item (iii) above. 

\begin{lemma}[Small ball estimates]\label{lem:small-ball-geom}
Let $(X,d,\mu)$ satisfy Assumption~\ref{hyp:mmspace}, and let $U\subset X$ be a bounded open set. We consider $x\in U$ and $\varepsilon>0$ such that $\overline{B(x,\varepsilon)}\subset U$. Write $V_{\varepsilon}:=B(x,\varepsilon)$, and let
\begin{equation}\label{eq:Teps-def}
T_{\varepsilon}:=\inf\bigl\{s\ge0:\,B^{x}_{s}\notin V_{\varepsilon}\bigr\}
\end{equation}
be the exit time from $V_{\varepsilon}$ of the Brownian motion $B^{x}$ started at $x$ (Definition~\ref{def:BM}). Moreover, let $\lambda_{1}(V_{\varepsilon})$ be the smallest Dirichlet eigenvalue of $V_{\varepsilon}$, that is the bottom of the spectrum of the operator $-\Delta_{D,V_{\varepsilon}}$ given by Proposition~\ref{prop:E-D}. Then the following holds:
\begin{enumerate}

\item The exit time $T_{\varepsilon}$ satisfies the spectral lower bound
\begin{equation}\label{eq:survival-spectral-geom}
\liminf_{t\to+\infty}\frac{\ln\Prob_{B}(T_{\varepsilon}>t)}{t}
\;\ge\;
-\lambda_{1}(V_{\varepsilon}).
\end{equation}

\item There exists a constant $c_{\lambda}>0$, depending only on the constants of Assumption~\ref{hyp:mmspace}, such that
\begin{equation}\label{eq:lambda1-ball-geom}
\lambda_{1}(V_{\varepsilon})
\;\le\;
\frac{c_{\lambda}}{\varepsilon^{d_{w}}}.
\end{equation}

\end{enumerate}
\end{lemma}

\begin{proof}
We focus on the proof of~\eqref{eq:survival-spectral-geom}. Indeed, the eigenvalue bound~\eqref{eq:lambda1-ball-geom} is a consequence of the sub-Gaussian estimate~\eqref{eq:subGauss}, and we refer to~\cite[Lem.~7.1]{mariano:22:spectral} for a proof.

In order to prove~\eqref{eq:survival-spectral-geom}, observe that the ball $V_{\varepsilon}$ is a bounded open set with compact closure. Therefore the constructions of Section~\ref{ssec:dirichlet-bc} apply with $V_{\varepsilon}$ in place of $U$. Namely Proposition~\ref{prop:E-D} yields the Dirichlet Laplacian $\Delta_{D,V_{\varepsilon}}$ and the associated semigroup $P^{D,V_{\varepsilon}}_{t}$, while the representation~\eqref{eq:PDU-prob}, applied on $V_{\varepsilon}$ with $f\equiv1$, reads
\begin{equation}\label{eq:survival-semigroup-geom}
\Prob_{B}\bigl(T_{\varepsilon}>t\bigr)
\;=\;
P^{D,V_{\varepsilon}}_{t}\mathbf{1}(x).
\end{equation}
In addition, exactly as in Notation~\ref{not:spectrum-D}, the operator $-\Delta_{D,V_{\varepsilon}}$ admits a discrete spectrum $0<\lambda_{1}(V_{\varepsilon})\le\lambda_{2}(V_{\varepsilon})\le\cdots$, together with an orthonormal basis $\{\Phi^{\varepsilon}_{j}\}_{j\ge1}$ of $L^{2}(V_{\varepsilon},\mu)$ made of eigenfunctions. Observe that $V_{\varepsilon}$ is a connected set. Indeed the space $(X,d)$ is geodesic (Definition~\ref{def:mmspace}), and for $y\in V_{\varepsilon}$ every point $w$ on a geodesic from $x$ to $y$ satisfies $d(x,w)\le d(x,y)<\varepsilon$. Hence, owing to the Perron-Frobenius argument recalled in Remark~\ref{rem:ground-state-positive}, the ground state $\Phi^{\varepsilon}_{1}$ can be chosen to be strictly positive on $V_{\varepsilon}$.

Next we claim that $\Phi^{\varepsilon}_{1}$ is a bounded function. To this aim, fix $t_{0}\in(0,\mathrm{diam}(X)^{d_{w}}/2)$. The eigenfunction relation gives
\begin{equation}\label{j21}
\Phi^{\varepsilon}_{1}=e^{\lambda_{1}(V_{\varepsilon})t_{0}}P^{D,V_{\varepsilon}}_{t_{0}}\Phi^{\varepsilon}_{1}.
\end{equation}
 Moreover, the semigroup $P^{D,V_{\varepsilon}}_{t_{0}}$ admits a kernel $p^{D,V_{\varepsilon}}_{t_{0}}$, given by the analogue of the spectral expansion~\eqref{eq:HK-spectral-D} on $V_{\varepsilon}$, and this kernel satisfies $p^{D,V_{\varepsilon}}_{t_{0}}\le p_{t_{0}}$ exactly as in Remark~\ref{rem:Phi1-positive}. Hence for every $y\in V_{\varepsilon}$ we get
\begin{equation}\label{j3}
\bigl|P^{D,V_{\varepsilon}}_{t_{0}}\Phi^{\varepsilon}_{1}(y)\bigr|
\;\le\;
\Bigl(\int_{V_{\varepsilon}}p^{D,V_{\varepsilon}}_{t_{0}}(y,z)^{2}\,d\mu(z)\Bigr)^{1/2}
\;=\;
p^{D,V_{\varepsilon}}_{2t_{0}}(y,y)^{1/2}
\;\le\;
p_{2t_{0}}(y,y)^{1/2}.
\end{equation}
Let us detail the three relations in~\eqref{j3}. The first inequality is obtained thanks to Cauchy-Schwarz and the normalization $\|\Phi^{\varepsilon}_{1}\|_{L^{2}(V_{\varepsilon},\mu)}=1$. The middle identity stems from the symmetry of the kernel $p^{D,V_{\varepsilon}}_{t_{0}}$ and the semigroup property, as in Corollary~\ref{cor:L2-bound-D}. The last inequality is due to the domination $p^{D,V_{\varepsilon}}_{2t_{0}}\le p_{2t_{0}}$ recalled above. Next, resorting to the on-diagonal upper bound in the sub-Gaussian estimate~\eqref{eq:subGauss}, we obtain
\begin{equation}\label{j4}
p_{2t_{0}}(y,y)^{1/2}
\;\le\;
c\,t_{0}^{-d_{h}/(2d_{w})}.
\end{equation}
Gathering~\eqref{j3} and~\eqref{j4} and recalling the eigenfunction relation~\eqref{j21}, we end up with
\begin{equation*}
\|\Phi^{\varepsilon}_{1}\|_{\infty}
\;\le\;
c\,e^{\lambda_{1}(V_{\varepsilon})t_{0}}\,t_{0}^{-d_{h}/(2d_{w})}
\;<\;
+\infty .
\end{equation*}
This shows that $\Phi^{\varepsilon}_{1}$ is a bounded function on $V_{\varepsilon}$, as claimed.

We can now conclude. Since $\Phi^{\varepsilon}_{1}$ is positive and bounded we have $\mathbf{1}\ge\Phi^{\varepsilon}_{1}/\|\Phi^{\varepsilon}_{1}\|_{\infty}$ on $V_{\varepsilon}$. Furthermore, $P^{D,V_{\varepsilon}}_{t}$ preserves positivity (which is clear from the representation~\eqref{eq:PDU-prob}). Hence~\eqref{eq:survival-semigroup-geom} yields
\begin{equation}\label{eq:survival-ground-state}
\Prob_{B}\bigl(T_{\varepsilon}>t\bigr)
\;\ge\;
\frac{P^{D,V_{\varepsilon}}_{t}\Phi^{\varepsilon}_{1}(x)}{\|\Phi^{\varepsilon}_{1}\|_{\infty}}
\;=\;
\frac{\Phi^{\varepsilon}_{1}(x)}{\|\Phi^{\varepsilon}_{1}\|_{\infty}}\,
e^{-\lambda_{1}(V_{\varepsilon})t}.
\end{equation}
Since $\Phi^{\varepsilon}_{1}(x)>0$, taking logarithms in~\eqref{eq:survival-ground-state}, dividing by $t$ and letting $t\to+\infty$ we obtain~\eqref{eq:survival-spectral-geom}.
\end{proof}

With Lemma~\ref{lem:small-ball-geom} in hand, we can turn to the main results of this subsection, namely the lower bounds on the moments of $u$. Since those results share a common setting, we first gather our standing assumptions.

\begin{assumption}[Setting for the moment lower bounds]\label{hyp:moment-lower-geom}
As in Theorem~\ref{thm:fk-moments-geom}, let $(X,d,\mu)$ satisfy Assumption~\ref{hyp:mmspace}, let $U\subset X$ be an inner uniform domain (Definition~\ref{def:inner-uniform}), let $\alpha>0$, $\beta>0$, and let $u$ be the unique solution of~\eqref{eq:pam-D-geom} given by Theorem~\ref{thm:existence-geom}. We assume that $u_{0}\in L^{\infty}(U,\mu)$ is nonnegative and bounded away from $0$ near a point $x\in U$. That is, we suppose that there exist $m_{0}>0$ and $\varepsilon_{0}>0$ such that
\begin{equation}\label{eq:u0-local-lower}
\overline{B(x,\varepsilon_{0})}\subset U,
\qquad\text{and}\qquad
u_{0}\ge m_{0}
\quad\text{$\mu$-a.e. on }B(x,\varepsilon_{0}).
\end{equation}
We also introduce the reduced radius
\begin{equation}\label{eq:eps1-def}
\varepsilon_{1}
:=
\frac{1}{4}\min(\varepsilon_{0},1),
\end{equation}
which will serve as a maximal confinement radius in the proofs below.
\end{assumption}

Our first lower bound is valid for every $\alpha>0$ and every noise intensity $\beta>0$. Its proof also contains the confinement procedure on which the other lower bounds of this subsection rely.

\begin{theorem}[Lower bound on moments, general case]\label{thm:moment-lower-gen-geom}
Let Assumption~\ref{hyp:moment-lower-geom} prevail. Then one can find a strictly positive constant $c_{1}$, depending on $x$, $\varepsilon_{0}$, $\alpha$, $U$ and the constants of Assumption~\ref{hyp:mmspace}, but not on $k$ or $\beta$, such that for every integer $k\ge2$ we have
\begin{equation}\label{eq:moment-lower-gen-geom}
\liminf_{t\to+\infty}\frac{\ln\E\bigl[u(t,x)^{k}\bigr]}{t}
\;\ge\;
c_{1}\,\beta^{2}k(k-1)-k\,\lambda_{1}\bigl(B(x,\varepsilon_{1})\bigr),
\end{equation}
where $\lambda_{1}(B(x,\varepsilon_{1}))$ is the first Dirichlet eigenvalue of the ball $B(x,\varepsilon_{1})$, as defined in Lemma~\ref{lem:small-ball-geom}, and $\varepsilon_{1}$ is the reduced radius~\eqref{eq:eps1-def}.
\end{theorem}

\begin{proof}
We divide this proof in three steps, which are the counterparts of Steps 1--3 in the proof of Theorem~\ref{thm:moment-lower}. Throughout the proof, $\varepsilon\in(0,\varepsilon_{1}]$ denotes a confinement radius, where $\varepsilon_{1}$ is the reduced radius~\eqref{eq:eps1-def}. Although the elementary choice $\varepsilon=\varepsilon_{1}$ will suffice for~\eqref{eq:moment-lower-gen-geom}, we carry out Steps 1 and 2 for a generic $\varepsilon\in(0,\varepsilon_{1}]$. In this way, the intermediate bound~\eqref{eq:liminf-eps-geom} below can be recycled for the optimization procedures leading to Theorems~\ref{thm:moment-lower-sub-geom} and~\ref{thm:moment-lower-crit-geom}. Observe that the reduction from $\varepsilon_{0}$ to $\varepsilon_{1}$ grants two properties which will be used repeatedly below. Namely for every $\varepsilon\in(0,\varepsilon_{1}]$ we have $2\varepsilon\le1/2$, as well as the inclusions
\begin{equation}\label{eq:3eps-in-U}
\overline{B(x,3\varepsilon)}
\subset
\overline{B(x,\varepsilon_{0})}
\subset
U .
\end{equation}

\smallskip

\noindent
{\it Step 1: Confinement.} For $1\le j\le k$ let
\begin{equation}\label{eq:Tj-eps-def}
T^{j}_{\varepsilon}
:=
\inf\bigl\{s\ge0:\,B^{x,j}_{s}\notin B(x,\varepsilon)\bigr\} ,
\end{equation}
where $B^{x,1},\ldots,B^{x,k}$ are the independent Brownian copies of Notation~\ref{not:BM-copies-geom}. Then we consider the event
\begin{equation}\label{eq:Aeps-event-geom}
\mathcal{A}_{\varepsilon,t}
:=
\bigcap_{j=1}^{k}\bigl\{t<T^{j}_{\varepsilon}\bigr\}.
\end{equation}
The event $\mathcal{A}_{\varepsilon,t}$ is the analogue of the small ball event~\eqref{eq:Beps-def}, with the coordinatewise conditions $\sup_{0\le s\le t}|B^{i,l}_{s}|\le\varepsilon/2$ replaced by exit time conditions. Since $B(x,\varepsilon)\subset U$, we have $T^{j}_{\varepsilon}\le T^{j}$, and thus $\mathbf{1}_{\{t<T^{j}\}}\ge\mathbf{1}_{\{t<T^{j}_{\varepsilon}\}}$ for every $j$. In addition, on the event $\mathcal{A}_{\varepsilon,t}$ we have $B^{x,j}_{t}\in B(x,\varepsilon)\subset B(x,\varepsilon_{0})$, so that $u_{0}(B^{x,j}_{t})\ge m_{0}$ almost surely thanks to~\eqref{eq:u0-local-lower}. Since $u_{0}$ is nonnegative, the integrand in the Feynman-Kac formula~\eqref{eq:fk-moment-formula-geom} is nonnegative, and restricting the expectation to $\mathcal{A}_{\varepsilon,t}$ we get the analogue of~\eqref{eq:fk-lower-start}:
\begin{equation}\label{eq:fk-lower-start-geom}
\E\bigl[u(t,x)^{k}\bigr]
\;\ge\;
m_{0}^{k}\,
\E_{B}\left[
\exp\left(
\beta^{2}\!\!\sum_{1\le i<j\le k}
\int_{0}^{t}G^{D,U}_{2\alpha}(B^{x,i}_{s},B^{x,j}_{s})\,ds
\right)\mathbf{1}_{\mathcal{A}_{\varepsilon,t}}
\right].
\end{equation}

We now bound the interaction term in~\eqref{eq:fk-lower-start-geom} from below on the event $\mathcal{A}_{\varepsilon,t}$. This is the place where the power lower bound of Assumption~\ref{hyp:Lambda-power} is replaced by Proposition~\ref{prop:RK-D-bounds}. Fix $s\in[0,t]$ and a pair $i<j$, and abbreviate $y:=B^{x,i}_{s}$, $z:=B^{x,j}_{s}$, so that $y,z\in B(x,\varepsilon)$ on $\mathcal{A}_{\varepsilon,t}$. Our first observation is that the inner distance between $y$ and $z$ is small. Indeed the space $(X,d)$ is geodesic (Definition~\ref{def:mmspace}), so $y$ and $z$ are joined by a curve $\gamma$ of length $d(y,z)\le2\varepsilon$. Every point $w$ on $\gamma$ satisfies $d(w,x)\le d(w,y)+d(y,x)<d(y,z)+\varepsilon\le3\varepsilon$, so that $\gamma\subset B(x,3\varepsilon)\subset U$ thanks to the inclusions~\eqref{eq:3eps-in-U}. According to the definition~\eqref{eq:inner-metric} of the inner metric, we thus get
\begin{equation}\label{eq:dU-confin}
d_{U}(y,z)
\;\le\;
\mathrm{length}(\gamma)
\;=\;
d(y,z)
\;\le\;
2\varepsilon .
\end{equation}
Our second observation is that, the set $\overline{B(x,\varepsilon_{0})}$ being a compact subset of $U$, the quantity
\begin{equation}\label{eq:m-hat-def}
\hat{m}
:=
\inf\bigl\{\Phi_{1}(w)\,;\;w\in\overline{B(x,\varepsilon_{0})}\bigr\}
\end{equation}
is strictly positive (see Remark~\ref{rem:Phi1-positive}). We now plug~\eqref{eq:dU-confin} and~\eqref{eq:m-hat-def} into the lower bounds~\eqref{eq:RK-D-small}, \eqref{eq:RK-D-crit} and~\eqref{eq:RK-D-large} of Proposition~\ref{prop:RK-D-bounds}, applied with $\alpha$ replaced by $2\alpha$ (the three cases there correspond to $a>0$, $a=0$ and $a<0$ respectively, for the parameter $a=d_{h}-2\alpha d_{w}$ introduced in~\eqref{eq:a-def}). Taking into account the fact that $r\mapsto r^{-a}$ is decreasing for $a>0$, as well as the bound $\max(1,|\ln r|)\ge\ln(1/r)\ge\ln(1/(2\varepsilon))$ valid for $0<r\le2\varepsilon\le1/2$ (recall that $2\varepsilon\le1/2$ thanks to our reduction $\varepsilon\le\varepsilon_{1}$), we obtain
\begin{equation}\label{eq:G-lower-ball}
G^{D,U}_{2\alpha}(y,z)
\;\ge\;
c\,\hat{m}^{2}\,g_{\alpha}(\varepsilon),
\qquad\text{where}\quad
g_{\alpha}(\varepsilon)
:=
\begin{cases}
(2\varepsilon)^{-a}, & \text{if } a>0\\[4pt]
\ln\bigl(\tfrac{1}{2\varepsilon}\bigr), & \text{if } a=0\\[4pt]
1, & \text{if } a<0 ,
\end{cases}
\end{equation}
and where the bound~\eqref{eq:G-lower-ball} trivially holds when $y=z$. Note for further use that $2\varepsilon\le1/2$ implies $g_{\alpha}(\varepsilon)\ge\ln2$ in all three cases. Reporting~\eqref{eq:G-lower-ball} in~\eqref{eq:fk-lower-start-geom}, and recalling that the sum over pairs contains $k(k-1)/2$ terms, we end up with the analogue of~\eqref{eq:fk-restricted}:
\begin{equation}\label{eq:fk-restricted-geom}
\E\bigl[u(t,x)^{k}\bigr]
\;\ge\;
m_{0}^{k}\,
\exp\Bigl(c_{0}\,\beta^{2}k(k-1)\,g_{\alpha}(\varepsilon)\,t\Bigr)\,
\Prob_{B}\bigl(\mathcal{A}_{\varepsilon,t}\bigr),
\qquad\text{where}\quad
c_{0}:=\frac{c\,\hat{m}^{2}}{2}.
\end{equation}

\smallskip

\noindent
{\it Step 2: Small ball estimate.} The paths $B^{x,1},\ldots,B^{x,k}$ being independent copies (Notation~\ref{not:BM-copies-geom}), we have, exactly as in~\eqref{f3},
\begin{equation}\label{eq:small-ball-product-geom}
\Prob_{B}\bigl(\mathcal{A}_{\varepsilon,t}\bigr)
=
\bigl[\Prob_{B}\bigl(T^{1}_{\varepsilon}>t\bigr)\bigr]^{k}.
\end{equation}
The role of the Gaussian small ball estimate~\eqref{eq:small-ball} is now played by Lemma~\ref{lem:small-ball-geom} (observe that $\overline{B(x,\varepsilon)}\subset\overline{B(x,\varepsilon_{0})}\subset U$ by~\eqref{eq:3eps-in-U}, so that the lemma indeed applies to $V_{\varepsilon}=B(x,\varepsilon)$). Namely, taking logarithms in~\eqref{eq:fk-restricted-geom}, dividing by $t$ and letting $t\to+\infty$, relations~\eqref{eq:small-ball-product-geom}, \eqref{eq:survival-spectral-geom} and~\eqref{eq:lambda1-ball-geom} yield
\begin{eqnarray}\label{eq:liminf-eps-geom}
\liminf_{t\to+\infty}\frac{\ln\E\bigl[u(t,x)^{k}\bigr]}{t}
&\ge&
c_{0}\,\beta^{2}k(k-1)\,g_{\alpha}(\varepsilon)-k\,\lambda_{1}\bigl(B(x,\varepsilon)\bigr) \notag\\
&\ge&
c_{0}\,\beta^{2}k(k-1)\,g_{\alpha}(\varepsilon)-\frac{c_{\lambda}\,k}{\varepsilon^{d_{w}}}.
\end{eqnarray}
Relation~\eqref{eq:liminf-eps-geom} is the counterpart of~\eqref{eq:fk-opt-prep}, with two differences. First, it is stated at the level of exponential rates rather than for every fixed $t$, since the survival estimate~\eqref{eq:survival-spectral-geom} is itself asymptotic. Second, the confinement cost is now $c_{\lambda}k\varepsilon^{-d_{w}}$ instead of $C_{0}dk\varepsilon^{-2}$, the walk dimension $d_{w}$ replacing the Euclidean exponent $2$.

\smallskip

\noindent
{\it Step 3: Choice of $\varepsilon$.} For the bound~\eqref{eq:moment-lower-gen-geom} no optimization is needed, and we simply take $\varepsilon=\varepsilon_{1}$ in~\eqref{eq:liminf-eps-geom}. Bounding $g_{\alpha}(\varepsilon_{1})\ge\ln2$ (as noted after~\eqref{eq:G-lower-ball}) and resorting to the first inequality in~\eqref{eq:liminf-eps-geom}, we obtain~\eqref{eq:moment-lower-gen-geom} with $c_{1}:=c_{0}\ln2$.
\end{proof}

When $(k-1)\beta^{2}$ is large, the general bound~\eqref{eq:moment-lower-gen-geom} can be improved by tuning the confinement radius $\varepsilon$ in the intermediate inequality~\eqref{eq:liminf-eps-geom}. Our next result implements this optimization in the subcritical case $\alpha<\frac{d_{h}}{2d_{w}}$. Recall that the parameter $a=d_{h}-2\alpha d_{w}$ has been defined in~\eqref{eq:a-def}, so that the subcritical case corresponds to $a>0$. Observe also that the exponent $\frac{d_{w}}{d_{w}-a}$ in~\eqref{eq:moment-lower-sub-geom} is strictly greater than $1$ when $a>0$. Hence Theorem~\ref{thm:moment-lower-sub-geom} improves on the growth given by~\eqref{eq:moment-lower-gen-geom} for large values of $(k-1)\beta^{2}$.

\begin{theorem}[Lower bound on moments, subcritical case]\label{thm:moment-lower-sub-geom}
Let Assumption~\ref{hyp:moment-lower-geom} prevail, and suppose that $0<\alpha<\frac{d_{h}}{2d_{w}}$ (equivalently $a>0$, where $a$ is the parameter defined in~\eqref{eq:a-def}). Then one can find strictly positive constants $c_{2},Q$, depending on $x$, $\varepsilon_{0}$, $\alpha$, $U$ and the constants of Assumption~\ref{hyp:mmspace}, but not on $k$ or $\beta$, such that for every integer $k\ge2$ and every $\beta$ satisfying $(k-1)\beta^{2}\ge Q$ we have
\begin{equation}\label{eq:moment-lower-sub-geom}
\liminf_{t\to+\infty}\frac{\ln\E\bigl[u(t,x)^{k}\bigr]}{t}
\;\ge\;
c_{2}\,k\,\bigl((k-1)\beta^{2}\bigr)^{\frac{d_{w}}{(1+2\alpha)d_{w}-d_{h}}}
\;=\;
c_{2}\,k\,\bigl((k-1)\beta^{2}\bigr)^{\frac{d_{w}}{d_{w}-a}}.
\end{equation}
\end{theorem}

\begin{proof}
Steps 1 and 2 in the proof of Theorem~\ref{thm:moment-lower-gen-geom} have been carried out for a generic confinement radius $\varepsilon\in(0,\varepsilon_{1}]$. We can thus start from the intermediate bound~\eqref{eq:liminf-eps-geom} obtained there, and our task is reduced to optimizing over $\varepsilon$. In the computations below, $c_{0}$ and $c_{\lambda}$ respectively denote the constants featured in~\eqref{eq:fk-restricted-geom} and~\eqref{eq:lambda1-ball-geom}.

We first note that $a<d_{w}$, since the definition~\eqref{eq:a-def} of $a$ gives $a<d_{h}$, and $d_{h}<d_{w}$ by the recurrence condition in Assumption~\ref{hyp:mmspace}. Owing to~\eqref{eq:liminf-eps-geom} and the definition of $g_{\alpha}$ in~\eqref{eq:G-lower-ball}, we are reduced to maximize over $\varepsilon\in(0,\varepsilon_{1}]$ the function
\begin{equation}\label{eq:h-eps-def}
h(\varepsilon)
:=
A\,\varepsilon^{-a}-c_{\lambda}\,k\,\varepsilon^{-d_{w}},
\qquad\text{where}\quad
A:=2^{-a}\,c_{0}\,\beta^{2}k(k-1).
\end{equation}
This is the same optimization as in Step 3 of the proof of Theorem~\ref{thm:moment-lower}, with the Euclidean parameters $(2,C_{0}dk)$ replaced by $(d_{w},c_{\lambda}k)$. Specifically, setting the derivative of $h$ to zero gives the following analogue of~\eqref{eq:eps-star}:
\begin{equation}\label{eq:eps-star-geom}
\varepsilon_{\ast}
=
\Bigl(\frac{d_{w}\,c_{\lambda}\,k}{a\,A}\Bigr)^{\frac{1}{d_{w}-a}}
=
\Bigl(\frac{2^{a}\,d_{w}\,c_{\lambda}}{a\,c_{0}}\Bigr)^{\frac{1}{d_{w}-a}}
\bigl((k-1)\beta^{2}\bigr)^{-\frac{1}{d_{w}-a}},
\end{equation}
where the second identity stems from the definition of $A$ in~\eqref{eq:h-eps-def}. From the expression~\eqref{eq:eps-star-geom} it is readily checked that $\varepsilon_{\ast}\le\varepsilon_{1}$ if and only if
\begin{equation}\label{eq:Q-sub-geom}
(k-1)\beta^{2}
\;\ge\;
\frac{2^{a}\,d_{w}\,c_{\lambda}}{a\,c_{0}\,\varepsilon_{1}^{d_{w}-a}},
\end{equation}
which is granted by our standing assumption $(k-1)\beta^{2}\ge Q$, provided we choose $Q$ to be the right-hand side of~\eqref{eq:Q-sub-geom}. At the point $\varepsilon_{\ast}$ the optimality condition reads $c_{\lambda}k\varepsilon_{\ast}^{-d_{w}}=\frac{a}{d_{w}}A\varepsilon_{\ast}^{-a}$, so that, exactly as in Step 3 of the proof of Theorem~\ref{thm:moment-lower},
\begin{equation}\label{eq:h-eps-star}
h(\varepsilon_{\ast})
=
\Bigl(1-\frac{a}{d_{w}}\Bigr)A\,\varepsilon_{\ast}^{-a}.
\end{equation}
Substituting the value~\eqref{eq:eps-star-geom} of $\varepsilon_{\ast}$ into~\eqref{eq:h-eps-star}, and resorting to the elementary identity $1+\frac{a}{d_{w}-a}=\frac{d_{w}}{d_{w}-a}$ (which replaces the Euclidean identity $2+\frac{a}{2-a}=\frac{4-a}{2-a}$ invoked for Theorem~\ref{thm:moment-lower}), we end up with
\begin{equation}\label{eq:h-eps-star-value}
h(\varepsilon_{\ast})
=
\Bigl(1-\frac{a}{d_{w}}\Bigr)2^{-a}c_{0}
\Bigl(\frac{a\,c_{0}}{2^{a}\,d_{w}\,c_{\lambda}}\Bigr)^{\frac{a}{d_{w}-a}}
\,k\,\bigl((k-1)\beta^{2}\bigr)^{\frac{d_{w}}{d_{w}-a}}
\equiv
c_{2}\,k\,\bigl((k-1)\beta^{2}\bigr)^{\frac{d_{w}}{d_{w}-a}},
\end{equation}
where $c_{2}>0$ thanks to $a<d_{w}$. Recalling that $d_{w}-a=(1+2\alpha)d_{w}-d_{h}$, relation~\eqref{eq:h-eps-star-value} proves~\eqref{eq:moment-lower-sub-geom}.
\end{proof}

We now turn to the critical case $\alpha=\frac{d_{h}}{2d_{w}}$, that is $a=0$, for which the improvement on the general bound~\eqref{eq:moment-lower-gen-geom} is of logarithmic order.

\begin{theorem}[Lower bound on moments, critical case]\label{thm:moment-lower-crit-geom}
Let Assumption~\ref{hyp:moment-lower-geom} prevail, and suppose that $\alpha=\frac{d_{h}}{2d_{w}}$ (equivalently $a=0$, where $a$ is the parameter defined in~\eqref{eq:a-def}). Then one can find strictly positive constants $c_{3},Q$, depending on $x$, $\varepsilon_{0}$, $\alpha$, $U$ and the constants of Assumption~\ref{hyp:mmspace}, but not on $k$ or $\beta$, such that for every integer $k\ge2$ and every $\beta$ satisfying $(k-1)\beta^{2}\ge Q$ we have
\begin{equation}\label{eq:moment-lower-crit-geom}
\liminf_{t\to+\infty}\frac{\ln\E\bigl[u(t,x)^{k}\bigr]}{t}
\;\ge\;
c_{3}\,k(k-1)\beta^{2}\,\ln\bigl((k-1)\beta^{2}\bigr).
\end{equation}
\end{theorem}

\begin{proof}
As in the proof of Theorem~\ref{thm:moment-lower-sub-geom}, we start from the bound~\eqref{eq:liminf-eps-geom}, established in Steps 1 and 2 of the proof of Theorem~\ref{thm:moment-lower-gen-geom} for a generic $\varepsilon\in(0,\varepsilon_{1}]$. We also keep the notation $c_{0}$ and $c_{\lambda}$ for the constants in~\eqref{eq:fk-restricted-geom} and~\eqref{eq:lambda1-ball-geom}. In the current situation we have $a=0$, and thus the function $g_{\alpha}$ in~\eqref{eq:G-lower-ball} is given by $g_{\alpha}(\varepsilon)=\ln(\frac{1}{2\varepsilon})$. The two terms in~\eqref{eq:liminf-eps-geom} are then balanced, up to the logarithmic factor, by the choice
\begin{equation}\label{eq:eps-crit-geom}
\varepsilon
:=
\bigl((k-1)\beta^{2}\bigr)^{-1/d_{w}},
\end{equation}
which satisfies $\varepsilon\le\varepsilon_{1}$ as soon as $(k-1)\beta^{2}\ge\varepsilon_{1}^{-d_{w}}$. With the choice~\eqref{eq:eps-crit-geom} we have $\varepsilon^{-d_{w}}=(k-1)\beta^{2}$ and $\ln(\frac{1}{2\varepsilon})=\frac{1}{d_{w}}\ln((k-1)\beta^{2})-\ln2$. Hence~\eqref{eq:liminf-eps-geom} becomes
\begin{equation}\label{eq:liminf-crit-geom}
\liminf_{t\to+\infty}\frac{\ln\E\bigl[u(t,x)^{k}\bigr]}{t}
\;\ge\;
k(k-1)\beta^{2}
\Bigl(\frac{c_{0}}{d_{w}}\ln\bigl((k-1)\beta^{2}\bigr)-c_{0}\ln2-c_{\lambda}\Bigr).
\end{equation}
If in addition $\ln((k-1)\beta^{2})\ge\frac{2d_{w}}{c_{0}}(c_{0}\ln2+c_{\lambda})$, the parenthesis in~\eqref{eq:liminf-crit-geom} is bounded below by $\frac{c_{0}}{2d_{w}}\ln((k-1)\beta^{2})$, which yields~\eqref{eq:moment-lower-crit-geom} with $c_{3}:=\frac{c_{0}}{2d_{w}}$. Summarizing, it suffices to choose the constant $Q$ in the statement so that $(k-1)\beta^{2}\ge Q$ implies both $(k-1)\beta^{2}\ge\varepsilon_{1}^{-d_{w}}$ and $\ln((k-1)\beta^{2})\ge\frac{2d_{w}}{c_{0}}(c_{0}\ln2+c_{\lambda})$. The proof is now complete.
\end{proof}

\begin{remark}[Comparison with the upper bounds]\label{rem:moment-lower-comp-geom}
Let us compare the lower bounds of Theorems~\ref{thm:moment-lower-gen-geom}, \ref{thm:moment-lower-sub-geom} and~\ref{thm:moment-lower-crit-geom} with the upper bounds of Theorem~\ref{thm:moment-upper-geom} regime by regime, in the spirit of Remark~\ref{rmk:moment-estimates}. In order to ease the reading, all the exponents below are expressed in terms of the parameter $a=d_{h}-2\alpha d_{w}$ introduced in~\eqref{eq:a-def}.

\begin{enumerate}

\item {\it Subcritical regime: matching rates.}
Let $0<\alpha<\frac{d_{h}}{2d_{w}}$, that is $a>0$. The large-$\beta$ asymptotics~\eqref{eq:theta-large-beta} of $\Theta_{\alpha}$ show that the upper bound in~\eqref{eq:moment-upper-geom} satisfies
\begin{equation}\label{eq:upper-rate-sub-geom}
\frac{k}{2}\,\Theta_{\alpha}\bigl(\sqrt{k-1}\,\beta\bigr)
\;\asymp\;
k\,\bigl((k-1)\beta^{2}\bigr)^{\frac{d_{w}}{d_{w}-a}},
\qquad\text{as}\quad (k-1)\beta^{2}\to+\infty,
\end{equation}
which is exactly the growth of the lower bound~\eqref{eq:moment-lower-sub-geom}. Since the dissipative terms ($-k\lambda_{1}$ in~\eqref{eq:moment-upper-geom}) only grow linearly in $k$, we obtain, for a fixed $\beta>0$ and $k\to\infty$,
\begin{equation}\label{eq:lyapunov-two-sided-geom}
c\,k^{\frac{2d_{w}-a}{d_{w}-a}}
\;\le\;
\liminf_{t\to+\infty}\frac{\ln\E\bigl[u(t,x)^{k}\bigr]}{t}
\;\le\;
\limsup_{t\to+\infty}\frac{\ln\E\bigl[u(t,x)^{k}\bigr]}{t}
\;\le\;
C\,k^{\frac{2d_{w}-a}{d_{w}-a}}
\end{equation}
for $k$ large enough. The super-linear growth exponent $\frac{2d_{w}-a}{d_{w}-a}$ in~\eqref{eq:lyapunov-two-sided-geom} is the geometric analogue of the Euclidean exponent $\frac{4-a}{2-a}$ obtained in Theorems~\ref{thm:moment-upper} and~\ref{thm:moment-lower}. Indeed, as detailed in Remark~\ref{rem:moment-growth-geom}, in the Euclidean normalization $d_{w}=2$ our parameter~\eqref{eq:a-def} coincides with the exponent $a$ of~\eqref{eq:Lambda-power}, and $\frac{2d_{w}-a}{d_{w}-a}$ then reduces to $\frac{4-a}{2-a}$. The solution $u$ therefore exhibits the full intermittency phenomenon described in Remark~\ref{rmk:moment-estimates}. Similarly, for a fixed $k$ and $\beta\to+\infty$, the growth $\beta^{\frac{2d_{w}}{d_{w}-a}}$ is matched on both sides. 

\item {\it Critical regime: a logarithmic gap.}
Let $\alpha=\frac{d_{h}}{2d_{w}}$, that is $a=0$. In this case the asymptotics~\eqref{eq:theta-crit-large-beta-2} yield an upper bound of order $k(k-1)\beta^{2}\,(\ln((k-1)\beta^{2}))^{2}$ in~\eqref{eq:moment-upper-geom}, while the lower bound~\eqref{eq:moment-lower-crit-geom} only exhibits one logarithmic factor. The two bounds thus differ by a factor $\ln((k-1)\beta^{2})$, and we are not aware of a result capturing the exact critical rate.

\end{enumerate}
\end{remark}

\begin{remark}[Phase transition in $\beta$]\label{rem:phase-transition-geom}
The combination of Theorems~\ref{thm:moment-upper-geom} and~\ref{thm:moment-lower-gen-geom} also reveals a phase transition in the noise intensity $\beta$. Namely fix $k\ge2$ and $\alpha>0$. For small $\beta$, the asymptotics $\Theta_{\alpha}(\beta')\sim C\beta'^{2}$ as $\beta'\to0$ (valid in the three regimes, see Theorem~\ref{thm:existence-geom}) show that the right-hand side of~\eqref{eq:moment-upper-geom} is negative for $\beta$ small enough, so that $\E[u(t,x)^{k}]$ decays exponentially fast (the dissipation induced by the spectral gap $\lambda_{1}$ wins over the noise). On the other hand, the general lower bound~\eqref{eq:moment-lower-gen-geom} is positive as soon as $c_{1}(k-1)\beta^{2}>\lambda_{1}(B(x,\varepsilon_{1}))$, in which case the moments grow exponentially fast (the noise wins). This weak/strong disorder transition is not a feature of the geometric setting, but rather of the Dirichlet boundary condition on the bounded domain $U$. Indeed, on the full space $\R^{d}$ the lower bound~\eqref{eq:moment-lower} of Theorem~\ref{thm:moment-lower} entails exponential growth for every $\beta>0$. In contrast, a transition of the same nature was exhibited in~\cite{foondun:15:behaviour} for the stochastic heat equation on a Euclidean interval with Dirichlet boundary condition, driven by a space-time white noise (in that setting the transition also disappears when the Dirichlet condition is replaced by a Neumann condition). We refer to~\cite[Cor.~3.19]{baudoin:24:parabolic} for a detailed study of the transition in our geometric framework.
\end{remark}

\begin{remark}[Lower bounds under Neumann boundary conditions]\label{rem:neumann-lower-geom}
As mentioned in Remark~\ref{rem:phase-transition-geom}, the phase transition disappears when the Dirichlet boundary condition is replaced by a Neumann condition. Although we will not include the corresponding computations in this book, let us summarize the results of~\cite[Section~4.3]{baudoin:24:parabolic}. In that reference, the parabolic Anderson model is considered on a uniform domain $U$ (Definition~\ref{def:uniform}), for a noise whose spatial covariance is defined through the Neumann fractional Laplacian of Definition~\ref{def:frac-lap-N}. The moments of the solution admit a Feynman-Kac representation analogous to~\eqref{eq:fk-moment-formula-geom}, in which the killed Brownian motions are replaced by reflected Brownian motions (recall Theorem~\ref{thm:hk-neum}) and the kernel $G^{D,U}_{2\alpha}$ is replaced by its Neumann counterpart $G^{N,U}_{2\alpha}$ of Proposition~\ref{prop:RK-N-bounds}. Starting from this representation, for $\alpha>0$ and an initial condition satisfying $u_{0}\ge m_{0}>0$ on $U$, it is proved in~\cite[Thm.~4.6]{baudoin:24:parabolic} that
\begin{equation}\label{eq:moment-lower-neumann}
\liminf_{t\to+\infty}\frac{\ln\E\bigl[u(t,x)^{k}\bigr]}{t}
\;\ge\;
c\,\beta^{2}k(k-1),
\qquad\text{for every }\beta>0 .
\end{equation}
The mechanism behind~\eqref{eq:moment-lower-neumann} is that the Neumann semigroup is conservative, that is $P^{N,U}_{t}\mathbf{1}=\mathbf{1}$ (equivalently, the principal Neumann eigenvalue vanishes, see Proposition~\ref{prop:spectrum-N}). Hence no confinement of the Brownian paths is needed, and no eigenvalue term like $-k\lambda_{1}(B(x,\varepsilon_{1}))$ arises in~\eqref{eq:moment-lower-neumann}. In particular the moments grow exponentially fast for every $\beta>0$, which confirms that the transition of Remark~\ref{rem:phase-transition-geom} is specific to the Dirichlet case. Let us also mention that the analogues of Theorems~\ref{thm:moment-lower-sub-geom} and~\ref{thm:moment-lower-crit-geom} hold in the Neumann setting. Namely as $\beta\to+\infty$ the rate in~\eqref{eq:moment-lower-neumann} can be improved to $c\,k(k-1)\,\beta^{\frac{2d_{w}}{d_{w}-a}}$ when $a>0$, and to $c\,k(k-1)\,\beta^{2}\ln\beta$ when $a=0$, where $a$ is the parameter of~\eqref{eq:a-def}. Finally, in the white noise case $\alpha=0$ the second moment satisfies a lower bound with rate $\beta^{\frac{2d_{w}}{d_{w}-d_{h}}}$, which is consistent with the substitution $a=d_{h}$ in the subcritical exponent.
\end{remark}

\begin{remark}[The white noise case $\alpha=0$]\label{rem:moment-lower-white-geom}
Theorems~\ref{thm:moment-lower-gen-geom}, \ref{thm:moment-lower-sub-geom} and~\ref{thm:moment-lower-crit-geom} leave out the white noise case $\alpha=0$, since their proofs rest on the Feynman-Kac formula of Theorem~\ref{thm:fk-moments-geom} (recall from Remark~\ref{rem:fk-intersection-geom} that the case $\alpha=0$ would require intersection local times on $X$). A lower bound for the second moment can still be obtained by a renewal argument performed directly on the mild equation~\eqref{eq:mild-pam-D-geom}, in the spirit of the Laplace transform computations in the proof of Lemma~\ref{lem:fk-bound-geom}. One gets
\begin{equation}\label{eq:L2-lower-white-geom}
\liminf_{t\to+\infty}\frac{\ln\E\bigl[u(t,x)^{2}\bigr]}{t}
\;\ge\;
\Upsilon(\beta^{2})-2\lambda_{1},
\end{equation}
where $\Upsilon$ is the inverse of the function $y\mapsto C_{1}\,y\,e^{C_{2}y^{1/d_{w}}}$ and $C_{1},C_{2}$ are positive constants. We refer to~\cite[Thm.~3.17]{baudoin:24:parabolic} for the details.
\end{remark}

\section*{Bibliographic notes}
\addcontentsline{toc}{section}{Bibliographic notes}
\markright{BIBLIOGRAPHIC NOTES}

\smallskip

\noindent
{\it Noises for PAM.} The whole of Section~\ref{sec:geom-noise} follows Baudoin
and his coauthors. The idea is to give up on the translation-invariant kernel
$\Lambda(x-y)$ and to let the inverse fractional Laplacian
$(-\Delta_{U})^{-2\alpha}$ produce the covariance instead. This is summarized in
Definition~\ref{def:geom-noise-D}, through the Riesz kernel
$G^{D,U}_{2\alpha}$. The fields obtained this way were
constructed and studied by Baudoin and Lacaux~\cite{baudoin:22:fractional} in
the Neumann setting, on the Sierpi\'nski gasket first and then on related
fractals such as the Sierpi\'nski carpet. The  Dirichlet setting is handled by Baudoin and
Chen~\cite{baudoin:23:dirichlet}, on the gasket alone,
where they are also identified as limits of discrete fractional fields on the
canonical approximating graphs. Baudoin and Chen are the source of the
Sobolev-inclusion results of Subsection~\ref{ssec:sobolev-inclusion}. What makes
that subsection sharper than its counterpart on a general mm-space is the
precise spectral asymptotics available on the gasket. This is worked out in Fukushima and
Shima~\cite{fukushima:92:spectral}, and has no analogue at the level of
generality of Assumption~\ref{hyp:mmspace}.

Taking a negative power of the Laplacian as a covariance is not
specific to fractals. The resulting objects are the \emph{fractional Gaussian
fields}. The survey of Lodhia, Sheffield, Sun and
Watson~\cite{lodhia:16:fractional} covers the family as a whole, on $\R^{d}$,
on manifolds and on graphs. The case $\alpha=1$ is the Gaussian free field,
which is the central object of Liouville quantum gravity (see Duplantier and
Sheffield~\cite{duplantier:11:liouville}). A reader used to that literature
should note that the two uses of the field are quite different. There the
field is a fixed random distribution and the difficulty lies in exponentiating
it (through Gaussian multiplicative chaos). Here the same kind of
field is white in time and enters an evolution equation as a multiplicative
potential.

\smallskip

\noindent
{\it Existence and uniqueness.} Section~\ref{sec:geom-existence} follows the
six-author paper~\cite{baudoin:24:parabolic}, which is the reference for this
chapter as a whole. As explained at the start of that section, the chaos route
is preferred here to the contraction route of
Chapter~\ref{chap:existence}, because it delivers the $L^{2}$-moment formula
which Section~\ref{sec:geom-moments} then needs. The point at which the
geometric proof really differs from the Euclidean one of
Section~\ref{sec:chaos} is the absence of Fourier analysis. On $\R^{d}$
translation invariance lets the spectral measure diagonalize all the
heat-kernel convolutions at once, whereas on $U$ the only available tool is a
bilinear identity which handles one pairing at a time, and this is what forces
the level-by-level renewal recursion. Remark~\ref{rem:renewal-vs-fourier}
records the comparison.  For the stochastic heat
equation on a bounded Euclidean domain, where the boundary is present but the
geometry is not, we refer to Foondun and Nualart~\cite{foondun:15:behaviour}.

\smallskip

\noindent
{\it Moment estimates.} Section~\ref{sec:geom-moments} again
follows~\cite{baudoin:24:parabolic} for its results, but not for its proof of
the moment formula. We mollify the noise in space only and then run It\^{o}'s
formula, exactly as in Section~\ref{sec:fk-moments}, whereas that paper
regularizes in space and in time, as do Hu and Nualart~\cite{hu:09:stochastic}
and~\cite[Section~3.2]{hu:15:stochastic} on $\R^{d}$. In this book a space-time
regularization appears only in Remark~\ref{rem:fk-geom}, where it produces the
heuristic formula for a single realization of $u(t,x)$ rather than for its
moments. The eigenvalue bound on small balls at the heart of the lower bound is
taken from Mariano and Wang~\cite{mariano:22:spectral}, and the Gaussian
regularity estimates from the entropy method of Adler and
Taylor~\cite{adler:07:random}.

Two features of this section have no Euclidean counterpart, and both come from
the boundary rather than from the fractal structure. The first is the phase
transition in the noise intensity of Remark~\ref{rem:phase-transition-geom},
which is driven by the competition between the noise and the principal
Dirichlet eigenvalue $\lambda_{1}$. This competition disappears under Neumann conditions, where the semigroup is
conservative, and the corresponding statements are
in~\cite[Section~4.3]{baudoin:24:parabolic}. The second feature is that the white noise
case $\alpha=0$ is excluded from the Feynman-Kac formula, for the reason given
in Remark~\ref{rem:fk-intersection-geom}. It would require intersection local
times on $(X,d,\mu)$, which are not available at this level of generality.
Finally, the reader should be aware that the upper and lower bounds of
Remark~\ref{rem:moment-lower-comp-geom} do not match in the critical regime, one
carrying a $(\ln)^{2}$ where the other carries a single $\ln$. Closing that gap
is open.

\smallskip

\noindent
{\it The parabolic Anderson model in geometric contexts.} The equation has been
studied on geometric structures well before the setting of this chapter. For instance 
Tindel and Viens obtained space-time
regularity for the stochastic heat equation on a Lie
group~\cite{tindel:99:on}, and then the almost sure exponential behavior of a
parabolic SPDE on a manifold~\cite{tindel:02:almost}. Peszat and
Tindel~\cite{peszat:10:stochastic} give necessary and sufficient conditions for
the existence of a function-valued solution on a Lie group, in terms of the
covariance kernel of the noise. That is the direct ancestor of
Section~\ref{sec:geom-noise}, with the Lie group structure playing the role
that the Dirichlet form plays here. The discrete counterpart of the present
chapter, namely directed polymers on infinite graphs with a proper geometric structure, is developed by Cosco,
Seroussi and Zeitouni~\cite{cosco:21:directed}.

The line of work this chapter belongs to has several ramifications. Besides~\cite{baudoin:24:parabolic} on fractals, it covers the
Heisenberg group in both the It\^{o} and the Stratonovich
readings~\cite{baudoin:23:parabolic-b,baudoin:23:parabolic-a}. It should be
stressed that the Heisenberg case is \emph{not} covered by
Assumption~\ref{hyp:mmspace}, since the recurrence condition $d_{h}<d_{w}$
fails there. The same program is being carried out on manifolds where the
geometry enters more intrinsically. Baudoin, Chen and
Ouyang~\cite{baudoin:24:moment} treat Cartan-Hadamard manifolds, Chen and
Ouyang~\cite{chen:25:global} the role of the global geometry in the
well-posedness problem, Chen, Neel and Ouyang~\cite{chen:26:sharp} the heat
kernel on the cut locus, and Geng and Ouyang~\cite{geng:25:parabolic} the
hyperbolic spaces, where a phase transition again appears. Several of these
are still preprints at the time of writing.

Finally, a different toolkit reaches equations that the $L^{2}$ methods of this
book do not. When the noise is white in space and \emph{constant} in time, the
product $u\dot{W}$ is analytically ill-defined in dimension two, It\^{o}
calculus is of no help, and a renormalization is needed. Dahlqvist, Diehl and
Driver~\cite{dahlqvist:19:parabolic} prove well-posedness in that regime on a
two-dimensional closed Riemannian manifold, which requires extending the theory
of regularity structures to a curved space. Mayorcas and
Singh~\cite{mayorcas:25:singular} carry the same theory over to homogeneous Lie
groups, the Heisenberg group among them. Neither setting is the one of this
chapter, where the noise is white in time and the solution is built by
It\^{o} theory. They are worth quoting because they show what the geometric
program looks like once the equation becomes singular.

\backmatter
\cleardoublepage
\phantomsection
\addcontentsline{toc}{chapter}{Bibliography}
\bibliographystyle{plainnat}
\bibliography{pam-course-biblio}

\end{document}